%% file: higher_braided_01.tex
\documentclass{amsbook}

\usepackage[utf8]{inputenc}\usepackage{tgpagella}\usepackage{mathpazo}

\usepackage{graphicx}\usepackage{amsmath,amsthm,leftidx}
\usepackage[all,color]{xy}\usepackage[pagebackref]{hyperref} \usepackage{MnSymbol}
\usepackage{bbm}
\usepackage[mathscr]{eucal}\usepackage{xcolor}\usepackage{enumitem}

\numberwithin{section}{chapter}
\numberwithin{subsection}{section}
\numberwithin{equation}{section}
\theoremstyle{plain}
        \newtheorem{theorem}[equation]{Theorem}
        \newtheorem{proposition}[equation]{Proposition}
        \newtheorem{lemma}[equation]{Lemma}
        
        \newtheorem{corollary}[equation]{Corollary}

\theoremstyle{definition}
        \newtheorem{assumption}[equation]{Assumption}
        
        \newtheorem{definition}[equation]{Definition}
        
        \newtheorem{remark}[equation]{Remark}
        \newtheorem{example}[equation]{Example}
        \newtheorem{ginterpretation}[equation]{Geometric Interpretation}
        \newtheorem{interpretation}[equation]{Interpretation}
        \newtheorem{convention}[equation]{Convention}
        \newtheorem{notation}[equation]{Notation}
        \newtheorem{motivation}[equation]{Motivation}

\usepackage{tikz}
\usepackage{pgfplots}
\pgfplotsset{compat=1.13}
\usepackage{braids}

\usetikzlibrary{arrows,decorations.pathmorphing}
\usetikzlibrary{backgrounds,positioning}
\usetikzlibrary{fit,petri,shapes.misc}
\usetikzlibrary{shapes.geometric}
\usetikzlibrary{decorations.markings}
\usetikzlibrary{calc}

\newlength{\defaultpgflinewidth}
\tikzset{auto}

\tikzset{empty/.style={circle,inner sep=0pt,minimum size=6mm}}
\tikzset{emptyvt/.style={circle,inner sep=0pt,minimum size=0mm}}

\tikzset{plain/.style={circle,draw,very thick,
inner sep=0pt,minimum size=6mm}}

\tikzset{xplain/.style={circle,draw,very thick,
inner sep=0pt,minimum size=8mm}}

\tikzset{smallplain/.style={circle,draw,very thick,
inner sep=0pt,minimum size=4mm}}

\tikzset{xsplain/.style={circle,draw, thick,
inner sep=0pt,minimum size=1.5mm}}

\tikzset{tinyplain/.style={circle,draw, thick,
inner sep=0pt,minimum size=1mm}}

\tikzset{smalldotted/.style={circle,draw,very thick, densely dotted,
inner sep=0pt,minimum size=3mm}}

\tikzset{dottedplain/.style={circle,draw,very thick, densely dotted,
inner sep=0pt,minimum size=6mm}}

\tikzset{normaldot/.style={circle,black,fill=black,inner sep=0pt,minimum size=2mm}}

\tikzset{rectplain/.style={rectangle,draw,very thick,minimum size=6mm}}
\tikzset{bigplain/.style={rectangle,draw,very thick,minimum size=1cm}}

\tikzset{triangular/.style={regular polygon, regular polygon sides=3, draw,very thick,
inner sep=0pt,minimum size=1.2cm}}

\tikzset{arrow/.style={->,thick}}
\tikzset{dashedarrow/.style={->,dashed,thick}}
\tikzset{dottedarrow/.style={->,dotted,thick}}
\tikzset{mapto/.style={|->,thick}}
\tikzset{->-/.style={decoration={markings, mark=at position #1 with {\arrow{>}}},postaction={decorate}}}

\tikzset{implies/.style={thick,double,double equal sign distance,-implies}} 

\tikzset{line/.style={thick}}
\tikzset{dottedline/.style={dotted,thick}}
\tikzset{dashedline/.style={dashed,thick}}

\tikzset{inputleg/.style={<-,thick}}
\tikzset{outputleg/.style={->,thick}}
\tikzset{dottedinput/.style={<-,dotted,thick}}

\newcommand{\nicearrow}{\SelectTips{cm}{10}}
\newcommand{\nicexy}{\nicearrow\xymatrix@C+5pt}
\newcommand{\narrowxy}{\nicearrow\xymatrix}
\newcommand{\adjoint}{\nicearrow\xymatrix{ \ar@<2pt>[r] & \ar@<2pt>[l]}}
\newcommand{\becomes}{\nicexy{{}\ar@{~>}[r]&{}}}

\renewcommand{\to}{\hspace{-.1cm}\nicearrow\xymatrix@C-.3cm{\ar[r]&}\hspace{-.1cm}}

\newcommand{\shortto}{\hspace{-.1cm}\nicearrow\xymatrix@C-.6cm{\ar[r]&}\hspace{-.1cm}}
\newcommand{\iso}{\hspace{-.1cm}\nicearrow\xymatrix@C-.2cm{\ar[r]^-{\cong}&}\hspace{-.1cm}}
\newcommand{\weq}{\hspace{-.1cm}\nicearrow\xymatrix@C-.2cm{\ar[r]^-{\sim}&}\hspace{-.1cm}}
\renewcommand{\mapsto}{\hspace{-.1cm}\nicexy@C-.2cm{\ar@{|->}[r]&}\hspace{-.1cm}}

\newcommand{\defn}{\overset{\mathrm{def}}{=\joinrel=}}

\newcommand{\fieldc}{\mathbbm{C}}
\newcommand{\fieldk}{\mathbbm{K}}
\newcommand{\bbN}{\mathbbm{N}}
\newcommand{\fieldr}{\mathbbm{R}}
\newcommand{\Rn}{\fieldr^n}

\newcommand{\bbZ}{\mathbbm{Z}}

\newcommand{\cala}{\mathcal{A}}
\newcommand{\calc}{\mathcal{C}}
\newcommand{\cald}{\mathcal{D}}
\newcommand{\Dn}{\cald^n}
\newcommand{\Dtwo}{\cald^2}
\newcommand{\opendtwo}{\Dtwo_{\circ}}
\newcommand{\calf}{\mathcal{F}}
\newcommand{\cali}{\mathcal{I}}
\newcommand{\calj}{\mathcal{J}}
\newcommand{\call}{\mathcal{L}}
\newcommand{\calr}{\mathcal{R}}
\newcommand{\calv}{\mathcal{V}}
\newcommand{\calw}{\mathcal{W}}
\newcommand{\calx}{\mathcal{X}}

\newcommand{\colorc}{\mathfrak{C}}

\newcommand{\colord}{\mathfrak{D}}
\newcommand{\colorcd}{\colorc\times\colord}
\newcommand{\colore}{\mathfrak{E}}

\newcommand{\sequence}{\mathsf{seq}}
\newcommand{\seq}{\mathsf{Seq}}
\newcommand{\seqc}{\seq^{\colorc}}
\newcommand{\seqcm}{\seqc(\M)}
\newcommand{\seqcn}{\seqc(\N)}

\newcommand{\seqd}{\seq^{\colord}}
\newcommand{\seqdm}{\seqd(\M)}

\newcommand{\bsubc}{\B_{\colorc}}
\newcommand{\bsubcop}{\bsubc^{\op}}
\newcommand{\bsubcopc}{\bsubcop \times \colorc}
\newcommand{\bseqcm}{\B\seqcm}

\newcommand{\Cell}{\mathsf{Cell}}

\newcommand{\smallm}{\raisebox{.1cm}{\tiny{$\M$}}}
\newcommand{\subm}[1]{{#1}_{{\raisebox{.1cm}{\tiny{$\M$}}}}}
\newcommand{\ptsubm}{\P^{T}_{\smallm}}

\newcommand{\gsubm}{\subm{\G}}
\newcommand{\gdashsubm}{\gsubm^-}
\newcommand{\gsubc}{\G_{\colorc}}
\newcommand{\gsubcop}{\gsubc^{\op}}
\newcommand{\gsubcopc}{\gsubcop \times \colorc}
\newcommand{\gseqcm}{\G\seqcm}
\newcommand{\gseqcn}{\G\seqcn}

\newcommand{\gsubd}{\G_{\colord}}
\newcommand{\gsubdop}{\gsubd^{\op}}
\newcommand{\gsubdopd}{\gsubdop \times \colord}
\newcommand{\gseqdm}{\G\seqdm}

\newcommand{\gammacac}{\gamma^{\Cac}}
\newcommand{\gammapcac}{\gamma^{\PCac}}
\newcommand{\gammas}{\gamma^{\S}}
\newcommand{\gammag}{\gamma^{\G}}
\newcommand{\Gone}{\G^1}
\newcommand{\Gonedash}{\G^{1-}}
\newcommand{\gonesubc}{\Gone_{\colorc}}
\newcommand{\gonesubcop}{(\gonesubc)^{\op}}
\newcommand{\gonesubcopc}{\gonesubcop \times \colorc}
\newcommand{\goneseqcm}{\Gone\seqcm}
\newcommand{\goneop}{\Gone\Operad}
\newcommand{\goneopcset}{\Gone\Operad^{\colorc}(\Set)}
\newcommand{\goneopcm}{\Gone\Operadcm}
\newcommand{\goneopcn}{\Gone\Operadcn}
\newcommand{\goneopm}{\Gone\Operadm}
\newcommand{\goneopset}{\Gone\Operadset}

\newcommand{\Gtwo}{\G^2}
\newcommand{\Gtwodash}{\G^{2-}}
\newcommand{\gtwosubc}{\Gtwo_{\colorc}}
\newcommand{\gtwosubcop}{(\gtwosubc)^{\op}}
\newcommand{\gtwosubcopc}{\gtwosubcop \times \colorc}
\newcommand{\gtwoseqcm}{\Gtwo\seqcm}
\newcommand{\gtwoop}{\Gtwo\Operad}
\newcommand{\gtwoopcset}{\Gtwo\Operad^{\colorc}(\Set)}
\newcommand{\gtwoopcm}{\Gtwo\Operadcm}
\newcommand{\gtwoopcn}{\Gtwo\Operadcn}
\newcommand{\gtwoopm}{\Gtwo\Operadm}
\newcommand{\gtwoopset}{\Gtwo\Operadset}
\newcommand{\Govero}{\G^-_{\O}}

\newcommand{\pseq}{\mathsf{PSeq}}
\newcommand{\pseqc}{\pseq^{\colorc}}
\newcommand{\pseqcm}{\pseqc(\M)}

\newcommand{\rsubc}{\R_{\colorc}}
\newcommand{\rsubcop}{\rsubc^{\op}}
\newcommand{\rsubcopc}{\rsubcop \times \colorc}
\newcommand{\rseqcm}{\R\seqcm}

\newcommand{\cacsubc}{\Cac_{\colorc}}
\newcommand{\cacsubcop}{\cacsubc^{\op}}
\newcommand{\cacsubcopc}{\cacsubcop \times \colorc}
\newcommand{\cacseqcm}{\Cac\seqcm}

\newcommand{\ssubc}{\S_{\colorc}}
\newcommand{\ssubcop}{\ssubc^{\op}}
\newcommand{\ssubcopc}{\ssubcop \times \colorc}
\newcommand{\sseqcm}{\mathsf{S}\seqcm}

\newcommand{\Cor}{\mathrm{Cor}} 
\newcommand{\Coru}{\underline{\Cor}}
 
\newcommand{\Ed}{\mathsf{Ed}}

\newcommand{\Int}{\mathsf{Int}}
\newcommand{\Intl}{\mathsf{Intl}}
\newcommand{\SCir}{\mathsf{SCir}}
\newcommand{\profscir}{\Prof(\SCir)}
\newcommand{\profscirscir}{\profscir\times\SCir}

\newcommand{\Iso}{\mathsf{Iso}}

\newcommand{\pr}{\mathsf{Pr}}
\newcommand{\Gpr}{\G^{\pr}}
\newcommand{\Prof}{\mathsf{Prof}}

\newcommand{\Profc}{\Prof(\colorc)}
\newcommand{\Profcc}{\Profc \times \colorc}
\newcommand{\Profd}{\Prof(\colord)}
\newcommand{\Profdd}{\Profd \times \colord}
\newcommand{\Profcd}{\Prof(\colorc\times\colord)}
\newcommand{\Profcdcd}{\Profcd\times\colorc\times\colord}

\newcommand{\profofk}{\Prof(K)}

\newcommand{\profoft}{\Prof(T)}
\newcommand{\profedt}{\Prof(\Ed(T))}
\newcommand{\profedtedt}{\profedt\times\Ed(T)}
\newcommand{\profedv}{\Prof(\Ed(V))}
\newcommand{\profedvedv}{\profedv\times\Ed(V)}

\newcommand{\profm}{\Prof(M)}
\newcommand{\profmm}{\profm\times M}
\newcommand{\profofu}{\Prof(u)}
\newcommand{\profofv}{\Prof(v)}
\newcommand{\profofw}{\Prof(w)}

\newcommand{\Vt}{\mathsf{Vt}}
\newcommand{\BVt}{\underline{\Vt}}
\newcommand{\st}{\mathsf{st}}
\newcommand{\nor}{\mathsf{nor}}

\newcommand{\op}{\mathsf{op}}

\newcommand{\B}{\mathsf{B}}
\newcommand{\C}{\mathsf{C}}
\newcommand{\Cop}{\C^{\op}}
\newcommand{\Cdiag}{\C^{\mathsf{diag}}}
\newcommand{\Ctilde}{\widetilde{\C}}
\newcommand{\Cac}{\mathsf{Cac}}
\newcommand{\PCac}{\mathsf{PCac}}

\newcommand{\D}{\mathsf{D}}
\newcommand{\Df}{\D^{\mathsf{fr}}}

\newcommand{\Dtilde}{\widetilde{\D}}

\newcommand{\Dftildetwo}{\widetilde{\Df_2}}

\newcommand{\G}{\mathsf{G}}
\newcommand{\I}{\mathsf{I}}
\newcommand{\J}{\mathsf{J}}

\newcommand{\M}{\mathsf{M}}

\newcommand{\Mtoc}{\M^{\colorc}}
\newcommand{\N}{\mathsf{N}}

\newcommand{\Nerve}{\mathsf{Ner}}
\newcommand{\Nerveg}{\Nerve^{\G}}
\newcommand{\Nervegone}{\Nerve^{\Gone}}
\newcommand{\Nervegtwo}{\Nerve^{\Gtwo}}
\newcommand{\Nervep}{\Nerve^{\P}}
\newcommand{\Nerves}{\Nerve^{\S}}
\newcommand{\Nerver}{\Nerve^{\R}}
\newcommand{\Nerveb}{\Nerve^{\B}}
\newcommand{\Nervepr}{\Nerve^{\PR}}
\newcommand{\Nervepb}{\Nerve^{\PB}}
\newcommand{\Nervecac}{\Nerve^{\Cac}}
\newcommand{\Nervepcac}{\Nerve^{\PCac}}

\newcommand{\Pa}{\mathsf{Pa}}
\newcommand{\PaG}{\Pa\G}
\newcommand{\PaP}{\Pa\P}
\newcommand{\PaS}{\Pa\S}
\newcommand{\PaR}{\Pa\R}
\newcommand{\PaPR}{\Pa\PR}
\newcommand{\PaB}{\Pa\B}
\newcommand{\PaPB}{\Pa\PB}
\newcommand{\PaCac}{\Pa\Cac}
\newcommand{\PaPCac}{\Pa\PCac}

\renewcommand{\O}{\mathsf{O}}

\renewcommand{\P}{\mathsf{P}}
\newcommand{\PB}{\mathsf{PB}}
\newcommand{\PR}{\mathsf{PR}}
\newcommand{\Q}{\mathsf{Q}}
\newcommand{\R}{\mathsf{R}}

\newcommand{\Real}{\mathsf{Re}}
\newcommand{\Realg}{\Real^{\G}}
\newcommand{\Realgone}{\Real^{\Gone}}
\newcommand{\Realgtwo}{\Real^{\Gtwo}}

\newcommand{\Realp}{\Real^{\P}}
\newcommand{\Reals}{\Real^{\S}}
\newcommand{\Realr}{\Real^{\R}}
\newcommand{\Realb}{\Real^{\B}}
\newcommand{\Realpr}{\Real^{\PR}}
\newcommand{\Realpb}{\Real^{\PB}}
\newcommand{\Realcac}{\Real^{\Cac}}
\newcommand{\Realpcac}{\Real^{\PCac}}

\newcommand{\Realhc}{\mathsf{C}\Real}
\newcommand{\Realghc}{\mathsf{C}\Realg}
\newcommand{\Realphc}{\mathsf{C}\Realp}

\newcommand{\Realshc}{\mathsf{C}\Reals}
\newcommand{\Realghcone}{\Realhc^{\Gone}}
\newcommand{\Realghctwo}{\Realhc^{\Gtwo}}

\newcommand{\Nervehc}{\mathsf{C}\Nerve}
\newcommand{\Nerveghc}{\mathsf{C}\Nerveg}
\newcommand{\Nervephc}{\mathsf{C}\Nervep}

\newcommand{\Nerveshc}{\mathsf{C}\Nerves}
\newcommand{\Nerveghcone}{\Nervehc^{\Gone}}
\newcommand{\Nerveghctwo}{\Nervehc^{\Gtwo}}

\renewcommand{\S}{\mathsf{S}}
\newcommand{\SO}{\mathsf{SO}}
\newcommand{\T}{\mathsf{T}}
\newcommand{\U}{\mathsf{U}}
\newcommand{\V}{\mathsf{V}}
\newcommand{\W}{\mathsf{W}}

\newcommand{\Yoneda}{\mathsf{Yon}}
\newcommand{\Yonedag}{\Yoneda^{\G}}
\newcommand{\Yonedagone}{\Yoneda^{\Gone}}
\newcommand{\Yonedagtwo}{\Yoneda^{\Gtwo}}
\newcommand{\Yonedap}{\Yoneda^{\P}}
\newcommand{\Yonedas}{\Yoneda^{\S}}
\newcommand{\Yonedar}{\Yoneda^{\R}}
\newcommand{\Yonedab}{\Yoneda^{\B}}
\newcommand{\Yonedapr}{\Yoneda^{\PR}}
\newcommand{\Yonedapb}{\Yoneda^{\PB}}
\newcommand{\Yonedacac}{\Yoneda^{\Cac}}
\newcommand{\Yonedapcac}{\Yoneda^{\PCac}}

\newcommand{\bopcm}{\mathsf{BOp}^{\colorc}_{\littlesub{\M}}}
\newcommand{\pbopcm}{\mathsf{PBOp}^{\colorc}_{\littlesub{\M}}}
\newcommand{\bopset}{\mathsf{BOp}(\Set)}
\newcommand{\pbopset}{\mathsf{PBOp}(\Set)}

\newcommand{\ropcm}{\mathsf{ROp}^{\colorc}_{\littlesub{\M}}}
\newcommand{\propcm}{\mathsf{PROp}^{\colorc}_{\littlesub{\M}}}

\newcommand{\ropset}{\mathsf{ROp}(\Set)}
\newcommand{\propset}{\mathsf{PROp}(\Set)}

\newcommand{\cacopcm}{\mathsf{CacOp}^{\colorc}_{\littlesub{\M}}}
\newcommand{\pcacopcm}{\mathsf{PCacOp}^{\colorc}_{\littlesub{\M}}}
\newcommand{\cacopset}{\mathsf{CacOp}(\Set)}
\newcommand{\pcacopset}{\mathsf{PCacOp}(\Set)}

\newcommand{\Aut}{\mathsf{Aut}}
\newcommand{\Cl}{\mathsf{Cl}}

\newcommand{\colim}{\mathsf{colim}}
\newcommand{\End}{\mathsf{End}}

\newcommand{\Ex}{\mathsf{Ex}}

\newcommand{\face}[1]{\partial^{#1}}
\newcommand{\horn}[1]{\Lambda^{#1}}
\newcommand{\Image}{\mathsf{Im}}
\newcommand{\Segal}[1]{\mathsf{Seg}^{#1}}
\newcommand{\Segalg}{\Segal{\G}}

\newcommand{\Hom}{\mathsf{Hom}}
\newcommand{\Homg}{\Hom^{\G}}

\newcommand{\Homm}{\Hom_{\M}}
\newcommand{\interior}[1]{{#1}^{\circ}}

\newcommand{\Id}{\mathrm{Id}}
\newcommand{\id}{\mathrm{id}}

\newcommand{\Ob}{\mathsf{Ob}}

\newcommand{\Reg}{\mathsf{Reg}}

\newcommand{\operadunit}{\mathsf{1}}
\newcommand{\operadunittilde}{\widetilde{\mathsf{1}}}

\newcommand{\tensorunit}{\mathbbm{1}}
\newcommand{\tensor}[1]{\otimes^{#1}}

\newcommand{\comp}{\circ}
\newcommand{\compi}{\circ_i}

\newcommand{\compv}{\circ_v}

\newcommand{\circb}{\overset{\raisebox{-.05cm}{\tiny{$\B$}}}{\circ}}
\newcommand{\circg}{\overset{\raisebox{-.05cm}{\tiny{$\G$}}}{\circ}}
\newcommand{\circp}{\overset{\raisebox{-.05cm}{\tiny{$\P$}}}{\circ}}
\newcommand{\circr}{\overset{\raisebox{-.05cm}{\tiny{$\R$}}}{\circ}}
\newcommand{\circs}{\overset{\raisebox{-.05cm}{\tiny{$\S$}}}{\circ}}
\newcommand{\circcac}{\overset{\raisebox{-.05cm}{\tiny{$\Cac$}}}{\circ}}
\newcommand{\circgsubm}{\circg_{\raisebox{.05cm}{\tiny{$\M$}}}}
\newcommand{\circgsubn}{\circg_{\raisebox{.05cm}{\tiny{$\N$}}}}

\newcommand{\dom}{\mathsf{dom}}
\newcommand{\codom}{\mathsf{cod}}

\newcommand{\colimover}[1]{\underset{#1}{\mathsf{colim}}}
\newcommand{\limover}[1]{\underset{#1}{\mathsf{lim}}}
\newcommand{\coprodover}[1]{\underset{#1}{\coprod}}
\newcommand{\sqcupover}[1]{\underset{#1}{\sqcup}}
\newcommand{\prodover}[1]{\underset{#1}{\prod}}

\newcommand{\bigtensorover}[1]{\underset{#1}{\bigotimes}}
\newcommand{\tensorg}{\overset{\raisebox{-.05cm}{\tiny{$\G$}}}{\otimes}}
\newcommand{\tensorgone}{\overset{\raisebox{-.05cm}{\tiny{$\Gone$}}}{\otimes}}
\newcommand{\tensorgtwo}{\overset{\raisebox{-.05cm}{\tiny{$\Gtwo$}}}{\otimes}}
\newcommand{\tensorp}{\overset{\raisebox{-.05cm}{\tiny{$\P$}}}{\otimes}}
\newcommand{\tensors}{\overset{\raisebox{-.05cm}{\tiny{$\S$}}}{\otimes}}
\newcommand{\tensorb}{\overset{\raisebox{-.05cm}{\tiny{$\B$}}}{\otimes}}
\newcommand{\tensorpb}{\overset{\raisebox{-.05cm}{\tiny{$\PB$}}}{\otimes}}
\newcommand{\tensorr}{\overset{\raisebox{-.05cm}{\tiny{$\R$}}}{\otimes}}
\newcommand{\tensorpr}{\overset{\raisebox{-.05cm}{\tiny{$\PR$}}}{\otimes}}
\newcommand{\tensorcac}{\overset{\raisebox{-.05cm}{\tiny{$\Cac$}}}{\otimes}}
\newcommand{\tensorpcac}{\overset{\raisebox{-.05cm}{\tiny{$\PCac$}}}{\otimes}}

\newcommand{\tensorsub}[1]{\otimes_{\raisebox{.05cm}{\tiny{$#1$}}}}
\newcommand{\Homsub}[1]{\Hom_{\raisebox{.05cm}{\tiny{$#1$}}}}

\newcommand{\smallsub}[1]{{\raisebox{.05cm}{\tiny{$#1$}}}}
\newcommand{\smallsubg}{\smallsub{\G}}
\newcommand{\smallsubs}{\smallsub{\S}}
\newcommand{\smallsubb}{\smallsub{\B}}
\newcommand{\smallsubr}{\smallsub{\R}}
\newcommand{\smallsubcac}{\smallsub{\Cac}}

\newcommand{\timesover}[1]{\underset{#1}{\times}}

\newcommand{\dotover}[1]{\underset{#1}{\cdot}}

\newcommand{\dqed}{\hfill$\diamond$}
\newcommand{\inv}[1]{{#1}^{-1}}
\newcommand{\sigmainv}{\inv{\sigma}}
\newcommand{\sigmabar}{\overline{\sigma}}
\newcommand{\sigmavbar}{\overline{\sigma_v}}

\newcommand{\taubar}{\overline{\tau}}
\newcommand{\rhobar}{\overline{\rho}}
\newcommand{\chibar}{\overline{\chi}}
\newcommand{\uchi}{\underline{\chi}}
\newcommand{\utau}{\underline{\tau}}

\newcommand{\bbar}{\overline{b}}
\newcommand{\Bbar}{\overline{\B}}
\newcommand{\gbar}{\overline{g}}
\newcommand{\Gbar}{\overline{\G}}
\newcommand{\Gonebar}{\overline{\Gone}}
\newcommand{\Gtwobar}{\overline{\Gtwo}}
\newcommand{\pbar}{\overline{p}}
\newcommand{\qbar}{\overline{q}}
\newcommand{\rbar}{\overline{r}}
\newcommand{\xbar}{\overline{x}}

\newcommand{\gammatilde}{\widetilde{\gamma}}
\newcommand{\gammao}{\gamma^{\O}}
\newcommand{\gammaq}{\gamma^{\Q}}

\newcommand{\gammazerofr}{\gamma^0_{\mathsf{fr}}}

\newcommand{\rtilde}{\widetilde{r}}
\newcommand{\stilde}{\widetilde{s}}

\newcommand{\As}{\mathsf{As}}
\newcommand{\Asc}{\As^{\C}}

\newcommand{\Cat}{\mathsf{Cat}}
\newcommand{\Chain}{\mathsf{Chain}}
\newcommand{\Chaink}{\mathsf{Chain}_{\fieldk}}
\newcommand{\Chainz}{\mathsf{Chain}_{\bbZ}}
\newcommand{\Com}{\mathsf{Com}}
\newcommand{\Comc}{\Com^{\C}}
\newcommand{\Commc}{\Com(\M)^{\C}}

\renewcommand{\emptyset}{\varnothing}

\newcommand{\Gpd}{\mathsf{Gpd}}
\newcommand{\GrOp}{\mathsf{GrOp}}

\newcommand{\Mag}{\mathsf{Mag}}
\newcommand{\Magma}{\mathsf{Magma}}
\newcommand{\Mon}{\mathsf{Mon}}
\newcommand{\Monm}{\Mon(\M)}
\newcommand{\Monmc}{\Monm^{\C}}

\newcommand{\Boperad}{\mathsf{BOp}}
\newcommand{\Boperadc}{\Boperad^{\colorc}}

\newcommand{\Boperadcm}{\Boperadc(\M)}
\newcommand{\PBoperadcm}{\mathsf{PBOp}^{\colorc}(\M)}

\newcommand{\gopc}{\mathsf{\G\Operadc}}
\newcommand{\littlesub}[1]{{\raisebox{.01cm}{\tiny{$#1$}}}}
\newcommand{\gopcm}{\mathsf{GOp}^{\colorc}_{\littlesub{\M}}}

\newcommand{\gopcset}{\mathsf{GOp}^{\colorc}(\Set)}
\newcommand{\gopccat}{\mathsf{GOp}^{\colorc}(\Cat)}
\newcommand{\gopedt}{\mathsf{GOp}^{\Ed(T)}}
\newcommand{\gopedtset}{\mathsf{GOp}^{\Ed(T)}(\Set)}
\newcommand{\gopset}{\mathsf{GOp}(\Set)}

\newcommand{\goneopcd}{\Gone\Operad^{\colorcd}}
\newcommand{\goneopcdset}{\goneopcd(\Set)}
\newcommand{\gopm}{\mathsf{GOp}(\M)}

\newcommand{\goperad}{\mathsf{GOp}}
\newcommand{\goperadc}{\goperad^{\colorc}}
\newcommand{\goperadcm}{\goperadc(\M)}
\newcommand{\goneoperadcm}{\Gone\Operadc(\M)}
\newcommand{\gtwooperadcm}{\Gtwo\Operadc(\M)}
\newcommand{\goneoperaddm}{\Gone\Operadd(\M)}
\newcommand{\gtwooperaddm}{\Gtwo\Operadd(\M)}

\newcommand{\goperadcn}{\goperadc(\N)}
\newcommand{\goperadcset}{\goperadc(\Set)}
\newcommand{\goperaddashm}{\goperad^{-}(\M)}

\newcommand{\goperadd}{\goperad^{\colord}}
\newcommand{\goperaddm}{\goperadd(\M)}
\newcommand{\goperaddset}{\goperadd(\Set)}

\newcommand{\gopcd}{\mathsf{GOp}^{\colorc\times\colord}}
\newcommand{\goperadcd}{\gopcd}
\newcommand{\goperadcdset}{\goperadcd(\Set)}

\newcommand{\iotaof}[1]{\iota^{\raisebox{-.05cm}{\tiny{$#1$}}}}
\newcommand{\iotag}{\iotaof{\G}} 
\newcommand{\Lambdaof}[1]{\Lambda^{\raisebox{-.05cm}{\tiny{$#1$}}}}
\newcommand{\Lambdag}{\Lambdaof{\G}} \newcommand{\Lambdap}{\Lambdaof{\P}}
\newcommand{\Lambdab}{\Lambdaof{\B}}
\newcommand{\Lambdas}{\Lambdaof{\S}}
\newcommand{\Lambdar}{\Lambdaof{\R}}
\newcommand{\Lambdacac}{\Lambdaof{\Cac}}

\newcommand{\MCat}{\mathsf{MCat}}
\newcommand{\MCatst}{\MCat_{\st}}
\newcommand{\MCatof}[1]{\mathsf{MCat}^{\raisebox{-.05cm}{\tiny{$#1$}}}}
\newcommand{\MCatg}{\MCatof{\G}}
\newcommand{\MCatgst}{\MCatg_{\mathsf{st}}}
 
\newcommand{\MCatgone}{\MCatof{\Gone}}
\newcommand{\MCatgtwo}{\MCatof{\Gtwo}}
\newcommand{\MCatgonest}{\MCatgone_{\mathsf{st}}}
\newcommand{\MCatgtwost}{\MCatgtwo_{\mathsf{st}}}
\newcommand{\MCatp}{\MCatof{\P}} 
\newcommand{\MCatpst}{\MCatp_{\mathsf{st}}}
\newcommand{\MCats}{\MCatof{\S}}
\newcommand{\MCatsst}{\MCats_{\mathsf{st}}}
\newcommand{\MCatb}{\MCatof{\B}}
\newcommand{\MCatbst}{\MCatb_{\mathsf{st}}}
\newcommand{\MCatr}{\MCatof{\R}}
\newcommand{\MCatrst}{\MCatr_{\mathsf{st}}}
\newcommand{\MCatcac}{\MCatof{\Cac}}
\newcommand{\MCatcacst}{\MCatcac_{\mathsf{st}}}

\newcommand{\MFun}{\mathsf{MFun}}
\newcommand{\MFunof}[1]{\mathsf{MFun}^{\raisebox{-.05cm}{\tiny{$#1$}}}}
\newcommand{\MFung}{\MFunof{\G}}

\newcommand{\MFuns}{\MFunof{\S}}
\newcommand{\MFunb}{\MFunof{\B}}
\newcommand{\MFunr}{\MFunof{\R}}
\newcommand{\MFuncac}{\MFunof{\Cac}}
\newcommand{\MFunstg}{\MFun_{\mathsf{stg}}}
\newcommand{\MFunstgof}[1]{\MFunstg^{\raisebox{-.05cm}{\tiny{$#1$}}}}
\newcommand{\MFungstg}{\MFung_{\mathsf{stg}}}
\newcommand{\MFunsstg}{\MFuns_{\mathsf{stg}}}
\newcommand{\MFunbstg}{\MFunb_{\mathsf{stg}}}
\newcommand{\MFunrstg}{\MFunr_{\mathsf{stg}}}
\newcommand{\MFuncacstg}{\MFuncac_{\mathsf{stg}}}

\newcommand{\agop}{\mathsf{Op}^{\mathsf{act}}}
\newcommand{\Operad}{\mathsf{Op}}
\newcommand{\Operadc}{\Operad^{\colorc}}
\newcommand{\Operadd}{\Operad^{\colord}}

\newcommand{\Operadset}{\Operad(\Set)}
\newcommand{\Operadcm}{\Operadc(\M)}
\newcommand{\Operadcn}{\Operadc(\N)}
\newcommand{\Operadm}{\Operad(\M)}

\newcommand{\popcm}{\mathsf{POp}^{\colorc}_{\littlesub{\M}}}
\newcommand{\popset}{\mathsf{POp}(\Set)}

\newcommand{\Poperad}{\mathsf{POp}}
\newcommand{\Poperadc}{\Poperad^{\colorc}}

\newcommand{\Poperadcm}{\Poperadc(\M)}
\newcommand{\Poperadm}{\Poperad(\M)}

\newcommand{\Roperad}{\mathsf{ROp}}
\newcommand{\Roperadc}{\Roperad^{\colorc}}

\newcommand{\Roperadcm}{\Roperadc(\M)}
\newcommand{\PRoperadcm}{\mathsf{PROp}^{\colorc}(\M)}

\newcommand{\Cacoperad}{\mathsf{CacOp}}
\newcommand{\Cacoperadc}{\Cacoperad^{\colorc}}
\newcommand{\Cacoperadcm}{\Cacoperadc(\M)}
\newcommand{\Cacoperadm}{\Cacoperad(\M)}
\newcommand{\Pcacoperadcm}{\mathsf{PCacOp}^{\colorc}(\M)}
\newcommand{\Pcacoperadm}{\mathsf{PCacOp}(\M)}

\newcommand{\sopcm}{\mathsf{SOp}^{\colorc}_{\littlesub{\M}}}
\newcommand{\sopset}{\mathsf{SOp}(\Set)}

\newcommand{\Soperad}{\mathsf{SOp}}
\newcommand{\Soperadc}{\Soperad^{\colorc}}

\newcommand{\Soperadcm}{\Soperadc(\M)}
\newcommand{\Soperadm}{\Soperad(\M)}

\newcommand{\Set}{\mathsf{Set}}
\newcommand{\Sset}{\mathsf{SSet}}

\newcommand{\CHau}{\mathsf{CHau}}
\newcommand{\Top}{\mathsf{Top}}

\newcommand{\Tree}{\mathsf{Tree}}
\newcommand{\Treec}{\Tree^{\colorc}}
\newcommand{\GTree}{\G\Tree}
\newcommand{\GTreec}{\G\Treec}

\newcommand{\GTreecduc}{\GTreec\duc}
\newcommand{\GTreecducop}{\GTreecduc^{\op}}
\newcommand{\Gonetreec}{\Gone\Treec}
\newcommand{\Gtwotreec}{\Gtwo\Treec}
\newcommand{\PTree}{\mathsf{P}\Tree}
\newcommand{\PTreec}{\mathsf{P}\Treec}
\newcommand{\RTreec}{\R\Treec}
\newcommand{\BTreec}{\B\Treec}
\newcommand{\PRTreec}{\PR\Treec}
\newcommand{\PBTreec}{\PB\Treec}
\newcommand{\STreec}{\mathsf{S}\Treec}
\newcommand{\CacTreec}{\Cac\Treec}
\newcommand{\PCacTreec}{\PCac\Treec}

\newcommand{\PTan}{\mathsf{PTan}}

\newcommand{\Treecat}{\Psi}
\newcommand{\Treecatg}{\Treecat^{\G}}
\newcommand{\Treecatgplus}{\Treecatg_+}
\newcommand{\Treecatgminus}{\Treecatg_-}
\newcommand{\Treecatgop}{(\Treecatg)^{\op}}
\newcommand{\Settotreecatgop}{\Set^{\Treecatgop}}
\newcommand{\Treecatgsub}{\overline{\Treecat}^{\G}}
\newcommand{\Treecatp}{\Treecat^{\P}}
\newcommand{\Treecatr}{\Treecat^{\R}}
\newcommand{\Treecatb}{\Treecat^{\B}}
\newcommand{\Treecatpr}{\Treecat^{\PR}}
\newcommand{\Treecatpb}{\Treecat^{\PB}}
\newcommand{\Treecats}{\Treecat^{\S}}
\newcommand{\Treecatcac}{\Treecat^{\Cac}}
\newcommand{\Treecatpcac}{\Treecat^{\PCac}}

\newcommand{\Treecatpop}{(\Treecatp)^{\op}}
\newcommand{\Treecatrop}{(\Treecatr)^{\op}}
\newcommand{\Treecatbop}{(\Treecatb)^{\op}}
\newcommand{\Treecatprop}{(\Treecatpr)^{\op}}
\newcommand{\Treecatpbop}{(\Treecatpb)^{\op}}
\newcommand{\Treecatsop}{(\Treecats)^{\op}}
\newcommand{\Treecatcacop}{(\Treecatcac)^{\op}}
\newcommand{\Treecatpcacop}{(\Treecatpcac)^{\op}}

\newcommand{\Settotreecatpop}{\Set^{\Treecatpop}}
\newcommand{\Settotreecatrop}{\Set^{\Treecatrop}}
\newcommand{\Settotreecatbop}{\Set^{\Treecatbop}}
\newcommand{\Settotreecatprop}{\Set^{\Treecatprop}}
\newcommand{\Settotreecatpbop}{\Set^{\Treecatpbop}}
\newcommand{\Settotreecatsop}{\Set^{\Treecatsop}}
\newcommand{\Settotreecatcacop}{\Set^{\Treecatcacop}}
\newcommand{\Settotreecatpcacop}{\Set^{\Treecatpcacop}}

\newcommand{\Treecatgone}{\Treecat^{\Gone}}
\newcommand{\Treecatgtwo}{\Treecat^{\Gtwo}}
\newcommand{\Treecatgoneop}{(\Treecatgone)^{\op}}
\newcommand{\Treecatgtwoop}{(\Treecatgtwo)^{\op}}
\newcommand{\Settotreecatgoneop}{\Set^{\Treecatgoneop}}
\newcommand{\Settotreecatgtwoop}{\Set^{\Treecatgtwoop}}

\newcommand{\Treecatu}{\underline{\Treecat}}
\newcommand{\Treecatug}{\Treecatu^{\G}}

\newcommand{\phihat}{\hat{\phi}}
\newcommand{\varphihat}{\hat{\varphi}}
\newcommand{\Khat}{\widehat{K}}
\newcommand{\That}{\widehat{T}}
\newcommand{\Uhat}{\widehat{U}}
\newcommand{\Vhat}{\widehat{V}}
\newcommand{\What}{\widehat{W}}
\newcommand{\Xhat}{\widehat{X}}

\newcommand{\Kunderline}{\underline{K}}
\newcommand{\Lunderline}{\underline{L}}
\newcommand{\Tunderline}{\underline{T}}
\newcommand{\Uunderline}{\underline{U}}
\newcommand{\Vunderline}{\underline{V}}
\newcommand{\Wunderline}{\underline{W}}
\newcommand{\Xunderline}{\underline{X}}

\newcommand{\Wg}{\W^{\G}}
\newcommand{\Wgone}{\W^{\Gone}}
\newcommand{\Wgtwo}{\W^{\Gtwo}}
\newcommand{\Wp}{\W^{\P}}
\newcommand{\Wr}{\W^{\R}}
\newcommand{\Wb}{\W^{\B}}
\newcommand{\Wpr}{\W^{\PR}}
\newcommand{\Wpb}{\W^{\PB}}
\newcommand{\Ws}{\W^{\S}}
\newcommand{\Wcac}{\W^{\Cac}}
\newcommand{\Wpcac}{\W^{\PCac}}
\newcommand{\Wgo}{\Wg\O}
\newcommand{\Wgq}{\Wg\Q}

\newcommand{\alg}{\mathsf{Alg}}
\newcommand{\algm}{\alg_{\littlesub{\M}}}

\newcommand{\algmo}{\algm(\O)}

\newcommand{\algmq}{\algm(\Q)}

\newcommand{\ua}{\underline a}
\newcommand{\ub}{\underline b}
\newcommand{\uc}{\underline c}

\newcommand{\ud}{\underline d}
\newcommand{\udfr}{\ud_{\mathsf{fr}}}
\newcommand{\udkfr}{\udfr^k}
\newcommand{\udnfr}{\udfr^n}

\newcommand{\ue}{\underline e}
\newcommand{\uf}{\underline f}
\newcommand{\ug}{\underline g}
\newcommand{\uh}{\underline h}
\newcommand{\uk}{\underline k}
\newcommand{\ul}{\underline l}
\newcommand{\um}{\underline m}
\newcommand{\ur}{\underline r}
\newcommand{\us}{\underline s}

\newcommand{\uz}{\underline z}
\newcommand{\uvarphi}{\underline \varphi}
\newcommand{\uphi}{\underline \phi}

\newcommand{\smallprof}[1]
{\raisebox{.05cm}{\scalebox{0.8}{#1}}}

\newcommand{\sbinom}[2]{\raisebox{.05cm}{\scalebox{0.8}{$\binom{#1}{#2}$}}}
\newcommand{\inout}[1]{\raisebox{.05cm}{\scalebox{0.8}{$\binom{\out(#1)}{\inp(#1)}$}}}
\newcommand{\inoutt}{\inout{t}}
\newcommand{\inoutphit}{\sbinom{\phi\,\out(t)}{\phi\,\inp(t)}}

\newcommand{\Binoutt}{\sbinom{B\,\out(t)}{B\,\inp(t)}}

\newcommand{\bjiuaji}{\smallprof{$\binom{b^j_i}{\ua^j_i}$}}

\newcommand{\ciuatwo}{\smallprof{$\binom{c_i}{\ua^2}$}}
\newcommand{\cjuaj}{\smallprof{$\binom{c_j}{\ua_j}$}}

\newcommand{\ciub}{\smallprof{$\binom{c_i}{\ub}$}}

\newcommand{\cjub}{\smallprof{$\binom{c_j}{\ub}$}}
\newcommand{\cjubj}{\smallprof{$\binom{c_j}{\ub_j}$}}
\newcommand{\fcjfubj}{\smallprof{$\binom{fc_j}{f\ub_j}$}}
\newcommand{\cjubjtauj}{\smallprof{$\binom{c_j}{\ub_j\tau_j}$}}
\newcommand{\cjubjtaubarj}{\smallprof{$\binom{c_j}{\ub_j\taubar_j}$}}
\newcommand{\csigmajubsigmaj}{\smallprof{$\binom{c_{\sigma(j)}}{\ub_{\sigma(j)}}$}}
\newcommand{\csigmabarjubsigmabarj}{\smallprof{$\binom{c_{\sigmabar(j)}}{\ub_{\sigmabar(j)}}$}}
\newcommand{\cc}{\smallprof{$\binom{c}{c}$}}
\newcommand{\ccsingle}{\smallprof{$\binom{c}{c}$}}
\newcommand{\cici}{\smallprof{$\binom{c_i}{c_i}$}}
\newcommand{\cjcj}{\smallprof{$\binom{c_j}{c_j}$}}

\newcommand{\cjudj}{\smallprof{$\binom{c_j}{\ud_j}$}}

\newcommand{\czerouc}{\smallprof{$\binom{c_0}{\uc}$}}
\newcommand{\czeroempty}{\smallprof{$\binom{c_0}{\varnothing}$}}

\newcommand{\dua}{\smallprof{$\binom{d}{\ua}$}}
\newcommand{\duaone}{\smallprof{$\binom{d}{\ua^1}$}}
\newcommand{\duaonecompsigmabariuatwo}{\smallprof{$\binom{d}{\ua^1 \comp_{\sigmabar(i)}\ua^2}$}}

\newcommand{\dub}{\smallprof{$\binom{d}{\ub}$}}

\newcommand{\duboneubn}{\smallprof{$\binom{d}{\ub_1,\ldots,\ub_n}$}}
\newcommand{\dubsigmaoneubsigman}{\smallprof{$\binom{d}{\ub_{\sigma(1)},\ldots,\ub_{\sigma(n)}}$}}
\newcommand{\dubsigmabaroneubsigmabarn}{\smallprof{$\binom{d}{\ub_{\sigmabar(1)},\ldots,\ub_{\sigmabar(n)}}$}}
\newcommand{\dubonetauoneubntaun}{\smallprof{$\binom{d}{\ub_1\tau_1,\ldots,\ub_n\tau_n}$}}
\newcommand{\dubonetaubaroneubntaubarn}{\smallprof{$\binom{d}{\ub_1\taubar_1,\ldots,\ub_n\taubar_n}$}}
\newcommand{\dc}{\smallprof{$\binom{d}{c}$}}

\newcommand{\dd}{\smallprof{$\binom{d}{d}$}}

\newcommand{\duc}{\smallprof{$\binom{d}{\uc}$}}
\newcommand{\dsigmauc}{\smallprof{$\binom{d}{\sigma\uc}$}}
\newcommand{\dsigmabaruc}{\smallprof{$\binom{d}{\sigmabar\uc}$}}
\newcommand{\duccompiub}{\smallprof{$\binom{d}{\uc\compi\ub}$}}
\newcommand{\dconecn}{\smallprof{$\binom{d}{c_1,\ldots,c_n}$}}
\newcommand{\ducsigma}{\smallprof{$\binom{d}{\uc\sigma}$}}
\newcommand{\ducsigmabar}{\smallprof{$\binom{d}{\uc\sigmabar}$}}

\newcommand{\dzeroud}{\smallprof{$\binom{d_0}{\ud}$}}
\newcommand{\dzeroempty}{\smallprof{$\binom{d_0}{\varnothing}$}}

\newcommand{\fdfub}{\smallprof{$\binom{fd}{f\ub}$}}
\newcommand{\fdufc}{\smallprof{$\binom{fd}{f\uc}$}}
\newcommand{\fdfuc}{\fdufc}

\newcommand{\fzerouf}{\smallprof{$\binom{f_0}{\uf}$}}
\newcommand{\fzeroufsigmabar}{\smallprof{$\binom{f_0}{\uf\sigmabar}$}}
\newcommand{\fzerocufc}{\smallprof{$\binom{f_0c}{\uf c}$}}

\newcommand{\siur}{\smallprof{$\binom{s_i}{\ur}$}}
\newcommand{\tus}{\smallprof{$\binom{t}{\us}$}}
\newcommand{\tustau}{\smallprof{$\binom{t}{\us\tau}$}}
\newcommand{\tsonesn}{\smallprof{$\binom{t}{s_1,\ldots,s_n}$}}
\newcommand{\ttsingle}{\smallprof{$\binom{t}{t}$}}

\newcommand{\inp}{\mathsf{in}}
\newcommand{\out}{\mathsf{out}}

\newcommand{\Phyl}{\mathsf{Phyl}}

\newcommand{\andspace}{\quad\text{and}\quad}

\newcommand{\ifspace}{\quad\text{if}\quad}

\newcommand{\forspace}{\quad\text{for}\quad}

\newcommand{\orspace}{\quad\text{or}\quad}

\newcommand{\withspace}{\quad\text{with}\quad}

\newcommand{\minusquad}{\mkern-18mu}

\makeindex

\begin{document}
\title{Infinity Operads and Monoidal Categories with Group Equivariance}
\author{Donald Yau}
\address{Department of Mathematics\\
	 The Ohio State University at Newark\\
	 1179 University Drive\\ 
	 Newark, OH 43055, USA}
\email{yau.22@osu.edu}


\begin{abstract}
This monograph provides a coherent development of operads, infinity operads, and monoidal categories, equipped with equivariant structures encoded by an action operad.  A group operad is a planar operad with an action operad equivariant structure.  In the first three parts of this monograph, we establish a foundation for group operads and for their higher coherent analogues called infinity group operads.  Examples include planar, symmetric, braided, ribbon, and cactus operads, and their infinity analogues.  For example, with the tools developed here, we observe that the coherent ribbon nerve of the universal cover of the framed little $2$-disc operad is an infinity ribbon operad.

In Part \ref{part:monoidal-cat} we define general monoidal categories equipped with an action operad equivariant structure, and provide a unifying treatment of coherence and strictification for them.  Examples of such monoidal categories include symmetric, braided, ribbon, and coboundary monoidal categories, which naturally arise in the representation theory of quantum groups and of coboundary Hopf algebras and in the theory of crystals of finite dimensional complex reductive Lie algebras.  

This is the first monograph that provides a detailed treatment of operads, infinity operads, and monoidal categories, with action operad equivariance.  Many illustrations and examples are included.  Assuming only basic category theory, this monograph is intended for graduate students and researchers.   In addition to being a coherent reference for the topics covered, this book is also suitable for a graduate student seminar and a reading course.
\end{abstract}

\subjclass[2000]{18A30, 18A40, 18D10, 18D50, 18G55, 20F36, 55P48, 81R50, 81T05}
\keywords{Action operads, group operads, symmetric operads, braided operads, ribbon operads, cactus operads, infinity operads, Boardman-Vogt construction, coherent nerve, monoidal categories, coherence, strictification.}

\date{\today}
\maketitle

\cleardoublepage
\thispagestyle{empty}
\vspace*{13.5pc}
\begin{center}
To Eun Soo and Jacqueline
\end{center}
\cleardoublepage

\setcounter{page}{6}
\tableofcontents


\mainmatter

\include{operads_with_group_equivariance}

\include{constructions_group_operads}

\include{infinity_group_operads}

\include{coherence_for_monoidal_categories}

\appendix
\backmatter


\chapter*{List of Notations}

\newcommand{\where}[1]{\> \> \pageref{#1} \> \>}
\newcommand{\blob}{\> \> \> \> \hspace{2em}}

\begin{tabbing}
\phantom{\textbf{Notation}} \= \hspace{1.5cm}\= \phantom{\textbf{Page}}\= \hspace{.5cm}\= \phantom{\textbf{Description}} \\
\textbf{Chapter \ref{ch:introduction}} \>\> \textbf{Page} \>\> \textbf{Description}\\
$|S|$ \where{notation:cardinality} cardinality of a set $S$\\
$\Ob(\C)$ \where{notation:object} class of objects in a category $\C$\\
$\C(X,Y)$, $\C(X;Y)$ \where{notation:morphism-set}  set of morphisms from $X$ to $Y$ in $\C$\\
$\Id_X$ \where{notation:identity-morphism} identity morphism of $X$\\
$gf$ or $g\comp f$ \where{notation:morphism-composition} composition of morphisms\\
$\C^{\op}$ \where{notation:opposite-category} opposite category of $\C$\\
$\C^{\D}$ \where{notation:diagram-category} diagram category of functors $\D \to \C$\\
$\colim\, F$, $\colimover{x\in \D}\, Fx$ \where{notation:colimit} colimit of $F : \D \to C$\\
$\varnothing$ or $\varnothing^{\C}$ \where{notation:initial-object} initial object\\
$\int^{c\in\C} F(c,c)$ \where{notation:coend} coend of $F : \Cop \times \C \to \D$\\
$\Set$ \where{notation:set} category of sets\\
$\Chaink$ \where{notation:chaink} category of chain complexes over $\fieldk$\\
$\CHau$ \where{notation:chau} category of c.g. weak Hausdorff spaces\\
$\Sset$ \where{notation:sset} category of simplicial sets\\
$\Cat$ \where{notation:cat} category of small categories\\
$x^{\otimes n}$ \where{notation:tensors} $x \otimes \cdots\otimes x$ with $n$ copies of $x$\\
$x^{\otimes 0}$ \where{notation:tensors} monoidal unit $\tensorunit$\\ 
$\bigotimes_{i=1}^n x_i$ \where{notation:tensors} $x_1 \otimes \cdots \otimes x_n$\\ 
$(X,\mu,\operadunit)$ \where{notation:monoid} a monoid\\
$\Mon(\C,\otimes,\tensorunit)$ \where{notation:monoid-category} category of monoids\\
\blob\\
\textbf{Chapter \ref{ch:planar-operad}}\blob \\
$\colorc$ \where{notation:color-set} a set of colors\\
$\uc$, $(c_1,\ldots,c_n)$ \where{notation:uc} a $\colorc$-profile\\
$|\uc|$ \where{notation:lengthuc} length of $\uc$\\
$\varnothing$ \where{notation:empty-prof} empty $\colorc$-profile\\
$\Profc$ \where{notation:profc} set of $\colorc$-profiles\\
$\pseqcm$ \where{notation:pseqcm} $\colorc$-colored planar sequences in $\M$\\
$X\duc$ \where{notation:xduc} an entry of $X \in \pseqcm$\\
$(\uc;d)$, $\duc$ \where{notation:duc} an element in $\Profcc$\\
$\M^{\colorc}$ \where{notation:mtoc} category of $\colorc$-colored objects in $\M$\\
$X_c$ \where{notation:xsubc} an entry of a $\colorc$-colored object $X$\\
$\circp$ \where{planar-circle-product} planar circle product\\
$\I$ \where{unit-operad} unit $\colorc$-colored planar sequence\\
$\Poperadcm$ \where{def:planar-operad} $\colorc$-colored planar operads in $\M$\\
$\O(n)$ \where{notation:oofn} $n$th entry of a $1$-colored planar operad\\
$\gamma$ \where{operadic-composition} operadic composition\\
$\operadunit_c$ \where{notation:colored-unit} $c$-colored unit\\
$\uc \compi \ub$ \where{def:compi} replace $c_i$ in $\uc$ with $\ub$\\
$\compi$ \where{operadic-compi} $\compi$-composition\\
$X_{\uc}$ \where{not:x-sub-c} $X_{c_1}\otimes \cdots \otimes X_{c_m}$\\
$\alg(T)$ \where{notation:algt} category of $T$-algebras\\
$\algmo$ \where{def:planar-operad-algebra} category of $\O$-algebras in $\M$\\
$\As$ \where{ex:as-planar-operad} associative operad\\
$\Cdiag$ \where{ex:diag-operad} planar operad for $\C$-diagrams\\
$\Mag$ \where{notation:mag} planar operad for magmas\\
$\Magma$ \where{notation:magma} category of magmas\\
\blob\\
\textbf{Chapter \ref{ch:symmetric-operad}}\blob \\
$\ssubc$ \where{notation:sigmac} groupoid of $\colorc$-profiles with left permutations\\
$\ssubcop$ \where{notation:sigmac} groupoid of $\colorc$-profiles with right permutations\\
$\sseqcm$ \where{notation:symseqcm} $\colorc$-colored symmetric sequences in $\M$\\
$\circs$ \where{symmetric-circle-product} symmetric circle product\\
$\Soperadcm$ \where{def:symmetric-operad} $\colorc$-colored symmetric operads in $\M$\\
$\sigma\langle k_1,\ldots,k_n\rangle$ \where{notation:block-permutation} block permutation induced by $\sigma$\\
$\tau_1\oplus\cdots\oplus\tau_n$ \where{notation:direct-sum} direct sum permutation\\
$\S$ \where{sn-operad} symmetric group operad\\
$\gamma^{\S}$ \where{sn-operad} operadic composition in $\S$\\
$\End(X)$ \where{ex:endomorphism-operad} endomorphism operad\\
$\Com$ \where{ex:operad-com} commutative operad\\
$\Asc$ \where{ex:diag-monoid-operad} operad for $\C$-diagrams of monoids\\
$\Comc$ \where{notation:comc} operad for $\C$-diagrams of commutative monoids\\
$\cali$ \where{notation:square} closed unit interval $[0,1]$\\
$\calj$ \where{notation:square} open unit interval $(0,1)$\\
$\C_2$ \where{notation:ctwo} little $2$-cube operad\\
$\C_n$ \where{ex:little-ncube} little $n$-cube operad\\
$\Dn$ \where{ex:little-n-disc} closed unit $n$-disc\\
$\D_n$ \where{notation:dn} little $n$-disc operad\\
$\SO(n)$ \where{notation:son} special orthogonal group on $\fieldr^n$\\
$\Df_n$ \where{notation:dfn} framed little $n$-disc operad\\
$\GrOp$ \where{notation:grop} operad of graph operations\\
$\inp(T)$, $\out(T)$ \where{def:tree-for-phyl} set of inputs and output of a tree $T$\\
$\Ed(T)$, $\Int(T)$  \where{notation:edge} sets of edges and internal edges in $T$\\
$\inp(v)$, $\out(v)$ \where{notation:inpe} set of inputs and output of a vertex $v$\\
$\Vt(T)$ \where{notation:inpe} set of vertices in $T$\\
$\Phyl$ \where{def:phylogenetic-operad} phylogenetic operad\\
$\fieldc$ \where{notation:fieldc} complex plane\\
$\Intl(C)$ \where{notation:interval-circle} set of intervals in a circle with marked points\\
$\interior{X}$, $\partial X$, $\overline{X}$ \where{notation:interior}
interior, boundary, and closure of $X \subseteq \fieldc$\\
$\Reg(T)$ \where{notation:regt} set of regions in a planar tangle $T$\\ 
$\PTan$ \where{def:planar-tangle-topology} planar tangle operad\\
\blob\\
\textbf{Chapter \ref{ch:group-operad}}\blob\\
$\G$, $\gamma^{\G}$, $\operadunit^{\G}$ \where{not:g} action operad, operadic composition, colored unit\\
$\omega$ \where{not:g} augmentation\\
$\id_n$ \where{not:g} multiplicative unit in $\G(n)$\\
$\gbar$ \where{not:g} underlying permutation of $g \in \G(n)$\\
$\P$ \where{ex:trivial-group-operad} planar group operad\\
$\gsubc$ \where{not:gsubc} groupoid of $\colorc$-profiles with $\G$-action\\
$\gseqcm$ \where{notation:gseqcm} $\colorc$-colored $\G$-sequences in $\M$\\
$\I^{\G}$ \where{g-unit} unit $\colorc$-colored $\G$-sequence\\
$\circg$ \where{notation:circg} $\G$-circle product\\
$\goperadcm$ \where{def:g-operad} $\colorc$-colored $\G$-operads in $\M$\\
$\algmo$ \where{def:g-operad-algebra} category of $\O$-algebras for a $\colorc$-colored $\G$-operad $\O$\\
$\PaG$ \where{not:pag} parenthesized $\G$-operad\\
$E_G$ \where{def:translation-category} translation category of a group $G$\\
$\Gpd$ \where{not:gpd} category of small groupoids\\
$E_{\G}$ \where{not:esubg} translation category $\G$-operad\\
$\Delta$ \where{ex:translation-g} category of finite ordinal numbers\\
$\Nerve$ \where{nerve-category} nerve functor\\
$E\G$ \where{not:eofg} nerve of translation category $\G$-operad\\
\blob\\
\textbf{Chapter \ref{ch:braided_operad}}\blob\\
$B_n$ \where{def:braid-group} the $n$th braid group\\
$s_i$ \where{def:braid-group} a generating braid\\
$PB_n$ \where{def:pure-braid-group} the $n$th pure braid group\\
$F(X;n)$ \where{not:fxn} configuration space\\
$\opendtwo$ \where{not:opendtwo} open unit disc in $\fieldr^2$\\
$\B$ \where{def:braid-group-operad} braid group operad\\
$\PB$ \where{def:pure-braid-group-operad} pure braid group operad\\
$\bseqcm$ \where{notation:bseqcm} $\colorc$-colored braided sequences in $\M$\\
$\circb$ \where{not:circb} braided circle product\\
$\Boperadcm$ \where{not:boperadcm} $\colorc$-colored braided operads in $\M$\\
$\PaB$ \where{ex:parenthesized-braid} parenthesized braided operad\\
$\Ctilde_2$ \where{not:ctilde} universal cover of the little $2$-cube operad\\
$\Dtilde_2$ \where{not:dtilde} universal cover of the little $2$-disc operad\\
\blob\\
\textbf{Chapter \ref{ch:ribbon_operad}}\blob\\
$R_n$ \where{def:ribbon-group} the $n$th ribbon group\\
$PR_n$ \where{not:pr} the $n$th pure ribbon group\\
$r_i$ \where{int:ribbon-group-generators} a generating ribbon\\
$\R$ \where{not:ribbongop} ribbon group operad\\
$\rseqcm$ \where{notation:rseqcm} $\colorc$-colored ribbon sequences in $\M$\\
$\circr$ \where{not:circr} ribbon circle product\\
$\Roperadcm$ \where{not:roperadcm} $\colorc$-colored ribbon operads in $\M$\\
$\PaR$ \where{not:par} parenthesized ribbon operad\\
$\Dftildetwo$ \where{not:dftilde} universal cover of the framed little $2$-disc operad\\
\blob\\ \blob\\
\textbf{Chapter \ref{ch:cactus}}\blob\\
$\rho^{(n)}_{p,q}$ \where{not:rhonpq} interval-reversing permutation\\
$Cac_n$ \where{def:cactus-group} the $n$th cactus group\\
$s^{(n)}_{p,q}$ \where{def:cactus-group} a generating cactus\\
$PCac_n$ \where{not:pcac} the $n$th pure cactus group\\
$k_{[i,j]}$ \where{not:kij} $k_i+\cdots+k_j$\\
$\Cac$ \where{thm:cactus-group-operad} cactus group operad\\
$\PCac$ \where{pure-cactus-action-operad} pure cactus group operad\\
$\cacseqcm$ \where{notation:cacseqcm} $\colorc$-colored cactus sequences in $\M$\\
$\circcac$ \where{not:circcac} cactus circle product\\
$\Cacoperadcm$ \where{not:cacoperadcm} $\colorc$-colored cactus operads in $\M$\\
\blob\\
\textbf{Chapter \ref{ch:naturality}}\blob\\
$\agop$ \where{not:actop} category of action operads\\
\blob\\
\textbf{Chapter \ref{ch:group-operad-category}}\blob\\
$\ell_v$ \where{not:ellv} planar structure at a vertex\\
$\ell_T$ \where{not:ellt} canonical ordering of the inputs in a planar tree\\
$\kappa : T \to \colorc$ \where{not:coloring} $\colorc$-coloring of a tree\\
$\profofu$ \where{not:profofu} profile of $u$\\
$\uparrow$ \where{ex:exceptional-edge} exceptional edge\\
$\Cor_n$ \where{ex:corolla} $n$-corolla\\
$T\bigl(\{\ub_j\};\uc;d\bigr)$ \where{not:tbcd} $2$-level tree\\
$T \comp_v T'$ \where{not:compv} tree substitution at a vertex $v$\\
$T(T_v)_{v\in \Vt(T)}$ \where{not:ttv} tree substitution\\
$(T,\sigma)$ \where{not:tsigma} a $\G$-tree with input equivariance $\sigma$\\
$\Prof(T,\sigma)$ \where{not:proftsigma} profile of a $\G$-tree\\
\begin{small}$(T,\sigma)\comp_v (T_v,\sigma_v)$\end{small} \where{def:treesub-g} $\G$-tree substitution at $v$\\
\begin{small}$(T,\sigma)(T_v,\sigma_v)_v$\end{small} \where{def:g-tree-sub} $\G$-tree substitution\\
$\gopc$ \where{def:gopc} symmetric operad encoding $\G$-operads in $\Set$\\
$\gopcm$ \where{def:gopcm} symmetric operad encoding $\G$-operads in $\M$\\
\blob\\
\textbf{Chapter \ref{ch:all-group-operads}}\blob\\
$\gopm$ \where{def:g-operads} category of all $\G$-operads in $\M$\\
$\sigma_{m,n}$ \where{def:special-permutation} special $(m,n)$-permutation\\
$\tensorg$ \where{not:tensorg} $\G$-tensor product for $\G$-operads in $\Set$\\
\blob\\
\textbf{Chapter \ref{ch:w-group}}\blob\\
$\GTreec$ \where{def:substitution-category} substitution category of $\G$-trees\\
$\O(v)$ \where{def:vertex-decoration} $\profofv$-entry of $\O$\\
$\O(T)$ \where{not:ooft} $\O$-decoration of $T$\\
$\gamma_{(T,\sigma)}$ \where{goperad-gammat} operadic composition parametrized by a $\G$-tree\\
$(J, \mu, 0, 1, \epsilon)$ \where{notation:segment} a segment\\
$\Delta^1$ \where{not:deltaone} simplicial interval $\Delta(-,[1])$\\
$\J(T)$ \where{not:joft} $J^{\otimes |\Int(T)|}$\\
$\Wg$ \where{wgoduc} $\G$-Boardman-Vogt construction\\
$\Ws$ \where{ex:bv-symmetric} symmetric Boardman-Vogt construction\\
$\xi$ \where{not:xi} standard section\\
$\Wp$ \where{not:wofp} planar Boardman-Vogt construction\\
$\Wb$ \where{not:wofb} braided Boardman-Vogt construction\\
$\Wr$ \where{not:wofb} ribbon Boardman-Vogt construction\\
$\Wcac$ \where{not:wofb} cactus Boardman-Vogt construction\\
\blob\\
\textbf{Chapter \ref{ch:group-tree-category}}\blob\\
$\Treecatg$ \where{def:treecatg} $\G$-tree category\\
$\Treecatp$ \where{ex:treecatp} $\P$-tree category\\
$\Omega_p$ \where{ex:treecatp} category of planar rooted trees\\
$\Treecats$ \where{ex:treecats} $\S$-tree category\\
$\Omega$ \where{ex:treecats} category of rooted trees\\
$\Treecatb$ \where{ex:treecatbr} braided tree category\\
$\Treecatr$ \where{ex:treecatbr} ribbon tree category\\
$\Treecatcac$ \where{ex:treecatcac} cactus tree category\\
$\Gbar^{(T,\sigma)}$ \where{not:gbartsigma} colored set generated by vertices in $T$\\
$\G^{(T,\sigma)}$ \where{not:gtsigma} free $\G$-operad generated by $(T,\sigma)$\\
$\G^-$ \where{gtree-gop-functor} functor from $\Treecatg$ to $\gopset$\\
$\P^-$ \where{pgs-treecat-operad} functor from $\Treecatp$ to $\popset$\\
$\S^-$ \where{pgs-treecat-operad} functor from $\Treecats$ to $\sopset$\\
$\B^-$ \where{not:rbdash} functor from $\Treecatb$ to $\bopset$\\
$\R^-$ \where{not:rbdash} functor from $\Treecatr$ to $\ropset$\\
$\Cac^-$ \where{not:cacdash} functor from $\Treecatcac$ to $\cacopset$\\
\blob\\
\textbf{Chapter \ref{ch:contractibility}}\blob\\
$\Tunderline^I$ \where{not:iclosure} $I$-closure of $T$\\
$\Tunderline$ \where{not:closureoft} closure of $T$\\
$\BVt(\Tunderline)$ \where{not:bvertex} set of bottom vertices in $\Tunderline$\\
$\Cl(T,\sigma)$ \where{not:closureoftsigma} closure of $(T,\sigma)$\\
$\Ex$ \where{not:outputex} output extension\\
$\That$ \where{not:outputex} underlying closed planar tree of $\Ex(T,\id_0)$\\
\blob\\
\textbf{Chapter \ref{ch:reedy}}\blob\\
$\deg$ \where{def:treecatg-reedy} number of vertices\\
$f \boxslash g$ \where{notation:fboxslashg} $f$ has the left lifting property w.r.t. $g$\\
$^{\boxslash}\!\cala$ \where{notation:boxslasha} morphisms with the left lifting property w.r.t. $\cala$\\
$\cala^{\boxslash}$ \where{notation:boxslasha} morphisms with the right lifting property w.r.t. $\cala$\\
$\Cell(\cali)$ \where{not:cell} relative $\cali$-cell complexes\\
\blob\\
\textbf{Chapter \ref{ch:group-nerve}}\blob\\
$\Yonedag$ \where{not:yoneda} Yoneda embedding of $\Treecatg$\\
$\Realg$, $\Nerveg$ \where{not:realg} $\G$-realization/nerve\\
$\Realb$, $\Nerveb$ \where{ex:braided-nerve} braided realization/nerve\\
$\Realr$, $\Nerver$ \where{ex:ribbon-nerve} ribbon realization/nerve\\
$\Realcac$, $\Nervecac$ \where{ex:cactus-nerve} cactus realization/nerve\\
$\Govero$ \where{def:goperad-colim} category of $\G^-$ over $\O$\\
$\Gpr$ \where{not:gpr} functor from $\Govero$ to $\gopset$\\
$\tensorsub{\G}$ \where{def:treecatg-presheaf-monoidal} monoidal product of $\Treecatg$-presheaves\\
$\Homsub{\G}$ \where{not:homsubg} internal hom of $\Treecatg$-presheaves\\
\blob\\
\textbf{Chapter \ref{ch:infinity-goperad}}\blob\\
$X(T,\sigma)_1$ \where{def:segal-condition} corolla base of $X(T,\sigma)$\\
$\chi_{(T,\sigma)}$ \where{g-segal-map} $\G$-Segal map\\
$\xi_t$ \where{vertex-inclusion} vertex inclusion at $t$\\
$\Segalg(T,\sigma)$ \where{def:segal-core} $\G$-Segal core\\
$\face{\xi_t}(T,\sigma)$ \where{not:facexit} image induced by $\xi_t$\\
$\iota_{(T,\sigma)}$ \where{not:segalcor-inclusion} $\G$-Segal core inclusion\\
$\face{\phi}(V,\sigma^V)$ \where{not:phiface} $\phi$-face of $\Yonedag(V,\sigma^V)$\\
$\horn{\phi}(V,\sigma^V)$ \where{not:phihorn} $\phi$-horn of $\Yonedag(V,\sigma^V)$\\
\blob\\
\textbf{Chapter \ref{ch:hc-nerve}}\blob\\
$\Realghc$, $\Nerveghc$ \where{hcr-hcn} coherent $\G$-realization and coherent $\G$-nerve\\
$\Realhc^{\P}$, $\Nervehc^{\P}$ \where{not:realphc} coherent $\P$-realization/nerve\\
$\Realhc^{\B}$, $\Nervehc^{\B}$ \where{ex:braided-hcrn} coherent braided realization/nerve\\
$\Realhc^{\R}$, $\Nervehc^{\R}$ \where{ex:ribbon-hcrn} coherent ribbon realization/nerve\\
$\Realhc^{\Cac}$, $\Nervehc^{\Cac}$ \where{ex:cactus-hcrn} coherent cactus realization/nerve\\
\blob\\
\textbf{Chapter \ref{ch:monoidal-cat}}\blob\\
$(F,F_2,F_0)$ \where{def:monoidal-functor} a monoidal functor\\
$\MCat$ \where{def:mcat-operad} monoidal category operad\\
$x_{[i,j]}$ \where{not:subscript-interval} $(x_i,\ldots,x_j)$\\
$p_{\st}$ \where{not:pst} strict non-associative monomial associated to $p$\\
$p_{\nor}$ \where{not:pnor} normal non-associative monomial associated to $p$\\
\textit{\textsf{PaP}} \where{not:fresse-pap} parenthesized permutation operad\\
\textit{\textsf{PaP}}$_+$ \where{not:fresse-papplus} unitary parenthesized permutation operad\\
$\MCatst$ \where{def:mcatst} strict monoidal category operad\\
\textit{\textsf{CoP}}$_+$ \where{not:copplus} unitary colored permutation operad\\
$\MFun$ \where{def:monoidal-functor-operad} monoidal functor operad\\
$\MFunstg$ \where{def:strong-monoidal-functor-operad} strong monoidal functor operad\\
\blob\\
\textbf{Chapter \ref{ch:gmonoidal-cat}}\blob\\
$\MCatg$ \where{def:mcatg-underlying-categories} $\G$-monoidal category operad\\
$\MCatp$ \where{ex:mcatpn} $\P$-monoidal category operad\\
$\MCats$ \where{ex:mcatsn} $\S$-monoidal category operad\\
$\MCatb$ \where{ex:mcatbn} $\B$-monoidal category operad\\
$\MCatr$ \where{ex:mcatrn} $\R$-monoidal category operad\\
$\MCatcac$ \where{ex:mcatcacn} $\Cac$-monoidal category operad\\
$\MCatgst$ \where{def:mcatgst-operad} strict $\G$-monoidal category operad\\
$\MFung$ \where{def:gmonoidal-functor-operad} $\G$-monoidal functor operad\\
$\MFungstg$ \where{def:strong-gmonoidal-functor-operad} strong $\G$-monoidal functor operad\\
\blob\\ \blob\\
\textbf{Chapter \ref{ch:coherence-gmonoidal}}\blob\\
$\MCatg(A)$ \where{mcatg-free-forget} free $\G$-monoidal category generated by $A$\\
$\MCatgst(A)$ \where{lambdag-a} free strict $\G$-monoidal category generated by $A$\\
\blob\\
\textbf{Chapter \ref{ch:braided-symmetric-monoidal-cat}}\blob\\
$\MCatb$ \where{not:mcatb-operad} $\B$-monoidal category operad\\
\textit{\textsf{PaB}} \where{rk:fresse-pab} parenthesized braid operad\\
\textit{\textsf{PaB}}$_+$ \where{not:pabplus} unitary parenthesized braid operad\\
$\MCatbst$ \where{strict-bmonoidal-braided-monoidal} strict $\B$-monoidal category operad\\
\textit{\textsf{CoB}}$_+$ \where{not:cobplus} unitary colored braid operad\\
$\MFunb$ \where{mfunb-algebra} $\B$-monoidal functor operad\\
$\MFunbstg$ \where{mfunbstg-algebra} strong $\B$-monoidal functor operad\\
$\MCatb(A)$ \where{free-braided-monoidal-strictification} free braided monoidal category generated by $A$\\
$\MCatbst(A)$ \where{free-braided-monoidal-strictification} free strict braided monoidal category of $A$\\
$\MCats$ \where{not:mcats} $\S$-monoidal category operad\\
\textit{\textsf{PaS}} \where{rk:fresse-pas} parenthesized symmetry operad\\
\textit{\textsf{PaS}}$_+$ \where{not:pasplus} unitary parenthesized symmetry operad\\
$\MCatsst$ \where{strict-smonoidal-symmetric-monoidal} strict $\S$-monoidal category operad\\
\textit{\textsf{CoS}}$_+$ \where{not:cosplus} unitary colored symmetry operad\\
$\MFuns$ \where{mfuns-algebra} $\S$-monoidal functor operad\\
$\MFunsstg$ \where{mfunsstg-algebra} strong $\S$-monoidal functor operad\\
$\MCats(A)$ \where{free-symmetric-monoidal-strictification} free symmetric monoidal category generated by $A$\\
$\MCatsst(A)$ \where{free-symmetric-monoidal-strictification} free strict symmetric monoidal category of $A$\\
\blob\\
\textbf{Chapter \ref{ch:ribbon-moncat}}\\
$\MCatr$ \where{not:mcatr} $\R$-monoidal category operad\\
$\MCatrst$ \where{strict-rmonoidal-ribbon-monoidal} strict $\R$-monoidal category operad\\
$\MFunr$ \where{mfunr-algebra} $\R$-monoidal functor operad\\
$\MFunrstg$ \where{mfunrstg-algebra} strong $\R$-monoidal functor operad\\
$\MCatr(A)$ \where{free-ribbon-monoidal-strictification} free ribbon monoidal category generated by $A$\\
$\MCatrst(A)$ \where{free-ribbon-monoidal-strictification} free strict ribbon monoidal category of $A$\\
\blob\\
\textbf{Chapter \ref{ch:coboundary-monoidal-category}}\blob\\
$(a_1\otimes \cdots \otimes a_n)_l$ \where{not:left-normalized} left normalized monoidal product\\
$(a_1\otimes \cdots \otimes a_n)_r$ \where{not:left-normalized} right normalized monoidal product\\
$c^{(k)}$ \where{def:coboundary-cat-cpq} interval-reversing isomorphism\\
$a_{[i,j]}$ \where{not:coboundary-monoidal} $a_ia_{i+1}\cdots a_j$ or $\tensorunit$\\
$a_{[j,i]\downarrow}$ \where{not:coboundary-monoidal} $a_ja_{j-1} \cdots a_i$ or $\tensorunit$\\ 
$\MCatcac$ \where{not:mcatcac} $\Cac$-monoidal category operad\\
$\MCatcacst$ \where{strict-cacmonoidal-coboundary-monoidal} strict $\Cac$-monoidal category operad\\
$\MFuncac$ \where{mfuncac-algebra} $\Cac$-monoidal functor operad\\
$\MFuncacstg$ \where{mfuncacstg-algebra} strong $\Cac$-monoidal functor operad\\
$\MCatcac(A)$ \where{free-coboundary-monoidal-strictification} free coboundary monoidal category generated by $A$\\
$\MCatcacst(A)$ \where{free-coboundary-monoidal-strictification} free strict coboundary monoidal category of $A$\\
\end{tabbing}
\printindex

\end{document}

%% file: operads_with_group_equivariance.tex
\chapter{Introduction}\label{ch:introduction}

\section{Overview}\label{sec:overview}

This monograph provides a coherent development of operads, infinity operads, and monoidal categories, equipped with equivariant structures from general groups.  The group actions in all three cases are encoded by an \emph{action operad}, which is a sequence of groups that forms a planar operad in sets that is compatible with the group structures and that maps to the symmetric group operad.  The precise formulation is in Definition \ref{def:augmented-group-operad}.  This monograph has two main objectives:
\begin{enumerate}
\item We establish a foundation for \emph{group operads} and for their infinity analogues called \emph{infinity group operads}.  Planar operads in the category $\Set$ are Lambek's multicategories \cite{lambek}.   A group operad is a planar operad, enriched in some symmetric monoidal category, with a compatible equivariant structure by an action operad.  Infinity group operads are higher coherent analogues of group operads.  Restricting to suitable action operads, we obtain planar, symmetric, braided, ribbon, and cactus operads, and their infinity analogues.
\item We define general monoidal categories equipped with an equivariant structure by an action operad on multiple monoidal products, and prove coherence and strictification results for them.  Examples of monoidal categories with action operad equivariance include symmetric, braided, ribbon, and coboundary monoidal categories.  We recover and unify known coherence results for these kinds of monoidal categories.
\end{enumerate}

Lambek \cite{lambek} originally introduced multicategories, which we call planar operads in $\Set$, in the context of proof theory.  A planar operad is an extension of a category in which the domain of each morphism is a finite sequence of objects.  Similar to category theory, operad theory is a powerful tool to organize operations with multiple inputs and one output, and has numerous applications throughout mathematics and physics.  For example, in homotopy theory, symmetric operads, which are planar operads equipped with symmetric group actions, are the primary tools for recognizing iterated loop spaces \cite{boardman-vogt,may}.  In mathematical physics, symmetric operads provide a convenient framework to study homotopical aspects of algebraic quantum field theory \cite{hk}, as done in \cite{bsw,bsw-haqft,bsw2,yau-hqft}.  Many other applications of operad theory are discussed in \cite{baez-otter,jones,markl08,mss,mendez,spivak,vallette,yau-wd}.

While operads usually appear as either planar operads or symmetric operads, there are some important situations where other types of equivariant structures arise naturally.  For example, the recognition of $E_2$-operads \cite{fiedorowicz,fresse-gt} uses the universal cover of the little $2$-cube operad, which is a planar operad equipped with braid group actions, called a braided operad.  Along the same lines,  the recognition of framed $E_2$-operads \cite{wahl} involves the universal cover of the framed little $2$-disc operad, which is a planar operad equipped with ribbon group actions, called a ribbon operad.  These examples will  be discussed in Chapter \ref{ch:braided_operad} and Chapter \ref{ch:ribbon_operad}.  Part \ref{part:operads-group-eq} and Part \ref{part:construction-group-operad} of this book are devoted to developing the basic theory of operads equipped with group equivariance, called group operads. 

Infinity group operads belong to the much larger landscape of higher categorical structures, including: 
\begin{itemize}\item higher categories \cite{batanin98,joyal,leinster-higher,lurie,riehl-categorical,simpson};\index{infinity-category@$\infty$-category}
\item higher symmetric operads \cite{batanin08,batanin10,cm11b,cm13,moerdijk,mt,mw07,mw09,weiss};\index{infinity-symmetric operad@$\infty$-symmetric operad}
\item $\infty$-cyclic operads \cite{hry-cyclic};\index{infinity-cyclic operad@$\infty$-cyclic operad}
\item $\infty$-properads and $\infty$-wheeled properads \cite{hry15,hry18a,hry18b}.\index{infinity-properad@$\infty$-properad}
\end{itemize}
Infinity categories and various flavors of $\infty$-operads are interesting in their own right, and they are also important in applications.  For example, $\infty$-categories play a central role in Lurie's work \cite{lurie-tft} in\index{topological field theory} topological field theories.  Speculatively, a higher coherent analogue of Fresse's work \cite{fresse-gt} on the\index{Grothendieck-Teichmuller group@Grothendieck-Teichm\"{u}ller group} Grothendieck-Teichm\"{u}ller groups would involve the\index{infinity-braided operad@$\infty$-braided operad} $\infty$-braided operads developed here.  Along the same lines, a higher coherent analogue of the work in \cite{dhr} would involve our\index{infinity-ribbon operad@$\infty$-ribbon operad} $\infty$-ribbon operads.  Furthermore, homotopical quantum field theories\index{quantum field theory}\index{homotopical quantum field theory} as developed in \cite{bsw-haqft,yau-hqft}, with a flavor of braids \cite{bhk,oeckl1,oeckl2,sasai}--i.e., homotopical braided quantum field theories--would likely involve $\infty$-braided operads or the\index{braided operad!Boardman-Vogt construction}\index{Boardman-Vogt construction!braided} braided Boardman-Vogt construction of braided operads in Definition \ref{def:braided-bv}.  Part \ref{part:infinity-group-operads} of this book is devoted to infinity group operads.

Monoidal categories and symmetric monoidal categories have long been an integral part of category theory \cite{benabou,kelly2,kelly-enriched,maclane-rice}.  A major part of their usage and understanding is a number of results collectively known as their Coherence Theorem.  Joyal and Street \cite{joyal-street} introduced\index{braided monoidal category}\index{monoidal category!braided} braided monoidal categories and\index{ribbon monoidal category}\index{monoidal category!ribbon} ribbon monoidal categories, and established their coherence.  Braided and ribbon monoidal categories arise naturally as the categories of representations of\index{quantum group} quantum groups \cite{drinfeld-quantum,chari,kassel,savage}.  A\index{coboundary monoidal category}\index{monoidal category!coboundary} coboundary monoidal category is a modification of a braided monoidal category, with an involutive braiding and with the Hexagon Axioms replaced by the Cactus Axiom.  Coboundary monoidal categories appear as the categories of representations of a\index{coboundary Hopf algebra} coboundary Hopf algebra \cite{drinfeld} and of\index{crystal} crystals of a finite dimensional complex reductive\index{Lie algebra} Lie algebra \cite{hen-kam}.  Strictification of coboundary monoidal categories was proved by Gurski \cite{gurski}.  In Part \ref{part:monoidal-cat} we will define $\G$-monoidal categories for an action operad $\G$, and prove several versions of coherence for $\G$-monoidal categories.  Restricting to suitable action operads, these coherence results restrict to and unify those for symmetric, braided, ribbon, and coboundary monoidal categories.

This book consists of the following four parts.
\begin{description}
\item[Part \ref{part:operads-group-eq}] Operads with Group Equivariance
\item[Part \ref{part:construction-group-operad}] Constructions of Group Operads
\item[Part \ref{part:infinity-group-operads}] Infinity Group Operads
\item[Part \ref{part:monoidal-cat}] Coherence for Monoidal Categories with Group Equivariance
\end{description}
Part \ref{part:infinity-group-operads} on infinity group operads depends on the first two parts.  In the first three parts, we work over a symmetric monoidal category with all small (co)limits in which the monoidal product commutes with small colimits on each side.  Part \ref{part:monoidal-cat} on monoidal categories with group equivariance is independent of Part \ref{part:infinity-group-operads} and almost all of Part \ref{part:construction-group-operad}, so it can be read right after Part \ref{part:operads-group-eq}.  

This monograph is aimed at graduate students and researchers with an interest in operads, $\infty$-operads, and monoidal categories.  We assume the reader is familiar with basic category theory, including (co)limits, adjoint functors, coends, and (symmetric) monoidal categories up to Mac Lane's Coherence Theorem.  Familiarity with other flavors of monoidal categories, such as braided and ribbon monoidal categories, is not assumed.  Moreover, no prior knowledge of operads and $\infty$-operads is needed.  In addition to being a coherent reference for the topics covered, this book is also suitable for a graduate student seminar and a reading course.  For example:
\begin{itemize}
\item The first three parts can be the basis of a graduate seminar on operads and $\infty$-operads.
\item Part \ref{part:operads-group-eq} and Part \ref{part:monoidal-cat} can be used for a reading course covering operads and monoidal categories.
\end{itemize}
Next we describe the content of each part in more detail.

\subsection*{Part \ref{part:operads-group-eq} : Operads with Group Equivariance}

The first part is a leisurely introduction to colored operads equipped with equivariant structures from general groups.  The main concept is that of a \emph{colored $\G$-operad} for an action operad $\G$.  In the rest of this Introduction, we will usually drop the word \emph{colored}, with the understanding that our treatment of $\G$-operads is entirely at the colored level.  Assuming no prior knowledge of operads, we carefully introduce operads with action operad equivariance, and provide many detailed examples.  The reader who already knows about planar and symmetric operads can go straight to Chapter \ref{ch:group-operad}, where action operads and $\G$-operads are defined.   

In Chapter \ref{ch:planar-operad} we set the stage for the rest of this monograph by recalling planar operads.  Our discussion generally goes from a more conceptual and compact form to a longer and more explicit form.  Planar operads are first defined conceptually as monoids with respect to a planar circle product.  We then unpack this definition and describe a planar operad in a more familiar form, in terms of a unit and an operadic composition.  A planar operad in sets is what Lambek \cite{lambek} called a \emph{multicategory}, which is like a category but with the domain of each multi-morphism a finite sequence of objects.  So a planar operad is a multicategory enriched in an ambient symmetric monoidal category.  Next we define algebras over a planar operad and end this chapter with a list of examples of planar operads and their algebras.

Chapter \ref{ch:symmetric-operad} is about symmetric operads, which are planar operads equipped with an equivariant structure by the symmetric groups.  They are also known as \emph{symmetric multicategories}, with the symmetric groups permuting the domain objects of each multi-morphism.  As in the planar case, we first define symmetric operads as monoids with respect to a symmetric circle product, which is due to Kelly \cite{kelly} in the one-colored case.  Then we unpack this definition to the one introduced by May \cite{may}.  After a brief discussion of algebras over a symmetric operad, we discuss the little cube operads  due to Boardman-Vogt \cite{bv-everything,boardman-vogt} and May \cite{may}, the little disc operads, and the framed little disc operads due to Getzler \cite{getzler}.  We end this chapter with three examples of symmetric operads from outside of homotopy theory and category theory.  They are (i) Male's operad of graph operations in non-commutative probability theory \cite{male}, (ii) the phylogenetic operad in evolutionary biology due to Baez and Otter \cite{baez-otter}, and (iii) Jones's planar tangle operad \cite{jones} that defines planar algebras.

In Chapter \ref{ch:group-operad} we first define action operads, which are essentially due to Wahl \cite{wahl} and are also used in Yoshida \cite{yoshida} and Zhang \cite{zhang}.  We adopted the \emph{action operad} terminology from Corner-Gurski \cite{corner-gurski} and Gurski \cite{gurski}.  An action operad is a sequence of groups that forms a planar operad in $\Set$ whose operadic composition is compatible with the group multiplications.  Moreover, an action operad is equipped with an augmentation to the symmetric group operad, which is made up of the symmetric groups.  For each action operad $\G$, we define \emph{$\G$-operads} as monoids with respect to a $\G$-circle product.  In unpacked form, a $\G$-operad is a planar operad with a compatible $\G$-equivariant structure.  We recover planar operads and symmetric operads if we take $\G$ to be the planar group operad $\P$ and the symmetric group operad $\S$, respectively.  We also discuss algebras over a $\G$-operad and some examples.  

Two remarks are in order.  We use the generic term \emph{group operads} to refer to $\G$-operads for some action operad $\G$.  This is different from the usage in \cite{yoshida,zhang}, in which an action operad $\G$ is called a group operad.  Moreover, the definition of an action operad requires that of a planar operad in Chapter \ref{ch:planar-operad}, but not that of a symmetric operad in Chapter \ref{ch:symmetric-operad}.  We discuss symmetric operads before group operads in order to give the reader with no prior knowledge of operads more concrete examples earlier.  

The last three chapters in Part \ref{part:operads-group-eq} are about three main examples of action operads and their corresponding group operads.

In Chapter \ref{ch:braided_operad} we discuss the action operad $\B$, called the braid group operad, that is made up of Artin's braid groups.  The unpacked form of a $\B$-operad, as in Chapter \ref{ch:group-operad} for the action operad $\B$, is usually called a \emph{braided operad}.  Originally introduced by Fiedorowicz \cite{fiedorowicz} in the one-colored case, a braided operad is a planar operad together with a compatible braid group action.  There is also an action operad $\PB$ consisting of the pure braid groups, the elements of which have trivial underlying permutations.  We end this chapter with Fiedorowicz's example of the level-wise universal cover of the little $2$-cube operad, which is a topological braided operad.  The universal cover of the little $2$-disc operad is also discussed.

Chapter \ref{ch:ribbon_operad} is the ribbon analogue of Chapter \ref{ch:braided_operad}.  Geometrically, a ribbon is a braid in which each string is replaced by a strip that can have an integer multiple of full $2\pi$ twists.  The ribbon groups constitute the ribbon group operad $\R$, and $\R$-operads are called \emph{ribbon operads}.  Originally introduced by Wahl \cite{wahl} in the one-colored case, a ribbon operad is equivalent to a planar operad together with a compatible ribbon group action.  There is also an action operad $\PR$ consisting of the pure ribbon groups, the elements of which have trivial underlying permutations.  We end this chapter with Wahl's example of the level-wise universal cover of the framed little $2$-disc operad, which is a topological ribbon operad.  

Chapter \ref{ch:cactus} is the cactus analogue of Chapter \ref{ch:braided_operad} and Chapter \ref{ch:ribbon_operad}.  As proved by Henriques and Kamnitzer \cite{hen-kam}, the cactus groups act on multiple monoidal products in Drinfel'd's coboundary monoidal categories \cite{drinfeld}.  The category of representations of a coboundary Hopf algebra \cite{drinfeld} and the category of crystals of a finite dimensional complex reductive Lie algebra \cite{hen-kam} are natural examples of coboundary monoidal categories.  Cactus groups are also related to the symmetry groups of closed connected manifolds with a cubical cell structure \cite{djs} and the geometry of moduli spaces of real, genus 0, stable curves with marked points \cite{dev}.  We define the cactus groups and show that they form an action operad $\Cac$, which was first observed by Gurski \cite{gurski}.  There is also an action operad $\PCac$ consisting of the pure cactus groups, the elements of which have trivial underlying permutations.  We end this chapter with the observation that the braid group operad $\B$ is incompatible with the cactus group operad $\Cac$ in a specific sense.  This finishes the first part.

\subsection*{Part \ref{part:construction-group-operad} : Constructions of Group Operads}

The second part of this monograph is about constructions and categorical properties of action operads and group operads.  In the rest of this Introduction, $\G$ denotes an action operad as in Chapter \ref{ch:group-operad}.

In Chapter \ref{ch:naturality} we consider naturality properties of $\G$-operads and their algebras.  First we study the induced adjunction between categories of $\G$-operads associated to a morphism of action operads.  While the right adjoint is easy to describe, we produce an explicit formula for the left adjoint.  Examples of this adjunction include those between ribbon operads and symmetric operads, between cactus operads and symmetric operads, and so forth.  Next we observe that the category of $\G$-operads is natural with respect to a symmetric monoidal functor on the underlying symmetric monoidal category, such as the singular chain functor and the nerve functor.  In the last section, we observe that each morphism of $\G$-operads has an induced adjunction between the categories of algebras.

In Chapter \ref{ch:group-operad-category} we construct an explicit symmetric operad whose algebras are exactly $\G$-operads.  There are two key benefits of this construction.  First, it provides another conceptual description of $\G$-operads, which are first defined in Chapter \ref{ch:group-operad} as monoids with respect to the $\G$-circle product.  Second, the symmetric operad that encodes $\G$-operads involves and motivates the combinatorial objects of \emph{$\G$-trees}, which play a crucial role in the $\G$-Boardman-Vogt construction in Chapter \ref{ch:w-group} and in infinity $\G$-operads in Part \ref{part:infinity-group-operads} .  A $\G$-tree consists of a planar tree and an element in $\G$, called the input equivariance, at the same level as the number of inputs of the planar tree.  We think of the input equivariance as a $\G$-equivariant operator that acts on the inputs of the planar tree.  The idea is that each vertex in a planar tree corresponds to the entry in a $\G$-operad with the same inputs-output combination, called the profile, while the input equivariance corresponds to the $\G$-equivariant structure.  Corresponding to the operadic composition in $\G$-operads is the operation of $\G$-tree substitution, which is a combination of planar tree substitution and the action operad structure in $\G$.  One consequence of identifying $\G$-operads as algebras over a symmetric operad is that the category of $\G$-operads has all small limits and colimits.

In Chapter \ref{ch:all-group-operads} we consider the category of $\G$-operads in which the color sets are allowed to vary.  We first show that this category of all $\G$-operads has all small limits and colimits, and is natural with respect to a change of the action operad.  When the underlying symmetric monoidal category is $\Set$, we equip the category of all $\G$-operads with a symmetric monoidal closed category structure.  For the symmetric group operad $\S$, this recovers the Boardman-Vogt tensor product on the category of all symmetric operads in $\Set$.  We observe that the functor between categories of all group operads induced by a morphism of action operads  is a symmetric monoidal functor.  However, it is usually not a strong symmetric monoidal functor in examples of interest.  In the last section, we show that, with a general underlying symmetric monoidal category that has all small limits and colimits, the category of all $\G$-operads is locally finitely presentable.

In Chapter \ref{ch:w-group} we define the $\G$-Boardman-Vogt construction for $\G$-operads.  This construction is needed in Part \ref{part:infinity-group-operads} when we define the coherent $\G$-nerve and the coherent $\G$-realization functors.  For each $\G$-operad $\O$, its $\G$-Boardman-Vogt construction $\Wg\O$ is a natural thickening of $\O$ that takes into account homotopical structure of the underlying symmetric monoidal category in the form of a commutative segment.  It is also a $\G$-operad and is equipped with an augmentation \[\alpha : \Wg\O \to \O\] that has nice naturality properties.  Restricting to the symmetric group operad $\S$ and the underlying category of spaces with the unit interval as the commutative segment, we recover the original Boardman-Vogt construction for symmetric operads in \cite{boardman-vogt}.  For a general underlying symmetric monoidal category, we recover the Boardman-Vogt construction for symmetric operads by Berger and Moerdijk \cite{berger-moerdijk-bv} using the symmetric group operad $\S$, and also the Boardman-Vogt construction for planar operads using the planar group operad $\P$.  We point out that in \cite{berger-moerdijk-bv} the Boardman-Vogt construction of a symmetric operad is inductively defined as a sequential colimit, in which each morphism is a pushout that depends on the previous inductive step.  On the other hand, our $\G$-Boardman-Vogt construction of a $\G$-operad is defined in one-step as a coend indexed by a combinatorially defined category of $\G$-trees, called the substitution category, in which morphisms are given by $\G$-tree substitution.  This finishes the second part.

\subsection*{Part \ref{part:infinity-group-operads} : Infinity Group Operads}

The third part contains the main discussion of infinity group operads.  A quasi-category \cite{joyal,lurie}, also known as an $(\infty,1)$-category and a weak Kan complex \cite{boardman-vogt}, is a simplicial set that satisfies the Kan lifting property for inner horns.  This is a way to relax the axioms of a category from  equalities to their up-to-homotopy analogues.  The theory of quasi-categories has been extended to $\infty$-symmetric operads, $\infty$-cyclic operads, $\infty$-properads, and $\infty$-wheeled properads, as mentioned earlier in this Introduction.  A major part of the development of these up-to-homotopy categorical objects is the understanding of the combinatorially defined categories that model the shapes of composition and higher associativity.  The first three chapters of Part \ref{part:infinity-group-operads} are devoted to constructing and understanding this category in the context of group operads.  The last three chapters are about the associated presheaf category.

In Chapter \ref{ch:group-tree-category} we construct the $\G$-tree category $\Treecatg$ with $\G$-trees as objects and two-sided $\G$-tree substitution as morphisms and composition.  For the planar group operad $\P$ and the symmetric group operad $\S$, we recover the planar dendroidal category $\Omega_p$ and the symmetric dendroidal category $\Omega$ due to Moerdijk and Weiss \cite{mw07,mw09,weiss}.  For the braid, ribbon, and cactus group operads, we obtain the braided, ribbon, and cactus tree categories.  Next we observe that the assignment sending an action operad $\G$ to the $\G$-tree category is natural and has nice categorical properties.  Analogous to the functor $\Delta \to \Cat$, sending a finite ordinal $[n]$ to the small category generated by it, is a functor \[\nicexy{\Treecatg \ar[r] & \gopset}\] from the $\G$-tree category to the category of all $\G$-operads in $\Set$ in Chapter \ref{ch:all-group-operads}.  We show that this functor is fully faithful and behaves nicely with respect to a change of the action operad.

In Chapter \ref{ch:contractibility}, for an action operad $\G$ with bottom level $\G(0)$ the trivial group, which holds in all examples of interest, we determine the homotopy type of the $\G$-tree category $\Treecatg$ by showing that its nerve is a contractible simplicial set.  This is analogous to the contractibility of the nerve of the finite ordinal category $\Delta$.  The main difference is that, on the one hand, the contractibility of $\Delta$ is the result of it having a terminal object $[0]$.  On the other hand, the $\G$-tree category has neither a terminal object nor an initial object.  So we need to work a bit harder to prove its contractibility.  The key point is that the $\G$-tree category has a full reflective subcategory $\Treecatug$ consisting of $\G$-trees with an empty set of inputs, called closed $\G$-trees.  As a result, the nerves of the categories $\Treecatg$ and $\Treecatug$ are homotopy equivalent.  We then show that the full reflective subcategory $\Treecatug$ has a contractible nerve by constructing a zig-zag of natural transformations between the identity functor and a constant functor.  In fact, our proof shows that $\Treecatug$ is a Grothendieck strict test category.  For the symmetric group operad $\S$,  we recover a result in Ara-Cisinski-Moerdijk \cite{acm} that says that the symmetric dendroidal category $\Omega$, which is equivalent to the $\S$-tree category $\Treecats$, has a contractible nerve.

In Chapter \ref{ch:reedy} the Reedy category structure on the finite ordinal category is extended to a dualizable generalized Reedy category structure on the $\G$-tree category $\Treecatg$.  We begin by defining $\G$-tree analogues of cofaces and codegeneracies as special cases of $\G$-tree substitution.  Then we show that each morphism in the $\G$-tree category factors into a composite of codegeneracies, isomorphisms, and cofaces.  Using this factorization, the dualizable generalized Reedy category structure on $\Treecatg$ is established next.  As a consequence, the diagram categories $\M^{\Treecatg}$ and $\M^{\Treecatgop}$ based on a cofibrantly generated model category $\M$, indexed by the $\G$-tree category $\Treecatg$ and its opposite category, both admit Reedy-type model category structures.  In the last section, we show that the $\G$-tree category admits an even tighter structure called an \emph{EZ-category}, with EZ standing for Eilenberg-Zilber, which is a particularly nice kind of a dualizable generalized Reedy category.

Chapter \ref{ch:group-nerve} is the first of three chapters in which we study the category of $\Treecatg$-presheaves, which is the $\G$-operad analogue of the category of simplicial sets.  The realization-nerve adjunction, between the categories of simplicial sets and of small categories, is the backbone of categorical homotopy theory.  In this chapter, we first construct the $\G$-realization-nerve adjunction
\[\nicexy@C+.4cm{\Settotreecatgop \ar@<2pt>[r]^-{\Realg} & \gopset \ar@<2pt>[l]^-{\Nerveg}}\] in which $\gopset$ is the category of all $\G$-operads in $\Set$.  For the symmetric group operad $\S$, we recover the dendroidal nerve in Moerdijk-Weiss \cite{mw09}.  We show that the counit of the $\G$-realization-nerve adjunction is a natural isomorphism, so the $\G$-nerve is a full and faithful functor.  Using the $\G$-nerve and the symmetric monoidal product of $\G$-operads in $\Set$ in Chapter \ref{ch:all-group-operads}, we equip the category of $\Treecatg$-presheaves with a symmetric monoidal closed structure.  With respect to these symmetric monoidal structures, we show that the $\G$-realization functor is strong symmetric monoidal, so the $\G$-nerve is symmetric monoidal.  In the last two sections, we observe that the $\G$-realization functor, the $\G$-nerve functor, and the category of $\Treecatg$-presheaves behave nicely with respect to a change of the action operad.

In Chapter \ref{ch:infinity-goperad} we establish the $\G$-operad analogue of the Nerve Theorem\index{Nerve Theorem for categories} for categories.  The Nerve Theorem is an important fact in categorical homotopy theory.  It characterizes nerves of small categories as (i) strict quasi-categories or (ii) simplicial sets that satisfy the Segal condition.  To prove the $\G$-operad version, we begin by defining the $\G$-operad analogues of the Segal condition, faces, horns, and (strict) quasi-categories.  A $\Treecatg$-presheaf is said to satisfy the \emph{$\G$-Segal condition} if it is determined by its one-dimensional restrictions.  An \emph{infinity $\G$-operad} is defined as a $\Treecatg$-presheaf that satisfies an inner horn lifting property that is analogous to the inner Kan lifting property for simplicial sets.  A \emph{strict infinity $\G$-operad} is an infinity $\G$-operad whose inner horn fillers are unique.  The main observation of this chapter is that a $\Treecatg$-presheaf is a strict infinity $\G$-operad if and only if it satisfies the $\G$-Segal condition, which in turn is equivalent to being isomorphic to the $\G$-nerve of some $\G$-operad in $\Set$.

In Chapter \ref{ch:hc-nerve} we discuss a homotopy coherent analogue of the $\G$-realization-nerve adjunction in Chapter \ref{ch:group-nerve}.  We define the coherent $\G$-realization-nerve adjunction
\[\nicexy@C+.4cm{\Settotreecatgop \ar@<2pt>[r]^-{\Realghc} & \gopm \ar@<2pt>[l]^-{\Nerveghc}}\] with $\gopm$ the category of all $\G$-operads in the underlying symmetric monoidal category $\M$, which is equipped with a commutative segment.  The $\G$-Boardman-Vogt construction in Chapter \ref{ch:w-group} plays a crucial role in the definitions of these functors.  For the symmetric group operad $\S$, we recover the homotopy coherent dendroidal nerve in \cite{mw09}.  We show that the coherent $\G$-realization of the $\G$-nerve is isomorphic to the $\G$-Boardman-Vogt construction.  This allows us to compute the morphisms and the internal hom object from the $\G$-nerve to the coherent $\G$-nerve.  The main result of this chapter says that, for a $\G$-operad whose entries are fibrant objects in $\M$, the coherent $\G$-nerve is an infinity $\G$-operad.  This observation ensures an ample supply of infinity $\G$-operads as the coherent $\G$-nerves of entrywise fibrant $\G$-operads.  For example, the coherent braided nerve of the universal cover of the little $2$-cube operad is an $\infty$-braided operad.  Similarly, the coherent ribbon nerve of the universal cover of the framed little $2$-disc operad is an $\infty$-ribbon operad.  This finishes the third part.

\subsection*{Part \ref{part:monoidal-cat} : Coherence for Monoidal Categories with Group Equivariance}

The last part of this monograph is about general small monoidal categories equipped with group actions on multiple monoidal products.  As mentioned above, this part is independent of Part \ref{part:infinity-group-operads} and most of Part \ref{part:construction-group-operad}.  With $\Cat$ denoting the category of small categories, we will mostly omit the word \emph{small} in the rest of this Introduction.  The two main objectives of Part \ref{part:monoidal-cat} are:
\begin{enumerate}[label=(\roman*)]
\item To construct general small monoidal categories equipped with group equivariance on multiple monoidal products in terms of action operads.
\item To establish their coherence results in a unifying manner, which restrict to coherence results for symmetric, braided, ribbon, and coboundary monoidal categories.  
\end{enumerate}

In Chapter \ref{ch:monoidal-cat} we set the stage by discussing monoidal categories in operadic terms.  First we construct the symmetric operad in $\Cat$ that encodes monoidal categories with general associativity and left/right unit isomorphisms.  Then we construct the related symmetric operads in $\Cat$ that encode (i) strict monoidal categories, (ii) general monoidal functors, and (iii) strong monoidal functors.  We will not reprove Mac Lane's Coherence Theorem for monoidal categories, which we consider as well-known and will use repeatedly.

In Chapter \ref{ch:gmonoidal-cat}, for an action operad $\G$, we construct a symmetric operad $\MCatg$ in $\Cat$, whose algebras are defined as \emph{$\G$-monoidal categories}.  The first main observation in this chapter is that a $\G$-monoidal category is exactly a monoidal category with general associativity and left/right unit isomorphisms, and a compatible action by $\G$ on multiple monoidal products.  In particular, for the planar group operad $\P$, a $\P$-monoidal category is precisely a monoidal category.  The identification of other types of monoidal categories with $\G$-equivariant structures will be given in later chapters.  We define a \emph{$\G$-monoidal functor} between $\G$-monoidal categories as a monoidal functor in the usual sense that is compatible with the $\G$-equivariant structures.  Then we discuss the related symmetric operads for (i) strict $\G$-monoidal categories, (ii) general $\G$-monoidal functors, and (iii) strong $\G$-monoidal functors.  We remark that strict $\G$-monoidal categories are also algebras over Gurski's categorical Borel construction \cite{gurski}.  Therefore, our $\G$-monoidal category operad $\MCatg$ may be regarded as a generalization of Gurski's categorical Borel construction that allows general monoidal categories with associativity and left/right unit isomorphisms.

In Chapter \ref{ch:coherence-gmonoidal} we establish several versions of coherence for $\G$-monoidal categories.  We first prove that each $\G$-monoidal category is adjoint equivalent to a strict $\G$-monoidal category via strong $\G$-monoidal functors.  Next we prove the free version of this coherence result.  It says that the free $\G$-monoidal category generated by a small category is equivalent to its free strict $\G$-monoidal category via a canonical strict $\G$-monoidal functor.  Then we prove that in the free $\G$-monoidal category generated by a set of objects, a diagram is commutative if and only if composites with the same (co)domain have the same image in the free strict $\G$-monoidal category generated by one object.  We also provide explicit descriptions of the free (strict) $\G$-monoidal category generated by a small category.

In the next three chapters, we apply the results in the previous two chapters to the braid group operad $\B$, the symmetric group operad $\S$, the ribbon group operad $\R$, and the cactus group operad $\Cac$.

In Chapter \ref{ch:braided-symmetric-monoidal-cat} we show that, for the braid group operad $\B$, $\B$-monoidal categories and (strong/strict) $\B$-monoidal functors are exactly braided monoidal categories and (strong/strict) braided monoidal functors, respectively.  As defined by Joyal and Street \cite{joyal-street}, a \emph{braided monoidal category} is a monoidal category equipped with a natural braiding \[\nicexy{X\otimes Y \ar[r]^-{\xi}_-{\cong} & Y \otimes X}\] that is compatible with the associativity and left/right unit isomorphisms.  With $\G=\B$ the coherence results for $\G$-monoidal categories in Chapter \ref{ch:coherence-gmonoidal} restrict to coherence results for braided monoidal categories and braided monoidal functors.  We remark that what Joyal and Street called a \emph{braided tensor functor} is what we call a strong braided monoidal functor, which is a braided monoidal functor with invertible structure morphisms.  Using the symmetric group operad $\S$ instead, we also recover symmetric monoidal categories, symmetric monoidal functors, and their coherence results.

In Chapter \ref{ch:ribbon-moncat}, using the ribbon group operad $\R$, we establish the ribbon analogues of the results in Chapter \ref{ch:braided-symmetric-monoidal-cat}.  A \emph{ribbon monoidal category}, which Joyal and Street \cite{joyal-street} called a \emph{balanced tensor category}, is a braided monoidal category together with a compatible natural isomorphism \[\nicexy{X \ar[r]^-{\theta}_-{\cong} & X}.\]  Just like braided monoidal categories, ribbon monoidal categories arise naturally as categories of representations of quantum groups \cite{chari,drinfeld-quantum,kassel,street}.  While a symmetric monoidal category is a braided monoidal category with an extra property, a ribbon monoidal category is a braided monoidal category with an extra structure $\theta$.  Therefore, the identification of $\R$-monoidal categories with ribbon monoidal categories requires a non-trivial modification of the proof in the braided case.

In Chapter \ref{ch:coboundary-monoidal-category}, using the cactus group operad $\Cac$, we establish the cactus analogues of the results in Chapter \ref{ch:braided-symmetric-monoidal-cat} and Chapter \ref{ch:ribbon-moncat}.  We show that $\Cac$-monoidal categories are exactly coboundary monoidal categories, and similarly for functors.  A \emph{coboundary monoidal category} is a monoidal category equipped with a natural isomorphism $\xi$ as in a braided monoidal category, but it satisfies different axioms.  Coboundary monoidal categories arise naturally as (i) the category of representations of a coboundary Hopf algebra \cite{drinfeld} and (ii) the category of crystals of a finite dimensional complex reductive Lie algebra \cite{hen-kam}.  The proof of the identification of $\Cac$-monoidal categories with coboundary monoidal categories follows the same outline as the braided case, but the details are more complicated and require more work.  This result provides another proof, as well as an operadic interpretation, of the statement that the cactus groups act on multiple monoidal products in coboundary monoidal categories, which was first observed in \cite{hen-kam}.  This finishes Part \ref{part:monoidal-cat}.

\section{Categorical Conventions}\label{sec:categorical-conventions}

Here we fix some notations and recall some basic categorical concepts.  All the basic category theory used in this book can be found in the references \cite{borceux1,borceux2,kelly-enriched,maclane}.  The more elementary books \cite{awodey,leinster,riehl} also have most of the category theory that we need.  Discussion of coends can also be found in the paper \cite{loregian}.

\subsection{Categories}\label{subsec:categories}

For a finite set $S$, its\index{cardinality} cardinality is denoted by\label{notation:cardinality} $|S|$.

A \index{groupoid}\emph{groupoid} is a category in which every morphism is an isomorphism.

For a\index{category} category $\C$, its collection of objects is denoted by\label{notation:object} $\Ob(\C)$.  To denote an\index{object} object $X$ in $\C$, we write either $X\in \Ob(\C)$ or $X\in\C$.  We call $\C$ \index{small category}\emph{small} if $\Ob(\C)$ is a set.  For two objects $X$ and $Y$ in $\C$, the set of morphisms from $X$ to $Y$ is denoted by either\label{notation:morphism-set} $\C(X,Y)$ or $\C(X;Y)$.  The identity morphism of $X$ is denoted by\label{notation:identity-morphism} $\Id_X$.  For composable morphisms $f : X \to Y$ and $g : Y \to Z$ in $\C$, their composite in $\C(X,Z)$ is denoted by either\label{notation:morphism-composition} $gf$ or $g \circ f$.

The \index{opposite category}\emph{opposite category} of a category $\C$ is denoted by\label{notation:opposite-category} $\C^{\op}$.  It has the same objects as $\C$ and morphism sets $\C^{\op}(X,Y) = \C(Y,X)$, with identity morphisms and composition inherited from $\C$. 

For a small category $\D$, the \index{diagram category}\emph{diagram category}\label{notation:diagram-category} $\C^{\D}$ has functors $\D \to \C$ as objects and natural transformations between such functors as morphisms.

A category $\C$ is \index{complete}\index{cocomplete}\emph{(co)complete} if it has all small (co)limits, i.e., (co)limits of functors whose domains are small categories.  For a functor $F : \D \to \C$, its colimit, if it exists, is denoted by either\label{notation:colimit} $\colim\, F$ or $\colim_{x\in \D}\, Fx$.
  
An \index{initial object}\emph{initial object} in $\C$, if it exists, is denoted by either\label{notation:initial-object} $\varnothing$ or $\varnothing^{\C}$.  It is characterized by the universal property that, for each object $X\in\C$, there exists a unique morphism $\varnothing \to X$.  

The categorical concept of a coend is central in this work.  Suppose given a functor\label{notation:coend} $F : \Cop \times \C \to \D$.
\begin{enumerate}
\item A \index{wedge}\emph{wedge} of $F$ is a pair $(X,\zeta)$ consisting of 
\begin{itemize}\item an object $X \in \D$ and 
\item a morphism $\zeta_c : F(c,c) \to X$ for each object $c \in \C$ 
\end{itemize}
such that the diagram \[\nicexy{F(c',c) \ar[d]_-{F(g,c)} \ar[r]^-{F(c',g)} & F(c',c') \ar[d]^-{\zeta_{c'}}\\ F(c,c) \ar[r]^-{\zeta_c} & X}\] is commutative for each morphism $g : c\to c' \in \C$.
\item A \index{coend}\emph{coend} of $F$ \cite{maclane} (IX.6) is an initial wedge $\left(\int^{c\in \C} F(c,c), \omega\right)$.  In other words, a coend of $F$ is a wedge of $F$ such that given any wedge $(X,\zeta)$ of $F$, there exists a unique morphism \[\nicexy{\dint^{c\in \C} F(c,c) \ar[r]^-{h} & X} \in \D\] such that the diagram
\[\nicexy{F(c,c) \ar[r]^-{\omega_c} \ar[dr]_-{\zeta_c} & \dint^{c\in \C} F(c,c) \ar[d]^-{h}\\ & X}\]
is commutative for each object $c \in \C$.  
\end{enumerate}
Coends are used in the circle product that defines each variant of operads in Part \ref{part:operads-group-eq}, in the $\G$-Boardman-Vogt construction for group operads in Chapter \ref{ch:w-group}, and in the coherent $\G$-nerve in Chapter \ref{ch:hc-nerve}.

\subsection{Monoidal Categories}\label{subsec:moncat}

Throughout Part \ref{part:operads-group-eq} to Part \ref{part:infinity-group-operads} of this book, \[(\M,\otimes,\tensorunit,\alpha,\lambda,\rho)\] denotes a complete and cocomplete symmetric monoidal category with monoidal product $\otimes$, monoidal unit $\tensorunit$, associativity isomorphism $\alpha$, left unit isomorphism $\lambda$, and right unit isomorphism $\rho$.  We always assume that the monoidal product commutes with small colimits on each side.  The definitions of a monoidal category and of a symmetric monoidal category are given in Definition \ref{def:monoidal-category} and Definition \ref{def:symmetric-monoidal-category}, respectively.  The commutation of the monoidal product with small colimits on each side means that, for each functor $F : \C \to \M$ from a small category $\C$ and each object $Y\in\M$, the natural morphisms
\[\colimover{x\in \C}\, (Fx \otimes Y) \iso \Bigl(\colimover{x\in\C}\, Fx\Bigr)\otimes Y \andspace \colimover{x\in \C}\, (Y \otimes Fx) \iso Y \otimes \Bigl(\colimover{x\in\C}\, Fx\Bigr)\] are isomorphisms in $\M$.  

Here are some basic symmetric monoidal categories that satisfy our assumptions above. 
\begin{itemize}
\item \label{notation:set} $(\Set, \times, *)$ : category of sets.\index{set}
\item \label{notation:chaink} $(\Chaink, \otimes, \fieldk)$ : category of chain complexes over a field $\fieldk$ \cite{weibel}.\label{chain complex}
\item \label{notation:chau} $(\CHau, \times, *)$ : category of compactly generated weak Hausdorff spaces \cite{gz,may-concise}.\index{space}
\item \label{notation:sset} $(\Sset, \times, *)$ : category of simplicial sets \cite{may67}.\index{simplicial set}
\item \label{notation:cat} $(\Cat, \times, *)$ : category of small categories.
\end{itemize}

The symmetric monoidal category $\M$ will serve as our ambient category in which various kinds of operads and constructions are defined.  Just for the definition of a $\G$-operad as in Proposition \ref{def:g-operad-generating} or Proposition \ref{prop:g-operad-compi}, we actually only need $\M$ to be a symmetric monoidal category.  The other assumptions on $\M$--namely, the existence of small (co)limits and the commutation of the monoidal product with small colimits on each side--allow us to give other descriptions of $\G$-operads and algebras over them.  

From Part \ref{part:operads-group-eq} to Part \ref{part:infinity-group-operads}, the \index{Coherence Theorem}Coherence Theorem for (symmetric) monoidal categories \cite{maclane} (XI.3) will be used tacitly to avoid parentheses for iterated monoidal products and the associativity isomorphism.  The exact statements of these Coherence Theorems are Theorem \ref{maclane-thm} and Theorem \ref{symmetric-monoidal-strictification}.  With this in mind, for objects $x, x_1,\ldots,x_n$ in a monoidal category $\M$, we sometimes use the notations\label{notation:tensors} \[x^{\otimes n} = \overbrace{x \otimes \cdots \otimes x}^{\text{$n$ copies of $x$}} \andspace \bigotimes_{i=1}^n x_i = x_1 \otimes \cdots \otimes x_n\] for multiple monoidal products with $x^{\otimes 0} = \tensorunit$, and similarly for (co)products.

A \emph{monoid}\index{monoid} in a monoidal category $(\C,\otimes,\tensorunit)$  \cite{maclane} (VII.3) is a triple\label{notation:monoid} $(X,\mu,\operadunit)$ with:
\begin{itemize}
\item $X$ an object in $\C$;
\item $\mu : X \otimes X \to X$ a morphism, called the \emph{multiplication};
\item $\operadunit : \tensorunit \to X$ a morphism, called the \emph{unit}.
\end{itemize}
This data is required to make the following associativity and unity diagrams commutative.
\[\nicexy{X\otimes X \otimes X \ar[d]_-{\mu\otimes \Id_X} \ar[r]^-{\Id_X\otimes \mu} & X \otimes X \ar[d]^-{\mu}\\ X \otimes X \ar[r]^-{\mu} & X}\qquad
\nicexy{\tensorunit \otimes X \ar[r]^-{\operadunit \otimes \Id_X} \ar[d]_-{\lambda}^-{\cong} & X \otimes X \ar[d]^-{\mu} & X \otimes \tensorunit \ar[l]_-{\Id_X \otimes \operadunit} \ar[d]^-{\rho}_-{\cong}\\ X \ar@{=}[r] & X & X \ar@{=}[l]}\]
A morphism of monoids \[f : (X,\mu^X,\operadunit^X) \to (Y,\mu^Y,\operadunit^Y)\] is a morphism $f : X \to Y$ in $\C$ that preserves the multiplications and the units in the obvious sense.  The category of monoids in a monoidal category $\C$ is denoted by\label{notation:monoid-category} $\Mon(\C)$ or $\Mon(\C,\otimes,\tensorunit)$.  Notice that to define monoids in $\C$, it does \emph{not} need to be a symmetric monoidal category.  This is important because all the variants of operads in this work are defined as monoids in some monoidal categories that are not symmetric.

\subsection{Adjoint Lifting Theorem}\label{subsec:adjoint-lifting}

For an adjunction \[\nicexy{\C \ar@<2pt>[r]^-{L} & \D \ar@<2pt>[l]^-{R}}\] with left adjoint $L$ and right adjoint $R$, we always display the left adjoint on top, pointing to the right.  If an adjunction is displayed vertically, then the left adjoint is written on the left-hand side.  We also write $L \dashv R$ or $(L,R)$ for such an adjoint pair of functors.

Suppose:
\begin{enumerate}[label=(\roman*)]
\item $L \dashv R$ is an adjoint pair as above.  
\item $T$ is a monad in the category $\C$, the precise formulation of which is recalled in Definition \ref{def:monad}.
\item $V$ is a monad in the category $\D$.
\item $R' : \alg(V) \to \alg(T)$ is a functor from the category of $V$-algebras to the category of $T$-algebras.
\item $\alg(V)$ has all coequalizers.
\item The diagram \[\nicexy{\alg(T) \ar[d]_-{U} & \alg(V) \ar[l]_-{R'} \ar[d]^-{U}\\ \C & \D \ar[l]_-{R}}\] is commutative, in which both vertical functors $U$ are the forgetful functors.
\end{enumerate}
Then the \index{Adjoint Lifting Theorem}\emph{Adjoint Lifting Theorem}--see \cite{borceux2} Theorem 4.5.6--says that $R'$ also admits a left adjoint \[L' : \alg(T) \to \alg(V).\]  Moreover, it follows from (i) the equality $UR' = RU$ of right adjoints and (ii) the uniqueness of left adjoints, that the diagram of left adjoints
\[\nicexy{\alg(T) \ar[r]^-{L'} & \alg(V)\\
\C \ar[u]^-{T} \ar[r]^-{L} & \D \ar[u]_-{V}}\] is commutative up to a natural isomorphism.  Although we will not use it in this book, we point out that there is a homotopical version of the Adjoint Lifting Theorem involving Quillen equivalences between model categories \cite{white-yau-halt}.

\part{Operads with Group Equivariance}\label{part:operads-group-eq}

\chapter{Planar Operads}\label{ch:planar-operad}

The purpose of this chapter is to discuss colored planar operads, their algebras, and some examples. Colored planar operads in the category of sets were originally defined by Lambek \cite{lambek} with the name \emph{multicategories}.  Colored symmetric operads are discussed in the next chapter.  A more gentle introduction to colored planar and symmetric operads is the book \cite{yau-operad}.  Other useful references for one-colored symmetric operads include the books by Fresse \cite{fresse-book}, Kriz-May \cite{krizmay}, Loday-Vallette \cite{lv}, and Markl-Shnider-Stasheff \cite{mss}, as well as the articles by Markl \cite{markl08}, May \cite{may97a,may97b}, and Vallette \cite{vallette}.

In Section \ref{sec:planar-operad-as-monoid} we define $\colorc$-colored planar operads as monoids with respect to a monoidal product called the $\colorc$-colored planar circle product.  In Section  \ref{sec:planar-generating} we provide more explicit descriptions of colored planar operads in terms of generating operations and generating axioms.  In Section \ref{sec:planar-operad-algebra} we discuss algebras over a colored planar operad, first conceptually and then in more concrete terms.  Section \ref{sec:planar-operad-example} contains some examples of colored planar operads and algebras.

Throughout Part \ref{part:operads-group-eq} to Part \ref{part:infinity-group-operads} of this book, $(\M,\otimes,\tensorunit)$ denotes a complete and cocomplete symmetric monoidal category whose monoidal product $\otimes$ commutes with small colimits on each side.  A simple example is the category $\Set$ of sets with Cartesian product as the monoidal product and a one-point set as the monoidal unit.

\section{Planar Operads as Monoids}\label{sec:planar-operad-as-monoid}

Colored planar operads are monoids in a monoidal category of colored planar sequences, which we first define.

\begin{definition}\label{def:colored-sequences}
Fix a set\label{notation:color-set} $\colorc$, whose elements are called \emph{colors}.
\begin{enumerate}
\item A \emph{$\colorc$-profile}\index{profile} is a finite, possibly empty, sequence\label{notation:uc} \[\uc=(c_1,\ldots,c_n)\] with each $c_i \in \colorc$.  In this case, we write\label{notation:lengthuc} $|\uc|=n$, called the \emph{length} of $\uc$.
\item The empty $\colorc$-profile\index{empty profile} is denoted by\label{notation:empty-prof} $\varnothing$.
\item The set\index{set of profiles} of all $\colorc$-profiles is denoted by\label{notation:profc} $\Profc$, which is also regarded as a discrete category.
\item The objects of the diagram category\label{notation:pseqcm} \[\pseqcm = \M^{\Profcc}\] are called \index{planar sequence}\emph{$\colorc$-colored planar sequences} in $\M$.  For an object $X$ in $\pseqcm$, we write\label{notation:xduc} \[X\duc \in \M\] for the value of $X$ at\label{notation:duc} $(\uc;d) = \duc \in \Profcc$, and call it an \emph{$m$-ary entry of $X$} if $|\uc|=m$.  We call $\uc$ the \index{input profile}\emph{input profile}, $c_i$ the \emph{$i$th input color}, and $d$ the \index{output color}\emph{output color}.
\item An object in the product category $\prod_{\colorc} \M = \M^{\colorc}$\label{notation:mtoc} is called a \index{colored object}\emph{$\colorc$-colored object in $\M$}, and similarly for a morphism of $\colorc$-colored objects.  A $\colorc$-colored object $X$ is also written as\label{notation:xsubc} $\{X_c\}$ with $X_c \in \M$ for each color $c \in \colorc$.
\item A $\colorc$-colored object $\{X_c\}_{c\in\colorc}$ is also regarded as a $\colorc$-colored planar sequence concentrated in $0$-ary entries:
\begin{equation}\label{colored-object-sm}
X\duc= \begin{cases}X_d & \text{ if $\uc=\varnothing$},\\ \varnothing^{\M} & \text{ if $\uc\not=\varnothing$}.\end{cases}
\end{equation}
\end{enumerate}\end{definition}

Next we define the monoidal product on the category of $\colorc$-colored planar sequences.  The following colored planar circle product is the colored planar version of the circle product defined by Kelly \cite{kelly}.

\begin{definition}\label{def:planar-circle-product}
Suppose $X,Y  \in \pseqcm$.
\begin{enumerate}
\item For each $\uc = (c_1,\ldots,c_m) \in \Profc$, define the object $Y^{\uc} \in \M^{\Profc}$ entrywise as the coend
\begin{equation}\label{ytensorc}
Y^{\uc}(\ub) = \int^{\{\ua_j\}\in\overset{m}{\prodover{j=1}} \Profc} \Profc\bigl(\ua_1,\ldots,\ua_m;\ub\bigr) \cdot \left[\bigotimes_{j=1}^m Y\cjuaj\right] \in \M
\end{equation}
for $\ub \in \Profc$, in which $(\ua_1,\ldots,\ua_m) \in \Profc$ is the concatenation.
\item The \index{circle product!planar}\index{planar circle product}\emph{$\colorc$-colored planar circle product} \[X \circp Y \in \pseqcm\] is defined entrywise as the coend
\begin{equation}\label{planar-circle-product}
(X \circp Y)\dub = \int^{\uc \in \Profc} X\duc \otimes Y^{\uc}(\ub)\in \M
\end{equation}
for $\dub \in \Profcc$.
\item Define the object\index{unit planar operad}\index{planar operad!unit} $\I \in \pseqcm$ by
\begin{equation}\label{unit-operad}
\I\duc = \begin{cases} \tensorunit & \text{ if $\uc=(d)$},\\
\emptyset^{\M}  & \text{ otherwise}\end{cases}
\end{equation}
for $\duc \in \Profcc$.
\end{enumerate}
\end{definition}

\begin{remark}\label{rk:planar-ycb-coprod}
Since $\Profc$ is a discrete category, each of its hom sets is either empty or contains a single element.  Therefore, the coend in \eqref{ytensorc} is the coproduct
\[Y^{\uc}(\ub) \cong \coprod_{(\ua_1,\ldots,\ua_m)=\ub}\, \bigotimes_{j=1}^m Y\cjuaj\]
indexed by the set of $\colorc$-profiles $\{\ua_j\}\in\prod_{j=1}^m \Profc$ whose concatenation $(\ua_1,\ldots,\ua_m)$ is $\ub$.  Similarly, a typical entry of the $\colorc$-colored planar circle product \eqref{planar-circle-product} is a coproduct
\[(X \circp Y)\dub \cong \coprod_{\uc\in\Profc} X\duc \otimes Y^{\uc}(\ub)\] indexed by the set of all $\colorc$-profiles.\dqed
\end{remark}

\begin{proposition}\label{planar-circle-product-monoidal}
$\bigl(\pseqcm, \circp, \I\bigr)$ is a \index{monoidal category!of planar sequences}monoidal category.
\end{proposition}

\begin{proof}
Suppose $X,Y,Z \in \pseqcm$.  To exhibit the associativity isomorphism, first note that for $\uc=(c_1,\ldots,c_m), \ub \in \Profc$, there exist canonical isomorphisms:
\begin{equation}\label{ycircz}
\begin{split}(Y \circp Z)^{\uc}(\ub) 
&= \int^{\ua_1,\ldots,\ua_m} \Profc\bigl(\ua_1,\ldots,\ua_m;\ub\bigr) \cdot \left[\bigotimes_{j=1}^m (Y\circ Z)\cjuaj\right]\\
&\cong \int^{\ua_1,\ldots,\ua_m} \int^{\ud_1,\ldots,\ud_m} \Profc\bigl(\ua;\ub\bigr) \cdot \bigotimes_{j=1}^m \Bigl[Y\cjudj \otimes Z^{\ud_j}(\ua_j)\Bigr]\\
&\cong \int^{\ud_1,\ldots,\ud_m} \left[\bigotimes_{j=1}^m Y\cjudj\right] \otimes Z^{\ud}(\ub)\\
&\cong \int^{\ud_1,\ldots,\ud_m} \int^{\ue} \Profc(\ud;\ue) \cdot \left[\bigotimes_{j=1}^m Y\cjudj\right] \otimes Z^{\ue}(\ub)\\
&\cong \int^{\ue} Y^{\uc}(\ue) \otimes Z^{\ue}(\ub)\end{split}
\end{equation}
in which $\ua=(\ua_1,\ldots,\ua_m)$ and $\ud=(\ud_1,\ldots,\ud_m)$.  Now we have the equalities and canonical isomorphism
\[\begin{split}
\bigl((X \circp Y) \circp Z\bigr)\dua 
&= \int^{\uc} (X \circp Y)\duc \otimes Z^{\uc}(\ua)\\
&= \int^{\uc} \int^{\ub} X\dub \otimes Y^{\ub}(\uc) \otimes Z^{\uc}(\ua)\\
&\cong \int^{\ub} X\dub \otimes (Y \circp Z)^{\ub}(\ua)\\
&= \bigl(X \circp (Y \circp Z)\bigr)\dua\end{split}\]
for $\dua \in \Profcc$, in which the isomorphism uses \eqref{ycircz}.  The rest of the axioms of a monoidal category are straightforward to check.
\end{proof}

\begin{definition}\label{def:planar-operad}
The category of \index{operad!planar}\index{planar operad}\emph{$\colorc$-colored planar operads in $\M$} is the category \[\Poperadcm = \Mon\bigl(\pseqcm, \circp, \I\bigr)\] of monoids in the monoidal category $\bigl(\pseqcm,\circp,\I\bigr)$.  If $\colorc$ has $n<\infty$ elements, we also refer to objects in $\Poperadcm$ as \emph{$n$-colored planar operads in $\M$}.
\end{definition}

Unfolding the above definition, we may express a $\colorc$-colored planar operad as follows.

\begin{corollary}\label{planar-operad-monoid}
A $\colorc$-colored planar operad in $\M$ is exactly a triple $(\O,\mu,\varepsilon)$ consisting of
\begin{enumerate}[label=(\roman*)]
\item an object $\O \in \pseqcm$, 
\item an operadic multiplication morphism $\mu : \O \circp \O \to \O \in \pseqcm$, and 
\item an operadic unit $\varepsilon : \I \to \O \in \pseqcm$ 
\end{enumerate}
such that the associativity and unity diagrams
\begin{equation}\label{planar-operad-monoid-axioms}
\nicexy{\O \circp \O \circp \O \ar[d]_-{(\mu,\Id_{\O})} \ar[r]^-{(\Id_{\O},\mu)} & \O \circp \O \ar[d]^-{\mu}\\ \O \circp \O \ar[r]^-{\mu} & \O}\qquad
\nicexy{\I \circp \O \ar[d]_-{\cong} \ar[r]^-{(\varepsilon, \Id_{\O})} & \O \circp \O \ar[d]^-{\mu} & \O \circp \I \ar[l]_-{(\Id_{\O},\varepsilon)} \ar[d]^-{\cong}\\ \O \ar@{=}[r] & \O & \O \ar@{=}[l]}
\end{equation}
in $\pseqcm$ are commutative.  A morphism of $\colorc$-colored planar operads is a morphism of the underlying $\colorc$-colored planar sequences that is compatible with the operadic multiplications and the operadic units.
\end{corollary}

\begin{notation}
If $\O$ is a $1$-colored planar operad with color set $\{*\}$, then we write\label{notation:oofn} \[\O(n) = \O\sbinom{*}{*,\ldots,*}\] for $\sbinom{*}{*,\ldots,*} \in \Prof(\{*\}) \times \{*\}$ in which the input profile has length $n$.  
\end{notation}

\section{Coherence for Planar Operads}\label{sec:planar-generating}

The purpose of this section is to unpack the definition of a colored planar operad, expressing it in terms of a few generating operations and generating axioms as follows.

\begin{definition}\label{def:planar-operad-generating}
A \emph{$\colorc$-colored planar operad} in $(\M,\otimes,\tensorunit)$ is a triple $\bigl(\O, \gamma, \operadunit\bigr)$
consisting of the following data.
\begin{enumerate}[label=(\roman*)]
\item $\O \in \pseqcm$.
\item For $\bigl(\uc = (c_1, \ldots , c_n); d\bigr) \in \Profcc$ with $n \geq 1$,  $\ub_j \in \Profc$ for $1 \leq j \leq n$, and $\ub = (\ub_1,\ldots,\ub_n)$ their concatenation, 
it is equipped with an \index{operadic composition}\emph{operadic composition} \label{notation:operadic-composition}
\begin{equation}\label{operadic-composition}
\nicexy{\O\duc \otimes \bigotimes\limits_{j=1}^n \O\cjubj \ar[r]^-{\gamma} & \O\dub \in \M.}
\end{equation}
\item For each $c \in \colorc$, it is equipped with a \index{operad!unit}\emph{$c$-colored unit}\label{notation:colored-unit}
\begin{equation}\label{c-colored-unit}
\nicexy{\tensorunit \ar[r]^-{\operadunit_c} & \O\cc \in \M.}
\end{equation}
\end{enumerate}
This data is required to satisfy the following associativity and unity axioms.
\begin{description}
\item[Associativity]
Suppose that:
\begin{itemize}
\item in \eqref{operadic-composition} $\ub_j = \bigl(b^j_1, \ldots , b^j_{k_j}\bigr) \in \Profc$ with at least one $k_j > 0$;
\item $\ua^j_i \in \Profc$ for each $1 \leq j \leq n$ and $1 \leq i \leq k_j$;
\item for each $1 \leq j \leq n$, 
\[\ua_j = \begin{cases}\bigl(\ua^j_1, \ldots , \ua^j_{k_j}\bigr) & \text{if $k_j > 0$},\\
\varnothing & \text{if $k_j = 0$}\end{cases}\]
with $\ua = (\ua_1,\ldots , \ua_n)$ their concatenation.
\end{itemize}
Then the \index{operad!associativity}\index{associativity!operad}\emph{associativity diagram}
\begin{equation}\label{operad-associativity}
\nicexy{\O\duc \otimes 
\left[\bigotimes\limits_{j=1}^n \O\cjubj\right] \otimes \bigotimes\limits_{j=1}^n 
\bigotimes\limits_{i=1}^{k_j} \O\bjiuaji \ar[r]^-{(\gamma, \Id)} \ar[d]_{\text{permute}}^-{\cong} & \O\dub \otimes \bigotimes\limits_{j=1}^{n}\bigotimes\limits_{i=1}^{k_j} \O\bjiuaji \ar[dd]^{\gamma} \\ \O\duc \otimes \bigotimes\limits_{j=1}^n \left[\O\cjubj 
\otimes \bigotimes\limits_{i=1}^{k_j} \O\bjiuaji\right] \ar[d]_{(\Id, \otimes_j \gamma)} &\\
\O\duc \otimes \bigotimes\limits_{j=1}^n \O\cjuaj \ar[r]^-{\gamma} & \O\dua}
\end{equation}
in $\M$ is commutative.
\item[Unity]
Suppose $d \in \colorc$.
\begin{enumerate}
\item For each $\uc = (c_1,\ldots,c_n) \in \Profc$ with $n \geq 1$, the \index{operad!unity}\index{unity!operad}\emph{right unity diagram}
\begin{equation}\label{right-unity}
\nicexy{\O\duc \otimes \tensorunit^{\otimes n}\ar[d]_-{(\Id, \otimes \operadunit_{c_j})} \ar[r]^-{\cong} & \O\duc \ar[d]^-{=} \\
\O\duc \otimes \bigotimes\limits_{j=1}^n \O\cjcj\ar[r]^-{\gamma}&\O\duc}
\end{equation}
in $\M$ is commutative.
\item For each $\ub \in \Profc$ the \emph{left unity diagram}
\begin{equation}\label{left-unity}
\nicexy{\tensorunit \otimes \O\dub \ar[d]_-{(\operadunit_d, \Id)} \ar[r]^-{\cong}& \O\dub \ar[d]^-{=}\\ \O\dd \otimes \O\dub\ar[r]^-{\gamma} & \O\dub}
\end{equation}
in $\M$ is commutative.
\end{enumerate}
\end{description}
A morphism of $\colorc$-colored planar operads is a morphism of the underlying $\colorc$-colored planar sequences that is compatible with the operadic compositions and the colored units in the obvious sense.
\end{definition}

\begin{proposition}\label{prop:planar-operad-same-def}
The\index{coherence!planar operad}\index{planar operad!coherence} definitions of a $\colorc$-colored planar operad in Definition \ref{def:planar-operad} and in Definition \ref{def:planar-operad-generating} are equivalent.
\end{proposition}

\begin{proof}
In Corollary \ref{planar-operad-monoid} the domain of the multiplication $\mu$ is $\O \circp \O$, which has entries
\begin{equation}\label{ocompo}\begin{split}
(\O\circp\O)\dub &= \int^{\uc} \O\duc \otimes \O^{\uc}(\ub)\\
&\cong \int^{\uc, \ua_1, \ldots, \ua_m} \Profc(\ua;\ub) \cdot \O\duc \otimes \bigotimes_{j=1}^m \O\cjuaj\end{split}
\end{equation}
for $\dub \in \Profcc$, where $\ua=(\ua_1,\ldots,\ua_m)$ is the concatenation.  Observe that $\Profc(\ua;\ub)$ is empty unless $\ua=\ub$, i.e., the concatenation of the $\ua_j$'s is $\ub$.  So the multiplication $\mu : \O\circp\O \to \O$ yields the entrywise operadic composition $\gamma$ \eqref{operadic-composition}.  The associativity diagram in \eqref{planar-operad-monoid-axioms} corresponds to the associativity diagram \eqref{operad-associativity}.  For each $c \in \colorc$, the $c$-colored unit $\operadunit_c$ in \eqref{c-colored-unit} corresponds to the $\ccsingle$-entry of the unit morphism $\varepsilon : \I \to \O$ in Corollary \ref{planar-operad-monoid}.  The unity diagram in \eqref{planar-operad-monoid-axioms} corresponds to the right unity diagram \eqref{right-unity} and the left unity diagram \eqref{left-unity}.
\end{proof}

Instead of the operadic composition $\gamma$, a colored planar operad can also be described in terms of binary operations as follows.

\begin{definition}\label{def:compi}
Suppose $\bigl(\uc = (c_1, \ldots , c_n); d\bigr) \in \Profcc$ with $n \geq 1$, $\ub \in \Profc$, and $1 \leq i \leq n$.  Define the\index{profile!composition} $\colorc$-profile
\[\uc \compi \ub = \bigl(\underbrace{c_1,\ldots,c_{i-1}}_{\varnothing~\mathrm{if}~i=1},\ub,\underbrace{c_{i+1},\ldots,c_n}_{\varnothing~\mathrm{if}~i=n}\bigr),\] which is obtained from $\uc$ by replacing the entry $c_i$ with $\ub$.
\end{definition}

The following is \cite{yau-operad} Definition 16.6.1.

\begin{definition}\label{def:planar-operad-compi}
A \emph{$\colorc$-colored planar operad} in $(\M,\otimes,\tensorunit)$ is a triple $(\O,\comp,\operadunit)$ consisting of the following data.
\begin{enumerate}[label=(\roman*)]
\item $\O\in \pseqcm$.
\item For $\bigl(\uc = (c_1, \ldots , c_n);d\bigr) \in \Profcc$ with $n \geq 1$, $1 \leq i \leq n$, and $\ub \in \Profc$, it is equipped with a\index{compi  composition@$\compi$-composition} morphism
\begin{equation}\label{operadic-compi}
\nicexy{\O\duc \otimes \O\ciub\ar[r]^-{\compi} & \O\sbinom{d}{\uc\compi\ub} \in \M}
\end{equation}
called the \label{notation:compi-operad}\emph{$\compi$-composition}. 
\item For each color $c \in \colorc$, it is equipped with a $c$-colored unit
\[\nicexy{\tensorunit \ar[r]^-{\operadunit_c} & \O\cc\in \M.}\]
\end{enumerate}
This data is required to satisfy the following associativity and unity axioms.  Suppose $d \in \colorc$, $\uc = (c_1, \ldots , c_n) \in \Profc$, $\ub \in \Profc$ with length $|\ub| = m$, and $\ua \in \Profc$ with length $|\ua| = l$.
\begin{description}
\item[Associativity]
There are two associativity axioms.
\begin{enumerate}
\item Suppose $n \geq 2$ and $1 \leq i < j \leq n$.  Then the \index{horizontal associativity}\index{operad!horizontal associativity}\emph{horizontal associativity diagram} in $\M$
\begin{equation}\label{compi-associativity}
\nicexy{\O\duc \otimes \O\sbinom{c_i}{\ua} \otimes \O\cjub \ar[d]_-{\mathrm{permute}}^-{\cong} 
\ar[r]^-{(\compi, \Id)} & \O\sbinom{d}{\uc\compi\ua} \otimes \O\cjub\ar[dd]^-{\comp_{j-1+l}}\\
\O\duc \otimes \O\cjub \otimes \O\sbinom{c_i}{\ua}\ar[d]_-{(\comp_j,\Id)} &\\
\O\sbinom{d}{\uc\comp_j\ub}\otimes \O\sbinom{c_i}{\ua}\ar[r]^-{\compi} & \O\sbinom{d}{(\uc\comp_j\ub)\compi\ua} = \O\sbinom{d}{(\uc\compi\ua)\comp_{j-1+l}\,\ub}}
\end{equation}
is commutative.
\item Suppose $n,m \geq 1$, $1 \leq i \leq n$, and $1 \leq j \leq m$.  Then the \index{vertical associativity}\index{operad!vertical associativity}\emph{vertical associativity diagram} in $\M$
\begin{equation}\label{compi-associativity-two}
\nicexy{\O\duc \otimes \O\ciub \otimes \O\sbinom{b_j}{\ua} \ar[d]_-{(\compi,\Id)} \ar[r]^-{(\Id,\comp_j)} &
\O\duc \otimes \O\sbinom{c_i}{\ub\comp_j\ua} \ar[d]^-{\compi}\\
\O\sbinom{d}{\uc\compi\ub} \otimes \O\sbinom{b_j}{\ua} \ar[r]^-{\comp_{i-1+j}} & 
\O\sbinom{d}{(\uc\compi\ub)\comp_{i-1+j}\,\ua} = \O\sbinom{d}{\uc\compi(\ub\comp_j\ua)}}
\end{equation}
is commutative.
\end{enumerate}
\item[Unity]
There are two unity axioms.
\begin{enumerate}
\item The \emph{left unity diagram} in $\M$\index{operad!unity}
\begin{equation}\label{compi-left-unity}
\nicexy{\tensorunit \otimes \O\duc \ar[dr]_-{\cong} \ar[r]^-{(\operadunit_d,\Id)}&
\O\dd \otimes \O\duc \ar[d]^-{\comp_1}\\ & \O\duc}
\end{equation}
is commutative.
\item If $n \geq 1$ and $1 \leq i \leq n$, then the \emph{right unity diagram} in $\M$
\begin{equation}\label{compi-right-unity}
\nicexy{\O\duc \otimes \tensorunit \ar[dr]_-{\cong} \ar[r]^-{(\Id,\operadunit_{c_i})}& 
\O\duc \otimes \O\cici\ar[d]^-{\compi}\\ & \O\duc}
\end{equation}
is commutative.
\end{enumerate}
\end{description}
A morphism of $\colorc$-colored planar operads is a morphism of the underlying $\colorc$-colored planar sequences that is compatible with the $\compi$-compositions and the colored units in the obvious sense.
\end{definition}

\begin{proposition}\label{prop:planar-operad-defs-equivalent}
The definition of a $\colorc$-colored planar operad in Definition \ref{def:planar-operad-generating} and in Definition \ref{def:planar-operad-compi} are equivalent.
\end{proposition}

\begin{proof}
The correspondence of structures is as follows.  Given a $\colorc$-colored planar operad $(\O,\gamma,\operadunit)$ in the sense of Definition \ref{def:planar-operad-generating}, the associated $\compi$-composition is the following composite:
\begin{equation}\label{compi-def}
\nicexy{\O\duc\otimes\O\ciub \ar[d]_-{\cong} \ar[r]^-{\compi} & \O\sbinom{d}{\uc\compi\ub}\\
\O\duc \otimes \tensorunit^{\otimes i-1} \otimes \O\ciub \otimes \tensorunit^{\otimes n-i} \ar[r]^-{\{\operadunit_{c_j}\}} &
\O\duc \otimes \Bigl[\bigotimes\limits_{j=1}^{i-1}\O\cjcj\Bigr] \otimes \O\ciub \otimes \Bigl[\bigotimes\limits_{j=i+1}^{n}\O\cjcj\Bigr]. \ar[u]_-{\gamma}}
\end{equation}
The bottom horizontal morphism is the monoidal product of the colored units $\operadunit_{c_j}$ for $1 \leq j\not=i \leq n$ with the identity morphisms of $\O\duc$ and $\O\ciub$.

Conversely, given a $\colorc$-colored operad $(\O,\comp,\operadunit)$ in the sense of Definition \ref{def:planar-operad-compi}, the operadic composition $\gamma$ is recovered as the composition\index{operadic composition!in terms of $\compi$}
\begin{equation}\label{compi-to-gamma}
\nicexy@C+.7cm{\O\duc \otimes\bigotimes\limits_{j=1}^n \O\cjubj \ar[r]^-{\gamma} \ar[d]_-{(\comp_1,\Id)} & \O\dub\\
\O\sbinom{d}{\uc\comp_1\ub_1}\otimes\bigotimes\limits_{j=2}^n\O\cjubj \ar[d]_-{(\comp_{k_1+1},\Id)} 
& \O\sbinom{d}{((\uc\comp_1\ub_1)\cdots)\comp_{k_1+\cdots+k_{n-1}+1}\ub_n} \ar@{=}[u]\\
\cdots \ar[r]^-{(\comp_{k_1+\cdots+k_{n-2}+1},\Id)} & \O\sbinom{d}{((\uc\comp_1\ub_1)\cdots)\comp_{k_1+\cdots+k_{n-2}+1}\ub_{n-1}} \otimes \O\sbinom{c_n}{\ub_n} \ar[u]_-{\comp_{k_1+\cdots+k_{n-1}+1}}}
\end{equation}
in which $k_j = |\ub_j|$.  
\end{proof}

\section{Algebras}\label{sec:planar-operad-algebra}

The purpose of this section is to define algebras over a colored planar operad.  The following notation will be useful.

\begin{notation}\label{not:x-sub-c}
For a\index{colored object} $\colorc$-colored object $X=\{X_c\}_{c\in\colorc}$ in $\M$ and $\uc=(c_1,\ldots,c_m) \in \Profc$, we will write \[X_{\uc} = X_{c_1} \otimes \cdots \otimes X_{c_m},\] which is the monoidal unit $\tensorunit$ if $\uc$ is the empty profile, i.e., $X_{\varnothing} = \tensorunit$.
\end{notation}

Recall the concepts of a monad and of an algebra over a monad \cite{maclane} VI.

\begin{definition}\label{def:monad}
A \index{monad}\emph{monad} in a category $\C$ is a triple $(T,\mu,\epsilon)$ in which
\begin{enumerate}[label=(\roman*)]
\item $T : \C \to \C$ is a functor and 
\item $\mu : T^2 \to T$, called the \emph{multiplication}, and $\epsilon : \Id_{\C} \to T$, called the \emph{unit}, are natural transformations
\end{enumerate} 
such that the \index{associativity!monad}associativity and \index{unity!monad}unity diagrams
\[\nicexy{T^3 \ar[d]_-{\mu T} \ar[r]^-{T\mu} & T^2 \ar[d]^-{\mu}\\ T^2 \ar[r]^-{\mu} & T}\qquad
\nicexy{\Id_{\C} \circ T \ar@{=}[d] \ar[r]^-{\epsilon T} & T^2 \ar[d]_-{\mu} & T\circ \Id_{\C} \ar[l]_-{T\epsilon} \ar@{=}[d]\\ T \ar@{=}[r] & T & T \ar@{=}[l]}\]
are commutative.
\end{definition}

\begin{definition}\label{def:monad-algebra}
Suppose $(T,\mu,\epsilon)$ is a monad in a category $\C$.  
\begin{enumerate}
\item A \index{algebra!of a monad}\emph{$T$-algebra} is a pair $(X,\lambda)$ with
\begin{enumerate}[label=(\roman*)]
\item $X$ an object in $\C$ and 
\item $\lambda : TX \to X$ a morphism, called the \emph{structure morphism}, 
\end{enumerate}
such that the associativity and unity diagrams
\[\nicexy{T^2X \ar[d]_-{\mu_X} \ar[r]^-{T\lambda} & TX \ar[d]^-{\lambda}\\ TX \ar[r]^-{\lambda} & X}\qquad \nicexy{X \ar[r]^-{\epsilon_X} \ar@{=}[dr] & TX \ar[d]^-{\lambda}\\ & X}\]
are commutative.
\item A \emph{morphism of $T$-algebras} \[f : (X,\lambda^X) \to (Y,\lambda^Y)\] is a morphism $f : X \to Y$ in $\C$ such that the diagram
\[\nicexy{TX \ar[r]^-{Tf} \ar[d]_-{\lambda^X} & TY \ar[d]^-{\lambda^Y}\\ X \ar[r]^-{f} & Y}\]
is commutative.
\item The category of $T$-algebras is denoted by\label{notation:algt} $\alg(T)$.
\end{enumerate}
\end{definition}

\begin{lemma}\label{lem:planar-operad-monad}
Suppose $\O$ is a $\colorc$-colored planar operad in $\M$.  Then it induces a \index{planar operad!induced monad}\index{monad!induced by a planar operad}monad whose functor is \[\O \circp - : \Mtoc \to \Mtoc\] and whose multiplication and unit are induced by those of $\O$ as in Corollary \ref{planar-operad-monoid}.
\end{lemma}

\begin{proof}
Suppose $Y=\{Y_c\}_{c\in\colorc}$ is a $\colorc$-colored object in $\M$, regarded as a $\colorc$-colored planar sequence as in \eqref{colored-object-sm}.  In \eqref{ytensorc} we have 
\[Y^{\uc}(\ub)= \begin{cases} Y_{\uc} & \text{if $\ub=\varnothing$},\\
\varnothing^{\M} & \text{if $\ub\not=\varnothing$}\end{cases}\] for $\ub,\uc\in\Profc$.  Substituting this into the definition \eqref{planar-circle-product} of the $\colorc$-colored planar circle product, we obtain
\begin{equation}\label{ocompyentries}
(\O\circp Y)\dub = \begin{cases} \dint^{\uc\in\Profc} \O\duc \otimes Y_{\uc} \cong \coprodover{\uc \in \Profc} \O\duc \otimes Y_{\uc} & \text{ if $\ub=\varnothing$},\\ \varnothing^{\M} & \text{ if $\ub\not=\varnothing$.}\end{cases}
\end{equation}
The functor \[\O \circp - : \pseqcm \to \pseqcm\] thus restricts to a functor $\Mtoc \to \Mtoc$.  Corollary \ref{planar-operad-monoid} now implies that $\O\circp -$ is a monad in $\Mtoc$.
\end{proof}

\begin{definition}\label{def:planar-operad-algebra}
Suppose $\O$ is a $\colorc$-colored planar operad in $\M$.  The category $\algmo$ \index{planar operad!algebra}\index{algebra!of a planar operad}of \emph{$\O$-algebras} is defined as the category of $(\O\circp -)$-algebras for the monad $\O\circp -$ in $\Mtoc$.
\end{definition}

Unwrapping the above definition, we obtain the following detailed description of an algebra over a colored planar operad.

\begin{definition}\label{def:planar-operad-algebra-generating}
Suppose $(\O,\gamma,\operadunit)$ is a $\colorc$-colored planar operad in $\M$.  An \emph{$\O$-algebra} is a pair $(X,\lambda)$ consisting of 
\begin{enumerate}[label=(\roman*)]
\item a $\colorc$-colored object $X=\{X_c\}_{c\in\colorc}$ and 
\item an \emph{$\O$-action structure morphism}
\begin{equation}\label{planar-operad-algebra-action}
\nicexy{\O \duc \otimes X_{\uc} \ar[r]^-{\lambda}& X_d \in \M}
\end{equation}
for each $\duc \in \Profcc$.  
\end{enumerate}
It is required that the following associativity and unity axioms hold.
\begin{description}
\item[Associativity]
For $\bigl(\uc = (c_1, \ldots , c_n);d\bigr) \in \Profcc$ with $n \geq 1$, $\ub_j \in \Profc$ for $1 \leq j \leq n$, and $\ub = (\ub_1,\ldots,\ub_n)$ their concatenation, the associativity diagram\index{associativity!planar operad  algebra}
\begin{equation}
\label{planar-operad-algebra-associativity}
\nicexy@C+15pt{\O\duc \otimes \left[\bigotimes\limits_{j=1}^n \O\cjubj\right]\otimes X_{\ub} 
\ar[r]^-{(\gamma, \Id)} \ar[d]_-{\mathrm{permute}}^-{\cong}& \O\dub \otimes X_{\ub} \ar[dd]^-{\lambda}\\ \O\duc \otimes \bigotimes\limits_{j=1}^n \left[\O \cjubj \otimes X_{\ub_j}\right] \ar[d]_-{(\Id,\otimes_j \lambda)} &\\ \O \duc \otimes X_{\uc} \ar[r]^-{\lambda}& X_d}
\end{equation}
in $\M$ is commutative.
\item[Unity] For each $c \in \colorc$, the unity diagram\index{unity!planar operad algebra}
\begin{equation}\label{planar-operad-algebra-unity}
\nicexy{\tensorunit \otimes X_c \ar[d]_-{(\operadunit_c, \Id)} \ar[r]^-{\cong}& X_c \ar@{=}[d]\\
\O\cc \otimes X_c \ar[r]^-{\lambda}& X_c}
\end{equation}
in $\M$ is commutative.
\end{description}
A \index{morphism!planar operad algebras}\emph{morphism of $\O$-algebras} \[f : (X,\lambda^X) \to (Y,\lambda^Y)\] is a morphism $f : X \to Y$ of $\colorc$-colored objects in $\M$ such that the diagram
\begin{equation}\label{algebra-map-compatibility}
\nicexy{\O\duc \otimes X_{\uc}\ar[d]_-{\lambda^X} \ar[r]^-{(\Id, \otimes f)}& \O\duc \otimes Y_{\uc} \ar[d]^-{\lambda^Y}\\ X_d \ar[r]^-{f} & Y_d}
\end{equation}
in $\M$ is commutative for all $\duc \in \Profcc$.
\end{definition}

\begin{proposition}\label{prop:planar-operad-algebra-defs}
Suppose\index{coherence!planar operad algebra}  $\O$ is a $\colorc$-colored planar operad in $\M$.  Then the definitions of an $\O$-algebra in Definition \ref{def:planar-operad-algebra} and in Definition \ref{def:planar-operad-algebra-generating} are equivalent.
\end{proposition}

\begin{proof}
For a $\colorc$-colored object $X = \{X_c\}_{c\in\colorc}$ in $\M$, the $\colorc$-colored object $\O \circp X$ has entries 
\[\begin{split}
(\O\circp X)_d &= \int^{\uc\in\Profc} \O\duc \otimes X_{\uc}\\
&\cong \coprod_{\uc\in\Profc} \O\duc \otimes X_{\uc}
\end{split}\] 
for $d \in \colorc$.  So an $(\O\circp -)$-algebra has $\O$-action structure morphisms as in \eqref{planar-operad-algebra-action}.  The associativity axiom \eqref{planar-operad-algebra-associativity} and the unity axiom \eqref{planar-operad-algebra-unity} correspond to those of an $(\O\circp -)$-algebra.
\end{proof}

\section{Examples of Planar Operads}\label{sec:planar-operad-example}

In this section, we provide some examples of planar operads and their algebras.

\begin{example}[Small Enriched Categories]\label{ex:enriched-cat}
Recall\index{enriched category}\index{category!enriched}\index{planar operad!underlying category} that the symbol $\varnothing$ denotes both the initial object in $\M$ and the empty $\colorc$-profile.  Using Definition \ref{def:planar-operad-compi}, a $\colorc$-colored planar operad $(\O,\gamma,\operadunit)$ in $\M$ concentrated in unary entries, i.e., \[\O\duc = \varnothing \ifspace |\uc|\not= 1,\] is precisely a small $\M$-enriched category with object set $\colorc$ in the sense of \cite{borceux2} Definition 6.2.1.  For two objects $c,d$ in $\O$, the hom object with domain $c$ and codomain $d$ is $\O\dc$.  The categorical composition is the operadic composition $\gamma$.  The unit morphism of each object $c$ in $\O$ is given by the $c$-colored unit $\operadunit_c$.  The categorical associativity and unity axioms are exactly those in Definition \ref{def:planar-operad-compi}.\dqed
\end{example}

\begin{example}[Monoids]\label{ex:monoids}
As a special case of Example \ref{ex:enriched-cat}, if $\colorc$ consists of a single element $*$, then a $1$-colored planar operad $(\O,\gamma,\operadunit)$ in $\M$ concentrated in unary entries is precisely a\index{monoid} monoid in $\M$.  In this case, an $\O$-algebra in the sense of Definition \ref{def:planar-operad-algebra-generating} is a pair $(X,\lambda)$ in which \[\lambda : \O(1) \otimes X \to X\] gives a \emph{left action} of the monoid $\O(1)$ on the object $X$ in the sense of \cite{maclane} VII.4.\dqed
\end{example}

\begin{example}[Associative Planar Operad]\label{ex:as-planar-operad}
There is a $1$-colored planar \index{associative operad}operad $\As$ in $(\M,\otimes,\tensorunit)$ with entries \[\As(n)=\tensorunit \forspace n \geq 0.\]  The operadic composition $\gamma$ is the natural isomorphism $\tensorunit^{\otimes n+1} \cong \tensorunit$, and the unit  $\operadunit = \Id_{\tensorunit}$ is the identity morphism of the monoidal unit $\tensorunit$.  The category of $\As$-algebras is exactly the category of monoids in $\M$.\dqed
\end{example}

\begin{example}[Planar Operad for Diagrams]\label{ex:diag-operad}
For a small category $\C$, recall that a \emph{$\C$-diagram in $\M$} is a functor $\C \to \M$, and a morphism between two $\C$-diagrams in $\M$ is a natural transformation.  Suppose $\colorc$ is the object set of $\C$.  There is a $\colorc$-colored planar operad\index{planar operad!diagrams}\index{diagram operad} $\Cdiag$ in $\M$ with entries \[\Cdiag\duc = \begin{cases} \coprodover{\C(c,d)}\tensorunit & \text{ if $\uc=c \in \colorc$},\\ \varnothing & \text{ if $|\uc|\not=1$}\end{cases}\]
for $\duc \in \Profcc$.  Its colored units correspond to the identity morphisms in $\C$.  Its operadic composition $\gamma$ corresponds to the categorical composition in $\C$.  The category of $\Cdiag$-algebras is precisely the category of $\C$-diagrams in $\M$.\dqed
\end{example}

\begin{example}[Magma Planar Operad]\label{ex:mag-planar-operad}
Recall that a\index{magma} \emph{magma} is a triple $(X,\mu,e)$ consisting of
\begin{enumerate}[label=(\roman*)]
\item a set $X$,
\item a \emph{multiplication} $\mu : X^{\times 2} \to X$, and 
\item a \emph{unit} $e \in X$ 
\end{enumerate}
such that \[\mu(x,e) = x = \mu(e,x) \forspace x \in X.\]  A morphism of magmas is a function between the underlying sets that preserves both the multiplications and the units.  There is a $1$-colored planar operad\label{notation:mag} $\Mag$\index{magma!operad}\index{planar operad!magma} in the category $(\Set, \times, *)$ of sets whose category of algebras is exactly the category of magmas.  

To describe the magma operad explicitly, we define \index{parenthesized word}\emph{parenthesized words} in one alphabet $\{*\}$ as follows.  
\begin{itemize}\item A parenthesized word of length $0$ is the empty sequence $\varnothing$.  
\item A parenthesized word of length $1$ is the sequence $*$.  
\item Inductively, for $n \geq 2$, a parenthesized word of length $n$ has the form $(u)(v)$, where $(u)$ and $(v)$ are parenthesized words of lengths $l_{(u)} \geq 1$ and $l_{(v)} \geq 1$, respectively, such that $l_{(u)}+l_{(v)}=n$.  
\end{itemize}
For example, the only parenthesized word of length $2$ is $**$.  The only two parenthesized words of length $3$ are $(**)*$ and $*(**)$.  There are five parenthesized words of length $4$: \[((**)*)*\quad (*(**))*\quad (**)(**)\quad *((**)*) \quad *(*(**)).\]  

For $n \geq 0$, suppose $\Mag(n)$ is the set of all parenthesized words of length $n$.  In terms of Definition \ref{def:planar-operad-compi}, for $1 \leq i \leq n$, $m \geq 0$, $(u) \in \Mag(n)$, and $(v) \in \Mag(m)$, the $\compi$-composition \[(u) \compi (v) \in\Mag(n+m-1)\] is the parenthesized word of length $n+m-1$ given by replacing the $i$th $*$ in $(u)$, counting from left to right, by $(v)$.  For example, we have \[\begin{split}
*((**)*) \comp_3 *(**) &= *\Bigl(\bigl(*(*(**))\bigr)*\Bigr),\\ *((**)*) \comp_3 \varnothing &= *(**).\end{split}\]  The unit in $\Mag(1)$ is the unique parenthesized word of length $1$.  With these structure maps, \[\Mag = \bigl\{\Mag(n)\bigr\}_{n\geq 0}\] is a $1$-colored planar operad in $\Set$.  The category of $\Mag$-algebras is equal to the category of magmas.  Also note that the underlying set of \[\Mag = \coprod_{n\geq 0}\Mag(n)\] is the underlying set of the free magma on one generator $*$.\dqed
\end{example}

\begin{example}[Free Magma on a Monoid]\label{ex:free-magma-monoid}
Generalizing the previous example, each monoid generates a colored planar operad as follows.  There is a free-forgetful adjunction \[\nicexy{\Set \ar@<2pt>[r]^-{F} & \Magma \ar@<2pt>[l]^-{U}}\] with\label{notation:magma} $\Magma$ the category of magmas and the left adjoint $F$ sending each set $X$ to the free magma generated by $X$.  More explicitly, the elements in the free magma $FX$ are the \emph{parenthesized words} in $X$ defined as follows.
\begin{itemize}\item A parenthesized word in $X$ of length $0$ is the empty sequence $\varnothing$.  
\item A parenthesized word in $X$ of length $1$ is an element $x\in X$.  
\item Inductively, for $n \geq 2$, a parenthesized word in $X$ of length $n$ has the form \[(u)(v),\] where $(u)$ and $(v)$ are parenthesized words in $X$ of lengths $l_{(u)} \geq 1$ and $l_{(v)} \geq 1$, respectively, such that $l_{(u)}+l_{(v)}=n$.  
\end{itemize}
The set of parenthesized words in $X$ of length $n$ is denoted by $\Mag^X(n)$.

The parenthesized words in Example \ref{ex:mag-planar-operad} are the parenthesized words in the set $\{*\}$.  Parenthesized words in $X$ of lengths $1$ and $2$ form the sets $X$ and $X^{\times 2}$, respectively.  Parenthesized words in $X$ of length $3$ have the forms $(xy)z$ or $x(yz)$ with $x,y,z\in X$.  Forgetting the parentheses, each parenthesized word $w$ in $X$ of length $n$ yields an underlying sequence \[\sequence(w)\in\Prof(X)\] of length $n$.  For example, \[\sequence\bigl((xy)z\bigr) = (x,y,z).\]

Suppose $(M,\mu^M,1^M)$ is a monoid.  There is an $M$-colored planar operad $\Mag^M$ in $\Set$ defined as follows.
\begin{description}
\item[Entries] For $\duc \in \profmm$, denote by $\mu^M(\uc)\in M$ the image of $\uc$ under the iterated multiplication $M^{\times n}\to M$.   Define the set \[\Mag^M\duc= \begin{cases} \varnothing & \text{if $d\not= \mu^M(\uc)$},\\
\bigl\{w\in\Mag^M(|\uc|) : \sequence(w)=\uc\bigl\} & \text{if $d=\mu^M(\uc)$}.
\end{cases}\]
\item[Colored Units] For each $c\in M$, the $c$-colored unit is the parenthesized word \[(c)\in \Mag^M\ccsingle\] of length $1$.
\item[Composition] The $\compi$-composition \[\nicexy{\Mag^M\duc \times \Mag^M\ciub \ar[r]^-{\compi} & \Mag^M\sbinom{d}{\uc\compi\ub}}\] sends $(w,w')$ to $w\compi w'$, which is obtained from the parenthesized word $w\in \Mag^M\duc$ by replacing the $i$th entry $c_i$ with the parenthesized word $w'\in\Mag^M\ciub$.
\end{description}
The $M$-colored planar operad axioms as in Definition \ref{def:planar-operad-compi} follow by inspection.  When $M$ is the terminal monoid $\{*\}$ with one element, the one-colored planar operad $\Mag^{\{*\}}$ in $\Set$ is precisely the magma operad $\Mag$ in Example \ref{ex:mag-planar-operad}.\dqed 
\end{example}

\chapter{Symmetric Operads}\label{ch:symmetric-operad}

In this chapter, we discuss colored symmetric operads, their algebras, and some examples.  One-colored symmetric operads were first introduced by May \cite{may} in the topological setting.  Colored symmetric operads are also known as symmetric multicategories.  A basic introduction to colored symmetric operads is the book \cite{yau-operad}.

In Section \ref{sec:symmetric_operad} we define colored symmetric operads as monoids with respect to the colored symmetric circle product on the category of colored symmetric sequences.  In Section \ref{sec:symmetric-operad-generating} we provide more explicit descriptions of colored symmetric operads in terms of generating operations and generating axioms.  In Section \ref{sec:symmetric-operad-algebra} we discuss algebras over colored operads and some basic examples.  All the assertions in these sections follow from direct inspection and are minor modifications of the planar case.  We will, therefore, omit the proofs.

The remaining sections contain examples of symmetric operads from different fields.  In Section \ref{sec:symmetric-operad-little-cube} we discuss the little $n$-cube operad, the little $n$-disc operad, and the framed little $n$-disc operad, all of which are important in homotopy theory.  In Section \ref{sec:symmetric-operad-graph-operation} we discuss the operad of graph operations in non-commutative probability theory.  In Section \ref{sec:symmetric-operad-phylogenetic} we discuss the phylogenetic operad, which is an organizational tool in genetics.  In Section \ref{sec:planar-tangle-operad} we discuss Jones's planar tangle operad, which parametrizes operations in shaded planar algebras.  For other examples of colored symmetric operads in applied contexts, the reader is referred to the books \cite{cg,mendez,spivak,yau-wd,yau-hqft}.

\section{Symmetric Operads as Monoids}\label{sec:symmetric_operad}

In this section, we define the colored symmetric circle product, the monoids of which are colored symmetric operads.  The \index{symmetric group}symmetric group on $n$ elements is denoted by $S_n$ with unit $\id_n$.  

\begin{definition}\label{def:symmetric-sequences}
Suppose $\colorc$ is a set, whose elements are called colors.
\begin{enumerate}
\item For $\colorc$-profiles $\ua = (a_1,\ldots,a_m)$ and $\ub$, a \index{left permutation}\emph{left permutation} $\sigma : \ua \to \ub$ is a permutation $\sigma \in S_{|\ua|}$ such that\label{notation:left-permutation}
\[\sigma\ua = (a_{\sigma^{-1}(1)}, \ldots , a_{\sigma^{-1}(m)}) = \ub\]
\item The \index{groupoid of profiles}\index{groupoid}\emph{groupoid of $\colorc$-profiles}, with $\Profc$ as object set and with left permutations as the isomorphisms, is denoted by $\ssubc$.\label{notation:sigmac}  The \index{opposite groupoid}opposite groupoid $\ssubcop$ is regarded as the groupoid of $\colorc$-profiles with \index{right permutation}\label{notation:right-permutation}\emph{right permutations}
\[\ua\sigma = (a_{\sigma(1)}, \ldots , a_{\sigma(m)})\]
as isomorphisms.
\item The objects of the diagram category\label{notation:symseqcm} \[\sseqcm = \M^{\ssubcopc}\] are called \index{symmetric sequence}\emph{$\colorc$-colored symmetric sequences} in $\M$.  
\end{enumerate}
\end{definition}

Forgetting the symmetric group actions, every $\colorc$-colored symmetric sequence has an underlying $\colorc$-colored planar sequence.  This forgetful functor will often be suppressed from the notation.  The symmetric version of the circle product is defined similarly to the planar case in \eqref{planar-circle-product}, with the groupoid $\ssubcop$, which parametrizes input profiles, replacing the discrete category $\Profc$ everywhere.  The following colored symmetric circle product is the colored version of the circle product defined by Kelly \cite{kelly}.

\begin{definition}\label{def:symmetric-circle-product}
Suppose $X,Y  \in \sseqcm$.
\begin{enumerate}
\item For each $\uc = (c_1,\ldots,c_m) \in \ssubc$, define the object $Y^{\uc}_{\smallsubs} \in \M^{\ssubcop}$ entrywise as the coend
\begin{equation}\label{symmetric-ytensorc}
Y^{\uc}_{\smallsubs} (\ub) = \int^{\{\ua_j\}\in\overset{m}{\prodover{j=1}} \ssubcop} \ssubcop\bigl(\ua_1,\ldots,\ua_m;\ub\bigr) \cdot \left[\bigotimes_{j=1}^m Y\cjuaj\right] \in \M
\end{equation}
for $\ub \in \ssubcop$, in which $(\ua_1,\ldots,\ua_m) \in \ssubcop$ is the concatenation.  The object $Y^{\uc}_{\smallsubs}$ is natural in $\uc \in \ssubc$ via left permutations of the tensor factors in $\bigotimes_{j=1}^m Y\cjuaj$.
\item The \index{circle product!symmetric}\index{symmetric circle product}\emph{$\colorc$-colored symmetric circle product} \[X \circs Y \in \sseqcm\] is defined entrywise as the coend
\begin{equation}\label{symmetric-circle-product}
(X \circs Y)\dub = \int^{\uc \in \ssubc} X\duc \otimes Y^{\uc}_{\smallsubs} (\ub)\in \M
\end{equation}
for $\dub \in \ssubcopc$.
\end{enumerate}
\end{definition}

The next observation is proved with the same argument as in the proof of Proposition \ref{planar-circle-product-monoidal}, with the groupoid $\ssubcop$ replacing the discrete category $\Profc$ everywhere.

\begin{proposition}\label{symmetric-circle-product-monoidal}
$\bigl(\sseqcm, \circs, \I\bigr)$ is a \index{monoidal category!of symmetric sequences}\index{symmetric sequence!monoidal category}monoidal category, where $\I$ is as in \eqref{unit-operad}.
\end{proposition}

\begin{definition}\label{def:symmetric-operad}
The category of \index{operad!symmetric}\index{symmetric operad}\emph{$\colorc$-colored symmetric operads in $\M$} is the category \[\Soperadcm = \Mon\bigl(\sseqcm, \circs, \I\bigr)\] of monoids in the monoidal category $\bigl(\sseqcm,\circs,\I\bigr)$.  If $\colorc$ has $n<\infty$ elements, we also refer to objects in $\Soperadcm$ as \emph{$n$-colored symmetric operads in $\M$}.
\end{definition}

\section{Coherence for Symmetric Operads}
\label{sec:symmetric-operad-generating}

As in the planar case, a symmetric operad can be described more explicitly in terms of generating operations.  First we need the following notations for block permutations and direct sum permutations.

\begin{definition}\label{def:block-direct-sum-permutation}
Suppose $n \geq 1$, $k_i\geq 0$ for $1\leq i \leq n$, and $k=k_1+\cdots+k_n$.
\begin{enumerate}
\item For $\sigma \in S_n$, the \index{block permutation}\emph{block permutation}\label{notation:block-permutation} \[\sigma\langle k_1,\ldots,k_n\rangle \in S_k\] is defined by
\[\bigl(\sigma\langle k_1,\ldots,k_n\rangle\bigr)(j) = \underbrace{k_{\sigmainv(1)}+\cdots+k_{\sigmainv(\sigma(i)-1)}}_{\text{$0$ if $\sigma(i)=1$}}+j\] for $1 \leq j \leq k_i$ and $1\leq i \leq n$.
\item For $\tau_i\in S_{k_i}$ for $1 \leq i \leq n$, the \index{direct sum!permutation}\emph{direct sum permutation}\label{notation:direct-sum} \[\tau_1\oplus\cdots\oplus\tau_n \in S_k\] is defined by \[\bigl(\tau_1\oplus\cdots\oplus\tau_n\bigr)(j) = \underbrace{k_1+\cdots+k_{i-1}}_{\text{$0$ if $i=1$}} +\tau_i(j)
\]for $1 \leq j \leq k_i$ and $1\leq i \leq n$.
\end{enumerate}
\end{definition}

\begin{interpretation}
The block permutation $\sigma\langle k_1,\ldots,k_n\rangle$ permutes $n$ consecutive blocks of lengths $k_1,\ldots,k_n$ as $\sigma$ permutes $\{1,\ldots,n\}$.  On the other hand, the direct sum $\tau_1\oplus\cdots\oplus\tau_n$ leaves the relative order of the $n$ consecutive blocks unchanged and only permutes within the $i$th block using $\tau_i$ for $1\leq i \leq n$.\dqed
\end{interpretation}

\begin{proposition}\label{def:symmetric-operad-generating}
A\index{coherence!symmetric operad}\index{symmetric operad!coherence} $\colorc$-colored symmetric operad in $\M$ is equivalent to a $\colorc$-colored planar operad $(\O,\gamma,\operadunit)$ as in Definition \ref{def:planar-operad-generating} together with a $\colorc$-colored symmetric sequence structure on $\O$ that satisfies the following equivariance axioms.  Suppose that in \eqref{operadic-composition} $|\ub_j| = k_j \geq 0$.
\begin{enumerate}
\item For each permutation $\sigma \in S_n$, the \index{equivariance!symmetric operad}\index{symmetric operad!equivariance}\emph{top equivariance diagram} 
\begin{equation}\label{operadic-eq-1}
\nicexy@C+.3cm{\O\duc \otimes \bigotimes\limits_{j=1}^n \O\cjubj 
\ar[d]_-{\gamma} \ar[r]^-{(\sigma, \sigma^{-1})}& \O\ducsigma \otimes 
\bigotimes\limits_{j=1}^n \O\csigmajubsigmaj \ar[d]^-{\gamma}\\
\O\duboneubn \ar[r]^-{\sigma\langle k_{\sigma(1)}, \ldots , k_{\sigma(n)}\rangle}& \O\dubsigmaoneubsigman}
\end{equation}
in $\M$ is commutative.  In the top horizontal morphism, $\sigma$ is the equivariant structure morphism of $\O$ corresponding to $\sigma \in S_n$, and $\inv{\sigma}$ permutes the $n$ tensor factors from the left.  The bottom horizontal morphism is the equivariant structure morphism of $\O$ corresponding to the block permutation 
\begin{equation}\label{block-permutation}
\sigma\langle k_{\sigma(1)}, \ldots , k_{\sigma(n)}\rangle \in S_{k_1+\cdots+k_n}
\end{equation} 
induced by $\sigma\in S_n$ that permutes $n$ consecutive blocks of lengths $k_{\sigma(1)}, \ldots , k_{\sigma(n)}$.  
\item Given permutations $\tau_j \in S_{k_j}$ for $1 \leq j \leq n$, the \emph{bottom equivariance diagram}
\begin{equation}\label{operadic-eq-2}
\nicexy@C+.3cm{\O\duc \otimes \bigotimes\limits_{j=1}^n \O\cjubj\ar[d]_-{\gamma} \ar[r]^-{(\Id, \otimes \tau_j)}& \O\duc \otimes \bigotimes\limits_{j=1}^n \O\cjubjtauj\ar[d]^-{\gamma} \\
\O\duboneubn \ar[r]^-{\tau_1 \oplus \cdots \oplus \tau_n}& \O\dubonetauoneubntaun}
\end{equation}
in $\M$ is commutative.  In the top horizontal morphism, each $\tau_j$ is the equivariant structure morphism of $\O$ corresponding to $\tau_j \in S_{k_j}$.  The bottom horizontal morphism is the equivariant structure morphism of $\O$ corresponding to the direct sum permutation 
\begin{equation}\label{direct-sum-permutation}
\tau_1 \oplus \cdots \oplus \tau_n \in S_{k_1+\cdots+k_n}
\end{equation} 
induced by the $\tau_j$'s.
\end{enumerate}
\end{proposition}

\begin{proof}
We reuse the proof of Proposition \ref{prop:planar-operad-same-def}.  The top equivariance diagram \eqref{operadic-eq-1} corresponds to the $\uc$-variable in the coend 
\[\begin{split}(\O\circs\O)\dub &= \int^{\uc \in \ssubc} \O\duc \otimes \O^{\uc}_{\smallsubs} (\ub)\\
&\cong \int^{\uc, \ua_1, \ldots, \ua_m} \ssubcop(\ua;\ub) \cdot \O\duc \otimes \bigotimes_{j=1}^m \O\cjuaj\end{split}\]
where $\ua = (\ua_1,\ldots,\ua_m)$, and the fact that the multiplication $\mu : \O \circs\O \to \O$ is a morphism of $\colorc$-colored symmetric sequences.  Similarly, the bottom equivariance diagram \eqref{operadic-eq-2} corresponds to the $\ua_j$ variables in the above coend and the fact that $\mu$ is a morphism of $\colorc$-colored symmetric sequences.
\end{proof}

As an example of Proposition \ref{def:symmetric-operad-generating}, we observe that the symmetric groups form a one-colored symmetric operad in $\Set$.  The following properties of permutations can be checked by a direct inspection, in which we use the notations in \eqref{block-permutation} and \eqref{direct-sum-permutation} for block permutations and direct sum permutations.

\begin{lemma}\label{symmetric-group-lemma}
For integers $n \geq 1$ and $k_1,\ldots,k_n \geq 0$ with $k=k_1+\cdots+k_n$, the following statements hold.
\begin{enumerate}
\item The equality \[(\sigma\tau)\langle k_1,\ldots,k_n\rangle = \sigma\langle k_{\tau^{-1}(1)},\ldots,k_{\tau^{-1}(n)}\rangle \cdot \tau\langle k_{1},\ldots,k_{n}\rangle \in S_k\] holds for $\sigma,\tau\in S_n$.
\item The equality \[(\sigma_1\tau_1) \oplus \cdots\oplus (\sigma_n\tau_n) = \bigl(\sigma_1\oplus\cdots\oplus\sigma_n\bigr) \cdot \bigl(\tau_1\oplus\cdots\oplus \tau_n\bigr) \in S_k\] holds for $\sigma_i,\tau_i\in S_{k_i}$ for $1\leq i \leq n$.
\item The equality \[\sigma\langle k_1,\ldots,k_n\rangle \cdot \bigl(\tau_1\oplus\cdots\oplus\tau_n\bigr) = \bigl(\tau_{\sigmainv(1)} \oplus\cdots\oplus \tau_{\sigmainv(n)}\bigr)\cdot \sigma\langle k_1,\ldots,k_n\rangle\] holds for $\sigma \in S_n$ and $\tau_i\in S_{k_i}$ for $1\leq i \leq n$.
\end{enumerate}
\end{lemma}

The following example of a one-colored symmetric operad follows from Proposition \ref{def:symmetric-operad-generating} and Lemma \ref{symmetric-group-lemma}.

\begin{proposition}\label{sn-operad}
There is a one-colored symmetric operad \[\S=\{S_n\}_{n\geq 0}\] in $\Set$ in which:
\begin{enumerate}[label=(\roman*)]
\item The operadic composition is given by \[\gamma^{\S}\bigl(\sigma;\tau_1,\ldots,\tau_n\bigr) = \sigma \langle k_1,\ldots,k_n\rangle \cdot \bigl(\tau_1\oplus\cdots\oplus\tau_n\bigr)\] for $\sigma \in S_n$ for $n\geq 1$ and $\tau_i \in S_{k_i}$ for $1 \leq i \leq n$ and $k_i \geq 0$.
\item The operadic unit is the identity permutation $\id_1\in S_1$.
\item The right $S_n$-action on $S_n$ is given by group multiplication.
\end{enumerate}
\end{proposition}

We call $\S$ in Proposition \ref{sn-operad} the \index{group operad!symmetric}\index{symmetric group operad}\emph{symmetric group operad}.  

\begin{remark}
In the symmetric group operad $\S$, block permutation and direct sum permutation can be expressed in terms of the operadic composition as
\[\begin{split}
\sigma \langle k_1,\ldots,k_n\rangle &= \gamma^{\S}\bigl(\sigma;\id_{k_1},\ldots,\id_{k_n}\bigr),\\
\bigl(\tau_1\oplus\cdots\oplus\tau_n\bigr) &= \gamma^{\S}\bigl(\id_n;\tau_1,\ldots,\tau_n\bigr).
\end{split}\]
It follows from Lemma \ref{symmetric-group-lemma}(1) that the operadic composition \[\nicexy{S_n \times S_{k_1}\times\cdots\times S_{k_n} \ar[r]^-{\gamma^{\S}} & S_{k_1+\cdots+k_n}}\] is \emph{not} a group homomorphism.  Therefore, $\S$ is not a one-colored symmetric operad in the symmetric monoidal category of groups.\dqed
\end{remark}

Next is the symmetric version of Proposition \ref{prop:planar-operad-defs-equivalent}, and it admits basically the same proof.

\begin{proposition}\label{prop:symmetric-operad-compi}
A $\colorc$-colored symmetric operad in $\M$ is equivalent to a $\colorc$-colored planar operad $(\O,\circ,\operadunit)$ as in Definition \ref{def:planar-operad-compi} together with a $\colorc$-colored symmetric sequence structure on $\O$ that satisfies the following equivariance axiom.  Suppose $|\uc| = n \geq 1$, $1 \leq i \leq n$, $\sigma \in S_n$, and $\tau \in S_m$.  Then the \emph{equivariance diagram} in $\M$\index{symmetric operad!equivariance}\index{equivariance!symmetric operad}
\begin{equation}\label{compi-eq}
\nicexy{\O\duc \otimes \O\sbinom{c_{\sigma(i)}}{\ub}\ar[d]_-{(\sigma,\tau)} \ar[r]^-{\comp_{\sigma(i)}}& 
\O\sbinom{d}{\uc\comp_{\sigma(i)}\ub} \ar[d]^-{\sigma\comp_{i}\tau}\\
\O\ducsigma \otimes \O\sbinom{c_{\sigma(i)}}{\ub\tau}\ar[r]^-{\compi} & 
\O\sbinom{d}{(\uc\sigma)\compi(\ub\tau)} 
=\O\sbinom{d}{(\uc\comp_{\sigma(i)}\ub)(\sigma\compi\tau)}}
\end{equation}
is commutative, where
\begin{equation}\label{symmetric-sigma-compi-tau}
\sigma \compi \tau = \overbrace{\sigma\langle \underbrace{1,\ldots,1}_{i-1}, m,\underbrace{1,\ldots,1}_{n-i}\rangle}^{\text{block permutation}}\cdot \overbrace{\bigl(\underbrace{\id_1 \oplus \cdots \oplus \id_1}_{i-1} \oplus \tau \oplus \underbrace{\id_1 \oplus \cdots \oplus \id_1}_{n-i}\bigr)}^{\text{direct sum permutation}} \in S_{n+m-1}
\end{equation}
is the product of a direct sum induced by $\tau$ with a block permutation induced by $\sigma$ that permutes consecutive blocks of the indicated lengths.
\end{proposition}

\begin{example}\label{ex:sn-compi}
For the symmetric group operad $\S$ in Proposition \ref{sn-operad}, the $\compi$-composition is given by $\sigma\compi\tau$.\dqed
\end{example}

\section{Algebras}\label{sec:symmetric-operad-algebra}

In this section, we discuss algebras over a colored symmetric operad, along with some basic examples.  Lemma \ref{lem:planar-operad-monad} has the following symmetric analogue.

\begin{lemma}\label{lem:symmetric-operad-monad}
Suppose $\O$ is a $\colorc$-colored symmetric operad in $\M$.  Then it induces a \index{symmetric operad!induced monad}\index{monad!induced by a symmetric operad}monad whose functor is \[\O \circs - : \Mtoc \to \Mtoc\] and whose multiplication and unit are induced by those of $\O$ as in Corollary \ref{planar-operad-monoid}.
\end{lemma}

\begin{definition}\label{def:symmetric-operad-algebra}
Suppose $\O$ is a $\colorc$-colored symmetric operad in $\M$.  The category $\algmo$ \index{symmetric operad!algebra}\index{algebra!of a symmetric operad}of \emph{$\O$-algebras} is defined as the category of $(\O\circs -)$-algebras for the monad $\O\circs -$ in $\Mtoc$.
\end{definition}

\begin{proposition}\label{prop:symmetric-operad-algebra-defs}
Suppose\index{coherence!symmetric operad algebra} $\O$ is a $\colorc$-colored symmetric operad in $\M$.  Then an $\O$-algebra is precisely a pair $(X,\lambda)$ as in Definition \ref{def:planar-operad-algebra-generating} that satisfies the following equivariance axiom.  For each $\duc \in \ssubcopc$ and each permutation $\sigma \in S_{|\uc|}$, the equivariance diagram\index{equivariance!symmetric operad algebra}
\begin{equation}\label{operad-algebra-eq}
\nicexy{\O\duc \otimes X_{\uc}\ar[d]_-{\lambda} \ar[r]^-{(\sigma,\inv{\sigma})}& \O \ducsigma \otimes X_{\uc\sigma}\ar[d]^-{\lambda}\\ X_d \ar@{=}[r]& X_d}
\end{equation}
in $\M$ is commutative.   In the top horizontal morphism, $\inv{\sigma}$ is the right permutation on the tensor factors in $X_{\uc}$ induced by $\sigma \in S_{|\uc|}$.
\end{proposition}

\begin{proof}
We reuse the proof of Proposition \ref{prop:planar-operad-algebra-defs} with the groupoid $\ssubc$ replacing the discrete category $\Profc$.  The equivariance diagram \eqref{operad-algebra-eq} comes from the $\uc$-variable in the coend \[(\O\circs X)_d = \int^{\uc\in\ssubc} \O\duc \otimes X_{\uc}\] for $d \in \colorc$. 
\end{proof}

Next we provide some basic examples of symmetric operads and their algebras.

\begin{example}[Colored Endomorphism Operads]\label{ex:endomorphism-operad}
Suppose that the underlying symmetric monoidal category $\M$ is \emph{closed}, which means that each functor \[-\otimes X : \M \to \M\] for $X\in \M$ admits a right adjoint $\Homm$, called the \emph{internal hom}.  For each $\colorc$-colored object $X=\{X_c\}_{c\in\colorc}$ in $\M$, there is a $\colorc$-colored symmetric operad $\End(X)$, called the \emph{endomorphism operad}\index{endomorphism operad!symmetric}\index{symmetric operad!endomorphism}, with entries
\[\End(X)\duc = \Homm(X_{\uc},X_d)\] for $\duc \in \ssubcopc$.  Its equivariant structure is induced by permutations of the tensor factors in $X_{\uc}$.  Its $d$-colored unit \[\tensorunit \to \Homm(X_d,X_d)\] is adjoint to the isomorphism  $\tensorunit \otimes X_d \cong X_d$.  Its operadic composition $\gamma$ is induced by the $\otimes$-$\Homm$-adjunction.  Via the $\otimes$-$\Homm$-adjunction, for a $\colorc$-colored symmetric operad $\O$ in $\M$, an $\O$-algebra structure $(X,\theta)$ is equivalent to a morphism \[\theta' : \O \to \End(X)\] of $\colorc$-colored symmetric operads.\dqed
\end{example}

\begin{example}[Commutative Operad]\label{ex:operad-com}
There is a $1$-colored symmetric operad $\Com$ in $\M$, called the \index{commutative operad}\emph{commutative operad}, with entries \[\Com(n) = \tensorunit \forspace n \geq 0,\] operadic composition induced by the isomorphism $\tensorunit \otimes \tensorunit \cong \tensorunit$, and unit the identity morphism.  The category of $\Com$-algebras is precisely the category of commutative monoids in $\M$.\dqed
\end{example}

\begin{example}[Operad for Diagrams of Monoids]\label{ex:diag-monoid-operad}
Suppose $\C$ is a small category with object set $\colorc$.  There is a $\colorc$-colored \index{operad!for diagrams of monoids}\index{monoid!operad for diagrams of}symmetric operad $\Asc$ in $\M$ with entries \[\Asc\duc = \coprodover{S_n \times \prod\limits_{j=1}^n \C(c_j,d)}\tensorunit \forspace \duc=\dconecn \in \ssubcopc.\]  A coproduct summand corresponding to an element $(\sigma,\uf) \in S_n \times \prod_j \C(c_j,d)$ is denoted by $\tensorunit_{(\sigma,\uf)}$.  We will describe the symmetric operad structure on $\Asc$ in terms of the subscripts.

Its equivariant structure sends $\tensorunit_{(\sigma,\uf)}$ to $\tensorunit_{(\sigma\tau,\uf\tau)}$ for $\tau \in S_{|\uc|}$.  Its $c$-colored unit corresponds to $\tensorunit_{(\id_1,\Id_c)}$.  Its operadic composition 
\[\nicexy{\Asc\duc \otimes \bigotimes\limits_{j=1}^n \Asc\cjubj \ar[r]^-{\gamma} & \Asc\dub}\]
corresponds to
\[\nicexy{\Bigl((\sigma,\uf);\bigl\{(\tau_j,\ug_j)\bigr\}_{j=1}^n\Bigr) \ar@{|->}[r] & \Bigl(\gamma^{\S}\bigl(\sigma;\tau_1,\ldots,\tau_n\bigr), \bigl(f_1\ug_1,\ldots,f_n\ug_n\bigr) \Bigr)}\] with \[f_j\ug_j=\bigl(f_jg_{j1},\ldots,f_jg_{jk_j}\bigr) \in \prod\limits_{i=1}^{k_j} \C(b_{ji},d) \ifspace \ug_j=\bigl(g_{j1},\ldots,g_{jk_j}\bigr) \in \prod\limits_{i=1}^{k_j} \C(b_{ji},c_j)\] and $\gamma^{\S}$ the operadic composition in the symmetric group operad in Proposition \ref{sn-operad}.

There is a canonical isomorphism \[\nicexy{\algm(\Asc) \ar[r]^-{\cong} & \Monmc}\] from the category of $\Asc$-algebras to the category of $\C$-diagrams in $\Monm$, the category of monoids in $\M$, defined as follows.  Each $\Asc$-algebra $(X,\lambda)$ has a restricted structure morphism \[\nicexy@C+.7cm{X_{\uc} \ar[r]^-{\lambda_{(\sigma,\uf)}} \ar[d]_-{(\sigma,\uf)~\mathrm{inclusion}} & X_d\\ \coprod\limits_{S_n\times\prod\limits_{j=1}^n \C(c_j,d)} X_{\uc} \ar[r]^-{\cong} & \Asc\duc \otimes X_{\uc} \ar[u]_-{\lambda}}\] for each $(\sigma,\uf) \in S_n \times \prod_{j=1}^n \C(c_j,d)$.  For a morphism $f : c \to d \in \C$, there is a restricted structure morphism \[\nicexy@C+.5cm{X_c \ar[r]^-{\lambda_{(\id_1,f)}} & X_d}\in \M.\]  The associativity and unity axioms of $(X,\lambda)$ imply that this is a $\C$-diagram in $\M$.  

For each $c \in \C$, the restricted structure morphisms \[\nicexy@C+1.2cm{X_c\otimes X_c \ar[r]^-{\lambda_{(\id_2,\{\Id_c,\Id_c\})}} & X_c} \andspace \nicexy@C+.4cm{\tensorunit \ar[r]^-{\lambda_{(\id_0,*)}} & X_c}\] give $X_c$ the structure of a monoid in $\M$, once again by the associativity and unity axioms of $(X,\lambda)$.  One can check that this gives a $\C$-diagram of monoids in $\M$; i.e., the morphisms $\lambda_{(\id_1,f)}$ are compatible with the entrywise monoid structures.  In summary, $\Asc$ is the $\colorc$-colored symmetric operad whose algebras are $\C$-diagrams of monoids in $\M$.\dqed
\end{example}

\begin{example}[Operad for Diagrams of Commutative Monoids]\label{ex:diag-com-operad}
Suppose $\C$ is a small category with object set $\colorc$.  There is a $\colorc$-colored \index{commutative monoid!operad for diagrams of}\index{operad!for diagrams of commutative monoids}symmetric operad \label{notation:comc}$\Comc$ in $\M$ with entries \[\Comc\duc = \coprodover{\prod\limits_{j=1}^n\C(c_j,d)} \tensorunit \forspace \duc=\dconecn \in \ssubcopc.\]  Its symmetric operad structure is defined as in Example \ref{ex:diag-monoid-operad} by ignoring the first component.  

Similar to Example \ref{ex:diag-monoid-operad}, there is a canonical isomorphism \[\nicexy{\algm(\Comc) \ar[r]^-{\cong} & \Commc}\] from the category of $\Comc$-algebras to the category of $\C$-diagrams of commutative monoids in $\M$.  To see that the monoid multiplication \[\nicexy@C+1.2cm{X_c\otimes X_c \ar[r]^-{\mu_c~=~\lambda_{\{\Id_c,\Id_c\}}} & X_c} \in \M\] is commutative, observe that the pair $\{\Id_c,\Id_c\}$ is fixed by the permutation $(1~2)$.  So the equivariance axiom \eqref{operad-algebra-eq} implies that $\mu_c$ is commutative.  In summary, $\Comc$ is the $\colorc$-colored symmetric operad whose algebras are $\C$-diagrams of commutative monoids in $\M$.\dqed  
\end{example}

\section{Little Cube and Little Disc Operads}\label{sec:symmetric-operad-little-cube}

The purpose of this section is to discuss the little $n$-cube operad, the little $n$-disc operad, and the framed little $n$-disc operad, all of which are important in homotopy theory.

\begin{example}[Little $2$-Cube Operad]\label{ex:little-2cube}
An important example of a $1$-colored symmetric operad is the\index{little $2$-cube operad}\index{operad!little $2$-cube} little $2$-cube operad $\C_2$, due to Boardman-Vogt \cite{boardman-vogt} (2.49) and May \cite{may}.  In this example we work over $\CHau$, the symmetric monoidal  category of compactly generated weak Hausdorff spaces.  

Suppose $\mathbb{R}$ is the topological space of real numbers.  Denote by \label{notation:interval} 
\begin{itemize}
\item $\cali$ the \index{closed unit interval}closed interval $[0,1]$ of real numbers, 
\item $\calj = (0,1)$ its \index{open unit interval}interior, 
\item$\cali^2 = [0,1] \times [0,1] \subseteq \mathbb{R}^2$ the \label{notation:square}\index{closed unit square}closed unit square, and 
\item $\calj^2 = (0,1) \times (0,1)$\label{notation:square-interior} its interior.  
\end{itemize}
A \index{little square}\emph{little square} is a function \[f=(f^1,f^2) : \cali^2 \to \cali^2\] such that each $f^i : \cali \to \cali$ is a linear function of the form \[f^i(t) = a_i + t(b_i-a_i) \forspace 0 \leq t \leq 1\] and some $0 \leq a_i < b_i \leq 1$.  For $n \geq 0$, an $n$-tuple $\uf = (f_1, \ldots, f_n)$ of little squares is said to have \index{pairwise disjoint interiors}\emph{pairwise disjoint interiors} if for any $1 \leq i < j \leq n$, the images of $f_i$ and $f_j$ have disjoint interiors, i.e., \[f_i(\calj^2) \cap f_j(\calj^2) = \varnothing.\] An $n$-tuple $\uf = (f_1,\ldots,f_n)$ of little squares with pairwise disjoint interiors is also regarded as a function \[\nicexy@C+.6cm{\coprod\limits_{i=1}^n \cali^2 \ar[r]^-{f_1 \amalg \cdots \amalg f_n} & \cali^2,}\] where if $n=0$ then this is regarded as the unique function $\varnothing \to \cali^2$.  The set of all continuous maps $\coprod_{i=1}^n \cali^2 \to \cali^2$ is a topological space with the compact-open topology.  The subset of $n$-tuples of little squares with pairwise disjoint interiors is given the subspace topology and is denoted by\label{notation:ctwo} $\C_2(n)$.  When $n=0$, $\C_2(0)$ is the one-point space containing only the function $\varnothing \to \cali^2$.

For example, here is a picture of an element
\begin{center}\begin{tikzpicture}[scale=2.5]
\draw [thick] (0,0) rectangle (1,1);
\draw [thick] (0.25,0.2) rectangle (0.7,0.4);
\draw [thick] (0.1,0.7) rectangle (0.4,0.95);
\draw [thick] (0.5,0.5) rectangle (0.9,0.85);
\node at (0.46,0.3) {1};
\node at (0.25,0.825) {2};
\node at (0.7,0.65) {3};
\node at (-.5,.5) {$A =$};
\node at (1.5,.5) {$\in \C_2(3)$};
\end{tikzpicture}\end{center}
with three pairwise disjoint little squares inside the unit square.

The $1$-colored symmetric operad structure on $\C_2$ is defined as follows.  
\begin{description}
\item[Unit] The unit is the identity map \[\bigl(\Id : \cali^2 \to \cali^2\bigr) \in \C_2(1).\]
\item[Equivariance]
Given $\uf = (f_1, \ldots, f_n) \in \C_2(n)$ and a permutation $\sigma  \in S_n$, define
\[f \sigma = \left(f_{\sigma(1)}, \ldots, f_{\sigma(n)}\right) \in \C_2(n),\] so the symmetric group $S_n$ acts on $\C_2(n)$ by permuting the labels of the $n$ little squares.
\item[Composition]
For $n \geq 1$, $1 \leq i \leq n$, and $m \geq 0$, the \emph{$\compi$-composition}
\[\nicexy{\C_2(n) \times \C_2(m) \ar[r]^-{\compi}& \C_2(n+m-1)}\]
is defined as follows.  Suppose $\uf = (f_1,\ldots,f_n) \in \C_2(n)$ and $\ug = (g_1,\ldots,g_m) \in \C_2(m)$.  Then \[\uf \compi \ug \in \C_2(n+m-1)\] is defined as the following composition.
\[\nicexy{\coprod\limits_{k=1}^{n+m-1} \cali^2 \ar[dr]_-{\Id^{i-1} \amalg \ug \amalg \Id^{n-i}} 
\ar[rr]^-{\uf \compi \ug} && \cali^2\\ & \coprod\limits_{k=1}^n \cali^2 \ar[ur]_-{\uf} &}\]
\end{description}
With the above structure, \[\C_2 = \bigl\{\C_2(n)\bigr\}_{n\geq 0}\] is a $1$-colored symmetric operad in $\CHau$, called the \emph{little $2$-cube operad}.

For example, with $A \in \C_2(3)$ and $B \in \C_2(2)$ as drawn on the left-hand side below, $A \comp_2 B \in \C_2(4)$ is the picture on the right below.
\begin{center}
\begin{tikzpicture}[scale=3]
\draw [thick] (0,0) rectangle (1,1); \node at (-.2,.5) {$B$};
\draw [thick, fill=yellow!20] (0.1,0.7) rectangle (0.4,0.9); \node at (0.25,0.8) {1};
\draw [thick, fill=green!20] (0.6,0.1) rectangle (0.8,0.8); \node at (0.7,0.45) {2};
\draw [thick] (0,-1.5) rectangle (1,-.5); \node at (-.2,-1) {$A $};
\draw [thick] (0.25,-1.3) rectangle (0.7,-1.1); \node at (0.46,-1.2) {1};
\draw [thick] (0.1,-.8) rectangle (0.4,-.55); \node at (0.25,-.675) {2};
\draw [thick] (0.5,-1) rectangle (0.9,-.65); \node at (0.7,-.85) {3};
\draw [->-=.5, gray!50, dotted, line width=2pt] (0,0) to (0.1,-.8);
\draw [->-=.5, gray!50, dotted, line width=2pt] (1,0) to (.4,-.8);
\node at (1.5,-.5) {$\becomes$};
\draw [thick] (2,-1) rectangle (3,0); \node at (2.5,.1) {$A \comp_2 B$};
\draw [thick] (2.25,-.8) rectangle (2.7,-.6); \node at (2.46,-.7) {1};
\draw [thick, fill=yellow!20] (2.13,-.125) rectangle (2.22,-.075); \node at (2.15,-.18) {\tiny{2}};
\draw [thick, fill=green!20] (2.28,-.275) rectangle (2.34,-.1); \node at (2.31,-.18) {\tiny{3}};
\draw [thick] (2.5,-.5) rectangle (2.9,-.15); \node at (2.7,-.35) {4};
\end{tikzpicture}
\end{center}
To compute $A \comp_2 B \in \C_2(4)$, we first scale $B$ down linearly to the dimensions of the little square labeled $2$ in $A$.  Then we replace this little square in $A$ by the linearly scaled down version of $B$, as indicated by the gray dotted lines above.\dqed
\end{example}

\begin{example}[Little $n$-Cube Operad]\label{ex:little-ncube}
In Example \ref{ex:little-2cube}, we can replace the square $\cali^2$ with the closed $n$-cube $\cali^n = [0,1]^{\times n}$ for any $1\leq n\leq \infty$ with similarly defined symmetric operad structure.  The result is the\index{little $n$-cube operad}\index{operad!little $n$-cube} \emph{little $n$-cube operad} $\C_n$.  For $k \geq 0$, the space $\C_n(k)$ is the space of $k$-tuples of little $n$-cubes with pairwise disjoint interiors.  The little $1$-cube operad $\C_1$ is also called an \emph{$A_\infty$-operad}, and the little $\infty$-cube operad $\C_\infty$ is also called an \emph{$E_\infty$-operad}.\dqed
\end{example}

\begin{example}[Little $n$-Disc Operad]\label{ex:little-n-disc}
For a positive integer $n$, a variant of the little $n$-cube operad is the little $n$-disc operad\index{little $n$-disc operad}\index{operad!little $n$-disc} $\D_n$ defined as follows.  Suppose \[\Dn = \Bigl\{(r_1,\ldots,r_n) \in \Rn :  r_1^2 +\cdots+r_n^2 \leq 1\Bigr\} \subseteq \Rn\] is the closed unit $n$-disc.  A \emph{little $n$-disc} is an embedding $f : \Dn \to \Dn$ of the form \[f(r_1,\ldots,r_n)= (c_1,\ldots,c_n) + \lambda(r_1,\ldots,r_n)\] for some $(c_1,\ldots,c_n) \in \Dn$ and $\lambda \in \fieldr$ such that \[0<\lambda \leq 1 - \sum_{i=1}^n c_i^2.\]  Geometrically, the function $f$ first dilates the closed unit $n$-disc down by the $\lambda$ factor, and then translates the result by adding the vector $(c_1,\ldots,c_n) \in \Rn$.  The set of $k$-tuples of little $n$-discs $(f_1,\ldots,f_k)$ with pairwise disjoint interiors is denoted by\label{notation:dn} $\D_n(k)$.  It is equipped with the subspace topology of the space of continuous maps \[\coprod_{i=1}^k \Dn \to \Dn,\] which itself has the compact-open topology, and $\D_n(0)$ is the one-point space.

The $1$-colored symmetric operad structure on \[\D_n = \bigl\{\D_n(k)\bigr\}_{k\geq 0}\] is defined similarly to the little $2$-cube operad $\C_2$ in Example \ref{ex:little-2cube}, with the closed unit $n$-disc $\Dn$ replacing the closed unit square $\cali^2$.  For example, with $A \in \D_2(3)$ and $B \in \D_2(2)$ as drawn below on the left, $A \comp_2 B \in \D_2(4)$ is the picture on the right below.
\begin{center}\begin{tikzpicture}[scale=1.2]
\draw [thick] (0,3) circle [radius=1];\node at (-1.3,3) {$A$};
\draw [thick] (-.4,3.4) circle [radius=.3]; \node at (-.4,3.4) {1};
\draw [thick] (0,2.5) circle [radius=.5]; \node at (0,2.5) {2};
\draw [thick] (.5,3.2) circle [radius=.25]; \node at (0.5,3.2) {3};
\draw [thick] (0,0) circle [radius=1]; \node at (0,1.2) {$B$};
\draw [thick, fill=green!20] (-.4,.4) circle [radius=.3]; \node at (-.4,.4) {2};
\draw [thick, fill=yellow!20] (.4,-.4) circle [radius=.25]; \node at (.4,-.4) {1};
\draw [->-=.5, gray!50, dotted, line width=2pt] (-1,0) to (-.5,2.5);
\draw [->-=.5, gray!50, dotted, line width=2pt] (1,0) to (.5,2.5);
\node at (1.3,1.5) {$\becomes$};
\draw [thick] (3,1.5) circle [radius=1]; \node at (3,2.7) {$A \comp_2 B$};
\draw [thick] (2.6,1.9) circle [radius=.3]; \node at (2.6,1.9) {1};
\draw [thick, fill=yellow!20] (3.2,.8) circle [radius=.125]; \node at (3.2,.8) {\scriptsize{2}};
\draw [thick, fill=green!20] (2.8,1.2) circle [radius=.15]; \node at (2.8,1.2) {\scriptsize{3}};
\draw [thick] (3.5,1.7) circle [radius=.25]; \node at (3.5,1.7) {4};
\end{tikzpicture}\end{center}
The two gray dotted lines are there to help visualize how the little $2$-disc labeled $2$ in $A$ is replaced by a scaled down version of $B$.  The $1$-colored symmetric operad $\D_n$ in $\CHau$ is called the \emph{little $n$-disc operad}.\dqed
\end{example}

\begin{example}[Framed Little $n$-Disc Operad]\label{ex:framed-little-n-disc}
An extension of the little $n$-disc operad in Example \ref{ex:little-n-disc} is the framed little $n$-disc operad\index{framed little $n$-disc operad}\index{operad!framed little $n$-disc} due to Getzler \cite{getzler}.  For $n \geq 2$, denote by\label{notation:son} $\SO(n)$ the\index{special orthogonal group} special orthogonal group of dimension $n$.  A \emph{framed little $n$-disc} is an embedding $f : \Dn \to \Dn$ of the form \[f(r_1,\ldots,r_n)= (c_1,\ldots,c_n) + \rho\Bigl[\lambda(r_1,\ldots,r_n)\Bigr]\] for some $\rho \in \SO(n)$, $(c_1,\ldots,c_n) \in \Dn$, and $\lambda \in \fieldr$ such that \[0<\lambda \leq 1 - \sum_{i=1}^n c_i^2.\]  So a framed little $n$-disc is determined by 
\begin{enumerate}[label=(\roman*)]
\item a dilation by factor $\lambda$,
\item a rotation by $\rho$, and 
\item a translation by $(c_1,\ldots,c_n)$.
\end{enumerate} 
The set of $k$-tuples of framed little $n$-discs with pairwise disjoint interiors is denoted by\label{notation:dfn} $\Df_n(k)$.  It is equipped with the product topology of \[\Df_n(k) = \D_n(k) \times \SO(n)^{\times k},\] where the $i$th copy of $\SO(n)$ corresponds to the rotation of the $i$th little $n$-disc.  The $1$-colored symmetric operad structure on \[\Df_n = \bigl\{\Df_n(k)\bigr\}_{k\geq 0}\] is defined similarly to the little $n$-disc operad $\D_n$ in Example \ref{ex:little-n-disc} with the rotations also taken into account.\dqed
\end{example}

\section{Operad in Non-Commutative Probability}
\label{sec:symmetric-operad-graph-operation}

The purpose of this section is to discuss the following one-colored symmetric operad from non-commutative probability theory due to Male \cite{male}.

\begin{example}[Operad of Graph Operations]\label{ex:graph-op}
By\index{operad!graph operations}\index{graph operation!operad} a \emph{finite directed graph}, or just a \index{graph}\emph{graph}, we mean a pair of finite sets $(V,E)$ with $V$ non-empty such that an element in $E$ is an ordered pair $(u,v) \in V^{\times 2}$, where each such ordered pair may appear in $E$ more than once.  Elements in $V$ and $E$ are called vertices and edges, respectively, and an edge $e = (u,v)$ is said to have initial vertex $u$ and terminal vertex $v$, denoted $e : u \to v$.  An edge of the form $(v,v)$ is called a \index{loop}loop at $v$.  An edge of the form $(u,v)$ or $(v,u)$ is said to \emph{connect} $u$ and $v$.  We say that a graph $(V,E)$ is \index{connected graph}\emph{connected} if for each pair of distinct vertices $u$ and $v$, there exist edges $e_i$ for $1 \leq i \leq n$ for some $n \geq 1$ such that each $e_i$ connects $v_{i-1}$ and $v_i$ with $v_0 = u$ and $v_n = v$.

A \index{graph operation}\emph{graph operation} is a connected graph $(V,E)$ equipped with
\begin{enumerate}[label=(\roman*)]
\item an ordering $\rho$ of the set $E$ of edges and
\item two possibly equal vertices $\inp$ and $\out$, called the \emph{input} and the \emph{output}.  
\end{enumerate}
An \emph{isomorphism} of graph operations is a pair of bijections $(V,E) \to (V',E')$ on vertices and edges that preserves the initial and the terminal vertices of each edge, the ordering on edges, and the input and the output.  We only consider graph operations up to isomorphisms.  That is, if there is an isomorphism
\[(V,E,\rho,\inp,\out) \iso (V',E',\rho',\inp',\out')\]
of graph operations, then we consider them to be the same.  For each $n \geq 0$, denote by \label{notation:grop} $\GrOp_n$ the set of graph operations with $n$ edges.  So $\GrOp_0$ contains only the graph with one vertex, which is both the input and the output, and no edges.  Here is an example of a graph operation with two vertices and four edges, two of which are loops:
\begin{center}
\begin{tikzpicture}
\matrix[row sep=1cm, column sep=1.5cm]{
\node [plain] (in) {$\inp$}; & \node [plain] (out) {$\out$};\\};
\draw [arrow, out=210, in=150, looseness=4] (in) to node{$1$} (in);
\draw [arrow, out=120, in=60, looseness=4] (in) to node{$2$} (in);
\draw [arrow, bend left] (in) to node{$3$} (out);
\draw [arrow, bend left] (out) to node{$4$} (in);
\end{tikzpicture}
\end{center}

There is a one-colored symmetric operad structure on graph operations given by \index{edge!substitution}\emph{edge substitution} as follows.  Suppose $G \in \GrOp_n$ with $n \geq 1$ and $G_i \in \GrOp_{k_i}$ for $1 \leq i \leq n$.  Then the operadic composition
\[G(G_1,\ldots,G_n) \in \GrOp_{k_1 + \cdots + k_n}\]
is obtained from $G$ by 
\begin{itemize}\item replacing the $i$th edge $e_i$ in $G$ by $G_i$ and 
\item identifying the initial (resp., terminal) vertex of $e_i$ with the input (resp., output) of $G_i$.
\end{itemize}
The edge ordering of the operadic composition is induced by those of $G$ and of the $G_i$'s.  The input and the output are inherited from $G$.  The symmetric group action on $\GrOp_n$ is given by permutation of the edge ordering.  The operadic unit is the graph operation $\inp \to \out$ with two vertices and one edge from the input to the output.

For example, suppose $G$, $H$, and $K$ are the following graph operations in $\GrOp_2$:
\begin{center}
\begin{tikzpicture}
\matrix[row sep=1cm, column sep=1cm]{
\node [plain] (ing) {$\inp$}; & \node [plain] (outg) {$\out$};
& \node [plain] (inh) {$\inp$}; & \node [plain] (outh) {$\out$};
& \node [plain] (ink) {$\inp$}; & \node [plain] (outk) {$\out$};\\};
\draw [arrow, out=120, in=60, looseness=4] (ing) to node{$1$} (ing);
\draw [arrow] (ing) to node{$2$} (outg);
\draw [arrow, bend left] (inh) to node{$1$} (outh);
\draw [arrow, bend right] (inh) to node[swap]{$2$} (outh);
\draw [arrow, bend left] (ink) to node{$1$} (outk);
\draw [arrow, bend left] (outk) to node{$2$} (ink);
\end{tikzpicture}
\end{center}
Then the operadic composition $G(H,K)$ is the graph operation with four edges above.\dqed
\end{example}

\section{Phylogenetic Operad}\label{sec:symmetric-operad-phylogenetic}

The purpose of this section is to discuss a one-colored symmetric operad due to Baez and Otter \cite{baez-otter}, called the phylogenetic operad, that has connections with biology.  Elements in the phylogenetic operad are trees, which we first define.  The set of non-negative real numbers is denoted by $[0,\infty)$, which is regarded as a subspace of the real line.  The following concept of trees is due to Weiss \cite{weiss}.

\begin{definition}\label{def:tree-for-phyl}
A \emph{tree}\index{tree} is a pair \[\bigl((T,\leq),\inp(T)\bigr)\] consisting of:
\begin{enumerate}[label=(\roman*)]
\item a finite partially ordered set $(T,\leq)$ that has a unique smallest element $\out(T)$, called the \emph{root}\index{root} or the\index{output} \emph{output}, such that the set \[\bigl\{f \in T : f \leq e\bigr\}\] is linearly ordered for each $e \in T$;  
\item a subset $\inp(T)\subseteq T$ of maximal elements (so $e \in \inp(T)$ implies there are no $f \in T$ such that $e< f$), called the\index{inputs} \emph{inputs} of $T$.
\end{enumerate}
We will usually abbreviate such a tree to $T$.  For a tree $T$:
\begin{enumerate}
\item An element in $T$ is called an\index{edge} \emph{edge}.  The set of edges in $T$ is also  denoted by\label{notation:edge} $\Ed(T)$.  An\index{internal edge} \emph{internal edge} is an edge that is neither the root nor an input of $T$.  The set of internal edges is denoted by $\Int(T)$.  
\item For $e,f \in T$, we write $e \prec f$ if (i) $e<f$ and (ii) there are no $g \in T$ such that $e < g < f$.  
\item For each $e \in T \setminus \inp(T)$, we define the set\label{notation:inpe} \[\inp(e) = \bigl\{f \in T : e \prec f\bigr\},\] and call \[v = \{e\} \cup \inp(e)\] the \index{vertex}\emph{vertex} with \emph{output} $\out(v)=e$ and \emph{inputs} $\inp(v) = \inp(e)$.  We call $v$ the\index{initial vertex} \emph{initial vertex} of $e$ and the\index{terminal vertex} \emph{terminal vertex} of $f \in \inp(e)$.  The set of vertices in $T$ is denoted by $\Vt(T)$.
\item An\index{edge!length function} \emph{edge length function} is a map \[L : \Ed(T) \to [0,\infty).\] We call $L(e)$ the \emph{length} of $e \in \Ed(T)$.
\item An\index{input ordering} \emph{input ordering} is a bijection \[\sigma : \{1,\ldots,|\inp(T)|\} \iso \inp(T).\]  For a given input ordering $\sigma$ of $T$, we call $\sigma(i)$ the $i$th input of $T$.
\item An\index{isomorphism of trees} \emph{isomorphism} $\psi : T \iso T'$ of trees is an isomorphism of partially ordered sets that preserves the inputs, as well as other structures (e.g., edge length functions and/or input orderings) that both $T$ and $T'$ possess.
\end{enumerate}
\end{definition}

Recall that $\CHau$ is the symmetric monoidal category of compactly generated weak Hausdorff spaces, with Cartesian product as the monoidal product.

\begin{definition}\label{def:phylogenetic-operad}
The\index{phylogenetic operad}\index{operad!phylogenetic} \emph{phylogenetic operad} $\Phyl$ is the one-colored symmetric operad in $\CHau$ defined as follows.
\begin{description}
\item[Entries] For $n \geq 0$ the space $\Phyl(n)$ is defined as the quotient \[\Phyl(n) = \Bigl(\coprodover{(T,\sigma)} [0,\infty)^{\times \Ed(T)}\Bigr)\Big/\sim\] in which the coproduct is indexed by the set of isomorphism classes of trees with input orderings $(T,\sigma)$ such that:
\begin{itemize}
\item $|\inp(T)|=n$.
\item $|\inp(v)|\not= 0,1$ for $v\in\Vt(T)$.
\end{itemize}
A point in $[0,\infty)^{\times \Ed(T)}$ is also regarded as an edge length function $L : \Ed(T) \to [0,\infty)$.  The identification $\sim$ is generated as follows.  Suppose $(T,\sigma,L)$ is such a triple, and $e\in \Int(T)$ has length $0$.  Then \[(T,\sigma,L) \sim (T',\sigma,L')\] in which:
\begin{itemize}
\item $T'$ is the tree obtained from $T$ by removing $e$, with the partial ordering induced by that of $T$.  So \[\Ed(T') = \Ed(T)\setminus \{e\} \andspace \inp(T')=\inp(T).\] 
\item The edge length function $L' : \Ed(T') \to [0,\infty)$ is the restriction of $L$.
\end{itemize}
\item[Unit] The unit \[(\uparrow,0) \in \Phyl(1) \cong [0,\infty)\] is the point $0\in [0,\infty)$, corresponding to the tree $\uparrow$ with only one edge, which has length $0$, and no vertices.
\item[Composition] Suppose $(T,\sigma,L)\in \Phyl(n)$ with $1 \leq i \leq n$ and $(T',\sigma',L')\in \Phyl(m)$.  Their $\compi$-composition \[(T,\sigma,L) \compi (T',\sigma',L') = \bigl(T\compi T', \sigma \compi \sigma', L \compi L'\bigr) \in \Phyl(n+m-1)\] is defined as follows.
\begin{itemize}
\item $T \compi T'$ is the tree \[T \compi T' = \frac{T \amalg T'}{\bigl(\sigma(i) =\out(T')\bigr)}\] obtained by identifying the $i$th input of $T$ with the output of $T'$, with inputs \[\inp(T\compi T') = \bigl[\inp(T) \setminus\{\sigma(i)\}\bigr] \amalg \inp(T').\]
\item The input ordering \[\nicexy{\bigl\{1,\ldots,|\inp(T)|-1+|\inp(T')|\bigr\} \ar[r]^-{\sigma\compi\sigma'}_-{\cong} & \inp(T\compi T')}\] is the bijection defined as\[(\sigma\compi\sigma')(j) = \begin{cases} \sigma(j) & \text{if $1 \leq j < i$},\\
\sigma'(j-i+1) & \text{if $i \leq j \leq i+|\inp(T')|-1$},\\
\sigma\bigl(j-|\inp(T')|+1\bigr) & \text{if $i+|\inp(T')|\leq j \leq |\inp(T)|-1+|\inp(T')|$}.\end{cases}\]
\item The edge length function \[\nicexy@C+.5cm{\Ed(T\compi T') \ar[r]^-{L\compi L'} & [0,\infty)}\] is defined as\[(L\compi L')(e) = \begin{cases} L(e) & \text{if $e\in \Ed(T) \setminus\{\sigma(i)\}$},\\
L'(e) & \text{if $e\in \Ed(T') \setminus\{\out(T')\}$},\\
L\bigl(\sigma(i)\bigr) + L'\bigl(\out(T')\bigr) & \text{if $e=\sigma(i)$}.
\end{cases}\]
\end{itemize}
\item[Equivariance]
For $(T,\sigma,L) \in \Phyl(n)$ and $\tau \in S_n$, the $\tau$-action is given by \[(T,\sigma,L)\tau = (T,\sigma\tau,L)\] with $\sigma\tau$ the composite \[\nicexy{\{1,\ldots,n\} \ar[r]^-{\tau}_-{\cong} & \{1,\ldots,n\} \ar[r]^-{\sigma}_-{\cong} & \inp(T) \ar@{<-}`u[ll]`[ll]_-{\sigma\tau}[ll].}\]
\end{description}
\end{definition}

It is an exercise to check that $\Phyl$ in Definition \ref{def:phylogenetic-operad} satisfies the axioms in Proposition \ref{prop:symmetric-operad-compi}, so it is indeed a one-colored symmetric operad in $\CHau$.

\begin{interpretation}
In each point $(T,\sigma,L)\in\Phyl(n)$ in the phylogenetic operad, each edge in the tree $T$ represents a species.  Each vertex $v$ represents a splitting of the species $\out(v)$ into the species in the set $\inp(v)$.  The edge lengths represent time.  The identification $\sim$ collapses length $0$ internal edges.  This means that, if a species splits and then immediately splits again, then we consider both splittings as happening at the same moment. 

The exclusion of $0$-ary vertices--i.e., those with $|\inp(v)|=0$--means that we are not considering extinct species.  The exclusion of unary vertices--i.e., those with $|\inp(v)|=1$--means that we only consider the situation when a species is split into at least two species.  We emphasize that the exclusion of $0$-ary and unary vertices is \emph{not} necessary for the mathematics to work.  There is an extended version of the phylogenetic operad without the requirement that $|\inp(v)|\not=0,1$ for $v\in\Vt(T)$.  The restriction on the size of $|\inp(v)|$ is implemented to reflect actual practice of biologists when they construct phylogenetic trees \cite{baum-smith}.\dqed
\end{interpretation}

\begin{example}
When we draw a tree, we draw\index{tree!drawing} each edge as either an arrow $\uparrow$ or simply as a line $\vert$.  A vertex is drawn as a dot $\bullet$ or a circle $\bigcircle$, with its output above it and its inputs below it, as follows.
\begin{center}\begin{tikzpicture}
\node[plain] (v) {$v$}; \node[below=.1cm of v] () {\scriptsize{$\inp(v)$}};
\draw[outputleg] (v) to node[at end]{\scriptsize{$\out(v)$}} +(0cm,.8cm);
\draw[thick] (v) to +(-.8cm,-.6cm);\draw[thick] (v) to +(.8cm,-.6cm);
\end{tikzpicture}\end{center}
The output of the tree is drawn at the top.  

On the left below is an example of a point $(T,\sigma,L)$ in the space $\Phyl(4)$.  It is identified under $\sim$ with the point $(T',\sigma,L')$ in $\Phyl(4)$ on the right.
\begin{center}\begin{tikzpicture}
\node[normaldot] (w) {}; \node[above=.8cm of w] (T) {$(T,\sigma,L)$};
\node[below left=1cm of w, normaldot] (u) {}; \node[below=.9cm of u] (u1){};
\node[below right=1cm of w, normaldot] (v) {}; \node[below=.9cm of v] (v1){}; 
\draw[outputleg] (w) to node[pos=.6]{\scriptsize{$1.5$}} +(0,.9cm);
\draw[thick] (u) to node{\scriptsize{$2$}} (w);
\draw[thick] (v) to node[swap]{\scriptsize{$0$}} (w);
\draw[thick] (u) to node[swap, outer sep=-2pt]{\scriptsize{$4$}} +(-.5cm,-.9cm);
\draw[thick] (u) to node[outer sep=-2pt]{\scriptsize{$0$}} +(.5cm,-.9cm);
\draw[thick] (v) to node[swap, outer sep=-2pt]{\scriptsize{$1$}} +(-.5cm,-.9cm);
\draw[thick] (v) to node[outer sep=-2pt]{\scriptsize{$0.5$}} +(.5cm,-.9cm);
\node[left=.2cm of u1] (){\scriptsize{$\sigma(2)$}}; 
\node[right=.05cm of u1] (){\scriptsize{$\sigma(4)$}};
\node[left=.05cm of v1] (){\scriptsize{$\sigma(3)$}}; 
\node[right=.2cm of v1] (){\scriptsize{$\sigma(1)$}};
\node[right=2.5cm of w] (s) {\Huge{$\sim$}};
\node[right=2.5cm of s,normaldot] (w') {}; \node[above=.8cm of w'] {$(T',\sigma,L')$};
\node[below left=1cm of w', normaldot] (u') {}; \node[below=.9cm of u'] (u1'){};
\draw[outputleg] (w') to node[pos=.6]{\scriptsize{$1.5$}} +(0,.9cm);
\draw[thick] (w') to node[swap, outer sep=-2pt]{\scriptsize{$1$}} +(0cm,-.9cm); 
\draw[thick] (w') to node[outer sep=-2pt]{\scriptsize{$0.5$}} +(.9cm,-.9cm);
\draw[thick] (u') to node{\scriptsize{$2$}} (w');
\draw[thick] (u') to node[swap, outer sep=-2pt]{\scriptsize{$4$}} +(-.5cm,-.9cm);
\draw[thick] (u') to node[outer sep=-2pt]{\scriptsize{$0$}} +(.5cm,-.9cm);
\node[left=.2cm of u1'](){\scriptsize{$\sigma(2)$}}; 
\node[right=.2cm of u1'](){\scriptsize{$\sigma(4)$}};
\node[below=.7cm of w'] (w1){\scriptsize{$\sigma(3)$}}; 
\node[right=.3cm of w1] (){\scriptsize{$\sigma(1)$}};
\end{tikzpicture}
\end{center}
There are seven edges in $T$, with the edge length of each edge written next to it.  By the identification $\sim$, the length $0$ internal edge in $T$ can be removed (i.e., shrunk away), resulting in $T'$.

An illustration of the $\comp_3$-composition \[\nicexy{\Phyl(4) \times \Phyl(3) \ar[r]^-{\comp_3} & \Phyl(6)}\] in the phylogenetic operad is in the picture below.
\begin{center}\begin{tikzpicture}
\node[normaldot] (w') {}; \node[left=1.5cm of w'] {$(T',\sigma,L')$};
\node[below left=1cm of w', normaldot] (u') {}; 
\node[below=.8cm of u'] (u1'){\scriptsize{$\sigma(4)$}};
\draw[outputleg] (w') to node[pos=.6]{\scriptsize{$1.5$}} +(0,.9cm);
\draw[thick] (w') to node[swap, outer sep=-2pt]{\scriptsize{$1$}} +(0cm,-.9cm); 
\draw[thick] (w') to node[outer sep=-2pt]{\scriptsize{$0.5$}} +(.9cm,-.9cm);
\draw[thick] (u') to node{\scriptsize{$2$}} (w');
\draw[thick] (u') to node[swap, outer sep=-2pt]{\scriptsize{$4$}} +(-.9cm,-.9cm);
\draw[thick] (u') to node[outer sep=-1pt]{\scriptsize{$0$}} +(0cm,-.9cm);
\node[left=.2cm of u1'](){\scriptsize{$\sigma(2)$}}; 
\node[below=.7cm of w'] (w1){\scriptsize{$\sigma(3)$}}; 
\node[right=.3cm of w1] (){\scriptsize{$\sigma(1)$}};
\node[below=.9cm of w'] (B) {}; \node[below=2.5cm of w'] (A) {};
\draw[gray!50, ->, dotted, line width=2pt] (A) to node[swap]{$\comp_3$} (B);
\node[below=3.5cm of w', normaldot] (w'') {}; \node[left=1.5cm of w''] () {$(T'',\tau,L'')$};
\node[below=.8cm of w''] (w1''){};
\node[below=1cm of w'', normaldot] (u'') {}; \node[below=.8cm of u''] (u1''){};
\draw[outputleg] (w'') to node[pos=.6]{\scriptsize{$5$}} +(0,.9cm);
\draw[thick] (w'') to node[outer sep=-2pt]{\scriptsize{$3$}} +(.9cm,-.9cm);
\draw[thick] (u'') to node{\scriptsize{$0.4$}} (w'');
\draw[thick] (u'') to node[pos=.6,swap, outer sep=-2pt]{\scriptsize{$1.2$}} +(-.5cm,-.9cm);
\draw[thick] (u'') to node[pos=.6, outer sep=-2pt]{\scriptsize{$0.5$}} +(.5cm,-.9cm);
\node[left=.2cm of u1''] (){\scriptsize{$\tau(1)$}}; 
\node[right=.1cm of u1''] (){\scriptsize{$\tau(3)$}};
\node[right=.5cm of w1''] (){\scriptsize{$\tau(2)$}}; 
\node[right=5cm of B,normaldot] (w) {}; 
\node[above=1cm of w] () {$(T',\sigma,L')\comp_3(T'',\tau,L'')$};
\node[below=.8cm of w] (w1){};
\node[below left=1cm of w, normaldot] (u) {}; 
\node[below=.7cm of u] (u1) {\scriptsize{$\pi(6)$}};
\node[below=1cm of w, normaldot] (v) {}; \node[below=.8cm of v] (v1) {};
\node[below=1cm of v, normaldot] (t){}; \node[below=.8cm of t] (t1) {};
\draw[outputleg] (w) to node[pos=.6]{\scriptsize{$1.5$}} +(0,.9cm);
\draw[thick] (w) to node[outer sep=-2pt]{\scriptsize{$0.5$}} +(.9cm,-.9cm);
\draw[thick] (u) to node{\scriptsize{$2$}} (w); 
\draw[thick] (v) to node[outer sep=-2pt]{\scriptsize{$6$}} (w);
\draw[thick] (t) to node[outer sep=-2pt]{\scriptsize{$0.4$}} (v);
\draw[thick] (u) to node[swap, outer sep=-2pt]{\scriptsize{$4$}} +(-.9cm,-.9cm);
\draw[thick] (u) to node[outer sep=-2pt]{\scriptsize{$0$}} +(0cm,-.9cm);
\draw[thick] (v) to node[outer sep=-2pt]{\scriptsize{$3$}} +(.9cm,-.9cm);
\draw[thick] (t) to node[pos=.6, swap, outer sep=-2pt]{\scriptsize{$1.2$}} +(-.5cm,-.9cm);
\draw[thick] (t) to node[pos=.6, outer sep=-2pt]{\scriptsize{$0.5$}} +(.5cm,-.9cm);
\node[right=.5cm of w1] () {\scriptsize{$\pi(1)$}};
\node[left=.2cm of u1]() {\scriptsize{$\pi(2)$}}; 
\node[right=.5cm of v1]() {\scriptsize{$\pi(4)$}};
\node[left=.2cm of t1]() {\scriptsize{$\pi(3)$}}; 
\node[right=.2cm of t1]() {\scriptsize{$\pi(5)$}};
\end{tikzpicture}
\end{center}
The $\comp_3$-composition \[(T',\sigma,L') \comp_3 (T'',\tau,L'') = \bigl(T' \comp_3 T'', \sigma\comp_3 \tau, L'\comp_3 L''\bigr)\in \Phyl(6)\]is on the right-hand side in the previous picture, in which $\pi = \sigma \comp_3 \tau \in S_6$.  The gray dotted arrow is there to help visualize the identification of the third input of $T'$ with the output of $T''$
\dqed
\end{example}

\begin{remark}There are several essentially equivalent formalisms\index{tree!formalisms} of trees in the literature, including those in \cite{baez-otter,fresse-gt,gk,hry-cyclic,kock,weiss,bluemonster,yau-operad}.  The interested reader can easily translate the discussion above as well as in Section \ref{sec:planar-tree} below about trees to other similar definitions of trees.\dqed
\end{remark}

\section{Planar Tangle Operad}\label{sec:planar-tangle-operad}

The purpose of this section is to discuss the colored symmetric operad called the planar tangle operad due to Jones \cite{jones}.  The significance of the planar tangle operad comes from the fact that it parametrizes the operations in shaded planar algebras.  We first define the set of colors in the planar tangle operad.  The complex plane is denoted by $\fieldc$.\label{notation:fieldc}

\begin{definition}\label{def:colored-circle}
For $k \geq 0$, a\index{shaded circle} \emph{shaded $k$-circle} is a tuple \[(C,P,S,\ast)\] as follows:
\begin{enumerate}[label=(\roman*)]
\item $C$ is a circle in $\fieldc$.
\item $P$ is a set of $2k$ distinct points on the circle $C$, called\index{marked point} \emph{marked points}.  Each connected component of $C \setminus P$ is called an\index{interval in a circle} \emph{interval} in $C$.  The set of intervals in $C$ is denoted by\label{notation:interval-circle} $\Intl(C)$.
\item $S : \Intl(C) \to \{b,w\}$ is a map, called the\index{shading} \emph{shading} of $C$, such that if the closures of two distinct intervals $I_1$ and $I_2$ have a non-empty intersection, then $S(I_1)\not=S(I_2)$.  For an interval $I$ in $C$, if $S(I)=b$ (resp., $S(I)=w$), then we say that $I$ is shaded black (resp., white).
\item $\ast \in \Intl(C)$, called the\index{distinguished interval} \emph{distinguished interval} in $C$, which may be shaded either black or white.
\end{enumerate}
We abbreviate such a shaded $k$-circle to $C$ if there is no danger of confusion.  A \emph{shaded circle} is a shaded $k$-circle for some $k$.  The set of shaded circles is denoted by $\SCir$.
\end{definition}

\begin{example}
Here is an example of a shaded $3$-circle.
\begin{center}\begin{tikzpicture}[scale=1]
\draw [thick] (0,0) circle [radius=1];
\coordinate (A) at (0:1); \node at (22.5:1.2) {\scriptsize{$b$}};
\coordinate (B) at (45:1); \node at (125:1.2) {\scriptsize{$w$}}; 
\coordinate (C) at (155:1); \node at (167.5:1.2) {\scriptsize{$b$}};
\coordinate (D) at (180:1); \node at (235:1.2) {\scriptsize{$w$}}; 
\coordinate (E) at (290:1); \node at (305:1.2) {\scriptsize{$b$}};
\coordinate (F) at (320:1); \node at (340:1.2) {\scriptsize{$w$}}; 
\foreach \X in {A,B,C,D,E,F} \fill [black] (\X) circle (2pt);
\node at (75:1.2) {$\ast$};
\end{tikzpicture}\end{center}
Each marked point is indicated by $\bullet$.  The shading is indicated by either $b$ or $w$ next to the interval.  The distinguished interval is indicated by $\ast$, which is shaded white.  Note that the shading of $C$ is uniquely determined by the shading of one interval, such as the distinguished interval, since the intervals are shaded in an alternating manner as one goes around the circle.\dqed
\end{example}

Recall from Definition \ref{def:colored-sequences} that $\Profc$ is the set of $\colorc$-profiles.  Below we will use this notation when the set $\colorc$ is the set $\SCir$ of shaded circles.  By a\index{closed disc} \emph{closed disc} we mean a subspace \[D = \bigl\{z \in \fieldc : |z-a|\leq r\bigr\} \subseteq \fieldc\] for some $a\in\fieldc$ and some positive real number $r$.   For a subspace $X \subseteq \fieldc$, its \index{interior}interior, \index{boundary}boundary, and \index{closure}closure are denoted by\label{notation:interior} $\interior{X}$, $\partial X$, and $\overline{X}$, respectively. 

\begin{definition}\label{def:admissible-scir-profile}
Suppose \[\sbinom{C_0}{C_1,\ldots,C_n} \in \profscirscir\] for some $n\geq 0$, so each tuple \[(C_i,P_i,S_i,\ast_i)\] is a shaded $k_i$-circle for some $k_i\geq 0$.  Suppose $D_i$ is the closed disc with boundary $\partial D_i=C_i$ for $0 \leq i \leq n$.  We say that $\sbinom{C_0}{C_1,\ldots,C_n}$ is \index{shaded circle!admissible}\index{admissible}\emph{admissible} if the closed discs $\{D_j\}_{j=1}^n$ are
\begin{enumerate}[label=(\roman*)]
\item pairwise disjoint and 
\item contained in the interior $D_0^{\circ}$, except when $(C_1,\ldots,C_n)=(C_0)$, in which case $D_0=D_1$.
\end{enumerate}
\end{definition}

\begin{example}
For each shaded circle $C$, \[\sbinom{C}{C}\in\profscirscir\] is admissible.  If $C' \not=C$ is another shaded circle, then \[\sbinom{C}{C'}\in\profscirscir\] is admissible if and only if the circle $C'$ is inside the circle $C$.\dqed
\end{example}

Next we recall the concept of a planar tangle from \cite{jones}.  We will reuse the notations in Definition \ref{def:admissible-scir-profile} below.

\begin{definition}\label{def:planar-tangle}
Suppose $\sbinom{C_0}{C_1,\ldots,C_n} \in \profscirscir$ is admissible for some $n\geq 0$.  A\index{planar tangle} \emph{planar $k_0$-tangle} with \emph{profile} $\sbinom{C_0}{C_1,\ldots,C_n}$ is a pair \[T = \bigl(St,Sh\bigr)\] as follows.
\begin{enumerate}[label=(\roman*)]
\item $St$ is a finite set of disjoint smooth curves in $D_0\setminus \coprod_{j=1}^n D_j^{\circ}$, called\index{string!in a planar tangle} \emph{strings} in $T$, such that:
\begin{itemize}
\item Each string is either a closed curve in $D_0^{\circ}\setminus \coprod_{j=1}^n D_j$ , or its two boundary points are distinct marked points of the shaded circles $\{C_i\}_{i=0}^n$.  
\item Each marked point in $\{C_i\}_{i=0}^n$ is the boundary point of precisely one string, which meets the corresponding circle transversally.  
\end{itemize}  
The connected components of \[\Bigl[D_0^{\circ} \setminus \coprod_{j=1}^n D_j\Bigr] \setminus St\] are called the\index{region in a planar tangle} \emph{regions} in $T$, the set of which is denoted by\label{notation:regt} $\Reg(T)$.
\item $Sh : \Reg(T) \to \{b,w\}$ is a map, called the\index{shading!of a planar tangle} \emph{shading} of $T$, such that:
\begin{itemize}
\item If $R_1$ and $R_2$ are distinct regions such that $\overline{R_1} \cap \overline{R_2} \not=\varnothing$, then \[Sh(R_1) \not= Sh(R_2).\]
\item If $R$ is a region and if $I$ is an interval in $C_i$ for some $0\leq i \leq n$ such that $\overline{R} \cap I \not=\varnothing$, then \[Sh(R)=S_i(I).\]
\end{itemize}
If $Sh(R) = b$ (resp., $Sh(R)=w$), then we say that $R$ is shaded black (resp., white).
\end{enumerate}
We call $C_0$ the\index{outside shaded circle} \emph{outside shaded circle} and $C_j$ for $1\leq j \leq n$ the $j$th\index{inside shaded circle} \emph{inside shaded circle} of $T$.  A \emph{planar tangle} is a planar $k$-tangle for some $k\geq 0$ with some admissible profile.  The set of planar tangles with profile $\sbinom{C_0}{C_1,\ldots,C_n}$ is denoted by $\PTan\sbinom{C_0}{C_1,\ldots,C_n}$.  If $\sbinom{C_0}{C_1,\ldots,C_n}$ is not admissible, we define $\PTan\sbinom{C_0}{C_1,\ldots,C_n} = \varnothing$. 
\end{definition}

\begin{example}
For each shaded circle $C$, $\PTan\sbinom{C}{C}$ is the one-point set, since the set of strings and the set of regions are both necessarily empty.\dqed
\end{example}

\begin{example}\label{ex:planar-2-tangle}
Here is an example of a planar $2$-tangle $T \in \PTan\sbinom{C_0}{C_1,C_2,C_3}$.
\def\ra{.5} \def\raplus{.58} 
\def\rb{.2} \def\rbplus{.27} 
\begin{center}\begin{tikzpicture}[scale=2]
\coordinate (c) at (0,0); \draw[thick] (c) circle [radius=1]; 
\node at ([shift={(110:1.15)}] c) {$\ast$}; 
\node at ([shift={(150:1.15)}] c) {\scriptsize{$w$}};
\coordinate (p01) at (80:1); \coordinate (p02) at (180:1);
\coordinate (p03) at (220:1); \coordinate (p04) at (330:1);
\coordinate (c1) at ([shift={(135:.4)}] c); \draw [thick] (c1) circle [radius=\ra]; 
\node at (c1) {\scriptsize{$D_1$}};
\draw [thick] (c1) circle [radius=\ra]; 
\coordinate (p11) at ([shift={(90:\ra)}] c1); 
\coordinate (p12) at ([shift={(270:\ra)}] c1);
\node at ([shift={(0:\raplus)}] c1) {$\ast$}; 
\node at ([shift={(30:\raplus)}] c1) {\scriptsize{$b$}}; 
\coordinate (c2) at ([shift={(-30:.6)}] c); \node at (c2) {\scriptsize{$D_2$}};
\draw [thick] (c2) circle [radius=\rb]; 
\coordinate (p21) at ([shift={(135:\rb)}] c2);  
\coordinate (p22) at ([shift={(15:\rb)}] c2); 
\coordinate (p23) at ([shift={(180:\rb)}] c2);  
\coordinate (p24) at ([shift={(270:\rb)}] c2); 
\node at ([shift={(100:\rbplus)}] c2) {$\ast$};
\node at ([shift={(60:\rbplus)}] c2) {\scriptsize{$b$}};
\coordinate (c3) at ([shift={(250:.7)}] c); \node at (c3) {\scriptsize{$D_3$}};
\draw [thick] (c3) circle [radius=.16]; 
\node at ([shift={(90:.24)}] c3) {$\ast$}; 
\node at ([shift={(150:.24)}] c3) {\scriptsize{$w$}};
\foreach \X in {p01,p02,p03,p04,p11,p12,p21,p22,p23,p24} \fill [black] (\X) circle (1pt);
\coordinate (e) at ([shift={(40:.7)}] c); \draw[thick,blue] (e) ellipse (.15 and .3);
\draw[thick,blue] (p01) to[out=-100, in=90] (p11);
\draw[thick,blue] (p02) to[out=0, in=40] (p03);
\draw[thick,blue] (p12) to[out=-90, in=145] (p21);
\draw[thick,blue] (p22) to[out=15, in=150] (p04);
\draw[thick,blue] (p23) to[out=180, in=270, looseness=4.5] (p24); 
\node at ([shift={(200:.9)}] c) {\scriptsize{$R_1$}}; 
\node at ([shift={(280:.9)}] c) {\scriptsize{$R_2$}};
\node at ([shift={(225:\rbplus)}] c2) {\scriptsize{$R_3$}};
\node at ([shift={(0:.9)}] c) {\scriptsize{$R_4$}}; 
\node at (e) {\scriptsize{$R_5$}};  
\end{tikzpicture}\end{center}
In this planar $2$-tangle $T$:
\begin{itemize}
\item $C_0$ is a shaded $2$-circle whose distinguished interval is shaded white.
\item $C_1$ is a shaded $1$-circle whose distinguished interval is shaded black.
\item $C_2$ is a shaded $2$-circle whose distinguished interval is shaded black.
\item $C_3$ is a shaded $0$-circle whose distinguished interval is shaded white.
\item There are six strings, all drawn in blue. 
\item There are five regions $\{R_1,\ldots,R_5\}$, with (i) $R_1$, $R_3$, and $R_4$ shaded black and (ii) $R_2$ and $R_5$ shaded white.
\end{itemize}
\dqed
\end{example}

Next we equip the set of planar tangles with a given profile with a suitable topology.

\begin{definition}\label{def:planar-tangle-topology}
The set $\PTan\sbinom{C_0}{C_1,\ldots,C_n}$ of planar tangles with a given profile $\sbinom{C_0}{C_1,\ldots,C_n} \in \profscirscir$ is equipped with the following topology.  A planar $k_0$-tangle $T\in\PTan\sbinom{C_0}{C_1,\ldots,C_n}$ is uniquely specified by its set of strings, which is a smooth map \[\nicexy{\coprodover{s} [0,1] \ar[r]^-{St } & D_0\setminus \overset{n}{\coprodover{j=1}} D_j^{\circ}}\] if there are $s$ strings.  The set $\PTan\sbinom{C_0}{C_1,\ldots,C_n}$ is topologized as a subspace of the space of continuous maps\[\nicexy{\coprodover{s \geq 0}\, \overset{s}{\coprodover{i=1}} [0,1] \ar[r] & D_0\setminus \overset{n}{\coprodover{j=1}} D_j^{\circ}},\] which itself has the compact-open topology.  
\end{definition}

The colored symmetric operad structure on planar tangles is defined next.  Recall that $\SCir$ is the set of shaded circles and that $\CHau$ is the symmetric monoidal category of compactly generated weak Hausdorff spaces.  In the following definition, $T\in\PTan\sbinom{C_0}{C_1,\ldots,C_n}$ is as in Definition \ref{def:planar-tangle}.

\begin{definition}\label{def:planar-tangle-operad}
The\index{planar tangle operad}\index{operad!planar tangle} \emph{planar tangle operad} is the $\SCir$-colored symmetric operad $\PTan$ in $\CHau$ defined as follows.
\begin{description}
\item[Entries] For each $\sbinom{C_0}{C_1,\ldots,C_n} \in\profscirscir$, the space $\PTan\sbinom{C_0}{C_1,\ldots,C_n}$ is as in Definition \ref{def:planar-tangle-topology}.
\item[Colored Units] For each $C\in\SCir$, the unique point in $\PTan\sbinom{C}{C}$ is the $C$-colored unit.
\item[Equivariance] The right $S_n$-action is defined as \[T\sigma = \bigl(St,Sh\bigr) \in \PTan\sbinom{C_0}{C_{\sigma(1)},\ldots,C_{\sigma(n)}}\] for $\sigma\in S_n$.  In other words, $\sigma$ permutes the ordering of the inside shaded circles, so the $j$th inside shaded circle in $T\sigma$ is $C_{\sigma(j)}$ for $1\leq j \leq n$. 
\item[Composition] Suppose $1 \leq i \leq n$, and \[T' = \bigl(St', Sh'\bigr) \in \PTan\sbinom{C_i}{C_1',\ldots,C_m'}\] is another planar tangle whose outside shaded circle $C_i$ is the $i$th inside shaded circle in $T$.  The $\compi$-composition is defined as the planar tangle
\[T\compi T' = \bigl(St'', Sh''\bigr)\in \PTan\sbinom{C_0}{C_1,\ldots,C_{i-1},C_1',\ldots,C_m',C_{i+1},\ldots,C_n}\] in which:
\begin{itemize}
\item The strings in $St''$ are the connected components of the topological union of the strings in $T$ and $T'$.
\item The shading \[Sh'' : \Reg(T\compi T')\to \{b,w\}\] is induced by the shadings in $T$ and $T'$ as \[Sh''(R'') = \begin{cases}Sh(R) & \text{if $R'' \cap R \not=\varnothing$ for some $R\in\Reg(T)$},\\
Sh'(R') & \text{if $R'' \cap R' \not=\varnothing$ for some $R'\in\Reg(T')$}.\end{cases}\]
\end{itemize}
\end{description}
\end{definition}

It is an exercise to check that $\PTan$ satisfies the axioms of a $\SCir$-colored symmetric operad in Proposition \ref{prop:symmetric-operad-compi}.

\begin{example}\label{ex:planar-tangle-compi}
Suppose $T \in \PTan\sbinom{C_0}{C_1,C_2,C_3}$ is the planar $2$-tangle in Example \ref{ex:planar-2-tangle}.  Suppose $T' \in \PTan\sbinom{C_1}{C_1',C_2'}$ is the planar $1$-tangle on the left-hand side in the picture below, whose outside shaded circle $C_1$ is the first inside shaded circle in $T$.  It has two inside shaded circles, four strings drawn in green, and three regions, two of which are shaded white.  The $\comp_1$-composition \[T\comp_1 T'\in \PTan\sbinom{C_0}{C_1',C_2',C_2,C_3}\] is depicted on the right-hand side below.
\def\ra{.5} \def\raplus{.58} 
\def\rb{.2} \def\rbplus{.27} 
\def\rc{.15} \def\rcplus{.22} 
\begin{center}\begin{tikzpicture}[scale=2]
\coordinate (c) at (0,0); \draw[thick] (c) circle [radius=1]; 
\node at ([shift={(70:1.15)}] c)  {$T$};
\node at ([shift={(110:1.15)}] c) {$\ast$}; 
\node at ([shift={(150:1.15)}] c) {\scriptsize{$w$}};
\coordinate (p01) at ([shift={(80:1)}] c); \coordinate (p02) at ([shift={(180:1)}] c);
\coordinate (p03) at ([shift={(220:1)}] c); \coordinate (p04) at ([shift={(330:1)}] c);
\coordinate (c1) at ([shift={(135:.4)}] c); \draw [thick] (c1) circle [radius=\ra]; 
\coordinate (p11) at ([shift={(90:\ra)}] c1); 
\coordinate (p12) at ([shift={(270:\ra)}] c1);
\node at ([shift={(0:\raplus)}] c1) {$\ast$}; 
\node at ([shift={(30:\raplus)}] c1) {\scriptsize{$b$}}; 
\node (d1) at ([shift={(225:.25)}] c1) {\scriptsize{$D_1'$}}; 
\draw[thick](d1) circle [radius=\rc];
\node at ([shift={(25:\rcplus)}] d1) {$\ast$};
\node at ([shift={(55:\rcplus)}] d1) {\scriptsize{$w$}};
\coordinate (q11) at ([shift={(90:\rc)}] d1); 
\coordinate (q12) at ([shift={(125:\rc)}] d1);
\coordinate (q13) at ([shift={(270:\rc)}] d1);
\coordinate (q14) at ([shift={(-30:\rc)}] d1);
\node (d2) at ([shift={(45:.25)}] c1) {\scriptsize{$D_2'$}}; 
\draw[thick](d2) circle [radius=\rc];
\node at ([shift={(300:\rcplus)}] d2) {$\ast$};
\node at ([shift={(330:\rcplus)}] d2) {\scriptsize{$b$}};
\coordinate (q21) at ([shift={(180:\rc)}] d2); 
\coordinate (q22) at ([shift={(270:\rc)}] d2);
\foreach \Y in {q11,q12,q13,q14,q21,q22} \fill [black] (\Y) circle (1pt);
\node (U) at ([shift={(180:1.5)}] c) {$T'$};
\coordinate (Ua) at ([shift={(65:.1)}] U); 
\coordinate (Ub) at ([shift={(180:\ra)}] c1); 
\draw [arrow, gray!50, line width=1.5pt, bend left] (Ua) to (Ub);
\coordinate (c2) at ([shift={(-30:.6)}] c); \node at (c2) {\scriptsize{$D_2$}};
\draw [thick] (c2) circle [radius=\rb]; 
\coordinate (p21) at ([shift={(135:\rb)}] c2);  
\coordinate (p22) at ([shift={(15:\rb)}] c2); 
\coordinate (p23) at ([shift={(180:\rb)}] c2);  
\coordinate (p24) at ([shift={(270:\rb)}] c2); 
\node at ([shift={(100:\rbplus)}] c2) {$\ast$};
\node at ([shift={(60:\rbplus)}] c2) {\scriptsize{$b$}};
\coordinate (c3) at ([shift={(250:.7)}] c); \node at (c3) {\scriptsize{$D_3$}};
\draw [thick] (c3) circle [radius=.16]; 
\node at ([shift={(90:.24)}] c3) {$\ast$}; 
\node at ([shift={(150:.24)}] c3) {\scriptsize{$w$}};
\foreach \X in {p01,p02,p03,p04,p11,p12,p21,p22,p23,p24} \fill [black] (\X) circle (1pt);
\coordinate (e) at ([shift={(40:.7)}] c); \draw[thick,blue] (e) ellipse (.15 and .3);
\draw[thick,blue] (p01) to[out=-100, in=90] (p11);
\draw[thick,blue] (p02) to[out=0, in=40] (p03);
\draw[thick,blue] (p12) to[out=-90, in=145] (p21);
\draw[thick,blue] (p22) to[out=15, in=150] (p04);
\draw[thick,blue] (p23) to[out=180, in=270, looseness=4.5] (p24); 
\draw[thick,green] (q11) to[out=90, in=180] (q21); 
\draw[thick,green] (q12) to[out=125, in=270] (p11); 
\draw[thick,green] (q13) to[out=270, in=90] (p12); 
\draw[thick,green] (q14) to[out=-30, in=270] (q22); 
\node at ([shift={(200:.9)}] c) {\scriptsize{$R_1$}}; 
\node at ([shift={(280:.9)}] c) {\scriptsize{$R_2$}};
\node at ([shift={(225:\rbplus)}] c2) {\scriptsize{$R_3$}};
\node at ([shift={(0:.9)}] c) {\scriptsize{$R_4$}}; 
\node at (e) {\scriptsize{$R_5$}};  
\node at ([shift={(0:1.5)}] c) {$\becomes$};
\coordinate (c') at ([shift={(0:3)}] c);  \draw[thick] (c') circle [radius=1]; 
\node at ([shift={(60:1.2)}] c')  {$T \comp_1 T'$};
\node at ([shift={(110:1.15)}] c') {$\ast$}; 
\node at ([shift={(150:1.15)}] c') {\scriptsize{$w$}};
\coordinate (p01') at ([shift={(80:1)}] c'); \coordinate (p02') at ([shift={(180:1)}] c');
\coordinate (p03') at ([shift={(220:1)}] c'); \coordinate (p04') at ([shift={(330:1)}] c');
\coordinate (c1') at ([shift={(135:.4)}] c');
\coordinate (p11') at ([shift={(90:\ra)}] c1'); 
\coordinate (p12') at ([shift={(270:\ra)}] c1');
\node (d1') at ([shift={(225:.25)}] c1') {\scriptsize{$D_1'$}}; 
\draw[thick](d1') circle [radius=\rc];
\node at ([shift={(25:\rcplus)}] d1') {$\ast$};
\node at ([shift={(55:\rcplus)}] d1') {\scriptsize{$w$}};
\coordinate (q11') at ([shift={(90:\rc)}] d1'); 
\coordinate (q12') at ([shift={(125:\rc)}] d1');
\coordinate (q13') at ([shift={(270:\rc)}] d1');
\coordinate (q14') at ([shift={(-30:\rc)}] d1');
\node (d2') at ([shift={(45:.25)}] c1') {\scriptsize{$D_2'$}}; 
\draw[thick](d2') circle [radius=\rc];
\node at ([shift={(300:\rcplus)}] d2') {$\ast$};
\node at ([shift={(330:\rcplus)}] d2') {\scriptsize{$b$}};
\coordinate (q21') at ([shift={(180:\rc)}] d2'); 
\coordinate (q22') at ([shift={(270:\rc)}] d2');
\foreach \Y in {q11',q12',q13',q14',q21',q22'} \fill [black] (\Y) circle (1pt);
\coordinate (c2') at ([shift={(-30:.6)}] c'); \node at (c2') {\scriptsize{$D_2$}};
\draw [thick] (c2') circle [radius=\rb]; 
\coordinate (p21') at ([shift={(135:\rb)}] c2');  
\coordinate (p22') at ([shift={(15:\rb)}] c2'); 
\coordinate (p23') at ([shift={(180:\rb)}] c2');  
\coordinate (p24') at ([shift={(270:\rb)}] c2'); 
\node at ([shift={(100:\rbplus)}] c2') {$\ast$};
\node at ([shift={(60:\rbplus)}] c2') {\scriptsize{$b$}};
\coordinate (c3') at ([shift={(250:.7)}] c'); \node at (c3') {\scriptsize{$D_3$}};
\draw [thick] (c3') circle [radius=.16]; 
\node at ([shift={(90:.24)}] c3') {$\ast$}; 
\node at ([shift={(150:.24)}] c3') {\scriptsize{$w$}};
\foreach \X in {p01',p02',p03',p04',p21',p22',p23',p24'} 
\fill [black] (\X) circle (1pt);
\coordinate (e') at ([shift={(40:.7)}] c'); \draw[thick,red] (e') ellipse (.15 and .3);
\draw[thick,red] (p01') to[out=-100, in=90] (p11');
\draw[thick,red] (p02') to[out=0, in=40] (p03');
\draw[thick,red] (p12') to[out=-90, in=145] (p21');
\draw[thick,red] (p22') to[out=15, in=150] (p04');
\draw[thick,red] (p23') to[out=180, in=270, looseness=4.5] (p24'); 
\draw[thick,red] (q11') to[out=90, in=180] (q21'); 
\draw[thick,red] (q12') to[out=125, in=270] (p11'); 
\draw[thick,red] (q13') to[out=270, in=90] (p12'); 
\draw[thick,red] (q14') to[out=-30, in=270] (q22'); 
\end{tikzpicture}\end{center}
In the planar $2$-tangle $T\comp_1 T'$:
\begin{itemize}
\item There are four inside shaded circles $\{C_1',C_2',C_2,C_3\}$.
\item There are eight strings, all drawn in red.
\item There are six regions, three of which are shaded white.\dqed
\end{itemize}
\end{example}

\chapter{Group Operads}\label{ch:group-operad}

As we have discussed in the previous chapters, a planar operad has no equivariant structure, while a symmetric operad has actions by the symmetric groups.  In this chapter, we discuss a more general version of colored operads, called colored $\G$-operads or group operads, whose equivariant structure comes from a suitable sequence $\G$ of groups.  We will call $\G$ an action operad.  For example, when $\G$ is the planar group operad $\P$ and the symmetric group operad $\S$, we recover colored planar operads and colored symmetric operads, respectively.   Furthermore, as we will discuss in Chapter \ref{ch:braided_operad}, Chapter \ref{ch:ribbon_operad}, and Chapter \ref{ch:cactus}, colored (pure) braided operads, colored (pure) ribbon operads, and colored (pure) cactus operads are also examples of $\G$-operads for suitable choices of $\G$.

Later chapters are all written at the level of colored $\G$-operads.  There are two main benefits to using the language of $\G$-operads.  First, it provides a unifying framework for proving results about operads with various kinds of equivariant structures.  In particular, everything that is true for $\G$-operads is automatically true for planar operads, symmetric operads, (pure) braided operads, (pure) ribbon operads, and (pure) cactus operads.  Second, results about changing the equivariant structures, say from the ribbon groups to the symmetric groups, can be conceptually proved by considering morphisms of action operads.

In Section \ref{sec:augmented-group-operad} we define the objects, called action operads, that can act on planar operads.  In Section \ref{sec:group-operad} we define colored $\G$-operads for an action operad $\G$ as monoids with respect to the colored $\G$-circle product in the category of colored $\G$-sequences, which are objects equipped with $\G$-equivariant structures.   In Section \ref{sec:goperad-generating} we unpack this conceptual definition and characterize colored $\G$-operads in terms of generating operations and generating axioms.  We also discuss algebras over a colored $\G$-operad.  

The remaining sections of this chapter are about examples of action operads and $\G$-operads.  In Section \ref{sec:parenthesized-goperad} we discuss $\G$-operads that are closely related to the planar operad structure of $\G$, called the canonical $\G$-operad and the parenthesized $\G$-operad.  In Section \ref{sec:goperad-example} we provide examples of $\G$-operads related to the translation category construction.

\section{Action Operads}\label{sec:augmented-group-operad}

Suppose $(\Set, \times, \{*\})$ is the symmetric monoidal category of sets with direct product as the monoidal product.  Recall that $S_n$ is the symmetric group on $n$ elements.  Also recall the symmetric group operad $\S$ in Proposition \ref{sn-operad}.  The purpose of this section is to define an action operad, which can act on planar operads as we will discuss in later sections.  For two elements $x,y$ in a group, we write their product as either $xy$ or $x\cdot y$.

\begin{definition}\label{def:augmented-group-operad}
An \emph{action operad}\index{action operad} is a tuple\label{not:g} \[\G = \bigl(\{\G(n)\}_{n\geq 0},\gammag,\operadunit^{\G},\omega\bigr)\] as follows.
\begin{enumerate}
\item Each $\G(n)$ is a group with multiplicative unit $\id_n$.
\item $\omega : \G(n)\to S_n$ is a group homomorphism for each $n$, called the \index{augmentation!of an action operad}\emph{augmentation}.  For an element $g \in \G(n)$, its image $\omega(g) \in S_n$ is called the \emph{underlying permutation}\index{underlying permutation} of $g$ and is also denoted by $\gbar$.
\item $\bigl(\{\G(n)\}_{n\geq 0},\gammag,\operadunit^{\G}\bigr)$ is a one-colored planar operad in $\Set$ as in Definition \ref{def:planar-operad-generating} in which the following statements hold.
\begin{enumerate}[label=(\roman*)]
\item The augmentation \[\omega : \bigl(\{\G(n)\}_{n\geq 0},\gammag,\operadunit^{\G}\bigr) \to \S\] is a morphism of one-colored planar operads in $\Set$.
\item For $\sigma,\rho \in \G(n)$ with $n \geq 1$, $\tau_i,\tau_i'\in \G(k_i)$ with $k_i \geq 0$ for $1\leq i \leq n$, and $k=k_1+\cdots+k_n$, the equality\index{action operad axiom}
\begin{equation}\label{group-operad-identites}
\begin{split} &\gammag\bigl(\sigma\rho;\tau_1\tau_1',\ldots,\tau_n\tau_n'\bigr)\\
&= \gammag\bigl(\sigma; \tau_{\rhobar^{-1}(1)},\ldots,\tau_{\rhobar^{-1}(n)}\bigr) \cdot \gammag\bigl(\rho; \tau_1',\ldots,\tau_n'\bigr)
\end{split}\end{equation}
holds in $\G(k)$, where $\rhobar = \omega(\rho)\in S_n$ is the underlying permutation of $\rho$.
\end{enumerate}
\end{enumerate}
We call \eqref{group-operad-identites} the \emph{action operad axiom}.  We denote such an action operad by $(\G,\omega)$ if we wish to emphasize the augmentation $\omega$.  We will also denote the one-colored planar operad $\bigl(\{\G(n)\}_{n\geq 0},\gammag,\operadunit^{\G}\bigr)$ in $\Set$ by $\G$. 
\end{definition}

The following result contains some consequences of the action operad axiom \eqref{group-operad-identites}.

\begin{lemma}\label{augmented-goperad-identity}
In any action operad $\G$:
\begin{enumerate}
\item The operadic unit $\operadunit^{\G}$ is the multiplicative unit $\id_1\in\G(1)$.
\item The equality \[\gammag\bigl(\id_n; \id_{k_1},\ldots,\id_{k_n}\bigr) = \id_k \in\G(k)\] holds for $n \geq 1$, $k_i\geq 0$ for $1\leq i \leq n$, and $k=k_1+\cdots+k_n$.
\item The operadic composition \[\gammag : \G(1)\times\G(1)\to\G(1)\] coincides with the group multiplication, and $\G(1)$ is an abelian group.
\item In the context of \eqref{group-operad-identites}, the equality 
\begin{equation}\label{gamma-sigma-tau}
\gammag\bigl(\sigma;\tau_1,\ldots,\tau_n\bigr) = \gammag\bigl(\sigma; \id_{k_1},\ldots,\id_{k_n}\bigr) \cdot \gammag\bigl(\id_n; \tau_1,\ldots,\tau_n\bigr)
\end{equation} 
holds in $\G(k)$.
\end{enumerate}
\end{lemma}

\begin{proof}
The first assertion follows from the following computation in $\G(1)$:
\[\begin{split} \id_1 &= \id_1 \cdot \id_1\\
&=\gammag\bigl(\id_1; \operadunit^{\G}\bigr) \cdot \gammag\bigl(\operadunit^{\G}; \id_1\bigr)\\
&= \gammag\bigl(\id_1 \cdot \operadunit^{\G}; \operadunit^{\G} \cdot \id_1\bigr)\\
&= \gammag(\operadunit^{\G}; \operadunit^{\G}\bigr)\\ &= \operadunit^{\G}.
\end{split}\]
The first and the fourth equalities use the fact that $\id_1$ is the multiplicative unit in $\G(1)$.  The second and the last equalities use the fact that $\operadunit^{\G}$ is the operadic unit in $\G$.  The third equality follows from the action operad axiom \eqref{group-operad-identites} and the fact that $\omega(\operadunit^{\G}) = \id_1\in S_1$ is the identity permutation.

The second assertion follows from the following computation in $\G(k)$:
\[\begin{split}
&\gammag\bigl(\id_n; \id_{k_1},\ldots,\id_{k_n}\bigr) \cdot \gammag\bigl(\id_n; \id_{k_1},\ldots,\id_{k_n}\bigr) \\
&= \gammag\bigl(\id_n \cdot \id_n; \id_{k_1} \cdot \id_{k_1},\ldots,\id_{k_n}\cdot\id_{k_n}\bigr)\\
&= \gammag\bigl(\id_n; \id_{k_1},\ldots,\id_{k_n}\bigr). 
\end{split}\]
The first equality follows from the action operad axiom \eqref{group-operad-identites} and the fact that $\omega(\id_n)=\id_n\in S_n$ is the identity permutation.  The second equality uses the fact that each $\id_j\in\G(j)$ is the multiplicative unit.

For the third assertion, observe that the action operad axiom \eqref{group-operad-identites} implies that the operadic composition \[\gammag : \G(1)\times\G(1)\to\G(1)\] preserves the group multiplication in $\G(1)$.  Since the operadic unit $\operadunit^{\G}$ is also the multiplicative unit $\id_1\in \G(1)$, the classical Eckmann-Hilton argument \cite{eh} establishes the assertion.

Finally, the equality \eqref{gamma-sigma-tau} is obtained from the action operad axiom \eqref{group-operad-identites} by replacing $\rho$, $\tau_i$, and $\tau_i'$ with $\id_n$, $\id_{k_i}$, and $\tau_i$, respectively.
\end{proof}

\begin{example}[Planar Group Operad]\label{ex:trivial-group-operad}
The \emph{planar group operad}\index{planar group operad}\index{group operad!planar} $\P$ is the action operad with $\P(n)=\{*\}$  the trivial group for $n \geq 0$.\dqed\end{example}

\begin{example}[Symmetric Group Operad]\label{ex:symmetric-group-operad}
The symmetric group operad $\S$ in Proposition \ref{sn-operad} is an action operad when equipped with the identity augmentation.\dqed\end{example}

More examples of action operads will be given in Chapters \ref{ch:braided_operad}, \ref{ch:ribbon_operad}, and \ref{ch:cactus}.  Several remarks regarding the definition of an action operad follow.

\begin{remark}[Some History]
Wahl \cite{wahl} (Section 1.2.0.2) considered a concept similar to an action operad, the difference being that Wahl required each augmentation $\omega : \G(n)\to S_n$ to be surjective.  What we call an action operad is what Yoshida \cite{yoshida} and Zhang \cite{zhang} called a \emph{group operad}.  However, we want to reserve the name group group for a $\G$-operad as we will define later.  Yoshida \cite{yoshida} established a fully faithful functor from the category of action operads to the category of \index{crossed interval group}crossed interval groups in the sense of Batanin-Markl \cite{batanin-markl}.  We follow Corner-Gurski \cite{corner-gurski} and Gurski \cite{gurski} in using the name action operad, and Lemma \ref{augmented-goperad-identity}(1)-(3) above is Lemma 1.21 in \cite{corner-gurski}.  It was also pointed out in \cite{gurski} that an action operad can act on a multicategory, which is what we call a colored $\G$-operad in $\Set$ below.\dqed
\end{remark}

\begin{remark}[Augmentation]
The augmentation $\omega : \G(1) \to S_1$ automatically preserves the operadic units.  So the condition that $\omega : \G\to\S$ is a morphism of one-colored planar operads in $\Set$ means that the diagram \[\nicexy{\G(n) \times \G(m) \ar[d]_-{(\omega,\omega)} \ar[r]^-{\compi} & \G(n+m-1) \ar[d]^-{\omega}\\ S_n \times S_m \ar[r]^-{\compi} & S_{n+m-1}}\] is commutative for $1 \leq i \leq n$ and $m \geq 0$.  Equivalently, the diagram \[\nicexy{\G(n) \times \G(k_1) \times \cdots \times \G(k_n) \ar[d]_-{(\omega,\ldots,\omega)} \ar[r]^-{\gamma^{\G}} & \G(k_1+\cdots+k_n) \ar[d]^-{\omega}\\ S_n \times S_{k_1} \times \cdots \times S_{k_n} \ar[r]^-{\gammas} & S_{k_1+\cdots+k_n}}\] is commutative for $n \geq 1$ and $k_1,\ldots,k_n \geq 0$.\dqed
\end{remark}

\begin{remark}
We are \emph{not} asking $\G$ to be a one-colored planar operad in the symmetric monoidal category of groups.  In other words, the operadic composition \[\nicexy{\G(n)\times \G(k_1)\times\cdots \times \G(k_n) \ar[r]^-{\gammag} & \G(k)}\] is a function, not a group homomorphism.  Indeed, choosing $\tau_i=\tau_i'=\id_{k_i}$, the action operad axiom \eqref{group-operad-identites} becomes \[\gammag\bigl(\sigma\rho; \id_{k_1},\ldots,\id_{k_n}\bigr) = \gammag\bigl(\sigma; \id_{k_{\rhobar^{-1}(1)}}, \ldots,\id_{k_{\rhobar^{-1}(n)}}\bigr) \cdot \gammag\bigl(\rho; \id_{k_1},\ldots,\id_{k_n}\bigr).\]  Therefore, the map \[\nicexy@C+1.5cm{\G(n)\ar[r]^-{\gammag(-;\id_{k_1},\ldots,\id_{k_n})} & \G(k)}\] is \emph{not} multiplicative.

On the other hand, the map \[\nicexy@C+1cm{\G(k_1)\times\cdots\times\G(k_n) \ar[r]^-{\gammag(\id_n;-,\ldots,-)} &\G(k)}\] is a group homomorphism.  This is a consequence of the action operad axiom \eqref{group-operad-identites} and Lemma \ref{augmented-goperad-identity}(2).\dqed
\end{remark}

\section{Group Operads as Monoids}\label{sec:group-operad}

Fix an action operad $(\G,\omega)$.  We now follow the development in Section \ref{sec:symmetric_operad} with the symmetric group operad $\S$ replaced by the action operad $\G$.  First we define the $\G$-analogue of symmetric sequences.  As before, $(\M,\otimes,\tensorunit)$ is our ambient symmetric monoidal category.  We always assume that $\M$ is complete and cocomplete and that the monoidal product commutes with colimits on each side.  An initial object in $\M$ is denoted by $\varnothing$.

\begin{definition}\label{def:g-sequences}
Suppose $\colorc$ is a set.
\begin{enumerate}
\item Define the groupoid\label{not:gsubc}\index{groupoid of profiles!$\G$-action} $\gsubc$ with object set $\Profc$.  For $\uc=(c_1,\ldots,c_m),\ud \in \Profc$, an isomorphism $\sigma : \uc \iso \ud \in \gsubc$ is an element $\sigma \in \G(m)$ such that \[\sigmabar\uc = \bigl(c_{\inv{\sigmabar}(1)},\ldots,c_{\inv{\sigmabar}(m)}\bigr) = \ud \in \Profc,\] where $\sigmabar = \omega(\sigma) \in S_{m}$ is the underlying permutation of $\sigma$.  Composition and identity morphisms are given by multiplication and units in the groups $\G(n)$. 
\item The \index{opposite groupoid}opposite groupoid $\gsubcop$ is regarded as the groupoid of $\colorc$-profiles in which elements in $\G$ act on $\colorc$-profiles from the right via their underlying permutations.
\item The objects of the diagram category\label{notation:gseqcm} \[\gseqcm = \M^{\gsubcopc}\] are called \index{G-sequence@$\G$-sequence}\emph{$\colorc$-colored $\G$-sequences} in $\M$.  
\item Define the $\colorc$-colored $\G$-sequence\index{unit G-operad@unit $\G$-operad} $\I^{\G}$ with entries
\begin{equation}\label{g-unit}
\I^{\G}\duc = \begin{cases}\varnothing & \text{if $\uc\not=(d)$},\\
\coprodover{\G(1)}\tensorunit & \text{if $\uc=(d)$}\end{cases}
\end{equation} 
for $\duc\in \Profcc$.  Its $\colorc$-colored $\G$-sequence structure is induced by the group multiplication in $\G(1)$.
\end{enumerate}
\end{definition}

\begin{example}
For the planar group operad $\P$, the groupoid $\P_{\colorc}$ is the discrete category with object set $\Profc$.  For the symmetric group operad $\S$, the groupoid $\S_{\colorc}$ is the groupoid of $\colorc$-profiles in Definition \ref{def:colored-sequences}.\dqed
\end{example}

Next we define the $\G$-version of the symmetric circle product with respect to which $\G$-operads are monoids.

\begin{definition}\label{def:g-circle-product}
Suppose $\G$ is an action operad, and $X,Y  \in \gseqcm$.
\begin{enumerate}
\item For each $\uc = (c_1,\ldots,c_m) \in \gsubc$, define the object $Y^{\uc}_{\smallsubg} \in \M^{\gsubcop}$ entrywise as the coend
\[Y^{\uc}_{\smallsubg}(\ub) = \int^{\{\ua_j\}\in\overset{m}{\prodover{j=1}} \gsubcop} \gsubcop\bigl(\ua_1,\ldots,\ua_m;\ub\bigr) \cdot \left[\bigotimes_{j=1}^m Y\cjuaj\right] \in \M\]
for $\ub \in \gsubcop$, in which $(\ua_1,\ldots,\ua_m) \in \gsubcop$ is the concatenation.  The object $Y^{\uc}_{\smallsubg}$ is natural in $\uc \in \gsubc$ via left permutations of the tensor factors in $\bigotimes_{j=1}^m Y\cjuaj$ by underlying permutations of elements in $\G(m)$.
\item The \index{circle product!G@$\G$}\index{G-circle product@$\G$-circle product}\emph{$\colorc$-colored $\G$-circle product}\label{notation:circg} \[X \circg Y \in \gseqcm\] is defined entrywise as the coend
\[(X \circg Y)\dub = \int^{\uc \in \gsubc} X\duc \otimes Y^{\uc}_{\smallsubg}(\ub)\in \M\]
for $\dub \in \gsubcopc$.
\end{enumerate}
\end{definition}

The following result is the $\G$-version of Proposition \ref{planar-circle-product-monoidal}, whose proof is readily adapted to the current case.

\begin{proposition}\label{g-circle-product-monoidal}
$\bigl(\gseqcm, \circg, \I^{\G}\bigr)$ is a \index{monoidal category!of $\G$-sequences}\index{G-sequence@$\G$-sequence!monoidal category}monoidal category.
\end{proposition}

We now define colored $\G$-operads as monoids in this monoidal category.

\begin{definition}\label{def:g-operad}
For an action operad $(\G,\omega)$, the category of \index{operad!G@$\G$-}\index{G-operad@$\G$-operad}\emph{$\colorc$-colored $\G$-operads in $\M$} is the category \[\goperadcm = \Mon\bigl(\gseqcm, \circg, \I^{\G}\bigr)\] of monoids in the monoidal category $\bigl(\gseqcm,\circg,\I^{\G}\bigr)$.  An object in $\goperadcm$ is also called a \index{group operad}\emph{group operad}.  If $\colorc$ has $n<\infty$ elements, we refer to objects in $\goperadcm$ as \emph{$n$-colored $\G$-operads in $\M$}.
\end{definition}

\begin{remark}[One-Colored Case]
In the one-colored case (i.e, when $\colorc=\{*\}$) with $\M=\Set$, the one-colored $\G$-circle product and Proposition \ref{g-circle-product-monoidal} above reduce to Definition 1.19 and Theorem 1.20 in Corner-Gurski \cite{corner-gurski}.  For our later discussion of infinity $\G$-operads, it is crucial that we be able to consider \emph{colored} $\G$-operads.\dqed
\end{remark}

\begin{remark}\label{rk:lms-goperad}
The reader should not confuse our $\colorc$-colored $\G$-operads in Definition \ref{def:g-operad} with the $G$-operads in \cite{gut-white,lms,salvatore-wahl}, which mean one-colored symmetric operads in $\Top^G$ for a group $G$.  Here $\Top$ is a suitable category of pointed spaces, and $\Top^G$ is the category of $G$-equivariant spaces and $G$-equivariant maps.\dqed 
\end{remark}

\begin{example}[Planar and Symmetric Operads]\label{ex:g-operad}
For the planar group operad $\P$, the category of $\colorc$-colored $\P$-operads is the category of $\colorc$-colored planar operads.  For the symmetric group operad $\S$ with identity augmentation,  the category of $\colorc$-colored $\S$-operads is the category of $\colorc$-colored symmetric operads.\dqed
\end{example}

\begin{example}[Braided, Ribbon, and Cactus Operads]\label{ex:braided-ribbon-operads}
As we will see in Chapter \ref{ch:braided_operad}, the (pure) braid groups form an action operad $\B$ ($\PB$) such that $\colorc$-colored $\B$-operads ($\PB$-operads) are $\colorc$-colored (pure) braided operads.  Similarly, in Chapter \ref{ch:ribbon_operad} we will show that the (pure) ribbon groups form an action operad $\R$ ($\PR)$ such that $\colorc$-colored $\R$-operads ($\PR$-operads) are $\colorc$-colored (pure) ribbon operads.  In Chapter \ref{ch:cactus} we will show that the (pure) cactus groups form an action operad $\Cac$ (resp., $\PCac$) such that $\colorc$-colored $\Cac$-operads ($\PCac$-operads) are $\colorc$-colored (pure) cactus operads.\dqed
\end{example}

\begin{remark}[Non-Examples]\label{rk:virtual-braid}
While the braid groups form an action operad, there are several sequences of groups related to the braid groups that do \emph{not} form action operads.  For instance,
\begin{enumerate}
\item the sequence of \index{virtual braid group}\index{braid group!virtual}virtual braid groups in \cite{kl}, 
\item the sequence of \index{welded braid group}\index{braid group!welded}welded braid groups, called the braid-permutation groups in \cite{fenn}, and
\item the sequence of \index{Grothendieck cartographical group}Grothendieck cartographical groups \cite{sv}, called the \index{twin group}twin groups in \cite{kho1,kho2},
\end{enumerate}
are not action operads because they are not one-colored planar operads in $\Set$.\dqed
\end{remark}

\section{Coherence for Group Operads}\label{sec:goperad-generating}

For a given action operad $\G$ as in Definition \ref{def:augmented-group-operad}, $\G$-operads can be described more explicitly as follows.  The proofs are minor modifications of those in the planar and symmetric cases.  The following description of $\G$-operads is the $\G$-version of Proposition \ref{prop:planar-operad-same-def} and Proposition \ref{def:symmetric-operad-generating}.

\begin{proposition}\label{def:g-operad-generating}
A\index{coherence!G-operad@$\G$-operad}\index{G-operad@$\G$-operad!coherence} $\colorc$-colored $\G$-operad in $\M$ is equivalent to a $\colorc$-colored planar operad $(\O,\gamma,\operadunit)$ as in Definition \ref{def:planar-operad-generating} together with a $\colorc$-colored $\G$-sequence structure on $\O$ that satisfies the following equivariance axioms.  Suppose that in \eqref{operadic-composition} $|\ub_j| = k_j \geq 0$.
\begin{enumerate}
\item For each element $\sigma \in \G(n)$, the \index{equivariance!G-operad@$\G$-operad}\index{G-operad@$\G$-operad!equivariance}\emph{top equivariance diagram} 
\[\nicexy@C+.3cm{\O\duc \otimes \bigotimes\limits_{j=1}^n \O\cjubj 
\ar[d]_-{\gamma} \ar[r]^-{(\sigma, \sigmabar^{-1})}& \O\ducsigmabar \otimes 
\bigotimes\limits_{j=1}^n \O\csigmabarjubsigmabarj \ar[d]^-{\gamma}\\
\O\duboneubn \ar[r]^-{\sigma\langle k_{\sigmabar(1)}, \ldots , k_{\sigmabar(n)}\rangle}& \O\dubsigmabaroneubsigmabarn}\]
in $\M$ is commutative.  In the top horizontal morphism, $\sigma$ is the $\G(n)$-equivariant structure morphism of $\O$ corresponding to $\sigma \in \G(n)$, and $\sigmabar \in S_n$ is the underlying permutation of $\sigma$.  The bottom horizontal morphism is the $\G(k)$-equivariant structure morphism of $\O$ corresponding to the element 
\begin{equation}\label{g-block-perm}
\sigma\langle k_{\sigmabar(1)},\ldots,k_{\sigmabar(n)}\rangle = \gamma^{\G}\bigl(\sigma; \id_{k_{\sigmabar(1)}},\ldots,\id_{k_{\sigmabar(n)}}\bigr) \in \G(k),
\end{equation} where $k=k_1+\cdots+k_n$.
\item Given elements $\tau_j \in \G(k_j)$ for $1 \leq j \leq n$, the \emph{bottom equivariance diagram}
\[\nicexy@C+.3cm{\O\duc \otimes \bigotimes\limits_{j=1}^n \O\cjubj\ar[d]_-{\gamma} \ar[r]^-{(\Id, \otimes \tau_j)}& \O\duc \otimes \bigotimes\limits_{j=1}^n \O\cjubjtaubarj\ar[d]^-{\gamma} \\
\O\duboneubn \ar[r]^-{\tau_1 \oplus \cdots \oplus \tau_n}& \O\dubonetaubaroneubntaubarn}\]
in $\M$ is commutative.  In the top horizontal morphism, each $\tau_j$ is the $\G(k_j)$-equivariant structure morphism of $\O$ corresponding to $\tau_j \in \G(k_j)$.  The bottom horizontal morphism is the $\G(k)$-equivariant structure morphism of $\O$ corresponding to the element 
\begin{equation}\label{g-block-sum}
\tau_1 \oplus \cdots \oplus \tau_n = \gamma^{\G}\bigl(\id_n; \tau_1,\ldots,\tau_n\bigr) \in \G(k).
\end{equation}
\end{enumerate}

Moreover, a morphism of $\colorc$-colored $\G$-operads in $\M$ is equivalent to a morphism of the underlying $\colorc$-colored planar operads that is also a morphism of the $\colorc$-colored $\G$-sequences.
\end{proposition}

\begin{remark}
In the one-colored case (i.e., when $\colorc = \{*\}$), the characterization in  Proposition \ref{def:g-operad-generating} is used as the definition of a one-colored $\G$-operad in Corner-Gurski \cite{corner-gurski}, Yoshida \cite{yoshida}, and Zhang \cite{zhang}.  In the general colored case, the characterization in  Proposition \ref{def:g-operad-generating} is used in Gurski \cite{gurski} as the definition of a colored $\G$-operad in $\Set$, called a \index{multicategory}$\Lambda$-multicategory.\dqed
\end{remark}

\begin{example}[Changing Color Sets]\label{ex:g-operad-change-colors}
Suppose\index{change of colors} $f : \colorc \to \colord$ is a function between sets.  It induces a functor \[f_* : \gsubc \to \gsubd,\] where $f$ is applied entrywise.  By composition it induces a monoidal functor \[\nicexy{\gseqcm = \M^{\gsubcopc} & \M^{\gsubdopd} = \gseqdm \ar[l]_-{f^*}}\] for the monoidal structure in Proposition \ref{g-circle-product-monoidal}.  Therefore, it induces a functor between their respective categories of monoids, i.e., a functor \[\nicexy{\goperadcm & \goperaddm \ar[l]_-{f^*}}\] from $\colord$-colored $\G$-operads to $\colorc$-colored $\G$-operads in $\M$.  Specifically, for a $\colord$-colored $\G$-operad $\P$, its image $f^*\P \in \gopcm$ has entries \[(f^*\P)\duc = \P\fdfuc\] for $\duc \in \Profcc$, where \[f\uc = (fc_1,\ldots,fc_n) \ifspace \uc = (c_1,\ldots,c_n)\in \Profc.\]  Its $\colorc$-colored $\G$-operad structure is directly induced by the $\colord$-colored $\G$-operad structure of $\P$.  For example, in the setting of Definition \ref{def:planar-operad-generating}, the operadic composition in $f^*\P$, \[\nicexy{(f^*\P)\duc \otimes \bigotimes\limits_{j=1}^n (f^*\P)\cjubj \ar[r]^-{\gamma^{f^*\P}} & (f^*\P)\dub\\ 
\P\fdfuc \otimes \bigotimes\limits_{j=1}^n \P\fcjfubj \ar@{=}[u] \ar[r]^-{\gamma^{\P}} & \P\fdfub \ar@{=}[u]}\] is the operadic composition $\gamma^{\P}$ in $\P$.  We call each functor $f^*$ a \emph{change-of-color functor}.\dqed
\end{example}

The following description of $\G$-operads is the $\G$-version of Proposition \ref{prop:planar-operad-defs-equivalent} and Proposition \ref{prop:symmetric-operad-compi}.

\begin{proposition}\label{prop:g-operad-compi}
A $\colorc$-colored $\G$-operad in $\M$ is equivalent to a $\colorc$-colored planar operad $(\O,\circ,\operadunit)$ as in Definition \ref{def:planar-operad-compi} together with a $\colorc$-colored $\G$-sequence structure on $\O$ that satisfies the following equivariance axiom.  Suppose $|\uc| = n \geq 1$, $1 \leq i \leq n$, $\sigma \in \G(n)$, and $\tau \in \G(m)$.  Then the \emph{equivariance diagram} in $\M$\index{G-operad@$\G$-operad!equivariance}\index{equivariance!G-operad@$\G$-operad}
\[\nicexy{\O\duc \otimes \O\sbinom{c_{\sigmabar(i)}}{\ub}\ar[d]_-{(\sigma,\tau)} \ar[r]^-{\comp_{\sigmabar(i)}}& 
\O\sbinom{d}{\uc\comp_{\sigmabar(i)}\ub} \ar[d]^-{\sigma\compi\tau}\\
\O\ducsigmabar \otimes \O\sbinom{c_{\sigmabar(i)}}{\ub\taubar}\ar[r]^-{\compi} & 
\O\sbinom{d}{(\uc\sigmabar)\compi(\ub\taubar)}
= \O\sbinom{d}{(\uc\comp_{\sigmabar(i)}\ub)(\sigmabar\compi\taubar)}}\]
is commutative.  The right vertical morphism is the $\G(n+m-1)$-equivariant structure morphism of $\O$ corresponding to the element
\begin{equation}\label{g-operad-compi}
\sigma\compi\tau = \gammag\bigl(\sigma; \underbrace{\id_1,\ldots,\id_1}_{i-1},\tau,\underbrace{\id_1,\ldots,\id_1}_{n-i}\bigr) \in \G(n+m-1).
\end{equation}

Moreover, a morphism of $\colorc$-colored $\G$-operads in $\M$ is equivalent to a morphism of the underlying $\colorc$-colored planar operads that is also a morphism of the $\colorc$-colored $\G$-sequences.
\end{proposition}

\begin{remark}[Further Characterizations of $\G$-Operads]
Later we will provide two more characterizations of colored $\G$-operads.
\begin{enumerate}
\item In Theorem \ref{gopcm-algebra} below, we will characterize $\colorc$-colored $\G$-operads as algebras over a $(\Profcc)$-colored symmetric operad $\gopcm$.  
\item In Corollary \ref{goperad-gammat} below, we will describe $\colorc$-colored $\G$-operads using structure morphisms $\gamma_{(T,\sigma)}$ that are parametrized by all $\colorc$-colored $\G$-trees.   This description will be important when we discuss the Boardman-Vogt construction of colored $\G$-operads.\dqed
\end{enumerate}
\end{remark}

Next we discuss algebras over a $\G$-operad.

\begin{lemma}\label{lem:g-operad-monad}
Suppose $\O$ is a $\colorc$-colored $\G$-operad in $\M$.  Then it induces a \index{G-operad@$\G$-operad!induced monad}\index{monad!induced by a $\G$-operad}monad whose functor is \[\O \circg - : \Mtoc \to \Mtoc\] and whose multiplication and unit are induced by those of $\O$ as in Corollary \ref{planar-operad-monoid}.
\end{lemma}

\begin{definition}\label{def:g-operad-algebra}
Suppose $\O$ is a $\colorc$-colored $\G$-operad in $\M$.  The category $\algmo$ \index{G-operad@$\G$-operad!algebra}\index{algebra!of a G-operad@of a $\G$-operad}of \emph{$\O$-algebras} is defined as the category of $(\O\circg -)$-algebras for the monad $\O\circg -$ in $\Mtoc$.
\end{definition}

\begin{remark}\label{rk:o-circg-x}
For a $\colorc$-colored $\G$-operad $\O$ in $\M$ and $X \in \Mtoc$, a typical entry of $\O\circg X \in \Mtoc$ is the coend \[(\O\circg X)_d = \int^{\uc\in\gsubc} \O\duc\otimes X_{\uc} \in\M\] for $d\in \colorc$.\dqed
\end{remark}

\begin{proposition}\label{prop:g-operad-algebra-defs}
Suppose\index{coherence!G-operad algebra@$\G$-operad algebra} $\O$ is a $\colorc$-colored $\G$-operad in $\M$.  Then an $\O$-algebra is precisely a pair $(X,\lambda)$ as in Definition \ref{def:planar-operad-algebra-generating} that satisfies the following equivariance axiom.  For each $\duc \in \gsubcopc$ and each element $\sigma \in \G(|\uc|)$, the equivariance diagram\index{equivariance!G-operad algebra@$\G$-operad algebra}
\[\nicexy{\O\duc \otimes X_{\uc}\ar[d]_-{\lambda} \ar[r]^-{(\sigma,\inv{\sigmabar})}& \O \ducsigmabar \otimes X_{\uc\sigmabar}\ar[d]^-{\lambda}\\ X_d \ar@{=}[r] & X_d}\]
in $\M$ is commutative.   In the top horizontal morphism, $\inv{\sigmabar}$ is the right permutation on the tensor factors in $X_{\uc}$ induced by the underlying permutation $\sigmabar \in S_{|\uc|}$.
\end{proposition}

\begin{example}[Symmetric Operads as $\G$-Operads]\label{ex:sop-gop}
Every\index{G-operad@$\G$-operad!induced by a symmetric operad}\index{symmetric operad!inducing a G-operad@inducing a $\G$-operad} $\colorc$-colored symmetric operad is naturally a $\colorc$-colored $\G$-operad, in which elements in $\G$ act via their underlying permutations.  The underlying $\colorc$-colored planar operad structure remains unchanged.\dqed\end{example}

\begin{example}[Colored Endomorphism $\G$-Operads]\label{ex:end-g}
For each $\colorc$-colored object $X=\{X_c\}_{c\in\colorc}$ in $\M$, recall from Example \ref{ex:endomorphism-operad} the endomorphism\index{G-operad@$\G$-operad!endomorphism}\index{endomorphism operad!G-@$\G$} symmetric operad $\End(X)$, which is a $\colorc$-colored symmetric operad.  By Example \ref{ex:sop-gop} we may regard $\End(X)$ as a $\colorc$-colored $\G$-operad.  Then for each $\colorc$-colored $\G$-operad $\O$ in $\M$, an $\O$-algebra structure on $X$ in the sense of Definition \ref{def:g-operad-algebra} is equivalent to a morphism $\O \to \End(X)$ of $\colorc$-colored $\G$-operads.\dqed\end{example}

\section{Parenthesized Group Operads}\label{sec:parenthesized-goperad}

In this section, we discuss examples of $\G$-operads that are related to or directly inherited from the planar operad structure of an action operad $\G$.

\begin{example}[Action Operads as $\G$-Operads in $\Set$]\label{ex:agop-goperad}
Each action operad $\G$ yields a one-colored $\G$-operad $\bigl(\{\G(n)\}_{n\geq 0},\gammag,\operadunit\bigr)$ in $\Set$ whose $\G$-sequence structure comes from the level-wise group multiplication on each $\G(n)$.  We call this one-colored $\G$-operad in $\Set$ the \index{G-operad@$\G$-operad!canonical}\index{canonical G-operad@canonical $\G$-operad}\emph{canonical $\G$-operad}, denoted also by $\G$.  This is Proposition 1.17 in \cite{gurski}.

Indeed, by definition $\G$ has an underlying one-colored planar operad in $\Set$.  Level-wise group multiplication defines a one-colored $\G$-sequence structure on $\{\G(n)\}$.  The top equivariance diagram in Proposition \ref{def:g-operad-generating} is the equality 
\[\begin{split}
&\gammag\bigl(\sigma;\tau_1,\ldots,\tau_n\bigr) \cdot \gammag\bigl(\rho; \id_{k_{\rhobar(1)}},\ldots,\id_{k_{\rhobar(n)}}\bigr)\\
&= \gammag\bigl(\sigma\rho; \tau_{\rhobar(1)},\ldots, \tau_{\rhobar(n)}\bigr) \in \G(k)
\end{split}\] 
for $\sigma,\rho\in \G(n)$ and $\tau_i\in\G(k_i)$ for $1\leq i \leq n$.  Similarly, the bottom equivariance diagram in Proposition \ref{def:g-operad-generating} is the equality 
\[\begin{split}
&\gammag\bigl(\sigma;\tau_1,\ldots,\tau_n\bigr) \cdot \gammag\bigl(\id_n; \rho_1,\ldots,\rho_n\bigr)\\
&= \gammag\bigl(\sigma; \tau_1\rho_1, \ldots, \tau_n\rho_n\bigr) \in \G(k)
\end{split}\] 
for $\sigma\in \G(n)$ and $\tau_i,\rho_i\in \G(k_i)$ for $1 \leq i \leq n$.  Moreover, the equivariance diagram in Proposition \ref{prop:g-operad-compi} is the equality 
\[\begin{split}
&\gammag\bigl(\sigma; \id_{1}^{\rhobar^{-1}(i)-1}, \tau, \id_{1}^{n-\rhobar^{-1}(i)}\bigr)
\cdot \gammag\bigl(\rho; \id_{1}^{i-1},\pi,\id_{1}^{n-i}\bigr)\\
&= \gammag\bigl(\sigma\rho; \id_{1}^{i-1},\tau\pi, \id_{1}^{n-i}\bigr)\in \G(n+m-1)
\end{split}\] 
for $\sigma,\rho\in \G(n)$ and $\tau,\pi\in\G(m)$, where \[\id_1^{j} = (\id_1,\ldots,\id_1)\in\G(1)^{\times j}.\]  These three equalities hold by the action operad axiom \eqref{group-operad-identites}.\dqed
\end{example}

\begin{example}[Parenthesized $\G$-Operad]\label{ex:parenthesized-goperad}
This example is an extension of both the magma operad $\Mag$ in Example \ref{ex:mag-planar-operad} and the canonical $\G$-operad $\G$ in $\Set$ in Example \ref{ex:agop-goperad}.  Suppose $\G$ is an action operad.  

There is a one-colored\index{parenthesized G-operad@parenthesized $\G$-operad} $\G$-operad\label{not:pag} $\PaG$ in $\Set$ defined as follows.
\begin{description}
\item[Entries] For $n\geq 0$, define the set \[\PaG(n) = \Mag(n) \times \G(n) \times \Mag(n),\] in which an element is written as $(w;\sigma;w')$ with $w,w'\in \Mag(n)$ and $\sigma \in \G(n)$.  We visualize such an element as\[\nicexy{w'\\ w \ar[u]^-{\sigma}}\] with $w$ the \emph{input parenthesized word} and $w'$ the \emph{output parenthesized word}, respectively.
\item[Unit] The operadic unit is \[(*;\id_1;*) \in \Mag(1)\times\G(1)\times\Mag(1).\]
\item[Composition] For $1 \leq i \leq n$ and $m \geq 0$, the $\compi$-composition \[\nicexy{\PaG(n) \times \PaG(m) \ar[r]^-{\compi} & \PaG(n+m-1)}\] is defined as \[(w;\sigma;w') \compi (x;\tau;x') = \bigl(w\compi x; \sigma\compi\tau; w' \comp_{\sigmabar(i)} x'\bigr)
\] for $w,w'\in \Mag(n)$, $x,x' \in \Mag(m)$, $\sigma\in\G(n)$, and $\tau\in \G(m)$.  Here \[w\compi x \andspace w' \comp_{\sigmabar(i)} x' \in \Mag(n+m-1)\] are the $\compi$-composition in the magma operad $\Mag$ with $\sigmabar \in S_n$ the underlying permutation of $\sigma$.  Likewise, \[\sigma\compi\tau\in \G(n+m-1)\] is the $\compi$-composition in the planar operad $\G$.
\item[$\G$-Equivariance] The one-colored $\G$-sequence structure \[\nicexy{\PaG(n)\times\G(n)\ar[r] & \PaG(n)}\] is level-wise defined as \[(w;\sigma;w') \tau = (w;\sigma\tau;w')\] for $\sigma, \tau\in \G(n)$.
\end{description}
That $\PaG$ is a one-colored $\G$-operad in $\Set$ follows from the fact that the magma operad $\M$ is a one-colored planar operad and that $\G$ is a one-colored $\G$-operad in $\Set$.  We call $\PaG$ the \emph{parenthesized $\G$-operad}.

Let us consider some special cases of the parenthesized $\G$-operad.  
\begin{enumerate}
\item For the planar group operad $\P$ in Example \ref{ex:trivial-group-operad}, there is an isomorphism \[\PaP \cong \Mag\times\Mag\] of one-colored planar operads in $\Set$ between the parenthesized planar operad and the product of two copies of the magma operad.  
\item
For the symmetric group operad $\S$ in Proposition \ref{sn-operad}, the elements in the parenthesized symmetric operad $\PaS$ are called \index{parenthesized permutation}\emph{parenthesized permutations} in Bar-Natan \cite{barnatan}, the set of which is denoted by \textbf{PaP} there.  On the other hand, our parenthesized symmetric operad $\PaS$ is different from Fresse's (Section 6.3 in \cite{fresse-gt}) parenthesized symmetry operad \textit{\textsf{PaS}}, which is a one-colored symmetric operad in $\Cat$ whose algebras are small symmetric monoidal categories without monoidal units.  See Remark \ref{rk:fresse-pas} for more discussion of Fresse's parenthesized symmetry operad.
\item We will discuss the braided and ribbon versions of the parenthesized $\G$-operad in Example \ref{ex:parenthesized-braid} and Example \ref{ex:parenthesized-ribbon} below.
\end{enumerate}

There is a morphism $\pi : \PaG \to \G$ of one-colored $\G$-operads in $\Set$, where the target $\G$ is the canonical $\G$-operad in Example \ref{ex:agop-goperad}, that is level-wise the projection \[\nicexy{\Mag(n)\times\G(n)\times\Mag(n) \ar[r]^-{\pi} & \G(n)}.\]\dqed
\end{example}

\section{Group Operads from Translation Categories}\label{sec:goperad-example}

In this section, we discuss some examples of $\G$-operads that are related to the translation category construction.

\begin{definition}\label{def:translation-category}
For a group $G$, its \index{translation category}\emph{translation category} $E_G$ is defined as the small groupoid with:
\begin{itemize}
\item object set the underlying set of $G$;
\item morphism set $E_G(g,h) = \{hg^{-1}\}$, a one-point set, for $g,h\in G$;
\item composition given by multiplication in $G$.  
\end{itemize}
As $G$ acts on itself from the right via group multiplication, $G$ also acts on the translation category $E_G$ from the right.  
\end{definition}

Denote by\label{not:gpd} $\Gpd$ the symmetric monoidal category with small groupoids\index{groupoid} as objects and functors as morphisms.  Its monoidal product is given by the Cartesian product.  Its monoidal unit is the groupoid $\{*\}$ with one object and no non-identity morphisms.

\begin{example}[Translation Categories of Action Operads]\label{ex:translation-cat-goperad}
Suppose $\G$ is an action operad.  Applying the translation category construction $E_?$ to the groups\index{G-operad@$\G$-operad!translation category} $\{\G(n)\}_{n\geq 0}$ and the one-colored planar operad structure in $\G$, we obtain a one-colored $\G$-operad\label{not:esubg} \[E_{\G} = \Bigl(\bigl\{E_{\G(n)}\bigr\}_{n\geq 0}, \gamma,\operadunit\Bigr)\] in the symmetric monoidal category $\Gpd$ of small groupoids.  In more details:
\begin{enumerate}[label=(\roman*)]
\item The operadic composition \[\nicexy{E_{\G(n)}\times E_{\G(k_1)} \times\cdots\times E_{\G(k_n)} \ar[r]^-{\gamma} & E_{\G(k)}},\] where $k=k_1+\cdots+k_n$, is the functor whose assignment on objects is the operadic composition \[\nicexy{\G(n)\times\G(k_1)\times \cdots\times \G(k_n) \ar[r]^-{\gammag} & \G(k)}\] in the action operad $\G$.  The assignment on morphisms is uniquely determined by the definition of the morphism sets in the translation category $E_{\G(k)}$.
\item The operadic unit is the functor \[\operadunit : \{*\} \to E_{\G(1)}\] determined by $\operadunit(*)=\id_1\in\G(1)$ on the unique object.
\item The one-colored $\G$-sequence structure is determined level-wise by the group multiplication \[\nicexy{\G(n)\times\G(n) \ar[r] & \G(n)}\] on $\G(n)$.
\end{enumerate}
This observation can be found in Fiedorowicz \cite{fiedorowicz} and also in Corollary 1.18 in Gurski \cite{gurski}.\dqed
\end{example}

\begin{example}[Nerves of Translation Categories of Action Operads]\label{ex:translation-g}
Suppose $\Delta$\index{finite ordinal category}\index{category!finite ordinal} denotes the category of finite ordinal numbers \[[n] = \{0 < 1 < \cdots < n\}\] for $n \geq 0$ with weakly order-preserving maps as morphisms.  For a small category $\C$, recall that its \emph{nerve}\index{nerve}
\begin{equation}\label{nerve-category}
\Nerve(\C) : \Delta^{\op} \to \Set,
\end{equation}
also known as its \index{classifying space}\emph{classifying space}, is the simplicial set whose set of $n$-simplicies is the set \[\Nerve(\C)_n = \Nerve(\C)([n]) = \Cat([n],\C)\] of functors from $[n]$ to $\C$.  See, for instance, Example 1.4 in \cite{goerss-jardine}.  The nerve functor \[\Nerve : \Cat \to \Sset,\] from the category of small categories to the category of simplicial sets, is product-preserving and takes each natural transformation to a simplicial homotopy.

Suppose $\G$ is an action operad.  Applying the nerve functor to the one-colored $\G$-operad $E_{\G}$ in groupoids, we obtain a one-colored $\G$-operad\label{not:eofg}\index{translation category!nerve of} \[E\G = \bigl\{E\G(n) = \Nerve(E_{\G(n)})\bigr\}_{n\geq 0}\] in the category of simplicial sets.  This one-colored $\G$-operad has the following two properties.
\begin{enumerate}
\item Each simplicial set $E\G(n)$ is contractible, since it is the nerve of a small connected groupoid.
\item The $\G(n)$-action \[E\G(n) \times \G(n) \to E\G(n)\] is free and proper.  
\end{enumerate}
We call a $\G$-operad with these two properties a\index{Ginfinity-operad@$\G_\infty$-operad}\index{operad!Ginfinity@$\G_\infty$} \emph{$\G_\infty$-operad}.  

When $\G$ is the symmetric group operad $\S$ in Example \ref{ex:symmetric-group-operad} with the identity augmentation, an $\S_\infty$-operad is usually called an\index{Einfinity-operad@$E_\infty$-operad}\index{operad!Einfinity@$E_\infty$} \emph{$E_\infty$-operad}.  In particular, \[E\S = \{ES_n\}_{n\geq 0}\] is a one-colored $E_\infty$-operad in the category of simplicial sets, called the\index{Barratt-Eccles operad}\index{operad!Barratt-Eccles} Barratt-Eccles operad \cite{be}.  For the action operad $\B=\{B_n\}$ defined by the braid groups, which we will discuss in details in Chapter \ref{ch:braided_operad}, the $\B_\infty$-operad terminology and the one-colored braided operad \[E\B=\{EB_n\}\] are due to Fiedorowicz \cite{fiedorowicz}.  At the generality of an action operad, this example is essentially due to Wahl \cite{wahl}.\dqed\end{example}

\chapter{Braided Operads}\label{ch:braided_operad}

The purpose of this chapter is to discuss group operads in which the action operad arises from the braid groups or the pure braid groups.  These group operads are called braided operads and pure braided operads.  They are planar operads with an appropriate action by the (pure) braid groups.  One-colored braided operads were originally introduced by Fiedorowicz \cite{fiedorowicz}.

In Section \ref{sec:braid-groups} we discuss (pure) braid groups.  In Section \ref{sec:braided-operad} we discuss the action operads $\B$ and $\PB$, called the (pure) braid group operad, which are defined by the (pure) braid groups.  The corresponding notions of $\B$-operads and $\PB$-operads are colored braided operads and colored pure braided operads, respectively.  In Section \ref{sec:braided-operad-ex} we provide several examples of braided operads.  In Section \ref{sec:braided-operad-examples} we discuss the one-colored braided operads obtained by taking  the universal covers of the spaces in the little $2$-cube operad and the little $2$-disc operad.

As in the previous chapter, $(\M,\otimes,\tensorunit)$ denotes a complete and cocomplete symmetric monoidal category whose monoidal product commutes with small colimits on each side.

\section{Braid Groups}\label{sec:braid-groups}

The purpose of this section is to recall the (pure) braid groups, which are needed to define colored (pure) braided operads.  For more detailed discussion of the braid groups, the reader is referred to \cite{artina,artinb,kassel-turaev}.  

\begin{definition}\label{def:braid-group}
The \index{braid group}\emph{$n$th braid group} $B_n$ is defined as the trivial group if $n=0,1$.  For $n \geq 2$ the braid group $B_n$ is the group generated by generators $s_1,\ldots,s_{n-1}$ that is subject to the \index{braid relations}\emph{braid relations}: \begin{align}s_is_j &= s_js_i & \text{ for $|i-j|\geq 2$ and $1\leq i,j \leq n-1$};\\ s_is_{i+1}s_i &= s_{i+1}s_is_{i+1} & \text{ for $1 \leq i \leq n-2$}.\label{braid-relations}\end{align}
An element in $B_n$ is called a \emph{braid}.
\end{definition}

\begin{example}\label{ex:B2}
The braid group $B_2$, with one generator $s_1$ and no relations, is the infinite cyclic group.\dqed
\end{example}

\begin{remark}\label{rk:pure-braid-group}
The \index{symmetric group!generators and relations}symmetric group $S_n$ admits a similar generator and relation description as the braid group $B_n$ with $s_i$ given by the adjacent transposition $(i,i+1)$ and with one extra relation: $s_is_i = \id$ for each $i$.\dqed\end{remark}  

\begin{definition}\label{def:pure-braid-group}
There is a natural projection \[\pi : B_n \to S_n\] that sends each generator $s_i \in B_n$ to the adjacent transposition $(i,i+1) \in S_n$.  The kernel of this projection is called the\index{pure braid group} \emph{pure braid group} and is denoted by $PB_n$.  The image of an element $b \in B_n$ under the projection $\pi$ is called the \index{underlying permutation!of a braid}\emph{underlying permutation} of $b$, denoted $\bbar$.
\end{definition}

\begin{ginterpretation}Elements in the braid group $B_n$ can be geometrically interpreted as follows.  A \index{topological interval}\emph{topological interval} is a topological space homeomorphic to the closed unit interval $\cali= [0,1] \subseteq \fieldr$.  A \index{geometric braid}\emph{geometric braid on $n$ strings} is a subset \[b=\coprod_{i=1}^n I_i \subseteq \fieldr^2 \times \cali\] such that the following statements hold:
\begin{enumerate}[label=(\roman*)]
\item Each $I_i$, called the\index{string!in a geometric braid} \emph{$i$th string of $b$}, is a topological interval via the homeomorphism \[I_i \subseteq \fieldr^2 \times \cali \to \cali,\] where the second map is the projection.
\item There are equalities 
\[\begin{split}
I_i \cap \bigl(\fieldr^2 \times \{0\}\bigr) &= \bigl\{(i,0,0)\bigr\} \forspace 1 \leq i \leq n,\\
\bigl(\coprod_{i=1}^n I_i\bigr) \cap \bigl(\fieldr^2 \times \{1\}\bigr) &= \bigl\{(1,0,1), (2,0,1),\ldots,(n,0,1)\bigr\}.
\end{split}\]  
\end{enumerate}
A \emph{braid on $n$ strings} is an isotopy class of geometric braids on $n$ strings.  We will use the same notation for a geometric braid and its isotopy class.  The \emph{underlying permutation $\bbar \in S_n$} of a braid on $n$ strings $b$ is the permutation given by $i \mapsto \bbar(i)$ for $1 \leq i \leq n$, where \[I_i \cap \bigl(\fieldr^2 \times \{1\}\bigr) = \bigl\{(\bbar(i),0,1)\bigr\}.\]

Similar to the fundamental group, using the $\cali$ coordinate, the set of braids on $n$ strings is a group under vertical composition.  This group is naturally identified with the braid group $B_n$, in which the generator $s_i \in B_n$ is identified with the braid on $n$ strings with the $(i+1)$st string crossing over the $i$th string when viewed from bottom to top.  We will identify a braid (i.e., an element in $B_n$) with the corresponding braid on $n$ strings.  Under this identification, the pure braid group consists of braids whose underlying permutation is the identity permutation.\dqed\end{ginterpretation} 

\begin{example}\label{ex1:braid} 
The generator $s_2 \in B_5$ and its inverse are the braids:
\begin{center}\begin{tikzpicture}[scale=.7]
\braid[number of strands=5, thick] (braid) at (0,0) a_2;\node at (2,.7) {$s_2 \in B_5$};
\draw[->] (-1,-1.5) -- (5,-1.5); \draw[->] (-1,-1.5) -- (-.5,-1.2); \draw[->] (-1,-1.5) -- (-1,1);
\foreach \x in {0,...,5} \node at (\x-1,-1.8) {\scriptsize{\x}};
\node at (-1.2,0) {\scriptsize{1}}; \draw[dotted] (-1,0) -- (5,0);
\braid[number of strands=5, thick] (braid-inv) at (8,0) a_2^{-1};\node at (10,.7) {$s_2^{-1} \in B_5$};
\draw[->] (7,-1.5) -- (13,-1.5); \draw[->] (7,-1.5) -- (7.5,-1.2); \draw[->] (7,-1.5) -- (7,1);
\foreach \x in {0,...,5} \node at (\x+7,-1.8) {\scriptsize{\x}};
\node at (6.8,0) {\scriptsize{1}};\draw[dotted] (7,0) -- (13,0);
\end{tikzpicture}\end{center}
Under this identification, the generator $s_i \in B_n$ has underlying permutation the adjacent transposition $(i,i+1)$.\dqed
\end{example}

\begin{example}\label{ex2:braid}
If $n=3$ and $i=1$, then the second braid relation $s_1s_2s_1=s_2s_1s_2$ in \eqref{braid-relations} is the equality
\begin{center}\begin{tikzpicture}[xscale=.7, yscale=.5]
\braid[number of strands=3, thick] (braid1) at (0,0) a_1 a_2 a_1;
\draw[dotted] (-.5,-1.2) -- (2.5,-1.2); 
\draw[dotted] (-.5,-2.2) -- (2.5,-2.2);
\node at (1,0.5) {$s_1s_2s_1 \in B_3$};
\draw[->] (-.5,-3.5) -- (2.5,-3.5); 
\foreach \x in {1,2,3} \node at (\x-1,-3.8) {\scriptsize{\x}};
\node at (4,-2) {\Huge{$=$}};
\braid[number of strands=3, thick] (braid2) at (6,0) a_2 a_1 a_2;
\draw[dotted] (5.5,-1.2) -- (8.5,-1.2); 
\draw[dotted] (5.5,-2.2) -- (8.5,-2.2);
\node at (7,0.5) {$s_2s_1s_2 \in B_3$};
\draw[->] (5.5,-3.5) -- (8.5,-3.5); 
\foreach \x in {1,2,3} \node at (\x+5,-3.8) {\scriptsize{\x}};
\end{tikzpicture}\end{center}
of braids on $3$ strings.  The horizontal dotted lines are there to help visualize the composition of braids.\dqed
\end{example}  

\begin{ginterpretation}
There is another interpretation of the braid groups in terms of configuration spaces.  For a topological space $X$ and an integer $n > 0$, denote by\label{not:fxn} \[F(X;n) = \Bigl\{(x_1,\ldots,x_n) \in X^{\times n} : x_i \not= x_j \ifspace i\not= j\Bigr\}\] the \index{configuration space}\emph{configuration space} of $n$ distinct ordered points in $X$ with the subspace topology.  The symmetric group $S_n$ acts on $F(X;n)$ by permuting the points.  Denote by\label{not:opendtwo} $\opendtwo$ the open unit disc in $\fieldr^2$, i.e., the interior of $\Dtwo$.  Then there is an identification\index{braid group!geometric interpretation} \[B_n = \pi_1\bigl(F(\opendtwo;n)/S_n, *\bigr)\] of the braid group with the fundamental group of the quotient space $F(\opendtwo;n)/S_n$ of the configuration space of $n$ distinct ordered points in the open unit disc.\dqed
\end{ginterpretation}

Recall the block permutation and the direct sum permutation in Definition \ref{def:block-direct-sum-permutation}.  Next we discuss their braid analogues.

\begin{example}[Direct Sum Braid]\label{ex:direct-sum-braid}
The direct sum homomorphism \[\nicexy{S_{k_1} \times \cdots \times S_{k_n} \ar[r] & S_{k_1+\cdots+k_n}}\] of permutations admits a lift to the braid groups, so there is a commutative diagram 
\begin{equation}\label{direct-sum-braid-lift}
\nicexy@C+.7cm{B_{k_1} \times \cdots \times B_{k_n} \ar[r]^-{\text{direct sum}} \ar[d]_-{(\pi,\ldots,\pi)} & B_{k_1+\cdots+k_n} \ar[d]^-{\pi}\\ S_{k_1} \times \cdots \times S_{k_n} \ar[r]^-{\text{direct sum}} & S_{k_1+\cdots+k_n}}
\end{equation}
in the category of groups.

Suppose $b_1 \in B_{k_1}, \ldots, b_n \in B_{k_n}$ are $n > 0$ braids.  Geometrically, the \index{direct sum!braid}\index{braid!direct sum}direct sum braid \[b_1 \oplus \cdots \oplus b_n \in B_{k_1+\cdots+k_n}\] is the braid obtained by placing $b_1,\ldots,b_n$ side-by-side from left to right.  For example, the direct sum braid of $s_1s_1 \in B_2$ and $s_1s_2s_1 \in B_3$,
\begin{center}\begin{tikzpicture}[xscale=.7, yscale=.5]
\braid[number of strands=2, thick] (braid0) at (-3,0) a_1 a_{-1} a_1;
\draw[->] (-3.5,-3.5) -- (-1.5,-3.5); \foreach \x in {1,2} \node at (\x-4,-3.8) {\scriptsize{\x}};
\node at (-2.5,0.5) {$s_1s_1 \in B_2$}; \node at (-1,-2) {\Huge{$\oplus$}};
\braid[number of strands=3, thick] (braid1) at (0,0) a_1 a_2 a_1;
\draw[->] (-.5,-3.5) -- (2.5,-3.5); \foreach \x in {1,2,3} \node at (\x-1,-3.8) {\scriptsize{\x}};
\node at (1,0.5) {$s_1s_2s_1 \in B_3$}; \node at (4,-2) {\Huge{$=$}};
\braid[number of strands=5, thick] (braid2) at (6,0) a_1-a_3 a_{-1}-a_4 a_1-a_3; 
\draw[->] (5.5,-3.5) -- (10.5,-3.5); \foreach \x in {1,...,5} \node at (\x+5,-3.8) {\scriptsize{\x}};
\node at (8,.5) {$s_1s_1 \oplus s_1s_2s_1 \in B_5$};
\end{tikzpicture}\end{center}
is the braid \[s_1s_1 \oplus s_1s_2s_1 = s_1 s_1 s_3 s_4 s_3 \in B_5.\]  

To define the direct sum braid algebraically, since the direct sum map is a group homomorphism, it is enough to define it on generating braids $s_{i_j} \in B_{k_j}$ for $1\leq j \leq n$ with each $1 \leq i_j \leq k_j-1$.  Algebraically, the direct sum braid is given by the product
\[\begin{split}
s_{i_1} \oplus \cdots \oplus s_{i_n} 
&= \prod_{j=1}^n s_{k_1+\cdots+k_{j-1}+i_j}\in B_{k_1+\cdots+k_n}\\
&= s_{i_1} s_{k_1+i_2} \cdots\, s_{k_1+\cdots+k_{n-1}+i_n} 
\end{split}\] 
of $n$ generating braids, where the index of each $s_{i_j}$ is shifted by $k_1+\cdots +k_{j-1}$ on the right-hand side.\dqed
\end{example}

\begin{example}[Block Braid]\label{ex:block-braid}
For $n>0$ and $k_1,\ldots,k_n\geq 0$, the block permutation map \[\nicexy{S_n \ar[r] & S_{k_1+\cdots+k_n}},\qquad S_n \ni \sigma \mapsto \sigma\langle k_1,\ldots,k_n\rangle\] admits a lift to the braid groups, so there is a commutative diagram\index{block braid}\index{braid!block} 
\begin{equation}\label{block-braid-lift}
\nicexy@C+1.5cm{B_n \ar[d]_-{\pi} \ar[r]^-{\text{block braid}} & B_{k_1+\cdots+k_n} \ar[d]^-{\pi}\\ S_n \ar[r]^-{\text{block permutation}} & S_{k_1+\cdots+k_n}}
\end{equation}
in $\Set$.  

Geometrically, for a braid $b \in B_n$, the block braid \[b\langle k_1,\ldots,k_n\rangle \in B_{k_1+\cdots+k_n}\] is obtained from $b$ by replacing its $i$th string by $k_i$ parallel strings for $1 \leq i \leq n$.  For example, for the generator $s_1 \in B_2$, $k_1 = 2$, and $k_2 = 3$, we have the block braid \[s_1\langle 2,3\rangle = s_3s_2s_1s_4s_3s_2 \in B_5,\] as illustrated in the picture below.
\begin{center}\begin{tikzpicture}[xscale=.7, yscale=.4]
\foreach \x in {0,...,4} {\coordinate (A\x) at (\x,-1.5); \coordinate (B\x) at (\x,5);}
\foreach \x in {0,1} {\pgfmathsetmacro\xplus{\x + 3}
\draw[thick] (A\x) to[out=90, in=270] (B\xplus);}
\foreach \x in {2,3,4} {\pgfmathsetmacro\xminus{\x -2}
\draw[line width=7pt, white] (A\x) to [out=90, in=270] (B\xminus);
\draw[thick] (A\x) to [out=90, in=270] (B\xminus);}
\draw[->] (-.5,-1.5) -- (4.5,-1.5); \foreach \x in {1,...,5} \node at (\x-1,-2) {\scriptsize{\x}};
\node at (2,5.7) {$s_1\langle 2,3\rangle \in B_5$}; \node at (6,2) {\Huge{$=$}};
\braid[number of strands=5, thick] (braid) at (8,5) a_3 a_2 a_1 a_4 a_3 a_2;
\draw[->] (7.5,-1.5) -- (12.5,-1.5); \foreach \x in {1,...,5} \node at (\x+7,-2) {\scriptsize{\x}};
\node at (10,5.7) {$s_3s_2s_1s_4s_3s_2 $};
\end{tikzpicture}\end{center}

Algebraically, for a generating braid $s_i \in B_n$ with $1\leq i \leq n-1$, the block braid is given by the product
\[s_i\langle k_1,\ldots,k_n\rangle = \sigma_1\cdots\,\sigma_{k_i} \in B_{k_1+\cdots+k_n}\] of $k_i$ braids, in which 
\[\begin{split}
\sigma_j &= \prod_{m=1}^{k_{i+1}} s_{k_1+\cdots+k_{i-1}+k_{i+1}+j-m} \in B_{k_1+\cdots+k_n}\\
&= s_{k_1+\cdots+k_{i-1}+k_{i+1}+j-1} \cdots\, s_{k_1+\cdots+k_{i-1}+j}
\end{split}\] for each $1\leq j \leq k_i$.  So the block braid $s_i\langle k_1,\ldots,k_n\rangle$ is a product of $k_ik_{i+1}$ generating braids in $B_{k_1+\cdots+k_n}$.  

To understand these formulas, note that the generating braid $s_i \in B_n$ has only one crossing, with the $(i+1)$st string crossing over the $i$th string.  In the induced block braid, each of the $k_{i+1}$ strings in the $(i+1)$st block crosses over each of the $k_i$ strings in the $i$th block.  The braid $\sigma_j$ above encodes the $k_{i+1}$ strings in the $(i+1)$st block crossing over the $j$th string in the $i$th block.  In particular, in the product formula for $\sigma_j$, the right-most generator represents the first string in the $(i+1)$st block crossing over the $j$th string in the $i$th block.  In general, the $l$th generator, counting from the right, in $\sigma_j$ represents the $l$th string in the $(i+1)$st block crossing over the $j$th string in the $i$th block.  

Inductively, for $\sigma,\tau\in B_n$, suppose the block braids $\sigma\langle \cdots\rangle$ and $\tau\langle \cdots\rangle\in B_{k_1+\cdots+k_n}$ have already been defined.  Then the block braid for the product $\sigma\tau$ is defined as the product
\[(\sigma\tau)\langle k_1,\ldots,k_n\rangle
= \sigma\langle k_{\taubar^{-1}(1)},\ldots,k_{\taubar^{-1}(n)} \rangle \tau\langle k_1,\ldots,k_n\rangle\]
of two block braids, in which $\taubar \in S_n$ is the underlying permutation of $\tau$.  An inspection of relevant pictures shows that this algebraic definition of block braid is well-defined; i.e., for each braid relation, the two sides define the same block braid.  Moreover, the algebraic definition agrees with the geometric definition above.  Note that the block permutation map and the block braid map are \emph{not} group homomorphisms.\dqed
\end{example}

\section{Braided Operads as Monoids}\label{sec:braided-operad}

The purpose of this section is to define (pure) braided operads using the (pure) braid groups.  To see that the (pure) braid groups have the required operadic composition, let us first consider the following example.

\begin{example}[Comp-$i$ of Braids]\label{ex:compi-braid}
Combining Example \ref{ex:direct-sum-braid} and Example \ref{ex:block-braid}, for braids $\sigma \in B_n$ and $\tau \in B_m$ and $1 \leq i \leq n$, we define\index{braid group operad!composition}
\begin{equation}\label{braided-compi}
\sigma \compi \tau = \overbrace{\sigma\langle \underbrace{1,\ldots,1}_{i-1}, m,\underbrace{1,\ldots,1}_{n-i}\rangle}^{\text{block braid}} \cdot \overbrace{\bigl(\underbrace{\id_1 \oplus \cdots \oplus \id_1}_{i-1} \oplus \tau \oplus \underbrace{\id_1 \oplus \cdots \oplus \id_1}_{n-i}\bigr)}^{\text{direct sum braid}} \in B_{n+m-1}
\end{equation}
as the product of a direct sum braid induced by $\tau$ with a block braid induced by $\sigma$.  Geometrically, $\sigma \compi \tau$ is the braid obtained from $\sigma$ by replacing its $i$th string with the braid $\tau$.  An inspection of the relevant pictures shows that \[\B = \Bigl(\{B_n\}_{n\geq 0}, \comp, \id_1\in B_1\Bigr)\] is a one-colored planar operad in $\Set$ as in Definition \ref{def:planar-operad-compi}.

For example, suppose $n=i=2$ and $m=3$.  Then we have 
\[\begin{split}
(s_1 \in B_2) \comp_2 (s_1^{-1} \in B_3) 
&= s_1\langle 1,3\rangle \cdot \bigl(\id_1 \oplus s_1^{-1}\bigr)\\
&= s_3 s_2 s_1 s_2^{-1} \in B_4,
\end{split}\] 
as illustrated in the picture below.
\begin{center}\begin{tikzpicture}[scale=.7]
\braid[number of strands=2, thick, style strands={2}{blue}] (braid0) at (0,0) a_1;
\draw[->] (-.5,-1.5) -- (1.5,-1.5); \foreach \x in {1,2} \node at (\x-1,-1.8) {\scriptsize{\x}};
\node at (-2,-.5) {$s_1 \in B_2$}; \node at (.5,-2.7) {$\comp_2$};
\braid[number of strands=3, thick, style strands={1,2,3}{magenta}] (braid1) at (0,-3.5) a_1^{-1};
\draw[->] (-.5,-5) -- (2.5,-5); \foreach \x in {1,2,3} \node at (\x-1,-5.3) {\scriptsize{\x}};
\node at (-2,-4) {$s_1^{-1} \in B_3$}; \node at (4,-2.7) {\Huge{$=$}};
\foreach \x in {1,...,4} {\coordinate (A\x) at (\x+5,-5); 
\coordinate (B\x) at (\x+5,-2.5); \coordinate (C\x) at (\x+5,0);}
\draw[thick, blue] (A1) to [out=90, in=270] (B1); 
\draw[thick, magenta] (A4) to [out=90, in=270] (B4);
\draw[thick, magenta] (A3) to [out=90, in=270] (B2);
\draw[line width=7pt, white] (A2) to [out=90, in=270] (B3);
\draw[thick, magenta] (A2) to [out=90, in=270] (B3);
\draw[thick, blue] (B1) to [out=90, in=270] (C4);
\foreach \x in {2,3,4} {\pgfmathsetmacro\xminus{\x -1}
\draw[line width=7pt, white] (B\x) to [out=90, in=270] (C\xminus);
\draw[thick, magenta] (B\x) to [out=90, in=270] (C\xminus);}
\draw[->] (5.5,-5) -- (9.5,-5); \foreach \x in {1,2,3,4} \node at (\x+5,-5.3) {\scriptsize{\x}};
\draw[dotted] (5.5,-2.5) -- (9.5,-2.5); \draw[dotted] (5.5,0) -- (9.5,0);
\end{tikzpicture}\end{center}

By an inspection of pictures and \eqref{compi-to-gamma}, the operadic composition \[\nicexy{B_n \times B_{k_1}\times \cdots \times B_{k_n} \ar[r]^-{\gamma} & B_{k_1+\cdots+k_n}}\] in the planar operad $\B$ is given by the product
\begin{equation}\label{gamma-braid-group-operad}
\gamma\bigl(\sigma; \tau_1,\ldots,\tau_n\bigr) = 
\overbrace{\sigma\langle k_1,\ldots,k_n\rangle}^{\text{block braid}} \cdot \overbrace{\bigl(\tau_1 \oplus \cdots \oplus \tau_n\bigr)}^{\text{direct sum braid}}
\end{equation}
in $B_{k_1+\cdots+k_n}$.  Geometrically, $\gamma\bigl(\sigma; \tau_1,\ldots,\tau_n\bigr)$ is the braid obtained from $\sigma$ by replacing its $i$th string with the braid $\tau_i$ for $1\leq i \leq n$.\dqed
\end{example}

The next observation about block braid and direct sum braid is also proved by an inspection of pictures.

\begin{lemma}\label{braid-lemma}
For $\sigma,\rho\in B_n$ and $\tau_i\in B_{k_i}$ with $1 \leq i \leq n$, the equalities
\[\begin{split}
(\sigma\rho)\langle k_1,\ldots,k_n\rangle &=\sigma\langle k_{\rhobar^{-1}(1)},\ldots,k_{\rhobar^{-1}(n)}\rangle \cdot \rho\langle k_1,\ldots,k_n\rangle,\\
\sigma\langle k_1,\ldots,k_n\rangle \cdot \bigl(\tau_1 \oplus\cdots\oplus \tau_n\bigr) &= \bigl(\tau_{\sigmabar^{-1}(1)} \oplus\cdots\oplus \tau_{\sigmabar^{-1}(n)}\bigr) \cdot \sigma\langle k_1,\ldots,k_n\rangle
\end{split}\] 
hold in $B_{k_1+\cdots+k_n}$, where $\sigmabar=\pi(\sigma)$ and $\rhobar = \pi(\rho)\in S_n$ are the underlying permutations of $\sigma$ and $\rho$, respectively.
\end{lemma}

\begin{proposition}\label{B-braid-group-operad}
The braid groups form an action operad \[\B = \Bigl(\{B_n\}_{n\geq 0}, \comp, \id_1\in B_1,\pi\Bigr)\] with:
\begin{itemize}
\item the $\compi$-composition as in \eqref{braided-compi};
\item the augmentation $\pi : \B \to \S$ given level-wise by the underlying permutation map $\pi : B_n \to S_n$.
\end{itemize}
\end{proposition}

\begin{proof}
We already noted in Example \ref{ex:compi-braid} that the braid groups form a one-colored planar operad in $\Set$ as in Definition \ref{def:planar-operad-compi}.  

To see that $\pi : \B \to \S$ is a morphism of one-colored planar operads in $\Set$, first note that $\pi$ sends the identity braid $\id_1\in B_1$ to the identity permutation in $S_1$.  That $\pi$ preserves the operadic compositions is a consequence of:
\begin{itemize}
\item the decompositions of the operadic composition $\gamma$ in the symmetric group operad in Proposition \ref{sn-operad}(i) and in $\B$ in \eqref{gamma-braid-group-operad};
\item the fact that level-wise $\pi : B_n \to S_n$ is a group homomorphism;
\item the commutative diagrams \eqref{direct-sum-braid-lift} and \eqref{block-braid-lift}.
\end{itemize}

It remains to prove the action operad axiom \eqref{group-operad-identites}.  So suppose $\sigma,\rho \in B_n$ with $n\geq 1$, and $\tau_i,\tau_i' \in B_{k_i}$ for $1\leq i \leq n$.  Then the action operad axiom for the braid groups follows from the following computation in $B_{k_1+\cdots+k_n}$:
\[\begin{split}
&\gamma\Bigl(\sigma\rho; \{\tau_i\tau_i'\}_{i=1}^n\Bigr)\\
&= (\sigma\rho)\langle k_1,\ldots,k_n\rangle \cdot \Bigl(\bigoplus_{i=1}^n \tau_i\tau_i'\Bigr)\\
&= \sigma\langle k_{\rhobar^{-1}(1)},\ldots,k_{\rhobar^{-1}(n)}\rangle \cdot \rho\langle k_1,\ldots,k_n\rangle \cdot \Bigl(\bigoplus_{i=1}^n \tau_i\Bigr) \cdot \Bigl(\bigoplus_{i=1}^n \tau_i'\Bigr)\\
&= \sigma\langle k_{\rhobar^{-1}(1)},\ldots,k_{\rhobar^{-1}(n)}\rangle \cdot \Bigl(\bigoplus_{i=1}^n \tau_{\rhobar^{-1}(i)}\Bigr) \cdot \rho\langle k_1,\ldots,k_n\rangle \cdot \Bigl(\bigoplus_{i=1}^n \tau_i'\Bigr)\\
&= \gamma\Bigl(\sigma; \{\tau_{\rhobar^{-1}(i)}\}_{i=1}^n\Bigr) \cdot \gamma\Bigl(\rho; \{\tau_i'\}_{i=1}^n\Bigr).
\end{split}\]
In the above computation, the first and the last equalities use \eqref{gamma-braid-group-operad}.  The second and the third equalities follow from Lemma \ref{braid-lemma} and the multiplicativity of direct sum braid.
\end{proof}

\begin{definition}\label{def:braid-group-operad}
The \index{braid group operad}\index{group operad!braid}\emph{braid group operad} $\B$ is the action operad in Proposition \ref{B-braid-group-operad}.
\end{definition}

Restricting the above discussion to the pure braid groups leads to the following action operad.

\begin{definition}\label{def:pure-braid-group-operad}
The \index{pure braid group operad}\index{braid group operad!pure}\emph{pure braid group operad} $\PB$ is the action operad with:
\begin{itemize}
\item $\PB(n) = PB_n$, the $n$th pure braid group in Definition \ref{def:pure-braid-group};
\item the $\compi$-composition as in \eqref{braided-compi};
\item the augmentation $\rho : \PB \to \S$ given level-wise by $\rho(b)=\id_n$ for $b\in PB_n$.
\end{itemize}
\end{definition}

Following the development in Section \ref{sec:group-operad}, for the action operad $\B$ (resp., $\PB$), we will use the adjective \emph{braided} (resp., \emph{pure braid}) instead of $\G$.

\begin{definition}\label{def:braided-operad}
Suppose $\colorc$ is a set.  With the action operad $\G=\B$:
\begin{enumerate}
\item The objects of the diagram category\label{notation:bseqcm} \[\bseqcm = \M^{\bsubcopc}\] are called \index{braided sequence}\emph{$\colorc$-colored braided sequences} in $\M$.  
\item For $X,Y  \in \bseqcm$, $X \circb Y \in \bseqcm$ is called the \index{circle product!braided}\index{braided circle product}\emph{$\colorc$-colored braided circle product}.  It is defined entrywise as the coend\label{not:circb} \[(X \circb Y)\dub = \int^{\uc \in \bsubc} X\duc \otimes Y^{\uc}_{\smallsubb}(\ub)\in \M\] for $\dub \in \bsubcopc$.  The object $Y^{\uc}_{\smallsubb} \in \M^{\bsubcop}$ is defined entrywise as the coend \[Y^{\uc}_{\smallsubb}(\ub) = \int^{\{\ua_j\}\in\prod_{j=1}^m \bsubcop} \bsubcop\bigl(\ua_1,\ldots,\ua_m;\ub\bigr) \cdot \left[\bigotimes_{j=1}^m Y\cjuaj\right] \in \M\] for $\ub \in \bsubcop$, in which $(\ua_1,\ldots,\ua_m) \in \bsubcop$ is the concatenation.  
\item The category of \index{operad!braided}\index{braided operad}\emph{$\colorc$-colored braided operads in $\M$} is the category\label{not:boperadcm} \[\Boperadcm = \Mon\bigl(\bseqcm, \circb, \I^{\B}\bigr)\] of monoids in the monoidal category $\bigl(\bseqcm,\circb,\I^{\B}\bigr)$.  If $\colorc$ has $n<\infty$ elements, we also refer to objects in $\Boperadcm$ as \emph{$n$-colored braided operads in $\M$}.
\item Suppose $\O$ is a $\colorc$-colored braided operad in $\M$.  The category $\algmo$ \index{braided operad!algebra}\index{algebra!of a braided operad}of \emph{$\O$-algebras} is defined as the category of $(\O\circb -)$-algebras for the monad $\O\circb -$ in $\Mtoc$.
\end{enumerate}
\end{definition}

\begin{remark}Note that the monoidal unit $\I^{\B}$ in $\bseqcm$ is equal to the monoidal unit $\I$ in $\pseqcm$ and in $\sseqcm$, since the first braid group $B_1$ is the trivial group.\dqed\end{remark}

The previous definition also works with the pure braid group operad $\PB$, which defines $\colorc$-colored \index{pure braided sequence}\index{braided sequence!pure}pure braided sequences, the $\colorc$-colored \index{pure braided circle product}\index{circle product!pure braided}pure braided circle product, $\colorc$-colored \index{pure braided operad}\index{operad!pure braided}pure braided operads, and algebras over them.  The explicit descriptions of $\G$-operads in Proposition \ref{def:g-operad-generating} and Proposition \ref{prop:g-operad-compi} now specialize to the braid group operad $\B$ and the pure braid group operad $\PB$.  Moreover, the description of algebras over a $\colorc$-colored $\G$-operad in Proposition \ref{prop:g-operad-algebra-defs} specializes to the (pure) braided case.

\section{Examples of Braided Operads}\label{sec:braided-operad-ex}

Some examples of (pure) braided operads are given in this section.

\begin{example}[Braided Operads as Pure Braided Operads]\label{ex:bop-pbop}
Each $\colorc$-colored braided operad is naturally a $\colorc$-colored pure braided operad by restriction of equivariant structure.\dqed
\end{example}

\begin{example}[Symmetric Operads as Braided Operads]\label{ex:sop-bop}
Every\index{braided operad!induced by a symmetric operad} $\colorc$-colored symmetric operad is naturally a $\colorc$-colored braided operad.  This is the $\G=\B$ special case of Example \ref{ex:sop-gop}.\dqed\end{example}

\begin{example}[Colored Endomorphism Braided Operad]\label{ex:end-braided}
For each $\colorc$-colored object $X=\{X_c\}_{c\in\colorc}$ in $\M$, the endomorphism operad\index{endomorphism operad!braided}\index{braided operad!endomorphism} $\End(X)$, which is a $\colorc$-colored symmetric operad, is naturally a $\colorc$-colored braided operad.  This is the $\G=\B$ special case of Example \ref{ex:end-g}.\dqed\end{example}

\begin{example}[Braid Group Operad as a Braided Operad]
The braid group operad $\B=\{B_n\}_{n\geq 0}$ yields a one-colored braided operad in $\Set$.  This is Example \ref{ex:agop-goperad} when $\G=\B$.  Similarly, the pure braid group operad $\PB = \{PB_n\}_{n\geq 0}$ yields a one-colored pure braided operad in $\Set$.\dqed
\end{example}

\begin{example}[Translation Category of Braid Group Operad]
Applying the \index{braid group operad!translation category}\index{translation category!braid group operad}translation category construction $E_?$ to the braid group operad $\B$, we obtain a one-colored braided operad \[E_{\B} = \Bigl(\bigl\{E_{B_n}\bigr\}_{n\geq 0}, \gamma,\operadunit\Bigr)\] in the symmetric monoidal category of small groupoids.  This is Example \ref{ex:translation-cat-goperad} when $\G=\B$.  Similarly, the pure braid group operad $\PB$ yields a one-colored pure braided operad \[E_{\PB} = \Bigl(\bigl\{E_{PB_n}\bigr\}_{n\geq 0}, \gamma,\operadunit\Bigr)\] in the symmetric monoidal category of small groupoids.\dqed
\end{example}

\begin{example}[Nerves of Translation Categories of Braid Groups]\label{ex:translation-braid}
Applying the nerve functor to the one-colored braided operad $E_{\B}$ in small groupoids, we obtain a one-colored braided operad\index{translation category!nerve of} \[E\B = \{EB_n\}_{n\geq 0}\] in the category of simplicial sets, where \[EB_n = \Nerve(E_{B_n})\] is the nerve of the translation category $E_{B_n}$ of the $n$th braid group $B_n$.  This is the $\G=\B$ special case of Example \ref{ex:translation-g}.  Furthermore, the one-colored braided operad $E\B$ is a\index{Binfinity-operad@$\B_\infty$-operad}\index{operad!Binfinity@$\B_\infty$} \emph{$\B_\infty$-operad} \cite{fiedorowicz}.  That is, each $EB_n$ is contractible, and the braid group action \[EB_n \times B_n \to EB_n\] is free and proper.  

Similarly, with the pure braid group operad $\PB$ in place of the braid group operad $\B$, there is a one-colored pure braided operad \[E\PB = \{EPB_n\}_{n\geq 0}\] in the category of simplicial sets, where \[EPB_n = \Nerve(E_{PB_n})\] is the nerve of the translation category $E_{PB_n}$ of the $n$th pure braid group $PB_n$.  Moreover, each $EPB_n$ is contractible, and the pure braid group action \[EPB_n \times PB_n \to EPB_n\] is free and proper.\dqed
\end{example}

\begin{example}[Parenthesized Braided Operad]\label{ex:parenthesized-braid}
The \index{parenthesized braided operad}\index{braided operad!parenthesized}\emph{parenthesized braided operad} $\PaB$ is the one-colored braided operad in $\Set$ as in Example \ref{ex:parenthesized-goperad} with $\G=\B$.  Its entries are \[\PaB(n) = \Mag(n) \times B_n \times \Mag(n)\] for $n\geq 0$.  For example, the element 
\[\Bigl(((\ast\ast)\ast)\ast; s_3s_2^{-1}s_1s_2^{-1}; (\ast\ast)(\ast\ast)\Bigr)\in \PaB(4)\]
may be visualized as follows.
\begin{center}\begin{tikzpicture}[yscale=.7]
\foreach \x in {1,...,4} {\coordinate (A\x) at (\x,0); 
\coordinate (B\x) at (\x,1.5); \coordinate (C\x) at (\x,5);}
\coordinate (D) at (2.6,3.15); \coordinate (E) at (2.4,3.35);
\foreach \x in {1,4} \draw[thick] (A\x) to [out=90, in=270] (B\x);
\draw[thick] (A3) to [out=90, in=270] (B2);
\draw[line width=7pt, white] (A2) to [out=90, in=270] (B3);
\draw[thick] (A2) to [out=90, in=270] (B3);
\draw[thick] (B1) to [out=90, in=270] (C4);
\foreach \x in {2,4} {\pgfmathsetmacro\xminus{\x -1}
\draw[line width=7pt, white] (B\x) to [out=90, in=270] (C\xminus);
\draw[thick] (B\x) to [out=90, in=270] (C\xminus);}
\draw[thick] (B3) to [out=90, in=315] (D); \draw[thick] (E) to [out=135, in=270] (C2);
\draw[->] (.5,0) -- (4.5,0);
\node at (1,-.3) {(($\ast$}; \node at (2,-.3) {$\ast$)}; \node at (3,-.3) {$\ast$)}; \node at (4,-.3) {$\ast$};
\node at (1,5.3) {($\ast$}; \node at (2,5.3) {$\ast$)}; \node at (3,5.3) {($\ast$}; \node at (4,5.3) {$\ast$)};
\end{tikzpicture}\end{center}

For $1 \leq i \leq n$ and $m \geq 0$, the $\compi$-composition \[\nicexy{\PaB(n) \times \PaB(m) \ar[r]^-{\compi} & \PaB(n+m-1)}\] is defined as \[(w;\sigma;w') \compi (x;\tau;x') = \bigl(w\compi x; \sigma\compi\tau; w' \comp_{\sigmabar(i)} x'\bigr)\] for $w,w'\in \Mag(n)$, $x,x' \in \Mag(m)$, $\sigma\in B_n$, and $\tau\in B_m$.  Here \[w\compi x \andspace w' \comp_{\sigmabar(i)} x' \in \Mag(n+m-1)\] are the $\compi$-composition in the magma operad $\Mag$ with $\sigmabar \in S_n$ the underlying permutation of $\sigma$.  Likewise, \[\sigma\compi\tau\in B_{n+m-1}\] is the $\compi$-composition of braids in \eqref{braided-compi}.  

The $\compi$-composition in $\PaB$ can be explained geometrically as follows.  To form the $\compi$-composition above, we replace the $i$th string in $(w;\sigma;w')$, counting from the left at the bottom, with the element $(x;\tau;x')$.  In doing so, the $i$th $\ast$ in the parenthesized word $w$ is replaced by $x$, which is the $\compi$-composition $w\compi x$.  In terms of the braid, the $i$th string in $\sigma$ is replaced by the braid $\tau$, forming the braid $\sigma \compi \tau$.  As one travels along the $i$th string in $\sigma$ from bottom to top, the $i$th $\ast$ in $w$ reaches the $\sigmabar(i)$th $\ast$ in the parenthesized word $w'$.  So the output parenthesized word of the $\compi$-composition is $w'\comp_{\sigmabar(i)} x'$.

Similarly, using the pure braid group operad $\G=\PB$ in Example \ref{ex:parenthesized-goperad}, we obtain the \index{parenthesized pure braided operad}\index{pure braided operad!parenthesized}\emph{parenthesized pure braided operad} $\PaPB$, which is a one-colored pure braided operad in $\Set$.

The elements in the parenthesized braided operad $\PaB$ are what Bar-Natan \cite{barnatan} called \index{parenthesized braids}\emph{parenthesized braids}, but the notation \textbf{PaB} there means something different from our $\PaB$.  Our parenthesized braided operad is also different from Fresse's (Section 6.2 in \cite{fresse-gt}) parenthesized braid operad \textit{\textsf{PaB}}, which is a one-colored symmetric operad in $\Cat$ whose algebras are small braided monoidal categories without monoidal units.  See Remark \ref{rk:fresse-pab} for more discussion of Fresse's parenthesized braid operad.\dqed
\end{example}

\section{Universal Cover of the Little $2$-Cube Operad}
\label{sec:braided-operad-examples} 

The purpose of this section is to discuss an important example of a one-colored braided operad that comes from universal covers of the spaces in the little $2$-cube operad.  We also discuss the variation where the little $2$-disc operad is used instead.

\begin{example}[Universal Cover of the Little $2$-Cube Operad]\label{ex:universal-cover-c2}
This example\index{little $2$-cube operad!universal cover}\index{operad!universal cover of little $2$-cube} provides another example of a $\B_\infty$-operad and is also due to Fiedorowicz \cite{fiedorowicz}.  Recall from Examples \ref{ex:little-2cube} and \ref{ex:little-ncube} the little $n$-cube operad $\C_n$ in $\CHau$.  We consider universal covers of the spaces $\C_2(n)$ for $n \geq 0$.  Although there are no consistent base points for the spaces $\{\C_2(n)\}_{n\geq 0}$ in the little $2$-cube operad $\C_2$, we can consider the $n$-tuple of little $2$-cubes $\uc^n \in \C_2(n)$ with pairwise disjoint interiors given by the picture:
\begin{center}\begin{tikzpicture}[scale=2]
\draw [thick] (0,0) rectangle (1,1);
\foreach \x in {0.2,0.4,0.8} \draw (\x,0) -- (\x,1);
\node at (0.1,.5) {1}; \node at (0.3,.5) {2}; 
\node at (.6,.5) {$\cdots$}; \node at (0.9,.5) {$n$};
\node at (-.7,.5) {$\uc^n \in \C_2(n)$};
\end{tikzpicture}\end{center}
In $\uc^n$ the $n$ little $2$-cubes partition the closed unit square into $n$ consecutive rectangles with height $1$ and equal width ordered from left to right.  For each $n \geq 0$, denote by\label{not:ctilde} $\Ctilde_2(n)$ the universal cover of $\C_2(n)$ at the base point $\uc^n$.  In other words, $\Ctilde_2(n)$ is the space of homotopy classes of paths \[\alpha : [0,1] \to \C_2(n)\] starting at $\alpha(0)=\uc^n$.  Denote by \[p : \Ctilde_2(n) \to \C_2(n)\] the universal covering map, which sends $\alpha \in \Ctilde_2(n)$ to its end point $\alpha(1) \in \C_2(n)$.

The spaces \[\Ctilde_2 = \bigl\{\Ctilde_2(n)\bigr\}_{n\geq 0}\] formed a one-colored braided operad in $\CHau$ with the following structure.  The operadic unit $\operadunittilde$ is the homotopy class of the constant path \[\operadunittilde : [0,1] \to \C_2(1)\] at the base point, so \[\operadunittilde(t) = \uc^1 \in \C_2(1)\] for $0 \leq t \leq 1$.  Note that \[p(\operadunittilde) = \operadunittilde(1) = \uc^1 \in \C_2(1)\] is the operadic unit of the little $2$-cube operad $\C_2$.

For $n \geq 1$, $k_1,\ldots,k_n \geq 0$, and $k = k_1+\cdots+k_n$, the operadic composition $\gammatilde$ in $\Ctilde_2$ is the lifting in the diagram \[\nicexy{\Ctilde_2(n) \times \Ctilde_2(k_1) \times \cdots \times \Ctilde_2(k_n) \ar[d]_-{(p,\ldots,p)} \ar@{-->}[r]^-{\gammatilde} & \Ctilde_2(k) \ar[d]^p\\
\C_2(n) \times \C_2(k_1) \times \cdots \times \C_2(k_n) \ar[r]^-{\gamma} & \C_2(k)}\] in $\CHau$ defined as follows.  Suppose given $\alpha \in \Ctilde_2(n)$ and $\alpha_i \in \Ctilde_2(k_i)$ for $1 \leq i \leq n$.  First we define a path \[\gamma^0 : [0,1] \to \C_2(k)\] such that: 
\begin{itemize}\item $\gamma^0(0)=\uc^k \in \C_2(k)$.
\item $\gamma^0(1)=\gamma\bigl(\alpha(0); \alpha_1(0),\ldots,\alpha_n(0)\bigr) = \gamma\bigl(\uc^n; \uc^{k_1},\ldots,\uc^{k_n}\bigr) \in \C_2(k)$.
\item For each $0 < t < 1$, $\gamma^0(t) \in \C_2(k)$ consists of $k$ little $2$-cubes that partition the closed unit square into $k$ consecutive rectangles with height $1$, but not necessarily equal width, ordered from left to right.
\end{itemize}
For example, if $k=5$, then a typical $\gamma^0(t) \in \C_2(5)$ may look like
\begin{center}\begin{tikzpicture}[scale=2]
\draw [thick] (0,0) rectangle (1,1);
\foreach \x in {0.1,0.2,0.4,0.7} \draw (\x,0) -- (\x,1);
\node at (.05,.5) {\scriptsize{1}}; \node at (.15,.5) {\scriptsize{2}};
\node at (.3,.5) {3}; \node at (.55,.5) {4}; \node at (.85,.5) {5};
\end{tikzpicture}\end{center}
with $5$ little $2$-cubes ordered from left to right.  Then we define a path \[\alpha(\alpha_1,\ldots,\alpha_n) : [0,1] \to \C_2(k)\] by \[\alpha(\alpha_1,\ldots,\alpha_n)(t) = \gamma\bigl(\alpha(t); \alpha_1(t),\ldots,\alpha_n(t)\bigr) \in \C_2(k) \forspace t\in [0,1].\]  Finally, we define the operadic composition as the homotopy class of the composition of paths \[\gammatilde\bigl(\alpha;\alpha_1,\ldots,\alpha_n\bigr) = \alpha(\alpha_1,\ldots,\alpha_n) \cdot \gamma^0 \in \Ctilde_2(k).\] 

The braid group action on $\Ctilde_2$ is defined as follows.  It is enough to define the action by the generator $s_i \in B_n$ for $1 \leq i \leq n-1$.  Denote by $\tau_i = (i,i+1) \in S_n$ the adjacent transposition that interchanges $i$ and $i+1$.   Suppose $\uc^n\tau_i$ is the same as $\uc^n$ but with the $i$th and the $(i+1)$st little $2$-cubes interchanged.  First we define the path \[\stilde_i : [0,1] \to \C_2(n)\] with \[\stilde_i(0) = \uc^n \andspace \stilde_i(1) = \uc^n\tau_i\] that is determined by the following sequence of pictures.
\begin{center}
\begin{tikzpicture}[scale=2.5]
\draw [thick] (0,0) rectangle (1,1); \foreach \x in {.3,.7} \draw [thick] (\x,0) -- (\x,1);
\draw [very thick, fill=yellow!20] (.3,0) rectangle (.5,1); 
\draw[->, shorten <=1pt, lightgray, line width=2pt] (.4,0) -- (.4,.25);
\draw [very thick, fill=green!20] (.5,0) rectangle (.7,1); 
\draw[->, shorten <=1pt, lightgray, line width=2pt] (.6,1) -- (.6,.75);
\node at (0.4,.5) {\scriptsize{$i$}}; \node at (0.6,.5) {\scriptsize{$i\!+\!1$}}; 
\node at (.15,.5) {$\cdots$}; \node at (.85,.5) {$\cdots$}; 
\node at (.5,1.2) {$\stilde_i(0) = \uc^n$};
\end{tikzpicture}\qquad
\begin{tikzpicture}[scale=2.5]
\draw [thick] (0,0) rectangle (1,1); \foreach \x in {.3,.7} \draw [thick] (\x,0) -- (\x,1);
\draw [very thick, fill=yellow!20] (.3,.5) rectangle (.5,1); 
\draw[->, shorten <=1pt, lightgray, line width=2pt] (.5,.75) -- (.65,.75);
\draw [very thick, fill=green!20] (.5,0) rectangle (.7,.5);
\draw[->, shorten <=1pt, lightgray, line width=2pt] (.5,.25) -- (.35,.25);
\node at (0.6,.25) {\scriptsize{$i\!+\!1$}}; \node at (0.4,.75) {\scriptsize{$i$}}; 
\node at (.15,.5) {$\cdots$}; \node at (.85,.5) {$\cdots$}; 
\node at (.5,1.2) {$\stilde_i(1/3)$};
\end{tikzpicture}\qquad
\begin{tikzpicture}[scale=2.5]
\draw [thick] (0,0) rectangle (1,1); \foreach \x in {.3,.7} \draw [thick] (\x,0) -- (\x,1);
\draw [very thick, fill=green!20] (.3,0) rectangle (.5,.5); 
\draw[->, shorten <=1pt, lightgray, line width=2pt] (.4,.5) -- (.4,.75);
\draw [very thick, fill=yellow!20] (.5,.5) rectangle (.7,1);
\draw[->, shorten <=1pt, lightgray, line width=2pt] (.6,.5) -- (.6,.25);
\node at (0.4,.25) {\scriptsize{$i\!+\!1$}}; \node at (0.6,.75) {\scriptsize{$i$}}; 
\node at (.15,.5) {$\cdots$}; \node at (.85,.5) {$\cdots$}; 
\node at (.5,1.2) {$\stilde_i(2/3)$};
\end{tikzpicture}\qquad
\begin{tikzpicture}[scale=2.5]
\draw [thick] (0,0) rectangle (1,1); \foreach \x in {.3,.7} \draw [thick] (\x,0) -- (\x,1);
\draw [very thick, fill=green!20] (.3,0) rectangle (.5,1); 
\draw [very thick, fill=yellow!20] (.5,0) rectangle (.7,1);
\node at (0.6,.5) {\scriptsize{$i$}}; \node at (0.4,.5) {\scriptsize{$i\!+\!1$}}; 
\node at (.15,.5) {$\cdots$}; \node at (.85,.5) {$\cdots$}; 
\node at (.5,1.2) {$\stilde_i(1) = \uc^n\tau_i$};
\end{tikzpicture}
\end{center}
For $0 \leq t \leq 1$ the little $2$-cubes in $\stilde(t)$ labeled by $1,\ldots,i-1,i+2,\ldots,n$ remain unchanged.  Now for an element $\alpha \in \Ctilde_2(n)$, we first define the path \[\alpha\tau_i : [0,1] \to \C_2(n)\] given by \[(\alpha\tau_i)(t) = \alpha(t)\tau_i \in \C_2(n) \forspace t \in [0,1].\]  Then we define the $s_i$-action $\alpha s_i \in \Ctilde_2(n)$ on $\alpha$ as the homotopy class of the composition of paths \[(\alpha\tau_i) \cdot \stilde_i : [0,1] \to \C_2(n).\] The one-colored braided operad axioms of $\Ctilde_2$ can be checked by inspecting pictures.  

The level-wise universal covering map \[p : \Ctilde_2 \to \C_2\] is a morphism of one-colored planar operads that respects the equivariant structure in the sense that \[p(\alpha b) = (\alpha b)(1) = \alpha(1)\bbar = p(\alpha)\bbar \] for $\alpha \in \Ctilde_2(n)$ and $b \in B_n$ with underlying permutation $\bbar \in S_n$.  In other words, $p : \Ctilde_2 \to \C_2$ is a morphism of one-colored braided operads when the one-colored symmetric operad $\C_2$ is regarded as a one-colored braided operad as in Example \ref{ex:sop-bop}.  Moreover, $\Ctilde_2$ is a $\B_\infty$-operad.  That is, each space $\Ctilde_2(n)$ is contractible, and the braid group action \[\Ctilde_2(n) \times B_n \to \Ctilde_2(n)\] is free and proper.\dqed
\end{example}

\begin{example}[Universal Cover of the Little $2$-Disc Operad]\label{ex:universal-cover-d2}
This example\index{little $2$-disc operad!universal cover}\index{operad!universal cover of little $2$-disc} provides yet another example of a $\B_\infty$-operad.  We proceed as in Example \ref{ex:universal-cover-c2} using the little $2$-disc operad $\D_2$ in Example \ref{ex:little-n-disc} instead of the little $2$-cube operad $\C_2$.  The result is a one-colored braided operad\label{not:dtilde} $\Dtilde_2$ in $\CHau$ that is also a $\B_\infty$-operad.  We will simply point out the modifications.  

Let $\Dtilde_2(n)$ be the universal cover of the space $\D_2(n)$ at the base point $\ud^n$ below.
\begin{center}\begin{tikzpicture}[scale=1.2,thick]
\draw (0cm,0cm) circle(1cm);
\foreach \x in {-.8,-.4,.8} \draw (\x cm,0cm) circle (.2cm);
\node at (-.8cm,0cm) {1}; \node at (-.4cm,0cm) {2}; \node at (.8cm,0cm) {$n$};
\node at (.2cm,0cm) {$\cdots$}; \node at (-2.3cm,0cm) {$\ud^n \in \D_2(n)$};
\end{tikzpicture}\end{center}
In $\ud^n$ the $n$ little $2$-discs are ordered from left to right and are centered at the horizontal axis.  Their horizontal diameters partition $[-1,1]$ into $n$ sub-intervals of equal length.  The operadic unit in $\Dtilde_2(1)$ is the constant path at the base point $\ud^1$.

To define the operadic composition $\gammatilde$ on $\Dtilde_2$, suppose $\alpha \in \Dtilde_2(n)$ with $n \geq 1$, $\alpha_i \in \Dtilde_2(k_i)$ with $1 \leq i \leq n$, and $k = k_1 + \cdots + k_n$.  We use a path \[\gamma^0 : [0,1] \to \D_2(k)\] such that:
\begin{itemize}
\item $\gamma^0(0) = \ud^k \in \D_2(k)$.
\item $\gamma^0(1) = \gamma\bigl(\alpha(0); \alpha_1(0), \ldots, \alpha_n(0)\bigr) = \gamma\bigl(\ud^n; \ud^{k_1},\ldots,\ud^{k_n}\bigr) \in \D_2(k)$.
\item For each $0 < t < 1$, $\gamma^0(t) \in \D_2(k)$ consists of $k$ little $2$-discs that are ordered from left to right and are centered at the horizontal axis.  Their diameters partition $[-1,1]$ into $k$ sub-intervals, not necessarily of equal length.  
\end{itemize}
For example, if $k=5$, then a typical $\gamma^0(t) \in \D_2(5)$ may look like
\begin{center}\begin{tikzpicture}[scale=1.5,thick]
\draw (0,0) circle(1);
\draw (-.8,0) circle (.2); \draw (-.3,0) circle (.3); \draw (.1,0) circle (.1);
\draw (.55,0) circle (.35); \draw (.95,0) circle (.05);
\end{tikzpicture}\end{center}
with $5$ little $2$-discs ordered from left to right.  Then we define the operadic composition as the homotopy class of the composition of paths \[\gammatilde\bigl(\alpha;\alpha_1,\ldots,\alpha_n\bigr) = \alpha(\alpha_1,\ldots,\alpha_n) \cdot \gamma^0 \in \Dtilde_2(k).\] 

To define the braid group action \[\Dtilde_2(n) \times B_n \to \Dtilde_2(n),\] it is again enough to define the action by the generator $s_i \in B_n$ for $1 \leq i \leq n-1$.  Here we use the path \[\stilde_i : [0,1] \to \D_2(n)\] with \[\stilde_i(0) = \ud^n \andspace \stilde_i(1) = \ud^n\tau_i\] that is determined by the following sequence of pictures.
\begin{center}
\begin{tikzpicture}[scale=1.2,thick]
\draw (0,0) circle (1); 
\draw[fill=yellow!20] (-.2,0) circle (.2); \draw[fill=green!20] (.2,0) circle (.2);
\draw[->, gray!50, line width=1.5 pt] (-.2,.21) to[out=90,in=135, looseness=2] (.18,.25);
\draw[->, gray!50, line width=1.5 pt] (.2,-.21) to[out=270,in=315, looseness=2] (-.18,-.25);
\node at (-.2,0) {\scriptsize{$i$}}; \node at (.2,0) {\scriptsize{$i\!\!+\!\!1$}};
\foreach \x in {-.7,.7} \node at (\x,0) {$\cdots$}; \node at (0,1.3) {$\stilde_i(0)$};
\end{tikzpicture}\qquad
\begin{tikzpicture}[scale=1.2,thick]
\draw (0,0) circle (1);
\draw[fill=green!20] (.141,-.141) circle (.2); \draw[fill=yellow!20] (-.141,.141) circle (.2);
\node at (-.141,.141) {\scriptsize{$i$}}; \node at (.141,-.141) {\scriptsize{$i\!\!+\!\!1$}};
\foreach \x in {-.7,.7} \node at (\x,0) {$\cdots$}; \node at (0,1.3) {$\stilde_i(1/3)$};
\end{tikzpicture}\qquad
\begin{tikzpicture}[scale=1.2,thick]
\draw (0,0) circle (1);
\draw[fill=yellow!20] (.141,.141) circle (.2); \draw[fill=green!20] (-.141,-.141) circle (.2);
\node at (.141,.141) {\scriptsize{$i$}}; \node at (-.141,-.141) {\scriptsize{$i\!\!+\!\!1$}};
\foreach \x in {-.7,.7} \node at (\x,0) {$\cdots$}; \node at (0,1.3) {$\stilde_i(2/3)$};
\end{tikzpicture}\qquad
\begin{tikzpicture}[scale=1.2,thick]
\draw (0,0) circle (1); 
\draw[fill=green!20] (-.2,0) circle (.2); \draw[fill=yellow!20] (.2,0) circle (.2);
\node at (-.2,0) {\scriptsize{$i\!\!+\!\!1$}}; \node at (.2,0) {\scriptsize{$i$}};
\foreach \x in {-.7,.7} \node at (\x,0) {$\cdots$}; \node at (0,1.3) {$\stilde_i(1)$};
\end{tikzpicture}\end{center}
For $0 \leq t \leq 1$ the little $2$-discs in $\stilde_i(t)$ labeled by $1,\ldots,i-1,i+2,\ldots,n$ remain unchanged.  Then we define the $s_i$-action $\alpha s_i \in \Dtilde_2(n)$ on $\alpha \in \Dtilde_2(n)$ as the homotopy class of the composition of paths \[(\alpha\tau_i) \cdot \stilde_i : [0,1] \to \D_2(n).\] 

When the little $2$-disc operad $\D_2$ is regarded as a one-colored braided operad, the level-wise universal covering map \[p : \Dtilde_2 \to \D_2\] is a morphism of one-colored braided operads.  Moreover, $\Dtilde_2$ is a $\B_\infty$-operad.  That is, each space $\Dtilde_2(n)$ is contractible, and the braid group action \[\Dtilde_2(n) \times B_n \to \Dtilde_2(n)\] is free and proper.\dqed
\end{example}

\chapter{Ribbon Operads}\label{ch:ribbon_operad}

The purpose of this chapter is to discuss group operads in which the action operad arises from the ribbon groups or the pure ribbon groups.  These group operads are called ribbon operads and pure ribbon operads.  They are planar operads with an appropriate action by the (pure) ribbon groups.  One-colored ribbon operads were originally introduced by Wahl \cite{wahl} and Salvatore-Wahl \cite{salvatore-wahl}.

In Section \ref{sec:ribbon-groups} we discuss the (pure) ribbon groups.  In Section \ref{sec:ribbon-operad} we discuss the action operads $\R$ and $\PR$, called the (pure) ribbon group operad, which are defined by the (pure) ribbon groups.  The corresponding notions of $\R$-operads and $\PR$-operads are colored ribbon operads and colored pure ribbon operads, respectively.  In Section \ref{sec:ribbon-ex} we discuss some examples of ribbon operads.  In Section \ref{sec:ribbon-operad-examples} we discuss the one-colored ribbon operad obtained by taking the universal covers of the spaces in the framed little $2$-disc operad.

As in the previous chapter, $(\M,\otimes,\tensorunit)$ denotes a complete and cocomplete symmetric monoidal category whose monoidal product commutes with small colimits on each side.

\section{Ribbon Groups}\label{sec:ribbon-groups}

The purpose of this section is to recall the (pure) ribbon groups, which are needed to define colored (pure) ribbon operads.  The additive group of integers is denoted by $\bbZ$.

\begin{definition}\label{def:ribbon-group}
For $n\geq 0$ the\index{ribbon group} \emph{$n$th ribbon group} $R_n$ is the semi-direct product \[R_n=B_n \rtimes \bbZ^n,\] where $B_n$ is the $n$th braid group.  An element in $R_n$ is called a \index{ribbon}\emph{ribbon}.
\begin{enumerate}
\item The image of a ribbon $r \in R_n$ under the projection map \[\pi : R_n \to B_n,\] whose kernel is the subgroup $\bbZ^n$,   is called the \index{underlying braid}\emph{underlying braid} of $r$.
\item Composing with the underlying permutation map $B_n \to S_n$, the projection map \[\pi : R_n \to S_n\] is also denoted by $\pi$.  The image of a ribbon $r \in R_n$ under this composite is called the\index{underlying permutation!of a ribbon} \emph{underlying permutation} of $r$, denoted $\rbar$.  
\item The kernel of the projection map $\pi : R_n \to S_n$ is called the \index{pure ribbon group}\emph{pure ribbon group} and is denoted by\label{not:pr} $PR_n$. 
\end{enumerate}
\end{definition}

\begin{ginterpretation}\label{gint:braid-on-ribbon}
There is a geometric interpretation of the ribbon group $R_n$ similar to the interpretation of the braid group $B_n$ as the group of braids on $n$ strings.  Instead of strings for $B_n$, each element in $R_n$ is interpreted as a braid on $n$ strips, each of which can have any finite number of full $2\pi$ twists in either direction.  The number of twists in the $i$th strip corresponds to the $i$th factor of the infinite cyclic group $\bbZ$ in $R_n = B_n \rtimes \bbZ^n$.  Each element in $R_n$ admits a unique expression as \[(b;k_1,\ldots,k_n)\] with 
\begin{itemize}
\item $b \in B_n$;
\item each $k_i$ representing the number of full $2\pi$ twists in the $i$th strip, counting from left to right at the bottom.
\end{itemize}
The projection $\pi : R_n \to B_n$ is given by \[\pi(b;k_1,\ldots,k_n) = (b;0,\ldots,0).\]  In other words, it forgets about the full $2\pi$ twists and only remembers the underlying braid.\dqed
\end{ginterpretation}

\begin{interpretation}\label{int:ribbon-group-generators}
In terms of\index{ribbon group!generators and relations} generators and relations, the ribbon group $R_n$ is the group generated by $r_1,\ldots,r_n$ such that $r_1,\ldots,r_{n-1}$ satisfy the braid relations in Definition \ref{def:braid-group} and that the additional relation 
\begin{equation}\label{ribbon-new-relation}
r_{n-1}r_nr_{n-1}r_n = r_nr_{n-1}r_nr_{n-1}
\end{equation}
is satisfied.  See, for example, \cite{joyal-street} Example 6.4 or \cite{street} Example 11.5.\dqed
\end{interpretation}

\begin{example}\label{ex:generator-ri}
The generator $r_i \in R_n$ for $1 \leq i \leq n-1$ is the braid on $n$ strips
\begin{center}\begin{tikzpicture}[xscale=1.2, yscale=.7, thick]
\foreach \x in {1,3,6,8} \draw (\x,0) rectangle (\x+.2,1);
\node at (2,.5) {$\cdots$}; \node at (7,.5) {$\cdots$};
\draw (4,0) -- (4.2,0) -- (5.2,1) -- (5,1) -- cycle;
\draw[shorten >=.2, shorten <=.2, line width=10pt, white] (5.1,0) -- (4.1,1);
\draw (5,0) -- (5.2,0) -- (4.2,1) -- (4,1) -- cycle;
\node at (1.1,-.4) {1}; \node at (3.1,-.4) {$i-1$}; \node at (4.1,-.4) {$i$};
\node at (5.1,-.4) {$i+1$}; \node at (6.1,-.4) {$i+2$}; \node at (8.1,-.4) {$n$};
\node at (-.2,.5) {$(r_i \in R_n) =$};
\end{tikzpicture}\end{center}
with the $(i+1)$st strip crossing over the $i$th strip when viewed from bottom to top. In terms of the semi-direct product $R_n = B_n \rtimes \bbZ^n$, we have \[(r_i \in R_n) = \bigl(s_i \in B_n; \underbrace{0,\ldots,0}_{n}\bigr),\] where the $n$ zeros indicate that none of the $n$ strips has any twists.\dqed
\end{example}

\begin{example}\label{ex:generator-rn}
The generator $r_n  \in R_n$ is the braid on $n$ strips
\begin{center}\begin{tikzpicture}[xscale=1.2, yscale=.7, thick]
\foreach \x in {1,3} \draw (\x,0) rectangle (\x+.2,1); \node at (2,.5) {$\cdots$}; 
\draw (4,0) to[out=90,in=270] (4.2,.5);
\draw[shorten >=1pt, shorten <=1pt, line width=5pt, white] (4.2,0) to[out=90,in=270] (4,.5);
\draw (4.2,0) to[out=90,in=270] (4,.5);
\draw (4,.5) to[out=90,in=270] (4.2,1);
\draw[shorten >=1pt, shorten <=1pt, line width=5pt, white] (4.2,.5) to[out=90,in=270] (4,1);
\draw (4.2,.5) to[out=90,in=270] (4,1);
\foreach \x in {0,1} \draw (4,\x) -- (4.2,\x);
\node at (1.1,-.4) {1}; \node at (3.1,-.4) {$n-1$}; \node at (4.1,-.4) {$n$};
\node at (-1.7,.2) {$(r_n \in R_n) = \bigl(\id_n \in B_n; \underbrace{0,\ldots,0}_{n-1},1\bigr)=$};
\end{tikzpicture}\end{center}
in which the $n$th strip has a full $2\pi$ twist.\dqed
\end{example}

\begin{example}\label{ex:ribbon-relation}
When $n=2$, the relation \[r_{1}r_2r_{1}r_2 = r_2r_{1}r_2r_{1}\] in \eqref{ribbon-new-relation} corresponds to the following equality of braids on $2$ strips.
\begin{center}
\begin{tikzpicture}[xscale=1.2, yscale=.6, thick]
\foreach \x in {1,1.2}{\draw (\x,0) -- (\x,1);} 
\foreach \y in {0,.5}{\draw (2,\y) to[out=90,in=270] (2.2,\y+.5);
\draw[shorten >=1pt, shorten <=1pt, line width=5pt, white] (2.2,\y) to[out=90,in=270] (2,\y+.5); \draw (2.2,\y) to[out=90,in=270] (2,\y+.5);}
\foreach \x in {1,1.2}{\draw (\x,1) -- (\x+1,2);} 
\draw[shorten >=.4cm, shorten <=.4cm, line width=8pt, white] (2.1,1) -- (1.1,2);
\foreach \x in {2,2.2}{\draw (\x,1) -- (\x-1,2);} 
\foreach \x in {1,1.2}{\draw (\x,2) -- (\x,3);} 
\foreach \y in {2,2.5}{\draw (2,\y) to[out=90,in=270] (2.2,\y+.5);
\draw[shorten >=1pt, shorten <=1pt, line width=5pt, white] (2.2,\y) to[out=90,in=270] (2,\y+.5); \draw (2.2,\y) to[out=90,in=270] (2,\y+.5);}
\foreach \x in {1,1.2}{\draw (\x,3) -- (\x+1,4);} 
\draw[shorten >=.4cm, shorten <=.4cm, line width=8pt, white] (2.1,3) -- (1.1,4);
\foreach \x in {2,2.2}{\draw (\x,3) -- (\x-1,4);}
\foreach \y in {0,4}{\draw (1,\y) -- (1.2,\y); \draw (2,\y) -- (2.2,\y);}
\node at (3,2) {\Huge{$=$}};
\end{tikzpicture}\qquad
\begin{tikzpicture}[xscale=1.2, yscale=.6, thick]
\foreach \x in {1,1.2}{\draw (\x,0) -- (\x+1,1);} 
\draw[shorten >=.4cm, shorten <=.4cm, line width=8pt, white] (2.1,0) -- (1.1,1);
\foreach \x in {2,2.2}{\draw (\x,0) -- (\x-1,1);} 
\foreach \x in {1,1.2}{\draw (\x,1) -- (\x,2);} 
\foreach \y in {1,1.5}{\draw (2,\y) to[out=90,in=270] (2.2,\y+.5);
\draw[shorten >=1pt, shorten <=1pt, line width=5pt, white] (2.2,\y) to[out=90,in=270] (2,\y+.5); \draw (2.2,\y) to[out=90,in=270] (2,\y+.5);}
\foreach \x in {1,1.2}{\draw (\x,2) -- (\x+1,3);} 
\draw[shorten >=.4cm, shorten <=.4cm, line width=8pt, white] (2.1,2) -- (1.1,3);
\foreach \x in {2,2.2}{\draw (\x,2) -- (\x-1,3);}
\foreach \x in {1,1.2}{\draw (\x,3) -- (\x,4);} 
\foreach \y in {3,3.5}{\draw (2,\y) to[out=90,in=270] (2.2,\y+.5);
\draw[shorten >=1pt, shorten <=1pt, line width=5pt, white] (2.2,\y) to[out=90,in=270] (2,\y+.5); \draw (2.2,\y) to[out=90,in=270] (2,\y+.5);}
\foreach \y in {0,4}{\draw (1,\y) -- (1.2,\y); \draw (2,\y) -- (2.2,\y);}
\end{tikzpicture}\end{center}
One can also verify this relation directly by the calculation \[r_1r_2 = (s_1;0,0) (\id_2;0,1)=(s_1;0,1).\]  So we have \[r_1r_2r_1r_2 = (s_1;0,1)(s_1;0,1) = (s_1s_1;1,1).\]  A similar calculation, starting with $r_2r_1=(s_1;1,0)$, yields the same result.\dqed
\end{example}

\begin{example}[Direct Sum Ribbon]\label{ex:direct-sum-ribbon}
Similar\index{direct sum!ribbon} to Example \ref{ex:direct-sum-braid}, the direct sum map on the symmetric groups admits a lift to the ribbon groups, so there is a commutative diagram 
\[\nicexy@C+.7cm{R_{k_1} \times \cdots \times R_{k_n} \ar[r]^-{\text{direct sum}} \ar[d]_-{(\pi,\ldots,\pi)} & R_{k_1+\cdots+k_n} \ar[d]^-{\pi}\\ S_{k_1} \times \cdots \times S_{k_n} \ar[r]^-{\text{direct sum}} & S_{k_1+\cdots+k_n}}\] 
in the category of groups.  Algebraically, in terms of the semi-direct product, given ribbons \[\sigma_i = (b_i;\uz_i) \in R_{k_i} = B_{k_i} \rtimes \bbZ^{k_i} \forspace 1 \leq i \leq n,\] the direct sum ribbon is \[\sigma_1 \oplus \cdots \oplus \sigma_n = \bigl(\overbrace{b_1 \oplus \cdots \oplus b_n}^{\text{direct sum braid}}; \overbrace{\uz_1,\ldots,\uz_n}^{\text{concatenation}} \bigr) \in B_{k_1+\cdots+k_n} \rtimes \bbZ^{k_1+\cdots+k_n}.\]  Geometrically, the direct sum ribbon $\sigma_1 \oplus \cdots \oplus \sigma_n \in R_{k_1+\cdots+k_n}$ is the ribbon obtained by placing $\sigma_1,\ldots,\sigma_n$ side-by-side from left to right.\dqed
\end{example}

\begin{example}[Block Ribbon]\label{ex:block-ribbon}
Similar\index{block ribbon} to Example \ref{ex:block-braid}, the block permutation map lifts to the ribbon groups, so there is a commutative diagram \[\nicexy@C+1.5cm{R_n \ar[d]_-{\pi} \ar[r]^-{\text{block ribbon}} & R_{k_1+\cdots+k_n} \ar[d]^-{\pi}\\ S_n \ar[r]^-{\text{block permutation}} & S_{k_1+\cdots+k_n}}\] in $\Set$ for $n \geq 1$ and $k_1,\ldots,k_n \geq 0$.  Geometrically, the block ribbon \[\sigma\langle k_1,\ldots,k_n\rangle \in R_{k_1+\cdots+k_n}\] is obtained from the ribbon $\sigma \in R_n$ by dividing its $i$th strip, for each $1 \leq i \leq n$, into $k_i$ strips, twisted individually and around each other according to the number of full $2\pi$ twists in the $i$th strip in $\sigma$.

For example, the generator 
\begin{center}\begin{tikzpicture}[xscale=1.2, yscale=.6, thick]
\foreach \y in {0,.5}{\draw (0,\y) to[out=90,in=270] (.2,\y+.5);
\draw[shorten >=1pt, shorten <=1pt, line width=5pt, white] (.2,\y) to[out=90,in=270] (0,\y+.5);
\draw (.2,\y) to[out=90,in=270] (0,\y+.5);}
\foreach \x in {0,1} \draw (0,\x) -- (.2,\x); \node at (-1,.5) {$(r_1 \in R_1) = $};
\end{tikzpicture}\end{center}
has associated block ribbon \[r_1\langle 2\rangle = r_1r_2r_1r_2 \in R_2,\] as illustrated in the picture below.
\begin{center}\begin{tikzpicture}[xscale=.7, yscale=.5, thick]
\draw [gray!50, dotted, line width=2pt] (1,0) -- (1,4);
\foreach \y in {0,2}{\draw (0,\y) to[out=90,in=270] (2,\y+2);
\draw[shorten >=1pt, shorten <=1pt, line width=5pt, white] (2,\y) to[out=90,in=270] (0,\y+2);
\draw (2,\y) to[out=90,in=270] (0,\y+2);}
\foreach \x in {0,4} \draw (0,\x) -- (2,\x); 
\node at (3.3,2) {\Huge{$=$}};
\end{tikzpicture}\qquad
\begin{tikzpicture}[xscale=1.2, yscale=.5, thick]
\foreach \x in {1,1.2}{\draw (\x,0) -- (\x,1);} 
\foreach \y in {0,.5}{\draw (2,\y) to[out=90,in=270] (2.2,\y+.5);
\draw[shorten >=1pt, shorten <=1pt, line width=5pt, white] (2.2,\y) to[out=90,in=270] (2,\y+.5); \draw (2.2,\y) to[out=90,in=270] (2,\y+.5);}
\foreach \x in {1,1.2}{\draw (\x,1) -- (\x+1,2);} 
\draw[shorten >=.4cm, shorten <=.4cm, line width=8pt, white] (2.1,1) -- (1.1,2);
\foreach \x in {2,2.2}{\draw (\x,1) -- (\x-1,2);} 
\foreach \x in {1,1.2}{\draw (\x,2) -- (\x,3);} 
\foreach \y in {2,2.5}{\draw (2,\y) to[out=90,in=270] (2.2,\y+.5);
\draw[shorten >=1pt, shorten <=1pt, line width=5pt, white] (2.2,\y) to[out=90,in=270] (2,\y+.5); \draw (2.2,\y) to[out=90,in=270] (2,\y+.5);}
\foreach \x in {1,1.2}{\draw (\x,3) -- (\x+1,4);} 
\draw[shorten >=.4cm, shorten <=.4cm, line width=8pt, white] (2.1,3) -- (1.1,4);
\foreach \x in {2,2.2}{\draw (\x,3) -- (\x-1,4);}
\foreach \y in {0,4}{\draw (1,\y) -- (1.2,\y); \draw (2,\y) -- (2.2,\y);}
\end{tikzpicture}\end{center}
Its underlying braid is $s_1s_1\in PB_2$, and each strip has one full $2\pi$ twist.

This example illustrates that the block braid map does \emph{not} lift to the block ribbon map.  In other words, the diagram \[\nicexy@C+1.2cm{R_n \ar[d]_-{\pi} \ar[r]^-{\text{block ribbon}} & R_{k_1+\cdots+k_n} \ar[d]^-{\pi}\\ B_n \ar[r]^-{\text{block braid}} & B_{k_1+\cdots+k_n}}\] in $\Set$ is not commutative in general.  For example, with $n=1$ and $k_1=2$, the diagram \[\nicexy@C+1.2cm{R_1 \ar[d]_-{\pi} \ar[r]^-{\text{block ribbon}} & R_{2} \ar[d]^-{\pi}\\ B_1 \ar[r]^-{\text{block braid}} & B_{2}}\] is not commutative because $\pi r_1 \in B_1$ is the identity braid, so its induced block braid is the identity braid in $B_2$.  On the other hand, we observed above that the underlying braid of the block ribbon $r_1\langle 2\rangle$ is the non-identity braid $s_1s_1\in PB_2$.\dqed
\end{example}

\begin{example}[Block Ribbon of Generators]
Writing $k = k_1+\cdots+k_n$, for $1\leq i \leq n-1$, the block ribbon induced by the generating ribbon $r_i \in R_n$ is given by \[r_i\langle k_1,\ldots,k_n\rangle = \bigl(\overbrace{s_i\langle k_1,\ldots,k_n\rangle}^{\text{block braid}}; \overbrace{0,\ldots,0}^{k}\bigr) \in B_k \rtimes \bbZ^{k} = R_k,\] as illustrated below.
\begin{center}\begin{tikzpicture}[xscale=1.2, yscale=1.6, thick]
\foreach \x in {1,2,7,8} \draw (\x,0) rectangle (\x+.2,1);
\node at (1.6,.5) {$\cdots$}; \node at (7.6,.5) {$\cdots$};
\foreach \x in {3,4}{\draw (\x,0) -- (\x+.2,0) -- (\x+2.2,1) -- (\x+2,1) -- cycle;}
\node at (3.6,0) {$\cdots$}; \node at (5.5,.95) {$\cdots$};
\draw[shorten >=.2, shorten <=.2, line width=1.1cm, white] (5.6,0) -- (3.6,1);
\foreach \x in {5,6}{\draw (\x,0) -- (\x+.2,0) -- (\x-1.8,1) -- (\x-2,1) -- cycle;}
\node at (4.6,.5) {$\cdots$};
\draw[thin] (1,-.1) -- (1,-.2) -- (2.2,-.2) -- (2.2,-.1);
\node at (1.6,-.3) {\scriptsize{$k_1+\cdots+k_{i-1}$}};
\draw[thin] (3,-.1) -- (3,-.2) -- (4.2,-.2) -- (4.2,-.1);
\node at (3.6,-.3) {\scriptsize{$k_{i}$}};
\draw[thin] (5,-.1) -- (5,-.2) -- (6.2,-.2) -- (6.2,-.1);
\node at (5.6,-.3) {\scriptsize{$k_{i+1}$}};
\draw[thin] (7,-.1) -- (7,-.2) -- (8.2,-.2) -- (8.2,-.1);
\node at (7.6,-.3) {\scriptsize{$k_{i+2}+\cdots+k_{n}$}};
\end{tikzpicture}\end{center}
So each of the $k_{i+1}$ strips in the $(i+1)$st block crosses over each of the $k_i$ strips in the $i$th block, and none of the strips has any twists.\dqed
\end{example}

\begin{example}[Block Ribbon of the Last Generator]
For the generating ribbon $r_n \in R_n$, writing $k^{n-1}=k_1+\cdots +k_{n-1}$, the block ribbon is \[r_n\langle k_1,\ldots,k_n\rangle = \bigl(\overbrace{\id_{k^{n-1}} \oplus b_{k_n}}^{\text{direct sum braid}}; \underbrace{0,\ldots,0}_{k^{n-1}}, \underbrace{1,\ldots,1}_{k_n}\bigr) \in B_k \rtimes \bbZ^k = R_k,\] with $\id_{k^{n-1}} \in B_{k^{n-1}}$ the identity braid and \[b_{k_n} = \Bigl[\bigl(s_{k_n-1}s_{k_n-2}\cdots s_1\bigr) \cdots \bigl(s_{k_n-1}s_{k_n-2}s_{k_n-3}\bigr) \bigl(s_{k_n-1}s_{k_n-2}\bigr)s_{k_n-1}\Bigr]^2 \in PB_{k_n},\] as illustrated below.
\begin{center}\begin{tikzpicture}[xscale=1, yscale=.6, thick]
\draw (1,0) to[out=90,in=270] (5,3);
\draw[shorten >=.2, shorten <=.2, line width=.2cm, white] (5,2) to[out=90,in=270] (4,3); \draw (2,0) to[out=90,in=270] (5,2) to[out=90,in=270] (4,3);
\draw[shorten >=.2, shorten <=.2, line width=.2cm, white] (5,1) to[out=90,in=270] (3,3); 
\draw (3,0) to[out=90,in=270] (5,1) to[out=90,in=270] (3,3);
\draw[shorten >=.2, shorten <=.2, line width=.2cm, white] (5,0) to[out=90,in=270] (1,3);
\draw (5,0) to[out=90,in=270] (1,3);
\draw (1,3) to[out=90,in=270] (5,6);
\draw[shorten >=.2, shorten <=.2, line width=.2cm, white] (5,5) to[out=90,in=270] (3,6); \draw (3,3) to[out=90,in=270] (5,5) to[out=90,in=270] (3,6);
\draw[shorten >=.2, shorten <=.2, line width=.2cm, white] (5,4) to[out=90,in=270] (2,6); \draw (4,3) to[out=90,in=270] (5,4) to[out=90,in=270] (2,6);
\draw[shorten >=.2, shorten <=.2, line width=.2cm, white] (5,3) to[out=90,in=270] (1,6); \draw (5,3) to[out=90,in=270] (1,6);
\node at (4,.1) {$\cdots$}; \node at (2,3) {$\cdots$}; \node at (4,5.9) {$\cdots$};
\foreach \x in {1,2,3}{\node at (\x,-.3) {\x};} \node at (5,-.3) {$k_n$};
\node at (-.5,3) {$b_{k_n} =$}; 
\end{tikzpicture}\end{center}
For instance, if $k_n=3$, then \[b_3=(s_2s_1)s_2(s_2s_1)s_2 \in PB_3.\] To visualize the pure braid $b_{k_n}$, one begins with two horizontal bars, one above the other, and places $k_n$ parallel vertical strings between them whose ends are tied to the bars.  Holding the top horizontal bar fixed, one rotates the bottom horizontal bar by one full round in the counter-clockwise direction as viewed from the top.  The resulting braid is $b_{k_n}$.\dqed
\end{example}

\section{Ribbon Operads as Monoids}\label{sec:ribbon-operad}

The purpose of this section is to define (pure) ribbon operads using the (pure) ribbon groups.  To see that the (pure) ribbon groups have the required operadic composition, let us first consider the following example.

\begin{example}[Comp-$i$ of Ribbons]\label{ex:compi-ribbon}
Combining\index{ribbon group operad!composition} Example \ref{ex:direct-sum-ribbon} and Example \ref{ex:block-ribbon}, for ribbons $\sigma \in R_n$ and $\tau \in R_m$ and $1 \leq i \leq n$, we define
\begin{equation}\label{ribbon-compi}
\sigma \compi \tau = \overbrace{\sigma\langle \underbrace{1,\ldots,1}_{i-1}, m,\underbrace{1,\ldots,1}_{n-i}\rangle}^{\text{block ribbon}} \cdot \overbrace{\bigl(\underbrace{\id_1 \oplus \cdots \oplus \id_1}_{i-1} \oplus \tau \oplus \underbrace{\id_1 \oplus \cdots \oplus \id_1}_{n-i}\bigr)}^{\text{direct sum ribbon}} \in R_{n+m-1}
\end{equation}
as the product of a direct sum ribbon induced by $\tau$ with a block ribbon induced by $\sigma$.  Geometrically, $\sigma \compi \tau$ is the ribbon obtained from $\sigma$ by replacing its $i$th strip by the ribbon $\tau$, twisted according to the number of full $2\pi$ twists in the $i$th strip in $\sigma$.  Similar to Section \ref{sec:braided-operad}, by an inspection of relevant pictures, one can see that the ribbon groups form an action operad.\dqed
\end{example}

\begin{definition}\label{def:ribbon-group-operad}
The\index{ribbon group operad}\index{group operad!ribbon} \emph{ribbon group operad}\label{not:ribbongop} $\R$ is the action operad with:
\begin{itemize}
\item $\R(n) = R_n$, the $n$th ribbon group;
\item the $\compi$-composition as in \eqref{ribbon-compi};
\item the augmentation $\pi : \R \to \S$ given level-wise by the underlying permutation map $\pi : R_n \to S_n$.
\end{itemize}
\end{definition}

\begin{definition}\label{def:pure-ribbon-group-operad}
The\index{pure ribbon group operad} \emph{pure ribbon group operad} $\PR$ is the action operad with:
\begin{itemize}
\item $\PR(n) = PR_n$, the $n$th pure ribbon group;
\item the $\compi$-composition as in \eqref{ribbon-compi};
\item the augmentation $\rho : \PR \to \S$ given level-wise by $\rho(r)=\id_n$ for $r\in PR_n$.
\end{itemize}
\end{definition}

We now restrict Section \ref{sec:group-operad} to the action operads $\R$ and $\PR$.  Following the planar, symmetric, and braided cases, for the action operad $\R$ (resp., $\PR$), we will use the adjective \emph{ribbon} (resp., \emph{pure ribbon}).

\begin{definition}\label{def:ribbon-operad}
Suppose $\colorc$ is a set.  With the action operad $\G=\R$:
\begin{enumerate}
\item The objects of the diagram category\label{notation:rseqcm} \[\rseqcm = \M^{\rsubcopc}\] are called \index{ribbon sequence}\emph{$\colorc$-colored ribbon sequences} in $\M$.  
\item For $X,Y  \in \rseqcm$, $X \circr Y \in \rseqcm$ is called the \index{circle product!ribbon}\index{ribbon circle product}\emph{$\colorc$-colored ribbon circle product}.  It is defined entrywise as the coend\label{not:circr} \[(X \circr Y)\dub = \int^{\uc \in \rsubc} X\duc \otimes Y^{\uc}_{\smallsubr}(\ub)\in \M\] for $\dub \in \rsubcopc$.  The object $Y^{\uc}_{\smallsubr} \in \M^{\rsubcop}$ is defined entrywise as the coend \[Y^{\uc}_{\smallsubr}(\ub) = \int^{\{\ua_j\}\in\prod_{j=1}^m \rsubcop} \rsubcop\bigl(\ua_1,\ldots,\ua_m;\ub\bigr) \cdot \left[\bigotimes_{j=1}^m Y\cjuaj\right] \in \M\] for $\ub \in \rsubcop$, in which $(\ua_1,\ldots,\ua_m) \in \rsubcop$ is the concatenation.  
\item The category of \index{operad!ribbon}\index{ribbon operad}\emph{$\colorc$-colored ribbon operads in $\M$} is the category\label{not:roperadcm} \[\Roperadcm = \Mon\bigl(\rseqcm, \circr, \I^{\R}\bigr)\] of monoids in the monoidal category $\bigl(\rseqcm,\circr,\I^{\R}\bigr)$.  If $\colorc$ has $n<\infty$ elements, we also refer to objects in $\Roperadcm$ as \emph{$n$-colored ribbon operads in $\M$}.
\item Suppose $\O$ is a $\colorc$-colored ribbon operad in $\M$.  The category $\algmo$ \index{ribbon operad!algebra}\index{algebra!of a ribbon operad}of \emph{$\O$-algebras} is defined as the category of $(\O\circr -)$-algebras for the monad $\O\circr -$ in $\Mtoc$.
\end{enumerate}
\end{definition}

The previous definition also works with the pure ribbon group operad $\PR$, which defines $\colorc$-colored \index{pure ribbon sequence}pure ribbon sequences, the $\colorc$-colored \index{circle product!pure ribbon}pure ribbon circle product, $\colorc$-colored \index{operad!pure ribbon}\index{pure ribbon operad}pure ribbon operads, and algebras over them.  The explicit descriptions of $\G$-operads in Proposition \ref{def:g-operad-generating} and Proposition \ref{prop:g-operad-compi} now specialize to the ribbon group operad $\R$ and the pure ribbon group operad $\PR$.  Moreover, the description of algebras over a $\colorc$-colored $\G$-operad in Proposition \ref{prop:g-operad-algebra-defs} specializes to the (pure) ribbon case.  Next we provide some examples of (pure) ribbon operads and their algebras.

\section{Examples of Ribbon Operads}\label{sec:ribbon-ex}

Some examples of (pure) ribbon operads are given in this section.

\begin{example}[Ribbon Operads as Pure Ribbon Operads]\label{ex:rop-prop}
Every $\colorc$-colored ribbon operad is naturally a $\colorc$-colored pure ribbon operad by restriction of the equivariant structure.\dqed
\end{example}

\begin{example}[Symmetric Operads as Ribbon Operads]\label{ex:sop-rop}
Every\index{ribbon operad!induced by a symmetric operad} $\colorc$-colored symmetric operad is naturally a $\colorc$-colored ribbon operad via the underlying permutation maps $\pi : R_n \to S_n$ for $n \geq 0$.\dqed\end{example}

\begin{example}[Colored Endomorphism Ribbon Operads]\label{ex:end-ribbon}
For each $\colorc$-colored object $X=\{X_c\}_{c\in\colorc}$ in $\M$, the endomorphism operad\index{endomorphism operad!ribbon}\index{ribbon operad!endomorphism} $\End(X)$, which is a $\colorc$-colored symmetric operad, is naturally a $\colorc$-colored ribbon operad.  This is the $\G=\R$ special case of Example \ref{ex:end-g}.\dqed\end{example}

\begin{example}[Ribbon Group Operad as a Ribbon Operad]
The ribbon group operad $\R=\{R_n\}_{n\geq 0}$ yields a one-colored ribbon operad in $\Set$.  This is Example \ref{ex:agop-goperad} when $\G=\R$.  Similarly, the pure ribbon group operad $\PR = \{PR_n\}_{n\geq 0}$ yields a one-colored pure ribbon operad in $\Set$.\dqed
\end{example}

\begin{example}[Translation Category of Ribbon Group Operad]
Applying the \index{ribbon group operad!translation category}\index{translation category!ribbon group operad}translation category construction $E_?$ to the ribbon group operad $\R$, we obtain a one-colored ribbon operad \[E_{\R} = \Bigl(\bigl\{E_{R_n}\bigr\}_{n\geq 0}, \gamma,\operadunit\Bigr)\] in the symmetric monoidal category of small groupoids.  This is Example \ref{ex:translation-cat-goperad} when $\G=\R$.  Similarly, the pure ribbon group operad $\PR$ yields a one-colored pure ribbon operad \[E_{\PR} = \Bigl(\bigl\{E_{PR_n}\bigr\}_{n\geq 0}, \gamma,\operadunit\Bigr)\] in the symmetric monoidal category of small groupoids.\dqed
\end{example}

\begin{example}[Nerves of Translation Categories of Ribbon Groups]\label{ex:translation-ribbon}
Applying the nerve functor to the one-colored ribbon operad $E_{\R}$ in small groupoids, we obtain a one-colored ribbon operad\index{translation category!nerve of} \[E\R = \{ER_n\}_{n\geq 0}\] in the category of simplicial sets, where \[ER_n = \Nerve(E_{R_n})\] is the nerve of the translation category $E_{R_n}$ of the $n$th ribbon group $R_n$.  This is the $\G=\R$ special case of Example \ref{ex:translation-g}.  Furthermore, the one-colored ribbon operad $E\R$ is an\index{Rinfinity-operad@$\R_\infty$-operad}\index{operad!Rinfinity@$\R_\infty$} \emph{$\R_\infty$-operad} \cite{wahl}.  That is, each $ER_n$ is contractible, and the ribbon group action \[ER_n \times R_n \to ER_n\] is free and proper.  

Similarly, there is a one-colored pure ribbon operad \[E\PR = \{EPR_n\}_{n\geq 0}\] in the category of simplicial sets, where \[EPR_n = \Nerve(E_{PR_n})\] is the nerve of the translation category $E_{PR_n}$ of the $n$th pure ribbon group $PR_n$. Moreover, each $EPR_n$ is contractible, and the pure ribbon group action \[EPR_n \times PR_n \to EPR_n\] is free and proper.\dqed\end{example}

\begin{example}[Parenthesized Ribbon Operad]\label{ex:parenthesized-ribbon}
The \index{parenthesized ribbon operad}\emph{parenthesized ribbon operad}\label{not:par} $\PaR$ is the one-colored ribbon operad in $\Set$ as in Example \ref{ex:parenthesized-goperad} with $\G=\R$.  Its entries are \[\PaR(n) = \Mag(n) \times R_n \times \Mag(n)\] for $n\geq 0$.  For $1 \leq i \leq n$ and $m \geq 0$, the $\compi$-composition \[\nicexy{\PaR(n) \times \PaR(m) \ar[r]^-{\compi} & \PaR(n+m-1)}\] is defined as \[(w;\sigma;w') \compi (x;\tau;x') = \bigl(w\compi x; \sigma\compi\tau; w' \comp_{\sigmabar(i)} x'\bigr)\] for $w,w'\in \Mag(n)$, $x,x' \in \Mag(m)$, $\sigma\in R_n$, and $\tau\in R_m$, where $\sigma\compi\tau\in R_{n+m-1}$ is the $\compi$-composition of ribbons in \eqref{ribbon-compi}.  The $\compi$-composition in $\PaR$ has a geometric interpretation similar to the one in the parenthesized braided operad $\PaB$ in Example \ref{ex:parenthesized-braid}.

Similarly, using the pure ribbon group operad $\G=\PR$ in Example \ref{ex:parenthesized-goperad}, we obtain the \index{parenthesized pure ribbon operad}\emph{parenthesized pure ribbon operad} $\PaPR$, which is a one-colored pure ribbon operad in $\Set$.\dqed
\end{example}

\section{Universal Cover of the Framed Little $2$-Disc Operad}\label{sec:ribbon-operad-examples}

The purpose of this section is to discuss the one-colored ribbon operad obtained by taking the universal covers of the spaces in the framed little $2$-disc operad.  This example is due to Wahl \cite{wahl}.

\begin{example}[Universal Cover of the Framed Little $2$-Disc Operad]\label{ex:universal-cover-fd2}
This\index{framed little $n$-disc operad!universal cover}\index{operad!universal cover of framed little $2$-disc}  example provides another example of an $\R_\infty$-operad.  We proceed as in Example \ref{ex:universal-cover-d2}, which in turn is an adaptation of Example \ref{ex:universal-cover-c2}, using the framed little $2$-disc operad $\Df_2$ in Example \ref{ex:framed-little-n-disc} instead of the little $2$-disc operad $\D_2$.  The result is a one-colored ribbon operad $\Dftildetwo$ in $\CHau$ that is also an $\R_\infty$-operad.  Let us point out the modifications that we need to make in Example \ref{ex:universal-cover-d2}.  

Let\label{not:dftilde} $\Dftildetwo(n)$ be the universal cover of the space $\Df_2(n)$ at the following base point.
\begin{center}\begin{tikzpicture}[scale=1.2, thick]
\draw (0,0) circle(1); \draw[very thick] (0,.9) -- (0,1.1);
\foreach \x in {-.8,-.4,.8} {\draw (\x,0) circle (.2); \draw[very thick] (\x,.15) -- (\x,.25);}
\node at (-.8,0) {1}; \node at (-.4,0) {2}; \node at (.8,0) {$n$};
\node at (.2,0) {$\cdots$}; \node at (-2.3,0) {$\udnfr \in \Df_2(n)$};
\end{tikzpicture}\end{center}
This is the image of $\ud^n\in \D_2(n)$ in $\Df_2(n)$.  In the picture above, we indicate the rotation of a (little) $2$-disc by marking the image of the point $(0,1)$.  In other words, $\udnfr$ is $\ud^n \in \D_2(n)$ in which the $n$ framed little $2$-discs all have identity rotation.  The operadic unit in $\Dftildetwo(1)$ is the constant path at the base point $\udfr^1$.

To define the operadic composition on $\Dftildetwo$, suppose $\alpha \in \Dftildetwo(n)$ with $n \geq 1$, $\alpha_i \in \Dftildetwo(k_i)$ with $1 \leq i \leq n$, and $k=k_1+\cdots+k_n$.   We modify the path $\gamma^0 : [0,1] \to \D_2(k)$ in Example \ref{ex:universal-cover-d2} to a path \[\gammazerofr : [0,1] \to \Df_2(k)\] such that:
\begin{itemize}
\item $\gammazerofr(0) = \udkfr \in \Df_2(k)$.
\item $\gammazerofr(1) = \gamma\bigl(\alpha(0); \alpha_1(0), \ldots, \alpha_n(0)\bigr) = \gamma\bigl(\udnfr; \udfr^{k_1},\ldots,\udfr^{k_n}\bigr) \in \Df_2(k)$.
\item For each $0 < t < 1$, $\gammazerofr(t) \in \Df_2(k)$ consists of $k$ framed little $2$-discs with identity rotation that are ordered from left to right and are centered at the horizontal axis.  Their diameters partition $[-1,1]$ into $k$ sub-intervals, not necessarily of equal length.  
\end{itemize}

To define the ribbon group action \[\Dftildetwo(n) \times R_n \to \Dftildetwo(n),\] it is enough to define the action by the generators $r_i \in R_n$ for $1 \leq i \leq n$.  For $1 \leq i \leq n-1$, the $r_i$-action is defined using a modification \[\rtilde_i : [0,1] \to \Df_2(n)\] of the path $\stilde_i : [0,1] \to \D_2(n)$ in Example \ref{ex:universal-cover-d2}.  It satisfies \[\rtilde_i(0) = \udnfr \andspace \rtilde_i(1) = \udnfr\tau_i,\] and the $n$ framed little $2$-discs in $\rtilde_i(t)$ all have identity rotation for $0 \leq t \leq 1$, as illustrated by the following sequence of pictures.
\begin{center}
\begin{tikzpicture}[scale=1.2,thick]
\draw (0,0) circle (1); \draw[very thick] (0,.9) -- (0,1.1);
\draw [fill=green!20] (.2,0) circle (.2); \draw [fill=yellow!20] (-.2,0) circle (.2);
\foreach \x in {-.2,.2} \draw[very thick] (\x,.15) -- (\x,.25);
\node at (-.2,0) {\scriptsize{$i$}}; \node at (.2,0) {\scriptsize{$i\!\!+\!\!1$}};
\foreach \x in {-.7,.7} \node at (\x,0) {$\cdots$}; \node at (0,1.3) {$\rtilde_i(0)$};
\end{tikzpicture}\qquad
\begin{tikzpicture}[scale=1.2,thick]
\draw (0,0) circle (1); \draw[very thick] (0,.9) -- (0,1.1);
\draw[fill=yellow!20] (-.141,.141) circle (.2); \draw[fill=green!20] (.141,-.141) circle (.2); 
\draw[very thick] (.141,0) -- (.141,.1); \draw[very thick] (-.141,.291) -- (-.141,.391);
\node at (-.141,.141) {\scriptsize{$i$}}; \node at (.141,-.141) {\scriptsize{$i\!\!+\!\!1$}};
\foreach \x in {-.7,.7} \node at (\x,0) {$\cdots$}; \node at (0,1.3) {$\rtilde_i(1/3)$};
\end{tikzpicture}\qquad
\begin{tikzpicture}[scale=1.2,thick]
\draw (0,0) circle (1); \draw[very thick] (0,.9) -- (0,1.1);
\draw[fill=yellow!20] (.141,.141) circle (.2); \draw[fill=green!20] (-.141,-.141) circle (.2); 
\draw[very thick] (-.141,0) -- (-.141,.1); \draw[very thick] (.141,.291) -- (.141,.391);
\node at (.141,.141) {\scriptsize{$i$}}; \node at (-.141,-.141) {\scriptsize{$i\!\!+\!\!1$}};
\foreach \x in {-.7,.7} \node at (\x,0) {$\cdots$}; \node at (0,1.3) {$\rtilde_i(2/3)$};
\end{tikzpicture}\qquad
\begin{tikzpicture}[scale=1.2,thick]
\draw (0,0) circle (1); \draw[very thick] (0,.9) -- (0,1.1);
\draw[fill=yellow!20] (.2,0) circle (.2); \draw[fill=green!20] (-.2,0) circle (.2);
\foreach \x in {-.2,.2} \draw[very thick] (\x,.15) -- (\x,.25);
\node at (-.2,0) {\scriptsize{$i\!\!+\!\!1$}}; \node at (.2,0) {\scriptsize{$i$}};
\foreach \x in {-.7,.7} \node at (\x,0) {$\cdots$}; \node at (0,1.3) {$\rtilde_i(1)$};
\end{tikzpicture}\end{center}
For $0 \leq t \leq 1$ the framed little $2$-discs in $\rtilde_i(t)$ labeled by $1,\ldots,i-1,i+2,\ldots,n$ remain unchanged.  

To define the $r_n$-action, we use the path \[\rtilde_n : [0,1] \to \Df_2(n)\] with \[\rtilde_n(0) = \udnfr = \rtilde_n(1)\] determined by the following sequence of pictures.
\begin{center}
\begin{tikzpicture}[scale=1.2,thick]
\draw (0,0) circle (1); \draw[very thick] (0,.9) -- (0,1.1);
\foreach \x in {-.8,0,.8} \draw (\x,0) circle (.2); 
\foreach \x in {-.8,0} \draw[very thick] (\x,.15) -- (\x,.25);
\draw[very thick, red] (.8,.15) -- (.8,.25);
\node at (-.8,0) {\scriptsize{$1$}}; \node at (0,0) {\scriptsize{$i$}}; \node at (.8,0) {\scriptsize{$n$}};
\foreach \x in {-.4,.4} \node at (\x,0) {$\cdots$}; \node at (0,1.3) {$\rtilde_n(0)$};
\end{tikzpicture}\qquad
\begin{tikzpicture}[scale=1.2,thick]
\draw (0,0) circle (1); \draw[very thick] (0,.9) -- (0,1.1);
\foreach \x in {-.8,0,.8} \draw (\x,0) circle (.2); 
\foreach \x in {-.8,0} \draw[very thick] (\x,.15) -- (\x,.25);
\draw[very thick, red] (.92,-.1) -- (1,-.175);
\node at (-.8,0) {\scriptsize{$1$}}; \node at (0,0) {\scriptsize{$i$}}; \node at (.8,0) {\scriptsize{$n$}};
\foreach \x in {-.4,.4} \node at (\x,0) {$\cdots$}; \node at (0,1.3) {$\rtilde_n(1/3)$};
\end{tikzpicture}\qquad
\begin{tikzpicture}[scale=1.2,thick]
\draw (0,0) circle (1); \draw[very thick] (0,.9) -- (0,1.1);
\foreach \x in {-.8,0,.8} \draw (\x,0) circle (.2); 
\foreach \x in {-.8,0} \draw[very thick] (\x,.15) -- (\x,.25);
\draw[very thick, red] (.68,-.1) -- (.6,-.175); 
\node at (-.8,0) {\scriptsize{$1$}}; \node at (0,0) {\scriptsize{$i$}}; \node at (.8,0) {\scriptsize{$n$}};
\foreach \x in {-.4,.4} \node at (\x,0) {$\cdots$}; \node at (0,1.3) {$\rtilde_n(2/3)$};
\end{tikzpicture}\qquad
\begin{tikzpicture}[scale=1.2,thick]
\draw (0,0) circle (1); \draw[very thick] (0,.9) -- (0,1.1);
\foreach \x in {-.8,0,.8} \draw (\x,0) circle (.2); 
\foreach \x in {-.8,0} \draw[very thick] (\x,.15) -- (\x,.25);
\draw[very thick, red] (.8,.15) -- (.8,.25);
\node at (-.8,0) {\scriptsize{$1$}}; \node at (0,0) {\scriptsize{$i$}}; \node at (.8,0) {\scriptsize{$n$}};
\foreach \x in {-.4,.4} \node at (\x,0) {$\cdots$}; \node at (0,1.3) {$\rtilde_n(1)$};
\end{tikzpicture}\end{center}
For $0 \leq t \leq 1$ the framed little $2$-discs in $\rtilde_n(t)$ labeled by $1,\ldots,n-1$ have identity rotation and remain unchanged.  The $n$th framed little $2$-disc in $\rtilde_n$ rotates exactly one round as $t$ goes from $0$ to $1$.  Then we define the $r_n$-action $\alpha r_n \in \Dftildetwo(n)$ on $\alpha \in \Dftildetwo(n)$ as the homotopy class of the composition of paths \[\alpha \cdot \rtilde_n : [0,1] \to \Df_2(n).\]

When the framed little $2$-disc operad $\Df_2$ is regarded as a one-colored ribbon operad, the level-wise universal covering map \[p : \Dftildetwo \to \Df_2\] is a morphism of one-colored ribbon operads.  Moreover, $\Dftildetwo$ is an $\R_\infty$-operad.  That is, each space $\Dftildetwo(n)$ is contractible, and the ribbon group action \[\Dftildetwo(n) \times R_n \to \Dftildetwo(n)\] is free and proper.\dqed
\end{example}

\chapter{Cactus Operads}\label{ch:cactus}

In this chapter, we discuss an example of an action operad that is made up of the cactus groups in \cite{djs,dev,hen-kam}.  In these references, the cactus groups arise in connection with:
\begin{itemize}
\item the geometry of\index{moduli space of stable curves} moduli spaces of real, genus 0, stable curves with marked points \cite{dev};
\item symmetry groups of closed connected manifolds with a cubical cell structure \cite{djs}; 
\item categorical aspects of \index{quantum group}quantum groups \cite{hen-kam}, namely, groups that can act on multiple monoidal products in a coboundary monoidal category in the sense of Drinfel'd \cite{drinfeld}.
\end{itemize}
We will define the planar operad structure on the cactus groups purely algebraically, without involving any geometric objects.  The action operad structure on the cactus groups was first made explicit by Gurski \cite{gurski}.

We define the cactus groups in Section \ref{sec:cactus-groups}.  As in the symmetric group operad $\S$, the braid group operad $\B$, and the ribbon group operad $\R$, the operadic composition on the cactus groups is a product of two auxiliary operations, called direct sum cactus and block cactus.  In Section \ref{sec:direct-sum-cacti} we define direct sum cactus and study their properties.  In Section \ref{sec:block-cactus} we define block cactus and prove some preliminary results about them that are needed to define the planar operad structure on the cactus groups.  In Section \ref{sec:cactus-action-operad} we define the action operad structure on the cactus groups and the corresponding cactus operads.  Some examples of cactus operads are given in Section \ref{sec:cactus-ex}.  In Section \ref{sec:relationship-braid-cactus} we observe that there are \emph{no} group homomorphisms from the cactus groups to the braid groups, or from the braid groups to the cactus groups, that are compatible with their operad structures.

\section{Cactus Groups}\label{sec:cactus-groups}

One way to think about the cactus groups is that they are generated by an abstraction of permutations that reverse an interval.   To motivate the generators and relations in the cactus groups, we first consider the relevant permutations, for which we need the following preliminary concepts from \cite{hen-kam}.
 
\begin{definition}\label{def:disjoint-contain}
Suppose given integers $1 \leq p < q \leq n$ and $1 \leq r < t \leq n$.
\begin{enumerate}
\item We say that $p < q$ and $r < t$ are \index{disjoint}\emph{disjoint} if either $q < r$ or $t < p$.
\item We say that $p < q$ \index{contain}\emph{contains} $r < t$ if $p \leq r < t \leq q$.
\end{enumerate}
\end{definition}

Recall that $S_n$ denotes the symmetric group on $n$ letters.

\begin{definition}\label{def:interval-reversal}
For $1 \leq p<q \leq n$, define the permutation\label{not:rhonpq} $\rho^{(n)}_{p,q}\in S_n$ by
\[\rho^{(n)}_{p,q}(i)  = \begin{cases} i & \text{if $1 \leq i < p$ or $q<i\leq n$},\\
p+q-i & \text{if $p \leq i \leq q$}.
\end{cases}\]
We call $\rho^{(n)}_{p,q}$ an \index{interval-reversing permutation}\emph{interval-reversing permutation}, and we say that it \emph{covers the interval} $[p,q]$.  By convention, we set $\rho^{(n)}_{r,t}=\id_n$ when $r\geq t$.
\end{definition}

\begin{example}
The interval-reversing permutation $\rho^{(5)}_{2,4}$ is given by \[\nicexy{(1,2,3,4,5) \ar@{|->}[r]^-{\rho^{(5)}_{2,4}} & (1,4,3,2,5)},\] so the interval $[2,4]$ is reversed.  In general, the interval-reversing permutation $\rho^{(n)}_{p,q}$ reverses the order of the integers in the interval $[p,q]$.\dqed
\end{example}

The following relations in the symmetric groups can be checked by a direct inspection.

\begin{lemma}\label{permuation-cactus-axioms}
Suppose $1 \leq p < q \leq n$ and $1 \leq r < t \leq n$.  In the permutation group $S_n$, the following relations hold.
\begin{description}
\item[Involution] Each interval-reversing permutation is an involution: \[\rho^{(n)}_{p,q} \rho^{(n)}_{p,q} = \id_n.\]
\item[Disjointness] If $p < q$ and $r < t$ are disjoint, then \[\rho^{(n)}_{p,q} \rho^{(n)}_{r,t} = \rho^{(n)}_{r,t}\rho^{(n)}_{p,q}.\]
\item[Containment] If $p<q$ contains $r<t$, then \[\rho^{(n)}_{p,q} \rho^{(n)}_{r,t} = \rho^{(n)}_{\hat{r},\hat{t}} \rho^{(n)}_{p,q},\] where $\hat{r}=\rho^{(n)}_{p,q}(t) = p+q-t$ and $\hat{t}=\rho^{(n)}_{p,q}(r) = p+q-r$.
\end{description}
\end{lemma}

The generators in the cactus groups are modeled after the interval-reversing permutations, with the relations in Lemma \ref{permuation-cactus-axioms}.

\begin{definition}\label{def:cactus-group}
For $n \geq 0$ the $n$th \index{cactus group}\emph{cactus group} $Cac_n$ is the group generated by the generators $s^{(n)}_{p,q}$ for $1 \leq p < q \leq n$ that satisfy the following three relations:
\begin{description}
\item[Involution] Each\index{involution axiom} generator is an involution: \[s^{(n)}_{p,q} s^{(n)}_{p,q} =\id_n,\] where $\id_n$ is the unit in $Cac_n$.
\item[Disjointness] If\index{disjointness axiom} $p < q$ and $r < t$ are disjoint, then \[s^{(n)}_{p,q} s^{(n)}_{r,t} = s^{(n)}_{r,t}s^{(n)}_{p,q}.\]
\item[Containment] If\index{containment axiom} $p<q$ contains $r<t$, then \[s^{(n)}_{p,q} s^{(n)}_{r,t} = s^{(n)}_{\hat{r},\hat{t}} s^{(n)}_{p,q},\] where $\hat{r}=\rho^{(n)}_{p,q}(t)$ and $\hat{t}=\rho^{(n)}_{p,q}(r)$.
\end{description}
Define the group homomorphism \[\pi : Cac_n \to S_n, \qquad \pi\bigl(s^{(n)}_{p,q}\bigr) = \rho^{(n)}_{p,q},\] whose kernel is denoted by\label{not:pcac} $PCac_n$ and is called the $n$th \index{pure cactus group}\emph{pure cactus group}.  Each element in $Cac_n$ (resp., $PCac_n$) is called a\index{cactus}\index{pure cactus} \emph{(pure) cactus}.  We also write $\xbar=\pi(x) \in S_n$ for the image of a cactus $x$ under $\pi$, called the \index{underlying permutation!of a cactus}\emph{underlying permutation} of $x$.  By convention, we set $s^{(n)}_{r,t}=\id_n$ when $r\geq t$.
\end{definition}

\begin{remark}Since the adjacent transpositions \[(i,i+1) = \rho^{(n)}_{i,i+1} = \pi\bigl(s^{(n)}_{i,i+1}\bigr)\] generate the symmetric group $S_n$, the map $\pi : Cac_n\to S_n$ is surjective.\dqed\end{remark}

\begin{remark}[Some History]
In Devadoss \cite{dev} the elements in the group $Cac_n$ are called \index{quasibraid}\emph{quasibraids}.  Henriques and Kamnitzer \cite{hen-kam} introduced the name \emph{cactus group}, which we adopted.  In \cite{dev,gurski,hen-kam} the cactus group was denoted by $J_n$.  Davis, Januszkiewicz, and Scott \cite{djs} used the name \index{mock reflection group}\emph{mock reflection group} for a group with relations similar to those in the cactus group.  The cactus group sits inside an exact sequence of groups \[\nicexy@C-.3cm{\{e\} \ar[r] & PCac_n=\pi_1\bigl(\widetilde{M}_0^{n+1}\bigr) \ar[r] & Cac_n \ar[r]^-{\pi} & S_n \ar[r] & \{e\}}.\] Here $\widetilde{M}_0^{n+1}$ is the total space of the $(n+1)$st tautological line bundle of $\overline{M}_0^{n+1}$, which is the Deligne-Knudson-Mumford moduli space of real, genus $0$, stable curves with $n+1$ marked points.  The group $Cac_n$ was introduced in \cite{dev} in this context.  In \cite{djs} the cactus groups are examples of their mock reflection groups, which are groups of symmetry on certain connected closed manifolds with a cubical cell structure.

On the other hand, in \cite{hen-kam} the cactus group $Cac_n$ was shown to be the group that acts on multiple monoidal products $x_1\otimes\cdots\otimes x_n$ in coboundary monoidal categories in the sense of Drinfel'd \cite{drinfeld}.  So the cactus groups play an important role in (i) quantum group theory, (ii) the geometry of moduli spaces of real, genus $0$, stable curves with marked points, and (iii) symmetries on connected closed manifolds with a cubical cell structure.\dqed
\end{remark}

\begin{example}\label{ex:second-cactus-group}
Here are the first few cactus groups.
\begin{enumerate}
\item The first two cactus groups $Cac_0$ and $Cac_1$ are both the trivial group.  
\item The second cactus group $Cac_2$ has one order-two generator $s^{(2)}_{1,2}$, so \[Cac_2 \cong \bbZ/2.\]  In particular, there are no group homomorphisms $Cac_2\to B_2 \cong\bbZ$ from the second cactus group to the second braid group, which is an infinite cyclic group.\dqed
\end{enumerate}  
\end{example}

\begin{proposition}\label{third-cactus-group}
There is a group isomorphism \[Cac_3 \cong \frac{F\{x,y\}}{(x^2,y^2)}\] with $F\{x,y\}$ the free group on two generators $x$ and $y$.
\end{proposition}

\begin{proof}
By definition the third cactus group $Cac_3$ has generators \[x=s^{(3)}_{1,2},\quad y=s^{(3)}_{1,3},\andspace z=s^{(3)}_{2,3}\] that satisfy the involution relation and the containment relation: \[x^2=y^2=z^2=\id_3,\quad yx=zy,\andspace yz=xy.\]  It follows from the last relation and $y^2=\id_3$ that \[z=y^2z = yxy.\]  But then the relation $yx=zy$ is simply a consequence of $y^2=\id_3$ because \[zy=yxyy=yx.\] Similarly, $z^2=\id_3$ is a consequence of $x^2=y^2=\id_3$ because \[z^2 = yxyyxy = yxxy = yy = \id_3.\]  So the cactus group $Cac_3$ is determined by two non-commuting generators $x$ and $y$, each with order $2$.
\end{proof}

\begin{example}\label{ex:cac-three}
The elements in $Cac_3$ are finite alternating sequences in $x$ and $y$.  The map \[\pi : Cac_3\to S_3\] sends $x$ to the adjacent transposition $(1,2)$ and $y$ to the transposition $(1,3)$.  For instance, the cacti $yxyxyx$ and $xyxyxy$ are both sent to the identity permutation in $S_3$, so belong to the pure cactus group $PCac_3$.\dqed
\end{example}

\section{Direct Sum Cacti}\label{sec:direct-sum-cacti}

The planar operad structure on the cactus groups is defined using two auxiliary operations that are analogous to the direct sum permutation and the block permutation in Definition \ref{def:block-direct-sum-permutation}.  The purpose of this section is to introduce the cactus version of direct sum.  The following notations will be used throughout the rest of this chapter.

\begin{notation}\label{not:cactus}
For integers $n \geq 1$ and $k_i \geq 0$ for $1 \leq i \leq n$, we will write
\begin{enumerate}[label=(\roman*)]
\item $k=k_1+\cdots+k_n$;
\item $k_{[i,j]} = k_i+\cdots+k_j$\label{not:kij} for $1\leq i\leq j\leq n$ (resp., $0$ if $i>j$);
\item $\uk=(k_1,\ldots,k_n)$.  
\end{enumerate}
\end{notation}

To motivate the definition of direct sum cactus, we first consider the interval-reversing permutations.  The following relation can be checked by a direct inspection.

\begin{lemma}\label{int-rev-direct-sum}
For $1 \leq p_i<q_i\leq k_i$ and $1\leq i \leq n$, the relation \[\rho^{(k_1)}_{p_1,q_1} \oplus\cdots\oplus \rho^{(k_n)}_{p_n,q_n} = \prod_{i=1}^n \rho^{(k)}_{k_{[1,i-1]}+p_i,\, k_{[1,i-1]}+q_i}\]
holds in $S_k$.
\end{lemma}

\begin{interpretation}
In Lemma \ref{int-rev-direct-sum}, the left-hand side is a direct sum of interval-reversing permutations.  The equality expresses it as a product of $n$ interval-reversing permutations, which commute with each other by Lemma \ref{permuation-cactus-axioms} because they cover disjoint intervals.  We think of $k$ as an interval that is partitioned into consecutive sub-intervals of lengths $k_1,\ldots,k_n$.  In the direct sum, the $i$th interval-reversing permutation $\rho^{(k_i)}_{p_i,q_i}$ acts on the $i$th sub-interval.  This is the reason for the shifting of the indices by \[k_{[1,i-1]}=k_1+\cdots+k_{i-1}\] on the right-hand side.\dqed
\end{interpretation}

\begin{example}
If $n=3$, then the equality\[\rho^{(k_1)}_{p_1,q_1} \oplus \rho^{(k_2)}_{p_2,q_2} \oplus \rho^{(k_3)}_{p_3,q_3} = \rho^{(k_1+k_2+k_3)}_{p_1,q_1} \cdot \rho^{(k_1+k_2+k_3)}_{k_1+p_2,\, k_1+q_2} \cdot \rho^{(k_1+k_2+k_3)}_{k_1+k_2+p_3,\, k_1+k_2+q_3}\] holds in $S_{k_1+k_2+k_3}$.  Moreover, the three interval-reversing permutations on the right commute with each other by the disjointness relation in Lemma \ref{permuation-cactus-axioms}.\dqed
\end{example}

Using Lemma \ref{int-rev-direct-sum} as a guide, we now define the cactus analogue of direct sum permutations.

\begin{definition}[Direct Sum Cactus]\label{def:cactus-direct-sum}
Define\index{direct sum!cactus} a group homomorphism \[\nicexy{Cac_{k_1}\times \cdots\times Cac_{k_n} \ar[r]^-{\bigoplus_{i=1}^n} & Cac_k}\] by \[\bigoplus_{i=1}^n s^{(k_i)}_{p_i,q_i} = s^{(k_1)}_{p_1,q_1} \oplus\cdots\oplus s^{(k_n)}_{p_n,q_n} = \prod_{i=1}^n s^{(k)}_{k_{[1,i-1]}+p_i,\, k_{[1,i-1]}+q_i} \in Cac_k\]
on generators for $1 \leq p_i<q_i\leq k_i$ and $1\leq i \leq n$, and extend multiplicatively.  We call $\bigoplus_{i=1}^n \sigma_i$ the \emph{direct sum cactus} induced by $\{\sigma_i\} \in Cac_{k_1}\times \cdots\times Cac_{k_n}$.
\end{definition}

\begin{remark}\label{rk:direct-sum-cactus-commute}
The $n$ generators in the big product above commute with each other by the disjointness relation in $Cac_k$.\dqed\end{remark}

\begin{lemma}\label{direct-sum-cactus-welldefined}
The group homomorphism in Definition \ref{def:cactus-direct-sum} is well-defined.
\end{lemma}

\begin{proof}
We need to check that the three relations in Definition \ref{def:cactus-group} for the cactus group are preserved.  

The preservation of the involution relation follows from the computation:
\[\begin{split}
\Bigl(\bigoplus_{i=1}^n s^{(k_i)}_{p_i,q_i}\Bigr)^2
&= \Bigl(\prod_{i=1}^n s^{(k)}_{k_{[1,i-1]}+p_i,\, k_{[1,i-1]}+q_i}\Bigr)^2\\
&= \prod_{i=1}^n \Bigl(s^{(k)}_{k_{[1,i-1]}+p_i,\, k_{[1,i-1]}+q_i}\Bigr)^2\\
&= \prod_{i=1}^n \id_k = \id_k.\end{split}\]
The second and the third equalities follow from the disjointness relation and the involution relation in $Cac_k$, respectively.

Similarly, the preservation of the disjointness (resp., containment) relation is a consequence of the disjointness relation (and the containment relation) in $Cac_k$.
\end{proof}

The direct sum constructions for permutations and for cacti are compatible in the following sense.

\begin{lemma}\label{direct-sum-cactus-permutation}
The diagram \[\nicexy{Cac_{k_1}\times \cdots\times Cac_{k_n} \ar[d]_-{(\pi,\ldots,\pi)} \ar[r]^-{\bigoplus_{i=1}^n} & Cac_k \ar[d]^-{\pi}\\ S_{k_1}\times \cdots \times S_{k_n} \ar[r]^-{\bigoplus_{i=1}^n} & S_k}\] is commutative, in which the bottom horizontal map is the direct sum of permutations.
\end{lemma}

\begin{proof}
All four maps in the diagram are group homomorphisms.  Therefore, it is enough to check the commutativity of the diagram on generators.  Since the map $\pi$ sends each generator $s^{(k_i)}_{p_i,q_i} \in Cac_{k_i}$ to the interval-reversing permutation $\rho^{(k_i)}_{p_i,q_i} \in S_{k_i}$, the commutativity of the diagram follows from Lemma \ref{int-rev-direct-sum} and Definition \ref{def:cactus-direct-sum}.
\end{proof}

The next observation is the associativity property of direct sum cacti.

\begin{lemma}\label{direct-sum-cacti-associativity}
Suppose $l_{i,j}\geq 0$ for $1\leq i \leq n$ and $1\leq j \leq k_i$ with $k_i = l_{i,1}+\cdots+l_{i,k_i}$.  Then the diagram \[\nicexy{\bigl(Cac_{l_{1,1}}\times \cdots\times Cac_{l_{1,k_1}}\bigr) \times\cdots\times \bigl(Cac_{l_{n,1}}\times \cdots\times Cac_{l_{n,k_n}}\bigr) \ar[d]_-{\bigl(\bigoplus_{i=1}^{k_1},\ldots,\bigoplus_{i=1}^{k_n}\bigr)} \ar[r]^-{\bigoplus_{i=1}^k} & Cac_k \ar@{=}[d]\\ 
Cac_{k_1}\times \cdots\times Cac_{k_n} \ar[r]^-{\bigoplus_{i=1}^n} & Cac_k}\]
is commutative.
\end{lemma}

\begin{proof}
The diagram involves only group homomorphisms, so it is enough to check its commutativity on generators, which follows from Definition \ref{def:cactus-direct-sum} and Remark \ref{rk:direct-sum-cactus-commute}.
\end{proof}

\section{Block Cacti}\label{sec:block-cactus}

The purpose of this section is to introduce the cactus version of block permutations, which we will use in the next section to define the action operad structure on the cactus groups.  To motivate its definition, we first consider block permutations induced by interval-reversing permutations in Definition \ref{def:interval-reversal}.  We continue to use Notation \ref{not:cactus}.  Also recall the notation \eqref{block-permutation} for block permutation.  The following relations can be checked by a direct inspection.  

\begin{lemma}\label{int-rev-block}
In the symmetric group $S_k$, the following relations hold.
\begin{enumerate}
\item Suppose $1\leq p<q\leq n$.  Then 
\[\rho^{(n)}_{p,q} \langle \uk\rangle = \Bigl(\prod_{j=p}^q \rho^{(k)}_{k_{[1,p-1]}+k_{[j+1,q]}+1,\, k_{[1,p-1]}+k_{[j,q]}}\Bigr) \cdot \rho^{(k)}_{k_{[1,p-1]}+1,\,k_{[1,q]}}\]
holds, in which $\rho^{(n)}_{p,q}\langle \uk\rangle$ is the block permutation induced by the interval-reversing permutation $\rho^{(n)}_{p,q} \in S_n$ that permutes $n$ consecutive blocks of lengths $k_1,\ldots,k_n$.
\item Suppose $\sigma,\tau\in S_k$.  Then \[(\sigma\tau)\langle \uk\rangle = \sigma \langle\tau\uk \rangle \cdot \tau\langle \uk \rangle\] with $\tau\uk = (k_{\tau^{-1}(1)}, \ldots, k_{\tau^{-1}(n)})$.
\end{enumerate}
\end{lemma}

\begin{interpretation}
In Lemma \ref{int-rev-block}(1), the block permutation $\rho^{(n)}_{p,q} \langle \uk\rangle$ 
\begin{enumerate}[label=(\roman*)]
\item reverses the order of the $p$th block up to and including the $q$th block and
\item keeps the relative order within each of the $q-p+1$ blocks unchanged.  
\end{enumerate}
The equality expresses this block permutation in terms of interval-reversing permutations. First the right-most interval-reversing permutation \[A=\rho^{(k)}_{k_{[1,p-1]}+1,\,k_{[1,q]}} \in S_k\] reverses the interval $[k_{[1,p-1]}+1,k_{[1,q]}]$, which covers the $p$th block to the $q$th block.  Then each interval-reversing permutation \[\rho^{(k)}_{k_{[1,p-1]}+k_{[j+1,q]}+1,\, k_{[1,p-1]}+k_{[j,q]}}\] in the big product reverses the order within the interval that corresponds to the $j$th block before the permutation $A$ was applied.  Notice that the $q-p+1$ interval-reversing permutations in the big product commute with each other by Lemma \ref{permuation-cactus-axioms} because they cover disjoint intervals. 

Lemma \ref{int-rev-block}(2) says that the block permutation induced by a product is equal to the product of two block permutations.  Note that taking the block permutation is not multiplicative.\dqed 
\end{interpretation}

Using Lemma \ref{int-rev-block} as a guide, we now define the cactus version of block permutations.

\begin{definition}[Block Cactus]\label{def:block-cactus}
Define\index{block cactus} a function \[\nicexy@C+.2cm{Cac_n \ar[r]^-{-\langle \uk \rangle} & Cac_k}\] inductively on the number of generators as follows.
\begin{enumerate}[label=(\roman*)]
\item Set $\id_n\langle \uk\rangle = \id_k$.
\item For $1 \leq p<q\leq n$, the generator $s^{(n)}_{p,q}\in Cac_n$ is sent to the product
\[s^{(n)}_{p,q} \langle \uk\rangle = \Bigl(\prod_{j=p}^q s^{(k)}_{k_{[1,p-1]}+k_{[j+1,q]}+1,\, k_{[1,p-1]}+k_{[j,q]}}\Bigr) \cdot s^{(k)}_{k_{[1,p-1]}+1,\,k_{[1,q]}} \in Cac_k.\]
\item Suppose $\sigma,\tau\in Cac_n$ such that $\sigma\langle \cdots \rangle$ and $\tau\langle \cdots \rangle$ have been defined.  Then $(\sigma\tau)\langle \uk\rangle$ is defined as the product
\begin{equation}\label{block-cactus-product}
(\sigma\tau)\langle \uk\rangle = \sigma \langle\tau\uk \rangle \cdot \tau\langle \uk \rangle \in Cac_k
\end{equation}
with \[\tau\uk = (k_{\taubar^{-1}(1)}, \ldots, k_{\taubar^{-1}(n)}),\] where $\taubar=\pi(\tau)\in S_n$ is the underlying permutation of $\tau$.
\end{enumerate}
We call $\sigma\langle \uk \rangle$ a \emph{block cactus} induced by $\sigma$.
\end{definition}

\begin{remark}In the definition of the block cactus $s^{(n)}_{p,q} \langle \uk\rangle$, the $q-p+1$ generators in the product $\prod_{j=p}^q$ commute with each other by the disjointness relation in $Cac_k$.\dqed\end{remark}

Below we will show that \eqref{block-cactus-product} is indeed well-defined.  The following result contains some preliminary calculation involving $s^{(n)}_{p,q} \langle \uk\rangle$.

\begin{lemma}\label{block-cactus-properties}
For $1 \leq p<q\leq n$, the following statements hold in $Cac_k$.
\begin{enumerate}
\item There are equalities 
\[\begin{split} &s^{(n)}_{p,q} \langle \uk\rangle
= s^{(k)}_{k_{[1,p-1]}+1,\,k_{[1,q]}} \cdot \Bigl(\prod_{j=p}^q s^{(k)}_{k_{[1,j-1]}+1,\, k_{[1,j]}}\Bigr)\\
&= s^{(k)}_{k_{[1,p-1]}+1,\,k_{[1,q]}} \cdot \Bigl(\id_{k_1}\oplus\cdots\oplus \id_{k_{p-1}} \oplus s^{(k_p)}_{1,k_p}\oplus\cdots\oplus s^{(k_q)}_{1,k_q} \oplus \id_{k_{q+1}} \oplus\cdots\oplus \id_{k_n}\Bigr).
\end{split}\]
\item The equality \[s^{(k)}_{k_{[1,p-1]}+1,\,k_{[1,q]}} = s^{(k)}_{k_{[1,p-1]}+1,\,k_{[1,p-1]} + k_{\xi(p)}+\cdots+k_{\xi(q)}}\] holds for each permutation $\xi$ of $\{p,\ldots,q\}$.
\item There is an equality 
\[\begin{split}
& s^{(n)}_{p,q}\langle k_1,\ldots,k_{p-1}, k_q,k_{q-1},\ldots,k_p,k_{q+1},\ldots,k_n\rangle\\
&= \Bigl(\prod_{j=p}^q s^{(k)}_{k_{[1,j-1]}+1,\, k_{[1,j]}}\Bigr) \cdot s^{(k)}_{k_{[1,p-1]}+1,\ k_{[1,q]}},\end{split}\] in which the $q-p+1$ generators in the product $\prod_{j=p}^q$ commute with each other. 
\item For $p\leq j \leq q$, there is an equality
\[\begin{split}
& s^{(k)}_{k_{[1,p-1]}+1,\,k_{[1,q]}} \cdot s^{(k)}_{k_{[1,p-1]}+k_{[j+1,q]}+1,\, k_{[1,p-1]}+k_{[j,q]}}\\
&= s^{(k)}_{k_{[1,j-1]}+1,\,k_{[1,j]}} \cdot s^{(k)}_{k_{[1,p-1]}+1,\,k_{[1,q]}}. \end{split}\]
\end{enumerate}
\end{lemma}

\begin{proof}
The first equality in (1) is a result of repeated applications of the containment relation in $Cac_k$ to bring the generator $s^{(k)}_{k_{[1,p-1]}+1,\,k_{[1,q]}}$ pass each of the $q-p+1$ generators in the product $\prod_{j=p}^q$.  The second equality in (1) follows from the definition of the direct sum cactus in Definition \ref{def:cactus-direct-sum}.

The equality in (2) holds because \[k_{[1,q]} = k_{[1,p-1]}+k_{[p,q]}= k_{[1,p-1]}+k_{\xi(p)}+\cdots+k_{\xi(q)}.\] 

The equality in (3) is the definition of $s^{(n)}_{p,q} \langle \uk\rangle$ with the indices $\{p,\ldots,q\}$ first reversed to $\{q,q-1,\ldots,p\}$.  The commutativity assertion holds by the disjointness relation in $Cac_k$.

The equality in (4) is true by the containment relation in $Cac_k$.
\end{proof}

\begin{lemma}\label{block-cactus-well-defined}
The function in Definition \ref{def:block-cactus} is well-defined.
\end{lemma}

\begin{proof}
We need to check that the three relations in Definition \ref{def:cactus-group} for the cactus group are preserved.  

For the involution relation, by definition we have \[\bigl(s^{(n)}_{p,q} s^{(n)}_{p,q}\bigr)\langle \uk\rangle = \underbrace{s^{(n)}_{p,q} \langle k_{\taubar^{-1}(1)}, \ldots, k_{\taubar^{-1}(n)}\rangle}_{X} \cdot \underbrace{s^{(n)}_{p,q} \langle \uk \rangle}_{Y}\] with $\taubar = \rho^{(n)}_{p,q} \in S_n$ and $X$ the block cactus in Lemma \ref{block-cactus-properties}(3).  Each of the block cacti $X$ and $Y$ on the right is a product of $q-p+2$ generators in $Cac_k$.  First we use the containment relation to bring the right-most generator in $X$ pass the first $q-p+1$ generators in $Y$.  In the resulting product, the two right-most generators are equal, so their product is $\id_k$ by the involution relation in $Cac_k$.  The remaining product of $2(q-p+1)$ generators is also equal to $\id_k$ by repeated applications of Lemma \ref{block-cactus-properties}(4), the disjointness relation, and the involution relation in $Cac_k$.  So \[\bigl(s^{(n)}_{p,q} s^{(n)}_{p,q}\bigr)\langle \uk\rangle = \id_k = \id_n\langle \uk\rangle.\]

For the disjointness relation, suppose $p<q$ and $r<t$ are disjoint.  Then in \[\bigl(s^{(n)}_{p,q} s^{(n)}_{r,t}\bigr)\langle \uk\rangle = \underbrace{s^{(n)}_{p,q}\langle k_{\chibar^{-1}(1)}, \ldots, k_{\chibar^{-1}(n)}\rangle}_{U} \cdot \underbrace{s^{(n)}_{r,t} \langle \uk \rangle}_{V}\] with $\chibar = \rho^{(n)}_{r,t}\in S_n$, observe that there is an equality\[U=s^{(n)}_{p,q}\langle \uk\rangle\] because $\chibar$ only reverses the interval $[r,t]$, which is disjoint from $[p,q]$.  Moreover, each of the $q-p+2$ generators in $U$ commutes with each of the $q-p+2$ generators in $V$ by the disjointness relation in $Cac_k$.  Therefore, the result is the same if we start with \[\bigl(s^{(n)}_{r,t} s^{(n)}_{p,q} \bigr)\langle\uk\rangle,\] which is therefore equal to $\bigl(s^{(n)}_{p,q} s^{(n)}_{r,t}\bigr)\langle \uk\rangle$.

For the containment relation, suppose $p<q$ contains $r<t$.  In \[\bigl(s^{(n)}_{p,q} s^{(n)}_{r,t}\bigr)\langle \uk\rangle = \underbrace{s^{(n)}_{p,q}\langle k_{\chibar^{-1}(1)}, \ldots, k_{\chibar^{-1}(n)}\rangle}_{Q} \cdot \underbrace{s^{(n)}_{r,t} \langle \uk \rangle}_{R}\] with $\chibar = \rho^{(n)}_{r,t}\in S_n$, observe that the right-most generator in $Q$ is equal to the right-most generator in $s^{(n)}_{p,q}\langle \uk\rangle$.  We first bring this generator in $Q$ pass each of the $q-p+2$ generators in $R$ using the containment relation in $Cac_k$.  Then we bring each of the $q-p+1$ remaining generators in $Q$ pass each of the generators that came from $R$ using the relations in $Cac_k$.  After these steps, the result is \[s^{(n)}_{\hat{r},\hat{t}}\langle k_{\taubar^{-1}(1)},\ldots, k_{\taubar^{-1}(n)} \rangle \cdot s^{(n)}_{p,q}\langle \uk\rangle = \bigl(s^{(n)}_{\hat{r},\hat{t}} s^{(n)}_{p,q} \bigr)\langle \uk\rangle\] with $\taubar = \rho^{(n)}_{p,q} \in S_n$, $\hat{r}=\taubar(t)$, and $\hat{t}=\taubar(r)$.    
\end{proof}

\begin{example}
For $1 \leq 2 <4\leq 5$, there is the block cactus
\[\begin{split}
& s^{(5)}_{2,4}\langle k_1,\ldots,k_5\rangle\\
&= s^{(k)}_{k_1+k_4+k_3+1,\, k_1+k_4+k_3+k_2}
\cdot s^{(k)}_{k_1+k_4+1,\, k_1+k_4+k_3} \cdot s^{(k)}_{k_1+1,\,k_1+k_4}
\cdot s^{(k)}_{k_1+1,\,k_1+k_2+k_3+k_4}\end{split}\]
in $Cac_{k_1+\cdots+k_5}$ for $k_i \geq 0$.\dqed
\end{example}

The next several observations will be used in the following section to show that the cactus groups form an action operad.  First, block permutations and block cacti are compatible in the following sense.

\begin{lemma}\label{block-perm-cacti-compatible}
The diagram \[\nicexy{Cac_n \ar[d]_-{\pi} \ar[r]^-{-\langle\uk\rangle} & Cac_k \ar[d]^-{\pi}\\ S_n \ar[r]^-{-\langle\uk\rangle} & S_k}\] is commutative.
\end{lemma}

\begin{proof}
For each generator in $Cac_n$, the above diagram is commutative by Lemma \ref{int-rev-block}(1) and that $\pi$ is a group homomorphism.  The general case--for a finite product of generators in $Cac_n$--follows by an induction on the number of generators in $Cac_n$, Lemma \ref{int-rev-block}(2), and \eqref{block-cactus-product}.
\end{proof}

The next observation is the multiplicative commutativity between block cacti and direct sum cacti.  

\begin{lemma}\label{direct-sum-block-cacti-commute}
For $\sigma\in Cac_n$ and $\tau_i\in Cac_{k_i}$ with $1 \leq i \leq n$, the equality
\[\sigma\langle\uk\rangle \cdot \bigl(\tau_1 \oplus\cdots\oplus \tau_n\bigr) = \bigl(\tau_{\sigmabar^{-1}(1)} \oplus\cdots\oplus \tau_{\sigmabar^{-1}(n)}\bigr) \cdot \sigma\langle\uk\rangle\] holds in $Cac_k$, where $\sigmabar=\pi(\sigma)\in S_n$ is the underlying permutation of $\sigma$.
\end{lemma}

\begin{proof}
First we consider the case when $\sigma$ and the $\tau_i$'s are all generators.  By the definition of the direct sum cactus, it is enough to prove the case when $\tau_j = \id_{k_j}$ for all $1 \leq j \leq n$ except for one index $i$.  In this case, with both $\sigma = s^{(n)}_{p,q}$ and $\tau_i = s^{(k_i)}_{p_i,q_i}$ generators, the desired equality becomes 
\begin{equation}\label{sum-block-cacti-commute-1}
s^{(n)}_{p,q}\langle\uk\rangle \cdot s^{(k)}_{k_{[1,i-1]}+p_i,\, k_{[1,i-1]}+q_i}
= s^{(k)}_{N+p_i,\,N+q_i} \cdot s^{(n)}_{p,q}\langle\uk\rangle,
\end{equation}
where \[N=k_{\sigmabar^{-1}(1)} +\cdots+ k_{\sigmabar^{-1}(\sigmabar(i)-1)}\] with $\sigmabar = \rho^{(n)}_{p,q}$.  If either $i<p$ or $q<i$, then \eqref{sum-block-cacti-commute-1} holds by the disjointness relation in $Cac_k$ and the definition of the block cactus $s^{(n)}_{p,q}\langle\uk\rangle$.  If $p\leq i \leq q$, then \eqref{sum-block-cacti-commute-1} holds by the disjointness relation and the containment relation in $Cac_k$.  

The general case--when both $\sigma$ and $\tau_i$ are finite products of generators--follows by an induction, using the multiplicativity of direct sum cacti and \eqref{block-cactus-product}.  In more details, suppose the assertion is true for $\sigma,\rho\in Cac_n$ and $\tau_i,\tau_i'\in Cac_{k_i}$ for $1\leq i \leq n$.  Then we compute as follows:
\[\begin{split} (\sigma\rho)\langle\uk\rangle \cdot \bigoplus_{i=1}^n (\tau_i\tau_i')
&= \sigma\langle\rho\uk\rangle \cdot \rho\langle\uk\rangle \cdot \Bigl(\bigoplus_{i=1}^n \tau_i\Bigr) \cdot \Bigl(\bigoplus_{i=1}^n \tau_i'\Bigr)\\
&= \sigma\langle\rho\uk\rangle \cdot \Bigl(\bigoplus_{i=1}^n \tau_{\rhobar^{-1}(i)}\Bigr) \cdot \rho\langle\uk\rangle \cdot \Bigl(\bigoplus_{i=1}^n \tau_i'\Bigr)\\
&= \Bigl(\bigoplus_{i=1}^n \tau_{\overline{\sigma\rho}^{-1}(i)}\Bigr) \cdot \sigma\langle\rho\uk\rangle \cdot \Bigl(\bigoplus_{i=1}^n \tau_{\rhobar^{-1}(i)}'\Bigr) \cdot \rho\langle\uk\rangle\\
&= \Bigl(\bigoplus_{i=1}^n \tau_{\overline{\sigma\rho}^{-1}(i)}\Bigr) \cdot \Bigl(\bigoplus_{i=1}^n \tau_{\overline{\sigma\rho}^{-1}(i)}'\Bigr) \cdot \sigma\langle\rho\uk\rangle \cdot \rho\langle\uk\rangle\\
&= \Bigl(\bigoplus_{i=1}^n \tau_{\overline{\sigma\rho}^{-1}(i)} \tau_{\overline{\sigma\rho}^{-1}(i)}'\Bigr) \cdot (\sigma\rho)\langle\uk\rangle.
\end{split}\]
This shows that the desired equality holds for $\sigma\rho$ and $\{\tau_i\tau_i'\}_{i=1}^n$.
\end{proof}

The next result is the compositional commutativity between block cacti and direct sum cacti.  

\begin{lemma}\label{direct-sum-block-cacti-commute2}
Suppose $l_i, m_{i,j} \geq 0$ for $1 \leq i \leq n$ and $1 \leq j \leq l_i.$  Define
\begin{itemize}
\item $\um_i = (m_{i,1},\ldots,m_{i,l_i})$ for $1\leq i \leq n$,
\item $\um = (\um_1,\ldots,\um_n)$,
\item $k_i = m_{i,1}+\cdots+m_{i,l_i}$ for $1 \leq i \leq n$,
\item $k=k_1+\cdots+k_n$, and $l=l_1+\cdots+l_n$.
\end{itemize}
Then the diagram\[\nicexy@C+1cm{Cac_{l_1}\times\cdots\times Cac_{l_n} \ar[d]_-{\bigl(-\langle\um_1\rangle,\ldots, -\langle \um_n\rangle\bigr)} \ar[r]^-{\bigoplus_{i=1}^n} & Cac_l \ar[d]^-{\langle \um\rangle}\\ Cac_{k_1}\times\cdots\times Cac_{k_n} \ar[r]^-{\bigoplus_{i=1}^n} & Cac_k}\]
is commutative.
\end{lemma}

\begin{proof}
By the multiplicativity of direct sum cacti and \eqref{block-cactus-product}, it is enough to check the commutativity of the above diagram on an element of the form \[\bigl(\id_{l_1},\ldots, \id_{l_{i-1}}, \sigma, \id_{l_{i+1}} \ldots,\id_{l_n}\bigr) \in Cac_{l_1}\times\cdots\times Cac_{l_n} \] with $\sigma \in Cac_{l_i}$ for some $1 \leq i \leq n$.  Furthermore, by an induction, we reduce to the case when $\sigma$ is a generator in $Cac_{l_i}$.  In this case, the desired equality follows from a direct inspection of the definitions of the direct sum cactus and the block cactus.
\end{proof}

The following result is the associativity property of block cacti.  It says that a block cactus of a block cactus is a block cactus.

\begin{lemma}\label{block-cacti-associativity}
Suppose $l_{i,j}\geq 0$ for $1\leq i \leq n$ and $1\leq j \leq k_i$. Define
\begin{itemize}
\item$\ul_i = (l_{i,1},\ldots,l_{i,k_i})$ for $1 \leq i \leq n$,
\item $l_i = l_{i,1}+\cdots+l_{i,k_i}$ for $1 \leq i \leq n$, and
\item $l=l_1+\cdots+l_n$.
\end{itemize}
Then the diagram \[\nicexy@C+1cm{Cac_n \ar[r]^-{\langle l_1,\ldots,l_n\rangle} \ar[d]_-{\langle k_1,\ldots,k_n\rangle} & Cac_l \ar@{=}[d]\\ Cac_k \ar[r]^-{\langle \ul_1,\ldots,\ul_n\rangle} & Cac_l}\] is commutative.
\end{lemma}

\begin{proof}
For each generator in $Cac_n$, the commutativity of the above diagram follows from the definition of the block cactus for a product \eqref{block-cactus-product} and the involution and disjointness relations in $Cac_l$.  In more details, suppose given a generator $s^{(n)}_{p,q}$ in $Cac_n$.  Writing $\ul=(\ul_1,\ldots,\ul_n)$ and using \eqref{block-cactus-product} and Lemma \ref{block-cactus-properties}(1) repeatedly, there are equalities
\[\begin{split} 
\bigl(s^{(n)}_{p,q}\langle\uk\rangle\bigr) \langle\ul\rangle 
&= \Bigl[s^{(k)}_{k_{[1,p-1]}+1,\, k_{[1,q]}} \cdot \prod_{j=p}^q s^{(k)}_{k_{[1,j-1]}+1,\, k_{[1,j]}} \Bigr] \langle\ul\rangle\\
&= \Bigl[s^{(k)}_{k_{[1,p-1]}+1,\, k_{[1,q]}} \langle \cdots\rangle\Bigr] \cdot \prod_{j=p}^q \Bigl[s^{(k)}_{k_{[1,j-1]}+1,\, k_{[1,j]}}\langle \cdots\rangle\Bigr]\\
&= \Bigl[\underbrace{s^{(l)}_{l_{[1,p-1]}+1,\, l_{[1,q]}}}_{Y} \cdot Z\Bigr] \cdot \prod_{j=p}^q \Bigl[\underbrace{s^{(l)}_{l_{[1,j-1]}+1,\, l_{[1,j]}}}_{Y_j} \cdot X_j\Bigr]
\end{split}\]
in $Cac_l$ in which
\[\begin{cases}
X_j = \prod_{m=1}^{k_j} s^{(l)}_{l_{[1,j-1]} + l_{j,1}+\cdots+ l_{j,m-1}+1,\, l_{[1,j-1]}+l_{j,1} +\cdots+ l_{j,m}},\\ Z = \prod_{j=p}^q Z_j \withspace Z_j = Y_j X_j Y_j.
\end{cases}\]
By the disjointness relation in $Cac_l$:
\begin{itemize}\item $Z_j$ commutes with $Y_i$ and $X_i$ if $p\leq i\not= j\leq q$.
\item The $k_j$ generators in the product $\prod_{m=1}^{k_j}$ in $X_j$ commute with each other.
\end{itemize}
Therefore, by repeated applications of the involution relation in $Cac_l$, we conclude that
\[\bigl(s^{(n)}_{p,q}\langle\uk\rangle\bigr) \langle\ul\rangle 
= Y\cdot \prod_{j=p}^q Y_j= s^{(n)}_{p,q}\langle\ul\rangle.\]
So the diagram is commutative for each generator in $Cac_n$.

The general case--for a finite product of generators in $Cac_n$--follows by an induction.  Indeed, suppose the diagram is commutative for $\sigma,\tau\in Cac_n$.  Then there are equalities
\[\begin{split} \Bigl[(\sigma\tau)\langle\uk\rangle\Bigr]\langle\ul\rangle
&= \Bigl[\sigma\langle\tau\uk\rangle \cdot \tau\langle\uk\rangle\Bigr] \langle\ul\rangle\\
&= \Bigl[\sigma\langle\tau\uk\rangle\Bigr]\langle \ul_{\taubar^{-1}(1)},\ldots,\ul_{\taubar^{-1}(n)}\rangle \cdot \Bigl[\tau\langle\uk\rangle\Bigr]\langle\ul\rangle \\
&= \sigma\langle l_{\taubar^{-1}(1)},\ldots,l_{\taubar^{-1}(n)}\rangle \cdot \tau\langle l_1,\ldots,l_n\rangle\\
&= (\sigma\tau)\langle l_1,\ldots,l_n\rangle.
\end{split}\]
The third equality follows from the assumption about $\sigma$ and $\tau$.  The other three equalities are from \eqref{block-cactus-product}.
\end{proof}

\section{Cactus Group Operad}\label{sec:cactus-action-operad}

In this section, we define the action operad structure on the cactus groups due to Gurski \cite{gurski}.  Then we define cactus operads as monoids in the monoidal category of cactus sequences equipped with the cactus circle product.  

\begin{motivation}
To motivate the operadic composition involving the cactus groups, recall from Proposition \ref{sn-operad} that in the symmetric group operad $\S=\{S_n\}_{n\geq 0}$, the operadic composition is given by a product \[\gamma^{\S}\bigl(\sigma;\tau_1,\ldots,\tau_n\bigr) = \sigma \langle k_1,\ldots,k_n\rangle \cdot \bigl(\tau_1\oplus\cdots\oplus\tau_n\bigr)\] of a block permutation induced by $\sigma$ and a direct sum permutation induced by the $\tau_i$'s.  The operadic composition in the braid group operad $\B$ in Definition \ref{def:braid-group-operad} and in the ribbon group operad $\R$ in Definition \ref{def:ribbon-group-operad} are defined in the same way using block braids and direct sum braids (resp., block ribbons and direct sum ribbons).  The operadic composition involving the cactus groups is defined analogously in terms of a block cactus and a direct sum cactus.\dqed
\end{motivation}

\begin{definition}[Operadic Composition of Cacti]\label{def:cactus-operad-comp}
Define a function \[\nicexy{Cac_n\times Cac_{k_1}\times\cdots\times Cac_{k_n}\ar[r]^-{\gammacac} & Cac_{k}}\] by \[\gammacac\bigl(\sigma;\tau_1,\ldots,\tau_n\bigr) = \sigma \langle k_1,\ldots,k_n\rangle \cdot \bigl(\tau_1\oplus\cdots\oplus\tau_n\bigr)\] for $\sigma \in Cac_n$ and $\tau_i\in Cac_{k_i}$, with:
\begin{itemize}
\item $\sigma\langle k_1,\ldots,k_n\rangle$ the block cactus in Definition \ref{def:block-cactus};
\item $\tau_1\oplus\cdots\oplus\tau_n$ the direct sum cactus in Definition \ref{def:cactus-direct-sum}.
\end{itemize}
\end{definition}

\begin{theorem}\label{thm:cactus-group-operad}
The sequence\index{cactus group operad}\index{group operad!cactus} \[\Cac = \bigl\{Cac_n\bigr\}_{n\geq 0}\] of cactus groups forms an action operad with the operadic composition $\gammacac$ in Definition \ref{def:cactus-operad-comp}.
\end{theorem}

\begin{proof}
To see that $\Cac$ is a one-colored planar operad in $\Set$ as in Definition \ref{def:planar-operad-compi}, first note that the operadic unity axioms hold by the definition of block cactus and direct sum cactus.  

Next we consider the operadic associativity axiom, which is the equality
\begin{equation}\label{cactus-operadic-associativity}
\gamma\Bigl(\sigma; \gamma(\tau_1;\uchi_1),\ldots,\gamma(\tau_n;\uchi_n)\Bigr) = 
\gamma\Bigl(\gamma(\sigma;\utau); \uchi_1,\ldots,\uchi_n\Bigr)
\end{equation}
in which
\begin{itemize}
\item $\gamma=\gammacac$, $\sigma\in Cac_n$, 
\item $\utau = (\tau_1,\ldots,\tau_n) \in Cac_{k_1}\times\cdots\times Cac_{k_n}$, and 
\item $\uchi_i = \bigl(\chi_{i,1},\ldots,\chi_{i,k_i}\bigr)\in Cac_{l_{i,1}}\times\cdots\times Cac_{l_{i,k_i}}$ for $1\leq i \leq n$.
\end{itemize}
After writing out both sides of \eqref{cactus-operadic-associativity} in terms of block cacti and direct sum cacti as in Definition \ref{def:cactus-operad-comp}, one checks that the two sides are equal using \eqref{block-cactus-product}, Lemma \ref{direct-sum-cacti-associativity}, Lemma \ref{direct-sum-block-cacti-commute2}, and Lemma \ref{block-cacti-associativity}.  In more details, suppose 
\begin{itemize}
\item $l_i=l_{i,1}+\cdots+l_{i,k_i}$, 
\item $\ul_i=(l_{i,1},\ldots,l_{i,k_i})$, 
\item $\ul=(\ul_1,\ldots,\ul_n)$, and 
\item $\tau_i\ul_i = \bigl(l_{i,\taubar_i^{-1}(1)}, \ldots, l_{i,\taubar_i^{-1}(k_i)}\bigr)$ with $\taubar_i=\pi(\tau_i)\in S_{k_i}$ the underlying permutation of $\tau_i$.
\end{itemize}
Using the results cited above, the operadic associativity \eqref{cactus-operadic-associativity} follows from the following computation in $Cac_{l_1+\cdots+l_n}$:
\[\begin{split} 
\gamma\Bigl(\sigma; \bigl\{\gamma(\tau_i;\uchi_i)\bigr\}_{i=1}^n\Bigr)
&= \sigma \langle l_1,\ldots,l_n\rangle \cdot \bigoplus_{i=1}^n \Bigl(\tau_i\langle\ul_i\rangle \cdot \bigoplus_{j=1}^{k_i} \chi_{i,j}\Bigr)\\
&= \sigma \langle l_1,\ldots,l_n\rangle \cdot \Bigl(\bigoplus_{i=1}^n \tau_i\langle\ul_i\rangle\Bigr) \cdot \Bigl(\bigoplus_{i,\,j} \chi_{i,j}\Bigr)\\
&= \bigl(\sigma\langle \uk\rangle\bigr) \langle \tau_1\ul_1,\ldots,\tau_n\ul_n\rangle \cdot \Bigl(\bigoplus_{i=1}^n \tau_i\Bigr)\langle \ul\rangle \cdot \Bigl(\bigoplus_{i,\,j} \chi_{i,j}\Bigr)\\
&= \Bigl(\sigma\langle \uk\rangle \cdot \bigoplus_{i=1}^n \tau_i\Bigr)\langle \ul\rangle \cdot \Bigl(\bigoplus_{i,\,j} \chi_{i,j}\Bigr)\\
&= \gamma\Bigl(\gamma(\sigma;\utau); \uchi_1,\ldots,\uchi_n\Bigr).
\end{split}\]

Next, Lemma \ref{direct-sum-cactus-permutation} and Lemma \ref{block-perm-cacti-compatible} ensure that $\pi : \Cac \to \S$ is a morphism of one-colored planar operads in $\Set$.

Finally, the action operad axiom \eqref{augmented-goperad-identity} follows from \eqref{block-cactus-product}, the multiplicativity of direct sum cacti, and Lemma \ref{direct-sum-block-cacti-commute}.  In more details, with the notations in \eqref{augmented-goperad-identity}, suppose \[\rho\uk=\bigl(k_{\rhobar^{-1}(1)},\ldots, k_{\rhobar^{-1}(n)}\bigr).\] Then the action operad axiom \eqref{augmented-goperad-identity} for the cactus groups follows from the following computation in $Cac_k$:
\[\begin{split}
\gamma\Bigl(\sigma\rho; \{\tau_i\tau_i'\}_{i=1}^n\Bigr) 
&= (\sigma\rho)\langle\uk\rangle \cdot \Bigl(\bigoplus_{i=1}^n \tau_i\tau_i'\Bigr)\\
&= \sigma\langle \rho\uk\rangle \cdot \rho\langle\uk\rangle \cdot \Bigl(\bigoplus_{i=1}^n \tau_i\Bigr) \cdot \Bigl(\bigoplus_{i=1}^n \tau_i'\Bigr)\\
&= \sigma\langle \rho\uk\rangle \cdot \Bigl(\bigoplus_{i=1}^n \tau_{\rhobar^{-1}(i)}\Bigr) \cdot \rho\langle \uk\rangle \cdot \Bigl(\bigoplus_{i=1}^n \tau_i'\Bigr)\\
&= \gamma\Bigl(\sigma; \{\tau_{\rhobar^{-1}(i)}\}_{i=1}^n\Bigr) \cdot \gamma\Bigl(\rho; \{\tau_i'\}_{i=1}^n\Bigr).
\end{split}\]
\end{proof}

\begin{definition}\label{def:cactus-group-operad}
The action operad $\Cac$ in Theorem \ref{thm:cactus-group-operad} is called the \emph{cactus group operad}.
\end{definition}

There is also an action operad that involves the pure cactus groups.

\begin{proposition}\label{pure-cactus-action-operad}
The sequence \[\PCac = \bigl\{PCac_n\bigr\}_{n\geq 0}\] of pure cactus groups forms an action operad with the operadic composition $\gammapcac$ defined as in Definition \ref{def:cactus-operad-comp}.
\end{proposition}

\begin{proof}
It is enough to show that the block cactus $\sigma\langle k_1,\ldots,k_n\rangle$ and the direct sum cactus $\tau_1\oplus\cdots\oplus\tau_n$ are both pure cacti when $\sigma$ and the $\tau_i$'s are pure cacti.  Recall that a pure cacti is an element in $Cac_n$ for some $n$ whose underlying permutation is the identity permutation.  If $\sigma$ is a pure cactus, then so is the block cactus $\sigma\langle k_1,\ldots,k_n\rangle$  by Lemma \ref{block-perm-cacti-compatible}.  Similarly, if the $\tau_i$'s are pure cacti, then the direct sum cactus $\tau_1\oplus\cdots\oplus\tau_n$ is a pure cactus by Lemma \ref{direct-sum-cactus-permutation}.
\end{proof}

\begin{definition}\label{def:pure-cactus-group-operad}
The action operad $\PCac$ in Proposition \ref{pure-cactus-action-operad} is called the \index{pure cactus group operad}\emph{pure cactus group operad}.
\end{definition}

We now restrict Section \ref{sec:group-operad} to the action operads $\Cac$ and $\PCac$.  Following the planar, symmetric, braided, and ribbon cases, for the action operad $\Cac$ (resp., $\PCac$), we will use the adjective \emph{cactus} (resp., \emph{pure cactus}).

\begin{definition}\label{def:cactus-operad}
Suppose $\colorc$ is a set.  With the action operad $\G=\Cac$:
\begin{enumerate}
\item The objects of the diagram category\label{notation:cacseqcm} \[\cacseqcm = \M^{\cacsubcopc}\] are called \index{cactus sequence}\emph{$\colorc$-colored cactus sequences} in $\M$.  
\item For $X,Y  \in \cacseqcm$, $X \circcac Y \in \cacseqcm$ is called the \index{circle product!cactus}\index{cactus circle product}\emph{$\colorc$-colored cactus circle product}.  It is defined entrywise as the coend\label{not:circcac} \[(X \circcac Y)\dub = \int^{\uc \in \cacsubc} X\duc \otimes Y^{\uc}_{\smallsubcac}(\ub)\in \M\] for $\dub \in \cacsubcopc$.  The object $Y^{\uc}_{\smallsubcac} \in \M^{\cacsubcop}$ is defined entrywise as the coend \[Y^{\uc}_{\smallsubcac}(\ub) = \int^{\{\ua_j\}\in\prod_{j=1}^m \cacsubcop} \cacsubcop\bigl(\ua_1,\ldots,\ua_m;\ub\bigr) \cdot \left[\bigotimes_{j=1}^m Y\cjuaj\right] \in \M\] for $\ub \in \cacsubcop$, in which $(\ua_1,\ldots,\ua_m) \in \cacsubcop$ is the concatenation.  
\item The category of \index{operad!cactus}\index{cactus operad}\emph{$\colorc$-colored cactus operads in $\M$} is the category\label{not:cacoperadcm} \[\Cacoperadcm = \Mon\bigl(\cacseqcm, \circcac, \I^{\Cac}\bigr)\] of monoids in the monoidal category $\bigl(\cacseqcm,\circcac,\I^{\Cac}\bigr)$.  If $\colorc$ has $n<\infty$ elements, we also refer to objects in $\Cacoperadcm$ as \emph{$n$-colored cactus operads in $\M$}.
\item Suppose $\O$ is a $\colorc$-colored cactus operad in $\M$.  The category $\algmo$ \index{cactus operad!algebra}\index{algebra!of a cactus operad}of \emph{$\O$-algebras} is defined as the category of $(\O\circcac -)$-algebras for the monad $\O\circcac -$ in $\Mtoc$.
\end{enumerate}
\end{definition}

The previous definition also works with the pure cactus group operad $\PCac$, which defines $\colorc$-colored \index{pure cactus sequence}pure cactus sequences, the $\colorc$-colored \index{pure cactus circle product}\index{circle product!pure cactus}pure cactus circle product, $\colorc$-colored \index{pure cactus operad}\index{operad!pure cactus}pure cactus operads, and algebras over them.  The explicit descriptions of $\G$-operads in Proposition \ref{def:g-operad-generating} and Proposition \ref{prop:g-operad-compi} now specialize to the cactus group operad $\Cac$ and the pure cactus group operad $\PCac$.  Moreover, the description of algebras over a $\colorc$-colored $\G$-operad in Proposition \ref{prop:g-operad-algebra-defs} specializes to the (pure) cactus case.  Next we provide some examples of (pure) cactus operads and their algebras.

\section{Examples of Cactus Operads}\label{sec:cactus-ex}

Some examples of (pure) cactus operads are given in this section.

\begin{example}[Cactus Operads as Pure Cactus Operads]\label{ex:cacop-pcacop}
Every $\colorc$-colored cactus operad is naturally a $\colorc$-colored pure cactus operad by restriction of the equivariant structure.\dqed\end{example}

\begin{example}[Colored Endomorphism Cactus Operads]\label{ex:end-cactus}
For each $\colorc$-colored object $X=\{X_c\}_{c\in\colorc}$ in $\M$, the endomorphism operad\index{endomorphism operad!cactus}\index{cactus operad!endomorphism} $\End(X)$, which is a $\colorc$-colored symmetric operad, is naturally a $\colorc$-colored cactus operad.  This is the $\G=\Cac$ special case of Example \ref{ex:end-g}.\dqed\end{example}

\begin{example}[Cactus Group Operad as a Cactus Operad]
The cactus group operad $\Cac=\{Cac_n\}_{n\geq 0}$ yields a one-colored cactus operad in $\Set$.  This is Example \ref{ex:agop-goperad} when $\G=\Cac$.  Similarly, the pure cactus group operad $\PCac = \{PCac_n\}_{n\geq 0}$ yields a one-colored pure cactus operad in $\Set$.\dqed
\end{example}

\begin{example}[Translation Category of Cactus Group Operad]
Applying the translation category\index{translation category!cactus group operad}\index{cactus group operad!translation category} construction $E_?$ to the cactus group operad $\Cac$, we obtain a one-colored cactus operad \[E_{\Cac} = \Bigl(\bigl\{E_{Cac_n}\bigr\}_{n\geq 0}, \gamma,\operadunit\Bigr)\] in the symmetric monoidal category of small groupoids.  This is Example \ref{ex:translation-cat-goperad} when $\G=\Cac$.  Similarly, the pure cactus group operad $\PCac$ yields a one-colored pure cactus operad \[E_{\PCac} = \Bigl(\bigl\{E_{PCac_n}\bigr\}_{n\geq 0}, \gamma,\operadunit\Bigr)\] in the symmetric monoidal category of small groupoids.\dqed
\end{example}

\begin{example}[Nerves of Translation Categories of Cactus Groups]\label{ex:translation-cactus}
Applying the nerve functor to the one-colored cactus operad $E_{\Cac}$ in small groupoids, we obtain a one-colored cactus operad\index{translation category!nerve of} \[E\Cac = \{ECac_n\}_{n\geq 0}\] in the category of simplicial sets, where \[ECac_n = \Nerve(E_{Cac_n})\] is the nerve of the translation category $E_{Cac_n}$ of the $n$th cactus group $Cac_n$.  This is the $\G=\Cac$ special case of Example \ref{ex:translation-g}.  Furthermore, the one-colored cactus operad $E\Cac$ is a\index{Cacinfinity-operad@$\Cac_\infty$-operad}\index{operad!cacinfinity@$\Cac_\infty$} \emph{$\Cac_\infty$-operad}.  That is, each $ECac_n$ is contractible, and the cactus group action \[ECac_n \times Cac_n \to ECac_n\] is free and proper.  

Similarly, there is a one-colored pure cactus operad \[E\PCac = \{EPCac_n\}_{n\geq 0}\] in the category of simplicial sets, where \[EPCac_n = \Nerve(E_{PCac_n})\] is the nerve of the translation category $E_{PCac_n}$ of the $n$th pure cactus group $PCac_n$. Moreover, each $EPCac_n$ is contractible, and the pure cactus group action \[EPCac_n \times PCac_n \to EPCac_n\] is free and proper.\dqed\end{example}

\begin{example}[Parenthesized Cactus Operad]\label{ex:parenthesized-cactus}
The \index{parenthesized cactus operad}\emph{parenthesized cactus operad} $\PaCac$ is the one-colored cactus operad in $\Set$ as in Example \ref{ex:parenthesized-goperad} with $\G=\Cac$.  Its entries are \[\PaCac(n) = \Mag(n) \times Cac_n \times \Mag(n)\] for $n\geq 0$.  For $1 \leq i \leq n$ and $m \geq 0$, the $\compi$-composition \[\nicexy{\PaCac(n) \times \PaCac(m) \ar[r]^-{\compi} & \PaCac(n+m-1)}\] is defined as \[(w;\sigma;w') \compi (x;\tau;x') = \bigl(w\compi x; \sigma\compi\tau; w' \comp_{\sigmabar(i)} x'\bigr)\] for $w,w'\in \Mag(n)$, $x,x' \in \Mag(m)$, $\sigma\in Cac_n$, and $\tau\in Cac_m$, where $\sigma\compi\tau\in Cac_{n+m-1}$ is the $\compi$-composition of cacti. 

Similarly, using the pure cactus group operad $\G=\PCac$ in Example \ref{ex:parenthesized-goperad}, we obtain the \index{parenthesized pure cactus operad}\emph{parenthesized pure cactus operad} $\PaPCac$, which is a one-colored pure cactus operad in $\Set$.\dqed
\end{example}

\section{Relationship with Braid Group Operad}
\label{sec:relationship-braid-cactus}

Near\index{braid group operad!incompatibility with cactus group operad}\index{cactus group operad!incompatibility with braid group operad} the end of the paper \cite{dev}, Devadoss suggested that there may be a close relationship between the braid groups $B_n$ and the cactus groups $Cac_n$.  In fact, that is one reason why cacti are called quasibraids in that paper.  While it may be the case that the cactus groups and the braid groups are related, in this section we observe that there are no morphisms from the cactus groups to the braid groups, or from the braid groups to the cactus groups, that are compatible with their planar operad structures.  In fact, we already observed in Example \ref{ex:second-cactus-group} that there are no group homomorphisms $Cac_2\to B_2$ from the second cactus group to the second braid group.  The following observation rules out the other direction as well.  Recall the braid group operad $\B$ in Definition \ref{def:braid-group-operad} and the cactus group operad $\Cac$ in Definition \ref{def:cactus-group-operad}.

\begin{theorem}\label{no-braid-to-cactus}
There does not exist any sequence of group homomorphisms \[\phi_n : B_n \to Cac_n \forspace n \geq 0\] that forms a morphism of planar operads in $\Set$ from the braid group operad to the cactus group operad.
\end{theorem}

\begin{proof}
Suppose such a sequence $\phi = \{\phi_n\}$ of group homomorphisms exists.  We will show that this supposition leads to a contradiction.  

First, there is only one possible group homomorphism \[\phi_2 : \bbZ \cong B_2 \to Cac_2 \cong \bbZ/2,\] which must send the generator $s_1\in B_2$ to the generator $s^{(2)}_{1,2} \in Cac_2$.  In the next paragraph, to avoid confusion, we will write the generators $s_1,s_2\in B_n$ as $s_1^{(n)}, s_2^{(n)}$.   We will also use Proposition \ref{third-cactus-group}, which identifies the cactus group $Cac_3$ with the quotient $F\{x,y\}/(x^2,y^2)$ in which $x=s^{(3)}_{1,2}$ and $y=s^{(3)}_{1,3}$.

Since $\phi$ is level-wise a group homomorphism and is compatible with the operadic composition, the map $\phi_3 : B_3 \to Cac_3$ must send 
\[\begin{cases}\nicexy{B_3 \ni s_1^{(3)} = \gamma^{\B}\bigl(\id_2; s_1^{(2)}, \id_1\Bigr) \ar@{|->}[r]^-{\phi_3} & \gammacac\bigl(\id_2; s^{(2)}_{1,2}, \id_1\bigr) = s^{(2)}_{1,2} \oplus \id_1 = s^{(3)}_{1,2} = x},\\
\nicexy{B_3 \ni s_2^{(3)} = \gamma^{\B}\bigl(\id_2; \id_1, s_1^{(2)}\bigr)  \ar@{|->}[r]^-{\phi_3} & \gammacac\bigl(\id_2; \id_1, s^{(2)}_{1,2}\bigr) = \id_1\oplus s^{(2)}_{1,2} = s^{(3)}_{2,3} = yxy}.
\end{cases}\]
The braid relation \[s_1^{(3)} s_2^{(3)} s_1^{(3)} = s_2^{(3)} s_1^{(3)} s_2^{(3)} \in B_3\] is preserved by the group homomorphism $\phi_3$.  This means that there is an equality \[x(yxy)x = (yxy)x(yxy) \in Cac_3.\]  Since $x^2=y^2=\id_3 \in Cac_3$, multiplying $xyxyx$ to the previous equality yields 
\[\begin{split}
\id_3 &= (xyxyx)(xyxyx)\\
&= (xyxyx)(yxyxyxy) \in Cac_3.
\end{split}\]  
But the right-most element in the previous line is not the unit in $Cac_3$ by Proposition \ref{third-cactus-group}, so the previous equality is actually a contradiction.
\end{proof}

%% file: constructions_group_operads.tex
\part{Constructions of Group Operads}
\label{part:construction-group-operad}

\chapter{Naturality}\label{ch:naturality}

In this chapter, we consider naturality properties of the category of $\colorc$-colored $\G$-operads and the category of algebras over a given $\colorc$-colored $\G$-operad.  As in previous chapters, $(\M,\otimes,\tensorunit)$ is a complete and cocomplete symmetric monoidal category in which the monoidal product commutes with small colimits on each side.

In Section \ref{sec:changing-group-operad} we consider the category of $\colorc$-colored $\G$-operads in $\M$ under a change of the action operad $\G$.  As a special case, the augmentation $\omega : \G \to\S$ from $\G$ to the symmetric group operad $\S$ induces an adjunction between the category of $\colorc$-colored $\G$-operads and the category of $\colorc$-colored symmetric operads.  Via this adjunction, each $\colorc$-colored $\G$-operad and its associated $\colorc$-colored symmetric operad have canonically isomorphic categories of algebras.

Section \ref{sec:symmetrization} contains examples of the results in Section \ref{sec:changing-group-operad}, including adjunctions between categories of planar operads, ribbon operads, symmetric operads, and so forth.

In Section \ref{sec:changing-base} we consider the category of $\colorc$-colored $\G$-operads when the ambient symmetric monoidal category is changed by a symmetric monoidal functor.  This applies to, for example, the singular chain functor from spaces to chain complexes and the nerve functor from small categories to simplicial sets.  We also consider the situation when the ambient symmetric monoidal category $\M$ and the action operad $\G$ are both changed.

In Section \ref{sec:change-algebra-category} we consider the category of $\O$-algebras when $\O$ is changed by a morphism of $\colorc$-colored $\G$-operads.  We observe that there is an induced adjunction between the algebra categories.

\section{Change of Action Operads}\label{sec:changing-group-operad} 

In this section, we observe that each morphism of action operads induces an adjunction between the respective categories of group operads.  As a consequence, for each action operad $\G$ and each $\colorc$-colored $\G$-operad $\O$, the category of $\O$-algebras can be identified with the category of algebras over an associated $\colorc$-colored symmetric operad.  See Corollary \ref{algebra-cat-bicomplete}.  We will also consider the situation when there are simultaneous changes of the action operad and the color set.  Specific examples will be given in Section \ref{sec:symmetrization}.

Recall from Definition \ref{def:augmented-group-operad} the concept of an action operad $(\G,\omega)$.

\begin{definition}\label{def:augmented-group-operad-morphism}
Suppose $(\Gone,\omega^1)$ and $(\Gtwo,\omega^2)$ are action operads.  A \emph{morphism}\index{action operad!morphism}\index{morphism!action operads} \[\varphi : (\Gone,\omega^1) \to (\Gtwo,\omega^2)\] of action operads is a morphism $\varphi : \Gone \to \Gtwo$ of one-colored planar operads in $\Set$ such that:
\begin{enumerate}[label=(\roman*)]
\item Level-wise $\varphi : \Gone(n) \to \Gtwo(n)$ is a group homomorphism.
\item The diagram \[\nicexy{\Gone \ar[r]^-{\varphi} \ar[dr]_-{\omega^1} & \Gtwo \ar[d]^-{\omega^2}\\ & \S}\] is commutative.
\end{enumerate}   
The category\index{category!action operads}\index{action operad!category} of action operads is denoted by\label{not:actop} $\agop$.
\end{definition}

The following statements are immediate from the definition.

\begin{lemma}\label{rk:varphi-augmentation}
In the context of Definition \ref{def:augmented-group-operad-morphism}:
\begin{enumerate}
\item The planar group operad $\P$ in Example \ref{ex:trivial-group-operad} is an initial object in the category $\agop$ of action operads.
\item The symmetric group operad $\S$ in Example \ref{ex:symmetric-group-operad} is a terminal object in the category $\agop$ of action operads.
\item For each element $g \in \Gone(n)$, the underlying permutations \[\gbar = \omega^1(g) = \omega^2(\varphi(g)) = \overline{\varphi(g)} \in S_n\] agree.  
\item $\varphi$ induces an object-preserving functor $\varphi : \gonesubc \to \gtwosubc$.
\item $\varphi$ induces an object-preserving functor 
\begin{equation}\label{gtwoseq-goneseq}
\nicexy{\goneseqcm = \M^{\gonesubcopc} & \M^{\gtwosubcopc} = \gtwoseqcm \ar[l]_-{\varphi^*}}.
\end{equation}
\item $\varphi$ induces a functor 
\begin{equation}\label{gtwoop-goneop}
\nicexy{\goneopcm & \gtwoopcm \ar[l]_-{\varphi^*}}
\end{equation}
that preserves the underlying $\colorc$-colored planar operads. 
\end{enumerate}
\end{lemma}

In the rest of this section, we will provide explicit descriptions of the left adjoints of the functors $\varphi^*$.  Fix a morphism $\varphi : (\Gone,\omega^1) \to (\Gtwo,\omega^2)$ of action operads.

\begin{definition}\label{def:phi-seq-left-adjoint}
Suppose $\O \in \goneseqcm$.  Define an object \[\varphi_!\O \in \gtwoseqcm = \M^{\gtwosubcopc}\] entrywise as the coend
\begin{equation}\label{phi-of-o}
(\varphi_!\O)\duc = \int^{\ua \in \gonesubc} \coprodover{\gtwosubc(\uc;\varphi\ua)} \O\dua \in \M
\end{equation}
for $\duc \in \gtwosubcopc$.  Its $\gtwosubcop$-equivariant structure comes from the $\uc$ variable in the set $\gtwosubc(\uc;\varphi\ua)$ in \eqref{phi-of-o}.  In other words, for $\sigma \in \Gtwo(|\uc|)$, the right $\sigma$-action is determined by the commutative diagram
\begin{equation}\label{sigma-phi-o}
\nicexy{(\varphi_!\O)\duc = \dint^{\ua \in \gonesubc} \coprodover{\gtwosubc(\uc;\varphi\ua)} \O\dua \ar[r]^-{\sigma} & \dint^{\ub \in \gonesubc} \coprodover{\gtwosubc(\uc\sigmabar;\varphi\ub)} \O\dub = (\varphi_!\O)\ducsigmabar\\
\coprodover{\gtwosubc(\uc;\varphi\ua)} \O\dua \ar[u]^-{\mathrm{natural}} & \coprodover{\gtwosubc(\uc\sigmabar;\varphi\taubar\uc)} \O\sbinom{d}{(\taubar\sigmabar)(\uc\sigmabar)} \ar[u]_-{\mathrm{natural}}\\
\O\sbinom{d}{\taubar\uc} \ar[u]^-{\tau}_-{\mathrm{summand}} \ar[r]^-{\Id} & \O\sbinom{d}{(\taubar\sigmabar)(\uc\sigmabar)} \ar[u]^-{\tau\sigma}_-{\mathrm{summand}}}
\end{equation}
for $\ua \in \gonesubc$ and $\tau \in \gtwosubc(\uc;\varphi\ua)$.
\end{definition}

\begin{interpretation} To understand the coend formula \eqref{phi-of-o} for $\varphi_!\O$, consider the situation where the morphism $\varphi : \Gone \to \Gtwo$ is level-wise a surjection, such as the projection $B_n \to S_n$ from the $n$th braid group to the $n$th symmetric group.  The coend in \eqref{phi-of-o} is a quotient indexed by $\gonesubc$.  The entire formula divides out part of the $\gonesubcop$-equivariant structure of $\O$ to the point where only the $\gtwosubcop$-equivariant structure remains.  On the other hand, if the morphism $\varphi$ is level-wise an inclusion, such as $\P(n) = \{\id\} \to B_n$, then the coend in \eqref{phi-of-o} gives $\O$ the partially free $\gtwosubcop$-equivariant structure with the original $\gonesubcop$-equivariant structure taken into account.\dqed
\end{interpretation}

\begin{remark}Since the functor $\varphi : \gonesubc \to \gtwosubc$ keeps the objects (i.e., $\colorc$-profiles) unchanged, the set $\gtwosubc(\uc;\varphi\ua)$ in \eqref{phi-of-o} is equal to \[\gtwosubc(\uc;\ua) = \Bigl\{\sigma \in \Gtwo(|\uc|) : \sigmabar\uc=\ua\Bigr\}\] for each object $\ua \in \gonesubc$ (i.e., for each $\colorc$-profile $\ua$).  To simplify the notation below, we will usually suppress $\varphi$ in the set $\gtwosubc(\uc;\varphi\ua)$ and just write $\gtwosubc(\uc;\ua)$.\dqed
\end{remark}

\begin{remark} In the diagram \eqref{sigma-phi-o}, the upper left natural morphism is for a generic object $\ua \in \gonesubc$ (i.e., a $\colorc$-profile).  The lower left vertical morphism is the coproduct summand inclusion corresponding to a generic element $\tau \in \gtwosubc(\uc;\ua)$, so $\taubar\uc=\ua$.  The lower right vertical morphism is the coproduct summand inclusion corresponding to $\tau\sigma \in \gtwosubc(\uc\sigmabar; \ub)$, so \[\taubar\sigmabar(\uc\sigmabar) = \taubar\uc = \ub.\]  The upper right natural morphism is for the object $\taubar\uc \in \gonesubc$.\dqed
\end{remark}

\begin{lemma}\label{phistar-seq-left-adjoint}\index{change of action operads!G-sequence@$\G$-sequence}
The\index{G-sequence@$\G$-sequence!change of action operads} construction $\varphi_!$ in \eqref{phi-of-o} defines a functor \[\varphi_! : \goneseqcm \to \gtwoseqcm\] that is left adjoint to $\varphi^*$ in \eqref{gtwoseq-goneseq}.
\end{lemma}

\begin{proof}
The functoriality of $\varphi_!$ is immediate from its definition.  That $\varphi_!$ defines a left adjoint of the functor $\varphi^*$ is also immediate from the coend definition in \eqref{phi-of-o}.
\end{proof}

For a $\colorc$-colored $\Gone$-operad $\O$, we now equip $\varphi_!\O$ with a $\colorc$-colored $\Gtwo$-operad structure.  To simplify the presentation below, for an action operad $\G$, the forgetful functor from $\goperadcm$ to $\gseqcm$, which forgets about the underlying $\colorc$-colored planar operad structure, will be suppressed from the notation.

\begin{definition}\label{def:phi-operad-left-adjoint}
Suppose $(\O,\comp,\operadunit^{\O}) \in \goneopcm$.  Consider the object $\varphi_!\O \in \gtwoseqcm$ in Definition \ref{def:phi-seq-left-adjoint}.
\begin{enumerate}
\item For each color $d \in \colorc$, equip $\varphi_!\O$ with the $d$-colored unit $\operadunit_d$ given by the following composite. \[\nicexy@C+.5cm{\tensorunit \ar[d]_-{\operadunit^{\O}_d} \ar[rr]^-{\operadunit_d} && (\varphi_!\O)\dd\\ \O\dd \ar[r]^-{\id_1\in \Gtwo(1)}_-{\mathrm{summand}} & \coprodover{\gtwosubc(d;d)} \O\dd \ar[r]^-{\mathrm{natural}} & \dint^{\ua \in \gonesubc} \coprodover{\gtwosubc(d;\ua)} \O\dua \ar@{=}[u]}\] Note that $\gtwosubc(d;d) = \Gtwo(1)$.  The bottom left horizontal morphism above is the coproduct summand inclusion corresponding to the group multiplication unit $\id_1$ in $\Gtwo(1)$.  The bottom right horizontal morphism is the natural morphism into the coend corresponding to the object $d \in \gonesubc$.
\item Suppose $d \in \colorc$, $\uc = (c_1,\ldots,c_n) \in \Profc$ with $n \geq 1$, $1 \leq i \leq n$, and $\ub=(b_1,\ldots,b_m) \in \Profc$.  Define $\compi$ for $\varphi_!\O$ via the  commutative diagram 
\begin{equation}\label{phi-o-compi}
\nicexy@C-.6cm{(\varphi_!\O)\duc \otimes (\varphi_!\O)\ciub \ar[r]^-{\compi} & (\varphi_!\O)\duccompiub\\
\dint^{(\ua^1,\,\ua^2) \in \gonesubc\times\gonesubc}\minusquad \coprodover{\gtwosubc(\uc;\ua^1) \times \gtwosubc(\ub;\ua^2)} \O\duaone \otimes \O\ciuatwo \ar[u]^-{\cong} & \dint^{\ua \in \gonesubc} \coprodover{\gtwosubc(\uc\compi\ub;\,\ua)} \O\dua \ar@{=}[u]\\
\coprodover{\gtwosubc(\uc;\ua^1) \times \gtwosubc(\ub;\ua^2)} \O\duaone \otimes \O\ciuatwo \ar[u]^-{\mathrm{natural}} & \coprodover{\gtwosubc(\uc\compi\ub;\, \ua^1 \comp_{\sigmabar(i)}\ua^2)} \O\duaonecompsigmabariuatwo \ar[u]_-{\mathrm{natural}}\\
\O\duaone \otimes \O\ciuatwo \ar[u]^-{(\sigma,\tau)}_-{\mathrm{summand}} \ar[r]^-{\comp_{\sigmabar(i)}} & \O\duaonecompsigmabariuatwo \ar[u]^-{\sigma\compi\tau}_-{\mathrm{summand}}}
\end{equation}
for $(\ua^1,\,\ua^2) \in \gonesubc\times\gonesubc$ and $(\sigma,\tau) \in \gtwosubc(\uc;\ua^1)\times \gtwosubc(\ub;\ua^2)$.
\end{enumerate}
\end{definition}

\begin{remark} In the diagram \eqref{phi-o-compi}, the upper left isomorphism is the following composite. \[\begin{split} (\varphi_!\O)\duc \otimes (\varphi_!\O)\ciub &=
\Bigl[\int^{\ua^1 \in \gonesubc} \coprodover{\gtwosubc(\uc;\ua^1)} \O\duaone\Bigr] \otimes 
\Bigl[\int^{\ua^2 \in \gonesubc} \coprodover{\gtwosubc(\ub;\ua^2)} \O\ciuatwo\Bigr]\\
&\cong \dint^{(\ua^1,\,\ua^2) \in \gonesubc\times\gonesubc} \Bigl[\coprodover{\gtwosubc(\uc;\ua^1)} \O\duaone\Bigr] \otimes \Bigl[\coprodover{\gtwosubc(\ub;\ua^2)} \O\ciuatwo\Bigr]\\  
&\cong \dint^{(\ua^1,\,\ua^2) \in \gonesubc\times\gonesubc}\minusquad \coprodover{\gtwosubc(\uc;\ua^1) \times \gtwosubc(\ub;\ua^2)} \O\duaone \otimes \O\ciuatwo 
\end{split}\]
The left vertical natural morphism in the diagram \eqref{phi-o-compi} is for a generic pair of objects \[(\ua^1, \ua^2) \in \gonesubc\times\gonesubc.\]  The right vertical natural morphism is for the object \[\ua^1 \comp_{\sigmabar(i)} \ua^2 \in \gonesubc.\]  For the bottom left vertical summand inclusion, note that $\sigma \in \Gtwo(n)$ with underlying permutation $\sigmabar\in S_n$ and $\tau \in \Gtwo(m)$ with underlying permutation $\taubar \in S_n$ such that \[\begin{split}\sigmabar\uc &= \bigl(c_{\sigmabar^{-1}(1)},\ldots,c_{\sigmabar^{-1}(n)}\bigr)= \bigl(a^1_1,\ldots,a^1_n) = \ua^1,\\
\taubar\ub &= \bigl(b_{\taubar^{-1}(1)},\ldots,b_{\taubar^{-1}(m)}\bigr)= \bigl(a^2_1,\ldots,a^2_m) = \ua^2.\end{split}\]  It follows that \[c_i = a^1_{\sigmabar(i)} \andspace \bigl(\sigma\compi\tau\bigr) \bigl(\uc\compi\ub\bigr) = \ua^1 \comp_{\sigmabar(i)} \ua^2.\] The bottom horizontal morphism is the $\comp_{\sigmabar(i)}$-composition in the $\colorc$-colored $\Gone$-operad $\O$.\dqed
\end{remark}

\begin{lemma}\label{phio-gtwo-operad}
For each $(\O,\comp,\operadunit^{\O}) \in \goneopcm$, the colored units and $\compi$ in Definition \ref{def:phi-operad-left-adjoint} give the object \[\varphi_!\O \in \gtwoseqcm\] the structure of a $\colorc$-colored $\Gtwo$-operad in $\M$.
\end{lemma}

\begin{proof}
Both the colored units and the $\compi$-composition in $\varphi_!\O$ are naturally defined by those in $\O$.  The $\colorc$-colored $\Gtwo$-operadic associativity and unity axioms for $\varphi_!\O$ thus follow from those of $\O$ as a $\colorc$-colored $\Gone$-operad.  In terms of the diagram \eqref{phi-o-compi}, the equivariance axiom for $\compi$ in $\varphi_!\O$ comes from the $\uc$ variable in $\gtwosubc(\uc;\ua^1)$ and the $\ub$ variable in $\gtwosubc(\ub;\ua^2)$. Via $\sigma$ and $\tau$, this translates into the equivariance axiom for $\comp_{\sigmabar(i)}$ in $\O$. 
\end{proof}

\begin{theorem}\label{phi-goneop-gtwoop}\index{change of action operads!G-operad@$\G$-operad}
For\index{G-operad@$\G$-operad!change of action operads} each morphism $\varphi : (\Gone,\omega^1) \to (\Gtwo,\omega^2)$ of action operads, Definitions \ref{def:phi-seq-left-adjoint} and \ref{def:phi-operad-left-adjoint} define a functor \[\varphi_! : \goneopcm \to \gtwoopcm\] that is left adjoint to the functor $\varphi^*$ in \eqref{gtwoop-goneop}.
\end{theorem}

\begin{proof}
As in Lemma \ref{phistar-seq-left-adjoint}, this is immediate from the definitions of $\varphi_!\O$ and its colored units and $\compi$-composition.
\end{proof}

Using the adjunction in Theorem \ref{phi-goneop-gtwoop}, we now observe that the category of algebras over a $\G$-operad is isomorphic to the category of algebras over an associated symmetric operad.

\begin{corollary}\label{algebra-cat-bicomplete}
Suppose\index{algebra!of a G-operad@of a $\G$-operad} $(\G,\omega)$ is an action operad, and $\O$ is a $\colorc$-colored $\G$-operad.  Then the following statements hold.
\begin{enumerate}
\item  There is a canonical isomorphism of categories \[\algmo \cong \algm(\omega_!\O),\] in which $\omega_!$ is the left adjoint in the adjunction  \[\nicexy{\goperadcm \ar@<2pt>[r]^-{\omega_!} & \Soperadcm \ar@<2pt>[l]^-{\omega^*}}\] in Theorem \ref{phi-goneop-gtwoop} induced by the augmentation $\omega : \G\to\S$.
\item The category $\algmo$ is complete and cocomplete.
\end{enumerate}
\end{corollary}

\begin{proof}
For the first assertion, by Example \ref{ex:end-g}, an $\O$-algebra structure on an object $X\in\Mtoc$ is equivalent to a morphism \[\O\to\omega^*\End(X) \in\goperadcm,\] where $\End(X)$ is the endomorphism symmetric operad of $X$ in Example \ref{ex:endomorphism-operad}.  By the adjunction $\omega_!\dashv\omega^*$, such an $\O$-algebra structure is equivalent to a morphism \[\omega_!\O \to \End(X)\in\Soperadcm.\]  So the category $\algmo$ of $\O$-algebras is isomorphic to the category $\algm(\omega_!\O)$ of $\omega_!\O$-algebras.  

For the second assertion, it is proved in Proposition 4.2.1 in \cite{white-yau} that for each colored symmetric operad, such as $\omega_!\O$, the category of algebras is complete and cocomplete.  Therefore, by the isomorphism in the first assertion, the category $\algmo$ is also complete and cocomplete.
\end{proof}

The following observation is about mixing changes in the color set and the action operad.  Recall the change-of-color functor in Example \ref{ex:g-operad-change-colors}

\begin{proposition}\label{change-color-goperad}\index{change of action operads!mixed with change of colors}
Suppose\index{G-operad@$\G$-operad!change of colors} $\varphi : \Gone\to\Gtwo$ is a morphism of action operads, and $f : \colorc \to \colord$ is a function between sets.  Then the diagram
\[\nicexy{\gtwooperadcm \ar[d]_-{\varphi^*} & \gtwooperaddm \ar[l]_-{f^*} \ar[d]^-{\varphi^*}\\
\goneoperadcm & \goneoperaddm \ar[l]_-{f^*}}\]
is commutative.
\end{proposition}

\begin{proof}
Suppose $\Q$ is a $\colord$-colored $\Gtwo$-operad in $\M$.  Then $f^*\varphi^*\Q$ is the $\colorc$-colored $\Gone$-operad obtained as follows.
\begin{enumerate}[label=(\roman*)]
\item $\varphi^*\Q$ is the $\colord$-colored $\Gone$-operad with the same entries and the same $\colord$-colored planar operad structure as $\Q$.  Its $\Gone$-equivariant structure is induced by the morphism $\varphi$, so $\sigma \in \Gone$ acts on $\varphi^*\Q$ as $\varphi\sigma$ acts on $\Q$. 
\item The bottom horizontal functor $f^*$ is the change-of-color functor \[\nicexy{\goneoperadcm & \goneoperaddm \ar[l]_-{f^*}}\] in Example \ref{ex:g-operad-change-colors} induced by the color change function $f : \colorc\to\colord$.  The object $f^*\varphi^*\Q$ is the $\colorc$-colored $\Gone$-operad with entries \[(f^*\varphi^*\Q)\duc = (\varphi^*\Q)\fdfuc = \Q\fdfuc\] for $\duc \in \Profcc$.  Its $\colorc$-colored $\Gone$-operad structure is directly induced by the $\colord$-colored $\Gone$-operad structure of $\varphi^*\Q$.
\end{enumerate}
A similar description of $\varphi^*f^*\Q$ reveals that it is equal to $f^*\varphi^*\Q$.
\end{proof}

\section{Symmetrization and Other Left Adjoints}\label{sec:symmetrization}

In this section, we discuss the adjunctions that go between various categories of colored operads, such as the one between colored ribbon operads and colored symmetric operads.  

\begin{example}[Planar Operads to Symmetric Operads]\label{ex:pop-sop}
Consider\index{planar operad!to symmetric operad}\index{symmetric operad!from planar operad} the morphism $\iota : \P \to \S$ of action operads, with $\P$ the planar group operad in Example \ref{ex:trivial-group-operad} and $\S$ the symmetric group operad in Example \ref{ex:symmetric-group-operad} with the identity augmentation.  Level-wise the morphism $\iota : \P(n) = \{\id\} \to S_n$ is the unit inclusion.  The functor \[\iota^* : \Soperadcm \to \Poperadcm\] is the forgetful functor from the category of $\colorc$-colored symmetric operads to the category of $\colorc$-colored planar operads that forgets about the equivariant structure.  Its left adjoint \[\iota_! : \Poperadcm \to \Soperadcm\] is given entrywise by \[(\iota_!\O)\duc = \int^{\ua \in \Profc} \coprodover{\ssubc(\uc;\ua)} \O\dua = \coprodover{\sigma \in S_{|\uc|}} \O\dsigmauc\] for $\O \in \Poperadcm$ and $\duc \in \ssubcopc$.  This entrywise expression and the rest of the $\colorc$-colored symmetric operad structure on $\iota_!\O$ agree with Definition 20.1.2 in \cite{yau-operad}, where the functor $\omega_!$ is called \index{symmetrization}\emph{symmetrization}.  In the one-colored case (i.e., $\colorc = \{*\}$), the above formula reduces to \[(\iota_!\O)(n) = \coprodover{\sigma \in S_n} \O(n) = \O(n) \cdot S_n\] for $n \geq 0$.\dqed
\end{example}

\begin{example}[Planar Operads to Ribbon Operads]\label{ex:pop-rop}
Consider\index{planar operad!to ribbon operad}\index{ribbon operad!from planar operad} the morphism $\iota : \P \to \R$ of action operads, where $\R$ is the ribbon group operad in Definition \ref{def:ribbon-group-operad}.  Once again $\iota^*$ is the forgetful functor.  Its left adjoint \[\iota_! : \Poperadcm \to \Roperadcm\] is given entrywise by \[(\iota_!\O)\duc = \int^{\ua \in \Profc} \coprodover{\rsubc(\uc;\ua)} \O\dua = \coprodover{\sigma \in R_{|\uc|}} \O\dsigmabaruc\] for $\O \in \Poperadcm$ and $\duc \in \rsubcopc$, where $\sigmabar \in S_{|\uc|}$ is the underlying permutation of $\sigma$.  This functor $\iota_!$ is called \index{ribbonization}\emph{ribbonization}.  In the one-colored case, the above formula reduces to \[(\iota_!\O)(n) = \coprodover{\sigma \in R_n} \O(n)= \O(n) \cdot R_n\] for $n \geq 0$.\dqed
\end{example}

\begin{example}[Planar Operads to Braided Operads]\label{ex:pop-bop} 
Similarly, consider the morphism\index{planar operad!to braided operad}\index{braided operad!from planar operad} $\iota : \P \to \B$ of action operads, where $\B$ is the braid group operad in Definition \ref{def:braid-group-operad}.  The left adjoint \[\iota_! : \Poperadcm \to \Boperadcm\] to the forgetful functor $\iota^*$ is given entrywise by \[(\iota_!\O)\duc = \int^{\ua \in \Profc} \coprodover{\bsubc(\uc;\ua)} \O\dua = \coprodover{\sigma \in B_{|\uc|}} \O\dsigmabaruc\] for $\O \in \Poperadcm$ and $\duc \in \bsubcopc$.  This functor $\iota_!$ is called\index{braidization} \emph{braidization}.  In the one-colored case, the above formula reduces to \[(\iota_!\O)(n) = \coprodover{\sigma \in B_n} \O(n)= \O(n) \cdot B_n\] for $n \geq 0$.\dqed
\end{example}

\begin{example}[Ribbon Operads to Symmetric Operads]\label{ex:rop-sop}
Consider\index{ribbon operad!to symmetric operad}\index{symmetric operad!from ribbon operad}  the morphism $\pi : \R \to \S$ of action operads, which is level-wise the underlying permutation map $\pi : R_n \to S_n$.  The right adjoint \[\pi^* : \Soperadcm \to \Roperadcm\] sends a $\colorc$-colored symmetric operad $\O$ to the $\colorc$-colored ribbon operad $\pi^*\O$ with the same underlying $\colorc$-colored planar operad whose $\rsubcop$-equivariant structure is induced by $\pi$.  Its left adjoint \[\pi_! : \Roperadcm \to \Soperadcm,\] called \index{symmetrization!of ribbon operad}\emph{symmetrization}, is given entrywise by \[(\pi_!\O)\duc = \int^{\ua \in \rsubc} \coprodover{\ssubc(\uc;\ua)} \O\dua\] for $\O \in \Roperadcm$ and $\duc \in \ssubcopc$.  In the one-colored case, this formula reduces to \[(\pi_!\O)(n) = \O(n) \dotover{R_n} S_n = \O(n)/PR_n\] for $n \geq 0$.  Here $PR_n$ denotes the $n$th pure ribbon group, which is the kernel of the underlying permutation map $\pi : R_n \to S_n$.  

For example, via the level-wise universal covering map \[p : \Dftildetwo \to \Df_2,\] the level-wise universal cover $\Dftildetwo$ in Example \ref{ex:universal-cover-fd2}, which is an $\R_\infty$-operad, has the framed little $2$-disc operad $\Df_2$ as its symmetrization \cite{wahl}.  By Corollary \ref{algebra-cat-bicomplete} the category of $\Dftildetwo$-algebras and the category of $\Df_2$-algebras are canonically isomorphic.\dqed
\end{example}

\begin{example}[Braided Operads to Symmetric Operads]\label{ex:bop-sop}
Consider\index{braided operad!to symmetric operad}\index{symmetric operad!from braided operad}  the morphism $\pi : \B \to \S$ of action operads, which is level-wise the underlying permutation map $\pi : B_n \to S_n$.  The right adjoint \[\pi^* : \Soperadcm \to \Boperadcm\] sends a $\colorc$-colored symmetric operad $\O$ to the $\colorc$-colored braided operad $\pi^*\O$ with the same underlying $\colorc$-colored planar operad whose $\bsubcop$-equivariant structure is induced by $\pi$.  Its left adjoint \[\pi_! : \Boperadcm \to \Soperadcm,\] called \index{symmetrization!of braided operad}\emph{symmetrization}, is given entrywise by \[(\pi_!\O)\duc = \int^{\ua \in \bsubc} \coprodover{\ssubc(\uc;\ua)} \O\dua\] for $\O \in \Boperadcm$ and $\duc \in \ssubcopc$.  In the one-colored case, this formula reduces to \[(\pi_!\O)(n) = \O(n) \dotover{B_n} S_n = \O(n)/PB_n\] for $n \geq 0$.  Here $PB_n$ denotes the $n$th pure braid group, which is the kernel of the underlying permutation map $\pi : B_n \to S_n$.  For example:
\begin{enumerate}
\item Via the level-wise universal covering map\index{little $2$-cube operad} \[p : \Ctilde_2 \to \C_2,\] the level-wise universal cover $\Ctilde_2$ in Example \ref{ex:universal-cover-c2}, which is a $\B_\infty$-operad, has the little $2$-cube operad $\C_2$ as its symmetrization \cite{fiedorowicz}.  By Corollary \ref{algebra-cat-bicomplete} the category of $\Ctilde_2$-algebras and the category of $\C_2$-algebras are canonically isomorphic.
\item Via the level-wise universal covering map\index{little $2$-disc operad} \[p : \Dtilde_2 \to \D_2,\] the level-wise universal cover $\Dtilde_2$ in Example \ref{ex:universal-cover-d2}, which is a $\B_\infty$-operad, has the little $2$-disc operad $\D_2$ as its symmetrization \cite{fresse-gt}.  By Corollary \ref{algebra-cat-bicomplete} the category of $\Dtilde_2$-algebras and the category of $\D_2$-algebras are canonically isomorphic.\dqed   
\end{enumerate}
\end{example}

\begin{example}[Planar Operads to Cactus Operads]\label{ex:pop-cacop}
Consider the morphism\index{planar operad!to cactus operad}\index{cactus operad!from planar operad} $\iota : \P \to \Cac$ of action operads from the planar group operad to the cactus group operad in Definition \ref{def:cactus-group-operad}, which is level-wise the unit inclusion $\iota : \P(n)=\{\id\}\to Cac_n$.  The functor \[\iota^* : \Cacoperadcm \to \Poperadcm\] is the forgetful functor from the category of $\colorc$-colored cactus operads to the category of $\colorc$-colored planar operads that forgets about the equivariant structure.  Its left adjoint \[\iota_! : \Poperadcm \to \Cacoperadcm\] is given entrywise by \[(\iota_!\O)\duc = \int^{\ua \in \Profc} \coprodover{\cacsubc(\uc;\ua)} \O\dua = \coprodover{\sigma \in Cac_{|\uc|}} \O\dsigmauc\] for $\O \in \Poperadcm$ and $\duc \in \cacsubcopc$.  This function $\iota_!$ is called \index{cactusization}\emph{cactusization}.  In the one-colored case, the above formula reduces to \[(\iota_!\O)(n) = \coprodover{\sigma \in Cac_n} \O(n) = \O(n) \cdot Cac_n\] for $n \geq 0$.\dqed
\end{example}

\begin{example}[Cactus Operads to Symmetric Operads]\label{ex:cacop-sop}
Consider the morphism\index{cactus operad!to symmetric operad}\index{symmetric operad!from cactus operad} $\pi : \Cac \to \S$ of action operads, which is level-wise the underlying permutation map $\pi : Cac_n \to S_n$.  The right adjoint \[\pi^* : \Soperadcm \to \Cacoperadcm\] sends a $\colorc$-colored symmetric operad $\O$ to the $\colorc$-colored cactus operad $\pi^*\O$ with the same underlying $\colorc$-colored planar operad whose $\cacsubcop$-equivariant structure is induced by $\pi$.  Its left adjoint \[\pi_! : \Cacoperadcm \to \Soperadcm,\] called \index{symmetrization!of cactus operad}\emph{symmetrization}, is given entrywise by \[(\pi_!\O)\duc = \int^{\ua \in \cacsubc} \coprodover{\ssubc(\uc;\ua)} \O\dua\] for $\O \in \Cacoperadcm$ and $\duc \in \ssubcopc$.  In the one-colored case, this formula reduces to \[(\pi_!\O)(n) = \O(n) \dotover{Cac_n} S_n = \O(n)/PCac_n\] for $n \geq 0$.  Here $PCac_n$ denotes the $n$th pure cactus group, which is the kernel of the underlying permutation map $\pi : Cac_n \to S_n$.\dqed
\end{example}

\begin{example}[Summary]\label{ex:all-augmented-group-operads}
Taking\index{action operad!main examples} into account the pure braid group operad $\PB$ in Definition \ref{def:pure-braid-group-operad}, the pure ribbon group operad $\PR$ in Definition \ref{def:pure-ribbon-group-operad}, and the pure cactus group operad $\PCac$ in Definition \ref{def:pure-cactus-group-operad}, there is a commutative diagram 
\[\nicexy{\P\ar@{=}[d] \ar[r]^-{\iota} & \PR \ar[r]^-{\iota} & \R \ar[r]^-{\pi} & \S \ar@{=}[d]\\
\P \ar[r]^-{\iota} & \PB \ar[r]^-{\iota} & \B \ar[r]^-{\pi} & \S\\
\P \ar[r]^-{\iota} \ar@{=}[u] & \PCac \ar[r]^-{\iota} & \Cac \ar[r]^-{\pi} & \S\ar@{=}[u]}\]
of morphisms of action operads.  Every morphism $\iota$ is level-wise an inclusion.  These morphisms induce adjunctions\index{group operad!main examples} \[\nicexy{\Poperadcm \ar@{=}[d] \ar@<2pt>[r]^-{\iota_!} & \PRoperadcm \ar@<2pt>[r]^-{\iota_!} \ar@<2pt>[l]^-{\iota^*} & \Roperadcm \ar@<2pt>[r]^-{\pi_!} \ar@<2pt>[l]^-{\iota^*} & \Soperadcm \ar@<2pt>[l]^-{\pi^*} \ar@{=}[d]\\ 
\Poperadcm \ar@<2pt>[r]^-{\iota_!} & \PBoperadcm \ar@<2pt>[r]^-{\iota_!} \ar@<2pt>[l]^-{\iota^*} & \Boperadcm \ar@<2pt>[r]^-{\pi_!} \ar@<2pt>[l]^-{\iota^*} & \Soperadcm \ar@<2pt>[l]^-{\pi^*}\\
\Poperadcm \ar@<2pt>[r]^-{\iota_!} \ar@{=}[u] & \Pcacoperadcm \ar@<2pt>[r]^-{\iota_!} \ar@<2pt>[l]^-{\iota^*} & \Cacoperadcm \ar@<2pt>[r]^-{\pi_!} \ar@<2pt>[l]^-{\iota^*} & \Soperadcm \ar@<2pt>[l]^-{\pi^*} \ar@{=}[u]}\] between $\colorc$-colored planar operads, $\colorc$-colored (pure) ribbon operads, $\colorc$-colored (pure) braided operads, $\colorc$-colored (pure) cactus operads, and $\colorc$-colored symmetric operads.\dqed
\end{example}

\begin{example}[Ribbon to Braid is a Non-Example]
The underlying braid group homomorphisms $\pi : R_n \to B_n$ from the ribbon groups do \emph{not} form a morphism $\pi : \R\to \B$ of action operads.  The reason is that this is not a morphism of one-colored planar operads in $\Set$.  In fact, these level-wise maps are not compatible with block ribbon and block braid, as we explained in Example \ref{ex:block-ribbon}.  So they are not compatible with the operadic compositions in the ribbon group operad $\R$ and in the braid group operad $\B$.\dqed
\end{example}


\section{Change of Base Categories}\label{sec:changing-base}

Fix an action operad $\G$.  The purpose of this section is to consider the category of $\colorc$-colored $\G$-operads under a change of the ambient symmetric monoidal category.   We observe that each symmetric monoidal functor induces a monoidal functor between the respective monoidal categories of $\colorc$-colored $\G$-sequences, where the monoidal products are the $\colorc$-colored $\G$-circle products.  Passing to their categories of monoids, there is an induced functor between the respective categories of $\colorc$-colored $\G$-operads.  We will also consider the situation when there are simultaneous changes of the ambient symmetric monoidal category and of the action operad.

Recall from Proposition \ref{g-circle-product-monoidal} that $\gseqcm$ is a monoidal category with the $\G$-circle product $\circg$ as the monoidal product and $\I$ as the monoidal unit.  The category $\goperadcm$ of $\colorc$-colored $\G$-operads in $\M$ is the category of monoids in $\gseqcm$.  We will denote the $\G$-circle product $\circg$ and the monoidal unit $\I$ by $\circgsubm$ and $\I^{\M}$, respectively, below to make the ambient category $\M$ explicit.

Recall that a \emph{monoidal functor}  \[F : (\M,\tensor{\M},\tensorunit^{\M}) \to (\N,\tensor{\N},\tensorunit^{\N})\] between two monoidal categories is a functor $F : \M \to \N$ equipped with 
\begin{enumerate}[label=(\roman*)]
\item a structure morphism $F_0 : \tensorunit^{\N} \to F\tensorunit^{\M} \in \N$ and 
\item a structure morphism \[F_2(x,y) : Fx \tensor{\N} Fy \to F(x\tensor{\M} y)\in \N\] that is natural in the objects $x$ and $y$, 
\end{enumerate}
satisfying suitable associativity and unity axioms.  The precise formulation is in Definition \ref{def:monoidal-functor}.  If, in addition, $\M$ and $\N$ are symmetric monoidal categories and if $F_2$ is compatible with the symmetry isomorphisms, then $F$ is called a \emph{symmetric monoidal functor}.  The precise formulation is in Definition \ref{def:symmetric-monoidal-functor}.  In the literature, a monoidal functor is sometimes called a \emph{lax} monoidal functor to emphasize that $F_0$ and $F_2$ are not necessarily invertible.

\begin{theorem}\label{change-base-category}
Suppose \[(F,F_2,F_0) : (\M,\tensor{\M},\tensorunit^{\M}) \to (\N,\tensor{\N},\tensorunit^{\N})\] is a symmetric monoidal functor between cocomplete symmetric monoidal categories in which the monoidal products commute with colimits on each side.  Then the following statements hold.
\begin{enumerate}
\item There is an induced monoidal functor\index{G-sequence@$\G$-sequence!change of base} \[\nicexy{\bigl(\gseqcm,\circgsubm,\I^{\M}\bigr) \ar[r]^-{F_*} & \bigl(\gseqcn,\circgsubn,\I^{\N}\bigr)}\] that sends each $X\in\gseqcm$ to the $\colorc$-colored $\G$-sequence in $\N$ given by the composite \[\nicexy{\gsubcopc \ar[r]^-{X} & \M \ar[r]^-{F} & \N \ar@{<-}`u[ll] `[ll]_-{F_*X} [ll]}.\]
\item The monoidal functor $F_*$ in the previous item induces a functor\index{G-operad@$\G$-operad!change of base} \[\nicexy{\goperadcm \ar[r]^-{F_*} & \goperadcn}.\] 
\end{enumerate}
\end{theorem}

\begin{proof}
For the first assertion, since $\gseqcm = \M^{\gsubcopc}$ and $\gseqcn = \N^{\gsubcopc}$ are diagram categories indexed by $\gsubcopc$, composition with $F$ induces a functor \[\nicexy{\gseqcm \ar[r]^-{F_*} & \gseqcn} \] as stated.  

To see that $F_*$ is a monoidal functor, first recall that the monoidal unit $\I^{\M}$ is the $\colorc$-colored $\G$-operad in \eqref{g-unit}.  The structure morphism \[\nicexy{\I^{\N} \ar[r]^-{(F_*)_0} & F_*\I^{\M}}\] is given as follows.
\begin{itemize}
\item For $d\in\colorc$ and an entry of the form $\dd$, it is the composite 
\[\nicexy{(\I^{\N})\dd \ar@{=}[d] \ar[rr]^-{F_*} && (F_*\I^{\M})\dd \ar@{=}[d]\\
\coprodover{\sigma \in \G(1)} \tensorunit^{\N}_{\sigma} \ar[r]^-{\coprod_{\sigma} F_0} & \coprodover{\sigma \in \G(1)} F\tensorunit^{\M}_{\sigma} \ar[r]^-{\{F\iota_{\sigma}\}}
& F\Bigl(\coprodover{\sigma \in \G(1)} \tensorunit^{\M}_{\sigma}\Bigr)}\]
in which:
\begin{itemize}
\item Each $\tensorunit^{\N}_{\sigma}$ is a copy of the monoidal unit $\tensorunit^{\N}$, and similarly for $\tensorunit^{\M}_{\sigma}$.
\item $F_0 : \tensorunit^{\N} \to F\tensorunit^{\M}$ is the structure morphism of $F$.
\item For each $\sigma \in \G(1)$, the morphism \[\nicexy{\tensorunit^{\M}_{\sigma} \ar[r]^-{\iota_{\sigma}} & \coprodover{\sigma \in \G(1)} \tensorunit^{\M}_{\sigma}}\] is the coproduct summand inclusion.  
\end{itemize}
\item For all other entries, it is the unique morphism $\varnothing^{\N}\to F(\varnothing^{\M})$.
\end{itemize}

For the other structure morphism, suppose given $X,Y\in \gseqcm$.  With the notations in Definition \ref{def:g-circle-product}, we first define the morphism $F_2^Y$ as the composite \[\nicexy{(F_*Y)^{\uc}_{\smallsubg}(\ub) \ar@{=}[d] \ar[r]^-{F_2^Y} & F\bigl(Y^{\uc}_{\smallsubg}(\ub)\bigr) \ar@{=}[d]\\
\dint^{\{\ua_j\}_{j=1}^m} \gsubcop(\ua;\ub) \cdot \Bigl[\overset{m}{\bigtensorover{j=1}} FY\cjuaj\Bigr] \ar[d]_-{F_2^{m-1}} & F\biggl[\dint^{\{\ua_j\}_{j=1}^m} \gsubcop(\ua;\ub) \cdot \Bigl[\overset{m}{\bigtensorover{j=1}} Y\cjuaj\Bigr] \biggr]\\
\dint^{\{\ua_j\}_{j=1}^m} \gsubcop(\ua;\ub) \cdot F\Bigl[\overset{m}{\bigtensorover{j=1}} Y\cjuaj\Bigr] \ar[r]^-{\text{natural}} & \dint^{\{\ua_j\}_{j=1}^m} F\biggl[\gsubcop(\ua;\ub) \cdot \Bigl[\overset{m}{\bigtensorover{j=1}} Y\cjuaj\Bigr]\biggr] \ar[u]_-{\text{natural}}}\]
in $\N$.  Here $\uc=(c_1,\ldots,c_m)$, $\ub$, and $\ua_j\in\gsubc$ with $\ua = \bigl(\ua_1,\ldots,\ua_m\bigr)$.  The vertical morphism $F_2^{m-1}$ is:
\begin{itemize}
\item the identity morphism if $m\leq 1$;
\item the $(m-1)$-fold iterate of the structure morphism $F_2$ if $m \geq 2$.  
\end{itemize}
The two natural morphisms come from the functoriality of $F$.

The structure morphism $(F_*)_2(X,Y)$ is defined entrywise as the composite 
\[\nicexy{\bigl(F_*X \circgsubn F_*Y\bigr)\dub \ar@{=}[d] \ar[r]^-{(F_*)_2(X,Y)} & F_*\bigl(X\circgsubm Y\bigr)\dub \ar@{=}[d]\\
\dint^{\uc\in\gsubc} (F_*X)\duc \tensor{\N} (F_*Y)^{\uc}_{\smallsubg}(\ub) \ar[d]_-{(F_2^Y)_*} & F\biggl[\dint^{\uc\in\gsubc} X\duc \tensor{\M} Y^{\uc}_{\smallsubg}(\ub)\biggr]\\
\dint^{\uc\in\gsubc} FX\duc \tensor{\N} F\bigl(Y^{\uc}_{\smallsubg}(\ub)\bigr) \ar[r]^-{(F_2)_*} & \dint^{\uc\in\gsubc} F\Bigl(X\duc \tensor{\M} Y^{\uc}_{\smallsubg}(\ub)\Bigr) \ar[u]_-{\text{natural}}}\]
in $\N$ for $\duc\in\Profcc$, where $(F_2^Y)_*$ and $(F_2)_*$ are induced by $F_2^Y$ and $F_2$, respectively.  The natural morphism comes from the functoriality of $F$.  Notice that $(F_*)_2(X,Y)$ is entirely constructed from the structure morphism $F_2$ and the naturality of $F$.  Now observe that
\begin{enumerate}
\item the fact that $(F_*)_2 \in \gseqcn$,
\item the naturality of $(F_*)_2$, and 
\item the required unity and associativity axioms for $\bigl(F_*,(F_*)_2,(F_*)_0\bigl)$ to be a monoidal functor
\end{enumerate}
all boil down to the assumption that $(F,F_2,F_0)$ is a symmetric monoidal functor.

For the second assertion, recall that the category $\goperadcm$ is the category of monoids in the monoidal category $\bigl(\gseqcm,\circgsubm,\I^{\M}\bigr)$.  The assertion now follows from the basic fact that each monoidal functor, such as $F_*$ above, induces a functor between the respective categories of monoids.  See, for example, Theorem 12.11 in \cite{bluemonster}.
\end{proof}

\begin{remark}The one-colored case of Theorem \ref{change-base-category} (i.e., when $\colorc = \{*\}$) is Theorem 1.16 in \cite{gurski}.\dqed\end{remark}

\begin{interpretation}\label{int:goperad-monoidal-functor}
Let us explain more explicitly the functor \[F_* : \goperadcm \to \goperadcn\] in Theorem \ref{change-base-category}.  Suppose $(\O,\comp^{\O},\operadunit^{\O})$ is a $\colorc$-colored $\G$-operad in $\M$ as in Proposition \ref{prop:g-operad-compi}.
\begin{description}
\item[Entries] The underlying entries of $F_*\O$ are obtained from those of $\O$ by applying $F$ entrywise, so \[(F_*\O)\duc = F\O\duc\] for $\duc\in \Profcc$.
\item[Equivariance] The $\G$-sequence structure of $F_*\O$ is obtained from that of $\O$ by applying $F$.  So for each $\sigma \in \G(|\uc|)$, the $\G$-equivariant structure morphism \[\nicexy{(F_*\O)\duc = F\O\duc \ar[r]^-{\sigma} & F\O\ducsigma = (F_*\O)\ducsigma}\] is $F$ applied to the $\G$-equivariant structure morphism $\sigma : \O\duc \to \O\ducsigma$.
\item[Colored Units] For each $d\in \colorc$, the $d$-colored unit of $F_*\O$ is the composite below.\[\nicexy{\tensorunit^{\N} \ar[d]_-{F_0} \ar[r]^-{\operadunit_d^{F_*\O}} & (F_*\O)\dd\\ F\tensorunit^{\M} \ar[r]^-{F(\operadunit_d^{\O})} & F\O\dd \ar@{=}[u]}\]
\item[Composition] For $\duc\in\Profcc$ with $1 \leq i \leq n=|\uc|$ and $\ub\in\Profc$, the $\compi$-composition in $F_*\O$ is the composite below.
\[\nicexy{F\O\duc \tensor{\N} F\O\ciub \ar[d]_-{F_2} \ar[r]^-{\compi} & (F_*\O)\sbinom{d}{\uc\compi\ub}\\
F\Bigl[\O\duc \tensor{\M} \O\ciub\Bigr] \ar[r]^-{F(\compi^{\O})} & F\O\sbinom{d}{\uc\compi\ub} \ar@{=}[u]}\]
\end{description}
\dqed
\end{interpretation}

Next we consider the situation when there are simultaneous changes in the ambient symmetric monoidal category and in the action operad.

\begin{theorem}\label{varphi-F}\index{change of action operads!mixed with change of categories}
Suppose $F : \M\to\N$ is a symmetric monoidal functor as in Theorem \ref{change-base-category}, and $\varphi : \Gone \to \Gtwo$ is a morphism of action operads.  Consider the diagram \[\nicexy{\goneopcm \ar[d]_-{F_*} \ar@<2pt>[r]^-{\varphi_!} & \gtwoopcm \ar@<2pt>[l]^-{\varphi^*} \ar[d]^-{F_*}\\
\goneopcn \ar@<2pt>[r]^-{\varphi_!} & \gtwoopcn \ar@<2pt>[l]^-{\varphi^*}}\] in which each pair $\varphi_! \dashv \varphi^*$ is the adjunction in Theorem \ref{phi-goneop-gtwoop}.  Then:
\begin{enumerate}
\item There is an equality \[\nicexy@C+1cm{\goneopcn & \gtwoopcm \ar[l]_-{F_*\varphi^* \,=\, \varphi^*F_*}}.\]
\item If $F$ commutes with colimits, then there is a natural isomorphism \[\nicexy@C+1cm{\goneopcm \ar[r]^-{\varphi_!F_* \,\cong\, F_*\varphi_!} & \gtwoopcn}.\]
\end{enumerate}
\end{theorem}

\begin{proof}
The first assertion is checked by a simple inspection of the definition of $\varphi^*$ in Lemma \ref{rk:varphi-augmentation} and the explicit description of $F_*$ in \ref{int:goperad-monoidal-functor}.

For the second assertion, suppose $\O$ is a $\colorc$-colored $\Gone$-operad in $\M$.  Recall the definition of $\varphi_!\O$ in \eqref{phi-of-o}.  For each $\duc \in \Profcc$, the assumed commutation of $F$ with colimits implies that there is a natural isomorphism
\[\begin{split}
(F_*\varphi_!\O)\duc &=F\bigl((\varphi_!\O)\duc\bigr)\\
&= F\biggl[\int^{\ua \in \gonesubc} \coprodover{\gtwosubc(\uc;\varphi\ua)} \O\dua\biggr]\\
&\cong \int^{\ua \in \gonesubc} \coprodover{\gtwosubc(\uc;\varphi\ua)} F\O\dua\\
&= \bigl(\varphi_!F_*\O\bigr)\duc \in\N.
\end{split}\]
An inspection of \eqref{sigma-phi-o}, Definition \ref{def:phi-operad-left-adjoint}, and \ref{int:goperad-monoidal-functor} implies that the $\colorc$-colored $\Gtwo$-operad structures of $F_*\varphi_!\O$ and $\varphi_!F_*\O$ are identified under the above natural isomorphism.
\end{proof}

\begin{example}[Singular Chain Functor]\label{ex:singular-chain-functor}
Theorem \ref{change-base-category} applies to all the action operads in Example \ref{ex:all-augmented-group-operads}.  For instance, consider the \index{singular chain functor}singular chain functor \[C_\bullet : \CHau \to \Chain\] from the category of compactly generated weak Hausdorff spaces to the category of chain complexes of abelian groups.  This is a symmetric monoidal functor, so it induces a functor from the category of $\colorc$-colored ribbon operads in spaces to the category of $\colorc$-colored ribbon operads in chain complexes, and similarly for other action operads.  For example:
\begin{enumerate}
\item The \index{operad!little $n$-cube}\index{little $n$-cube operad}little $n$-cube operad $\C_n$ in Example \ref{ex:little-ncube}, the \index{operad!little $n$-disc}\index{little $n$-disc operad}little $n$-disc operad $\D_n$ in Example \ref{ex:little-n-disc}, the \index{operad!framed little $n$-disc}\index{framed little $n$-disc operad}framed little $n$-disc operad $\Df_n$ in Example \ref{ex:framed-little-n-disc}, and the \index{operad!phylogenetic}\index{phylogenetic operad}phylogenetic operad $\Phyl$ in Section \ref{sec:symmetric-operad-phylogenetic} are one-colored symmetric operads in spaces.  Applying the singular chain functor $C_\bullet$ level-wise, we obtain one-colored symmetric operads $C_{\bullet}\C_n$, $C_{\bullet}\D_n$, $C_{\bullet}\Df_n$, and $C_{\bullet}\Phyl$ in chain complexes.
\item The \index{planar tangle operad}\index{operad!planar tangle}planar tangle operad $\PTan$ in Definition \ref{def:planar-tangle-operad} is an $\SCir$-colored symmetric operad in spaces.  Applying the singular chain functor $C_\bullet$ entry-wise, we obtain an $\SCir$-colored symmetric operad $C_{\bullet}\PTan$ in chain complexes.
\item The universal cover\index{operad!universal cover of little $2$-cube} $\Ctilde_2$ of the little $2$-cube operad in Example \ref{ex:universal-cover-c2} and the universal cover $\Dtilde_2$ of the little $2$-disc operad in Example \ref{ex:universal-cover-d2} are one-colored braided operads in spaces.  Applying the singular chain functor $C_\bullet$ level-wise, we obtain one-colored braided operads $C_\bullet \Ctilde_2$ and $C_\bullet \Dtilde_2$ in chain complexes.
\item The universal cover\index{operad!universal cover of framed little $2$-disc} $\Dftildetwo$ of the framed little $2$-disc operad in Example \ref{ex:universal-cover-fd2} is a one-colored ribbon operad in spaces.  Applying the singular chain functor $C_\bullet$ level-wise, we obtain a one-colored ribbon operad $C_\bullet\Dftildetwo$ in chain complexes.\dqed
\end{enumerate}
\end{example}

\begin{example}[Nerve Functor]\label{ex:nerve-goperad}
For a group $G$, recall from Definition \ref{def:translation-category} that its translation category $E_G$ is the small groupoid whose object set is the underlying set of $G$, with a unique morphism between any two objects.  The translation category $E_G$ has a natural action by the group $G$.  For each action operad $\G$, applying the translation category construction level-wise, we obtain a one-colored $\G$-operad \[E_{\G}=\{E_{\G(n)}\}_{n\geq 0}\] in $\Cat$, the symmetric monoidal category of small categories.  The \index{nerve}nerve functor \[\Nerve : \Cat\to\Sset\] is a strong symmetric monoidal functor.  Applying the nerve functor level-wise, by Theorem \ref{changing-goperad} $\Nerve(E_{\G})$ is a one-colored $\G$-operad in the category of simplicial sets.  We already discussed this fact in Example \ref{ex:translation-g}.\dqed
\end{example}

\begin{example}\label{ex:set-to-m}
For each symmetric monoidal category $(\M,\otimes,\tensorunit)$, there is a strong symmetric monoidal functor $\iota : \Set \to \M$ defined by \[\iota(X) = \coprodover{X} \tensorunit.\] By Theorem \ref{change-base-category} there is an induced functor \[\nicexy{\gopcset \ar[r]^-{\iota_*} & \goperadcm}\] that sends each $\colorc$-colored $\G$-operad in $\Set$ to a $\colorc$-colored $\G$-operad in $\M$.  

Moreover, $\iota$ commutes with colimits because it is left adjoint to the\index{underlying set} functor \[\M(\tensorunit,-) : \M \to \Set.\]  For each morphism $\varphi : \Gone\to\Gtwo$ of action operads, by Theorem \ref{varphi-F}, in the diagram \[\nicexy{\goneopcset \ar[d]_-{\iota_*} \ar@<2pt>[r]^-{\varphi_!} & \gtwoopcset \ar@<2pt>[l]^-{\varphi^*} \ar[d]^-{\iota_*}\\
\goneopcm \ar@<2pt>[r]^-{\varphi_!} & \gtwoopcm \ar@<2pt>[l]^-{\varphi^*}}\] 
there are (i) an equality \[\iota_*\varphi^*=\varphi^*\iota_*\] and (ii) a natural isomorphism \[\varphi_!\iota_* \cong \iota_*\varphi_!.\]\dqed
\end{example}

\begin{example}[Underlying Objects]\label{ex:underlying-object-operad}
The \index{underlying object}underlying object set functor \[\Ob : \Cat \to \Set\] is strong symmetric monoidal.  By Theorem \ref{change-base-category} there is an induced functor \[\nicexy{\gopccat \ar[r]^-{\Ob_*} & \gopcset}\] that sends each $\colorc$-colored $\G$-operad in $\Cat$ to its underlying object $\colorc$-colored $\G$-operad in $\Set$.\dqed
\end{example}

\section{Change of Algebra Categories}\label{sec:change-algebra-category}

In this section, for a fixed action operad $\G$, we observe that each morphism of $\colorc$-colored $\G$-operads induces a morphism between the respective monads.  As a consequence, there is an induced adjunction between the respective categories of algebras.  

We recalled the concept of a monad in Definition \ref{def:monad}.  Let us now recall from Definition 4.5.8 in \cite{borceux2} the concept of a morphism of monads.

\begin{definition}\label{def:monad-morphism}
Suppose $(T,\mu,\epsilon)$ and $(T',\mu',\epsilon')$ are two monads in a category $\C$.  A \emph{morphism of monads}\index{monad!morphism}\index{morphism!monads} \[\theta : (T,\mu,\epsilon) \to (T',\mu',\epsilon')\] consists of a natural transformation $\theta : T \to T'$ such that the diagrams\[\nicexy{\Id_{\C} \ar[r]^-{\epsilon} \ar[dr]_-{\epsilon'} & T \ar[d]^-{\theta}\\ & T'} \qquad \nicexy{TT \ar[d]_-{\mu} \ar[r]^-{T\theta} & TT' \ar[r]^-{\theta_{T'}} & T'T' \ar[d]^-{\mu'}\\ T \ar[rr]^-{\theta} && T'}\]
are commutative.
\end{definition}

For each $\colorc$-colored $\G$-operad $\O$, recall from Definition \ref{def:g-operad-algebra} that the category $\algmo$ of $\O$-algebras is defined as the category of algebras over the monad in $\Mtoc$ with underlying functor $\O\circg -$, multiplication $\mu^{\O}$, and unit $\epsilon^{\O}$.

\begin{theorem}\label{changing-goperad}\index{morphism!G-operads@$\G$-operads}
Suppose\index{algebra!change of $\G$-operads} $f : \O\to\Q$ is a morphism of $\colorc$-colored $\G$-operads in $\M$.  Then the following statements hold.
\begin{enumerate}
\item There is an induced morphism \[\nicexy{\bigl(\O\circg -,\mu^{\O},\epsilon^{\O}\bigr) \ar[r]^-{f_*} & \bigl(\Q\circg -,\mu^{\Q},\epsilon^{\Q}\bigr)}\] of monads in $\Mtoc$ defined by the morphism
\begin{equation}\label{goperad-monad-morphism}
\nicexy@C+.5cm{(\O\circg X)_d = \dint^{\uc\in\gsubc} \O\duc\otimes X_{\uc} \ar[r]^-{\int^{\uc} f \otimes X_{\uc}} & \dint^{\uc\in\gsubc} \Q\duc\otimes X_{\uc} = (\Q\circg X)_d}\in\M
\end{equation}
for $X\in\Mtoc$ and $d\in\colorc$.
\item The morphism of monads $f_*$ in the previous item induces a functor \[\nicexy{\algmo & \algmq\ar[l]_-{f^*}}.\]
\item The functor $f^*$ in the previous item admits a left adjoint \[\nicexy{\algmo \ar[r]^-{f_!}& \algmq}.\]
\end{enumerate}
\end{theorem}

\begin{proof}
For the first assertion, note that \eqref{goperad-monad-morphism} defines a natural transformation $f_*$ from the functor $\O\circg-$ to the functor $\Q\circg-$.  To check the commutativity of the two diagrams in Definition \ref{def:monad-morphism} for $f_*$, note that the unit \[\epsilon^{\O}_X : X \to \O\circg X\] is defined entrywise as the composite \[\nicexy@C+.5cm{X_d \ar[r]^-{(\epsilon^{\O}_X)_d} \ar[d]_-{\cong} & (\O\circg X)_d\\ \tensorunit \otimes X_d \ar[r]^-{(\operadunit^{\O}_d,\Id)} & \O\dd\otimes X_d \ar[u]_-{\text{natural}}}\] for each $d\in\colorc$, where $\operadunit^{\O}_d$ is the $d$-colored unit of $\O$.  The commutativity of the first diagram in Definition \ref{def:monad-morphism} follows from the commutativity of the diagram \[\nicexy{\tensorunit \ar[r]^-{\operadunit^{\O}_d} \ar[dr]_-{\operadunit^{\Q}_d} & \O\dd \ar[d]^-{f}\\ & \Q\dd}\] which holds because $f$ is compatible with the colored units.

Similarly, the commutativity of the second diagram in Definition \ref{def:monad-morphism} follows from the commutativity of the diagram \[\nicexy{\O\circg\O \ar[d]_-{\mu^{\O}} \ar[r]^-{f \circg f} & \Q \circg \Q \ar[d]^-{\mu^{\Q}}\\ \O \ar[r]^-{f} & \Q}\] which holds because $f$ is a morphism of $\colorc$-colored $\G$-operads.

The second assertion follows from the first assertion and the fact that each morphism of monads induces a functor between the respective algebra categories in the opposite direction.  See, for example, Proposition 4.5.9 in \cite{borceux2}.

For the last assertion, consider the solid-arrow diagram \[\nicexy@C+.5cm{\algmo \ar@<2pt>[d]^-{U^{\O}} \ar@{-->}@<2pt>[r]^-{f_!} & \algmq \ar@<2pt>[l]^-{f^*} \ar@<2pt>[d]^-{U^{\Q}} \\
\Mtoc \ar@<2pt>[u]^-{\O\circg-} \ar@{=}[r] & \Mtoc \ar@<2pt>[u]^-{\Q\circg-}}\] with vertical free-forgetful adjunctions and \[U^{\O}f^*=U^{\Q}.\]  We observed in Corollary \ref{algebra-cat-bicomplete} that $\algmq$ is cocomplete.  By the Adjoint Lifting Theorem--see Section \ref{subsec:adjoint-lifting}--a left adjoint $f_!$ of $f^*$ exists.
\end{proof}

\chapter{Group Operads as Algebras}\label{ch:group-operad-category}

As before $(\M,\otimes,\tensorunit)$ is a fixed complete and cocomplete symmetric monoidal category whose monoidal product commutes with small colimits on each side.  Fix an action operad $(\G,\omega)$ as in Section \ref{sec:group-operad}.  The category $\goperadcm$ of $\colorc$-colored $\G$-operads in $\M$ was defined as the category of monoids in the monoidal category of $\colorc$-colored $\G$-sequences in $\M$ with respect to the $\colorc$-colored $\G$-circle product $\circg$.  There are two main purposes of this chapter:
\begin{enumerate}
\item We introduce the language of trees that is suitable for the study of $\G$-operads and $\infty$-$\G$-operads.  A $\G$-tree is a planar tree equipped with an element in $\G$ that acts as the input equivariance.  We describe the substitution procedure for $\G$-trees that will be used not only in this chapter but also in several subsequent chapters.  Indeed, the entire presheaf theory of $\G$-operads is based on the $\G$-tree category $\Treecatg$ in Chapter \ref{ch:group-tree-category}, whose objects are $\G$-trees and whose morphisms and composition involve $\G$-tree substitution.
\item Using the language of $\G$-trees, we give an explicit construction of the colored symmetric operad $\gopcm$ whose category of algebras is the category of $\colorc$-colored $\G$-operads in $\M$.
\end{enumerate}

Planar trees and their substitution are discussed in Section \ref{sec:planar-tree}.  In Section \ref{sec:g-tree} we introduce $\G$-trees and discuss their substitution.  In Section \ref{sec:symmetric-g-operad} we construct the colored symmetric operad $\gopcm$ in $\M$ whose category of algebras is the category $\goperadcm$.

\section{Planar Trees}\label{sec:planar-tree}

In this section, we discuss planar trees, which will serve two main purposes.  First, in Section \ref{sec:g-tree} we will equip planar trees with input equivariance coming from the action operad $\G$.  The symmetric operad whose algebras are $\colorc$-colored $\G$-operads is constructed from these $\G$-trees.  Second, when equipped with a suitable notion of morphisms, these $\G$-trees form a category, the presheaf category of which will give combinatorial models of $\infty$-$\G$-operads, as we will discuss in later chapters.  All the assertions in this section follow from simple inspection of definitions, so we will omit most of them.

Recall the concept of a tree in Definition \ref{def:tree-for-phyl}.

\begin{definition}\label{def:tree}
Suppose $T$ is a tree.
\begin{enumerate}
\item A \emph{planar structure} on $T$ is a choice of a bijection\label{not:ellv} \[\ell_v : \{1,\ldots,|\inp(v)|\} \iso \inp(v)\] for each vertex $v$ in $T$.  A planar structure on $T$ induces a \emph{canonical ordering}\label{not:ellt} \[\ell_T : \{1,\ldots,|\inp(T)|\} \iso \inp(T)\] on the set of inputs of $T$.  A \index{planar tree}\index{tree!planar}\emph{planar tree} is a tree with a chosen planar structure.  In this case, we regard $\inp(v)$ and $\inp(T)$ as ordered sets.
\item For a set $\colorc$, a \emph{$\colorc$-coloring}\index{tree!coloring} of $T$ is a function\label{not:coloring} \[\kappa : T \to \colorc,\] and $\kappa(e) \in \colorc$ is called the color of $e \in T$.  A \emph{$\colorc$-colored tree}\index{colored tree} is a tree with a chosen $\colorc$-coloring.
\item Suppose $T$ is equipped with a planar structure and a $\colorc$-coloring $\kappa : T \to \colorc$.  For each $u \in \{T\} \amalg \Vt(T)$, we call
\begin{itemize}
\item $\kappa\,\out(u) \in \colorc$ the \index{output color}\emph{output color},
\item $\ell_u(j) \in \inp(u)$ the\index{input} \emph{$j$th input},
\item $\kappa\,\inp(u) \in \Profc$ the\index{input profile} \emph{input profile}, and
\item the pair\label{not:profofu} \[\profofu=\sbinom{\kappa\,\out(u)}{\kappa\,\inp(u)} \in \Profcc\] the \emph{profile}\index{profile} of $u$. 
\end{itemize}
\item A \emph{vertex ordering}\index{vertex ordering} on $T$ is a bijection \[\rho : \{1,\ldots,|\Vt(T)|\} \iso \Vt(T),\] with respect to which we regard $\Vt(T)$ as an ordered set and call $\rho(j)$ the \emph{$j$th vertex}.  A planar structure on $T$ induces a \emph{canonical vertex ordering}\index{canonical vertex ordering} on $T$ in which the initial vertex of the root $r$ is ordered first, followed by the initial vertices of $f \in \inp(r)$, followed by the initial vertices of $g \in \inp(f)$ with $f \in \inp(r)$, etc.
\item An \emph{isomorphism}\index{tree!isomorphism} $\psi : T \iso T'$ of trees is an isomorphism of partially ordered sets that preserves the inputs and any extra structure (e.g., planar structure, $\colorc$-coloring, and/or vertex ordering) that $T$ and $T'$ both possess.
\end{enumerate}
\end{definition}

\begin{remark}In what follows, we will mostly be concerned with isomorphism classes of trees.  We will use the same symbol to denote a tree and its isomorphism class and omit mentioning isomorphism classes.  It will sometimes be convenient to regard a tree as a $\{*\}$-colored tree, in which every edge is assigned the same color $*$.\dqed
\end{remark}

Next are several examples of basic planar trees.

\begin{example}[Exceptional Edges]\label{ex:exceptional-edge}
The \index{exceptional edge}\emph{exceptional edge}, denoted $\uparrow$, is the planar tree with a single edge that is both the root and an input and that has no vertices.  If it is equipped with a $\colorc$-coloring with $\kappa(\uparrow)=c \in \colorc$, then we call it the \emph{$c$-colored exceptional edge}, denoted $\uparrow_c$.  Its profile is $\cc \in \Profcc$.\dqed
\end{example}

\begin{example}[Corollas]\label{ex:corolla}
For $n\geq 0$, the\index{corolla} \emph{$n$-corolla} \[\Cor_n=\{r,i_1,\ldots,i_n\}\]is the planar tree with $n+1$ edges, root $r$, $n$ inputs \[\inp(\Cor_n) = \{i_1,\ldots,i_n\},\] and a unique vertex $v$ containing all $n+1$ edges.  The planar structure is given by \[\ell_v(j) = i_j \forspace 1 \leq j \leq n.\]  The $n$-corolla has no internal edges.  We visualize the $n$-corolla as 
\begin{center}\begin{tikzpicture}
\node[plain] (v) {$v$}; \node[below=.1cm of v] () {$\cdots$};
\draw[outputleg] (v) to node[at end]{\scriptsize{$r$}} +(0cm,.8cm);
\draw[thick] (v) to node[swap, at end, outer sep=-2pt]{\scriptsize{$i_1$}} +(-.8cm,-.6cm);
\draw[thick] (v) to node[at end, outer sep=-2pt]{\scriptsize{$i_n$}} +(.8cm,-.6cm);
\end{tikzpicture}\end{center}
with the inputs drawn from left to right according to their orders.  If the $n$-corolla is equipped with a $\colorc$-coloring with \[\kappa(r) = d \andspace \kappa(i_j) = c_j \in \colorc,\] then its profile is $\duc \in \Profcc$, and we call it the \emph{$\duc$-corolla}, denoted $\Cor_{(\uc;d)}$, where \[\uc=(c_1,\ldots,c_n) \in \Profc.\]  Observe that we draw the root as an arrow, each non-root edge as a line, and each vertex as a circle.  We also orient the pictures from bottom to top, so inputs are at the bottom, and the root is at the top.\dqed
\end{example}

\begin{example}[$2$-Level Trees]\label{ex:2-level-tree}
For $n \geq 1$ and $k_1,\ldots,k_n \geq 0$, there is a planar tree $T$, called a \index{tree!two-level}\emph{$2$-level tree}, that can be pictorially represented as follows.
\begin{center}\begin{tikzpicture}
\matrix[row sep=.1cm, column sep=1.2cm]{
& \node [plain, label=below:$...$] (v) {$v$}; &\\
\node [plain, label=below:$...$] (u1) {$u_1$}; &&
\node [plain, label=below:$...$] (um) {$u_n$};\\};
\draw [outputleg] (v) to node[at end]{\scriptsize{$r$}} +(0,.8cm);
\draw [thick] (u1) to node{\scriptsize{$t_1$}} (v);
\draw [thick] (um) to node[swap]{\scriptsize{$t_n$}} (v);
\draw [thick] (u1) to node[below left=.2cm]{\scriptsize{$i_{1,1}$}} +(-.8cm,-.6cm);
\draw [thick] (u1) to node[below right=.2cm]{\scriptsize{$i_{1,k_1}$}} +(.8cm,-.6cm);
\draw [thick] (um) to node[below left=.2cm]{\scriptsize{$i_{n,1}$}} +(-.8cm,-.6cm);
\draw [thick] (um) to node[below right=.2cm]{\scriptsize{$i_{n,k_n}$}} +(.8cm,-.6cm);
\end{tikzpicture}
\end{center}
Formally $T$ is defined as follows.
\begin{itemize}
\item $T = \Bigl\{r,\{t_j\}_{1\leq j \leq n}, \{i_{j,l}\}_{1\leq j \leq n}^{1 \leq l \leq k_j}\Bigr\}$ with root $r$, inputs $\inp(T) = \Bigl\{\{i_{j,l}\}_{1\leq j \leq n}^{1 \leq l \leq k_j}\Bigr\}$, and $n$ internal edges $\{t_j\}_{1 \leq j \leq n}$.
\item There are $n+1$ vertices:
\[v=\bigl\{r,t_1,\ldots,t_n\bigr\} \andspace 
u_j=\bigl\{t_j,i_{j,1},\ldots,i_{j,k_j}\bigr\} \forspace 1 \leq j \leq n.\]
\item The planar structure is defined as \[\ell_v(j)=t_j \andspace \ell_{u_j}(l) = i_{j,l}\] for $1 \leq j \leq n$ and $1 \leq l \leq k_j$.
\end{itemize}
Under the canonical vertex ordering, the vertex $v$ is ordered first, followed by $u_1,\ldots,u_n$.  

Suppose $T$ has a $\colorc$-coloring such that $\kappa(r)=d$, $\kappa(t_j) = c_j$, and $\kappa(i_{j,l}) = b_{j,l}$ in $\colorc$.  Then we denote $T$ by\label{not:tbcd} \[T\bigl(\{\ub_j\};\uc;d\bigr)\] with $\uc = (c_1,\ldots,c_n)$ and $\ub_j=(b_{j,1},\ldots,b_{j,k_j})$ in $\Profc$.  Its profile is $\dub \in \Profcc$, where $\ub = (\ub_1,\ldots,\ub_n)$ is the concatenation of the $\ub_j$'s.\dqed
\end{example}

Next we define tree substitution at a vertex, which gives rise to $\compi$-composition of planar trees.

\begin{definition}[Tree Substitution at a Vertex]\label{def:planar-tree-compv}
Suppose $T$ and $T'$ are planar $\colorc$-colored trees, and $v \in \Vt(T)$ such that $\profofv = \Prof(T')$ as $\colorc$-profiles.  The \index{tree substitution!at a vertex}\emph{tree substitution at $v$}, denoted by $T\comp_v T'$, is defined as the planar $\colorc$-colored tree\label{not:compv} \[T\comp_v T' = \frac{T \amalg T'}{\bigl\{\out(v)\equiv \out(T'),\, \inp(v) \equiv\inp(T')\bigr\}}.\]  In this quotient, the outgoing edge $\out(v)$ of $v$ is identified with the root $\out(T')$ of $T'$, and the incoming edges of $v$ are identified with the inputs of $T'$.  Its planar $\colorc$-colored tree structure is uniquely induced by those in $T$ and $T'$.  In this case, we also say that the tree substitution \emph{occurs at $v$}.  If both $T$ and $T'$ are equipped with vertex orderings, then the tree substitution at $v$ is given the induced vertex ordering.
\end{definition}

\begin{lemma}\label{lem:planar-tree-sub-vertex}
In the context of Definition \ref{def:planar-tree-compv}, the following statements hold.
\begin{enumerate}
\item The root and the inputs of $T \comp_v T'$ are those of $T$, and \[\Prof(T \comp_v T') = \profoft.\]
\item $\Vt(T \comp_v T') = \bigl[\Vt(T) \setminus \{v\}\bigr] \amalg \Vt(T')$.
\item Internal edges in $T'$ remain internal edges in $T \comp_v T'$.
\item Left Unity: If $T = \Cor_{\profoft}$, then \[\Cor_{\profoft} \comp_v T' = T'.\]
\item Right Unity: If $T' = \Cor_{\profofv}$, then \[T \comp_v \Cor_{\profofv} = T.\]
\item Vertical Associativity: Suppose $u \in \Vt(T')$ and $T''$ is a planar $\colorc$-colored tree with $\profofu = \Prof(T'')$.  Then \[\bigl(T\comp_v T'\bigr) \comp_u T'' = T\comp_v \bigl(T' \comp_u T''\bigr).\]
\item Horizontal Associativity: Suppose $v\not= w \in \Vt(T)$ and $T''$ is a planar $\colorc$-colored tree with $\profofw = \Prof(T'')$.  Then \[\bigl(T\comp_v T'\bigr) \comp_w T'' = \bigl(T\comp_w T''\bigr) \comp_v T'.\]
\end{enumerate}
\end{lemma}

\begin{proof}
As an example, let us prove the horizontal associativity.  By definition both sides are equal to the planar $\colorc$-colored tree \[\frac{T \amalg T' \amalg T''}{\bigl\{\sbinom{\out(v)}{\inp(v)} \equiv \sbinom{\out(T')}{\inp(T')},\, \sbinom{\out(w)}{\inp(w)} \equiv \sbinom{\out(T'')}{\inp(T'')}\bigr\}}\] with structures induced by those in $T$, $T'$, and $T''$.
\end{proof}

Instead of using only one vertex at a time, we can also substitute planar trees into all the vertices at the same time.  This yields the following operation on planar trees.

\begin{definition}[Tree Substitution]\label{def:planar-tree-sub}
Suppose $T$ is a planar $\colorc$-colored tree, and for each vertex $v \in \Vt(T)$, $T_v$ is a planar $\colorc$-colored tree such that $\profofv = \Prof(T_v)$.  Define the \index{tree substitution}\emph{tree substitution}\label{not:ttv} \[T(T_v)_{v\in \Vt(T)} =  \Bigl((T \comp_{v_1} T_{v_1}) \comp_{v_2} \cdots \Bigr) \comp_{v_n} T_{v_n},\] where $(v_1,\ldots,v_n)$ is an arbitrary choice of an ordering of $\Vt(T)$.  To shorten the notation, we will sometimes abbreviate the tree substitution to $T(T_v)_{v}$ or $T(T_v)$.
\end{definition}

The following statements follow immediately from Lemma \ref{lem:planar-tree-sub-vertex}.

\begin{lemma}\label{planar-tree-sub-associative} 
In the context of Definition \ref{def:planar-tree-sub}, the following statements hold.
\begin{enumerate}
\item The tree substitution is independent of the choice of an ordering of $\Vt(T)$.
\item $\Prof\bigl(T(T_v)_{v\in \Vt(T)} \bigr) = \profoft$.\index{profile!tree substitution}
\item $\Vt\bigl(T(T_v)_{v\in\Vt(T)}\bigr) = \coprod_{v\in\Vt(T)} \Vt(T_v)$.\index{vertex!tree substitution}
\item Internal edges in each $T_v$ remain internal edges in the tree substitution.
\item Unity:\index{unity!tree substitution} $T(\Cor_{\profofv})_{v\in\Vt(T)} = T$.
\item Associativity:\index{associativity!tree substitution}  Suppose for each $v \in \Vt(T)$ and each $u \in \Vt(T_v)$, $T_v^u$ is a planar $\colorc$-colored tree such that $\Prof(u) = \Prof(T_v^u)$.  Then there is an equality \[\Bigl[T(T_v)_{v\in\Vt(T)}\Bigr](T_v^u)_{v\in\Vt(T),u\in\Vt(T_v)} = T\Bigl[T_v(T_v^u)_{u\in\Vt(T_v)}\Bigr]_{v\in\Vt(T)}.\]  
\end{enumerate}
\end{lemma}

\begin{proof}
For example, the first assertion follows from the horizontal associativity property in Lemma \ref{lem:planar-tree-sub-vertex} by an induction.  Similarly, the last assertion about associativity of tree substitution follows from the vertical associativity and the horizontal associativity properties in Lemma \ref{lem:planar-tree-sub-vertex}.
\end{proof}

\begin{example}\label{ex:treesub}
Consider the tree substitution $K = T(T_u,T_v,T_w)$ as in the following picture.
\begin{center}\begin{tikzpicture}
\node[plain] (w) {$w$}; \node[above=.7cm of w] (T) {$T$};
\node[below right=.7cm of w, plain] (v) {$v$}; \node[below left=.7cm of w, plain] (u) {$u$};
\draw[outputleg] (w) to node[pos=.9]{\scriptsize{$e$}} +(0,.9cm);
\draw[thick] (u) to node{\scriptsize{$c$}} (w);
\draw[thick] (v) to node[swap]{\scriptsize{$d$}} (w);
\draw[thick] (u) to node[swap, outer sep=-2pt]{\scriptsize{$a$}} +(-.6cm,-.9cm);
\draw[thick] (u) to node[pos=.9, outer sep=-2pt, swap]{\scriptsize{$b$}} +(0cm,-.9cm);
\draw[thick] (v) to node{\scriptsize{$d$}} +(0cm,-.9cm);
\node[left=3cm of T,plain] (w2) {$w_2$}; \node[above=.5cm of w2,plain] (w1) {$w_1$};
\node[left=.6cm of w2] (Hw) {$T_w$};
\draw[outputleg] (w1) to node{\scriptsize{$e$}} +(0,.7cm);
\draw[thick] (w1) to node[swap, outer sep=-2pt]{\scriptsize{$c$}} +(-.7cm,-.5cm);
\draw[thick] (w2) to node[swap]{\scriptsize{$g$}} (w1);
\draw[thick] (w2) to node[swap]{\scriptsize{$d$}} +(0,-.7cm);
\node[draw=gray!50,inner sep=15pt,ultra thick,rounded corners,
fit=(w1) (w2)] (Hwbox) {};
\draw [gray!50, ->, dotted, line width=2pt, shorten >=0.03cm, bend left=10] (Hwbox) to (w);
\node[left=2cm of u,plain] (u1) {$u_1$}; \node[below=.5cm of u1,plain] (u2) {$u_2$};
\node[left=.6cm of u1] {$T_u$};
\draw[outputleg] (u1) to node{\scriptsize{$c$}} +(0,.7cm);
\draw[thick] (u1) to node[swap, outer sep=-2pt]{\scriptsize{$a$}} +(-.7cm,-.3cm);
\draw[thick] (u1) to node[outer sep=-2pt]{\scriptsize{$b$}} +(-.5cm,-.6cm);
\draw[thick] (u2) to node[swap]{\scriptsize{$f$}} (u1);
\node[draw=gray!50,inner sep=15pt,ultra thick,rounded corners,
fit=(u1) (u2)] (Hubox) {};
\draw [gray!50, ->, dotted, line width=2pt, shorten >=0.03cm, bend left=10] (Hubox) to (u);
\node[below=1.5cm of w] (Hv) {$\uparrow_d$}; \node[below=.2cm of Hv] () {$T_v$}; 
\node[draw=gray!50,inner sep=4pt,ultra thick,rounded corners,
fit=(Hv)] (Hvbox) {};
\draw [gray!50, ->, dotted, line width=2pt, shorten >=0.03cm, bend left=15] (Hvbox) to (v);
\node[right=3.5cm of w,plain] (w1') {$w_1$}; \node[above=.7cm of w1'] {$K$};
\node[below right=.7cm of w1', plain] (w2') {$w_2$}; 
\node[below left=.7cm of w1', plain] (u1') {$u_1$};
\node[below right=.7cm of u1',plain] (u2') {$u_2$};
\draw[outputleg] (w1') to node{\scriptsize{$e$}} +(0,.9cm);
\draw[thick] (w2') to node[swap]{\scriptsize{$g$}} (w1');
\draw[thick] (u1') to node{\scriptsize{$c$}} (w1');
\draw[thick] (u2') to node[swap]{\scriptsize{$f$}} (u1');
\draw[thick] (u1') to node[swap, outer sep=-2pt]{\scriptsize{$a$}} +(-.6cm,-.9cm);
\draw[thick] (u1') to node[pos=.9, outer sep=-2pt, swap]{\scriptsize{$b$}} +(0cm,-.9cm);\draw[thick] (w2') to node{\scriptsize{$d$}} +(0,-.9cm);
\end{tikzpicture}
\end{center}
Each gray dotted arrow represents a tree substitution at the indicated vertex.  Observe that 
\[\begin{split}\profofu&=\sbinom{c}{a,b}=\Prof(T_u),\quad \profofv=\dd=\Prof(T_v),\\
\profofw&=\sbinom{e}{c,d}=\Prof(T_w),\andspace \profoft=\sbinom{e}{a,b,d}=\profofk.\end{split}\]
Internal edges in $T_u$, $T_v$, and $T_w$--namely, $f$ and $g$--become internal edges in the tree substitution $K$.  The $d$-colored internal edge in $T$ is no longer an internal edge in $K$ because $T_v$ is the $d$-colored exceptional edge.

Suppose given the vertex orderings: $(u,v,w)$ for $T$, $(u_1,u_2)$ for $T_u$, and $(w_1,w_2)$ for $T_w$.  Then the induced vertex ordering for $K$ is $(u_1,u_2,w_1,w_2)$.\dqed
\end{example}

\section{Group Trees}\label{sec:g-tree}

Fix an action operad $(\G,\omega)$ as in Definition \ref{def:augmented-group-operad}.  In this section, we consider $\G$-versions of planar trees and tree substitution.  Recall that for an element $g \in \G(n)$, its underlying permutation $\omega(g) \in S_n$ is denoted by $\gbar$.  

\begin{definition}\label{def:g-tree}
Suppose $\colorc$ is a set.
\begin{enumerate}
\item A\index{G-tree@$\G$-tree}\index{tree!G@$\G$-} \emph{$\colorc$-colored $\G$-tree} is a pair\label{not:tsigma} $(T,\sigma)$ consisting of 
\begin{enumerate}[label=(\roman*)]
\item a planar $\colorc$-colored tree $T$, called the \index{underlying planar  tree}\emph{underlying planar $\colorc$-colored tree}, and 
\item an element $\sigma \in \G(|\inp(T)|)$, called the \index{input equivariance}\emph{input equivariance}.
\end{enumerate}
We call 
\begin{itemize}
\item $\bigl(\kappa\,\inp(T)\bigr)\sigmabar \in \Profc$ the \emph{input profile}, and
\item the pair\label{not:proftsigma} \[\Prof(T,\sigma) = \sbinom{\kappa\,\out(T)}{(\kappa\,\inp(T))\sigmabar} \in \Profcc\] the \emph{profile}\index{profile!G-tree@$\G$-tree} of $(T,\sigma)$.
\end{itemize}
We will often omit mentioning $\colorc$ and just call $(T,\sigma)$ a $\G$-tree.  A \emph{group tree}\index{group tree}\index{tree!group} is a $\G$-tree for some action operad $\G$.
\item An \index{isomorphism of G-trees@isomorphism of $\G$-trees}\emph{isomorphism} $\psi : (T,\sigma) \iso (T',\sigma)$ of $\G$-trees is an isomorphism of planar $\colorc$-colored trees $T \iso T'$ with the same input equivariance $\sigma \in \G(|\inp(T)|)$.
\item If $\G$ is the symmetric group operad $\S$ in Example \ref{ex:symmetric-group-operad},  the braid group operad $\B$ in Definition \ref{def:braid-group-operad}, the ribbon group operad $\R$ in Definition \ref{def:ribbon-group-operad}, or the cactus group operad $\Cac$ in Definition \ref{def:cactus-group-operad}, then a $\G$-tree is called a \index{symmetric tree}\emph{symmetric tree}, a \index{braided tree}\emph{braided tree}, a \index{ribbon tree}\emph{ribbon tree}, or a \index{cactus tree}\emph{cactus tree}, respectively.
\end{enumerate}
\end{definition}

\begin{remark} For the planar group operad $\P$ in Example \ref{ex:trivial-group-operad}, a $\P$-tree is exactly a planar $\colorc$-colored tree.  As before, we will mostly be working with isomorphism classes of $\G$-tree, so we will omit mentioning isomorphism classes from now on.\dqed\end{remark}

\begin{ginterpretation} For a $\G$-tree $(T,\sigma)$, $\inp(T)$ is equipped with the canonical ordering coming from the planar structure of $T$.  We think of the input equivariance $\sigma \in \G(|\inp(T)|)$ as acting on the inputs of $T$.  The input profile \[(\kappa\,\inp(T))\sigmabar \in \Profc\] of $(T,\sigma)$ is the result of the underlying permutation $\sigmabar \in S_{|\inp(T)|}$ acting on the input profile $\kappa\,\inp(T) \in \Profc$ of the planar tree $T$.  We can visualize a $\G$-tree $(T,\sigma)$ as 
\begin{center}\begin{tikzpicture}
\node[triangular] (T) {$T$}; \node[below=.15cm of T] () {\scriptsize{$\inp(T)$}};
\draw[outputleg] (T) to node[at end]{} +(0cm,1.1cm);
\node[below =1cm of T, bigplain] (s) {$\sigma$}; 
\node[above=.1cm of s] () {$\cdots$};
\node[below=.04cm of s] () {$\cdots$};
\node[below=.25cm of s] () {\scriptsize{$\inp(T)\sigmabar$}};
\draw[thick] ([xshift=.05cm]s.north west) -- ([xshift=-.1cm]T.south west); 
\draw[thick] ([xshift=-.05cm]s.north east) -- ([xshift=.1cm]T.south east);
\draw[thick] ([xshift=.05cm]s.south west) to +(0cm,-.6cm); 
\draw[thick] ([xshift=-.05cm]s.south east) to +(0cm,-.6cm);
\end{tikzpicture}\end{center}
in which the planar tree $T$ is represented as a triangle and the input equivariance $\sigma \in \G(|\inp(T)|)$ is the rectangular box under the inputs of $T$.\dqed
\end{ginterpretation}

\begin{example}[Symmetric Trees]\label{ex:symmetric-tree}
For the symmetric group operad $\S$ in Example \ref{ex:symmetric-group-operad}, a symmetric tree is a pair $(T,\sigma)$ with $T$ a planar $\colorc$-colored tree and $\sigma \in S_{|\inp(T)|}$ a permutation.  For instance, for the $\duc$-corolla $\Cor_{(\uc;d)}$ in Example \ref{ex:corolla} and a permutation $\sigma \in S_{|\uc|}$, the symmetric tree $(\Cor_{(\uc;d)}, \sigma)$ is called a \index{permuted corolla}\index{corolla!permuted}\emph{permuted corolla} in \cite{yau-operad,bluemonster}.    For another example, consider the symmetric tree $(T,\sigma)$
\begin{center}\begin{tikzpicture}[scale=1.3]
\node[plain] (w) {$w$}; \node[left=1 of w] (T) {$T$}; \node[below=1.25 of T] (sigma) {$\sigma$};
\node[below=2 of w] (c1) {};
\node[left=.7 of c1] (c2) {};\node[right=.7 of c1] (c3) {};
\node[above=1 of c2, plain] (u) {$u$}; \node[above=1 of c3, plain] (v) {$v$}; 
\draw[outputleg] (w) to +(0,.5);
\draw[thick] (u) to (w); \draw[thick] (v) to (w);
\draw[thick] (c2) to[out=90, in=335] (u);
\draw[thick] (c1) to[out=90, in=270] (u);
\draw[thick] (c3) to (v);
\end{tikzpicture}\end{center}
in which $T$ is the planar tree in Example \ref{ex:treesub} and the input equivariance is the adjacent transposition $\sigma = (1~2) \in S_3$.\dqed
\end{example}

\begin{example}[Braided Trees]\label{ex:braided-tree}
For the braid group operad $\B$ in Definition \ref{def:braid-group-operad} (resp., the pure braid group operad $\PB$ in Definition \ref{def:pure-braid-group-operad}), a (pure) braided tree is a pair $(T,\sigma)$ with $T$ a planar $\colorc$-colored tree and $\sigma \in B_{|\inp(T)|}$ a braid (resp., $\sigma \in PB_{|\inp(T)|}$).  For the $\duc$-corolla $\Cor_{(\uc;d)}$ and a braid $\sigma \in B_{|\uc|}$, the braided tree $(\Cor_{(\uc;d)}, \sigma)$ is called a \index{braided corolla}\index{corolla!braided}\emph{braided corolla}.  For another example, consider the braided tree $(T,\sigma)$
\begin{center}\begin{tikzpicture}[scale=1.3]
\node[plain] (w) {$w$}; \node[left=1.5 of w] (T) {$T$}; \node[below=1.75 of T] (sigma) {$\sigma$};
\node[below=3 of w] (c1) {};\node[left=.7 of c1] (c2) {};\node[right=.7 of c1] (c3) {};
\node[above=.7 of c1] (i2) {};\node[above=.7 of c3] (i3) {};
\node[above=2 of c3, plain] (v) {$v$}; \node[above=2 of c2, plain] (u) {$u$};
\draw[outputleg] (w) to node[pos=.8]{\scriptsize{$e$}} +(0,.7);
\draw[thick] (u) to node{\scriptsize{$c$}} (w); 
\draw[thick] (v) to node[swap]{\scriptsize{$d$}} (w);
\draw[thick] (c1) to[out=120, in=270] node[pos=.9]{\scriptsize{$a$}} (u);
\draw[line width=6pt, white, shorten <=.5pt, shorten >=.5pt] (c2) to[out=100, in=270] (i3);
\draw[thick, shorten >=-8pt] (c2) to[out=80, in=270] (i3);
\draw[line width=4pt, white, shorten <=.5pt, shorten >=.5pt] (c3) to[out=80, in=270] (i2);
\draw[thick, shorten >=-8pt] (c3) to[out=100, in=270] (i2);
\draw[thick] (i2) to[out=90, in=270] node[swap, pos=.95]{\scriptsize{$d$}} (v);
\draw[line width=6pt, white, shorten <=.5pt, shorten >=.5pt] (i3) to[out=90, in=315] (u);
\draw[thick] (i3) to[out=90, in=315] node[swap, pos=.9]{\scriptsize{$b$}} (u);
\end{tikzpicture}\end{center}
in which $T$ is the planar tree in Example \ref{ex:treesub} and the input equivariance is the braid \[\sigma = s_2s_2s_1^{-1} \in B_3.\]  Note that the underlying permutation of $\sigma$ is the adjacent transposition $\sigmabar = (1~2)$.\dqed
\end{example}

\begin{example}[Ribbon Trees]\label{ex:ribbon-tree}
For the ribbon group operad $\R$ in Definition \ref{def:ribbon-group-operad} (resp., the pure ribbon group operad $\PR$ in Definition \ref{def:pure-ribbon-group-operad}), a (pure) ribbon tree is a pair $(T,r)$ with $T$ a planar $\colorc$-colored tree and $r \in R_{|\inp(T)|}$ a ribbon (resp., $r\in PR_{|\inp(T)|}$).  For the $\duc$-corolla $\Cor_{(\uc;d)}$ and a ribbon $r \in R_{|\uc|}$, the ribbon tree $(\Cor_{(\uc;d)},r)$ is called a \index{ribbon corolla}\index{corolla!ribbon}\emph{ribbon corolla}.\dqed
\end{example}

\begin{example}[Cactus Trees]\label{ex:cactus-tree}
For the cactus group operad $\Cac$ in Definition \ref{def:cactus-group-operad} (resp., the pure cactus group operad $\PCac$ in Definition \ref{def:pure-cactus-group-operad}), a (pure) cactus tree is a pair $(T,r)$ with $T$ a planar $\colorc$-colored tree and $r \in Cac_{|\inp(T)|}$ a cactus (resp., $r\in PCac_{|\inp(T)|}$).  For the $\duc$-corolla $\Cor_{(\uc;d)}$ and a cactus $r \in Cac_{|\uc|}$, the cactus tree $(\Cor_{(\uc;d)},r)$ is called a \index{cactus corolla}\index{corolla!cactus}\emph{cactus corolla}.\dqed
\end{example}

Next we discuss tree substitution for $\G$-trees.

\begin{definition}[Group Tree Substitution at a Vertex]\label{def:treesub-g}
Suppose $(T,\sigma)$ is a $\G$-tree with $v \in \Vt(T)$, and $(T_v,\sigma_v)$ is a $\G$-tree with $\profofv = \Prof(T_v,\sigma_v)$.  The \index{G-tree@$\G$-tree!substitution at a vertex}\emph{$\G$-tree substitution at $v$}, denoted \[(T,\sigma) \comp_v (T_v,\sigma_v) = \bigl(T' \comp_{v'} T_v, \sigma_v \comp^v \sigma\bigr),\] is defined as follows.
\begin{enumerate}[label=(\roman*)]
\item Suppose $T'$ is the planar $\colorc$-colored tree obtained from $T$ by changing the ordering of $\inp(v)$ to the composite \[\nicexy{\{1,\ldots,|\inp(v)|\} \ar[r]^-{\sigmavbar^{-1}}_-{\cong} 
& \{1,\ldots,|\inp(v)|\} \ar[r]^-{\ell_v}_-{\cong} 
& \inp(v) \ar@{<-}`u[ll] `[ll]_-{\ell_{v'}} [ll]}.\]  To avoid confusion, the vertex $v$ in $T'$ with this new ordering $\ell_{v'}$ is denoted by $v'$.  Its input profile is 
\[\begin{split}
\kappa\,\inp(v') &= \bigl(\kappa\,\inp(v)\bigr)\sigmavbar^{-1}\\
&= \bigl(\kappa\, \inp(T_v,\sigma_v)\bigr)\sigmavbar^{-1}\\ 
&= \bigl(\kappa\, \inp(T_v)\sigmavbar\bigr)\sigmavbar^{-1}\\
&= \kappa\, \inp(T_v).\end{split}\] 
Since $\Prof(v') = \Prof(T_v)$, we can form the tree substitution $T' \comp_{v'} T_v$ at $v'$ as in Definition \ref{def:planar-tree-compv}.  This is the planar $\colorc$-colored tree for the desired $\G$-tree substitution at $v$.
\item To define the input equivariance \[\sigma_v \comp^v \sigma \in \G(|\inp(T)|),\] for each $1 \leq j \leq |\inp(v)|$, first define \[k_v^j = \left|\Bigl\{e \in \inp(T) : \ell_v\bigl(\sigmavbar^{-1}(j)\bigr) \leq e\Bigr\}\right|\] as the number of inputs of $T$ that are equal to (in which case $k_v^j=1$) or descended from the $\sigmavbar^{-1}(j)$th input of $v$.  Then define 
\[\begin{split}
k_v &= k_v^1 + \cdots + k_v^{|\inp(v)|},\\
k_v^- &= \left|\biggl\{e \in \inp(T) : \parbox{6cm}{$\exists v_e \in \Vt(T), f_1,f_2 \in \inp(v_e)$ such that $\ell_{v_e}(f_1) < \ell_{v_e}(f_2)$, $f_1 \leq e$, $f_2 \leq \out(v)$}\biggr\}\right|,\\
k_v^+ &= \left|\biggl\{e \in \inp(T) : \parbox{6cm}{$\exists v_e \in \Vt(T), f_1,f_2 \in \inp(v_e)$ such that $\ell_{v_e}(f_1) < \ell_{v_e}(f_2)$, $f_2 \leq e$, $f_1 \leq \out(v)$}\biggr\}\right|.
\end{split}\] So $k_v^-$ and $k_v^+$ are the numbers of inputs of $T$ that are ordered, under the canonical ordering, before and after those descended from $\out(v)$, respectively.  Similarly, $k_v$ is the number of inputs of $T$ that are descended from $\out(v)$.  Note that \[|\inp(T)| = k_v^- + k_v^+  + k_v.\] Using the notation in \eqref{g-block-perm} and \eqref{g-block-sum}, we define
\begin{equation}\label{sigmav-compv-sigma}
\sigma_v \comp^v \sigma = \Bigl(\id_1^{\oplus {k_v^-}} \oplus \sigma_v\langle k_v^1,\ldots,k_v^{|\inp(v)|}\rangle \oplus \id_1^{\oplus {k_v^+}}\Bigr) \cdot \sigma \in \G(|\inp(T)|)
\end{equation}
with the product taken in $\G(|\inp(T)|)$.  This is the input equivariance for the desired $\G$-tree substitution at $v$.
\end{enumerate}
\end{definition}

\begin{interpretation} In order to substitute $T_v$ into $T$ at $v$, first we need to make sure that the profiles of $T_v$ and $v$ agree, which they do up to the underlying permutation $\sigmabar_v$.  So we (i) permute the inputs of $v$ from the left using the underlying permutation of $\sigma_v$ and (ii) form the tree substitution $T' \comp_{v'} T_v$ at $v'$.  In terms of $\G$, the effect on the inputs of $T$ is the element 
\[\begin{split}
&\gammag\Bigl(\id_{k_v^- + k_v^+ + 1}; \overbrace{\id_1,\ldots,\id_1}^{k_v^-}, \gammag\bigl(\sigma_v; \id_{k_v^1},\ldots, \id_{k_v^{|\inp(v)|}}\bigr), \overbrace{\id_1,\ldots, \id_1}^{k_v^+} \Bigr)\\
&= \id_1^{\oplus {k_v^-}} \oplus \sigma_v\langle k_v^1,\ldots,k_v^{|\inp(v)|}\rangle \oplus \id_1^{\oplus {k_v^+}} \in \G(|\inp(T)|)\end{split}\] 
that appeared in \eqref{sigmav-compv-sigma}.\dqed
\end{interpretation}

\begin{example}
Suppose $\G$ is the symmetric group operad $\S$ in Example \ref{ex:symmetric-group-operad}.  Consider the planar tree $T$
\begin{center}\begin{tikzpicture}
\matrix[row sep=.5cm, column sep=.9cm]{
& \node [plain] (w) {$y$}; &&\\
\node [plain] (v1) {$w$}; && \node [plain] (v2) {$x$}; &\\
& \node [plain] (u1) {$v$}; && \\
& \node [plain] (u2) {$u$}; &&\\};
\draw [outputleg] (w) to +(0,.8cm); \draw [thick] (w) to +(0cm,-.9cm);
\draw [thick] (v1) to (w); \draw [thick] (v2) to (w);
\draw [thick] (u1) to node[inner sep=-1pt]{\scriptsize{$\out(v)$}} (v2); \draw [thick] (u2) to (u1);
\draw [thick] (v1) to +(-.3cm,-.9cm); \draw [thick] (v1) to +(.3cm,-.9cm);
\draw [thick] (v2) to  +(0cm,-.9cm); \draw [thick] (v2) to +(.3cm,-.9cm);
\draw [thick] (u1) to +(-.5cm,-.7cm); \draw [thick] (u1) to +(.5cm,-.7cm);
\draw [thick] (u2) to +(-.3cm,-.9cm); \draw [thick] (u2) to +(0cm,-.9cm);
\draw [thick] (u2) to +(.3cm,-.9cm);
\end{tikzpicture}\end{center}
and the cyclic permutation $\sigma_v = (1~3~2) \in S_3$.  Then we have 
\[\begin{split}
k_v^1 &= 3, k_v^2=1=k_v^3, k_v^-=3, k_v^+=2,\\
\sigma_v\langle 3,1,1\rangle &= (1~3~5~2~4) \in S_5.
\end{split}\]
On the other hand, consider the vertex $x$ and the transposition $\sigma_x = (1~3) \in S_3$.  Then we have 
\[\begin{split}k_x^1 &= 1=k_x^2, k_x^3 = 5, k_x^- = 3, k_x^+ = 0,\\ \sigma_x\langle 1,1,5\rangle &= (1~7~5~3)(4~2~6) \in S_7.
\end{split}\]\dqed
\end{example}

\begin{example}
Suppose $\G$ is the braid group operad $\B$ in Definition \ref{def:braid-group-operad}.  Suppose $(T,\sigma)$ is a braided tree and $(\Cor_{(\uc;d)}, \sigma')$ is a braided corolla with unique vertex $v$ and $\profofv = \duc = \Prof(T,\sigma)$.  Then the braided tree substitution at $v$ is \[\bigl(\Cor_{(\uc;d)}, \sigma'\bigr) \comp_v (T,\sigma) = (T,\sigma\sigma'),\] where $\sigma\sigma' \in B_{|\uc|}$ is the product in the braid group.  For example, if $(T,\sigma)$ is the braided tree in Example \ref{ex:braided-tree} and $\sigma' = s_2 \in B_3$, then the braided tree substitution is $(T,s_2s_2s_1^{-1}s_2)$, as depicted below.
\begin{center}\begin{tikzpicture}[xscale=1.7, yscale=1.3]
\node[plain] (w) {$w$}; \node[left=1.5 of w] (T) {$T$}; \node[below=1.75 of T] (b) {$\sigma$}; \node[below=1 of b] () {$\sigma'$};
\node[below=3 of w, inner sep=0] (c1) {};\node[left=.7 of c1, inner sep=0] (c2) {};\node[right=.7 of c1, inner sep=0] (c3) {};
\node[below=1 of c1] (d2) {}; \node[below=1 of c2] (d1) {}; \node[below=1 of c3] (d3) {};
\node[above=1 of c1] (i2) {};\node[above=1 of c3] (i3) {};
\node[above=2 of c3, plain] (v) {$v$}; \node[above=2 of c2, plain] (u) {$u$};
\draw[outputleg] (w) to +(0,.5);
\draw[thick] (u) to (w); \draw[thick] (v) to (w);
\draw[thick] (c1) to [out=90, in=270] (u);
\draw[line width=6pt, white, shorten <=.5pt, shorten >=.5pt] (c2) to[out=100, in=270] (i3);
\draw[thick, shorten >=-8pt] (c2) to[out=80, in=270] (i3);
\draw[line width=6pt, white, shorten <=.5pt, shorten >=.5pt] (c3) to[out=80, in=270] (i2);
\draw[thick, shorten >=-8pt] (c3) to[out=100, in=270] (i2);
\draw[thick] (i2) to[out=90, in=270] (v);
\draw[line width=6pt, white, shorten <=.5pt, shorten >=.5pt] (i3) to[out=90, in=315] (u);
\draw[thick] (i3) to[out=90, in=315] (u);
\draw[thick, shorten >=-2pt] (d1) to (c2);
\draw[thick, shorten >=-2pt] (d2) to[out=80, in=270] (c3);
\draw[line width=6pt, white, shorten <=.5pt, shorten >=.5pt] (d3) to[out=135, in=270]  (c1);
\draw[thick, shorten >=-2pt] (d3) to[out=135, in=270] (c1);
\end{tikzpicture}\end{center}
\dqed
\end{example}

Next is the $\G$-analogue of Lemma \ref{lem:planar-tree-sub-vertex}.

\begin{lemma}\label{lem:g-tree-sub-vertex}
In the context of Definition \ref{def:treesub-g}, the following statements hold.
\begin{enumerate}
\item $\Prof(T') = \Prof(T) \Bigl(\id_{k_v^-} \oplus \sigmavbar\langle k_v^1,\ldots,k_v^n\rangle^{-1} \oplus \id_{k_v^+}\Bigr)$.
\item $\Prof\bigl((T,\sigma) \comp_v (T_v,\sigma_v)\bigr)= \Prof(T,\sigma)$. 
\item $\Vt\bigl((T,\sigma) \comp_v (T_v,\sigma_v)\bigr)= \bigl[\Vt(T) \setminus \{v\}\bigr] \amalg \Vt(T_v)$.
\item Internal edges in $T_v$ remain internal edges in $(T,\sigma) \comp_v (T_v,\sigma_v)$.
\item Left Unity: If \[(T,\sigma) = (\Cor_{\profoft},\id_{|\inp(T)|}),\] then \[(\Cor_{\profoft},\id_{|\inp(T)|}) \comp_v (T_v,\sigma_v) = (T_v,\sigma_v).\]
\item Right Unity: If \[(T_v,\sigma_v) = (\Cor_{\profofv},\id_{|\inp(v)|}),\] then \[(T,\sigma) \comp_v (\Cor_{\profofv},\id_{|\inp(v)|}) = (T,\sigma).\]
\item Vertical Associativity: Suppose $u \in \Vt(T_v)$ and $(T_u,\sigma_u)$ is a $\G$-tree such that $\profofu = \Prof(T_u,\sigma_u)$.  Then \[\bigl((T,\sigma)\comp_v (T_v,\sigma_v)\bigr) \comp_u (T_u,\sigma_u) = (T,\sigma) \comp_v \bigl((T_v,\sigma_v) \comp_u (T_u,\sigma_u)\bigr).\]
\item Horizontal Associativity: Suppose $v\not= w \in \Vt(T)$ and $(T_w,\sigma_w)$ is a $\G$-tree such that $\profofw = \Prof(T_w,\sigma_w)$.  Then \[\bigl((T,\sigma)\comp_v (T_v,\sigma_v)\bigr) \comp_w (T_w,\sigma_w) = \bigl((T,\sigma)\comp_w (T_w,\sigma_w)\bigr) \comp_v (T_v,\sigma_v).\]
\item Corolla Decomposition: With $v$ denoting the unique vertex in the $\profoft$-corolla $\Cor_{\profoft}$, there is an equality \[(T,\sigma) = (\Cor_{\profoft},\sigma) \compv (T,\id_{|\inp(T)|}).\]
\end{enumerate}
\end{lemma}

\begin{proof}
All of these assertions are straightforward to check.  As an example, consider the horizontal associativity assertion.  The equality of the planar $\colorc$-colored trees follows from the horizontal associativity of tree substitution at a vertex in Lemma \ref{lem:planar-tree-sub-vertex}.  

For the equality of input equivariances, there are three cases.  First suppose neither $\out(v) < \out(w)$ nor $\out(w) < \out(v)$.  In this case, both input equivariances are equal to \[\Bigl(\id_1^{\oplus {k_v^-}} \oplus \sigma_v\langle k_v^1,\ldots,k_v^{|\inp(v)|}\rangle \oplus \id_1^{\oplus k_w^--k_v^--k_v} \oplus \sigma_w\langle k_w^1,\ldots,k_w^{|\inp(w)|}\rangle \oplus \id_1^{\oplus {k_w^+}}\Bigr) \cdot \sigma\] or \[\Bigl(\id_1^{\oplus {k_w^-}} \oplus \sigma_w\langle k_w^1,\ldots,k_w^{|\inp(w)|}\rangle \oplus \id_1^{\oplus k_v^- -k_w^- -k_w} \oplus \sigma_v\langle k_v^1,\ldots,k_v^{|\inp(v)|}\rangle \oplus \id_1^{\oplus {k_v^+}}\Bigr) \cdot \sigma.\] The other two cases, when either $\out(v)<\out(w)$ or $\out(w)<\out(v)$, follow from the action operad axiom \eqref{group-operad-identites}.
\end{proof}

We can also substitute $\G$-trees using all the vertices at the same time.

\begin{definition}[Group Tree Substitution]\label{def:g-tree-sub}
Suppose $(T,\sigma)$ is a $\G$-tree, and for each vertex $v \in \Vt(T)$, $(T_v,\sigma_v)$ is a $\G$-tree such that $\profofv = \Prof(T_v,\sigma_v)$.  Define the \emph{$\G$-tree substitution}\index{G-tree@$\G$-tree!substitution} \[(T,\sigma)(T_v,\sigma_v)_{v\in \Vt(T)} = \Bigl(\bigl((T,\sigma) \comp_{v_1} (T_{v_1},\sigma_{v_1})\bigr) \comp_{v_2} \cdots \Bigr) \comp_{v_n} (T_{v_n},\sigma_{v_n}),\] where $(v_1,\ldots,v_n)$ is an arbitrary choice of an ordering of $\Vt(T)$.  To simplify the notation, we will sometimes abbreviate the $\G$-tree substitution to $(T,\sigma)(T_v,\sigma_v)_v$ or $(T,\sigma)(T_v,\sigma_v)$.
\end{definition}

The following statements follow immediately from Lemma \ref{lem:g-tree-sub-vertex}.

\begin{lemma}\label{lem:g-tree-sub}
In the context of Definition \ref{def:g-tree-sub}, the following statements hold.
\begin{enumerate}
\item The $\G$-tree substitution is independent of the choice of an ordering of $\Vt(T)$.
\item $\Prof\bigl((T,\sigma)(T_v,\sigma_v)\bigr) = \Prof(T,\sigma)$.\index{profile!G-tree substitution@$\G$-tree substitution}
\item $\Vt\bigl((T,\sigma)(T_v,\sigma_v)\bigr) = \coprod_{v\in\Vt(T)} \Vt(T_v)$.
\item Internal edges in each $T_v$ remain internal edges in the $\G$-tree substitution.
\item Unity: $(T,\sigma)\bigl(\Cor_{\profofv},\id_{|\inp(v)|}\bigr) = (T,\sigma)$.\index{unity!G-tree substitution@$\G$-tree substitution}
\item Associativity: Suppose for each $v \in \Vt(T)$ and each $u \in \Vt(T_v)$, $(T_v^u,\sigma_v^u)$ is a $\G$-tree such that $\profofu = \Prof(T_v^u,\sigma_v^u)$.  Then there is an equality\index{associativity!G-tree substitution@$\G$-tree substitution} 
\[\begin{split}
&\Bigl[(T,\sigma)(T_v,\sigma_v)_{v\in\Vt(T)}\Bigr](T_v^u,\sigma_v^u)_{v\in\Vt(T),u\in\Vt(T_v)}\\
&= (T,\sigma)\Bigl[(T_v,\sigma_v)(T_v^u,\sigma_v^u)_{u\in\Vt(T_v)}\Bigr]_{v\in\Vt(T)}.
\end{split}\]
\end{enumerate}
\end{lemma}

\begin{proof}
All of these assertions are straightforward to check.  For example, the associativity assertion follows from the vertical associativity and the horizontal associativity properties in Lemma \ref{lem:g-tree-sub-vertex} by an induction.
\end{proof}

\section{Symmetric Operad for Group Operads}
\label{sec:symmetric-g-operad}

In this section, we describe the symmetric operad whose algebras, in the sense of Definition \ref{def:symmetric-operad-algebra}, are $\G$-operads as in Definition \ref{def:g-operad}.

\begin{definition}\label{def:gopc}
Fix a set $\colorc$ and an action operad $(\G,\omega)$.  Define the $(\Profcc)$-colored symmetric operad\index{symmetric operad!for G-operads@for $\G$-operads}\index{G-operad@$\G$-operad!symmetric operad encoding} $\gopc$ in $\Set$ as follows.
\begin{description}
\item[Entries] For $t, s_1,\ldots,s_n \in \Profcc$, define the set \[\gopc\tsonesn = \Bigl\{(T,\sigma,\rho) : \Prof(T,\sigma)=t,\, \Prof(\rho(j)) = s_j \text{ for $1\leq j \leq n$}\Bigr\}\] in which:
\begin{itemize}
\item $(T,\sigma)$ is a $\colorc$-colored $\G$-tree with profile $t$.
\item $\rho$ is a vertex ordering of $T$ with respect to which the $j$th vertex has profile $s_j$.
\end{itemize}
\item[Units] For $t \in \Profcc$, the $t$-colored unit is \[\Bigl(\Cor_t, \id_{|\inp(T)|}, \Id\Bigr) \in \gopc\ttsingle\] with:
\begin{itemize}
\item $\Cor_t$  the $t$-corolla in Example \ref{ex:corolla};
\item $\id_{|\inp(T)|} \in \G(|\inp(T)|)$ the group multiplication unit;
\item $\Id$ the trivial vertex ordering for the $t$-corolla.
\end{itemize}
\item[Composition] Suppose $t$ and $\us = (s_1,\ldots,s_n)$ are as above with $1 \leq i \leq n$, $\ur = (r_1,\ldots,r_m)$ with each $r_j \in \Profcc$.  The $\compi$-composition \[\nicexy{\gopc\tus \times \gopc\siur \ar[r]^-{\compi} & \gopc\sbinom{t}{\us\compi\ur}}\] is given by the $\G$-tree substitution \[(T,\sigma,\rho) \compi (T',\sigma',\rho') = \Bigl((T,\sigma) \comp_{\rho(i)} (T',\sigma'), \rho\compi\rho'\Bigr)\] at $\rho(i) \in \Vt(T)$ for $(T,\sigma,\rho) \in \gopc\tus$ and $(T',\sigma',\rho') \in \gopc\siur$.  The vertex ordering \[\nicexy{\{1,\ldots,n+m-1\} \ar[r]^-{\rho \compi \rho'}_-{\cong} & \bigl(\Vt(T) \setminus\{\rho(i)\}\bigr) \amalg \Vt(T')}\] is given by 
\begin{equation}\label{rho-compi-rho}
(\rho\compi\rho')(j) = \begin{cases}\rho(j) \in \Vt(T)\setminus\{\rho(i)\} & \text{ if $j < i$},\\  \rho'(j-i+1) \in \Vt(T')& \text{ if $i \leq j \leq i+m-1$},\\
 \rho(j-m+1) \in \Vt(T) \setminus \{\rho(i)\}& \text{ if $i+m \leq j \leq n+m-1$}.\end{cases}
\end{equation}
\item[Symmetry] With $t$ and $\us$ as above and $\tau \in S_n$, the $\colorc$-colored symmetric sequence structure \[\nicexy{\gopc\tus \ar[r]^-{\tau} & \gopc\tustau}\] is given by \[(T,\sigma,\rho)\tau = (T,\sigma,\rho\tau) \in \gopc\tustau,\] in which the vertex ordering $\rho\tau$ is the composite \[\nicexy{\{1,\ldots,n\} \ar[r]^-{\tau}_-{\cong} & \{1,\ldots,n\} \ar[r]^-{\rho}_-{\cong} & \Vt(T) \ar@{<-}`u[ll]`[ll]_-{\rho\tau}[ll].}\]
\end{description}
\end{definition}

\begin{lemma}\label{lem:gopc-symmetric-operad}
$\gopc$ is a $(\Profcc)$-colored symmetric operad in $\Set$.
\end{lemma}

\begin{proof}
Using the description of a $\colorc$-colored symmetric operad in Proposition \ref{prop:symmetric-operad-compi}, the associativity and unity axioms for $\gopc$ follow from Lemma \ref{lem:g-tree-sub-vertex}.  The equivariance axiom follows directly from the definition of $\rho \comp_i \rho'$ in \eqref{rho-compi-rho}. 
\end{proof}

Recall the strong symmetric monoidal functor $\iota : \Set \to \M$ that sends each set $X$ to the $X$-indexed coproduct $\coprod_X \tensorunit \in \M$.  As discussed in Example \ref{ex:set-to-m}, it induces a functor \[\nicexy{\gopcset \ar[r]^-{\iota_*} & \goperadcm}\] from $\colorc$-colored $\G$-operads in $\Set$ to $\colorc$-colored $\G$-operads in $\M$.  For $\O\in\gopcset$, we call $\iota_*\O \in \goperadcm$ its \index{image}\emph{image} in $\M$.

\begin{definition}\label{def:gopcm}
For the $(\Profcc)$-colored symmetric operad $\gopc$ in Lemma \ref{lem:gopc-symmetric-operad}, its image in $\M$ is denoted by $\gopcm$.
\end{definition}

\begin{remark} A generic entry of $\gopcm$ has the form\[\gopcm\tus = \coprodover{\gopc\tus} \tensorunit \in \M\]  Its $(\Profcc)$-colored symmetric operad structure is uniquely determined by that in $\gopc$.\dqed
\end{remark}

In the next result, we use the description in Proposition \ref{prop:symmetric-operad-algebra-defs} of an algebra over a colored symmetric operad.

\begin{theorem}\label{gopcm-algebra}
For each action operad $\G$, there is a canonical isomorphism \[\algm\bigl(\gopcm\bigr) \cong \goperadcm\] between the category of $\gopcm$-algebras and the category of $\colorc$-colored $\G$-operads in $\M$.
\end{theorem}

\begin{proof}
Since the symmetric operad $\gopcm$ is $(\Profcc)$-colored, a $\gopcm$-algebra is, first of all, a $(\Profcc)$-colored object \[\O = \bigl\{\O\duc \in \M: \duc \in \Profcc\bigr\}.\]  The rest of its $\colorc$-colored $\G$-operad structure is given as follows.  For $t,s_1,\ldots,s_n \in \Profcc$ with $\us=(s_1,\ldots,s_n)$, the $\gopcm$-action structure morphism \eqref{planar-operad-algebra-action} \[\nicexy{\gopcm\tus \otimes \O_{\us} \cong \coprodover{(T,\sigma,\rho) \in \gopc\tus} \,\,\overset{n}{\bigtensorover{j=1}} \O\bigl(\Prof(\rho(j))\bigr) \ar[r]^-{\lambda} & \O(t)}\] is determined by the restriction 
\begin{equation}\label{goperad-restricted-structure}
\nicexy@C+.4cm{\O(T,\sigma,\rho)\defn \overset{n}{\bigtensorover{j=1}} \O\bigl(\Prof(\rho(j))\bigr) \ar[r]^-{\lambda_{(T,\sigma,\rho)}} & \O(t) = \O\bigl(\Prof(T,\sigma)\bigr)}
\end{equation}
for $(T,\sigma,\rho) \in \gopc\tus$.  

For the operadic composition, consider 
\begin{itemize}
\item $d \in \colorc$, $\uc=(c_1,\ldots,c_n),\ub_1,\ldots,\ub_n \in \Profc$,
\item $\ub = (\ub_1,\ldots,\ub_n)$,
\item $k=|\ub_1|+\cdots+|\ub_n|$ as in Definition \ref{def:planar-operad-generating},
\item the $2$-level tree $T = T(\{\ub_j\};\uc;d)$ in Example \ref{ex:2-level-tree}, and
\item the canonical vertex ordering \[\nicexy{\{1,\ldots,n+1\} \ar[r]^-{\rho}_-{\cong} & \Vt(T)}\] given by \[\rho(i) =\begin{cases}v & \text{ if $i=1$},\\ u_{i-1} & \text{ if $2 \leq i \leq n+1$}. \end{cases}\]  
\end{itemize}
The restricted $\gopcm$-action structure morphism 
\begin{equation}\label{goperad-restricted-gamma}
\nicexy@C+1cm{\O(T,\id_k,\rho) = \O\duc \otimes \overset{n}{\bigtensorover{j=1}} \O\cjubj \ar[r]^-{\gamma\,=\,\lambda_{(T,\id_k,\rho)}} & \O\dub}
\end{equation}
is the operadic composition in \eqref{operadic-composition}.

For each color $c\in \colorc$, consider the $c$-colored exceptional edge $\uparrow_c$ in Example \ref{ex:exceptional-edge} with the trivial vertex ordering $\rho$.  The restricted $\gopcm$-action structure morphism 
\begin{equation}\label{goperad-restricted-unit}
\nicexy@C+1cm{\O(\uparrow_c,\id_1,\rho) = \tensorunit \ar[r]^-{\operadunit_c\,=\,\lambda_{(\uparrow_c,\id_1,\rho)}} & \O\cc}
\end{equation}
is the $c$-colored unit in \eqref{c-colored-unit}.

For the $\colorc$-colored $\G$-sequence structure, consider
\begin{itemize}
\item $\duc \in \Profcc$ as above with $|\uc|=n$, 
\item $\sigma \in \G(n)$ with underlying permutation $\sigmabar \in S_n$, and
\item the $\duc$-corolla $\Cor_{(\uc;d)}$ in Example \ref{ex:corolla} with the trivial vertex ordering $\rho$.
\end{itemize}
The restricted $\gopcm$-action structure morphism 
\begin{equation}\label{goperad-restricted-equivariance}
\nicexy@C+1.5cm{\O\bigl(\Cor_{(\uc;d)}, \sigma,\rho\bigr)=\O\duc \ar[r]^-{\sigma\,=\,\lambda_{(\Cor_{(\uc;d)}, \sigma,\rho)}} & \O\ducsigmabar}
\end{equation}
is the $\colorc$-colored $\G$-sequence structure in Definition \ref{def:g-sequences}.  A direct inspection shows that the above structure specifies a $\colorc$-colored $\G$-operad structure on $(\O,\gamma,\operadunit)$ as in Proposition \ref{def:g-operad-generating}.
 
For the other direction, suppose $(\O,\gamma,\operadunit)$ is a $\colorc$-colored $\G$-operad in $\M$.  The underlying entries of $\O$ is a $(\Profcc)$-colored object.  To define the restricted $\gopcm$-action structure morphism $\lambda_{(T,\sigma,\rho)}$ in \eqref{goperad-restricted-structure}, first note that it is enough to define it when $\rho$ is the canonical vertex ordering by the equivariance axiom of a $\gopcm$-algebra \eqref{operad-algebra-eq}.  For two vertices $u$ and $v$ in $T$, we write $u \prec v$ if there exists an internal edge $e \in T$ with initial vertex $v$ and terminal vertex $u$.  We define the restricted structure morphisms $\lambda_{(T,\sigma,\rho)}$ with $\rho$ the canonical vertex ordering by induction on the number \[N = \max\bigl\{n :  v_1 \prec v_2 \prec \cdots \prec v_n \text{ exists in $\Vt(T)$}\bigr\}.\] Note that if the chain $v_1 \prec \cdots \prec v_N$ realizes the maximum, then $v_1$ must be the initial vertex of the output of $T$.
\begin{itemize}
\item If $N=0$, then $\Vt(T)$ is empty, and $T$ is the $c$-colored exceptional edge $\uparrow_c$ for some $c \in \colorc$.  In this case, we define the restricted structure morphism $\lambda_{(\uparrow_c,\id_1,\rho)}$ as in \eqref{goperad-restricted-unit} using the $c$-colored unit of $\O$.
\item If $N=1$, then $T$ has exactly one vertex, so it is the $\duc$-corolla $\Cor_{(\uc;d)}$ for some $\duc \in \Profcc$.  In this case, we define the restricted structure morphism $\lambda_{(\Cor_{(\uc;d)}, \sigma,\rho)}$ as in \eqref{goperad-restricted-equivariance} using the $\colorc$-colored $\G$-sequence structure of $\O$.
\item If $N=2$, then we define the restricted structure morphism using \eqref{goperad-restricted-gamma}, \eqref{goperad-restricted-unit}, and \eqref{goperad-restricted-equivariance}.
\item Inductively, if $N>2$, then we define the restricted structure morphism using the induction hypothesis and the case $N=2$.
\end{itemize}
One can check that this defines a $\gopcm$-algebra structure on $\O$ and that the two constructions above are functorial and are inverses of each other.
\end{proof}

\begin{example}
In the context of Theorem \ref{gopcm-algebra}:
\begin{enumerate}
\item If $\G$ is the planar group operad\index{planar operad!symmetric operad encoding} $\P$ in Example \ref{ex:trivial-group-operad}, then there is a canonical isomorphism \[\algm\bigl(\popcm\bigr) \cong \Poperadcm\] between the category of $\popcm$-algebras and the category of $\colorc$-colored planar operads in $\M$. 
\item If $\G$ is the symmetric group operad\index{symmetric operad!symmetric operad encoding} $\S$ in Example \ref{ex:symmetric-group-operad}, then there is a canonical isomorphism \[\algm\bigl(\sopcm\bigr) \cong \Soperadcm\] between the category of $\sopcm$-algebras and the category of $\colorc$-colored symmetric operads in $\M$.  This is also a special case of Theorem 14.1 in \cite{bluemonster}.
\item If $\G$ is the braid group operad\index{braided operad!symmetric operad encoding} $\B$ in Definition \ref{def:braid-group-operad}, then there is a canonical isomorphism \[\algm\bigl(\bopcm\bigr) \cong \Boperadcm\] between the category of $\bopcm$-algebras and the category of $\colorc$-colored braided operads in $\M$. 
\item If $\G$ is the pure braid group\index{pure braided operad!symmetric operad encoding} operad $\PB$ in Definition \ref{def:pure-braid-group-operad}, then there is a canonical isomorphism \[\algm\bigl(\pbopcm\bigr) \cong \PBoperadcm\] between the category of $\pbopcm$-algebras and the category of $\colorc$-colored pure braided operads in $\M$. 
\item If $\G$ is the ribbon group operad\index{ribbon operad!symmetric operad encoding} $\R$ in Definition \ref{def:ribbon-group-operad}, then there is a canonical isomorphism \[\algm\bigl(\ropcm\bigr) \cong \Roperadcm\] between the category of $\ropcm$-algebras and the category of $\colorc$-colored ribbon operads in $\M$.
\item If $\G$ is the pure ribbon group operad\index{pure ribbon operad!symmetric operad encoding} $\PR$ in Definition \ref{def:pure-ribbon-group-operad}, then there is a canonical isomorphism \[\algm\bigl(\propcm\bigr) \cong \PRoperadcm\] between the category of $\propcm$-algebras and the category of $\colorc$-colored pure ribbon operads in $\M$. 
\item If $\G$ is the cactus group operad\index{cactus operad!symmetric operad encoding} $\Cac$ in Definition \ref{def:cactus-group-operad}, then there is a canonical isomorphism \[\algm\bigl(\cacopcm\bigr) \cong \Cacoperadcm\] between the category of $\cacopcm$-algebras and the category of $\colorc$-colored cactus operads in $\M$.
\item If $\G$ is the pure cactus group operad\index{pure cactus operad!symmetric operad encoding} $\PCac$ in Definition \ref{def:pure-cactus-group-operad}, then there is a canonical isomorphism \[\algm\bigl(\pcacopcm\bigr) \cong \Pcacoperadcm\] between the category of $\pcacopcm$-algebras and the category of $\colorc$-colored pure cactus operads in $\M$. \dqed
\end{enumerate}
\end{example}

Recall that our ambient symmetric monoidal category $\M$ is assumed to have all small limits and colimits.

\begin{corollary}\label{goperad-bicomplete}
The category $\goperadcm$ of $\colorc$-colored $\G$-operads in $\M$ has all small\index{G-operad@$\G$-operad!bicomplete} limits and colimits, with filtered colimits, reflexive coequalizers, and limits created and preserved by the forgetful functor \[\goperadcm \to \M^{\Profcc}.\]
\end{corollary}

\begin{proof}
This is a consequence of Theorem \ref{gopcm-algebra} above and Proposition 4.2.1 in \cite{white-yau}, which is the same assertion for the category of algebras over a general colored symmetric operad, including $\gopcm$.
\end{proof}

\chapter{Group Operads with Varying Colors}\label{ch:all-group-operads}

Fix an action operad $\G$ and a complete and cocomplete symmetric monoidal category $(\M,\otimes,\tensorunit)$ whose monoidal product commutes with small colimits on each side.  The purpose of this chapter is to study the category of all $\G$-operads, in which the color sets may vary.  

In Section \ref{sec:category-all-goperad} we introduce the category $\gopm$ of all $\G$-operads in $\M$ and observe that it has all small limits and colimits.  We also observe that each morphism of action operads induces an adjunction between the respective categories of all group operads.

The category of all symmetric operads in $\Set$ has a symmetric monoidal structure called the Boardman-Vogt tensor product.  In Section \ref{sec:symmetric-monoidal} we construct a symmetric monoidal structure on the category of all $\G$-operads in $\Set$ that recovers the Boardman-Vogt tensor product when $\G$ is the symmetric group operad $\S$.   In Section \ref{sec:goperad-monoidal-closed} we extend this symmetric monoidal structure to a symmetric monoidal \emph{closed} structure.

In Section \ref{sec:comparing-sms} we show that these symmetric monoidal structures are compatible with changing the action operad.  In Section \ref{sec:non-strong} we observe that for the morphisms among the planar, symmetric, (pure) ribbon, (pure) braided, and (pure) cactus group operads, the induced symmetric monoidal functors are not strong monoidal.  In Section \ref{sec:lfp} we observe that the category of all $\G$-operads is locally finitely presentable.

\section{Category of All Group Operads}\label{sec:category-all-goperad}

First we specify the category of all $\G$-operads.

\begin{definition}\label{def:g-operads}
Given an action operad $(\G,\omega)$, define the category\index{category!all G-operads@all $\G$-operads}\index{G-operad@$\G$-operad!category of} \[\gopm\] of $\G$-operads in $\M$ as follows.
\begin{description}
\item[Objects] An object in $\gopm$ is a pair $(\colorc,\O)$ consisting of a set $\colorc$ and a $\colorc$-colored $\G$-operad $\O$ in $\M$.
\item[Morphisms] A morphism \[f : (\colorc,\O) \to (\colord,\P)\] in $\gopm$ consists of
\begin{enumerate}[label=(\roman*)]
\item a function $f_0 : \colorc \to \colord$ and
\item a morphism $f_1 : \O \to f_0^*\P$ of $\colorc$-colored $\G$-operads in $\M$, where \[f_0^* : \goperaddm \to \goperadcm\] is the change-of-color functor in Example \ref{ex:g-operad-change-colors}.
\end{enumerate}
\end{description}
Identity morphisms and composition are defined in the obvious way.  We will often write both $f_0$ and $f_1$ simply as $f$ below.
\end{definition}

\begin{remark}When $\M$ is the category $\Set$, the category $\gopm$ is referred to as the category of \index{G-multicategory@$\G$-multicategory}$\G$-multicategories in Gurski \cite{gurski}.\dqed\end{remark}

Recall that our ambient category $\M$ is assumed to have all small limits and colimits.

\begin{proposition}\label{gopm-bicomplete}
The category $\gopm$ has all\index{G-operad@$\G$-operad!bicomplete} small limits and colimits.
\end{proposition}

\begin{proof}
As is usually the case, small limits in $\gopm$ are computed entrywise in $\M$.  Indeed, it is enough to observe that $\gopm$ has small products and equalizers.  For products of colored $\G$-operads, the color set is the product of the given color sets, and each entry is the product in $\M$ of the corresponding entries in the given set of colored $\G$-operads.  For equalizers of colored $\G$-operads, the color set is the equalizer of the given color sets in $\Set$, and each entry is the entrywise equalizer in $\M$.

Next we consider small colimits in $\gopm$.  Writing $\Cat$ for the category of small categories, for a functor $F : \C^{\op} \to \Cat$ for some category $\C$, recall that the \index{Grothendieck construction}\emph{Grothendieck construction} for $F$ is the category \[\int_{\C^{\op}} F\] in which:
\begin{itemize}
\item An object is a pair of objects $\bigl(c \in \C; x \in Fc\bigr)$.
\item A morphism $f : (c;x) \to (d;y)$ is a pair of morphisms \[\Bigl(f_0 : c \to d \in \C; f_1 : x \to (Ff_0^{\op})(y) \in Fc\Bigr).\]
\end{itemize}
Identity morphisms and composition are defined in the obvious way.  The category $\gopm$ is the Grothendieck construction $\int_{\Set^{\op}} \goperaddashm$ for the functor \[\nicexy@C+.5cm{\Set^{\op} \ar[r]^-{\goperaddashm} & \Cat}, \quad \Set \ni \colorc \mapsto \goperadcm\] that sends each set $\colorc$ to the category $\goperadcm$ of $\colorc$-colored $\G$-operads in $\M$.  On functions between sets, this functor is given by the change-of-color functor in Example \ref{ex:g-operad-change-colors}.  Since each category $\goperadcm$ is cocomplete by Corollary \ref{goperad-bicomplete}, the Grothendieck construction \[\int_{\Set^{\op}} \goperaddashm = \gopm\] is also cocomplete by Proposition 2.4.4 in Harpaz-Prasma \cite{hp}.
\end{proof}

\begin{remark}In the previous proof, if the category $\Cat$ is replaced by the category $\Set$ with $f_1$ an equality of elements, then the Grothendieck construction is usually called the \emph{category of elements}.  See, for example, \cite{maclane-moerdijk} (pages 41-43).\dqed\end{remark}

\begin{theorem}\label{goneopm-gtwoopm}\index{change of action operads!all G-operads@all $\G$-operads}
Suppose $\varphi : \Gone\to\Gtwo$ is a morphism of action operads.  Then there is an induced adjunction \[\nicexy{\goneopm \ar@<2pt>[r]^-{\varphi_!} & \gtwoopm \ar@<2pt>[l]^-{\varphi^*}}\] in which:
\begin{enumerate}[label=(\roman*)]
\item The left adjoint $\varphi_!$ restricts to the fixed color set version in Theorem \ref{phi-goneop-gtwoop} for each colored $\Gone$-operad.
\item The right adjoint $\varphi^*$ restricts to the fixed color set version in \eqref{gtwoop-goneop} for each colored $\Gtwo$-operad.
\end{enumerate}
\end{theorem}

\begin{proof}
Suppose $\O$ is a $\colorc$-colored $\Gone$-operad in $\M$, and $\Q$ is a $\colord$-colored $\Gtwo$-operad in $\M$.  A morphism \[\nicexy{(\colorc,\O) \ar[r]^-{f} & \varphi^*(\colord,\Q) = (\colord,\varphi^*\Q)} \in \goneopm\] consists of:
\begin{enumerate}[label=(\roman*)]
\item a function $f : \colorc \to \colord$;
\item a morphism $f : \O \to f^*\varphi^*\Q$ of $\colorc$-colored $\Gone$-operads.
\end{enumerate}
By Proposition \ref{change-color-goperad} $f^*\varphi^*\Q$ is equal to $\varphi^*f^*\Q$, in which $f^*\Q$ is the $\colorc$-colored $\Gtwo$-operad in $\M$ obtained from $\Q$ by applying the change-of-color functor $f^*$ in Example \ref{ex:g-operad-change-colors}.  By the adjunction \[\nicexy{\goneoperadcm \ar@<2pt>[r]^-{\varphi_!} & \gtwooperadcm \ar@<2pt>[l]^-{\varphi^*}}\] in Theorem \ref{phi-goneop-gtwoop}, the morphism \[f : \O\to\varphi^*f^*\Q\] corresponds to a morphism \[\varphi_!\O\to f^*\Q\] of $\colorc$-colored $\Gtwo$-operads.  Therefore, the original morphism $f$ corresponds to a morphism \[\nicexy{\varphi_!(\colorc,\O) = (\colorc,\varphi_!\O) \ar[r] & (\colord,\Q) \in \gtwoopm},\] which establishes the desired adjunction.
\end{proof}

\section{Symmetric Monoidal Structure}\label{sec:symmetric-monoidal}

In this section, we construct a symmetric monoidal closed structure on the category of all $\G$-operads in $\Set$ that generalizes the Boardman-Vogt tensor product of symmetric operads in \cite{boardman-vogt}.   The desired symmetric monoidal product on $\gopset$ involves certain special permutations, which we define next.

\begin{definition}\label{def:special-permutation}
For $m,n \geq 0$, the \index{special permutation}\emph{special $(m,n)$-permutation} $\sigma_{m,n} \in S_{mn}$ is given by \[\sigma_{m,n}\bigl((k-1)n+j\bigr) = (j-1)m+k\] for $1 \leq k \leq m$ and $1 \leq j \leq n$, and by \[\sigma_{m,0} = \sigma_{0,n}= \id_0 \in S_0\] if either $m=0$ or $n=0$.
\end{definition}

\begin{interpretation}
The special $(m,n)$-permutation $\sigma_{m,n}$ redistributes $m$ group of $n$ objects into $n$ groups of $m$ objects in such a way that the $j$th object in the $k$th group is sent to the $k$th object in the $j$th group.  In particular, we have \[\sigma_{m,n}^{-1}=\sigma_{n,m}.\]  Special permutations were used by Markl \cite{markl-bialgebra} to define certain special elements in props. \dqed
\end{interpretation}

\begin{example}\label{ex:id-special-permutation}
For each $n \geq 0$, $\sigma_{1,n} = \sigma_{n,1}= \id_n \in S_n$.\dqed
\end{example}

\begin{example}
The special $(2,3)$-permutation $\sigma_{2,3} \in S_6$ is the cyclic permutation $(2~3~5~4)$, i.e., $\sigma_{2,3}(1)=1$, $\sigma_{2,3}(2)=3$, $\sigma_{2,3}(3)=5$, $\sigma_{2,3}(4)=2$, $\sigma_{2,3}(5)=4$, and $\sigma_{2,3}(6)=6$.  
\begin{center}\begin{tikzpicture}
\draw (-.5,.5) node{$\sigma_{2,3}\,=$};
\foreach \x in {1,...,6} {\draw (\x,-.2) node {\x}; \draw (\x,1.2) node {\x};}
\draw[->] (1,0)--(1,1);  \draw[->] (2,0)--(3,1); \draw[->] (3,0)--(5,1); 
\draw[->] (4,0)--(2,1); \draw[->] (5,0)--(4,1); \draw[->] (6,0)--(6,1);
\end{tikzpicture}\end{center}
The special $(2,2)$-permutation $\sigma_{2,2} \in S_4$ is the adjacent transposition $(2~3)$.
\begin{center}\begin{tikzpicture}
\draw (-.5,.5) node{$\sigma_{2,2}\,=$};
\foreach \x in {1,...,4} {\draw (\x,-.2) node {\x}; \draw (\x,1.2) node {\x};}
\draw[->] (1,0)--(1,1);  \draw[->] (2,0)--(3,1); \draw[->] (3,0)--(2,1); \draw[->] (4,0)--(4,1);
\end{tikzpicture}\end{center}
So we have $\sigma_{2,2}(1)=1$, $\sigma_{2,2}(2)=3$, $\sigma_{2,2}(3)=2$, and $\sigma_{2,2}(4)=4$.\dqed
\end{example}

\begin{definition}
For $k \geq 0$, an element $g \in \G(k)$ is\index{G-operad@$\G$-operad!special element} \emph{special} if its underlying permutation $\omega(g)=\gbar \in S_k$ is a special $(m,n)$-permutation $\sigma_{m,n}$ for some $m,n \geq 0$ with $mn=k$. 
\end{definition}

\begin{example}
For the symmetric group operad $\S$ with the identity augmentation, the special elements are exactly the special $(m,n)$-permutations.\dqed
\end{example}

\begin{example}\label{ex:special-planar}
For each $n\geq 0$, the group multiplication unit $\id_n \in \G(n)$ is special because its underlying permutation is the identity  permutation $\id_n = \sigma_{1,n} \in S_n$.  Similarly, in
\begin{itemize}
\item the planar group operad $\P$ in Example \ref{ex:trivial-group-operad}, 
\item the pure braid group operad $\PB$ in Definition \ref{def:pure-braid-group-operad}, 
\item the pure ribbon group operad $\PR$ in Definition \ref{def:pure-ribbon-group-operad}, and
\item the pure cactus group operad $\PCac$ in Definition \ref{def:pure-cactus-group-operad},
\end{itemize}
every element is special because the underlying permutation is the identity permutation.\dqed
\end{example}

\begin{example}[Special Braids]
For the braid group operad\index{braid!special} $\B$ in Definition \ref{def:braid-group-operad}:
\begin{itemize}
\item The generator $s_2 \in B_4$ is special because its underlying permutation is the special $(2,2)$-permutation $\sigma_{2,2}\in S_4$.
\item Every element in the pure braid group $PB_n$ is special.
\item More generally, if $\sigma \in B_n$ is a special element and if $p\in PB_n$, then $p\sigma$ and $\sigma p \in B_n$ are both special.\dqed
\end{itemize}
\end{example}

\begin{example}[Special Ribbons]
For the ribbon group operad\index{ribbon!special} $\R$ in Definition \ref{def:ribbon-group-operad}:
\begin{itemize}
\item The generator $r_2 \in R_4$ in Example \ref{ex:generator-ri} is special because its underlying permutation is the special $(2,2)$-permutation $\sigma_{2,2}\in S_4$.
\item Every element in the pure ribbon group $PR_n$ is special.  In particular, the generator $r_n \in R_n$ in Example \ref{ex:generator-rn} is special.
\item More generally, if $r \in R_n$ is a special element and if $p\in PR_n$, then $pr$ and $rp \in R_n$ are both special.\dqed
\end{itemize}
\end{example}

\begin{example}[Special Cacti]
For the cactus group operad\index{cactus!special} $\Cac$ in Definition \ref{def:cactus-group-operad}:
\begin{itemize}
\item The generator $s^{(4)}_{2,3} \in Cac_4$ is special because its underlying permutation is the interval-reversing permutation $\rho^{(4)}_{2,3}=\sigma_{2,2}\in S_4$.
\item Every element in the pure cactus group $PCac_n$ is special.  For example, in the notation of Proposition \ref{third-cactus-group} and Example \ref{ex:cac-three}, the elements $yxyxyx$ and $xyxyxy$ are in $PCac_3$, so they are special.
\item More generally, if $c \in Cac_n$ is a special element and if $p\in PCac_n$, then $pc$ and $cp \in Cac_n$ are both special.\dqed
\end{itemize}
\end{example}

\begin{example}\label{ex:special-element-morphism}
For each morphism $\varphi : \Gone \to \Gtwo$ of action operads as in Definition \ref{def:augmented-group-operad-morphism}, an element $g \in \Gone$ is special if and only if $\varphi(g) \in \Gtwo$ is special.\dqed 
\end{example}

We now define the monoidal product of two colored $\G$-operads in $\Set$.

\begin{definition}\label{def:goperad-tensor}
Suppose
\begin{itemize}
\item $\T \in \goperadcset$ is a $\colorc$-colored $\G$-operad and 
\item $\Q \in \goperaddset$ is a $\colord$-colored $\G$-operad in $\Set$.  
\end{itemize}
Define the $(\colorcd)$-colored $\G$-operad\label{not:tensorg} $\T\tensorg\Q$, called the \emph{$\G$-tensor product}\index{G-tensor product@$\G$-tensor product}\index{monoidal product!of G-operads@of $\G$-operads} of $\T$ and $\Q$, as follows.
\begin{enumerate}
\item First we define the \index{smash product}\emph{smash product}\[\T\wedge\Q \in \Set^{\Profcdcd}\] as the graded set consisting of the following two types of elements:
\begin{enumerate}[label=(\roman*)]
\item For $\czerouc \in \Profcc$, $t \in \T\czerouc$, and $d \in \colord$, it contains the element \[t \otimes d \in (\T\wedge\Q)\sbinom{(c_0;d)}{(\uc;d)},\] where \[(\uc;d)=\bigl((c_1;d),\ldots,(c_n;d)\bigr) \in\Profcd\] if $\uc=(c_1,\ldots,c_n) \in \Profc$.
\item For $\dzeroud \in \Profdd$, $q \in \Q\dzeroud$, and $c \in \colorc$, it contains the element \[c \otimes q \in (\T\wedge\Q)\sbinom{(c;d_0)}{(c;\ud)},\] where \[(c;\ud) = \bigl((c;d_1),\ldots,(c;d_m)\bigr) \in \Profcd\] if $\ud=(d_1,\ldots,d_m) \in \Profd$.
\end{enumerate}
\item Then we define the quotient \[\T\tensorg\Q = \frac{\gopcd\circs(\T\wedge\Q)}{\sim} \in \goperadcdset\] in which:
\begin{itemize}
\item $\gopcd$ is the $(\Profcdcd)$-colored symmetric operad in $\Set$ in Lemma \ref{lem:gopc-symmetric-operad}.
\item $\circs$ is the $(\colorcd)$-colored symmetric circle product in \eqref{symmetric-circle-product}.  
\end{itemize}
So $\gopcd\circs(\T\wedge\Q)$ is the free $(\colorcd)$-colored $\G$-operad in $\Set$ generated by $\T\wedge\Q$.  The quotient consists of the following three types of relations $\sim$:
\begin{enumerate}
\item For each color $c\in \colorc$, the assignment \[\nicexy{\Q \ar[r]^-{c \otimes -}& \T\tensorg\Q} \qquad \begin{cases} \colord \ni d \mapsto (c;d) \in \colorcd, &\\
\Q\dzeroud \ni q \mapsto c\otimes q \in (\T\tensorg\Q)\sbinom{(c;d_0)}{(c;\ud)}\end{cases}\] defines a morphism in $\gopset$, the category of $\G$-operads in $\Set$ in Definition \ref{def:g-operads}.
\item For each color $d \in \colord$, the assignment \[\nicexy{\T \ar[r]^-{- \otimes d}& \T\tensorg\Q} \qquad \begin{cases} \colorc \ni c \mapsto (c;d) \in \colorcd, &\\
\T\czerouc \ni t \mapsto t\otimes d \in (\T\tensorg\Q)\sbinom{(c_0;d)}{(\uc;d)}\end{cases}\] defines a morphism in $\gopset$.
\item Suppose $t \in \T\czerouc$ and $q \in \Q\dzeroud$, and $g \in \G(mn)$ is a special element with underlying special $(m,n)$-permutation $\gbar=\sigma_{m,n} \in S_{mn}$, where $m=|\ud|$ and $n=|\uc|$.  Then $\sim$ identifies the elements \[\gamma^{\T\tensorg\Q}\bigl(t\otimes d_0; c_1\otimes q, \cdots, c_n\otimes q\bigr)g \in (\T\tensorg\Q)\sbinom{(c_0;d_0)}{(\uc;d_1),\ldots,(\uc;d_m)}\] and \[\gamma^{\T\tensorg\Q}\bigl(c_0\otimes q; t\otimes d_1,\ldots, t\otimes d_m\bigr) \in (\T\tensorg\Q)\sbinom{(c_0;d_0)}{(\uc;d_1),\ldots,(\uc;d_m)}.\]  This is called the \emph{interchange relation}.
\end{enumerate}
\end{enumerate}
\end{definition}

\begin{ginterpretation}
In the interchange relation, the element \[\gamma^{\T\tensorg\Q}\bigl(t\otimes d_0; c_1\otimes q, \cdots, c_n\otimes q\bigr) \in (\T\tensorg\Q)\sbinom{(c_0;d_0)}{(c_1;\ud),\ldots,(c_n;\ud)}\] may be visualized as follows.
\begin{center}\begin{tikzpicture}
\matrix[row sep=.1cm, column sep=1.2cm]{
& \node [plain, label=below:$...$] (v) {\footnotesize{$td_0$}}; &\\
\node [plain, label=below:$...$] (u1) {\footnotesize{$c_1q$}}; &&
\node [plain, label=below:$...$] (um) {\footnotesize{$c_nq$}};\\};
\draw [outputleg] (v) to node[at end]{\scriptsize{$(c_0;d_0)$}} +(0,.8cm);
\draw [thick] (u1) to node{\scriptsize{$(c_1;d_0)$}} (v);
\draw [thick] (um) to node[swap]{\scriptsize{$(c_n;d_0)$}} (v);
\draw [thick] (u1) to node[below left=.2cm]{\scriptsize{$(c_1;d_1)$}} +(-.8cm,-.6cm);
\draw [thick] (u1) to node[below right=.2cm]{\scriptsize{$(c_1;d_m)$}} +(.8cm,-.6cm);
\draw [thick] (um) to node[below left=.2cm]{\scriptsize{$(c_n;d_1)$}} +(-.8cm,-.6cm);
\draw [thick] (um) to node[below right=.2cm]{\scriptsize{$(c_n;d_m)$}} +(.8cm,-.6cm);
\end{tikzpicture}
\end{center}
We abbreviated $t\otimes d_0$ to $td_0$ and $c_i \otimes q$ to $c_iq$ to save space.  This is the $2$-level tree in Example \ref{ex:2-level-tree} with the top vertex decorated by $t \otimes d_0$ and the bottom vertices decorated by the $c_i \otimes q$'s.  Similarly, one can visualize the other element \[\gamma^{\T\tensorg\Q}\bigl(c_0\otimes q; t\otimes d_1,\ldots, t\otimes d_m\bigr) \in (\T\tensorg\Q)\sbinom{(c_0;d_0)}{(\uc;d_1),\ldots,(\uc;d_m)}\] that appears in the interchange relation as the $2$-level tree 
\begin{center}\begin{tikzpicture}
\matrix[row sep=.1cm, column sep=1.2cm]{
& \node [plain, label=below:$...$] (v) {\footnotesize{$c_0q$}}; &\\
\node [plain, label=below:$...$] (u1) {\footnotesize{$td_1$}}; &&
\node [plain, label=below:$...$] (um) {\footnotesize{$td_m$}};\\};
\draw [outputleg] (v) to node[at end]{\scriptsize{$(c_0;d_0)$}} +(0,.8cm);
\draw [thick] (u1) to node{\scriptsize{$(c_0;d_1)$}} (v);
\draw [thick] (um) to node[swap]{\scriptsize{$(c_0;d_m)$}} (v);
\draw [thick] (u1) to node[below left=.2cm]{\scriptsize{$(c_1;d_1)$}} +(-.8cm,-.6cm);
\draw [thick] (u1) to node[below right=.2cm]{\scriptsize{$(c_n;d_1)$}} +(.8cm,-.6cm);
\draw [thick] (um) to node[below left=.2cm]{\scriptsize{$(c_1;d_m)$}} +(-.8cm,-.6cm);
\draw [thick] (um) to node[below right=.2cm]{\scriptsize{$(c_n;d_m)$}} +(.8cm,-.6cm);
\end{tikzpicture}
\end{center}
with $c_0\otimes q$ on top and the $t\otimes d_j$'s at the bottom.  Since the underlying permutation $\gbar$ is the special $(m,n)$-permutation $\sigma_{m,n}$, we have \[\bigl((c_1;\ud),\ldots,(c_n;\ud)\bigr)\gbar = \bigl((\uc;d_1),\ldots,(\uc;d_m)\bigr).\]  In other words, $\gbar$ takes the input profile of the first element above to the input profile of the other element.  This accounts for the profiles in the interchange relation.\dqed
\end{ginterpretation}

\begin{example}[$0$-ary Interchange]
If $|\ud|=0$ (so $q \in \Q\dzeroempty$ and $g \in \G(0)$), then the interchange relation becomes \[\gamma^{\T\tensorg\Q}\bigl(t\otimes d_0; c_1\otimes q, \cdots, c_n\otimes q\bigr)g \sim c_0\otimes q \in (\T\tensorg\Q)\sbinom{(c_0;d_0)}{\varnothing}.\]  Similarly, if $|\uc|=0$ (so $t \in \T\czeroempty$ and $g \in \G(0)$), then the interchange relation becomes \[(t\otimes d_0)g \sim \gamma^{\T\tensorg\Q}\bigl(c_0\otimes q; t\otimes d_1,\ldots, t\otimes d_m\bigr) \in (\T\tensorg\Q)\sbinom{(c_0;d_0)}{\varnothing}.\]  In particular, if $|\uc|=0=|\ud|$, then the interchange relation becomes \[(t\otimes d_0)g \sim c_0 \otimes q \in (\T\tensorg\Q)\sbinom{(c_0;d_0)}{\varnothing}.\]  For example, if $\G$ is the planar group operad $\P$, the symmetric group operad $\S$, the braid group operad $\B$, the ribbon group operad $\R$, or the cactus group operad $\Cac$, then $\G(0)$ is the trivial group.  So the last interchange relation becomes $t \otimes d_0 \sim c_0 \otimes q$.\dqed
\end{example}

\begin{example}[One-Colored Interchange]\label{ex:interchange-one-color}
To understand the interchange relation better, let us consider the one-colored case.  Suppose $\T$ and $\Q$ are both one-colored $\G$-operads in $\Set$, so the $\G$-tensor product $\T\tensorg\Q$ is also a one-colored $\G$-operad.  By the definition of the quotient $\sim$, the $\G$-tensor product is equipped with morphisms \[\nicexy{\T \ar[r] & \T\tensorg\Q & \Q \ar[l]}\]of one-colored $\G$-operads.  By composition, each $(\T\tensorg\Q)$-algebra $X$ has an induced $\T$-algebra structure \[\nicexy{X^{\times n} \ar[r]^-{t} & X} \qquad (t \in \T(n),\,n\geq 0)\] and an induced $\Q$-algebra structure \[\nicexy{X^{\times m} \ar[r]^-{q} & X} \qquad (q \in \Q(m),\, m\geq 0).\]  The interchange relation implies that these two structures commute in the following sense.  For any $t \in \T(n)$, $q\in\Q(m)$, and $g\in \G(mn)$ a special element with underlying special permutation $\gbar=\sigma_{m,n}$, as in the interchange relation in the one-colored case, the diagram
\[\nicexy{(X^{\times n})^{\times m} \ar[d]_-{\sigma_{m,n}} \ar[rr]^-{(t)^{\times m}} && X^{\times m} \ar[d]^-{q}\\ (X^{\times m})^{\times n} \ar[r]^-{(q)^{\times n}} & X^{\times n} \ar[r]^-{t} & X}\]
is commutative.  The special $(m,n)$-permutation $\sigma_{m,n} \in S_{mn}$ permutes the factors in $(X^{\times n})^{\times m}$ from the left.  In other words, up to the special permutation $\sigma_{m,n}$, the action by $t$ and the action by $q$ commute in a suitable sense.\dqed
\end{example}

\begin{theorem}\label{goperad-symmetric-monoidal}
For each action operad $(\G,\omega)$,\[\bigl(\gopset, \tensorg, \I^{\G}\bigr)\] is a\index{G-operad@$\G$-operad!symmetric monoidal category}\index{symmetric monoidal category!of G-operads@of $\G$-operads} symmetric monoidal category, where $\I^{\G}$ is the one-colored $\G$-operad in $\Set$ in \eqref{g-unit}.
\end{theorem}

\begin{proof}
For a $\colorc$-colored $\G$-operad $\U$, a $\colord$-colored $\G$-operad $\Q$, and an $\colore$-colored $\G$-operad $\T$, the associativity isomorphism \[(\U\tensorg\Q)\tensorg\T \cong \U\tensorg(\Q\tensorg\T)\] of $(\colorc\times\colord\times\colore)$-colored $\G$-operads is induced at the smash product level by the assignments \[\begin{cases}
(u\otimes d) \otimes e \mapsto u\otimes (d;e) & \text{for $u \in \U$, $d \in \colord$, $e \in \colore$},\\
(c\otimes q)\otimes e \mapsto c\otimes (q\otimes e) & \text{for $c \in \colorc$, $q \in \Q$, $e\in\colore$},\\
(c;d)\otimes t \mapsto c \otimes (d\otimes t) & \text{for $c\in\colorc$, $d\in\colord$, $t \in T$.}\end{cases}\]
The commutativity isomorphism \[\U\tensorg\Q \cong \Q\tensorg\U\] is induced at the smash product level by the assignments \[\begin{cases}(c;d) \mapsto (d;c) & \text{for $c\in\colorc$, $d\in\colord$},\\
u\otimes d \mapsto d \otimes u & \text{for $u\in\U$, $d \in \colord$},\\ c\otimes q \mapsto q\otimes c & \text{for $c\in\colorc$, $q\in \Q$}.\end{cases}\]

For the left unit isomorphism $\I^{\G}\tensorg\U\cong\U$, recall that the $\{*\}$-colored $\G$-operad $\I^{\G}$ is $\G(1)$ at level $1$ and is empty at all other levels.  The left unit isomorphism is induced at the smash product level by the assignments 
\[\begin{cases} (*;c) \mapsto c & \text{for $c\in\colorc$},\\
*\otimes u \mapsto u & \text{for $u\in\U$},\\
\sigma \otimes c \mapsto \operadunit_c^{\U} \cdot \sigmainv\in \U\cc & \text{for $\sigma \in \G(1)$ and $c\in\colorc$},\end{cases}\]
where $\operadunit_c^{\U} \in \U\cc$ is the $c$-colored unit in $\U$.  The right unit isomorphism is defined similarly.
\end{proof}

\begin{example}[Planar Tensor Product]\label{ex:planar-tensor}
Suppose $\G=\P$ is the planar group operad in Example \ref{ex:trivial-group-operad}.  For a $\colorc$-colored planar operad $\T$ and a $\colord$-colored planar operad $\Q$, the $\P$-tensor product $\T\tensorp\Q$ is called the \index{planar tensor product}\index{planar operad!monoidal product}\emph{planar tensor product} of $\T$ and $\Q$.  The category $\popset$ of planar operads in $\Set$ is symmetric monoidal with respect to the planar tensor product.\dqed
\end{example}

\begin{example}[Boardman-Vogt Tensor Product of Symmetric Operads]\label{ex:bv-tensor}
Suppose $\G=\S$ is the symmetric group operad with the identity augmentation over itself.  For a $\colorc$-colored symmetric operad $\T$ and a $\colord$-colored symmetric operad $\Q$, the $\S$-tensor product $\T\tensors\Q$ is exactly the \index{Boardman-Vogt tensor product}Boardman-Vogt tensor product of $\T$ and $\Q$ \cite{boardman-vogt} (Chapter II.3).  So our $\G$-tensor product is a generalization of the Boardman-Vogt tensor product of symmetric operads to $\G$-operads.  The fact that the category $\sopset$ of symmetric operads in $\Set$ is symmetric monoidal with respect to the Boardman-Vogt tensor product $\tensors$ was essentially proved in \cite{boardman-vogt} (Proposition 2.18).  See also \cite{mw07} (Section 5.1).\dqed
\end{example}

\begin{example}[Braided Tensor Product]\label{ex:braided-tensor}
Suppose $\G=\B$ is the braid group operad in Definition \ref{def:braid-group-operad}.  For a $\colorc$-colored braided operad $\T$ and a $\colord$-colored braided operad $\Q$, the $\B$-tensor product $\T\tensorb\Q$ is called the \index{braided tensor product}\index{braided operad!monoidal product}\emph{braided tensor product} of $\T$ and $\Q$.  The category $\bopset$ of braided operads in $\Set$ is symmetric monoidal with respect to the braided tensor product.\dqed
\end{example}

\begin{example}[Pure Braided Tensor Product]\label{ex:pure-braided-tensor}
Suppose $\G=\PB$ is the pure braid group operad in Definition \ref{def:pure-braid-group-operad}.  For a $\colorc$-colored pure braided operad $\T$ and a $\colord$-colored pure braided operad $\Q$, the $\PB$-tensor product $\T\tensorpb\Q$ is called the\index{pure braided tensor product}\index{pure braided operad!monoidal product} \emph{pure braided tensor product} of $\T$ and $\Q$.  The category $\pbopset$ of pure braided operads in $\Set$ is symmetric monoidal with respect to the pure braided tensor product.\dqed
\end{example}

\begin{example}[Ribbon Tensor Product]\label{ex:ribbon-tensor}
Suppose $\G=\R$ is the ribbon group operad in Definition \ref{def:ribbon-group-operad}.  For a $\colorc$-colored ribbon operad $\T$ and a $\colord$-colored ribbon operad $\Q$, the $\R$-tensor product $\T\tensorr\Q$ is called the \index{ribbon tensor product}\index{ribbon operad!monoidal product}\emph{ribbon tensor product} of $\T$ and $\Q$.  The category $\ropset$ of ribbon operads in $\Set$ is symmetric monoidal with respect to the ribbon tensor product.\dqed
\end{example}

\begin{example}[Pure Ribbon Tensor Product]\label{ex:pure-ribbon-tensor}
Suppose $\G=\PR$ is the pure ribbon group operad in Definition \ref{def:pure-ribbon-group-operad}.  For a $\colorc$-colored pure ribbon operad $\T$ and a $\colord$-colored pure ribbon operad $\Q$, the $\PR$-tensor product $\T\tensorpr\Q$ is called the\index{pure ribbon tensor product}\index{pure ribbon operad!monoidal product} \emph{pure ribbon tensor product} of $\T$ and $\Q$.  The category $\propset$ of pure ribbon operads in $\Set$ is symmetric monoidal with respect to the pure ribbon tensor product.\dqed
\end{example}

\begin{example}[Cactus Tensor Product]\label{ex:cactus-tensor}
Suppose $\G=\Cac$ is the cactus group operad in Definition \ref{def:cactus-group-operad}.  For a $\colorc$-colored cactus operad $\T$ and a $\colord$-colored cactus operad $\Q$, the $\Cac$-tensor product $\T\tensorcac\Q$ is called the\index{cactus tensor product}\index{cactus operad!monoidal product} \emph{cactus tensor product} of $\T$ and $\Q$.  The category $\cacopset$ of cactus operads in $\Set$ is symmetric monoidal with respect to the cactus tensor product.\dqed
\end{example}

\begin{example}[Pure Cactus Tensor Product]\label{ex:pure-cactus-tensor}
Suppose $\G=\PCac$ is the pure cactus group operad in Definition \ref{def:pure-cactus-group-operad}.  For a $\colorc$-colored pure cactus operad $\T$ and a $\colord$-colored pure cactus operad $\Q$, the $\PCac$-tensor product $\T\tensorpcac\Q$ is called the\index{pure cactus tensor product}\index{pure cactus operad!monoidal product}  \emph{pure cactus tensor product} of $\T$ and $\Q$.  The category $\pcacopset$ of pure cactus operads in $\Set$ is symmetric monoidal with respect to the pure cactus tensor product.\dqed
\end{example}

\section{Closed Structure}\label{sec:goperad-monoidal-closed}

The purpose of this section is to extend the symmetric monoidal structure on the category $\gopset$ of all $\G$-operads in $\Set$ to a symmetric monoidal \emph{closed} structure.  Recall that a monoidal category $\C$ is \index{monoidal category!closed}\emph{closed} if for each object $X$, the functor \[-\otimes X : \C \to \C\] admits a right adjoint.  First we define the necessary objects that will show that the category of $\G$-operads in $\Set$ is a symmetric monoidal closed category.

\begin{definition}\label{def:g-natural-transformation}
Suppose\index{G-operad@$\G$-operad!internal hom}
\begin{itemize}
\item $(\T,\gamma^{\T},\operadunit^{\T})$ is a $\colorc$-colored $\G$-operad and 
\item $(\Q,\gamma^{\Q},\operadunit^{\Q})$ is a $\colord$-colored $\G$-operad in $\Set$.  
\end{itemize}
Define the colored $\G$-operad $\Homg(\T,\Q)$ in $\Set$ as follows.
\begin{description}
\item[Colors] Its set of colors is $\gopset(\T,\Q)$, i.e., the set of morphisms $\T \to \Q$ of $\G$-operads in $\Set$.
\item[Entries] For $f_0,f_1,\ldots,f_n \in \gopset(\T,\Q)$ with $n\geq 0$, an element \[\theta \in \Homg(\T,\Q)\fzerouf \withspace \uf=(f_1,\ldots,f_n)\] is an assignment 
\[\theta = \Bigl\{\colorc \ni c \mapsto \theta_c \in \Q\fzerocufc\Bigr\},\] where \[\uf c = (f_1c,\ldots,f_nc) \in \Profd,\] that is \emph{natural} in the following sense.  For each 
\begin{itemize}
\item $\czerouc \in \Profcc$ with $\uc=(c_1,\ldots,c_m) \in \Profc$, 
\item $t \in \T\czerouc$, and 
\item special element $g \in \G(nm)$ with underlying permutation $\gbar=\sigma_{n,m} \in S_{nm}$, 
\end{itemize}
the naturality condition
\begin{equation}\label{goperad-hom-naturality}
\gamma^{\D}\bigl(f_0t; \theta_{\uc}\bigr)g = \gamma^{\D}\bigl(\theta_{c_0}; \uf t\bigr)\in \Q\sbinom{f_0c_0}{f_1\uc,\ldots,f_n\uc}
\end{equation}
holds, where \[\begin{split}f_it &\in \Q\sbinom{f_ic_0}{f_i\uc} \forspace 0\leq i \leq n,\\
\uf t &= \bigl(f_1t, \ldots,f_nt\bigr) \in \Q\sbinom{f_1c_0}{f_1\uc}\times\cdots\times\Q\sbinom{f_nc_0}{f_n\uc},\\
\theta_{\uc} &= \bigl(\theta_{c_1},\ldots,\theta_{c_m}\bigr)\in \Q\sbinom{f_0c_1}{\uf c_1}\times\cdots\times\Q\sbinom{f_0c_m}{\uf c_m}.\end{split}\]
\item[Units] The $f_0$-colored unit is the element \[\operadunit_{f_0} = \Bigl\{\colorc \ni c \mapsto \operadunit^{\Q}_{f_0c}\in \Q\sbinom{f_0c}{f_0c}\Bigr\} \in \Homg(\T,\Q)\sbinom{f_0}{f_0},\] in which $ \operadunit^{\Q}_{f_0c}$ is the $f_0c$-colored unit in $\Q$.
\item[Composition] Suppose 
\begin{itemize}
\item $\theta \in \Homg(\T,\Q)\fzerouf$ is as above with $1 \leq i \leq n$ and
\item $\psi \in \Homg(\T,\Q)\sbinom{f_i}{\uh}$ for some $\uh \in \Prof\bigl(\gopset(\T,\Q)\bigr)$.  
\end{itemize}
The $\compi$-composition \[\nicexy{\Homg(\T,\Q)\fzerouf \times \Homg(\T,\Q)\sbinom{f_i}{\uh} \ar[r]^-{\compi} & \Homg(\T,\Q)\sbinom{f_0}{\uf \compi \uh}}\] sends the pair $(\theta,\psi)$ to the element \[\theta\compi\psi = \Bigl\{\colorc \ni c \mapsto \theta_c \compi \psi_c \in \Q\sbinom{f_0c}{(\uf\compi\uh)c}\Bigr\}\] in which $\theta_c \compi \psi_c$ is the $\compi$-composition in $\Q$.
\item[$\G$-Equivariance] Suppose $\theta$ is as above and $\sigma \in \G(n)$ with underlying permutation $\sigmabar \in S_n$.  The $\gopset(\T,\Q)$-colored $\G$-sequence structure \[\nicexy{\Homg(\T,\Q)\fzerouf \ar[r]^-{\sigma} & \Homg(\T,\Q)\fzeroufsigmabar}\] sends $\theta$ to the element \[\theta\sigma = \Bigl\{\colorc \ni c \mapsto \theta_c\sigma \in \Q\sbinom{f_0c}{(\uf\sigmabar)c}\Bigr\},\] where the $\sigma$ in $\theta_c\sigma$ indicates the $\colord$-colored $\G$-sequence structure in $\Q$.
\end{description}
\end{definition}  

\begin{remark}
In the naturality condition \eqref{goperad-hom-naturality}, the element $\gamma^{\D}\bigl(f_0t; \theta_{\uc}\bigr)$ belongs to $\Q\sbinom{f_0c_0}{\uf c_1,\ldots,\uf c_m}$.  The special element $g \in \G(nm)$ yields an isomorphism \[\nicexy{\Q\sbinom{f_0c_0}{\uf c_1,\ldots,\uf c_m} \ar[r]^-{g}_-{\cong} & \Q\sbinom{f_0c_0}{f_1\uc,\ldots,f_n\uc}},\] which accounts for the profiles in the naturality condition.\dqed
\end{remark}

\begin{lemma}\label{hompq-goperad}
$\Homg(\T,\Q)$ is a $\gopset(\T,\Q)$-colored $\G$-operad in $\Set$.
\end{lemma}

\begin{proof}
The naturality \eqref{goperad-hom-naturality} of the $f_0$-colored unit $\operadunit_{f_0}$, the $\compi$-composition $\theta\compi\psi$, and $\theta\sigma$ follows from the unity axiom, the associativity axiom, and the equivariance axioms in the $\colord$-colored $\G$-operad $\Q$.  Similarly,  since the colored units, the $\compi$-composition, and the $\G$-equivariant structure in $\Homg(\T,\Q)$ are defined in terms of those in $\Q$, the required $\gopset(\T,\Q)$-colored $\G$-operad axioms follow from the $\colord$-colored $\G$-operad axioms in $\Q$.
\end{proof}

\begin{theorem}\label{goperad-monoidal-closed}
The category \[\bigl(\gopset, \tensorg, \I^{\G}\bigr)\] is a\index{G-operad@$\G$-operad!symmetric monoidal closed structure} symmetric monoidal closed category in which the right adjoint of $-\tensorg\T$ is $\Homg(\T,-)$ for each $\T \in \gopset$.  In particular, $\tensorg$ commutes with colimits in both variables separately.
\end{theorem}

\begin{proof}
We observed in Theorem \ref{goperad-symmetric-monoidal} that the category $\bigl(\gopset, \tensorg, \I^{\G}\bigr)$ is symmetric monoidal.  To check that it is closed, suppose $\U$ is a $\colorc$-colored $\G$-operad, $\Q$ is a $\colord$-colored $\G$-operad, and $\T$ is an $\colore$-colored $\G$-operad in $\Set$.  The desired natural bijection
\begin{equation}\label{gopset-enriched-iso}
\gopset\bigl(\U\tensorg\Q,\T\bigr) \cong \gopset\bigl(\U,\Homg(\Q,\T)\bigr)
\end{equation}
is given as follows.

Suppose $\varphi \in \gopset\bigl(\U\tensorg\Q,\T\bigr)$, so $\varphi : \U\tensorg\Q \to \T$ is a morphism of $\G$-operads.  The associated morphism \[\phi : \U \to \Homg(\Q,\T)\] of $\G$-operads is defined as follows.
\begin{itemize}
\item For each color $c \in \colorc$, $\phi(c) : \Q \to \T$ is the composite morphism \[\nicexy{\Q \ar[r]^-{c \otimes -} & \U\tensorg\Q \ar[r]^-{\varphi} & \T \ar@{<-}`u[ll]`[ll]_-{\phi(c)}[ll]}.\]
\item For each element $u \in \U\czerouc$ with $\uc=(c_1,\ldots,c_m) \in \Profc$, the element \[\phi(u) \in \Homg(\Q,\T)\sbinom{\phi c_0}{\phi\uc}\] is the assignment
\[\colord \ni d \mapsto \phi(u)(d) = \varphi(u \otimes d)
\in \T\sbinom{(\phi c_0)(d)}{(\phi\uc)(d)} = \T\sbinom{\varphi(c_0;d)}{\varphi(c_1;d),\ldots,\varphi(c_m;d)}.\]
\end{itemize}
One can check that $\phi$ indeed defines a morphism $\U \to \Homg(\Q,\T)$.

Conversely, suppose $\phi : \U \to \Homg(\Q,\T)$ is a morphism of $\G$-operads. The associated morphism \[\varphi : \U\tensorg\Q \to \T\] is defined as follows.
\begin{itemize}
\item For each pair of colors $(c;d) \in \colorcd$, we define \[\varphi(c;d) = \phi(c)(d) \in \colore.\]  This makes sense because $\phi(c) : \Q \to \T$ is a morphism of $\G$-operads, so on color sets it is a function $\phi(c) : \colord \to \colore$.
\item For $u \in \U\czerouc$ with $\uc=(c_1,\ldots,c_n)\in\Profc$ and $d \in \colorc$, we define \[\varphi(u \otimes d)=\phi(u)(d) \in \T\sbinom{\phi(c_0)(d)}{(\phi\uc)(d)} = \T\sbinom{\varphi(c_0;d)}{\varphi(c_1;d),\ldots,\varphi(c_m;d)}.\]
\item For $c\in\colorc$ and $q\in \Q\dzeroud$ with $\ud=(d_1,\ldots,d_n)\in\Profd$, we define \[\varphi(c\otimes q)=\phi(c)(q) \in \T\sbinom{\phi(c)(d_0)}{\phi(c)(\ud)} = \T\sbinom{\varphi(c;d_0)}{\varphi(c;d_1),\ldots,\varphi(c;d_n)}\]
\end{itemize}
One can check that $\varphi$ indeed defines a morphism $\U\tensorg\Q \to \T$ and that the two constructions above are inverses of each other.
\end{proof}

\begin{example}
When $\G=\S$ is the symmetric group operad, the fact that the category $\sopset$ of symmetric operads in $\Set$ is symmetric monoidal closed was proved in Theorem 4.1.3 in \cite{mt}.\dqed
\end{example}

\begin{corollary}
There is a natural isomorphism \[\Homg\bigl(\U\tensorg\Q,\T\bigr) \cong \Homg\bigl(\U,\Homg(\Q,\T)\bigr)\] in $\gopset$ whose restriction to color sets is the natural bijection \eqref{gopset-enriched-iso}.
\end{corollary}

\begin{proof}
This is a purely formal consequence of the adjunction between $-\tensorg\U$ and $\Homg(\U,-)$.  More explicitly, suppose \[\varphi_0,\, \varphi_1,\ldots,\varphi_n \in \gopset\bigl(\U\tensorg\Q,\T\bigr)\] and $\uvarphi = (\varphi_1,\ldots,\varphi_n)$. 
Using the notation in the proof of Theorem \ref{goperad-monoidal-closed}, their images across the bijection \eqref{gopset-enriched-iso} are \[\phi_0,\,\phi_1,\ldots,\phi_n \in \gopset\bigl(\U,\Homg(\Q,\T)\bigr)\] with $\uphi = (\phi_1,\ldots,\phi_n)$.  For an element \[\theta = \Bigl\{\colorcd\ni (c;d) \mapsto \theta_{(c;d)} \in \T\sbinom{\varphi_0(c;d)}{\uvarphi(c;d)}\Bigr\} \in \Homg\bigl(\U\tensorg\Q,\T\bigr)\sbinom{\varphi_0}{\uvarphi},\]
the corresponding element across the desired isomorphism is \[\theta' = \Biggl\{\colorc \ni c \mapsto \theta'_c = \Bigl\{\colord \ni d \mapsto \theta_{(c;d)}\Bigr\} \in \Homg(\Q,\T)\sbinom{\phi_0c}{\uphi c}\Biggr\}\] in $\Homg\bigl(\U,\Homg(\Q,\T)\bigr)\sbinom{\phi_0}{\uphi}$.
\end{proof}

\section{Comparing Symmetric Monoidal Structures}\label{sec:comparing-sms}

In this section, we compare the symmetric monoidal categories of $\G$-operads for different action operads $\G$.  Recall from Theorem \ref{goneopm-gtwoopm} that each morphism $\varphi : \Gone \to \Gtwo$ of action operads induces a functor \[\nicexy{\goneopm & \gtwoopm \ar[l]_-{\varphi^*}}\] that preserves the underlying colored planar operads.

\begin{theorem}\label{goperad-monoidal-functor}\index{change of action operads!symmetric monoidal functor}
For each morphism\index{G-operad@$\G$-operad!change of action operads} $\varphi : \Gone \to \Gtwo$ of action operads, the induced functor \[\nicexy{\Bigl(\goneopset,\tensorgone,\I^{\Gone}\Bigr) & \Bigl(\gtwoopset,\tensorgtwo,\I^{\Gtwo}\Bigr) \ar[l]_-{\varphi^*}}\] is a \index{symmetric monoidal functor!between G-operad categories@between $\G$-operad categories}symmetric monoidal functor.\end{theorem}

\begin{proof}
The image $\varphi^*\I^{\Gtwo}$ is the $\{*\}$-colored $\Gone$-operad that is $\Gtwo(1)$ in level $1$ and is empty in all other entries, with $\Gone(1)$ acting on the right via the morphism $\varphi$.  The structure morphism \[\nicexy{\I^{\Gone} \ar[r]^-{\varphi^*_0} & \varphi^*\I^{\Gtwo}}\] is determined by the morphism $\varphi : \Gone(1) \to \Gtwo(1)$ in level $1$.

Suppose $\T$ is a $\colorc$-colored $\Gtwo$-operad and $\Q$ is a $\colord$-colored $\Gtwo$-operad in $\Set$.  For the functor $\varphi^*$ to be a symmetric monoidal functor, it must be equipped with a morphism
\begin{equation}\label{psipq}
\nicexy{\varphi^*\T \tensorgone \varphi^*\Q \ar[r]^-{\psi} & \varphi^*\bigl(\T\tensorgtwo\Q\bigr)} \in \goneopset
\end{equation}
that is natural in $\T$ and $\Q$.  To define this morphism $\psi$, first observe that at the smash product level there is an equality of graded sets\[\varphi^*\T \wedge \varphi^*\Q = \T\wedge\Q \in \Set^{\Profcdcd}\] because $\varphi^*\T$ has the same underlying sets as $\T$ and similarly for $\varphi^*\Q$.  So there is a morphism of graded sets \[\nicexy{\varphi^*\T \wedge\varphi^*\Q \ar[r] & \varphi^*\bigl(\T\tensorgtwo\Q\bigr)} \in \Set^{\Profcdcd}\] that sends each element $t \otimes d$ or $c \otimes q$ (with $t\in \T$, $q \in \Q$, $c\in\colorc$, and $d \in\colord$) to a generator with the same name on the right-hand side.  

Recall from Theorem \ref{gopcm-algebra} with $\M=\Set$ that $\goneopcd$ is the $(\Profcdcd$)-colored symmetric operad in $\Set$ whose category of algebras is the category $\goneopcdset$ of $(\colorcd)$-colored $\Gone$-operads in $\Set$.  In other words, \[\goneopcd \circs-\] is the free $(\colorcd)$-colored $\Gone$-operad functor, where $\circs$ is the $(\colorcd)$-colored symmetric circle product in \eqref{symmetric-circle-product}.  So the previous morphism induces a morphism of $(\colorcd)$-colored $\Gone$-operads\[\nicexy{\goneopcd \circs\bigl(\varphi^*\T \wedge\varphi^*\Q\bigr) \ar[r] & \varphi^*\bigl(\T\tensorgtwo\Q\bigr)} \in \goneopcdset.\] This morphism respects the relations $\sim$ that define the $\Gone$-tensor product $\tensorgone$ because:
\begin{itemize}\item The $\Gone$-equivariant structure in $\varphi^*\T$ and $\varphi^*\Q$ acts through $\varphi$.
\item $\varphi$ respects the augmentations of $\Gone$ and $\Gtwo$ in the sense of Lemma \ref{rk:varphi-augmentation}.
\end{itemize}
Therefore, the previous morphism induces a morphism $\psi$ that makes the following diagram commutative. \[\nicexy@C-.5cm{\goneopcd \circs\bigl(\varphi^*\T \wedge\varphi^*\Q\bigr) \ar[r] \ar[d]_-{\mathrm{quotient}} & \varphi^*\bigl(\T\tensorgtwo\Q\bigr)\\
\frac{\goneopcd \circs\bigl(\varphi^*\T \wedge\varphi^*\Q\bigr)}{\sim} \ar@{=}[r] &  \varphi^*\T \tensorgone \varphi^*\Q \ar@{-->}[u]_-{\psi} &}\] 
One can check that this morphism $\psi$ satisfies the conditions for $\varphi^*$ to be a symmetric monoidal functor.
\end{proof}

Recall that a \index{comonoidal functor}\index{monoidal category!comonoidal functor}\emph{comonoidal functor}, sometimes called a \emph{lax} comonoidal functor or an \index{oplax monoidal functor}\emph{oplax} monoidal functor in the literature, \[F : (\C,\otimes,\tensorunit) \to (\C',\otimes',\tensorunit')\] between two monoidal categories is a functor $F : \C \to \C'$ equipped with 
\begin{enumerate}[label=(\roman*)]
\item a morphism $F^0 : F\tensorunit \to \tensorunit' \in \C'$ and 
\item a morphism \[F^2(x,y) : F(x\otimes y)\to Fx \otimes' Fy\in \C'\] that is natural in the objects $x$ and $y$, 
\end{enumerate}
satisfying suitable coassociativity and counity axioms that are dual to those of a monoidal functor in Definition \ref{def:monoidal-functor}.  A \index{symmetric comonoidal functor}\emph{symmetric comonoidal functor} between symmetric monoidal categories is a comonoidal functor in which $F^2$ is compatible with the symmetry isomorphisms in the obvious sense.

\begin{corollary}\label{goperad-comonoidal}\index{change of action operads!symmetric comonoidal functor}
For each morphism $\varphi : \Gone \to \Gtwo$ of action operads, the induced functor \[\nicexy{\Bigl(\goneopset,\tensorgone,\I^{\Gone}\Bigr) \ar[r]^-{\varphi_!} & \Bigl(\gtwoopset,\tensorgtwo,\I^{\Gtwo}\Bigr)}\] in Theorem \ref{goneopm-gtwoopm} is a symmetric comonoidal functor.
\end{corollary}

\begin{proof}
The functor $\varphi^*$ is symmetric monoidal by Theorem \ref{goperad-monoidal-functor}.  So its left adjoint $\varphi_!$ is symmetric comonoidal by Theorem 1.2 in \cite{kelly-adjunction}.
\end{proof}

\begin{remark}
For $\Gone$-operads $\O$ and $\Q$, the comonoidal structure morphism of $\varphi_!$ is the composite \[\nicexy{\varphi_!\bigl(\O\tensorgone\Q\bigr) \ar[d]_-{\varphi_!(\epsilon_{\O},\epsilon_{\Q})} \ar[r]^-{\varphi_!^2} & \varphi_!\O\tensorgtwo \varphi_!\Q\\ 
\varphi_!\bigl(\varphi^*\varphi_!\O \tensorgone \varphi^*\varphi_!\Q\bigr) \ar[r]^-{\varphi_!(\varphi^*_2)} & \varphi_!\varphi^*\bigl(\varphi_!\O\tensorgtwo \varphi_!\Q\bigr) \ar[u]_-{\eta}}\] with:
\begin{itemize}
\item each $\epsilon$ the unit of the adjunction $\varphi_! \dashv \varphi^*$;
\item $\varphi^*_2$ the monoidal structure morphism of $\varphi^*$;
\item $\eta$ the counit of the adjunction $\varphi_! \dashv \varphi^*$.  
\end{itemize}
The compatibility of $\varphi_!^2$ with the symmetry isomorphisms follows from the same property of $\varphi^*_2$.\dqed
\end{remark}

\begin{example}[Comparison with the Boardman-Vogt Tensor Product]\label{ex:bv-tensor-compare}
As discussed in Example \ref{ex:all-augmented-group-operads}, there is a commutative diagram 
\[\nicexy{\P\ar@{=}[d] \ar[r]^-{\iota} & \PR \ar[r]^-{\iota} & \R \ar[r]^-{\pi} & \S \ar@{=}[d]\\
\P \ar[r]^-{\iota} & \PB \ar[r]^-{\iota} & \B \ar[r]^-{\pi} & \S\\
\P \ar[r]^-{\iota} \ar@{=}[u] & \PCac \ar[r]^-{\iota} & \Cac \ar[r]^-{\pi} & \S\ar@{=}[u]}\]
of morphisms of action operads.  These morphisms\index{Boardman-Vogt tensor product} induce adjunctions 
\[\begin{small}\nicexy@C-.3cm{\bigl(\popset,\tensorp\bigr) \ar@{=}[d]
\ar@<2pt>[r]^-{\iota_!} & \bigl(\propset,\tensorpr\bigr) \ar@<2pt>[r]^-{\iota_!} \ar@<2pt>[l]^-{\iota^*} & \bigl(\ropset,\tensorr\bigr) \ar@<2pt>[r]^-{\pi_!} \ar@<2pt>[l]^-{\iota^*} & \bigl(\sopset,\tensors\bigr) \ar@<2pt>[l]^-{\pi^*} \ar@{=}[d]\\ 
\bigl(\popset,\tensorp\bigr) \ar@<2pt>[r]^-{\iota_!} & 
\bigl(\pbopset,\tensorpb\bigr) \ar@<2pt>[r]^-{\iota_!} \ar@<2pt>[l]^-{\iota^*} & \bigl(\bopset,\tensorb\bigr) \ar@<2pt>[r]^-{\pi_!} \ar@<2pt>[l]^-{\iota^*} & \bigl(\sopset,\tensors\bigr) \ar@<2pt>[l]^-{\pi^*}\\
\bigl(\popset,\tensorp\bigr) \ar@{=}[u] \ar@<2pt>[r]^-{\iota_!} & \bigl(\pcacopset,\tensorpcac\bigr)  \ar@<2pt>[r]^-{\iota_!} \ar@<2pt>[l]^-{\iota^*} & \bigl(\cacopset,\tensorcac\bigr) \ar@<2pt>[r]^-{\pi_!} \ar@<2pt>[l]^-{\iota^*} & \bigl(\sopset,\tensors\bigr) \ar@<2pt>[l]^-{\pi^*} \ar@{=}[u]}
\end{small}\]
in which $\tensors$ is the Boardman-Vogt tensor product of colored symmetric operads.  By Theorem \ref{goperad-monoidal-functor} every right adjoint in this diagram--i.e., the functors denoted by $(-)^*$--is a symmetric monoidal functor. By Corollary \ref{goperad-comonoidal} every left adjoint in this diagram--i.e., the functors denoted by $(-)_!$--is a symmetric comonoidal functor.\dqed
\end{example}

\section{Non-Strong Monoidality}\label{sec:non-strong}

Recall that a \index{strong monoidal functor}\index{monoidal functor!strong}\emph{strong} monoidal functor is a monoidal functor $(F,F_2,F_0)$ such that the structure morphisms $F_0$ and $F_2$ are isomorphisms.  It is sometimes called a \emph{tensor functor} in the literature \cite{joyal-street}.  In this section, we observe that the symmetric monoidal functors in Example \ref{ex:bv-tensor-compare} are not strong monoidal.  

\begin{proposition}\label{tensorp-tensors}
The symmetric monoidal functor \[\nicexy{\bigl(\popset,\tensorp\bigr) & \bigl(\sopset,\tensors\bigr) \ar[l]_-{\iota^*}}\] induced by the morphism $\iota : \P \to \S$ is \underline{not} strong monoidal.
\end{proposition}

\begin{proof}
It suffices to show that the morphism \[\nicexy{\iota^*\T \tensorp \iota^*\Q \ar[r]^-{\psi} & \iota^*\bigl(\T\tensors\Q\bigr)} \in \popset\] in \eqref{psipq} is not an isomorphism in general.  So suppose both $\T$ and $\Q$ are one-colored symmetric operads with elements $t \in \T(2)$ and $q \in \Q(2)$.  Consider the element \[\alpha = \gamma^{\iota^*(\T\tensors\Q)}\bigl(t\otimes *;*\otimes q,*\otimes q\bigr) \in \iota^*\bigl(\T\tensors\Q\bigr)(4),\] which may be visualized as the following decorated $2$-level tree.
\begin{center}\begin{tikzpicture}
\matrix[row sep=.1cm, column sep=.5cm]{
& \node [plain] (v) {$t$}; &\\ \node [plain] (u1) {$q$}; && \node [plain] (um) {$q$};\\};
\draw [outputleg] (v) to +(0,.8cm);
\draw [thick] (u1) to (v); \draw [thick] (um) to (v);
\draw [thick] (u1) to +(-.8cm,-.6cm); \draw [thick] (u1) to +(.8cm,-.6cm);
\draw [thick] (um) to +(-.8cm,-.6cm); \draw [thick] (um) to +(.8cm,-.6cm);
\end{tikzpicture}
\end{center}
Although $\alpha$ is in the image of the morphism $\psi$, the element \[\alpha(1~3) \in \iota^*\bigl(\T\tensors\Q\bigr)(4)\] is not in the image of $\psi$ for the  transposition $(1~3)\in S_4$.   So $\psi$ is not an isomorphism, and $\iota^*$ is not a strong monoidal functor.  
\end{proof}

\begin{proposition}\label{tensorp-tensorr}
The symmetric monoidal functors \[\nicexy{& \bigl(\cacopset,\tensorcac\bigr) \ar[d]^-{\iota^*}&\\ 
\bigl(\ropset,\tensorr\bigr) \ar[r]^-{\iota^*} & \bigl(\popset,\tensorp\bigr) & \bigl(\bopset,\tensorb\bigr) \ar[l]_-{\iota^*}}\] are not strong monoidal.
\end{proposition}

\begin{proof}
We slightly modify the proof of Proposition \ref{tensorp-tensors} by first taking $\T$ and $\Q$ to be one-colored braided, ribbon, or cactus operads.  Then we replace the permutation $(1~3)$ with a braid, a ribbon, or a cactus whose underlying permutation is $(1~3)$.
\end{proof}

The proofs of Proposition \ref{tensorp-tensors} and Proposition \ref{tensorp-tensorr} above depend on the fact that the morphisms of action operads\[\P \to \S, \quad \P\to\R, \quad \P\to\B, \andspace \P\to\Cac\] are not surjective.  A similar argument establishes the following result.

\begin{proposition}
The symmetric monoidal functors \[\nicexy{& \bigl(\propset,\tensorpr\bigr) \ar`l/3pt[dl]_-{\iota^*} [dl] & \bigl(\ropset,\tensorr\bigr) \ar[l]_-{\iota^*}\\
\bigl(\popset,\tensorp\bigr) \ar@{<-}`d/3pt[dr][dr]^-{\iota^*}
& \bigl(\pbopset,\tensorpb\bigr) \ar[l]_-{\iota^*} & \bigl(\bopset,\tensorb\bigr) \ar[l]_-{\iota^*}\\
& \bigl(\pcacopset,\tensorpcac\bigr)  & \bigl(\cacopset,\tensorcac\bigr) \ar[l]_-{\iota^*}}\] in Example \ref{ex:bv-tensor-compare} are \underline{not} strong monoidal.
\end{proposition}

\begin{proposition}\label{tensors-tensorb}
The symmetric monoidal functors \[\nicexy{\bigl(\pbopset,\tensorpb\bigr) &\\
\bigl(\bopset,\tensorb\bigr) & \bigl(\sopset,\tensors\bigr) \ar[l]_-{\pi^*} \ar`u/3pt[ul][ul]_-{\rho^*}}\] induced by the morphisms $\pi : \B \to \S$ and $\rho : \PB \to \S$ are \underline{not} strong monoidal.  
\end{proposition}

\begin{proof}
In the braided case, it suffices to show that the morphism \[\nicexy{\pi^*\T \tensorb \pi^*\Q \ar[r]^-{\psi} & \pi^*\bigl(\T\tensors\Q\bigr)} \in \bopset\] in \eqref{psipq} is not an isomorphism in general.  Suppose both $\T$ and $\Q$ are one-colored symmetric operads with elements $t \in \T(2)$ and $q\not= q' \in \Q(2)$.  Consider the element \[\beta = \gamma^{\pi^*\T\tensorb\pi^*\Q}\bigl(t\otimes *;*\otimes q,*\otimes q'\bigr) \in \bigl(\pi^*\T\tensorb\pi^*\Q\bigr)(4),\] which may be visualized as the following decorated $2$-level tree.
\begin{center}\begin{tikzpicture}
\matrix[row sep=.1cm, column sep=.5cm]{
& \node [plain] (v) {$t$}; &\\ \node [plain] (u1) {$q$}; && \node [plain] (um) {$q'$};\\};
\draw [outputleg] (v) to +(0,.8cm);
\draw [thick] (u1) to (v); \draw [thick] (um) to (v);
\draw [thick] (u1) to +(-.8cm,-.6cm); \draw [thick] (u1) to +(.8cm,-.6cm);
\draw [thick] (um) to +(-.8cm,-.6cm); \draw [thick] (um) to +(.8cm,-.6cm);
\end{tikzpicture}
\end{center}  
Suppose $b\in PB_4$ is the braid 
\begin{center}\begin{tikzpicture}[xscale=.7, yscale=.4]
\braid[number of strands=4, thick] (braid) at (0,0) a_2 a_2;
\node at (-2,-1) {$b = s_2s_2 =$};
\end{tikzpicture}\end{center}
whose underlying permutation is the identity permutation $\id_4 \in S_4$.  The elements $\beta$ and $\beta b$ are distinct in the braided tensor product $\bigl(\pi^*\T\tensorb\pi^*\Q\bigr)(4)$.

On the other hand, recall that the braid groups act on $\pi^*\bigl(\T\tensors\Q\bigr)$ via underlying permutations.  So there are equalities
\[\begin{split}\psi(\beta b) &= \psi(\beta)b \\
&= \psi(\beta)\bbar\\
& = \psi(\beta)\id_4\\
&= \psi(\beta) \in \pi^*\bigl(\T\tensors\Q\bigr)(4).\end{split}\]  
This shows that the morphism $\psi$ is not injective, and $\pi^*$ is not a strong monoidal functor.

The proof above also establishes the pure braided case because $b$ is a pure braid.
\end{proof}

\begin{proposition}\label{tensors-tensorcac}
The symmetric monoidal functors 
\[\nicexy{\bigl(\pcacopset,\tensorpcac\bigr) &\\
\bigl(\cacopset,\tensorcac\bigr) & \bigl(\sopset,\tensors\bigr) \ar[l]_-{\pi^*} \ar`u/3pt[ul][ul]_-{(\pi\iota)^*}}\] induced by the morphisms $\pi : \Cac \to \S$ and $\pi\iota : \PCac\to\S$ are \underline{not} strong monoidal.  
\end{proposition}

\begin{proof}
We slightly modify the proof of Proposition \ref{tensors-tensorb} by replacing the pure braid $b \in PB_4$ with the pure cactus $s^{(4)}_{2,3} s^{(4)}_{2,3} \in PCac_4$.
\end{proof}

Similarly, we have the following ribbon analogue.

\begin{proposition}\label{tensors-tensorr}
The symmetric monoidal functors \[\nicexy{\bigl(\propset,\tensorpr\bigr) &\\
\bigl(\ropset,\tensorr\bigr) & \bigl(\sopset,\tensors\bigr) \ar[l]_-{\pi^*} \ar`u/3pt[ul][ul]_-{\rho^*}}\] induced by the morphisms $\pi : \R \to \S$ and $\rho : \PR \to \S$ are \underline{not} strong monoidal.  
\end{proposition}

\section{Local Finite Presentability}\label{sec:lfp}

The purpose of this section is to observe that the category of $\G$-operads is locally finitely presentable.  Our categorical references in this section are \cite{adamek-rosicky,borceux1,borceux2}.  Let us first recall some relevant definitions.

\begin{definition}
Suppose $\C$ is a category.
\begin{enumerate}
\item A \index{category!generator}\index{generator}\emph{generator} of $\C$ is a set of objects $\calx = \{X_i\}_{i\in I}$ in $\C$ such that, for each pair of distinct morphisms $f,g : A \to B$ in $\C$, there exists a morphism \[h_i : X_i \to A\] for some index $i \in I$ satisfying \[fh_i \not= gh_i.\]
\item A \index{strong generator}\emph{strong generator} of $\C$ is a generator $\calx = \{X_i\}_{i\in I}$ such that, for each object $B$ and each proper sub-object $A$ of $B$, there exists a morphism \[h_i : X_i \to B\] for some index $i \in I$ that does not factor through $A$.
\item An object $X$ in $\C$ is a \emph{finitely presentable object}\index{finitely presentable object} if the representable functor \[\C(X,-) : \C \to \Set\] preserves $\aleph_0$-filtered colimits.
\item $\C$ is said to be \index{category!locally finitely presentable}\index{locally finitely presentable}\emph{locally finitely presentable} if it is cocomplete (i.e., has all small colimits) and has a set $\calx$ of finitely presentable objects such that every object in $\C$ is a directed colimit of objects in $\calx$.
\end{enumerate}
\end{definition}

\begin{definition}Suppose $(\G,\omega)$ is an action operad.  Fix a countable sequence of distinct colors $c_0,c_1,c_2,\ldots$.
\begin{enumerate}
\item For each $n \geq 0$, define $\V_n$ as the $\{c_0,\ldots,c_n\}$-colored free $\G$-operad in $\Set$ generated by one element $v_n \in \V_n\sbinom{c_0}{c_1,\ldots,c_n}$.  
\item Denote by $\calv$ the set of $\G$-operads $\{\V_n\}_{n\geq 0}$.
\end{enumerate}
\end{definition}

\begin{remark}
Over $\Set$ the only elements in the free $\G$-operad $\V_n$ are:
\begin{enumerate}[label=(\roman*)]
\item The $c_j$-colored units for $0\leq j \leq n$ and their images under the $\G(1)$-action.
\item $v_ng \in \V_n\sbinom{c_0}{c_{\gbar(1)},\ldots,c_{\gbar(n)}}$ for $g \in \G(n)$ with underlying permutation $\gbar \in S_n$.
\dqed
\end{enumerate}
\end{remark}

\begin{lemma}\label{vn-finitely-presentable}
For each $n \geq 0$, the $\G$-operad\index{G-operad@$\G$-operad!strong generator} $\V_n \in \gopset$ is characterized by the following universal property.
\begin{quote}For each $\G$-operad $\O \in \gopset$ and each element $x\in \O\dzeroud$ with $\dzeroud \in \Profdd$ and $\ud = (d_1,\ldots,d_n)$, where $\colord$ is the set of colors of $\O$, there exists a unique morphism \[f : \V_n \to \O \in \gopset\] such that 
\[\begin{cases} f(c_i) = d_i \forspace 0 \leq i \leq n,\\
f(v_n)=x.\end{cases}\]
\end{quote}
\end{lemma}

\begin{proof}
Since $\V_n$ is freely generated as a $\G$-operad by the single element $v_n$, by adjunction a morphism $f : \V_n \to \O$ of $\G$-operads is uniquely determined by the image \[f(v_n) \in \O\sbinom{fc_0}{fc_1,\ldots,fc_n}\] of the generator $v_n$.
\end{proof}

\begin{theorem}\label{goperad-locally-finitely-presentable}
The category $\gopm$ of $\G$-operads in $\M$ is locally finitely presentable.
\end{theorem}

\begin{proof}
By Theorem 1.11 in \cite{adamek-rosicky}, a category is locally finitely presentable if and only if it is cocomplete and has a strong generator consisting of finitely presentable objects.  We observed in Proposition \ref{gopm-bicomplete} that $\gopm$ is complete and cocomplete.  It remains to show that $\gopm$ has a strong generator consisting of finitely presentable objects.  To simplify notations, we will consider the case when the ambient category is $\M=\Set$.  The general case is proved similarly.  

It follows from Lemma \ref{vn-finitely-presentable} that \[\calv = \{\V_n\}_{n\geq 0}\] is a generator of $\gopset$ with each $\V_n$ a finitely presentable object.  Indeed, if $f,g : \Q \to \T$ are distinct morphisms in $\gopset$, then there exists an element $q \in \Q$ such that \[f(q) \not= g(q).\]  By Lemma \ref{vn-finitely-presentable} there exists a morphism $h_n : \V_n \to \Q$ such that \[h_n(v_n) = q,\] so \[fh_n(v_n) = f(q) \not= g(q) = gh_n(v_n).\]  Similarly, $\calv$ is a strong generator because $\Q$ being a proper sub-object of $\T$ implies that the underlying graded set of $\Q$ is a graded subset of the underlying graded set of $\T$.  So if $t \in \T \setminus \Q$, then by Lemma \ref{vn-finitely-presentable} there exists a morphism $h_m : \V_m \to \T$ such that \[h_m(v_m) = t.\]  This morphism $h_m$ does not factor through $\Q$ because $t \not\in \Q$.
\end{proof}

\begin{remark} Without going into details, Theorem \ref{goperad-locally-finitely-presentable} is equivalent to the assertion that the category $\gopm$ is categorically equivalent to the category of models of an \index{essentially finitely algebraic theory}essentially finitely algebraic theory.  We refer the reader to Chapter 3.D in \cite{adamek-rosicky} for discussion of essentially algebraic theories.\dqed
\end{remark}

\chapter{Boardman-Vogt Construction for Group Operads}\label{ch:w-group}

Fix an action operad $(\G,\omega)$ as in Definition \ref{def:augmented-group-operad} and a set $\colorc$ of colors.  Recall that $(\M,\otimes,\tensorunit)$ is an arbitrary but fixed complete and cocomplete symmetric monoidal category in which the monoidal product commutes with small colimits on each side.  It serves as our ambient category.  The Boardman-Vogt construction of a colored topological symmetric operad was introduced by Boardman and Vogt in \cite{boardman-vogt} as a particular resolution of the given symmetric operad.  This construction was extended from the ambient category of topological spaces to general symmetric monoidal categories equipped with a suitable interval for one-colored symmetric operads in \cite{berger-moerdijk-bv} and for general colored symmetric operads in \cite{berger-moerdijk-resolution}.  

In this chapter, we extend the Boardman-Vogt construction to $\G$-operads for an action operad $\G$ in a symmetric monoidal category equipped with a suitable interval.  When $\G$ is the symmetric group operad, we recover the Boardman-Vogt construction for symmetric operads.  When $\G$ is the braid group operad $\B$ or the ribbon group operad $\R$, we obtain the braided or the ribbon Boardman-Vogt construction.  The $\G$-Boardman-Vogt construction for $\G$-operads will be important when we discuss the coherent $\G$-nerve of $\G$-operads in Chapter \ref{ch:hc-nerve}.

Our construction of the $\G$-Boardman-Vogt construction is quite different from the one given by Berger and Moerdijk in \cite{berger-moerdijk-bv,berger-moerdijk-resolution} and is much closer in spirit to the original construction given by Boardman and Vogt in \cite{boardman-vogt}.  The Boardman-Vogt construction for symmetric operads given by Berger and Moerdijk is an inductively defined sequential colimit with each map a pushout that involves the previous inductive stages.  On the other hand, the original Boardman-Vogt construction in \cite{boardman-vogt} is given in one step as a quotient of a coproduct.  Our $\G$-Boardman-Vogt construction for $\G$-operads is also given in one step as a coend, which is a categorical version of a quotient.  When applied to topological symmetric operads, our construction is exactly the original Boardman-Vogt construction.

Since our $\G$-Boardman-Vogt construction is defined as a coend, our first task is to define the indexing category for this coend.  This is done in Section \ref{sec:sub-cat-group-tree}.  In Section \ref{sec:vertex-decoration} and Section \ref{sec:internal-edge-decoration}, we define the two functors, one covariant and one contravariant, from which the coend in the $\G$-Boardman-Vogt construction is taken.  The $\G$-Boardman-Vogt construction itself is defined in Section \ref{sec:boardman-vogt}.  In Section \ref{sec:w-augmentation} we observe that the $\G$-Boardman-Vogt construction is equipped with a natural augmentation over the original $\G$-operad.  In Section \ref{sec:wcommute_leftadjoint} we observe that the $\G$-Boardman-Vogt construction behaves well with respect to morphisms of action operads.

\section{Substitution Category}\label{sec:sub-cat-group-tree}

The purpose of this section is to discuss the category that serves as the indexing category of the coend that defines our $\G$-Boardman-Vogt construction of $\G$-operads.  Recall from Definition \ref{def:g-tree} that a $\colorc$-colored $\G$-tree is a pair $(T,\sigma)$ consisting of a planar $\colorc$-colored tree $T$ and an input equivariance $\sigma \in \G(|\inp(T)|)$.  We will also use the concept of $\G$-tree substitution in Definition \ref{def:g-tree-sub}.  As before, when the set $\colorc$ of colors is clear, we will omit mentioning it.

\begin{definition}\label{def:substitution-category}
Define the \index{substitution category}\index{G-tree@$\G$-tree!substitution category}\emph{substitution category of $\G$-trees}, denoted $\GTreec$, as follows.
\begin{description}
\item[Objects] The objects in $\GTreec$ are\index{G-tree@$\G$-tree} $\G$-trees.
\item[Morphisms] Each morphism in $\GTreec$ has the form \[\nicexy@C+1cm{(T,\sigma)(T_v,\sigma_v)_{v\in\Vt(T)} \ar[r]^-{(T_v,\sigma_v)_{v\in\Vt(T)}} & (T,\sigma)}\] in which the left-hand side is a\index{G-tree@$\G$-tree!substitution} $\G$-tree substitution.
\item[Identity] The identity morphism of a $\G$-tree $(T,\sigma)$ is \[\nicexy@C+1.5cm{(T,\sigma) = (T,\sigma)\bigl(\Cor_{\profofv},\id_{|\inp(v)|}\bigr)_{v} \ar[r]^-{(\Cor_{\profofv},\id_{|\inp(v)|})_{v}} & (T,\sigma)},\] in which:
\begin{itemize}
\item $v$ runs through $\Vt(T)$.
\item $\Cor_{\profofv}$ is the $\profofv$-corolla in Example \ref{ex:corolla}.
\end{itemize}
\item[Composition] The composition of the morphisms \[\nicexy@C+.6cm{\Bigl[(T,\sigma)(T_v,\sigma_v)_v\Bigr](T_v^u,\sigma_v^u)_{v,u} \ar[r]^-{(T_v^u,\sigma_v^u)_{v,u}} & (T,\sigma)(T_v,\sigma_v)_v \ar[r]^-{(T_v,\sigma_v)_v} & (T,\sigma)}\] is \[\nicexy@C+1.5cm{(T,\sigma)\Bigl[(T_v,\sigma_v)(T_v^u,\sigma_v^u)_{u}\Bigr]_v = \Bigl[(T,\sigma)(T_v,\sigma_v)_v\Bigr](T_v^u,\sigma_v^u)_{v,u} \ar[r]^-{\bigl[(T_v,\sigma_v)(T_v^u,\sigma_v^u)_u\bigr]_{v}} & (T,\sigma)}\] in which $v$ and $u$ run through $\Vt(T)$ and $\Vt(T_v)$, respectively.  The equality on the left is the associativity of $\G$-tree substitution in Lemma \ref{lem:g-tree-sub}.
\end{description}
For each $\duc \in \Profcc$, the \emph{substitution category of $\G$-trees with profiles $\duc$}, denoted by $\GTreecduc$, is the full subcategory of $\GTreec$ consisting of $\G$-trees $(T,\sigma)$ such that $\Prof(T,\sigma)=\duc$.
\end{definition}

\begin{remark}
The fact that the substitution category $\GTreec$ is actually a category (i.e., composition is unital and associative) is a consequence of the unity and the associativity of $\G$-tree substitution in Lemma \ref{lem:g-tree-sub} and the left unity in Lemma \ref{lem:g-tree-sub-vertex}.  That the full subcategory $\GTreecduc$ is well-defined is a consequence of Lemma \ref{lem:g-tree-sub}(2).\dqed
\end{remark}

\begin{example}
When $\G$ is the symmetric group operad $\S$, the substitution category $\STreec\duc$ of symmetric trees is the substitution category in Definition 3.2.11 in \cite{yau-hqft}.\dqed
\end{example}

Recall from Definition \ref{def:augmented-group-operad-morphism} the concept of a morphism $\varphi : \Gone\to\Gtwo$ of action operads.

\begin{proposition}\label{gtreec-functor}\index{change of action operads!substitution category}
Suppose $\varphi : \Gone \to \Gtwo$ is a morphism of action operads. 
\begin{enumerate}
\item Suppose $(T,\sigma) \compv (T_v,\sigma_v)$ is a $\Gone$-tree substitution at $v$ as in Definition \ref{def:treesub-g}.  Then we have \[\bigl(T,\varphi\sigma\bigr) \compv \bigl(T_v,\varphi\sigma_v\bigr) = \Bigl(T'\comp_{v'} T_v, \varphi\bigl(\sigma_v\comp^v\sigma\bigr)\Bigr).\]
\item Suppose $(T,\sigma)(T_v,\sigma_v)_v$ is a $\Gone$-tree substitution as in Definition \ref{def:g-tree-sub}.  Then the underlying planar $\colorc$-colored tree of the $\Gtwo$-tree substitution \[\bigl(T,\varphi\sigma\bigr)\bigl(T_v,\varphi\sigma_v\bigr)_{v}\] is equal to the underlying planar $\colorc$-colored tree of $(T,\sigma)(T_v,\sigma_v)_v$.  Furthermore, the input equivariance of the $\Gtwo$-tree substitution is equal to $\varphi$ applied to the input equivariance of the $\Gone$-tree substitution.
\item There is an induced functor\index{substitution category!change of action operads} \[\varphi_* : \Gonetreec \to \Gtwotreec\] that sends a $\Gone$-tree $(T,\sigma)$ to the $\Gtwo$-tree $\bigl(T,\varphi\sigma\bigr)$.
\end{enumerate}
\end{proposition}

\begin{proof}
For the first assertion, to see that the $\compv$ makes sense, observe that 
\[\begin{split}
\Prof\bigl(T_v,\varphi\sigma_v\bigr) &= \Prof(T_v)\overline{\varphi\sigma_v}\\
&= \Prof(T_v)\overline{\sigma_v}\\
&= \Prof(T_v,\sigma_v)\\
&= \profofv,
\end{split}\] 
since the underlying permutation of $\varphi\sigma_v$ agrees with that of $\sigma_v$ by Lemma \ref{rk:varphi-augmentation}.  The underlying planar $\colorc$-colored tree of $\bigl(T,\varphi\sigma\bigr) \compv\bigl(T_v,\varphi\sigma_v\bigr)$ is equal to $T'\comp_{v'} T_v$, which is the underlying planar $\colorc$-colored tree of $(T,\sigma) \compv (T_v,\sigma_v)$, because $T'$ in Definition \ref{def:treesub-g} remains unchanged if $\sigma_v$ is replaced by $\varphi\sigma_v$, since they have the same underlying permutation.  For the input equivariance, we have 
\[\begin{split}(\varphi\sigma_v) \comp^v (\varphi\sigma) &= 
\Bigl(\id_{k_v^-} \oplus (\varphi\sigma_v)\langle k_v^1,\ldots,k_v^n\rangle \oplus \id_{k_v^+}\Bigr) \cdot (\varphi\sigma)\\ &=\varphi\bigl(\sigma_v \comp^v\sigma\bigr).\end{split}\] Here we use the facts that $\varphi$ is a morphism of one-colored planar operads $\Gone \to \Gtwo$ and that level-wise it is a group homomorphism.

The second assertion follows from the first assertion because $\G$-tree substitution is by definition an iteration of $\G$-tree substitution at a single vertex.

For the third assertion, in the context of Definition \ref{def:substitution-category}, $\varphi_*$ sends the morphism \[\nicexy@C+.5cm{(T,\sigma)(T_v,\sigma_v)_{v} \ar[r]^-{(T_v,\sigma_v)_{v}} & (T,\sigma)}\]  in $\Gonetreec$ to the morphism \[\nicexy@C+.7cm{\bigl(T,\varphi\sigma\bigr)\bigl(T_v,\varphi\sigma_v\bigr)_{v} \ar[r]^-{(T_v,\varphi\sigma_v)_{v}} & \bigl(T,\varphi\sigma\bigr)}\] in $\Gtwotreec$.  The domain of the second morphism is correct by assertion (2).  The assignment $\varphi_*$ preserves the identity morphism of a typical $\Gone$-tree $(T,\sigma)$ because level-wise $\varphi$ preserves the group unit.  Similarly, by assertion (2),  $\varphi_*$ sends the composition $\bigl[(T_v,\sigma_v)(T_v^u,\sigma_v^u)_u\bigr]_{v}$ in $\Gonetreec$ as in Definition \ref{def:substitution-category} to the morphism \[\bigl[\bigl(T_v,\varphi\sigma_v\bigr)\bigl(T_v^u,\varphi\sigma_v^u\bigr)_u\bigr]_{v},\] which is the composition of $\bigl(T_v,\varphi\sigma_v\bigr)_v$ and $\bigl(T_v^u,\varphi\sigma_v^u\bigr)_{v,u}$ in $\Gtwotreec$.
\end{proof}

\begin{example}\label{ex:sub-tree-morphism}
The morphisms 
\[\nicexy{\P\ar@{=}[d] \ar[r]^-{\iota} & \PR \ar[r]^-{\iota} & \R \ar[r]^-{\pi} & \S \ar@{=}[d]\\
\P \ar[r]^-{\iota} & \PB \ar[r]^-{\iota} & \B \ar[r]^-{\pi} & \S\\
\P \ar[r]^-{\iota} \ar@{=}[u] & \PCac \ar[r]^-{\iota} & \Cac \ar[r]^-{\pi} & \S\ar@{=}[u]}\] 
of action operads in Example \ref{ex:all-augmented-group-operads} induce the functors \[\nicexy{\PTreec \ar@{=}[d] \ar[r]^-{\iota_*} & \PRTreec \ar[r]^-{\iota_*} & \RTreec\ar[r]^-{\pi_*} & \STreec\ar@{=}[d]\\
\PTreec \ar[r]^-{\iota_*} & \PBTreec \ar[r]^-{\iota_*} & \BTreec \ar[r]^-{\pi_*}  & \STreec\\
\PTreec \ar[r]^-{\iota_*} \ar@{=}[u] & \PCacTreec \ar[r]^-{\iota_*} & \CacTreec \ar[r]^-{\pi_*}  & \STreec \ar@{=}[u].}\] 
More explicitly:
\begin{itemize}
\item For a planar tree $T$, $\iota_*(T)$ is the (pure) ribbon/cactus tree $(T,\id_{|\inp(T)|})$ whose input equivariance is the identity.  
\item For a pure ribbon/braided/cactus tree $(T,\sigma)$, the ribbon/braided/cactus tree $\iota_*(T,\sigma)$ is $(T,\sigma)$.
\item For a ribbon/braided/cactus tree $(T,\sigma)$, the symmetric tree $\pi_*(T,\sigma)$ is $(T,\pi(\sigma))$, where $\pi(\sigma) \in S_{|\inp(T)|}$ is the underlying permutation of the ribbon/braid/cactus $\sigma$.  For example, the ribbon/braided/cactus corolla $(\Cor_{(\uc;d)};\sigma)$ in Examples \ref{ex:ribbon-tree}, \ref{ex:braided-tree}, and \ref{ex:cactus-tree} is sent to the permuted corolla in Example \ref{ex:symmetric-tree}.\dqed
\end{itemize}
\end{example}

\section{Vertex Decoration}\label{sec:vertex-decoration}

Our Boardman-Vogt construction of a $\G$-operad will be defined as a coend, which involves two functors, the first of which is discussed in this section.  Recall the notion of a $\colorc$-colored $\G$-operad in $\M$ in Definition \ref{def:g-operad}.  We now define vertex decoration of $\G$-trees by a $\colorc$-colored $\G$-operad.  The key observation here is that vertex decoration defines a functor from the substitution category of $\G$-trees to the ambient symmetric monoidal category $\M$.

\begin{definition}\label{def:vertex-decoration}
Suppose $\O \in \M^{\Profcc}$, and $(T,\sigma)$ is a $\G$-tree.
\begin{enumerate}
\item For a vertex $v \in T$, define \[\O(v)=\O\bigl(\profofv\bigr) \in \M\] as the entry of $\O$ corresponding to the profile of $v$.  We call this object the\index{vertex decoration}\index{decoration!vertex} \emph{$\O$-decoration of $v$}.
\item Define the unordered tensor product\label{not:ooft} \[\O(T)=\O(T,\sigma)=\bigtensorover{v\in\Vt(T)}\O(v)\in\M\] as the $S_{|\Vt(T)|}$-coinvariants \[\Bigl[\coprod_{\rho}\O(v_{\rho(1)})\otimes\cdots\otimes\O(v_{\rho(|\Vt(T)|)}) \Bigr]_{S_{|\Vt(T)|}} \in\M\] of the ordered tensor products $\bigotimes_{j=1}^{|\Vt(T)|} \O(v_{\rho(j)})$ indexed by all the vertex orderings \[\nicexy{\{1,\ldots,|\Vt(T)|\} \ar[r]^-{\rho}_-{\cong} & \Vt(T)}\] of $T$.  We call this object the\index{tree!decoration}\index{G-tree@$\G$-tree!decoration} \emph{$\O$-decoration of $T$}, or of $(T,\sigma)$.
\end{enumerate}
\end{definition}

\begin{example}
If $T$ is the $c$-colored exceptional edge $\uparrow_c$ in Example \ref{ex:exceptional-edge}, then $\O(\uparrow_c)=\tensorunit$ because $\uparrow_c$ has an empty set of vertices.\dqed
\end{example}

\begin{example}
If $T$ is the $\duc$-corolla in Example \ref{ex:corolla}, then \[\O\bigl(\Cor_{(\uc;d)},\sigma\bigr)=\O\duc.\]\dqed
\end{example}

\begin{example}\label{ex1:vertex-decoration}
For the planar trees in Example \ref{ex:treesub}, there are isomorphisms \[\begin{split}
\O(T_u) &\cong \O\sbinom{c}{a,b,f} \otimes \O\sbinom{f}{\varnothing},\quad \O(T_v) \cong \tensorunit, \quad \O(T_w) \cong \O\sbinom{e}{c,g} \otimes \O\sbinom{g}{d},\\
\O(T) &\cong \O\sbinom{e}{c,d} \otimes \O\sbinom{c}{a,b} \otimes \O\sbinom{d}{d}.\\
\O(K) &\cong \O\sbinom{e}{c,d} \otimes \O\sbinom{c}{a,b,f}\otimes \O\sbinom{f}{\varnothing}\otimes \O\sbinom{g}{d}.\end{split}\]\dqed
\end{example}

\begin{proposition}\label{goperad-functor}
Suppose $(\O,\gamma)$ is a $\colorc$-colored $\G$-operad in $\M$.
\begin{enumerate}
\item For a $\G$-tree substitution $(T,\sigma)(T_v,\sigma_v)_{v\in\Vt(T)}$, there is an isomorphism \[\O\Bigl((T,\sigma)(T_v,\sigma_v)\Bigr) \cong \bigtensorover{v\in\Vt(T)} \O(T_v,\sigma_v).\]
\item Each morphism \[(T_v,\sigma_v) : (T,\sigma)(T_v,\sigma_v) \to (T,\sigma)\] in the substitution category $\GTreec$ as in Definition \ref{def:substitution-category} yields a morphism \[\nicexy@C+.5cm{\O\Bigl((T,\sigma)(T_v,\sigma_v)\Bigr) \cong \bigtensorover{v\in\Vt(T)} \O(T_v,\sigma_v) \ar[r]^-{\bigtensorover{v} \gamma_{(T_v,\sigma_v)}} & \bigtensorover{v\in\Vt(T)} \O(v)=\O(T,\sigma)}\] such that, for each vertex $v$ in $T$, the structure morphism \[\nicexy{\O(T_v,\sigma_v) = \bigtensorover{u\in\Vt(T_v)}\O(u) \ar[r]^-{\gamma_{(T_v,\sigma_v)}} & \O(v)=\O\bigl(\Prof(T_v,\sigma_v)\bigr)},\] when restricted to the ordered tensor product $\bigotimes_{j=1}^{|\Vt(T_v)|} \O(u_{\rho_v(j)})$ for each vertex ordering $\rho_v$ of $T_v$, is the morphism $\lambda_{(T_v,\sigma_v,\rho_v)}$ in \eqref{goperad-restricted-structure}.
\item $\O$ defines a functor\index{functor!induced by a G-operad@induced by a $\G$-operad}
\begin{equation}\label{o-functor}
\O : \GTreec \to \M
\end{equation}
that sends
\begin{itemize}
\item each $\G$-tree $(T,\sigma)$ to the object $\O(T,\sigma)\in\M$ and \item each morphism $(T_v,\sigma_v)$ in $\GTreec$ to the morphism $\bigotimes_{v} \gamma_{(T_v,\sigma_v)}$.
\end{itemize}
For $\duc \in \Profcc$, the restriction of this functor to the full subcategory $\GTreecduc$ is also denoted by $\O$ and is called the \emph{vertex decoration functor}.
\end{enumerate}
\end{proposition}

\begin{proof}
First the first assertion, the isomorphism follows from the decomposition \[\Vt\bigl((T,\sigma)(T_v,\sigma_v)\bigr) = \coprodover{v\in\Vt(T)} \Vt(T_v)\] in Lemma \ref{lem:g-tree-sub}.

For the second assertion, recall from Theorem \ref{gopcm-algebra} that $\gopcm$ is the $(\Profcc$)-colored symmetric operad in $\M$ whose category of algebras is precisely the category of $\colorc$-colored $\G$-operads in $\M$.  It is the image in $\M$, under the strong symmetric monoidal functor $\Set \to \M$ that sends a set $X$ to $\coprod_{X}\tensorunit$, of the $(\Profcc)$-colored symmetric operad $\gopc$ in Definition \ref{def:gopc}.  The morphism $\gamma_{(T_v,\sigma_v)}$ is well-defined by the equivariance axiom \eqref{operad-algebra-eq} of $\O$ as a $\gopcm$-algebra and the definition of the equivariant structure in $\gopcm$.

For the last assertion, $\O$ sends an identity morphism \[\nicexy@C+1.5cm{(T,\sigma)\bigl(\Cor_{\profofv}, \id_{|\inp(v)|}\bigr) = (T,\sigma) \ar[r]^-{(\Cor_{\profofv},\id_{|\inp(v)|})} & (T,\sigma)}\] to the identity morphism of $\O(T,\sigma)$ because \[\gamma_{(\Cor_{(\uc;d)},\id_{|\uc|})} = \lambda_{(\Cor_{(\uc;d)},\id_{|\uc|},\rho)}\] is the identity morphism for each corolla $\Cor_{(\uc;d)}$, which has a unique trivial vertex ordering $\rho$, by the unity axiom \eqref{planar-operad-algebra-unity} of $\O$ as a $\gopcm$-algebra and the definition of the colored units in $\gopcm$.

Finally, a composition $\bigl[(T_v,\sigma_v)(T_v^u,\sigma_v^u)_u\bigr]_{v}$ in $\GTreec$ as in Definition \ref{def:substitution-category} is sent by $\O$ to the morphism \[\begin{split}\bigtensorover{v\in\Vt(T)} \gamma_{[(T_v,\sigma_v)(T_v^u,\sigma_v^u)_u]_{v}} &= \bigtensorover{v\in\Vt(T)} \Bigl(\gamma_{(T_v,\sigma_v)} \circ \bigtensorover{u\in\Vt(T_v)}\gamma_{(T_v^u,\sigma_v^u)}\Bigr)\\
&= \Bigl(\bigtensorover{v\in\Vt(T)} \gamma_{(T_v,\sigma_v)}\Bigr) \circ \Bigl(\bigtensorover{v\in\Vt(T)} \bigtensorover{u\in\Vt(T_v)}\gamma_{(T_v^u,\sigma_v^u)}\Bigr).\end{split}\] The first equality above holds by the associativity axiom \eqref{planar-operad-algebra-associativity} of $\O$ as a $\gopcm$-algebra and the definition of the operadic composition in $\gopcm$.  The last morphism above is the composition of the images under $\O$ of the morphisms $(T_v^u,\sigma_v^u)_{v,u}$ and $(T_v,\sigma_v)_v$.   Therefore, $\O$ preserves categorical composition.
\end{proof}

The previous proof leads to another useful description of $\G$-operads that we will use below to define the $\G$-Boardman-Vogt construction of $\G$-operads.

\begin{corollary}\label{goperad-gammat}
Each $\colorc$-colored $\G$-operad in $\M$ is equivalent to an object $\O \in \M^{\Profcc}$ equipped with a structure morphism\index{G-operad@$\G$-operad!coherence}\index{coherence!G-operad@$\G$-operad} \[\nicexy{\O(T,\sigma)\ar[r]^-{\gamma_{(T,\sigma)}} & \O\bigl(\Prof(T,\sigma)\bigr) \in \M}\] for each $\G$-tree $(T,\sigma)$ that satisfies the following two axioms.
\begin{description}
\item[Unity] $\gamma_{(\Cor_{(\uc;d)},\id_{|\uc|})}$ is the identity morphism of $\O\duc$ for each $\duc \in \Profcc$, where $\Cor_{(\uc;d)}$ is the $\duc$-corolla in Example \ref{ex:corolla}.
\item[Associativity] For each $\G$-tree substitution $(T,\sigma)(T_v,\sigma_v)$, the diagram \[\nicexy@C+.5cm{\bigtensorover{v\in\Vt(T)} \O(T_v,\sigma_v) \ar[r]^-{\bigtensorover{v} \gamma_{(T_v,\sigma_v)}} & \bigtensorover{v\in\Vt(T)} \O(v)=\O(T,\sigma) \ar[d]^-{\gamma_{(T,\sigma)}}\\ \O\bigl((T,\sigma)(T_v,\sigma_v)\bigr) \ar[u]^-{\cong} \ar[r]^-{\gamma_{(T,\sigma)(T_v,\sigma_v)}} & \O\bigl(\Prof(T,\sigma)\bigr)}\] in $\M$ is commutative.
\end{description}
Each morphism \[\varphi : (\O,\gammao) \to (\Q,\gammaq)\] of $\colorc$-colored $\G$-operads in $\M$ is equivalent to a morphism $\varphi : \O \to \Q$ in $\M^{\Profcc}$ such that the diagram \[\nicexy{\O(T,\sigma) \ar[d]_-{\gammao_{(T,\sigma)}} \ar[r]^-{\bigtensorover{v\in\Vt(T)}\varphi} & \Q(T,\sigma)\ar[d]^-{\gammaq_{(T,\sigma)}}\\ \O\bigl(\Prof(T,\sigma)\bigr) \ar[r]^-{\varphi} & \Q\bigl(\Prof(T,\sigma)\bigr)}\] is commutative for each $\G$-tree $(T,\sigma)$.
\end{corollary}

\begin{proof}
This is a consequence of Theorem \ref{gopcm-algebra} that describes $\colorc$-colored $\G$-operads in $\M$ as algebras over the $(\Profcc)$-colored symmetric operad $\gopcm$.  As in Proposition \ref{goperad-functor}(2), the structure morphism $\gamma_{(T,\sigma)}$ is defined by the requirement that, when restricted to the ordered tensor product $\bigotimes_{j=1}^{|\Vt(T)|} \O(v_{\rho(j)})$ for each vertex ordering $\rho$ of $T$, it is the structure morphism $\lambda_{(T,\sigma,\rho)}$ in \eqref{goperad-restricted-structure}.  As $\rho$ runs through all the vertex orderings of $T$, this uniquely determines $\gamma_{(T,\sigma)}$ by the equivariance axiom \eqref{operad-algebra-eq} of $\O$ as a $\gopcm$-algebra.  The above unity and associativity axioms are exactly those of $\O$ as a $\gopcm$-algebra in Definition \ref{def:planar-operad-algebra-generating}.
\end{proof}

\begin{example}
Consider the braided tree $(T,\sigma)$ in Example \ref{ex:braided-tree}.  For each braided operad $(\O,\gamma)$ in $\M$, the structure morphism $\gamma_{(T,\sigma)}$ takes the form \[\nicexy{\O(T,\sigma)\cong\O\sbinom{e}{c,d} \otimes\O\sbinom{c}{a,b} \otimes\O\sbinom{d}{d} \ar[r]^-{\gamma_{(T,\sigma)}} & \O \sbinom{e}{b,a,d}=\O\bigl(\Prof(T,\sigma)\bigr)}.\] As in Example \ref{ex:gamma-tsigma}, this morphism factors as the composite \[\nicexy{\O(T,\sigma)\cong\O\sbinom{e}{c,d} \otimes\O\sbinom{c}{a,b} \otimes\O\sbinom{d}{d} \ar[r]^-{\gamma_{(T,\sigma)}} \ar[d]_-{\gamma_{(T,\id_{|\inp(T)|})}} & \O \sbinom{e}{b,a,d}=\O\bigl(\Prof(T,\sigma)\bigr)\\
\O\bigl(\Cor_{\profoft},\sigma\bigr) \ar@{=}[r] & \O\sbinom{e}{a,b,d}\ar[u]_-{\sigma}
}\] in which the right vertical morphism is the braided equivariant structure morphism of $\O$ corresponding to the braid $\sigma \in B_3$.\dqed
\end{example}

\begin{example}\label{ex:gamma-tsigma}
Recall from Lemma \ref{lem:g-tree-sub-vertex} the corolla decomposition \[(T,\sigma) = \bigl(\Cor_{\profoft},\sigma\bigr)\compv \bigl(T,\id_{|\inp(T)|}\bigr)\] for each $\G$-tree $(T,\sigma)$.  Suppose $(\O,\gamma)$ is a $\colorc$-colored $\G$-operad in $\M$.  By associativity there is a corresponding decomposition of the structure morphism $\gamma_{(T,\sigma)}$ as the composite \[\nicexy@C+.6cm{\O(T,\sigma) \ar[r]^-{\gamma_{(T,\sigma)}} \ar@{=}[d]& \O\bigl(\Prof(T,\sigma)\bigr) = \O\bigl(\profoft\sigmabar\bigr)\\
\O(T,\id_{|\inp(T)|}) \ar[r]^-{\gamma_{(T,\id_{|\inp(T)|})}} & \O\bigl(\profoft\bigr) = \O\bigl(\Cor_{\profoft},\sigma\bigr). \ar[u]^-{\sigma\,=}_-{\gamma_{(\Cor_{\profoft},\sigma)}}}\] Here the right vertical morphism is the $\G$-equivariant structure morphism of $\O$ corresponding to the element $\sigma \in \G$.\dqed
\end{example}

\begin{example}\label{ex:phi-o-gamma}
Suppose $\varphi : \Gone \to \Gtwo$ is a morphism of action operads, and $(\O,\gammao)$ is a $\colorc$-colored $\Gone$-operad in $\M$.  Recall the left adjoint \[\varphi_! : \goneopcm \to \gtwoopcm\] in Theorem \ref{phi-goneop-gtwoop}, which is defined in \eqref{phi-of-o}.  Suppose $(T,\sigma) \in \Gtwotreec\duc$.  To describe the structure morphism $\gamma_{(T,\sigma)}$ for the $\Gtwo$-operad $\varphi_!\O$, we first compute its domain as follows.
\begin{equation}\label{phi-o-tsigma}\begin{split}
(\varphi_!\O)(T,\sigma) &= \bigtensorover{v\in\Vt(T)} (\varphi_!\O)(v)\\
&= \bigtensorover{v\in\Vt(T)}\int^{\ua^v \in\gonesubc} \coprodover{\sigma_v\in\gtwosubc(\inp(v);\ua^v)} \O(\sigmavbar v)\\
&\cong \int^{\{\ua^v\}\in\prodover{v\in\Vt(T)}\gonesubc} \bigtensorover{v\in\Vt(T)} \Bigl(\coprodover{\sigma_v\in\gtwosubc(\inp(v);\ua^v)} \O(\sigmavbar v)\Bigr)\\
&\cong \int^{\{\ua^v\}\in\prodover{v\in\Vt(T)}\gonesubc} \coprodover{\{\sigma_v\}\in\prodover{v\in\Vt(T)} \gtwosubc(\inp(v);\ua^v)} \Bigl(\bigtensorover{v\in\Vt(T)} \O(\sigmavbar v)\Bigr)
\end{split}
\end{equation}
In the above computation, for each $v\in\Vt(T)$, $\sigmavbar \in S_{|\inp(v)|}$ is the underlying permutation of $\sigma_v$, and $\sigmavbar\inp(v)=\ua^v$.  We used the abbreviation \[\O(\sigmavbar v) = \O\bigl(\sigmavbar\profofv\bigr) = \O\sbinom{\out(v)}{\sigmavbar\inp(v)}.\]

Define the $\Gtwo$-tree substitution \[(T',\tau) \defn (T,\sigma)\bigl(\Cor_{\sigmavbar v}, \sigma_v\bigr) \in \Gtwotreec\duc,\] where $\Cor_{\sigmavbar v}$ is the $\sbinom{\out(v)}{\sigmavbar\inp(v)}$-corolla in Example \ref{ex:corolla}.  Note that \[(T',\id_{|\uc|}) \in \Gonetreec\sbinom{d}{\taubar\uc},\] since $T'$ is a $\colorc$-colored planar tree such that \[\Prof(T')\taubar = \Prof(T,\sigma)=\duc.\] More explicitly, $T'$ is obtained from $T$ by changing the ordering \[\ell_v : \{1,\ldots,|\inp(v)|\} \iso \inp(v)\] at each vertex $v$ to $\ell_v\sigmavbar^{-1}$.  The structure morphism $\gamma_{(T,\sigma)}$ for $\varphi_!\O$ is determined by the commutative diagram \[\nicexy{(\varphi_!\O)(T,\sigma) \ar[r]^-{\gamma_{(T,\sigma)}} & \varphi_!\O\duc = \dint^{\ua\in\gonesubc}\coprodover{\gtwosubc(\uc;\ua)}\O\dua \\
\bigtensorover{v\in\Vt(T)} \O(\sigmavbar v) \ar[u]^-{\{\ua^v\},\,\{\sigma_v\}}_-{\text{natural}} & \coprodover{\gtwosubc(\uc;\taubar\uc)} \O\sbinom{d}{\taubar\uc} \ar[u]^-{\taubar\uc}_-{\text{natural}}\\
\O(T',\id_{|\uc|}) \ar@{=}[u] \ar[r]^-{\gammao_{(T',\id_{|\uc|})}} & \O\sbinom{d}{\taubar\uc} \ar[u]^-{\tau}_-{\text{summand}}}\] for $\{\ua^v\} \in \prod_{v} \gonesubc$ and $\{\sigma_v\} \in \prod_v\gtwosubc\bigl(\inp(v),\ua^v\bigr)$.\dqed
\end{example}

\section{Internal Edge Decoration}\label{sec:internal-edge-decoration}

In this section, we discuss the second functor that we need to define our $\G$-Boardman-Vogt construction of a $\G$-operad.  The key observation here is that, when the ambient symmetric monoidal category $(\M,\otimes,\tensorunit)$ is equipped with a suitable notion of an interval, internal edge decoration of $\G$-trees by this interval defines a contravariant functor from the substitution category of $\G$-trees to the ambient category $\M$.

The following concept is due to Berger and Moerdijk \cite{berger-moerdijk-bv}.

\begin{definition}\label{def:segment}
A \index{segment}\emph{segment} in $\M$ is a tuple \label{notation:segment} $(J, \mu, 0, 1, \epsilon)$ in which:
\begin{itemize}\item $(J, \mu, 0)$ is a monoid in $\M$.
\item $1 : \tensorunit \to J$ is an \index{absorbing element}absorbing element.
\item $\epsilon : J \to \tensorunit$ is a \index{counit}counit.  
\end{itemize}
A \index{commutative segment}\emph{commutative segment} is a segment whose multiplication $\mu$ is commutative. \end{definition}

\begin{remark}\label{rk:absorbing}
In a segment, $J$ is a monoid with unit $0 : \tensorunit \to J$ and multiplication $\mu : J \otimes J \to J$, which is commutative if the segment is commutative.  That $1 : \tensorunit \to J$ is an absorbing element means that the diagram
\[\nicexy{\tensorunit \otimes J \ar[d]_-{(\Id,\epsilon)} \ar[rr]^-{(1,\Id)} && J \otimes J \ar[d]_-{\mu} && J \otimes \tensorunit \ar[ll]_-{(\Id,1)} \ar[d]^-{(\epsilon,\Id)}\\
\tensorunit \otimes \tensorunit \ar[r]^-{\cong} & \tensorunit \ar[r]^-{1} & J & \tensorunit \ar[l]_-{1} & \tensorunit \otimes \tensorunit \ar[l]_-{\cong}}\]
is commutative.  The counity of $\epsilon$ means the diagrams
\[\nicexy{J \otimes J \ar[d]_-{\mu} \ar[r]^-{(\epsilon,\epsilon)} & \tensorunit\otimes\tensorunit \ar[d]^-{\cong}\\ J \ar[r]^-{\epsilon} & \tensorunit}\qquad
\nicexy{\tensorunit \ar[d]_-{1} \ar[r]^-{0} \ar@{=}[dr] & J \ar[d]^-{\epsilon}\\ J \ar[r]^-{\epsilon} & \tensorunit}\]
are commutative.  A commutative segment provides a concept of homotopy from the $0$-end $0 : \tensorunit \to J$ to the $1$-end $1 : \tensorunit \to J$.\dqed\end{remark}  

\begin{example}\label{ex:trivial-segment}
The \index{trivial commutative segment}\emph{trivial commutative segment} is the monoidal unit $\tensorunit$ with $0,1,\epsilon=\Id_{\tensorunit}$ and $\mu : \tensorunit\otimes\tensorunit\cong\tensorunit$ the canonical isomorphism.\dqed
\end{example}

\begin{example}\label{ex:com-segment}
Here are some examples of non-trivial commutative segments.
\begin{enumerate}
\item In the category $\CHau$ of compactly generated weak Hausdorff spaces, the \index{unit interval}unit interval $[0,1]$ equipped with the multiplication \[\mu(a,b) = \max\{a,b\}\] is a commutative segment.
\item In the category $\Sset$ of simplicial sets, the simplicial interval\index{simplicial set!segment}\label{not:deltaone} $\Delta^1=\Delta(-,[1])$ is a commutative segment with the multiplication induced by the maximum operation.
\item In the category $\Chaink$ of chain complexes over a field $\fieldk$ of characteristic $0$, the \index{normalized chain complex}\index{chain complex!segment}normalized chain complex $J=N\Delta^1$ of $\Delta^1$ is a commutative segment whose structure is uniquely determined by that on the simplicial interval $\Delta^1$ and the monoidal structure of the normalized chain functor \cite{weibel} (8.3.6 page 265).
\item In the category $\Cat$ of small categories, the \index{small category!segment}groupoid \[J = \bigl\{\nicexy@C-.5cm{0 \ar@{<->}[r]^-{\cong} & 1}\bigr\}\] with two objects $\{0,1\}$ and a unique isomorphism from $0$ to $1$ is a commutative segment with the multiplication induced by the maximum operation.\dqed
\end{enumerate}
\end{example}

For a tree $T$, recall that $\Int(T)$ denotes its set of internal edges, so $|\Int(T)|$ is the number of internal edges in $T$.  Also recall the $\colorc$-colored exceptional edge $\uparrow_c$ in Example \ref{ex:exceptional-edge}.

\begin{definition}\label{functor-J}
Suppose $(J,\mu,0,1,\epsilon)$ is a commutative segment in $\M$, and $(T,\sigma)$ is a $\G$-tree.\index{internal edge decoration}\index{decoration!internal edge}
\begin{enumerate}
\item Define the object\label{not:joft} \[\J(T)=\J(T,\sigma) = \bigotimes_{e \in\Int(T)} J = J^{\otimes |\Int(T)|}\in\M.\]
\item Suppose $v$ is a vertex in $T$, and $(T_v,\sigma_v)$ is a $\G$-tree with $\Prof(T_v,\sigma_v) = \profofv$.  Define the morphism \[\iota_v : \J(T,\sigma) \to \J\bigl((T,\sigma)\compv (T_v,\sigma_v)\bigr) \in \M\] as induced by:
\begin{itemize}
\item a copy of $0 : \tensorunit \to J$ for each internal edge in $T_v$ (if one exists);
\item the multiplication $\mu : J \otimes J \to J$ if $T_v=\,\uparrow_c$ and if $v$ is the initial vertex of an internal edge and also the terminal vertex of another internal edge;
\item the counit $\epsilon : J \to\tensorunit$ if $T_v=\,\uparrow_c$, if $T$ is not the $\ccsingle$-corolla, and if $v$ is either the initial vertex of the root of $T$ or the terminal vertex of an input of $T$;
\item the identity morphism of $\tensorunit$ if $T_v =\,\uparrow_c$ and if $T$ is the $\ccsingle$-corolla.
\end{itemize}
\end{enumerate}
\end{definition}

The following is immediate from the axioms of a commutative segment.

\begin{lemma}\label{j-is-functor}
For each $\duc \in \Profcc$, Definition \ref{functor-J} defines a functor \[\J : \GTreecducop \to \M\] that sends
\begin{itemize}
\item each $\G$-tree $(T,\sigma)$ to $\J(T,\sigma)$ and
\item each morphism \[(T_v,\sigma_v) : (T,\sigma)(T_v,\sigma_v) \to (T,\sigma)\in \GTreecduc\] to the composite of the $\iota_v$'s as $v$ runs through $\Vt(T)$.
\end{itemize}
This is called the \emph{internal edge decoration functor}\label{notation:functorj}\index{functor!induced by a segment}.
\end{lemma}

Fix a commutative segment $(J,\mu,0,1,\epsilon)$ in $\M$.  The following observation will be needed to define the $\G$-operad structure on the $\G$-Boardman-Vogt construction of a $\G$-operad.  We now use the morphism $1 : \tensorunit \to J$ of the commutative segment.

\begin{lemma}\label{lem:morphism-pi}
Suppose $(T,\sigma)(T_v,\sigma_v)$ is a $\G$-tree substitution, and $I$ is its set of internal edges that are not in any of the $T_v$'s.  Then there is a morphism \[\nicexy{\bigtensorover{v\in\Vt(T)} \J(T_v,\sigma_v) \ar[r]^-{\pi} & \J\bigl((T,\sigma)(T_v,\sigma_v)\bigr)}\] of the form $\bigl(\bigotimes_I 1\bigr)\otimes \Id_{\bigotimes_{\sqcup_{v\in\Vt(T)} \Int(T_v)}J}$ up to isomorphism.
\end{lemma}

\begin{proof}
Each internal edge in each $T_v$ becomes a unique internal edge in the $\G$-tree substitution $(T,\sigma)(T_v,\sigma_v)$.  There is a decomposition \[\Int\bigl((T,\sigma)(T_v,\sigma_v)\bigr) = I \sqcup \coprodover{v\in\Vt(T)} \Int(T_v).\]  The morphism $\pi$ is the composite
\[\nicexy{\bigtensorover{v\in\Vt(T)} \J(T_v,\sigma_v) \ar[r]^-{\pi} \ar[d]_-{\cong} & \J\bigl((T,\sigma)(T_v,\sigma_v)\bigr)\\ 
\bigl(\bigtensorover{I} \tensorunit\bigr) \otimes \bigl(\bigtensorover{\coprodover{v\in\Vt(T)} \Int(T_v)} J\bigr) \ar[r]^-{(\bigotimes_I 1,\Id)} & \bigtensorover{\Int((T,\sigma)(T_v,\sigma_v))} J \ar@{=}[u]}\]
in which $1 : \tensorunit \to J$ is part of the commutative segment.
\end{proof}

\begin{interpretation} The morphism $\pi$ in Lemma \ref{lem:morphism-pi} assigns length $1$ to each new internal edge in the $\G$-tree substitution that is not in any of the $T_v$'s.\dqed\end{interpretation}

\begin{example}\label{ex:segment-j} 
Consider the morphism \[(T_u,T_v,T_w) : K=T(T_u,T_v,T_w) \to T \in \PTreec\] in Example \ref{ex:treesub}.  Counting the number of internal edges, we have
\[\begin{split} \J(T) &\cong J_c\otimes J_d, \qquad \J(K) \cong J_c\otimes J_f\otimes J_g,\\
\J(T_u) &= J_f, \qquad \J(T_v)=\tensorunit, \andspace \J(T_w)=J_g,\end{split}\]
in which we use $J_c$ to denote a copy of the segment $J$ corresponding to a $c$-colored internal edge.  The morphism $\J(T) \to \J(K)$ is the composite in the following diagram.
\[\nicexy@C+.4cm{\J(T)\cong J_c\otimes J_d \ar[d]_-{\cong} \ar[r] & J_c\otimes J_f\otimes J_g\cong \J(K)\\ 
J_c\otimes\tensorunit\otimes\tensorunit\otimes J_d \ar[r]^-{(\Id,0,0,\epsilon)} & J_c\otimes J_f\otimes J_g \otimes\tensorunit \ar[u]_-{\cong}}\]
Each of $T_u$ and $T_w$ has one internal edge.  This accounts for the morphisms \[0 : \tensorunit \to J_f \andspace 0 : \tensorunit \to J_g.\]  The exceptional edge $T_v=\,\uparrow_d$ accounts for the counit $\epsilon : J_d \to \tensorunit$.  

Since $K=T(T_u,T_v,T_w)$, the morphism $\pi$ in Lemma \ref{lem:morphism-pi} is \[\nicexy{\J(T_u)\otimes\J(T_v)\otimes\J(T_w) \cong \tensorunit \otimes J_f\otimes J_g \ar[r]^-{(1,\Id)} & J_c\otimes J_f \otimes J_g \cong \J(K)}\] with $1 : \tensorunit \to J_c$ corresponding to the $c$-colored internal edge in $K$.\dqed
\end{example}

\section{Boardman-Vogt Construction}\label{sec:boardman-vogt}

Fix a commutative segment $(J,\mu,0,1,\epsilon)$ in $(\M,\otimes,\tensorunit)$ for the rest of this chapter.  In this section, we define the $\G$-Boardman-Vogt construction of $\G$-operads for an action operad $\G$.  Intuitively, the $\G$-Boardman-Vogt construction of a $\G$-operad $\O$ involves (i) the decoration of vertices of $\G$-trees by elements in $\O$ and (ii) the decoration of internal edges of $\G$-trees by a commutative segment.  The substitution category $\GTreecduc$ is used to parametrize the relations among these decorated $\G$-trees.  We first define the entries of the $\G$-Boardman-Vogt construction and then define its $\G$-operad structure.

\begin{definition}\label{def:wg-entries}
Suppose $\O$ is a $\colorc$-colored $\G$-operad in $\M$ for an action operad $\G$, and $\duc \in \Profcc$.  Define the coend\index{Boardman-Vogt construction!G-@$\G$-}\index{G-operad@$\G$-operad!Boardman-Vogt construction}
\begin{equation}\label{wgoduc}
\Wg\O\duc=\int^{(T,\sigma)\in\GTreecduc} \J(T,\sigma)\otimes\O(T,\sigma) \in \M,
\end{equation} in which:
\begin{itemize}
\item $\GTreecduc$ is the substitution category of $\G$-trees with profiles $\duc$ in Definition \ref{def:substitution-category}.
\item $\J: \GTreecducop \to \M$ is the internal edge decoration functor in Lemma \ref{j-is-functor}.
\item $\O : \GTreecduc \to \M$ is the vertex decoration functor in \eqref{o-functor}.
\end{itemize}
For each $\G$-tree $(T,\sigma)$, we denote by \[\nicexy{\J(T,\sigma)\otimes\O(T,\sigma) \ar[r]^-{\eta_{(T,\sigma)}} & \Wg\O\duc} \in\M\] the natural morphism.
\end{definition}

\begin{notation}\label{not:j-on-morphism}
To simplify the notation: 
\begin{enumerate}
\item The natural morphism $\eta_{(T,\sigma)}$ is sometimes abbreviated to $\eta$ when its domain is clear.
\item The image of each morphism \[\nicexy{(T,\sigma)(T_v,\sigma_v) \ar[r]^-{(T_v,\sigma_v)} & (T,\sigma)\in \GTreecduc}\] under the internal edge decoration functor $\J$ is denoted by \[\nicexy{\J(T,\sigma) \ar[r]^-{\J} & \J\bigl((T,\sigma)(T_v,\sigma_v)\bigr) \in \M}.\]  So we use the name of the functor to denote its action on a typical morphism.
\item If $f : A \to B$ is a morphism and if $\Id_Y$ is the identity morphism of an object $Y$ in the same category, then we sometimes abbreviate both morphisms \[\nicexy{A\otimes Y \ar[r]^-{f\otimes\Id_Y} & B\otimes Y} \andspace \nicexy{Y \otimes A \ar[r]^-{\Id_Y\otimes f} & Y\otimes B}\] to $f$.\dqed 
\end{enumerate}
\end{notation}

The following characterization of $\Wg\O\duc$ is simply the definition of a coend as a universal wedge, as recalled in Section \ref{subsec:categories}.

\begin{lemma}\label{wg-universal-property}
For each $\colorc$-colored $\G$-operad $(\O,\gamma)$ in $\M$, the object $\Wg\O\duc$ is characterized by the following\index{Boardman-Vogt construction!universal property} universal property: Suppose $X \in \M$ is an object, and \[f_{(T,\sigma)} : \J(T,\sigma)\otimes\O(T,\sigma) \to X\in\M\] is a morphism for each $\G$-tree $(T,\sigma) \in \GTreecduc$ such that the solid-arrow diagram
\[\nicexy{\J(T,\sigma) \otimes \O\bigl((T,\sigma)(T_v,\sigma_v)\bigr) \ar[r]^-{\bigtensorover{v} \gamma_{(T_v,\sigma_v)}} \ar[d]_-{\J} & \J(T,\sigma)\otimes\O(T,\sigma) \ar[d]_-{\eta_{(T,\sigma)}} \ar`r/3pt[ddr][ddr]_(.4){f_{(T,\sigma)}} &\\
\J\bigl((T,\sigma)(T_v,\sigma_v)\bigr)\otimes \O\bigl((T,\sigma)(T_v,\sigma_v)\bigr) \ar[r]^-{\eta_{(T,\sigma)(T_v,\sigma_v)}} \ar`d/3pt[drr][drr]^-{f_{(T,\sigma)(T_v,\sigma_v)}} & \Wg\O\duc \ar@{.>}[dr]|-{f} &\\ && X}\] is commutative for each morphism \[(T_v,\sigma_v) : (T,\sigma)(T_v,\sigma_v) \to (T,\sigma)\] in $\GTreecduc$.  Then there exists a unique morphism \[f : \Wg\O\duc \to X\in\M\] that makes the diagram commutative.
\end{lemma}

Next we define the $\colorc$-colored $\G$-operad structure on the objects $\Wg\O\duc$ as $\duc$ runs though $\Profcc$.  We will use Corollary \ref{goperad-gammat}, which characterizes $\colorc$-colored $\G$-operads in $\M$ in terms of structure morphisms $\gamma_{(T,\sigma)}$ as $(T,\sigma)$ runs through all the $\G$-trees.

\begin{notation}\label{not:gtreecv}
For a vertex $v$ in a $\G$-tree, we will abbreviate the substitution category $\GTreec(\profofv)$ of $\G$-trees with profile $\profofv$ to $\GTreec(v)$.\dqed
\end{notation}

\begin{definition}\label{def:wg-goperad-structure}
Suppose $(\O,\gamma,\operadunit)$ is a $\colorc$-colored $\G$-operad in $\M$.  For each $\G$-tree $(T,\sigma)$ with profile $\duc$, define the morphism $\gamma_{(T,\sigma)}$ via the commutative diagram 
\[\nicexy@C-.6cm{\bigtensorover{v\in\Vt(T)} \Wg\O(v) = \Wg\O(T,\sigma) \ar[r]^-{\gamma_{(T,\sigma)}} & \Wg\O\duc\\
\bigtensorover{v\in\Vt(T)} \dint^{(T_v,\sigma_v)\in\GTreec(v)} \J(T_v,\sigma_v)\otimes \O(T_v,\sigma_v) \ar@{=}[u] &\\
\dint^{\{(T_v,\sigma_v)\}\in\prodover{v} \GTreec(v)} \Bigl[\bigtensorover{v} \J(T_v,\sigma_v)\Bigr] \otimes \Bigl[\bigtensorover{v} \O(T_v,\sigma_v)\Bigr] \ar[u]^-{\cong} &\\
\Bigl[\bigtensorover{v} \J(T_v,\sigma_v)\Bigr] \otimes \Bigl[\bigtensorover{v} \O(T_v,\sigma_v)\Bigr] \ar[u]^-{\{\eta_{(T_v,\sigma_v)}\}_{v}} \ar[r]^-{(\pi,\cong)} & \J(K)\otimes \O(K) \ar[uuu]_-{\eta_K}}\]
for $\bigl\{(T_v,\sigma_v)\bigr\} \in\prod_{v} \GTreec(v)$, where \[K=(T,\sigma)(T_v,\sigma_v) \in \GTreecduc\] is the $\G$-tree substitution.  In the bottom horizontal morphism, $\pi$ is the morphism in Lemma \ref{lem:morphism-pi}, and the  isomorphism is from Proposition \ref{goperad-functor}(1).
\end{definition}

\begin{ginterpretation}
Each entry\index{Boardman-Vogt construction!geometric interpretation} $\Wg\O\duc$ is a coend constructed from the vertex decoration functor \[\O : \GTreecduc \to \M\] and the internal edge decoration functor \[\J : \GTreecducop \to \M.\]  Each node \[\J(T,\sigma)\otimes\O(T,\sigma)\] in it should be thought of as the $\G$-tree $(T,\sigma)$ in which:
\begin{itemize}\item Each internal edge is decorated by the commutative segment $J$, i.e., assigned a length set-theoretically.
\item Each vertex $v$ is decorated by the entry $\O(v)$, or set-theoretically an element in $\O(v)$.  
\end{itemize}
The coend that defines $\Wg\O\duc$ is the categorical way of first taking the coproduct sum of these decorated $\G$-trees and then taking a quotient using the $\G$-operad structure of $\O$ and the commutative segment structure of $J$.

In Definition \ref{def:wg-goperad-structure}, the structure morphism $\gamma_{(T,\sigma)}$ for $\Wg\O$ says that, given a decorated $\G$-tree $(T_v,\sigma_v)$ for each vertex $v$ in a $\G$-tree $(T,\sigma)$, we $\G$-tree substitute them into $(T,\sigma)$.  All the vertices in the $\G$-tree substitution are from the $T_v$'s, so the vertex decoration is already determined.  This corresponds to the bottom horizontal isomorphism in the big diagram in Definition \ref{def:wg-goperad-structure}.  Internal edges in the $\G$-tree substitution that are not already decorated (i.e., not in any of the $T_v$'s) are given length $1$.  This corresponds to the bottom horizontal morphism $\pi$ in the big diagram.\dqed
\end{ginterpretation}  

\begin{proposition}\label{wg-is-goperad}
In the context of Definition \ref{def:wg-goperad-structure}, $\Wg\O$ is a $\colorc$-colored $\G$-operad in $\M$.
\end{proposition}

\begin{proof}
We need to verify the unity and associativity axioms of a $\colorc$-colored $\G$-operads in Corollary \ref{goperad-gammat}.  If $(T,\sigma) = (\Cor_{(\uc;d)}, \id_{|\uc|})$, then it has a unique vertex, and \[K=\bigl(\Cor_{(\uc;d)}, \id_{|\uc|}\bigr)(T_v,\sigma_v) = (T_v,\sigma_v)\] provided $\Prof(T_v,\sigma_v)=\duc$.  In this case, the morphism $\pi$ is the identity morphism, and so is $\gamma_{(\Cor_{(\uc;d)}, \id_{|\uc|})}$.

For the associativity axiom in Corollary \ref{goperad-gammat}, suppose $(T,\sigma)(T_v,\sigma_v)$ is a $\G$-tree substitution, and $(T_v^u,\sigma_v^u)$ is a $\G$-tree with $\Prof(T_v^u,\sigma_v^u) =\profofu$ for each vertex $u$ in each $T_v$ with $v\in\Vt(T)$.  It suffices to show that the diagram \[\nicexy{\Bigl[\bigtensorover{v,u} \J(T_v^u,\sigma_v^u)\Bigr]\ar`d/3pt[dr]_(.8){\pi}[dr] \ar[r]^-{\bigtensorover{v}\pi} & \Bigl[\bigtensorover{v} \J\bigl((T_v,\sigma_v)(T_v^u,\sigma_v^u)\bigr)\Bigr] \ar[d]^-{\pi}\\
& \J\bigl((T,\sigma)(T_v,\sigma_v)(T_v^u,\sigma_v^u)\bigr)}\] is commutative.  The commutativity of this diagram is immediate from the definition of $\pi$ in Lemma \ref{lem:morphism-pi} and the fact that $1 : \tensorunit \to J$ is an absorbing element, as explained in Remark \ref{rk:absorbing}.
\end{proof}

\begin{definition}\label{def:g-bv-construction}
The $\colorc$-colored $\G$-operad $\Wgo$ in Proposition \ref{wg-is-goperad} is called the \emph{$\G$-Boardman-Vogt construction}, or the \emph{$\Wg$-construction}, of $\O$.  
\end{definition}

\begin{example}[Units and Equivariance]\label{ex:wg-units}
In the $\colorc$-colored $\G$-operad $\Wgo$, for each $c\in\colorc$, the $c$-colored unit is the composite \[\nicexy{\tensorunit \ar[d]_-{\cong} \ar[r]^-{\operadunit_c\,=\,\gamma_{(\uparrow_c,\id_1)}} & \Wgo\ccsingle\\ \tensorunit\otimes\tensorunit \ar@{=}[r] & \J(\uparrow_c,\id_1)\otimes\O(\uparrow_c,\id_1) \ar[u]_-{\eta_{(\uparrow_c,\id_1)}}}\] in which $\uparrow_c$ is the $c$-colored exceptional edge in Example \ref{ex:exceptional-edge}.  The $\colorc$-colored $\G$-sequence structure on $\Wgo$ is given as follows.  For $\duc\in\Profcc$ and $\sigma\in\G(|\uc|)$, the $\G$-sequence structure is determined by the commutative diagram \[\nicexy@C+.5cm{\Wgo\duc \ar[r]^-{\sigma\,=\,\gamma_{(\Cor_{(\uc;d)},\sigma)}} & \Wgo\ducsigmabar\\ \J(T,\tau)\otimes\O(T,\tau) \ar[u]^-{\eta_{(T,\tau)}} \ar@{=}[r] & \J(T,\tau\sigma)\otimes\O(\tau\sigma) \ar[u]_-{\eta_{(T,\tau\sigma)}}}\]
for $\G$-trees $(T,\tau)$.  In particular, these are \emph{not} induced by the colored units and the $\G$-sequence structure of $\O$.\dqed
\end{example}

Recall the action operads in Example \ref{ex:all-augmented-group-operads}.

\begin{definition}\label{def:braided-bv}
For the action operad $\G=\P$ (resp., $\S$, $\R$, $\PR$, $\B$, $\PB$, $\Cac$, or $\PCac$), the $\G$-Boardman-Vogt construction $\Wg\O$ of a $\G$-operad $\O$ is called the \index{planar operad!Boardman-Vogt construction}\emph{planar} (resp., \index{symmetric operad!Boardman-Vogt construction}\emph{symmetric}, \index{ribbon operad!Boardman-Vogt construction}\emph{ribbon}, \index{pure ribbon operad!Boardman-Vogt construction}\emph{pure ribbon}, \index{braided operad!Boardman-Vogt construction}\emph{braided}, \index{pure braided operad!Boardman-Vogt construction}\emph{pure braided}, \index{cactus operad!Boardman-Vogt construction}\emph{cactus}, or \index{pure cactus operad!Boardman-Vogt construction}\emph{pure cactus}) \emph{Boardman-Vogt construction} of $\O$.
\end{definition}

\begin{example}[Symmetric and Planar Boardman-Vogt Constructions]\label{ex:bv-symmetric}
For the symmetric group operad $\S$, our symmetric Boardman-Vogt construction $\Ws\O$ for a $\colorc$-colored symmetric operad $\O$ is isomorphic to the one given by Berger and Moerdijk \cite{berger-moerdijk-bv,berger-moerdijk-resolution}.  The difference is that ours is given in one step as a coend indexed by a combinatorially defined substitution category.  On the other hand, the Boardman-Vogt construction\index{Boardman-Vogt construction!by Berger-Moerdijk} defined by Berger and Moerdijk is an inductively defined sequential colimit with each map a pushout that depends on the previous inductive stage.  For the proof of the identification of our definition and the one by Berger and Moerdijk, we refer the reader to Section 6.6 in \cite{yau-hqft}.  Berger and Moerdijk also mentioned a planar version of their Boardman-Vogt construction in Remark 4.4(a) in \cite{berger-moerdijk-bv}.  This is isomorphic to our planar Boardman-Vogt construction.\dqed
\end{example}

\begin{example}[Topological Boardman-Vogt Construction]\label{ex:bv-top}
Over the ambient category $\CHau$ of compactly generated weak Hausdorff spaces with the commutative segment given by the unit interval $[0,1]$, we can write the topological coend\index{Boardman-Vogt construction!topological} $\Wg\O\duc$ as a quotient \[\Wg\O\duc = \left.\biggl[\coprodover{(T,\sigma)\in\GTreecduc} \Bigl(\prodover{\Int(T)} [0,1]\times \prodover{v\in\Vt(T)}\O(v)\Bigr)\biggr] \middle/\sim \right.\] of a coproduct.  So each point in the space $\Wgo\duc$ is represented by a pair \[\Bigl(\{t_e\}_{e\in\Int(T)}, \{o_v\}_{v\in\Vt(T)}\Bigr) \in \prodover{\Int(T)} [0,1]\times \prodover{v\in\Vt(T)}\O(v)\] for a $\G$-tree $(T,\sigma)\in\GTreecduc$.  We say the $\G$-tree $(T,\sigma)$ is decorated, with each internal edge $e$ decorated by the length $t_e \in [0,1]$ and each vertex $v$ decorated by the element $o_v\in\O(v)$.

Suppose $\G$ is the symmetric group operad $\S$, and $\O$ is a $\colorc$-colored symmetric operad in $\CHau$.  Then this quotient of a coproduct agrees with the original definition of the Boardman-Vogt construction in Chapter III.1 in \cite{boardman-vogt} and in Definition 2.6 in \cite{vogt}.  Furthermore, the symmetric operad structure also agrees with the one in Definition \ref{def:wg-goperad-structure}.\dqed
\end{example}

\section{Augmentation}\label{sec:w-augmentation}

The purpose of this section is to observe that the $\G$-Boardman-Vogt construction is equipped with a natural augmentation, which can be used to compare with the original $\G$-operad.  As before $(J,\mu,0,1,\epsilon)$ is a fixed commutative segment in the ambient symmetric monoidal category $(\M,\otimes,\tensorunit)$.  We will use the description of $\colorc$-colored $\G$-operads in Corollary \ref{goperad-gammat}.

\begin{theorem}\label{wg-augmentation}
Suppose $\G$ is an action operad.  Then the following statements hold.
\begin{enumerate}
\item $\Wg : \goperadcm \to \goperadcm$ is a functor.
\item For each $\colorc$-colored $\G$-operad $\O$ in $\M$, there is a morphism \[\alpha : \Wgo \to \O\] of $\colorc$-colored $\G$-operads, called the\index{Boardman-Vogt construction!augmentation}\index{augmentation!of the Boardman-Vogt construction} \emph{augmentation}, that is entrywise determined by the commutative diagram \[\nicexy@C+.8cm{\Wgo\duc \ar[r]^-{\alpha} & \O\duc\\ 
& \O(T,\sigma)\ar[u]_-{\gamma_{(T,\sigma)}}\\
\J(T,\sigma)\otimes\O(T,\sigma) \ar[uu]^-{\eta_{(T,\sigma)}} \ar[r]^-{\epsilon^{\otimes|\Int(T)|}} & \tensorunit^{\otimes|\Int(T)|}\otimes\O(T,\sigma) \ar[u]_-{\cong}}\] for $\duc \in \Profcc$ and $(T,\sigma)\in\GTreecduc$.
\item $\alpha : \Wg \to \Id$ is a natural transformation of endofunctors on $\goperadcm$.
\end{enumerate}
\end{theorem}

\begin{proof}
For the first assertion, given a morphism $\varphi : \O \to \Q$ of $\colorc$-colored $\G$-operads in $\M$, the morphism \[\Wg\varphi : \Wgo \to \Wgq\] is defined entrywise by the commutative diagram \[\nicexy{\Wgo\duc \ar[r]^-{\Wg\varphi} & \Wgq\duc\\ \J(T,\sigma)\otimes\O(T,\sigma) \ar[u]^-{\eta_{(T,\sigma)}} \ar[r]^-{\bigtensorover{v\in\Vt(T)}\varphi} & \J(T,\sigma)\otimes\Q(T,\sigma) \ar[u]_-{\eta_{(T,\sigma)}}}\] for $\duc \in \Profcc$ and $(T,\sigma)\in\GTreecduc$.  That $\Wg\varphi$ is entrywise well-defined follows from (i) the coend definitions of $\Wgo\duc$ and $\Wgq\duc$ and (ii) the compatibility of $\varphi$ with the structure morphisms $\gamma$ in $\O$ and $\Q$.  

To see that $\Wg\varphi$ is a morphism of $\colorc$-colored $\G$-operads, we need to show that, in the context of Definition \ref{def:wg-goperad-structure}, the diagram \[\nicexy@C+.8cm{\Bigl[\bigtensorover{v} \J(T_v,\sigma_v)\Bigr] \otimes \Bigl[\bigtensorover{v} \O(T_v,\sigma_v)\Bigr] \ar[d]_-{(\pi,\cong)} \ar[r]^-{\bigtensorover{v}\bigtensorover{u\in\Vt(T_v)}\varphi} & \Bigl[\bigtensorover{v} \J(T_v,\sigma_v)\Bigr] \otimes \Bigl[\bigtensorover{v} \Q(T_v,\sigma_v)\Bigr] \ar[d]^-{(\pi,\cong)}\\
\J(K)\otimes\O(K) \ar[r]^-{\bigtensorover{u\in\Vt(K)}\varphi} & \J(K)\otimes\Q(K)}\] is commutative.  This diagram is commutative by construction and  Lemma \ref{lem:g-tree-sub}(3).  The functoriality of the construction $\Wg$ is immediate from the definition of $\Wg\varphi$.

For the second assertion, to see that $\alpha$ is entrywise well-defined, we use the characterization of $\Wgo\duc$ in Lemma \ref{wg-universal-property}.  We need to show that, for each $\G$-tree substitution $K=(T,\sigma)(T_v,\sigma_v)$ in $\GTreecduc$, the diagram \[\nicexy{\J(T,\sigma) \otimes \O\bigl((T,\sigma)(T_v,\sigma_v)\bigr) \ar[r]^-{\bigtensorover{v} \gamma_{(T_v,\sigma_v)}} \ar[d]_-{\J} & \J(T,\sigma)\otimes\O(T,\sigma) \ar[d]^-{(\cong)\epsilon^{|\Int(T)|}} \\
\J\bigl((T,\sigma)(T_v,\sigma_v)\bigr)\otimes \O\bigl((T,\sigma)(T_v,\sigma_v)\bigr) \ar[d]_-{(\cong)\epsilon^{|\Int(K)|}} & \O(T,\sigma) \ar[d]^-{\gamma_{(T,\sigma)}}\\ 
\O\bigl((T,\sigma)(T_v,\sigma_v)\bigr) \ar[r]^-{\gamma_{(T,\sigma)(T_v,\sigma_v)}}
& \O\duc}\] is commutative.  The commutativity of this diagram follows from (i) the associativity of the structure morphism $\gamma$ in $\O$ in Corollary \ref{goperad-gammat} and (ii) the fact that $\epsilon$ is the counit of the commutative segment.  Similarly, the compatibility of $\alpha$ with the $\G$-operad structure morphisms $\gamma$ in $\Wgo$ and $\O$ boils down to the diagram \[\nicexy@C+1.5cm{\Bigl[\bigtensorover{v} \J(T_v,\sigma_v)\Bigr] \otimes \Bigl[\bigtensorover{v} \O(T_v,\sigma_v)\Bigr] \ar[d]_-{(\cong)(\otimes\epsilon)\pi} \ar[r]^-{(\bigtensorover{v}\gamma_{(T_v,\sigma_v)})(\cong)(\otimes\epsilon)} & \bigtensorover{v}\O(v)=\O(T,\sigma)\ar[d]^-{\gamma_{(T,\sigma)}}\\ 
\O(K)\ar[r]^-{\gamma_K} & \O\duc}\] that is commutative for the same reasons.

For the last assertion, suppose given a morphism $\varphi : \O \to \Q$ of $\colorc$-colored $\G$-operads in $\M$.  The commutativity of the diagram \[\nicexy{\Wgo\ar[d]_-{\alpha^{\O}} \ar[r]^-{\Wg\varphi} & \Wgq\ar[d]^-{\alpha^{\Q}}\\ \O\ar[r]^-{\varphi} & \Q}\] of $\colorc$-colored $\G$-operads boils down to the diagram \[\nicexy{\J(T,\sigma)\otimes\O(T,\sigma) \ar[r]^-{\bigtensorover{v} \varphi} \ar[d]_-{\gamma_{(T,\sigma)}(\cong)(\otimes\epsilon)} & \J(T,\sigma)\otimes\Q(T,\sigma) \ar[d]^-{\gamma_{(T,\sigma)}(\cong)(\otimes\epsilon)}\\ \O\duc \ar[r]^-{\varphi} & \Q\duc}\] for $\duc\in\Profcc$ and $(T,\sigma)\in\GTreecduc$.  The commutativity of this diagram follows from the compatibility of $\varphi$ with the structure morphisms $\gamma$ in $\O$ and $\Q$, as in Corollary \ref{goperad-gammat}.\end{proof}

\begin{interpretation}
The augmentation $\alpha : \Wgo \to \O$ goes from the $\G$-Boardman-Vogt construction to the given $\G$-operad.  Intuitively, the augmentation consists of (i) forgetting the length of each internal edge via the counit $\epsilon$ and (ii) composing in the $\G$-operad $\O$.  There is also an entrywise morphism that goes in the other direction, from a $\G$-operad to its $\G$-Boardman-Vogt construction.  Set-theoretically, this morphism sends an element in $\O$ to the corolla decorated by this element together with the identity input equivariance.\dqed
\end{interpretation}

\begin{proposition}\label{standard-section}
Suppose $\O$ is a $\colorc$-colored $\G$-operad $\O$ in $\M$.
\begin{enumerate}
\item There is a natural morphism\label{not:xi} \[\xi : \O \to \Wgo \in \M^{\Profcc},\] called the \index{standard section}\emph{standard section}, that is entrywise defined as the composite \[\nicexy@C-.5cm{\O\duc \ar[r]^-{\cong} & \J\bigl(\Cor_{(\uc;d)},\id_{|\uc|}\bigr) \otimes\O\bigl(\Cor_{(\uc;d)},\id_{|\uc|}\bigr) \ar[r]^-{\eta} & \Wgo\duc \ar@{<-}`u[ll] `[ll]_-{\xi} [ll]}\] for $\duc\in\Profcc$, where $\Cor_{(\uc;d)}$ is the $\duc$-corolla in Example \ref{ex:corolla}. 
\item There is an equality \[\alpha\circ\xi = \Id : \O \to \O \in\M^{\Profcc}.\]
\item When $\M$ is the category $\CHau$ of compactly generated weak Hausdorff spaces, the augmentation $\alpha : \Wgo\to\O$ is an entrywise homotopy equivalence with the standard section $\xi$ as its homotopy inverse.
\end{enumerate}
\end{proposition}

\begin{proof}
The first assertion is immediate from the definition.  The equality in the second assertion holds because (i) a corolla has no internal edges and (ii) $\gamma_{(\Cor_{(\uc;d)},\id_{|\uc|})}$ is the identity morphism on $\O\duc$ by the unity axiom of a $\colorc$-colored $\G$-operad in Corollary \ref{goperad-gammat}.  For the last assertion, we use the notation in Example \ref{ex:bv-top}.  An entrywise homotopy \[H=\{H_p\}_{0\leq p \leq 1} : \xi\circ\alpha \simeq \Id : \Wgo\duc\to\Wgo\duc\] is given by \[H_p\bigl(\{t_e\}_{e\in\Int(T)}, \{o_v\}_{v\in\Vt(T)}\bigr) = \bigl(\{\min\{p,t_e\}\}_{e\in\Int(T)}, \{o_v\}_{v\in\Vt(T)}\bigr).\]  Intuitively, the map $H_p$ replaces each internal edge length $t_e \in [0,1]$ for $e\in\Int(T)$ in a representing decorated tree by the minimum of $p$ and $t_e$.
\end{proof}

\begin{remark}
The standard section $\xi \in \M^{\Profcc}$ is \emph{not} a morphism of $\G$-operads in general.  In other words, $\xi$ is not in general compatible with the structure morphisms $\gamma$ in $\O$ and $\Wgo$.  One can see this by writing out the details of the diagram \[\nicexy{\O(T,\sigma) \ar[d]_-{\gamma_{(T,\sigma)}} \ar[r]^-{\bigtensorover{v\in\Vt(T)}\xi} & \Wgo(T,\sigma) \ar[d]^-{\gamma_{(T,\sigma)}}\\ \O\duc \ar[r]^-{\xi} & \Wgo\duc}\] for $(T,\sigma) \in \GTreecduc$ and checking that it cannot be commutative in general as soon as $T$ has a non-empty set of internal edges.\dqed 
\end{remark}

\begin{remark}
When $\G$ is the symmetric group operad $\S$, Proposition \ref{standard-section} is from \cite{boardman-vogt,vogt}.  For symmetric operads in a nice enough monoidal model category, the augmentation $\alpha : \Wgo \to \O$ is an entrywise weak equivalence.  This is proved in \cite{berger-moerdijk-bv}, although they used an inductive definition of the symmetric Boardman-Vogt construction.  Using our one-step coend definition of the symmetric Boardman-Vogt construction, this is proved in Section 6.6 in \cite{yau-hqft}.  A similar argument, with the substitution category $\GTreecduc$ replacing the symmetric substitution category $\STreec\duc$, shows that the augmentation is an entrywise weak equivalence for a $\G$-operad in a nice enough monoidal model category.  Since we will not use this result in this work, we refer the reader to Section 6.6 in \cite{yau-hqft} for details.\dqed
\end{remark}

\section{Change of Action Operads}\label{sec:wcommute_leftadjoint}

Recall from Definition \ref{def:augmented-group-operad-morphism} that a morphism of action operads is level-wise a group homomorphism and is a morphism of one-colored planar operads in $\Set$ that is compatible with the augmentation over the symmetric group operad.  In this section, we observe that the $\G$-Boardman-Vogt construction is well-behaved with respect to the change-of-category functor associated to a morphism of action operads.  As before, a commutative segment $(J,\mu,0,1,\epsilon)$ in $(\M,\otimes,\tensorunit)$ has been fixed, with respect to which the $\G$-Boardman-Vogt construction is defined.

\begin{theorem}\label{thm:w-left-adjoint}\index{change of action operads!Boardman-Vogt construction}
Suppose\index{Boardman-Vogt construction!change of action operads} $\varphi : (\Gone,\omega^1) \to (\Gtwo,\omega^2)$ is a morphism of action operads.  Then the diagram \[\nicexy{\goneopcm \ar[d]_-{\Wgone} \ar[r]^-{\varphi_!} & \gtwoopcm \ar[d]^-{\Wgtwo}\\ \goneopcm \ar[r]^-{\varphi_!} & \gtwoopcm}\] is commutative up to a natural isomorphism, in which $\varphi_!$ is the left adjoint in Theorem \ref{phi-goneop-gtwoop}.
\end{theorem}

To improve readability, we break the proof of Theorem \ref{thm:w-left-adjoint} into several steps.  First we compute each of the two composites.

\begin{lemma}\label{phi-w}
Suppose $(\O,\gamma)\in\goneopcm$ and $\duc \in \Profcc$.
\begin{enumerate}
\item There is an equality \[\varphi_!\Wgone\O\duc = \int^{\ua\in\gonesubc} \coprodover{\zeta\in\gtwosubc(\uc;\ua)} \int^{(U,\chi)\in\Gonetreec\dua} \J(U,\chi)\otimes\O(U,\chi).\]
\item There is a natural isomorphism \[\Wgtwo\varphi_!\O\duc \cong\int^{(T,\sigma)\in\Gtwotreec\duc} \int^{\{\ua^v\}\in\prodover{v\in\Vt(T)}\gonesubc} \mkern-35mu \coprodover{\{\sigma_v\}\in\prodover{v\in\Vt(T)}\gtwosubc(\inp(v);\ua^v)} \mkern-30mu \J(T',\tau)\otimes\O(T',\tau),\] where $(T',\tau)$ is the $\Gtwo$-tree substitution \[(T',\tau) = (T,\sigma)\bigl(\Cor_{\sigmavbar v}, \sigma_v\bigr) \in \Gtwotreec\duc\] as in Example \ref{ex:phi-o-gamma}. 
\end{enumerate}
\end{lemma}

\begin{proof}
For the first assertion, by the definition \eqref{phi-of-o} of $\varphi_!$, we have \[\varphi_!\Wgone\O\duc = \int^{\ua\in\gonesubc} \coprodover{\gtwosubc(\uc;\ua)} \Wgone\O\dua.\] Now we substitute in the coend definition \eqref{wgoduc} of $\Wgone\O\dua$ to obtain the desired expression.

For the second assertion, we start with the coend definition \eqref{wgoduc} of the $\Gtwo$-Boardman-Vogt construction $\Wgtwo\varphi_!\O\duc$.  We have the following natural isomorphisms, the first of which is \eqref{phi-o-tsigma}:
\[\begin{split}
&\Wgtwo\varphi_!\O\duc\\
&= \int^{(T,\sigma)\in\Gtwotreec\duc} \J(T,\sigma)\otimes (\varphi_!\O)(T,\sigma)\\
&\cong\int^{(T,\sigma)\in\Gtwotreec\duc} \J(T,\sigma)\otimes \Bigl[\int^{\{\ua^v\}\in\prodover{v\in\Vt(T)}\gonesubc} \mkern-30mu  \coprodover{\{\sigma_v\}\in \prodover{v\in\Vt(T)}\gtwosubc(\inp(v);\ua^v)} \bigl(\bigtensorover{v\in\Vt(T)} \O(\sigmavbar v)\bigr)\Bigr]\\
&\cong \int^{(T,\sigma)\in\Gtwotreec\duc} \int^{\{\ua^v\}\in\prodover{v\in\Vt(T)}\gonesubc} \J(T,\sigma)\otimes\Bigl[\coprodover{\{\sigma_v\}} \bigl(\bigtensorover{v\in\Vt(T)} \O(\sigmavbar v)\bigr)\Bigr]\\
&\cong \int^{(T,\sigma)\in\Gtwotreec\duc} \int^{\{\ua^v\}\in\prodover{v\in\Vt(T)}\gonesubc} \coprodover{\{\sigma_v\}} \J(T',\tau)\otimes\O(T',\tau).
\end{split}\]
In the last isomorphism, we use the fact that $T'$ is obtained from $T$ by changing the ordering $\ell_v$ of $\inp(v)$ for each $v \in\Vt(T)$ to $\ell_v\sigmavbar^{-1}$.  In particular, the sets of internal edges in $T$ and $T'$ are equal, so $\J(T,\sigma) =\J(T',\tau)$.
\end{proof}

Using Lemma \ref{phi-w}, we now define maps between $\varphi_!\Wgone\O$ and $\Wgtwo\varphi_!\O$ that will be shown to be mutual inverses.

\begin{lemma}\label{phiw-to-wphi}
Suppose $(\O,\gamma)\in\goneopcm$.  There is a morphism \[\nicexy{\varphi_!\Wgone\O \ar[r]^-{\theta} & \Wgtwo\varphi_!\O} \in \gtwoopcm\] determined entrywise by the commutative diagram \[\nicexy{\varphi_!\Wgone\O\duc \ar[r]^-{\theta} & \Wgtwo\varphi_!\O\duc\\ 
\J(U,\chi)\otimes\O(U,\chi) \ar[u]^-{\ua,\,\zeta}_-{\mathrm{natural}} \ar@{=}[r] & \J\bigl(U,\varphi(\chi)\zeta\bigr)\otimes \O\bigl(U,\varphi(\chi)\zeta\bigr) \ar[u]^-{\{\id_{|\inp(v)|}\}}_-{\mathrm{natural}}}\] for
\begin{itemize}
\item $\duc\in\Profcc$, $\ua\in\gonesubc$, $\zeta\in\gtwosubc(\uc;\ua)$, and
\item $(U,\chi) \in \Gonetreec\dua$.
\end{itemize} 
\end{lemma}

\begin{proof}
By the coend definitions of the Boardman-Vogt constructions $\Wgone$ and $\Wgtwo$, for each $\duc\in\Profcc$, there is a well-defined morphism $\theta'$ in $\M$ determined by the commutative diagram \[\nicexy@C-.8cm{\Wgone\O\duc = \dint^{(U,\chi)\in\Gonetreec\duc} \J(U,\chi)\otimes\O(U,\chi) \ar[r]^-{\theta'} & \varphi^*\Wgtwo\varphi_!\O\duc\\ 
\J(U,\chi)\otimes\O(U,\chi) \ar[u]^-{\eta_{(U,\chi)}} \ar@{=}[r] & \J\bigl(U,\varphi(\chi)\bigr)\otimes \O\bigl(U,\varphi(\chi)\bigr) \ar[u]^-{\{\id_{|\inp(v)|}\}}_-{\text{natural}}}\] for $(U,\chi) \in \Gonetreec\duc$.  When $\duc$ runs through $\Profcc$, the entrywise morphisms $\theta'$ define a morphism of $\colorc$-colored $\Gone$-operads because they are compatible with the $\Gone$-operad structure morphisms $\gamma$ in $\Wgone\O$ and $\varphi^*\Wgtwo\varphi_!\O$, in the sense of Corollary \ref{goperad-gammat}.  By the adjunction $\varphi_! \dashv \varphi^*$ in Theorem \ref{phi-goneop-gtwoop}, $\theta'$ correspons to a morphism \[\nicexy{\varphi_!\Wgone\O \ar[r]^-{\theta} & \Wgtwo\varphi_!\O} \in \gtwoopcm\] that is the desired morphism.
\end{proof}

\begin{lemma}\label{wphi-to-phiw}
Suppose $(\O,\gamma)\in\goneopcm$ and $\duc\in\Profcc$.  There is a morphism $\pi$ in $\M$ determined by the commutative diagram \[\nicexy@C-.5cm{\Wgtwo\varphi_!\O\duc \ar[r]^-{\pi} & \varphi_!\Wgone\O\duc\\
\J(T',\tau)\otimes\O(T',\tau) \ar[u]^-{(T,\sigma),\,\{\sigma_v\}}_-{\mathrm{natural}} \ar@{=}[d] & \coprodover{\gtwosubc(\uc;\taubar\uc)} \dint^{(U,\chi)\in\Gonetreec\sbinom{d}{\taubar\uc}} \J(U,\chi)\otimes\O(U,\chi) \ar[u]_-{\eta_{\taubar\uc}} \\
\J\bigl(T',\id_{|\inp(T)|}\bigr)\otimes \O\bigl(T',\id_{|\inp(T)|}\bigr) \ar[r]^-{\eta_{(T',\id)}} & \dint^{(U,\chi)\in\Gonetreec\sbinom{d}{\taubar\uc}} \J(U,\chi)\otimes\O(U,\chi) \ar[u]^-{\tau}_-{\mathrm{summand}}}\] for 
\begin{itemize}\item $(T,\sigma) \in \Gtwotreec\duc$ and  
\item $\{\sigma_v\}\in \prod_{v\in\Vt(T)} \gtwosubc\bigl(\inp(v),\sigmavbar\inp(v)\bigr)$.
\end{itemize}
\end{lemma}

\begin{proof}
The lower-left vertical equality comes from the fact that $T'$ is obtained from $T$ by changing the ordering of $\inp(v)$ by $\sigmavbar^{-1}$ for $v\in\Vt(T)$, so its set of internal edges is the same as that of $T$.  In the lower-right corner, observe that \[\Prof(T')\taubar = \profoft\sigmabar=\duc\] and that $\Prof(T') = \sbinom{d}{\taubar\uc}$, so the bottom horizontal natural morphism is defined.  The lower-right vertical morphism is the coproduct summand inclusion corresponding to $\tau\in\gtwosubc\bigl(\uc;\taubar\uc\bigr)$.  That $\pi$ is well-defined follows from the coend definitions of $\Wgone$ and $\Wgtwo$.
\end{proof}

\begin{proof}[Proof of Theorem \ref{thm:w-left-adjoint}]
It is enough to show that the entrywise morphism $\pi$ in Lemma \ref{wphi-to-phiw} is an entrywise inverse of the morphism $\theta$ in Lemma \ref{phiw-to-wphi}.  Indeed, the composite $\theta\pi$ is determined entrywise by the commutative diagram \[\nicexy@C-.6cm{\Wgtwo\varphi_!\O\duc \ar[r]^-{\pi} & \varphi_!\Wgone\O\duc \ar[r]^-{\theta} & \Wgtwo\varphi_!\O\duc\\
\J(T',\tau)\otimes\O(T',\tau) \ar[u]^-{(T,\sigma),\,\{\sigma_v\}}_-{\mathrm{natural}} \ar@{=}[r] & \J(T',\id)\otimes\O(T',\id)\ar@{=}[r] & \J(T',\tau)\otimes\O(T',\tau) \ar[u]^-{(T',\tau),\,\{\id_{|\inp(v)|}\}}_-{\mathrm{natural}}}\] for
\begin{itemize}
\item $(T,\sigma) \in \Gtwotreec\duc$ and  
\item $\{\sigma_v\}\in \prod_{v\in\Vt(T)} \gtwosubc\bigl(\inp(v),\sigmavbar\inp(v)\bigr)$.
\end{itemize}  
This composite coincides with the identity morphism of $\Wgtwo\varphi_!\O\duc$ by Lemma \ref{phi-w}(2).  

Similarly, the composite $\pi\theta$ is determined entrywise by the commutative diagram \[\nicexy@C-.6cm{\varphi_!\Wgone\O\duc \ar[r]^-{\theta} & \Wgtwo\varphi_!\O\duc \ar[r]^-{\pi} & \varphi_!\Wgone\O\duc\\
\J(U,\chi)\otimes\O(U,\chi) \ar[u]^-{\ua,\,\zeta}_-{\mathrm{natural}} \ar@{=}[r] & \J\bigl(U,\varphi(\chi)\zeta\bigr)\otimes \O\bigl(U,\varphi(\chi)\zeta\bigr) \ar@{=}[r] & \J(U,\id)\otimes\O(U,\id) \ar[u]^-{\inp(U),\,\varphi(\chi)\zeta}_-{\mathrm{natural}}}\] for
\begin{itemize}
\item $\ua\in\gonesubc$, $\zeta\in\gtwosubc(\uc;\ua)$, and 
\item $(U,\chi) \in \Gonetreec\dua$.
\end{itemize}  
This composite coincides with the identity morphism of $\varphi_!\Wgone\O\duc$ by Lemma \ref{phi-w}(1).  So \[\theta : \varphi_!\Wgone\O \iso \Wgtwo\varphi_!\O\] is actually an isomorphism of $\colorc$-colored $\Gtwo$-operads in $\M$.  The naturality of this isomorphism with respect to $\O$ again follows from a direct inspection.
\end{proof}

\begin{example}[Planar, Symmetric, and $\G$-Boardman-Vogt Constructions]\label{ex:psw}
For\index{planar operad!Boardman-Vogt construction}\index{symmetric operad!Boardman-Vogt construction} each action operad $(\G,\omega)$, there are morphisms \[\nicexy{\P \ar[r]^-{\iota} & \G \ar[r]^-{\omega} & \S}\] of action operads, in which $\iota$ is level-wise the unit inclusion $\P(n)=\{\id_n\} \to \G(n)$.  The induced diagram\label{not:wofp} \[\nicexy{\Poperadcm \ar[d]_{\Wp} \ar[r]^-{\iota_!} & \goperadcm \ar[d]_-{\Wg}\ar[r]^-{\omega_!} & \Soperadcm \ar[d]^-{\Ws}\\
\Poperadcm \ar[r]^-{\iota_!} & \goperadcm \ar[r]^-{\omega_!} & \Soperadcm}\]is commutative up to natural isomorphisms by Theorem \ref{thm:w-left-adjoint}.  In the one-colored case, the commutativity of the outer-most diagram \[\nicexy{\Poperadcm \ar[d]_{\Wp} \ar[r]^-{\iota_!} & \Soperadcm \ar[d]^-{\Ws}\\ \Poperadcm \ar[r]^-{\iota_!} & \Soperadcm}\] was noted in Remark 4.4(a) in \cite{berger-moerdijk-bv}.  However, the reader is reminded that our $\G$-Boardman-Vogt construction is defined in one step as a coend, whereas the planar/symmetric Boardman-Vogt construction in \cite{berger-moerdijk-bv} is an inductive sequential colimit.\dqed
\end{example}

\begin{example}[Ribbon, Braided, and Cactus Boardman-Vogt Construction]\label{ex:rbw}
Theorem\index{ribbon operad!Boardman-Vogt construction}\index{braided operad!Boardman-Vogt construction} \ref{thm:w-left-adjoint} applies to the morphisms of action operads in Example \ref{ex:all-augmented-group-operads}.  For example, the morphisms \[\nicexy{\PR \ar[r]^-{\iota} & \R},\qquad \nicexy{\PB\ar[r]^-{\iota} & \B}, \andspace \nicexy{\PCac \ar[r]^-{\iota} & \Cac}\] of action operads induce the diagrams\label{not:wofb} 
\[\nicexy{\PRoperadcm \ar[d]_{\Wpr} \ar[r]^-{\iota_!} & \Roperadcm \ar[d]^-{\Wr} \\ \PRoperadcm \ar[r]^-{\iota_!} & \Roperadcm}\qquad
\nicexy{\PBoperadcm \ar[d]_{\Wpb} \ar[r]^-{\iota_!} & \Boperadcm \ar[d]^-{\Wb} \\ \PBoperadcm \ar[r]^-{\iota_!} & \Boperadcm}\]
\[\nicexy{\Pcacoperadcm \ar[d]_{\Wpcac} \ar[r]^-{\iota_!} & \Cacoperadcm \ar[d]^-{\Wcac}\\ \Pcacoperadcm \ar[r]^-{\iota_!} & \Cacoperadcm}\] 
that are commutative up to natural isomorphisms.\dqed
\end{example}

%% file: infinity_group_operads.tex
\part{Infinity Group Operads}\label{part:infinity-group-operads}

\chapter{Category of Group Trees}\label{ch:group-tree-category}

Throughout this chapter, unless otherwise specified, the color set is the one-point set $\{*\}$.  The main purpose of this chapter is to study categorical properties of the $\G$-tree category $\Treecatg$ for an action operad $\G$, whose objects are $\G$-trees and whose morphisms and composition are given by $\G$-tree substitution.

The $\G$-Tree category $\Treecatg$ is defined in Section \ref{sec:category-g-tree}.  The planar dendroidal category $\Omega_p$ of Moerdijk-Weiss \cite{mw07,mw09} is isomorphic to the $\P$-tree category $\Treecatp$ with $\P$ the planar group operad in Example \ref{ex:trivial-group-operad}.  Moreover, the Moerdijk-Weiss symmetric dendroidal category $\Omega$ is equivalent to the $\S$-tree category $\Treecats$, where $\S$ is the symmetric group operad in Example \ref{ex:symmetric-group-operad}.  Similar to the (planar) dendroidal case, presheaves on $\Treecatg$ are models of $\infty$-$\G$-operads, as we will see in later chapters.  

In Section \ref{sec:naturality-group-tree} we observe that the $\G$-tree category construction defines a functor from the category of action operads to the category of small categories.  Furthermore, this functor has nice categorical properties.

In Section \ref{sec:mapping} we observe that there is a fully faithful functor \[\nicexy{\Treecatg \ar[r]^-{\G^-} & \gopset}\] from the $\G$-tree category $\Treecatg$ to the category of $\G$-operads in $\Set$.  This functor sends each $\G$-tree $(T,\sigma)$ to the $\G$-operad freely generated by the vertices in $T$.  This is analogous to the functor \[\Delta \to \Cat,\] from the \index{finite ordinal category}finite ordinal category to the category of small categories, that sends each object $n\in \Delta$ to the category $\{0 \to 1 \to \cdots \to n\}$.

When considering the action of the functor $\G^-$ on morphisms of $\G$-trees, there is a fundamental difference between 
\begin{enumerate}[label=(\roman*)]
\item the planar and symmetric dendroidal cases, when $\G$ is the planar group operad $\P$ or the symmetric group operad $\S$, and
\item the braided, ribbon, and cactus cases, when $\G$ is the braid group operad $\B$, the ribbon group operad $\R$, or the cactus group operad $\Cac$.
\end{enumerate}  
In the dendroidal cases, the action of the functor $\G^-$ on morphisms is determined by its action on edge sets.  This unique determination property does not hold in the braided, ribbon, and cactus cases, as we will see in Examples \ref{ex:associated-roperad}, \ref{ex:associated-boperad}, and \ref{ex:associated-cacoperad}.  This is ultimately due to the non-triviality of the pure braid groups, the pure ribbon groups, and the pure cactus groups.

In Section \ref{sec:free-gop-changeg} we show that the functor $\G^-$ is well-behaved with respect to a change of the action operad $\G$.  This observation will be important in Chapter \ref{ch:hc-nerve} when we discuss the coherent $\G$-nerve.

\section{Morphisms of Group Trees}\label{sec:category-g-tree}

Since we are now using the one-point set as our color set, a $\G$-tree as in Definition \ref{def:g-tree} is a pair $(T,\sigma)$ consisting of
\begin{enumerate}[label=(\roman*)]
\item a planar tree $T$ without any edge colorings and
\item an input equivariance $\sigma \in \G(|\inp(T)|)$.
\end{enumerate}  
In this section, we use $\G$-trees to define a category $\Treecatg$ whose presheaves provide a combinatorial model of infinity $\G$-operads.  We also observe that the assignment $\G \mapsto \Treecatg$ is natural and well-behaved.

When $\G$ is the planar group operad $\P$ or the symmetric group operad $\S$, we recover from the category $\Treecatg$ the planar/symmetric dendroidal categories of Moerdijk and Weiss \cite{mt,mw07,mw09}.  Below we will use $\G$-tree substitution at a vertex as in Definition \ref{def:treesub-g} and $\G$-tree substitution as in Definition \ref{def:g-tree-sub}.

\begin{definition}\label{def:treecatg}
Define the\index{G-tree category@$\G$-tree category}\index{category!G-tree@$\G$-tree}\index{morphism!G-trees@$\G$-trees} \emph{$\G$-tree category} $\Treecatg$ as follows.
\begin{description}
\item[Objects] The objects in $\Treecatg$ are $\G$-trees.
\item[Morphisms] Each morphism in $\Treecatg$ has the form \[\nicexy@C+3.5cm{(T,\sigma) \ar[r]^-{\phi\,=\,\bigl(((U,\sigma^U),u); (T_t,\sigma_t)_{t\in\Vt(T)}\bigr)} & (V,\sigma^V)}\] such that: 
\begin{itemize}
\item $(U,\sigma^U)$ is a $\G$-tree, and $u\in \Vt(U)$ with $\profofu=\Prof(T,\sigma)$.
\item $(T_t,\sigma_t)$ is a $\G$-tree with $\Prof(t)=\Prof(T_t,\sigma_t)$ for each $t\in\Vt(T)$.
\item $(V,\sigma^V)$ is equal to the $\G$-tree substitution
\begin{equation}\label{vut}
(V,\sigma^V) = (U,\sigma^U)\comp_u \Bigl[(T,\sigma)(T_t,\sigma_t)_{t\in\Vt(T)}\Bigr].
\end{equation}
\end{itemize}
For such a morphism\index{G-tree@$\G$-tree!morphism} $\phi$ in $\Treecatg$, we call
\begin{itemize}
\item the $(T_t,\sigma_t)$ the \emph{inside $\G$-trees},\index{inside G-tree@inside $\G$-tree}
\item $(U,\sigma^U)$ the\index{outside G-tree@outside $\G$-tree} \emph{outside $\G$-tree}, and
\item $u\in \Vt(U)$ the\index{substitution vertex} \emph{substitution vertex}.
\end{itemize}
We also say that the morphism $\phi$ \emph{corresponds to the decomposition} \eqref{vut}.
\item[Identity] The identity morphism of each $\G$-tree $(T,\sigma)$ is \[\Id_{(T,\sigma)} = \Bigl(\bigl((\Cor_{\profoft},\id_{|\inp(T)|}),u\bigr); \bigl(\Cor_{\profofv},\id_{|\inp(v)|}\bigr)_{v\in\Vt(T)}\Bigr)\] in which each $\Cor_?$ is a corolla as in Example \ref{ex:corolla} and $u$ is the unique vertex in $\Cor_{\profoft}$.
\item[Composition] Given a morphism $\phi : (T,\sigma) \to (V,\sigma^V)$ as above and another morphism \[\nicexy@C+3.5cm{(V,\sigma^V) \ar[r]^-{\varphi\,=\,\bigl(((W,\sigma^W),w); (V_v,\sigma_v)_{v\in\Vt(V)}\bigr)} & (X,\sigma^X)}\] in $\Treecatg$, their composite $\varphi\phi$ is given by\[\nicexy@C+2.5cm{(T,\sigma) \ar[r]^-{\varphi\phi\,=\,\bigl((K,u); (H_t)_{t\in\Vt(T)}\bigr)} & (X,\sigma^X)}\] in which:
\[\begin{split}
K &= (W,\sigma^W) \comp_w\Bigl[(U,\sigma^U)\bigl((V_v,\sigma_v)_{v\in \Vt(U)\setminus\{u\}}, (\Cor_{\profofu},\id_{|\inp(u)|})\bigr)\Bigr],\\
H_t &= (T_t,\sigma_t)\bigl(V_v,\sigma_v\bigr)_{v\in\Vt(T_t)}.
\end{split}\]
\end{description}
\end{definition}

\begin{interpretation}
In Definition \ref{def:treecatg} the target $(V,\sigma^V)$ of the morphism \[\phi : (T,\sigma)\to(V,\sigma^V)\] is obtained from $(T,\sigma)$ by first $\G$-tree substituting the inside $\G$-trees $(T_t,\sigma_t)$ at $t\in\Vt(T)$.  The result has the same profile as $(T,\sigma)$, namely $\profofu$, and is $\G$-tree substituted at the substitution vertex $u\in \Vt(U)$ to form $(V,\sigma^V)$.  One can visualize such a morphism $\phi \in \Treecatg$ as follows, where the input equivariances are omitted to simplify the picture.
\begin{center}\begin{tikzpicture}
\node[triangular, rounded corners] (T1) {$T_1$};
\node[below=.3cm of T1] (d) {$\vdots$};
\node[triangular, rounded corners, below=.3cm of d] (Tn) {$T_n$};
\node[plain, right=3cm of T1] (t1) {$t_1$}; \node[below=.5cm of t1] (d2) {$\vdots$};
\node[plain, below=.5cm of d2] (tn) {$t_n$};
\node[above=.5cm of t1] () {$T$};
\node[draw=gray!50,inner sep=10pt,ultra thick,rounded corners,
fit=(t1) (tn)] (Tbox) {};
\draw [gray!50, ->, dotted, line width=2pt, shorten >=0.03cm, bend left=10] (T1) to (t1);
\draw [gray!50, ->, dotted, line width=2pt, shorten >=0.03cm, bend left=10] (Tn) to (tn);
\node[plain, right=3cm of d2] (u) {$u$}; \node[above right=.5cm of u] (du1) {$\vdots$}; \node[right=.5cm of u] (du2) {$\vdots$}; 
\node[above=.5cm of du1] () {$U$};
\node[draw=gray!50,inner sep=10pt,ultra thick,rounded corners,
fit=(u) (du1) (du2)] (Ubox) {};
\draw [gray!50, ->, dotted, line width=2pt, shorten >=0.03cm, bend left=10] (Tbox) to (u);
\end{tikzpicture}\end{center}
By associativity of $\G$-tree substitution, we may also express the target as\[(V,\sigma^V)= \bigl[(U,\sigma^U) \comp_u (T,\sigma)\bigr]\Bigl\{(T_t,\sigma_t)_{t\in\Vt(T)}, (\Cor_{\Prof(u')},\id_{|\inp(u')|})_{u'\in\Vt(U)\setminus\{u\}}\Bigr\}.\]  Here we are using the decomposition\[\Vt\bigl((U,\sigma^U)\comp_u (T,\sigma)\bigr) =\Vt(T) \amalg \bigl[\Vt(U)\setminus\{u\}\bigr].\] In other words, $(V,\sigma^V)$ can also be obtained by first $\G$-tree substituting $(T,\sigma)$ at the substitution vertex $u$ and then the inside $\G$-trees $(T_t,\sigma_t)$ at $t \in \Vt(T)$.\dqed
\end{interpretation}

\begin{lemma}\label{lem:gf-well-defined}
In Definition \ref{def:treecatg} the composite \[\varphi\phi : (T,\sigma) \to (X,\sigma^X)\] is well-defined.
\end{lemma}

\begin{proof}
To simplify the notation, we suppress the input equivariance from the notation, so the $\G$-trees $(T,\sigma)$ and $(V,\sigma^V)$ are written as $T$ and $V$, etc.  We abbreviate $(\Cor_{\profofu},\id_{|\inp(u)|})$ to $\Cor_u$.  Furthermore, we abbreviate $v \in \Vt(V)$ to $v\in V$, and similarly for other $\G$-trees.  We have the following $\G$-tree  decompositions of $X$:
\[\begin{split} X &= W\comp_w \bigl[V(V_v)_{v\in V}\bigr]\\
&= W\comp_w \Bigl[\bigl[U\comp_u\bigl(T(T_t)_{t\in T}\bigr)\bigr](V_v)_{v\in V}\Bigr]\\ 
&= W\comp_w\Bigl[\bigl[U\bigl((V_v)_{u\not= v\in U},\Cor_u\bigr)\bigr]\comp_u \bigl[T\bigl(T_t(V_v)_{v\in T_t}\bigr)_{t\in T}\bigr] \Bigr]\\
&= \Bigl[W\comp_w \bigl[U\bigl((V_v)_{u\not= v\in U},\Cor_u\bigr)\bigr]\Bigr] \comp_u \bigl[T\bigl(T_t(V_v)_{v\in T_t}\bigr)_{t\in T}\bigr]\\
&= K\comp_u\bigl[T(H_t)_{t\in T}\bigr].\end{split}\]
The last equality is from the definitions of $K$ and $H_t$.  The first two equalities are from the definitions of the morphisms $\varphi : V\to X$ and $\phi : T \to V$ in $\Treecatg$.  The third and the fourth equalities hold by the associativity of $\G$-tree substitution and the decomposition \[\Vt(V) = \Bigl(\coprodover{t\in T} \Vt(T_t)\Bigr) \amalg [\Vt(U)\setminus\{u\}]\] in Lemma \ref{lem:g-tree-sub-vertex}.  The equality \[X=K\comp_u\bigl[T(H_t)_{t\in T}\bigr]\] ensures that the composite \[\varphi\phi = \bigl((K,u); (H_t)_{t\in T}\bigr) : T \to X\] is well-defined.
\end{proof}

\begin{proposition}
For each action operad $\G$, $\Treecatg$ is a category.
\end{proposition}

\begin{proof}
That the identity morphisms are actually identity for the composition is a consequence of the unity of $\G$-tree substitution in Lemma \ref{lem:g-tree-sub-vertex}.  

For the associativity of the composition, suppose \[\nicexy@C+3.5cm{(X,\sigma^X) \ar[r]^-{\psi\,=\,\bigl(((Y,\sigma^Y),y); (X_x,\sigma_x)_{x\in\Vt(X)}\bigr)} & (Z,\sigma^Z)}\] is a morphism in $\Treecatg$.   Using the abbreviations in the proof of Lemma \ref{lem:gf-well-defined}, both composites $(\psi\varphi)\phi$ and $\psi(\varphi\phi)$ are equal to the morphism \[\nicexy@C+2cm{(T,\sigma) \ar[r]^-{\bigl((A,u); (B_t)_{t\in\Vt(T)}\bigr)} & (Z,\sigma^Z)}\] in $\Treecatg$ in which
\[\begin{split}
A &= Y \comp_y\Bigl[W\bigl((X_x)_{w\not= x\in W}, U'\bigr)\Bigr],\\
U' &= U\Bigl(\bigl(V_v(X_x)_{x\in V_v}\bigr)_{u\not= v\in U}, \Cor_u\Bigr),\\
B_t &= T_t\bigl(V_v(X_x)_{x\in V_v}\bigr)_{v\in T_t}.\end{split}\]
To see that the $\G$-trees $A$ and $B_t$ make sense, one uses the decomposition \[\Vt(X) = [\Vt(W)\setminus\{w\}] \amalg \Bigl[\coprodover{u\not=v \in U} \Vt(V_v)\Bigr] \amalg \Bigl[\coprodover{t\in T,\, v\in T_t} \Vt(V_v)\Bigr],\] which follows from the second equality \[X = W\comp_w \Bigl[\bigl[U\comp_u\bigl(T(T_t)_{t\in T}\bigr)\bigr](V_v)_{v\in V}\Bigr]\]in the proof of Lemma \ref{lem:gf-well-defined}.
\end{proof}

\begin{example}[Relation to the Substitution Category]\label{ex:gtreec-treecatg}
If $\colorc$ is the one-point set, we write the substitution category $\GTreec$ in Definition \ref{def:substitution-category} as $\GTree$.  Both the one-colored substitution category $\GTree$ and the $\G$-tree category $\Treecatg$ in Definition \ref{def:treecatg} have one-colored $\G$-trees as objects.  These two categories are related as follows.\index{G-tree@$\G$-tree!relation to substitution category}

Let $\Treecatgsub$ be the non-full subcategory of the $\G$-tree category $\Treecatg$ consisting of all the objects (i.e., all $\G$-trees) and only the morphisms of the form \[\nicexy@C+5cm{(T,\sigma) \ar[r]^-{\phi\,=\,\bigl(((\Cor_{\profoft},\id_{|\inp(T)|}),u); (T_t,\sigma_t)_{t\in\Vt(T)}\bigr)} & (V,\sigma^V)}\] in which $u$ is the unique vertex in the $\profoft$-corolla $\Cor_{\profoft}$.  In other words, a morphism in $\Treecatg$ belongs to the subcategory $\Treecatgsub$ if and only if the outside $\G$-tree is a corolla with the identity input equivariance.  In this case, the target of the morphism is the $\G$-tree substitution\[(V,\sigma^V) = (T,\sigma)(T_t,\sigma_t)_{t\in\Vt(T)}.\]  The one-colored substitution category is \[\GTree=(\Treecatgsub)^{\op},\] the opposite category of $\Treecatgsub$, as one can check by a direct inspection of Definitions \ref{def:treecatg} and \ref{def:substitution-category}.\dqed
\end{example}

\begin{example}[Planar Dendroidal Category]\label{ex:treecatp}
When $\G$ is the planar group operad $\P$ in Example \ref{ex:trivial-group-operad}, there is an isomorphism of categories \[\Treecatp\cong\Omega_p,\] in which $\Omega_p$ is the \emph{category of planar rooted trees} in Definition 2.2.1 in \cite{mt}, also known as the\index{planar dendroidal category} \emph{planar dendroidal category}.  The point-set level difference between the categories $\Treecatp$ and $\Omega_p$ is that, on the one hand, our $\P$-tree category $\Treecatp$ is defined completely combinatorially in terms of planar trees and their tree substitution.  It is independent of the concept of a planar operad.  On the other hand, $\Omega_p$ is defined as the full subcategory of the category $\popset$ of planar operads in $\Set$ whose objects are colored planar operads freely generated by the vertices of planar trees, with color sets given by the sets of edges.  In particular, morphisms in $\Omega_p$ are morphisms between freely generated planar operads in $\Set$. So the definition of the category $\Omega_p$ depends on that of a planar operad.\dqed
\end{example}

\begin{example}[Symmetric Dendroidal Category]\label{ex:treecats}
When $\G$ is the symmetric group operad $\S$ in Example \ref{ex:symmetric-group-operad}, there is an equivalence  of categories \[\Treecats\simeq\Omega,\] in which $\Omega$ is the \emph{category of rooted trees} in Definition 2.3.1 in \cite{mt} and Section 3 in \cite{mw07}, also known as the\index{symmetric dendroidal category}\index{dendroidal category} \emph{symmetric dendroidal category}.  Once again, our $\S$-tree category $\Treecats$ is defined completely combinatorially in terms of symmetric trees and their tree substitution.  It does not involve any other categories.  On the other hand, $\Omega$ is defined as the full subcategory of the category $\sopset$ of symmetric operads in $\Set$ whose objects are colored symmetric operads freely generated by the vertices of trees without a planar structure, with color sets given by the sets of edges.  In particular, morphisms in $\Omega$ are morphisms between freely generated symmetric operads in $\Set$.  So the definition of the category $\Omega$ depends on that of a symmetric operad.\dqed
\end{example}

\begin{example}[Braided and Ribbon Tree Categories]\label{ex:treecatbr}
When $\G$ is the braid group operad $\B$ in Definition \ref{def:braid-group-operad} (resp., pure braid group operad $\PB$ in Definition \ref{def:pure-braid-group-operad}), the $\B$-tree category $\Treecatb$ (resp., $\PB$-tree category $\Treecatpb$) is called the\index{braided tree category}\index{pure braided tree category} \emph{(pure) braided tree category}.  Its objects are (pure) braided trees, and its morphisms are defined in terms of (pure) braided tree substitution.  Similarly, when $\G$ is the ribbon group operad $\R$ in Definition \ref{def:ribbon-group-operad} (resp., pure ribbon group operad $\PR$ in Definition \ref{def:pure-ribbon-group-operad}), the $\R$-tree category $\Treecatr$ (resp., $\PR$-tree category $\Treecatpr$) is called the\index{ribbon tree category}\index{pure ribbon tree category} \emph{(pure) ribbon tree category}.  Its objects are (pure) ribbon trees, and its morphisms are defined in terms of (pure) ribbon tree substitution.\dqed
\end{example}

\begin{example}[Cactus Tree Categories]\label{ex:treecatcac}
When $\G$ is the cactus group operad $\Cac$ in Definition \ref{def:cactus-group-operad} (resp., pure cactus group operad $\PCac$ in Definition \ref{def:pure-cactus-group-operad}), the $\Cac$-tree category $\Treecatcac$ (resp., $\PCac$-tree category $\Treecatpcac$) is called the\index{cactus tree category}\index{pure cactus tree category} \emph{(pure) cactus tree category}.  Its objects are (pure) cactus trees, and its morphisms are defined in terms of (pure) cactus tree substitution.\dqed
\end{example}

\section{Naturality of Group Trees}\label{sec:naturality-group-tree}

Recall from Definition \ref{def:augmented-group-operad-morphism} that the category of action operads is denoted by $\agop$.  We now observe that the assignment $\G \mapsto \Treecatg$, which sends an action operad $\G$ to the $\G$-tree category $\Treecatg$, is natural and has nice categorical properties.  The category of small categories is denoted by $\Cat$.

\begin{theorem}\label{treecatg-functor}\index{change of action operads!G-tree category@$\G$-tree category}
The\index{G-tree category@$\G$-tree category!naturality} $\G$-tree category construction defines a functor \[\Treecat : \agop \to \Cat,\qquad \G \mapsto \Treecatg.\]  Moreover, suppose $\varphi : \Gone\to\Gtwo$ is a morphism of action operads with induced functor $\Treecat^{\varphi} : \Treecatgone\to\Treecatgtwo$.
\begin{enumerate}
\item If $\varphi$ is level-wise injective, then $\Treecat^{\varphi}$ is faithful and injective on objects.
\item If $\varphi$ is level-wise surjective, then $\Treecat^{\varphi}$ is full and surjective on objects.
\end{enumerate}
\end{theorem}

\begin{proof}
For a morphism $\varphi : \Gone \to \Gtwo$ of action operads, the functor \[\Treecat^{\varphi} : \Treecatgone \to \Treecatgtwo\] on objects sends each $\Gone$-tree $(T,\sigma)$ to the $\Gtwo$-tree $(T,\varphi\sigma)$.  It sends a morphism $\phi : (T,\sigma) \to (V,\sigma^V)$ in $\Treecatgone$ as in Definition \ref{def:treecatg} to the morphism \[\nicexy@C+3.5cm{(T,\varphi\sigma) \ar[r]^-{\bigl(((U,\varphi\sigma^U),u); (T_t,\varphi\sigma_t)_{t\in\Vt(T)}\bigr)} & (V,\varphi\sigma^V)}\] in $\Treecatgtwo$. That this is a well-defined morphism in $\Treecatgtwo$ and that $\Treecat^{\varphi}$ defines a functor are both consequences of the assumption that $\varphi$ is a morphism of action operads.  Since the functoriality of the assignment $\Treecat^?$ is immediate from its definition, this proves the first assertion.  The other two assertions follow from a simple inspection, using the fact that, for each action operad $\G$, two $\G$-trees $(T,\sigma)$ and $(T',\sigma')$ are equal if and only if $T=T'$ as planar trees and $\sigma=\sigma'\in \G(|\inp(T)|)$.
\end{proof}

\begin{example}[Planar Dendroidal Subcategory]\label{ex:treecatpr}
For each action operad $\G$, consider the morphism \[\iota : \P \to\G\] from the planar group operad $\P$ to $\G$ given level-wise by the unit inclusion $\{\id_n\}=\P(n) \to \G(n)$.  We assume $\G \not= \P$.  The morphism $\iota : \P\to\R$ from $\P$ to the ribbon group operad $\R$ in Example \ref{ex:pop-rop} is one such example.  By Proposition \ref{treecatg-functor} there is a functor \[\nicexy{\Omega_p \cong \Treecatp \ar[r]^-{\Treecat^{\iota}} & \Treecatg}\] between their corresponding tree categories that is faithful and injective on objects because the morphism $\iota : \P \to\G$ is level-wise injective.  The functor $\Treecat^{\iota}$ is not surjective on objects because \[\Treecat^{\iota}(T,\id_{|\inp(T)|})=(T,\id_{|\inp(T)|}),\] so the image of $\Treecat^{\iota}$ does not contain any $\G$-tree $(T,\sigma)$ with $\sigma\not=\id_{|\inp(T)|}$.

Moreover, the functor $\Treecat^{\iota}$ is not full.  For example, with $\Cor_n$ denoting the $n$-corolla in Example \ref{ex:corolla}, we have the trivial morphism set \[\Treecatp\bigl((\Cor_n,\id_n), (\Cor_n,\id_n)\bigr) = \P(n) = \{\id_n\}\] in the planar tree category.  On the other hand, we have the morphism set \[\Treecatg\bigl((\Cor_n,\id_n), (\Cor_n,\id_n)\bigr) \cong \G(n)\] in the $\G$-tree category.  As long as $\G(n)$ is not the trivial group, the last morphism set is not the one-element set.  For instance,  the $n$th ribbon group $\R(n)=R_n$ in Definition \ref{def:ribbon-group} is non-trivial for each $n\geq 1$.\dqed
\end{example}

\begin{example}[Pure Ribbon, Pure Braided, and Pure Cactus Tree Categories]\label{ex:treecatpure}
Consider the morphisms \[\iota : \PR\to\R, \quad \iota : \PB\to\B, \andspace \iota : \PCac\to\Cac\] from the pure ribbon/braid/cactus group operads to the ribbon/braid/cactus group operads in Example \ref{ex:all-augmented-group-operads}.  By Proposition \ref{treecatg-functor} the functors \[\nicexy{\Treecatpr \ar[r]^-{\Treecat^{\iota}} & \Treecatr}, \qquad \nicexy{\Treecatpb \ar[r]^-{\Treecat^{\iota}} & \Treecatb}, \andspace \nicexy{\Treecatpcac \ar[r]^-{\Treecat^{\iota}} & \Treecatcac}\] between their corresponding tree categories are faithful and injective on objects because the morphisms $\iota$ are level-wise injective.  The functors $\Treecat^{\iota}$ are not surjective on objects and are not full.\dqed
\end{example}

\begin{example}[Ribbon and Symmetric Tree Categories]\label{ex:treecatrs}
Consider the morphism \[\pi : \R\to\S\] from the ribbon group operad to the symmetric group operad in Example \ref{ex:rop-sop}.  By Proposition \ref{treecatg-functor} there is a functor \[\nicexy{\Treecatr \ar[r]^-{\Treecat^{\pi}} & \Treecats}\] between their corresponding tree categories that is full and surjective on objects because the morphism $\pi : \R\to\S$ is level-wise surjective.  The functor $\Treecat^{\pi}$ is not injective on objects, since the underlying permutation morphism $R_n\to S_n$ is not injective for $n\geq 1$.  For the same reason, the functor $\Treecat^{\pi}$ is not faithful.\dqed
\end{example}

\begin{example}[Braided and Symmetric Tree Categories]\label{ex:treecatbs}
Consider the morphism \[\pi : \B\to\S\] from the braid group operad to the symmetric group operad in Example \ref{ex:bop-sop}.  By Proposition \ref{treecatg-functor} there is a functor \[\nicexy{\Treecatb \ar[r]^-{\Treecat^{\pi}} & \Treecats}\] between their corresponding tree categories that is full and surjective on objects because the morphism $\pi : \B\to\S$ is level-wise surjective.  The functor $\Treecat^{\pi}$ is not injective on objects, since the underlying permutation morphism $B_n\to S_n$ is not injective for $n\geq 2$.  For the same reason, the functor $\Treecat^{\pi}$ is not faithful.\dqed
\end{example}

\begin{example}[Cactus and Symmetric Tree Categories]\label{ex:treecatcacs}
Consider the morphism \[\pi : \Cac\to\S\] from the cactus group operad to the symmetric group operad in Example \ref{ex:cacop-sop}.  By Proposition \ref{treecatg-functor} there is a functor \[\nicexy{\Treecatcac \ar[r]^-{\Treecat^{\pi}} & \Treecats}\] between their corresponding tree categories that is full and surjective on objects because the morphism $\pi : \Cac\to\S$ is level-wise surjective.  The functor $\Treecat^{\pi}$ is not injective on objects, since the underlying permutation morphism $Cac_n\to S_n$ is not injective for $n\geq 3$.  For the same reason, the functor $\Treecat^{\pi}$ is not faithful.\dqed
\end{example}

\section{From Group Trees to Group Operads}\label{sec:mapping}

In this section, we construct a full and faithful functor from the $\G$-tree category $\Treecatg$ in Definition \ref{def:treecatg} to the category $\gopset$ of all $\G$-operads in $\Set$ in Definition \ref{def:g-operads}.  Later in this work, we will use this functor to define the $\G$-nerve and the coherent $\G$-nerve of a $\G$-operad.  

Recall the following:
\begin{enumerate}[label=(\roman*)]
\item A $\G$-tree $(T,\sigma)$ in this chapter consists of a planar tree $T$ without edge coloring and an input equivariance $\sigma \in \G(|\inp(T)|)$.
\item From Definition \ref{def:tree}, the set of edges in $T$ is denoted by $\Ed(T)$.  We will regard a planar tree $T$ as $\Ed(T)$-colored with each edge colored by itself.  In this case, for each vertex $v \in T$, its profile is \[\profofv = \sbinom{\out(v)}{\inp(v)} \in\profedtedt.\]  
\item For each set $\colorc$, the category of $\colorc$-colored $\G$-operads in $\M$ is denoted by $\goperadcm$, so $\gopedtset$ is the category of $\Ed(T)$-colored $\G$-operads in $\Set$.  
\item From Theorem \ref{gopcm-algebra}, $\gopedt$ is the $(\profedtedt)$-colored symmetric operad in $\Set$ whose category of algebras is isomorphic to the category $\gopedtset$.  In particular, there is a free-forgetful adjunction \[\nicexy@C+1cm{\Set^{\profedtedt} \ar@<2pt>[r]^-{\gopedt \circs -} & \gopedtset \ar@<2pt>[l]^-{\mathrm{forget}}},\] where $\circs$ is the symmetric circle product \eqref{symmetric-circle-product} for the color set $\Ed(T)$.
\end{enumerate}
We first define the $\G$-operad generated by the vertices in a $\G$-tree.  

\begin{definition}\label{def:gtree-goperad}
Suppose $(T,\sigma)$ is a $\G$-tree.
\begin{enumerate}
\item Define the\index{G-tree@$\G$-tree!generates a graded set} graded set\label{not:gbartsigma} \[\Gbar^{(T,\sigma)} \in \Set^{\profedtedt}\] by \[\Gbar^{(T,\sigma)}(t)=\begin{cases} \{v\} & \text{if $v\in\Vt(T)$ with $\profofv=t$},\\
\varnothing & \text{otherwise},\end{cases}\] for $t\in\profedtedt$.
\item Define the free\index{G-operad@$\G$-operad!generated by a G-tree@generated by a $\G$-tree}\index{G-tree@$\G$-tree!generates a G-operad@generates a $\G$-operad} $\Ed(T)$-colored $\G$-operad in\label{not:gtsigma} $\Set$,\[\G^{(T,\sigma)} = \gopedt \circs \Gbar^{(T,\sigma)} \in \gopedtset,\] generated by $\Gbar^{(T,\sigma)}$.
\end{enumerate}
\end{definition}

\begin{ginterpretation}\label{int:gtsigma}
Let us try to understand the free $\Ed(T)$-colored $\G$-operad $\G^{(T,\sigma)}$ better.  For
\begin{itemize}
\item each $\G$-tree $(T,\sigma)$ and 
\item each $t\in\profedtedt$,
\end{itemize}
by the definition \eqref{symmetric-circle-product} of the $\Ed(T)$-colored symmetric circle product $\circs$, we have the set \[\begin{split}
\G^{(T,\sigma)}(t) &= \bigl[\gopedt \circs \Gbar^{(T,\sigma)}\bigr](t)\\
&= \int^{\us\in\S_{\profedtedt}} \gopedt\tus\times \prod_{i=1}^n \Gbar^{(T,\sigma)}(s_i)
\end{split}\] with $\us=(s_1,\ldots,s_n)$, $n\geq 0$, and each $s_i \in\profedtedt$.  Recall from Definition \ref{def:gopc} that $\gopedt\tus$ is the set containing triples $(U,\tau,\rho)$ in which:
\begin{itemize}
\item $(U,\tau)$ is an $\Ed(T)$-colored $\G$-tree with profile $t$.
\item $\rho$ is a vertex ordering of $U$ such that the $j$th vertex has profile $s_j$ for each $1\leq j \leq |\Vt(U)|=n$.
\end{itemize}
Each set $\Gbar^{(T,\sigma)}(s_i)$ is either empty or, if $s_i=\profofu$ for some necessarily unique vertex $u$ in $T$, the set $\{u\}$.  The symmetric groups act on $\gopedt$ by reordering the vertex orderings.  Therefore, from the above coend, we conclude that $\G^{(T,\sigma)}(t)$ is the set of $\Ed(T)$-colored $\G$-trees $(U,\tau)$ with profile $t$ in which each vertex has the same profile as some vertex in $T$.  It is in this sense that $\G^{(T,\sigma)}$ is the free $\Ed(T)$-colored $\G$-operad generated by the vertices in $T$.\dqed
\end{ginterpretation}

\begin{lemma}\label{ptree-poperad}
For each planar tree $T$, each entry of the free $\Ed(T)$-colored planar operad $\P^T$ is either empty or contains a single element.
\end{lemma}

\begin{proof}
This follows from the description of the entries of $\P^T$ in Geometric Interpretation \ref{int:gtsigma}.
\end{proof}

\begin{example}[Braided Operad Generated by a Braided Tree]\label{ex:bop-btree}
With the braid group operad $\G=\B$ in Definition \ref{def:braid-group-operad}, consider the braided tree $(T,\sigma)$ 
\begin{center}\begin{tikzpicture}[scale=1.3]
\node[plain] (w) {$w$}; \node[left=1.5 of w] (T) {$T$}; \node[below=1.75 of T] (sigma) {$\sigma$};
\node[below=3 of w] (c1) {};\node[left=.7 of c1] (c2) {};\node[right=.7 of c1] (c3) {};
\node[above=.7 of c1] (i2) {};\node[above=.7 of c3] (i3) {};
\node[above=2 of c3, plain] (v) {$v$}; \node[above=2 of c2, plain] (u) {$u$};
\draw[outputleg] (w) to node[pos=.8]{\scriptsize{$f$}} +(0,.7);
\draw[thick] (u) to node{\scriptsize{$d$}} (w); 
\draw[thick] (v) to node[swap]{\scriptsize{$e$}} (w);
\draw[thick] (c1) to[out=120, in=270] node[pos=.9]{\scriptsize{$a$}} (u);
\draw[line width=6pt, white, shorten <=.5pt, shorten >=.5pt] (c2) to[out=100, in=270] (i3);
\draw[thick, shorten >=-8pt] (c2) to[out=80, in=270] (i3);
\draw[line width=4pt, white, shorten <=.5pt, shorten >=.5pt] (c3) to[out=80, in=270] (i2);
\draw[thick, shorten >=-8pt] (c3) to[out=100, in=270] (i2);
\draw[thick] (i2) to[out=90, in=270] node[swap, pos=.95]{\scriptsize{$c$}} (v);
\draw[line width=6pt, white, shorten <=.5pt, shorten >=.5pt] (i3) to[out=90, in=315] (u);
\draw[thick] (i3) to[out=90, in=315] node[swap, pos=.9]{\scriptsize{$b$}} (u);
\end{tikzpicture}\end{center}
with
\begin{itemize}
\item edges $\Ed(T)=\{a,b,c,d,e,f\}$, 
\item vertices $\Vt(T)=\{u,v,w\}$, and
\item input equivariance the braid $\sigma = s_2s_2s_1^{-1} \in B_3$.   
\end{itemize}
The graded set \[\Bbar^{(T,\sigma)} \in \Set^{\profedtedt}\] has only three non-empty components, each corresponding to a vertex in $T$:
\[\Bbar^{(T,\sigma)}\sbinom{d}{a,b}=\{u\},\quad 
\Bbar^{(T,\sigma)}\sbinom{e}{c}=\{v\},\andspace 
\Bbar^{(T,\sigma)}\sbinom{f}{d,e}=\{w\}.\] 
These three vertices freely generate the $\{a,\ldots,f\}$-colored braided operad $\B^{(T,\sigma)}$ in $\Set$.

Let us list the elements in $\B^{(T,\sigma)}$.
\begin{enumerate}
\item It contains the six colored exceptional edges \[\uparrow_x~ \in\B^{(T,\sigma)}\sbinom{x}{x}\] for $x\in\{a,\ldots,f\}$.
\item Corresponding to the vertex $u$, it contains the braided corollas \[\bigl(\Cor_{(a,b;d)};\tau\bigr)\in \B^{(T,\sigma)}\sbinom{d}{(a,b)\taubar}\] for $\tau\in B_2$.
\item Corresponding to the vertex $v$, it contains the braided corolla \[\bigl(\Cor_{(c;e)};\id_1\bigr)\in \B^{(T,\sigma)}\sbinom{e}{c}.\] 
\item Corresponding to the vertex $w$, it contains the braided corollas \[\bigl(\Cor_{(d,e;f)};\tau\bigr)\in \B^{(T,\sigma)}\sbinom{f}{(d,e)\taubar}\] for $\tau\in B_2$.
\item Corresponding to the $d$-colored internal edge, it contains the $\{a,b,d,e,f\}$-colored braided trees \[\bigl(U_d,\tau\bigr)\in \B^{(T,\sigma)}\sbinom{f}{(a,b,e)\taubar}\] for $\tau\in B_3$, where $U_d$ is the following $\{a,b,d,e,f\}$-colored planar tree.
\begin{center}\begin{tikzpicture}[scale=1.3]
\node[plain] (w) {$w$}; \node[below left=1, plain] (u) {$u$};
\draw[outputleg] (w) to node[pos=.8]{\scriptsize{$f$}} +(0,.7);
\draw[thick] (w) to node[pos=.7]{\scriptsize{$e$}} +(.6,-.6);
\draw[thick] (u) to node{\scriptsize{$d$}} (w); 
\draw[thick] (u) to node[swap, pos=.7]{\scriptsize{$a$}} +(0,-.7);
\draw[thick] (u) to node[pos=.7]{\scriptsize{$b$}} +(.6,-.6);
\end{tikzpicture}\end{center}
\item Corresponding to the $e$-colored internal edge, it contains the $\{c,d,e,f\}$-colored braided trees \[\bigl(U_e,\tau\bigr)\in \B^{(T,\sigma)}\sbinom{f}{(d,c)\taubar}\] for $\tau\in B_2$, where $U_e$ is the following $\{c,d,e,f\}$-colored planar tree.
\begin{center}\begin{tikzpicture}[scale=1.3]
\node[plain] (w) {$w$}; \node[below right=1, plain] (v) {$v$};
\draw[outputleg] (w) to node[pos=.8]{\scriptsize{$f$}} +(0,.7);
\draw[thick] (w) to node[swap, pos=.7]{\scriptsize{$d$}} +(-.6,-.6);
\draw[thick] (v) to node[swap]{\scriptsize{$e$}} (w); 
\draw[thick] (v) to node[pos=.7]{\scriptsize{$c$}} +(0,-.7);
\end{tikzpicture}\end{center}
\item Finally, corresponding to $T$ itself, it contains the $\{a,\ldots,f\}$-colored braided trees \[\bigl(T,\tau\bigr)\in \B^{(T,\sigma)}\sbinom{f}{(a,b,c)\taubar}\] for $\tau\in B_3$.
\end{enumerate}
This is a complete description of all the elements in the braided operad $\B^{(T,\sigma)}$.\dqed
\end{example}

We now define the action of $\G^-$ on morphisms of $\G$-trees.

\begin{definition}\label{def:Gbar-morphism}
Suppose given a morphism\[\nicexy@C+3.5cm{(T,\sigma) \ar[r]^-{\phi\,=\,\bigl(((U,\sigma^U),u); (T_t,\sigma_t)_{t\in\Vt(T)}\bigr)} & (V,\sigma^V)}\in\Treecatg\] as in Definition \ref{def:treecatg}.  Define a morphism \[\nicexy{\bigl(\Ed(T),\G^{(T,\sigma)}\bigr) \ar[r]^-{\G^{\phi}} & \bigl(\Ed(V),\G^{(V,\sigma^V)}\bigr)}\in \gopset\] of $\G$-operads in $\Set$ as follows.
\begin{enumerate}[label=(\roman*)]
\item Since there is a $\G$-tree decomposition \[(V,\sigma^V) = (U,\sigma^U)\comp_u \Bigl[(T,\sigma)(T_t,\sigma_t)_{t\in\Vt(T)}\Bigr],\] each edge in $T$ has a canonical image in $\Ed(V)$.  This defines the map \[\G^{\phi}: \Ed(T) \to \Ed(V)\] between color sets.
\item Since the $\Ed(T)$-colored $\G$-operad $\G^{(T,\sigma)}$ is freely generated by the graded set $\Gbar^{(T,\sigma)}$, it suffices to define the images of the vertices in $T$.  For each $t\in\Vt(T)$, we define \[\G^{\phi}(t)= (T_t',\sigma_t) \in\G^{(V,\sigma^V)}\sbinom{\G^{\phi}\out(t)}{\G^{\phi}\inp(t)},\] where $T_t'$ is $T_t$ with its edges regarded as edges in $V$ via the above decomposition of $(V,\sigma^V)$. 
\end{enumerate}
\end{definition}

\begin{theorem}\label{gtree-gop-functor}
Definitions \ref{def:gtree-goperad} and \ref{def:Gbar-morphism} define a full and faithful functor\index{G-tree category@$\G$-tree category!to G-operads@to $\G$-operads}\index{G-operad@$\G$-operad!from G-tree@from $\G$-tree} \[\nicexy{\Treecatg\ar[r]^-{\G^-} & \gopset},\qquad \begin{cases}(T,\sigma) \mapsto \G^{(T,\sigma)},& \\ \phi \mapsto \G^{\phi}.\end{cases}\]
\end{theorem}

\begin{proof}
The fact that $\G^-$ defines a functor follows from a direct inspection.  To see that $\G^-$ preserves composition of morphisms, the key is to observe that the element $\G^{\phi}(t)=(T_t',\sigma)$ is the image of an operadic composition in $\G^{(V,\sigma^V)}$, namely $\gamma_{(T_t',\sigma)}$ applied to the vertices in $T_t'$, which are all generators in the $\Ed(V)$-colored $\G$-operad $\G^{(V,\sigma^V)}$.

To see that the functor $\G^-$ is injective on hom-sets, observe that in Definition \ref{def:Gbar-morphism} the morphism $\G^{\phi} \in \gopset$ determines the inside $\G$-trees $\bigl\{(T_t,\sigma_t)\bigr\}_{t\in\Vt(T)}$ of the morphism $\phi \in \Treecatg$.  Furthermore, since $\G^{\phi}$ also determines the images of the edges of $T$ in $V$, the outside $\G$-tree $(U,\sigma^U)$ and the substitution vertex $u$ are also determined by $\G^{\phi}$.

To see that the functor $\G^-$ is surjective on hom-sets, suppose given a morphism \[\Phi : \bigl(\Ed(T),\G^{(T,\sigma)}\bigr) \to \bigl(\Ed(V),\G^{(V,\sigma^V)}\bigr)\] of $\G$-operads.  To construct a morphism \[\phi : (T,\sigma) \to (V,\sigma^V) \in \Treecatg\] such that $\G^{\phi}=\Phi$, we first define $\{\Phi(t)\}_{t\in\Vt(T)}$ as the inside $\G$-trees of $\phi$.  We now use 
\begin{itemize}
\item the map $\Phi : \Ed(T) \to \Ed(V)$ between edge sets, 
\item the $\G$-tree substitution $(T,\sigma)\bigl(\Phi(t)\bigr)_{t\in\Vt(T)}$, and 
\item the target $(V,\sigma^V)$ 
\end{itemize}
to construct the outside $\G$-tree $(U,\sigma^U)$ and the substitution vertex $u$ of $\phi$.  To obtain the input equivariance $\sigma^U$, it is crucial that each level $\G(n)$ is a group, so each element in it has a unique inverse.
\end{proof}

\begin{proposition}\label{gdash-identifies-objects}
For each planar tree $T$ and elements $\sigma,\tau\in\G(|\inp(T)|)$, there is an equality\[\G^{(T,\sigma)}=\G^{(T,\tau)}\] of $\Ed(T)$-colored $\G$-operads in $\Set$.
\end{proposition}

\begin{proof}
There is an equality of graded sets\[\Gbar^{(T,\sigma)}=\Gbar^{(T,\tau)}\] because each of these graded sets only depends on the set $\Vt(T)$ of vertices in $T$.  So the $\Ed(T)$-colored $\G$-operads freely generated by them are also equal.
\end{proof}

\begin{corollary}\label{gdash-not-embedding}
If $\G(n)$ is non-trivial for some $n\geq 0$, then the functor $\G^-$ in Theorem \ref{gtree-gop-functor} is \underline{not} injective on objects.
\end{corollary}

\begin{example}[Unique Determination on Edge Sets in Planar and Symmetric Cases]\label{ex:associated-psoperad}
If $\G$ is either the planar group operad $\P$ in Example \ref{ex:trivial-group-operad} or the symmetric group operad $\S$ in Example \ref{ex:symmetric-group-operad}, then the morphism \[\G^{\phi} : \bigl(\Ed(T),\G^{(T,\sigma)}\bigr) \to \bigl(\Ed(V),\G^{(V,\sigma^V)}\bigr)\] of $\G$-operads in Definition \ref{def:Gbar-morphism} is uniquely determined by the induced map \[\phi : \Ed(T)\to\Ed(V)\] between edge sets.  This is clear for the planar group operad.  For the symmetric group operad, this is true because for each planar tree $U$, each ordering of $\inp(U)$ is obtained from the canonical ordering $\ell_U$ as in Definition \ref{def:tree} by composing with a unique permutation.\dqed
\end{example}

\begin{example}[Ribbon Case]\label{ex:associated-roperad}
In contrast to Example \ref{ex:associated-psoperad}, if $\G$ is the ribbon group operad $\R$ in Definition \ref{def:ribbon-group-operad}, then the morphism $\R^{\phi}$ of ribbon operads is \emph{not} determined by its action on edge sets.  For instance, consider the non-identity morphism \[\nicexy@C+3cm{(\Cor_1,\id_1) \ar[r]^-{\phi\,=\,\bigl(((\Cor_1,r_1^{-1}),u); (\Cor_1,r_1)\bigr)} & (\Cor_1,\id_1)}\in \Treecatr\] with
\begin{itemize}\item $\Cor_1$ the $1$-corolla in Example \ref{ex:corolla}, 
\item $u$ its unique vertex, 
\item $\id_1\in R_1$ the identity ribbon on one strip, and 
\item $r_1\in PR_1$ the pure ribbon 
\begin{center}\begin{tikzpicture}[xscale=1.5, yscale=.7, thick]
\draw (4,0) to[out=90,in=270] (4.2,.5);
\draw[shorten >=1pt, shorten <=1pt, line width=5pt, white] (4.2,0) to[out=90,in=270] (4,.5);
\draw (4.2,0) to[out=90,in=270] (4,.5);
\draw (4,.5) to[out=90,in=270] (4.2,1);
\draw[shorten >=1pt, shorten <=1pt, line width=5pt, white] (4.2,.5) to[out=90,in=270] (4,1);
\draw (4.2,.5) to[out=90,in=270] (4,1);
\foreach \x in {0,1} \draw (4,\x) -- (4.2,\x);
\end{tikzpicture}\end{center}
on one strip with one full $2\pi$ twist in Example \ref{ex:generator-rn}. 
\end{itemize} 
Then $\R^{\phi}$ is the identity map on the edge set of $\Cor_1$.  So the action on edge sets does not determine the morphism $\R^{\phi}$ of ribbon operads in general.  A more general example replaces $(\Cor_1,\id_1)$ with $(\Cor_n,\id_n)$ and $r_1\in PR_1$ with any non-identity element that belongs to the $n$th pure ribbon group $PR_n$.\dqed
\end{example}

\begin{example}[Braided Case]\label{ex:associated-boperad}
For the braid group operad $\B$ in Definition \ref{def:braid-group-operad}, the morphism $\B^{\phi}$ of braided operads is \emph{not} determined by its action on edge sets.  For instance, consider the non-identity morphism \[\nicexy@C+4cm{(\Cor_2,\id_2) \ar[r]^-{\phi\,=\,\bigl(((\Cor_2,s_1^{-1}s_1^{-1}),u); (\Cor_2,s_1s_1)\bigr)} & (\Cor_2,\id_2)}\in \Treecatb\] with
\begin{itemize}\item $\Cor_2$ the $2$-corolla, 
\item $u$ its unique vertex, and 
\item $s_1s_1\in PB_2$ the pure braid 
\begin{center}\begin{tikzpicture}[xscale=1.5, yscale=.7, thick]
\draw (4,0) to[out=90,in=270] (4.2,.5);
\draw[shorten >=1pt, shorten <=1pt, line width=5pt, white] (4.2,0) to[out=90,in=270] (4,.5);
\draw (4.2,0) to[out=90,in=270] (4,.5);
\draw (4,.5) to[out=90,in=270] (4.2,1);
\draw[shorten >=1pt, shorten <=1pt, line width=5pt, white] (4.2,.5) to[out=90,in=270] (4,1);
\draw (4.2,.5) to[out=90,in=270] (4,1);
\end{tikzpicture}\end{center}
on two strings.
\end{itemize}
Then $\B^{\phi}$ is the identity map on the edge set of $\Cor_2$, since the underlying permutation of $s_1s_1$ is the identity permutation in $S_2$.   So the action on edge sets does not determine the morphism $\B^{\phi}$ of braided operads in general.  A more general example replaces $(\Cor_2,\id_2)$ with $(\Cor_n,\id_n)$, $n\geq 2$, and $s_1s_1\in PB_2$ with any non-identity element that belongs to the $n$th pure braid group $PB_n$.\dqed
\end{example}

\begin{example}[Cactus Case]\label{ex:associated-cacoperad}
For the cactus group operad $\Cac$ in Definition \ref{def:cactus-group-operad}, the morphism $\Cac^{\phi}$ of cactus operads is \emph{not} determined by its action on edge sets.  For instance, consider the non-identity morphism \[\nicexy@C+3.2cm{(\Cor_3,\id_3) \ar[r]^-{\phi\,=\,\bigl(((\Cor_3,q^{-1}),u); (\Cor_3,q)\bigr)} & (\Cor_3,\id_3)}\in \Treecatcac\] with
\begin{itemize}\item $\Cor_3$ the $3$-corolla, 
\item $u$ its unique vertex, and 
\item $\id_3\not= q\in PCac_3$ any pure cactus (e.g., the pure cacti in Example \ref{ex:cac-three}).
\end{itemize}
Then $\Cac^{\phi}$ is the identity map on the edge set of $\Cor_3$, since the underlying permutation of $q$ is the identity permutation in $S_3$.   So the action on edge sets does not determine the morphism $\Cac^{\phi}$ of cactus operads in general.  A more general example replaces $(\Cor_3,\id_3)$ with $(\Cor_n,\id_n)$, $n\geq 3$, and $q$ with any non-identity element that belongs to the $n$th pure cactus group $PCac_n$.\dqed
\end{example}

\section{Change of Action Operads}\label{sec:free-gop-changeg}

In this section, we observe that the functor \[\G^- : \Treecatg \to \gopset\] in Theorem \ref{gtree-gop-functor} is well-behaved with respect to changing the action operad $\G$.

\begin{theorem}\label{gdash-left-adjoint}\index{change of action operads!G-trees to G-operads@$\G$-trees to $\G$-operads}
Suppose\index{G-operad@$\G$-operad!change of action operads}\index{G-tree category@$\G$-tree category!change of action operads} $\varphi : \Gone \to \Gtwo$ is a morphism of action operads.  Then the diagram \[\nicexy{\Treecatgone \ar[d]_-{\Gonedash} \ar[r]^-{\Treecat^{\varphi}} & \Treecatgtwo \ar[d]^-{\Gtwodash}\\
\goneopset \ar[r]^-{\varphi_!} & \gtwoopset}\]is commutative up to a natural isomorphism, in which:
\begin{itemize}\item $\varphi_!$ is the left adjoint in Theorem \ref{goneopm-gtwoopm}.
\item $\Treecat^{\varphi}$ is the functor in Proposition \ref{treecatg-functor}.
\item $\Gonedash$ and $\Gtwodash$ are the functors in Theorem \ref{gtree-gop-functor}.
\end{itemize}
\end{theorem}

\begin{proof}
Suppose $(T,\sigma)$ is a $\Gone$-tree.  First note that there is an equality \[(\Gonebar)^{(T,\sigma)} = (\Gtwobar)^{(T,\varphi\sigma)} \in \Set^{\profedtedt}\] because each side only depends on the vertices in the planar tree $T$ and not on the input equivariances $\sigma \in\Gone(|\inp(T)|)$ and $\varphi\sigma\in\Gtwo(|\inp(T)|)$.  There is a factorization
\[\nicexy@C-1cm{\Set^{\profedtedt} \ar@{<-}`d/3pt[dr]_-{\text{forget}} [dr] && \gtwoop^{\Ed(T)}(\Set) \ar[ll]_-{\text{forget}}\\
& \goneop^{\Ed(T)}(\Set) \ar@{<-}`r/3pt[ur][ur]_-{\varphi^*} &},\] 
so the diagram of their left adjoints 
\[\nicexy@C-1cm{\Set^{\profedtedt}  \ar`d/3pt[dr]_-{\goneop^{\Ed(T)}\circs-} [dr] \ar[rr]^-{\gtwoop^{\Ed(T)}\circs-} && \gtwoop^{\Ed(T)}(\Set)\\ & \goneop^{\Ed(T)}(\Set) \ar`r/3pt[ur][ur]_-{\varphi_!}}\]
is commutative up to a natural isomorphism.  Now we have the isomorphism
\[\begin{split} \varphi_!(\Gone)^{(T,\sigma)} &= \varphi_!\Bigl[\goneop^{\Ed(T)} \circs (\Gonebar)^{(T,\sigma)}\Bigr]\\
&\cong \gtwoop^{\Ed(T)}\circs (\Gtwobar)^{(T,\varphi\sigma)}\\
&= (\Gtwo)^{\Treecat^{\varphi}(T,\sigma)}
\end{split}\]
in which the last equality follows from the definition \[\Treecat^{\varphi}(T,\sigma) = (T,\varphi\sigma) \in \Treecatgtwo.\]  One can similarly check by a direct inspection that these isomorphisms are compatible with morphisms in $\Treecatgone$.
\end{proof}

The next observation is a special case of Theorem \ref{gdash-left-adjoint}.

\begin{corollary}\label{pgs-treecat-operad}
For each action operad $(\G,\omega)$, the diagram
\[\nicexy{\Treecatp \ar[d]_-{\P^-} \ar[r]^-{\Treecat^{\iota}} & \Treecatg \ar[d]_-{\G^-} \ar[r]^-{\Treecat^{\omega}} & \Treecats \ar[d]^-{\S^-}\\
\popset \ar[r]^-{\iota_!} & \gopset \ar[r]^-{\omega_!} & \sopset}\]
is commutative up to natural isomorphisms, in which $\iota : \P \to \G$ is the unique morphism from the planar group operad.
\end{corollary}

\begin{corollary}\label{ptree-free-goperad}
For each action operad $\G$, there is a natural isomorphism \[\G^{(T,\sigma)} \cong \iota_!\P^{(T,\id_{|\inp(T)|})} \in \gopset\] for $\G$-trees $(T,\sigma)$.
\end{corollary}

\begin{proof}
Given a $\G$-tree $(T,\sigma)$ with $n=|\inp(T)|$, there are equalities and a natural isomorphism
\[\begin{split}
\G^{(T,\sigma)} &= \G^{(T,\id_n)}\\
&= \G^{\Treecat^{\iota}(T,\id_n)}\\
&\cong \iota_!\P^{(T,\id_n)}.
\end{split}\]
The first equality is from Proposition \ref{gdash-identifies-objects}.  The second equality holds by the definition of the morphism $\iota$.  The isomorphism holds by Corollary \ref{pgs-treecat-operad}.
\end{proof}

\begin{example}[Ribbon, Braided, and Cactus Cases]
Consider the morphisms 
\[\nicexy{\P\ar@{=}[d] \ar[r]^-{\iota} & \PR \ar[r]^-{\iota} & \R \ar[r]^-{\pi} & \S \ar@{=}[d]\\
\P \ar[r]^-{\iota} & \PB \ar[r]^-{\iota} & \B \ar[r]^-{\pi} & \S\\
\P \ar[r]^-{\iota} \ar@{=}[u] & \PCac \ar[r]^-{\iota} & \Cac \ar[r]^-{\pi} & \S\ar@{=}[u]}\] 
of action operads in Example \ref{ex:all-augmented-group-operads}.  By Theorem \ref{gdash-left-adjoint} the diagrams\label{not:rbdash}\label{not:cacdash}
\[\nicexy{\Treecatp \ar[r]^-{\Treecat^{\iota}} \ar[d]_-{\P^-} &\Treecatpr \ar[d]^-{\PR^-} \ar[r]^-{\Treecat^{\iota}} & \Treecatr \ar[d]^-{\R^-} \ar[r]^-{\Treecat^{\pi}} & \Treecats \ar[d]^-{\S^-}&\\
\popset \ar[r]^-{\iota_!} & \propset \ar[r]^-{\iota_!} & \ropset \ar[r]^-{\pi_!} & \sopset}\]
\[\nicexy{\Treecatp \ar[r]^-{\Treecat^{\iota}} \ar[d]_-{\P^-} &\Treecatpb \ar[d]^-{\PB^-} \ar[r]^-{\Treecat^{\iota}} & \Treecatb \ar[d]^-{\B^-} \ar[r]^-{\Treecat^{\pi}} & \Treecats \ar[d]^-{\S^-}&\\
\popset \ar[r]^-{\iota_!} & \pbopset \ar[r]^-{\iota_!} & \bopset \ar[r]^-{\pi_!} & \sopset}\]
\[\nicexy{\Treecatp \ar[r]^-{\Treecat^{\iota}} \ar[d]_-{\P^-} &\Treecatpcac \ar[d]^-{\PCac^-} \ar[r]^-{\Treecat^{\iota}} & \Treecatcac \ar[d]^-{\Cac^-} \ar[r]^-{\Treecat^{\pi}} & \Treecats \ar[d]^-{\S^-}&\\
\popset \ar[r]^-{\iota_!} & \pcacopset \ar[r]^-{\iota_!} & \cacopset \ar[r]^-{\pi_!} & \sopset}\]
are commutative up to natural isomorphisms.\dqed
\end{example}

\chapter{Contractibility of Group Tree Category}\label{ch:contractibility}

The purpose of this chapter is to show that for an action operad $\G$ with $\G(0)$ the trivial group, the $\G$-tree category $\Treecatg$ is contractible.  In other words, the nerve of the category $\Treecatg$ is a contractible simplicial set.  To prove this, in Section \ref{sec:closed-trees} we first construct a full subcategory $\Treecatug$ of $\Treecatg$ whose objects are $\G$-trees with an empty set of inputs, called closed $\G$-trees.  It is equipped with a functor \[\Cl : \Treecatg\to \Treecatug,\] called the closure functor.  In Section \ref{sec:closure-reflection} we show that $\Treecatug$ is a full reflective subcategory of $\Treecatg$ with the closure functor as the reflection.  This implies that the $\G$-tree category $\Treecatg$ is contractible if and only if the full subcategory $\Treecatug$ is.  

To show that the category $\Treecatug$ is contractible, in Section \ref{sec:output-extension} we show that $\Treecatug$ has an endofunctor, called the output extension functor.  Finally, in Section \ref{sec:contractibility} we prove that $\Treecatug$ is contractible by relating the identity functor and a constant functor via the output extension functor.  In fact, the proof shows that $\Treecatug$ is a strict test category in the sense of Grothendieck; see Remark \ref{rk:strict-test-cat}.  When $\G$ is the symmetric group operad $\S$, we recover Theorem 4.10 in \cite{acm}, which says that the symmetric dendroidal category $\Omega \simeq\Treecats$ is contractible.

\section{Closed Group Trees}\label{sec:closed-trees}

Throughout this chapter, and this chapter only, we assume that $\G$ is an action operad with $\G(0)$ the trivial group.  This is a mild assumption.  For example, 
\begin{itemize}
\item the planar group operad $\P$ in Example \ref{ex:trivial-group-operad}, 
\item the symmetric group operad $\S$ in Example \ref{ex:symmetric-group-operad}, 
\item the braid group operad $\B$ in Definition \ref{def:braid-group-operad}, 
\item the pure braid group operad $\PB$ in Definition \ref{def:pure-braid-group-operad}, 
\item the ribbon group operad $\R$ in Definition \ref{def:ribbon-group-operad}, 
\item the pure ribbon group operad $\PR$ in Definition \ref{def:pure-ribbon-group-operad}, 
\item the cactus group operad $\Cac$ in Definition \ref{def:cactus-group-operad}, and 
\item the pure cactus group operad $\PCac$ in Definition \ref{def:pure-cactus-group-operad} 
\end{itemize}
all satisfy $\G(0)=\{\id_0\}$.

The purpose of this section is to show that there is a functor \[\Cl : \Treecatg \to \Treecatug,\] called the closure functor, from the category $\Treecatg$ of $\G$-trees to the full subcategory $\Treecatug$ consisting of $\G$-trees with an empty set of inputs, called closed $\G$-trees.  In the next section, we will show that $\Treecatug$ is a full reflective subcategory of $\Treecatg$ with the closure functor as reflection.  In Section \ref{sec:output-extension} and Section \ref{sec:contractibility}, we will use the full reflective subcategory $\Treecatug$ to observe that the $\G$-tree category $\Treecatg$ is contractible.  Recall the definition of a planar tree in Definition \ref{def:tree}.  We begin with the definition of the full subcategory $\Treecatug$.

\begin{definition}\label{def:closed-tree}
Suppose $T$ is a planar tree.
\begin{enumerate}
\item An \emph{$n$-ary vertex} $v \in \Vt(T)$ is a vertex such that $|\inp(v)|=n$.  A $0$-ary vertex is also called a\index{closed vertex} \emph{closed vertex}.
\item $T$ is said to be\index{tree!closed}\index{planar tree!closed} \emph{closed} if $\inp(T)=\varnothing$.  A $\G$-tree is\index{G-tree@$\G$-tree!closed} \emph{closed} if its underlying planar tree is closed.  
\item Define $\Treecatug$ as the full subcategory of the $\G$-tree category $\Treecatg$ consisting of closed $\G$-trees.\index{closed G-tree category@closed $\G$-tree category}\index{category!closed G-trees@closed $\G$-trees}
\item For a subset $I \subseteq \inp(T)$, define the\index{tree!closure}\index{closure} \emph{$I$-closure of $T$}, denoted\label{not:iclosure} $\Tunderline^I$, to be the planar tree obtained from $T$ by changing the set of inputs to \[\inp(\Tunderline^I) = \inp(T) \setminus I.\]
\item The \emph{closure of $T$}, denoted\label{not:closureoft} $\Tunderline$, is defined as the $\inp(T)$-closure of $T$, so \[\inp(\Tunderline)=\varnothing.\]
\item For each input $e$ in $T$, the vertex $v_e\in\Vt(\Tunderline)$ with $\out(v_e)=e$ is called the\index{bottom vertex} \emph{bottom vertex induced by $e$}.  
\item The set of bottom vertices in $\Tunderline$ is denoted by\label{not:bvertex} $\BVt(\Tunderline)$.
\item For an element $\sigma \in \G(|\inp(T)|)$, the\index{G-tree@$\G$-tree!closure} \emph{closure of $(T,\sigma)$} is defined as the closed $\G$-tree\label{not:closureoftsigma} \[\Cl(T,\sigma)=(\Tunderline,\id_0),\] where $\id_0\in\G(0)$ is the multiplicative unit.
\end{enumerate}
\end{definition}

\begin{example}
There is a canonical bijection \[\inp(T) \cong \BVt(\Tunderline)\] that sends an input $e \in \inp(T)$ to the bottom vertex $v_e \in \BVt(\Tunderline)$, which is the initial vertex of $e$ in the closure $\Tunderline$.\dqed
\end{example}

\begin{example}\label{ex:closure-exedge}
The closure of the\index{exceptional edge!closure} exceptional edge $\uparrow$ is the $0$-corolla \[\Cl(\uparrow)=\Cor_0 = 
\begin{tikzpicture}\node[xsplain] (v) {}; \draw[outputleg] (v) to +(0,.3); \end{tikzpicture}\] in Example \ref{ex:corolla}.\dqed
\end{example}

\begin{example}[Closure of Corollas]\label{ex:closed-corolla}
Consider\index{corolla!closure} the $4$-corolla and the subset $\{2,3\}\subseteq \inp\bigl(\Cor_4\bigr)$ consisting of the second and the third inputs.  The $\{2,3\}$-closure $\Coru_4^{\{2,3\}}$ of the $4$-corolla is the planar tree
\begin{center}\begin{tikzpicture}
\node[plain] (v) {$v$}; \node[below=.5 of v] (u) {};
\node[plain, left=.1 of u] (u2) {$u_2$}; \node[plain, right=.1 of u] (u3) {$u_3$};
\draw[outputleg] (v) to +(0,.8); \draw[thick] (u3) to (v); \draw[thick] (u2) to (v);
\draw[thick] (v) to +(1,-.6); \draw[thick] (v) to +(-1,-.6);
\end{tikzpicture}\end{center}
with three vertices, two of which are closed, two internal edges, and two inputs.  The closure $\Coru_4$ of the $4$-corolla is the planar tree
\begin{center}\begin{tikzpicture}
\node[plain] (v) {$v$}; \node[below=.5 of v] (u) {}; 
\node[plain, left=.1 of u] (u2) {$u_2$}; \node[plain, left=.1 of u2] (u1) {$u_1$};
\node[plain, right=.1 of u] (u3) {$u_3$}; \node[plain, right=.1 of u3] (u4) {$u_4$};
\draw[outputleg] (v) to +(0,.8); 
\foreach \x in {1,2,3,4} \draw[thick] (u\x) to (v); 
\end{tikzpicture}\end{center}
with five vertices, four of which are bottom vertices, and four internal edges.  In general, the closure $\Coru_n$ of the $n$-corolla has $n+1$ vertices, $n$ of which are bottom vertices, and $n$ internal edges.\dqed
\end{example}

\begin{example}[Closure of Symmetric Trees]\label{ex:closed-symmetric-tree}
Consider\index{symmetric tree!closure} the symmetric tree $(T,\sigma)$ in Example \ref{ex:symmetric-tree}.  Its closure is the closed symmetric tree $(\Tunderline,\id_0)$ with $\Tunderline$ the closed planar tree
\begin{center}\begin{tikzpicture}[scale=1.3]
\node[plain] (w) {$w$}; \node[left=1 of w] (T) {$\Tunderline$}; 
\node[plain, below=2 of w] (c1) {$t_2$};
\node[plain, left=.7 of c1] (c2) {$t_1$}; \node[plain, right=.7 of c1] (c3) {$t_3$};
\node[above=1 of c2, plain] (u) {$u$}; \node[above=1 of c3, plain] (v) {$v$}; 
\draw[outputleg] (w) to +(0,.6);
\draw[thick] (u) to (w); \draw[thick] (v) to (w);
\draw[thick] (c2) to[out=90, in=335] (u); \draw[thick] (c1) to[out=90, in=270] (u);
\draw[thick] (c3) to (v);
\end{tikzpicture}\end{center}
with six vertices, three of which are closed, and five internal edges.\dqed
\end{example}

The next observation says that the closure of a $\G$-tree can be obtained by $\G$-tree substitution into the closure of a corolla.  We also observe that taking the closure and $\G$-tree substitution are commuting operations in a suitable sense.

\begin{proposition}\label{closure-substitution}
Suppose $(T,\sigma)$ is a $\G$-tree with $|\inp(T)|=n$.  Then: 
\begin{enumerate}
\item Its closure is\index{closure!as tree substitution} \[\Cl(T,\sigma) = \bigl(\Coru_n,\id_0\bigr) \comp_v (T,\sigma),\] where:
\begin{itemize}
\item $v$ is the unique vertex in the $n$-corolla $\Cor_n$, regarded as a vertex in the closure $\Coru_n$.
\item $\comp_v$ is the $\G$-tree substitution at $v$.
\end{itemize}
\item $\Vt(\Tunderline) = \Vt(T) \amalg \BVt(\Tunderline)$.
\item Suppose the $\G$-tree substitution\index{closure!of tree substitution} $(U,\sigma^U) \comp_u (T,\sigma)$ at $u\in\Vt(U)$ is defined.  Then its closure is \[\Cl\Bigl((U,\sigma^U) \comp_u (T,\sigma)\Bigr) = \Cl(U,\sigma^U) \comp_u (T,\sigma).\]
\end{enumerate}
\end{proposition}

\begin{proof}
For the first assertion, the equality 
\begin{equation}\label{closure-decomp}
\Tunderline = \Coru_n \comp_v T
\end{equation} of planar trees is immediate from Definition \ref{def:planar-tree-compv}.  That the input equivariance of the $\G$-tree substitution \[\bigl(\Coru_n,\id_0\bigr) \comp_v (T,\sigma)\] is $\id_0\in\G(0)$ follows from our assumption that $\G(0)$ is the trivial group.  

The second assertion follows from the decomposition \eqref{closure-decomp}, Lemma \ref{lem:planar-tree-sub-vertex}(2), and that $\Vt(\Coru_n)$ consists of $v$ and $n$ bottom vertices corresponding to the $n$ inputs of the $n$-corolla.

For the last assertion, we have:
\[\begin{split} \Cl\Bigl((U,\sigma^U) \comp_u (T,\sigma)\Bigr) 
&= \bigl(\Coru_{|\inp(U)|}, \id_0\bigr) \comp_v \Bigl[(U,\sigma^U) \comp_u (T,\sigma)\Bigr]\\
&= \Bigl[\bigl(\Coru_{|\inp(U)|},\id_0\bigr)\comp_v (U,\sigma^U)\Bigr] \comp_u (T,\sigma)\\
&= \Cl(U,\sigma^U) \comp_u (T,\sigma).\end{split}\]
The first and the last equalities follow from the first assertion. The second equality follows from the associativity of $\G$-tree substitution at a vertex in Lemma \ref{lem:planar-tree-sub-vertex}(6).
\end{proof}

\begin{proposition}\label{closure-functor}
Closure defines a functor\index{closure!functoriality} \[\Cl : \Treecatg \to \Treecatug.\]
\end{proposition}

\begin{proof}
We need to define the action of $\Cl$ on morphisms.  Suppose given a morphism\[\nicexy@C+3.5cm{(T,\sigma) \ar[r]^-{\phi\,=\,\bigl(((U,\sigma^U),u); (T_t,\sigma_t)_{t\in\Vt(T)}\bigr)} & (V,\sigma^V)}\in\Treecatg\] as in Definition \ref{def:treecatg}.  To simplify the presentation in the rest of this proof, we will omit writing the input equivariances when there is no danger of confusion.  

To define the action of $\Cl$ on $\phi$, we first make the following definitions:
\begin{enumerate}[label=(\roman*)]
\item Define the planar tree 
\begin{equation}\label{uu}
\begin{split}
U^u &= \frac{U}{\bigl\{f \in \Ed(U) : \out(u)\leq f\bigr\}},\\ \inp(U^u) &= \bigl\{e\in\inp(U) : \out(u) \not< e\bigr\}
\end{split}
\end{equation}
in which the image of $\out(u)$ is \emph{not} an input.  Equivalently, $U^u$ is obtained from $U$ by removing all the edges descended from $\out(u)$, so the image of the vertex $u$ in $U^u$, which will be denoted simply by $u$, is a closed vertex.
\item For each $1 \leq i \leq |\inp(u)|=|\inp(T)|$, define the planar tree 
\begin{equation}\label{uui}
U^u_i=\bigl\{f \in \Ed(U) : \ell_u(i) \leq f\bigr\}
\end{equation}
consisting of
\begin{itemize}
\item the $i$th input $\ell_u(i)$ of $u\in\Vt(U)$ and
\item all the edges descended from it, 
\end{itemize}
with the planar structure induced by $U$.
\end{enumerate}
Then we have 
\begin{equation}\label{vclosure-tclosure}
\begin{split} \Vunderline &= \underline{U\comp_u \bigl[T(T_t)_{t\in\Vt(T)}\bigr]}\\
&= \underline{U} \comp_u \bigl[T(T_t)_{t\in\Vt(T)}\bigr] \\
&= \Uunderline^u \comp_u \Bigl[\Tunderline\bigl\{(\Uunderline_i^u)_{i=1}^{|\inp(T)|}, (T_t)_{t\in \Vt(T)}\bigr\}\Bigr].
\end{split}
\end{equation}
Here $\Uunderline^u$ is the closure of $U^u$.  For each $1 \leq i \leq |\inp(u)|=|\inp(T)|$, the closure $\Uunderline_i^u$ of $U^u_i$ is substituted into the bottom vertex in the closure $\Tunderline$ of $T$ induced by the $\sigmabar(i)$-th input of $T$.  The first equality follows from the decomposition \[V=U\comp_u \bigl[T(T_t)_{t\in\Vt(T)}\bigr].\]  The second equality follows from Proposition \ref{closure-substitution}(3).  The last equality follows from inspection.

We now define $\Cl(\phi)$ as the morphism 
\begin{equation}\label{closure-morphism}
\nicexy@C+3.5cm{\Cl(T,\sigma)=\Tunderline \ar[r]^{\bigl((\Uunderline^u,u); (\Uunderline^u_i)_{i=1}^{|\inp(T)|}, (T_t)_{t\in \Vt(T)}\bigr)} & \Vunderline = \Cl(V,\sigma^V)} \in \Treecatug
\end{equation}
corresponding to \eqref{vclosure-tclosure}.

To see that $\Cl$ preserves identity morphisms, note that for an identity morphism $\Id_{(T,\sigma)}$, the inside and the outside $\G$-trees are all corollas, 
\begin{equation}\label{inoutside-corolla}
\begin{split}
(T_t,\sigma_t) &= (\Cor_t,\id_t) \forspace t \in \Vt(T),\\
(U,\sigma^U) &= (\Cor_T,\id_T),
\end{split}
\end{equation}
in which the subscripts $t$ and $T$ on the right-hand side denote $|\inp(T_t)|$ and $|\inp(T)|$, respectively.  With $u$ denoting the unique vertex in $\Cor_T$, the planar tree $U^u$ in \eqref{uu} consists of a single edge with no inputs.  Therefore, $U^u$ is the $0$-corolla, and so is its closure.  Furthermore, each planar tree $U^u_i$ in \eqref{uui} consists of only the $i$th input of $u$, and it is an input of $U^u_i$.  So $U^u_i$ is the exceptional edge $\uparrow$, and its closure is the $0$-corolla.  It follows that the morphism $\Cl\bigl(\Id_{(T,\sigma)}\bigr)$ is the identity morphism of the closure $\Tunderline$.

To check that $\Cl$ preserves compositions in $\Treecatg$, suppose \[\nicexy@C+3.5cm{(V,\sigma^V) \ar[r]^-{\varphi\,=\,\bigl(((W,\sigma^W),w); (V_v,\sigma_v)_{v\in\Vt(V)}\bigr)} & (X,\sigma^X)} \in \Treecatg\] is another morphism.  Recall from Definition \ref{def:treecatg} that the composition is \[\nicexy@C+2.5cm{(T,\sigma) \ar[r]^-{\varphi\phi\,=\,\bigl((K,u); (H_t)_{t\in\Vt(T)}\bigr)} & (X,\sigma^X)}\] in which:
\[\begin{split}
K &= W \comp_w\Bigl[U\bigl((V_v)_{u\not=v\in \Vt(U)}, \Cor_{\profofu}\bigr)\Bigr],\\
H_t &= T_t(V_v)_{v\in\Vt(T_t)}.\end{split}\]
By definition \eqref{closure-morphism}, $\Cl(\varphi\phi)$ is the morphism\[\nicexy@C+3.5cm{\Cl(T,\sigma)=\Tunderline \ar[r]^{\bigl((\Kunderline^u,u); (\Kunderline^u_i)_{i=1}^{|\inp(T)|}, (H_t)_{t\in \Vt(T)}\bigr)} & \Xunderline = \Cl(X,\sigma^X)} \in \Treecatug\]
in which:
\begin{itemize}
\item $\Kunderline^u$ is the closure of the planar tree\[K^u=\frac{K}{\bigl\{f\in\Ed(K) : \out(u)\leq f\bigr\}}\] defined as in \eqref{uu}.
\item $\Kunderline^u_i$ is the closure of the planar tree\[K^u_i=\bigl\{f\in\Ed(K) : \ell_u(i)\leq f\bigr\}\] defined as in \eqref{uui}.
\end{itemize}

Moreover, by definition \eqref{closure-morphism}, \[\Cl(\varphi) : \Vunderline \to \Xunderline\in\Treecatug\] is the morphism corresponding to the decomposition \[\Xunderline = \Wunderline^w \comp_w \Bigl[\Vunderline\Bigl((\Wunderline^w_j)_{j=1}^{|\inp(V)|}, (V_v)_{v\in\Vt(V)}\Bigr)\Bigr]\] in which:
\begin{itemize}
\item $\Wunderline^w$ is the closure of the planar tree \[W^w=\frac{W}{\bigl\{f\in\Ed(W) : \out(w)\leq f\bigr\}}\] defined as in \eqref{uu}.
\item $\Wunderline^w_j$ is the closure of the planar tree\[W^w_j=\bigl\{f\in\Ed(W) : \ell_w(j)\leq f\bigr\}\] defined as in \eqref{uui}.
\end{itemize}
By Definition \ref{def:treecatg} the composition \[\Cl(\varphi)\Cl(\phi) : \Cl(T,\sigma) \to \Cl(X,\sigma^X)\in\Treecatug\] corresponds to the decomposition\[\Xunderline = A \comp_u \Bigl[\Tunderline\Bigl((D_i)_{i=1}^{|\inp(T)|}, (H_t)_{t\in\Vt(T)}\Bigr)\Bigr]\] 
in which \[\begin{split} A &= \Wunderline^w \comp_w \Bigl[\Uunderline^u \Bigl((\Wunderline^w_j)^{j\in\inp(U)}_{\out(u)\not\leq j}, (V_v)^{v\in\Vt(U)}_{\out(u)\not\leq\out(v)}, \Cor_0\Bigr)\Bigr],\\
D_i&=\Uunderline^u_i\Bigl((\Wunderline_j^w)^{j\in\inp(U)}_{\ell_u(i)\leq j}, (V_v)^{v\in\Vt(U)}_{\ell_u(i)\leq\out(v)}\Bigr)\end{split}\] for $1 \leq i \leq |\inp(T)|$.  By a direct inspection we have\[A=\Kunderline^u \andspace D_i=\Kunderline^u_i\] for each $i$.  So the composition $\Cl(\varphi)\Cl(\phi)$ is equal to $\Cl(\varphi\phi)$.
\end{proof}

\section{Closure as Reflection}\label{sec:closure-reflection}

The purpose of this section is to show that the closure functor is left adjoint to the full subcategory inclusion from $\Treecatug$ to $\Treecatg$.  We continue to assume that $\G(0)$ is the trivial group in our action operad $\G$.  To show that the closure functor is a left adjoint, we will use the following concept.

\begin{definition}
Suppose $(T,\sigma)$ is a $\G$-tree with $|\inp(T)|=n$.  The morphism
\begin{equation}\label{closure-inclusion}
\iota_{(T,\sigma)} : \nicexy@C+4.5cm{(T,\sigma) \ar[r]^{\bigl(((\Coru_n,\id_0),v); (\Cor_{\Prof(t)},\id_{|\inp(t)|})_{t\in\Vt(T)}\bigr)} & \Cl(T,\sigma)}\in \Treecatg,
\end{equation}
corresponding to the $\G$-tree substitution decomposition \[\Cl(T,\sigma) = \bigl(\Coru_n,\id_0\bigr) \comp_v (T,\sigma)\] in Proposition \ref{closure-substitution}(1), is called the\index{canonical inclusion} \emph{canonical inclusion} of $(T,\sigma)$ into its closure.  Here $v$ is the unique vertex in the $n$-corolla $\Cor_n$, regarded as a vertex in the closure $\Coru_n$.
\end{definition}

Canonical inclusions are compatible with morphisms of $\G$-trees in the following sense.

\begin{proposition}\label{closure-coaugmented}
For each morphism\[\nicexy@C+3.5cm{(T,\sigma) \ar[r]^-{\phi\,=\,\bigl(((U,\sigma^U),u); (T_t,\sigma_t)_{t\in\Vt(T)}\bigr)} & (V,\sigma^V)}\in\Treecatg\] as in Definition \ref{def:treecatg}, the diagram
\begin{equation}\label{iota-phi}
\nicexy{(T,\sigma) \ar[r]^-{\phi} \ar[d]_-{\iota_{(T,\sigma)}} & (V,\sigma^V) \ar[d]^-{\iota_{(V,\sigma^V)}}\\
\Cl(T,\sigma) \ar[r]^-{\Cl(\phi)} & \Cl(V,\sigma^V)}
\end{equation}
in $\Treecatg$ is commutative.
\end{proposition}

\begin{proof}
As above, we will omit writing the input equivariances if there is no danger of confusion.  Furthermore, to save space below, we will write $\Cor_{\profofv}$ as $\Cor_v$ for each vertex $v$.   By the definition \eqref{closure-morphism} of $\Cl(\phi)$ and Definition \ref{def:treecatg}, the composition $\Cl(\phi)\iota_{(T,\sigma)}$ corresponds to the $\G$-tree substitution decomposition:
\[\begin{split}
\Cl(V,\sigma^V) &= \Bigl[\Uunderline^u \comp_u \Coru_{|\inp(T)|} \bigl((\Uunderline^u_i)_{i=1}^{|\inp(T)|}, \Cor_{u}\bigr)\Bigr] \comp_u T\bigl(\Cor_{t}(T_t)\bigr)_{t\in\Vt(T)}\\
&= \Uunderline \comp_u T(T_t)_{t\in\Vt(T)}.
\end{split}\]
The other composition $\iota_{(V,\sigma^V)}\phi$ corresponds to the $\G$-tree substitution decomposition:
\[\begin{split}
&\Cl(V,\sigma^V) \\
&= \Bigl[\Coru_{|\inp(U)|} \comp_w U\bigl((\Cor_{v})_{u\not=v \in\Vt(U)}, \Cor_{u}\bigr)\Bigr] \comp_u T\Bigl[T_t(\Cor_{v})_{v\in\Vt(T_t)}\Bigr]_{t\in\Vt(T)}\\
&= \Uunderline \comp_u T(T_t)_{t\in\Vt(T)}.
\end{split}\]
Here $w$ is the unique vertex in the corolla $\Cor_{|\inp(U)|}$, regarded as a vertex in its closure $\Coru_{|\inp(U)|}$.  Therefore, the morphisms $\Cl(\phi)\iota_{(T,\sigma)}$ and $\iota_{(V,\sigma^V)}\phi$ are equal.
\end{proof}

\begin{example}[Closure of Symmetric Trees]\label{ex:closure-stree}
If $\G$ is the\index{closure!symmetric tree}\index{symmetric tree!closure} symmetric group operad $\S$ in Example \ref{ex:symmetric-group-operad}, then the $\S$-tree category $\Treecats$ is equivalent to the dendroidal category $\Omega$ of Moerdijk-Weiss, as we mentioned in Example \ref{ex:treecats}.  In the setting of the dendroidal category $\Omega$, closed symmetric trees are discussed in \cite{acm,moerdijk}.  The existence of the induced morphism $\Cl(\phi)$ in \eqref{closure-morphism}, the functor $\Cl$, and the commutative diagram \eqref{iota-phi} are claimed in 4.4 in \cite{acm}.  However, the detailed formula for the morphism $\Cl(\phi)$ and the detailed proofs of the functoriality of $\Cl$ and of the commutativity of the diagram \eqref{iota-phi}, are not given there.\dqed
\end{example}

Next we observe that $\Treecatug$ is a full reflective subcategory of $\Treecatg$ with the closure functor $\Cl : \Treecatg \to \Treecatug$ as the reflection.

\begin{theorem}\label{treecatug-full-reflective}
Suppose\index{closure!as reflection} $\G$ is an action operad with $\G(0)$ the trivial group.  Then there is an adjunction \[\nicexy{\Treecatg \ar@<2pt>[r]^-{\Cl} & \Treecatug \ar@<2pt>[l]^-{\eta}}\] with
\begin{itemize}
\item $\eta$ the full subcategory inclusion and 
\item $\Cl$ the left adjoint.
\end{itemize}
\end{theorem}

\begin{proof}
Suppose given $(T,\sigma)\in\Treecatg$, $(V,\id_0) \in \Treecatug$, and a morphism \[\nicexy@C+3.5cm{(T,\sigma) \ar[r]^-{\phi\,=\,\bigl(((U,\id_0),u); (T_t,\sigma_t)_{t\in\Vt(T)}\bigr)} & (V,\id_0)}\in\Treecatg\] as in Definition \ref{def:treecatg}.  Since $(V,\id_0)$ is already closed, it is equal to its closure, and its canonical inclusion is the identity morphism.  By Proposition \ref{closure-coaugmented} the diagram \[\nicexy{(T,\sigma) \ar[r]^-{\phi} \ar[d]_-{\iota_{(T,\sigma)}} & (V,\id_0) \ar@{=}[d]\\
\Cl(T,\sigma) \ar[r]^-{\Cl(\phi)} & \Cl(V,\id_0)}\] in $\Treecatg$ is commutative.  The morphism $\Cl(\phi)$ is the unique morphism that completes the above diagram for the following reasons.
\begin{enumerate}
\item $T$ and its closure have the same set of edges.
\item $\Vt(\Tunderline)$ consists of $\Vt(T)$ and the bottom vertices, one for each input of $T$.
\item $(V,\id_0)$ is closed.
\item $V = U \comp_u T(T_t)_{t\in\Vt(T)}$, so $U$ is also closed.
\end{enumerate}
This establishes the adjunction $\Cl \dashv \eta$.
\end{proof}

\section{Output Extension}\label{sec:output-extension}

We continue to assume that $\G(0)$ is the trivial group in our action operad $\G$.  The purpose of this section is to develop the key tool, called the output extension functor, that we will use to show that the category $\Treecatug$ of closed $\G$-trees is contractible.

\begin{definition}\label{def:output-extension}
Suppose $\Lunderline_1 \in \Treecatug$ is the closed $\G$-tree
\begin{center}\begin{tikzpicture}
\node[plain] (z) {$z$}; \node[plain, below=.5 of z] (y) {$y$};
\draw[outputleg] (z) to node[pos=.8]{$r$} +(0,.8); 
\draw[thick] (y) to (z);
\end{tikzpicture}\end{center}
with root $r$, \index{top vertex}\emph{top vertex} $z$, and bottom vertex $y$.  For each closed $\G$-tree $(T,\id_0) \in \Treecatug$, define its\index{output extension} \emph{output extension}\label{not:outputex} \[\Ex(T,\id_0) = (\That,\id_0) \in \Treecatug\] with underlying closed planar tree \[\That=\Lunderline_1 \comp_y T.\]  We also call $\That$ the \emph{output extension of $T$}.
\end{definition}

\begin{ginterpretation}
The output extension can be equivalently defined as\[\Ed(\That) = \{r\} \amalg \Ed(T),\qquad r=\out(\That) < \out(T)\] with the rest of the planar tree structure inherited from $T$.  We can visualize $\That$ as follows.
\begin{center}\begin{tikzpicture}
\node[plain] (z) {$z$}; \node[triangular, below=.5 of z] (t) {$T$};
\draw[outputleg] (z) to node[pos=.8]{$r$} +(0,.8); 
\draw[thick] (t) to (z);
\end{tikzpicture}\end{center}
In $\That$ there is a new root $r$ that is not an input and whose initial vertex is the top vertex $z$ from $\Lunderline_1$.  Observe that\[\Vt(\That) = \{z\} \amalg \Vt(T)\] with the input of $z$ identified with the output of $T$.  \dqed
\end{ginterpretation}

\begin{example}[Output Extension of Closed Corollas]\label{ex:output-ex-closed-corolla}
The output extension $\widehat{\Coru_4}$ of the closed $4$-corolla $\Coru_4$ in Example \ref{ex:closed-corolla} is the closed planar tree
\begin{center}\begin{tikzpicture}
\node[plain] (z) {$z$}; \node[plain, below=.5 of z] (v) {$v$}; 
\node[below=.5 of v] (u) {}; 
\node[plain, left=.1 of u] (u2) {$u_2$}; \node[plain, left=.1 of u2] (u1) {$u_1$};
\node[plain, right=.1 of u] (u3) {$u_3$}; \node[plain, right=.1 of u3] (u4) {$u_4$};
\draw[outputleg] (z) to node[pos=.8]{\footnotesize{$r$}} +(0,.8); 
\draw[thick] (v) to node{\footnotesize{$e$}} (z); 
\foreach \x in {2,3,4} \draw[thick] (u\x) to (v);
\draw[thick] (u1) to node{\footnotesize{$f$}} (v); 
\end{tikzpicture}\end{center}
with six vertices.  The edges $e$ and $f$ will be used in another example below.\dqed
\end{example}

We will use the output extension to connect the identity functor on the category $\Treecatug$ and a constant functor.  First we need to make sure that output extension is functorial, for which we need the following concept.

\begin{definition}\label{def:truncation}
Suppose $U$ is a closed planar tree, and $e$ is an edge in $U$.  Define the planar tree $U_e$ with 
\[\begin{split}
\Ed(U_e) &=\bigl\{f \in\Ed(U) : e \not< f\bigr\},\\
\inp(U_e) &= \{e\},\end{split}\] 
and with the rest of the planar tree structure inherited from $U$.  We call $U_e$ the\index{truncation} \emph{truncation of $U$ below $e$}.
\end{definition}

\begin{ginterpretation}
The truncation $U_e$ below $e$ is obtained from $U$ by (i) removing the edges descended from $e$ and (ii) redefining $e$ as an input.  In particular, $U_e$ is \emph{not} closed.\dqed
\end{ginterpretation}

\begin{example}
Consider the closed planar tree $U=\widehat{\Coru_4}$ in Example \ref{ex:output-ex-closed-corolla}.  Then the truncation $U_r$ below $r$ is the exceptional edge $\uparrow$.  The truncation $U_e$ below $e$ is the $1$-corolla:
\begin{center}\begin{tikzpicture}
\node[plain] (z) {$z$}; 
\draw[outputleg] (z) to node[pos=.8]{\footnotesize{$r$}} +(0,.8); 
\draw[thick] (z) to node[swap]{\footnotesize{$e$}} +(0,-.8);
\end{tikzpicture}\end{center}
The truncation $U_f$ below $f$ is the planar tree
\begin{center}\begin{tikzpicture}
\node[plain] (z) {$z$}; \node[plain, below=.5 of z] (v) {$v$}; 
\node[below=.5 of v] (u) {}; 
\node[plain, left=.1 of u] (u2) {$u_2$}; \node[inner sep=0pt,minimum size=6mm, left=.1 of u2] (u1) {};
\node[plain, right=.1 of u] (u3) {$u_3$}; \node[plain, right=.1 of u3] (u4) {$u_4$};
\draw[outputleg] (z) to node[pos=.8]{\footnotesize{$r$}} +(0,.8); 
\draw[thick] (v) to node{\footnotesize{$e$}} (z); 
\foreach \x in {2,3,4} \draw[thick] (u\x) to (v);
\draw[thick] (u1) to node{\footnotesize{$f$}} (v); 
\end{tikzpicture}\end{center}
whose only input is $f$.\dqed
\end{example}

\begin{lemma}\label{sub-into-output-extension}
Suppose $U$ and $T$ are closed planar trees, and $u$ is a closed vertex in $U$.  Then there are equalities \[\widehat{U\comp_u T} = \Uhat \comp_u T = \That \comp_z \Uhat_{\out(u)}\] in which:
\begin{itemize}
\item $\Uhat_{\out(u)}$ is the truncation of the output extension $\Uhat$ below $\out(u)$.
\item $z$ is the top vertex in $\Lunderline_1$ in Definition \ref{def:output-extension}, regarded as a vertex in the output extension $\That$.
\end{itemize}
\end{lemma}

\begin{proof}
By the vertical associativity of planar tree substitution at a vertex in Lemma \ref{lem:planar-tree-sub-vertex}(6), the first equality follows from the computation: \[\begin{split}
\widehat{U\comp_u T} &= \Lunderline_1 \comp_y \bigl(U\comp_u T\bigr)\\
&= \bigl(\Lunderline_1 \comp_y U\bigr) \comp_u T\\
&= \Uhat \comp_u T.
\end{split}\]
For the other equality, by the definition of $\comp_u$ in Definition \ref{def:planar-tree-sub}, we have \[\Ed\bigl(\Uhat \comp_u T\bigr) = \frac{\{r\}\amalg \Ed(U)\amalg\Ed(T)}{\bigl(\out(u)=\out(T)\bigr)}.\] Since $u$ is a closed vertex in $U$, it is still a closed vertex in the output extension.  It follows that 
\[\begin{split}
\Ed\bigl(\Uhat_{\out(u)}\bigr) &= \bigl\{f\in\{r\}\amalg\Ed(U) : \out(u)\not<f\bigr\}\\ &= \{r\} \amalg\Ed(U).
\end{split}\]
This implies that
\[\begin{split}
\Ed\bigl(\That \comp_z \Uhat_{\out(u)}\bigr)
&= \frac{\{r'\}\amalg\Ed(T) \amalg \{r\}\amalg\Ed(U)}{\bigl(\out(z)=r'=r,\,\inp(z)=\out(u)\bigr)}\\ &= \Ed\bigl(\Uhat \comp_u T\bigr)
\end{split}\]
because in the output extension $\That$, the output of $T$ is identified with  the input of the top vertex $z$.  Moreover, $\Uhat \comp_u T$ and $\That \comp_z \Uhat_{\out(u)}$ are both closed because $\Uhat$ and $\That$ are closed.  Their planar structure also coincide because they are both inherited from $U$ and $T$.
\end{proof}

\begin{ginterpretation}
Writing $e$ for $\out(u)$, one can visualize the output extension $\widehat{U \comp_u T}$ as follows.
\begin{center}\begin{tikzpicture}
\node[plain] (z) {$z$}; \node[triangular, below=.5 of z] (uhat) {$U_{e}$};
\node[triangular, below=.5 of uhat] (T) {$T$};
\draw[outputleg] (z) to node[pos=.8]{\footnotesize{$r$}} +(0,.8); 
\draw[thick] (uhat) to (z); 
\draw[thick] (T) to node{\footnotesize{$e$}} (uhat);
\end{tikzpicture}\end{center}
Here $U_{e}$ is the truncation of $U$ below $e=\out(u)$, which is equal to $U$ with $\out(u)$ redefined as an input.\dqed
\end{ginterpretation}

\begin{proposition}\label{output-ext-functor}
Output extension defines a functor\index{output extension!functoriality} \[\Ex : \Treecatug \to \Treecatug.\]
\end{proposition}

\begin{proof}
We need to define the action of $\Ex$ on a morphism.  Suppose \[\nicexy@C+3.5cm{(T,\id_0) \ar[r]^-{\phi\,=\,\bigl(((U,\id_0),u); (T_t,\sigma_t)_{t\in\Vt(T)}\bigr)} & (V,\id_0)}\in\Treecatug\] is a morphism as in Definition \ref{def:treecatg} with both $T$ and $V$ (and hence also $U$) closed.  To simplify the presentation below, we will omit writing the input equivariances when there is no danger of confusion.  By Lemma \ref{sub-into-output-extension} and the horizontal associativity of planar tree substitution at a vertex in Lemma \ref{lem:planar-tree-sub-vertex}(7), there are equalities:
\begin{equation}\label{phihat-decomp}
\begin{split} \Vhat &= \Ex\Bigl(U \comp_u \bigl(T(T_t)_{t\in\Vt(T)}\bigr)\Bigr)\\
&= \Ex\bigl(T(T_t)_{t\in\Vt(T)}\bigr) \comp_z \Uhat_{\out(u)}\\
&= \bigl[\That(T_t)_{t\in\Vt(T)}\bigr] \comp_z \Uhat_{\out(u)}\\
&= \That\Bigl((T_t)_{t\in\Vt(T)}, \Uhat_{\out(u)}\Bigr).
\end{split}
\end{equation}
In the previous line, $\Uhat_{\out(u)}$ is substituted into the top vertex $z$ in the output extension $\That$, and its input equivariance is $\id_1 \in \G(1)$.  This makes sense because $\Uhat_{\out(u)}$ has only one input, which is created by the truncation of $\Uhat$ below $\out(u)$.

We now define \[\Ex(\phi)=\phihat\] as the morphism 
\begin{equation}\label{output-ext-phi}
\nicexy@C+4cm{\That \ar[r]^-{\phihat\,=\,\bigl((\Cor_0,u); (T_t,\sigma_t)_{t\in\Vt(T)}, \Uhat_{\out(u)}\bigr)} & \Vhat}\in\Treecatug
\end{equation}
corresponding to \eqref{phihat-decomp}.

To see that $\Ex$ preserves identity morphisms, note that for the identity morphism $\Id_{(T,\id_0)}$, the inside and the outside $\G$-trees are all corollas, as in \eqref{inoutside-corolla}.  Since $T$ is closed, $U$ is the $0$-corolla, and its output extension $\Uhat$ is $\Lunderline_1$.  The truncation $\Uhat_{\out(u)}$ below $\out(u)$ is therefore the $1$-corolla.  So $\Ex\bigl(\Id_{(T,\id_0)}\bigr)$ is the identity morphism of the output extension $\That$.

To see that $\Ex$ preserves composition, suppose given a morphism \[\nicexy@C+3.5cm{(V,\id_0) \ar[r]^-{\varphi\,=\,\bigl(((W,\id_0),w); (V_v,\sigma_v)_{v\in\Vt(V)}\bigr)} & (X,\id_0)}\in\Treecatug\] with $X$ (and hence also $W$) closed.  By Definition \ref{def:treecatg} the composite $\varphi\phi$ is \[\nicexy@C+2.5cm{T \ar[r]^-{\varphi\phi\,=\,\bigl((K,u); (H_t)_{t\in\Vt(T)}\bigr)} & X}\in\Treecatug\] in which:
\[\begin{split}
K &= W \comp_w \Bigl[U\bigl((V_v)_{v\in \Vt(U)\setminus\{u\}}, \Cor_0\bigr)\Bigr],\\
H_t &= T_t(V_v)_{v\in\Vt(T_t)}.
\end{split}\]
By \eqref{output-ext-phi} $\Ex(\varphi\phi)$ is the morphism\[\nicexy@C+3.5cm{\That \ar[r]^-{\bigl((\Cor_0,u); (H_t)_{t\in\Vt(T)}, \Khat_{\out(u)}\bigr)} & \Xhat}\in\Treecatug\] in which 
\[\Khat_{\out(u)} = \Bigl\{\What \comp_w \Bigl[U\bigl((V_v)_{v\in \Vt(U)\setminus\{u\}}, \Cor_0\bigr)\Bigr]\Bigr\}_{\out(u)}\] 
by Lemma \ref{sub-into-output-extension}.

On the other hand, by \eqref{output-ext-phi} $\Ex(\varphi) = \varphihat$ is the morphism \[\nicexy@C+4cm{\Vhat \ar[r]^-{\varphihat\,=\,\bigl((\Cor_0,w); (V_v,\sigma_v)_{v\in\Vt(V)}, \What_{\out(w)}\bigr)} & \Xhat}\in\Treecatug.\]  By Definition \ref{def:treecatg} the composite $\Ex(\varphi)\Ex(\phi)$ is given by \[\nicexy@C+3cm{\That \ar[r]^-{\varphihat\phihat\,=\,\bigl((A,u); (H_t)_{t\in\Vt(T)}, B\bigr)} & \Xhat}\in\Treecatug\] in which:
\[\begin{split} A &= \Cor_0 \comp_w \bigl[\Cor_0(\Cor_0)\bigr] = \Cor_0,\\
B &= \Uhat_{\out(u)} \Bigl((V_v)_{v\in \Vt(U)\setminus\{u\}}, \What_{\out(w)}\Bigr)\\
&= \Bigl\{\What \comp_w \Bigl[U\bigl((V_v)_{v\in \Vt(U)\setminus\{u\}}, \Cor_0\bigr)\Bigr]\Bigr\}_{\out(u)}\\
&= \Khat_{\out(u)}.
\end{split}\]
Therefore, the morphisms $\Ex(\varphi\phi)$ and $\Ex(\varphi)\Ex(\phi)$ are equal.
\end{proof}

\begin{remark}
There are no analogues of Proposition \ref{output-ext-functor} for the $\G$-tree category $\Treecatg$.  For a general $\G$-tree $(T,\sigma)$ with $n=|\inp(T)|$, it is possible to define its output extension as the $\G$-tree substitution \[\Ex(T,\sigma)= (L^n,\id_n) \comp_y (T,\sigma) \in \Treecatg\] in which $L^n$ is the planar tree
\begin{center}\begin{tikzpicture}
\node[plain] (z) {$z$}; \node[plain, below=.5 of z] (y) {$y$}; 
\node[below=.1 of y] () {$\cdots$};
\draw[outputleg] (z) to node[pos=.8]{$r$} +(0,.8); 
\draw[thick] (y) to (z);
\draw[thick] (y) to node[swap]{\footnotesize{$1$}} +(-.6,-.6); 
\draw[thick] (y) to node{\footnotesize{$n$}} +(.6,-.6);
\end{tikzpicture}\end{center}
with two vertices and $n$ inputs, all adjacent to the vertex $y$.  If $n=0$ (i.e., $T$ is closed), then $L^0=\Lunderline_1$.  So this definition agrees with Definition \ref{def:output-extension} when $T$ is closed.  This might lead one to think that the proofs above go through for general $\G$-trees that are not necessarily closed.

However, if $T$ and $U$ are not closed, then the equality \[\Uhat\comp_u T= \That\comp_z \Uhat_{\out(u)}\] in Lemma \ref{sub-into-output-extension}, which we used in \eqref{phihat-decomp}, does not hold in general for two reasons.  
\begin{enumerate}\item If $(U,\sigma^U)$ is not a closed $\G$-tree and if $u$ is not a closed vertex, then it is not clear how one should define the input equivariance for the truncation $\Uhat_{\out(u)}$ below $\out(u)$. 
\item The left-hand side $\Uhat \comp_u T$ has the same number of inputs as $U$, while the right-hand side $\That\comp_z \Uhat_{\out(u)}$ has the same number of inputs as $T$.  But in general $|\inp(U)|$ and $|\inp(T)|$ are different.
\end{enumerate}
As a result, there are no reasonable analogues of the morphism $\phihat$ in \eqref{output-ext-phi}, and $\Ex$ does not define a functor on the category $\Treecatg$.  This is the reason why we have to consider the full reflective subcategory $\Treecatug$ of closed $\G$-trees.\dqed
\end{remark}

\section{Contractibility}\label{sec:contractibility}

As in the previous section, we continue to assume that $\G(0)$ is the trivial group.  The main purpose of this section is to show that the $\G$-tree category $\Treecatg$ is contractible; i.e., its nerve is a contractible simplicial set.  We prove this by showing that the full reflective subcategory $\Treecatug$ of closed $\G$-trees is contractible.  This is done by relating the identity functor and a constant functor via the output extension functor $\Ex$.  We will continue to omit writing input equivariances when there is no danger of confusion.

\begin{proposition}\label{output-ex-identity}
Definition \ref{def:output-extension} defines a natural transformation\index{output extension!morphism from identity functor} $\theta$, \[\nicexy@C-.3cm{\Treecatug \ar@/^1pc/[rr]^-{\Id} \ar@/_1pc/[rr]_-{\Ex} & \theta\, \Downarrow & \Treecatug}\]
from the identity functor to $\Ex$ on $\Treecatug$.
\end{proposition}

\begin{proof}
For each closed $\G$-tree $(T,\id_0)$, the morphism \[\theta_{T} : \nicexy@C+2.5cm{T \ar[r]^-{\bigl((\Lunderline_1,y); (\Cor_t)_{t\in\Vt(T)}\bigr)} & \Ex(T)=\That=\Lunderline_1\comp_y T}\in\Treecatug\] is from Definition \ref{def:output-extension}, where $\Cor_t$ is the $|\inp(t)|$-corolla in Example \ref{ex:corolla}.  Suppose \[\nicexy@C+3.5cm{(T,\id_0) \ar[r]^-{\phi\,=\,\bigl(((U,\id_0),u); (T_t,\sigma_t)_{t\in\Vt(T)}\bigr)} & (V,\id_0)}\in\Treecatug\] is a morphism as in Definition \ref{def:treecatg} with both $T$ and $V$ (and hence also $U$) closed.  We must show that the diagram \[\nicexy{T \ar[d]_-{\phi} \ar[r]^-{\theta_T} & \That \ar[d]^{\phihat}\\ V \ar[r]^-{\theta_V} & \Vhat}\] in $\Treecatug$ is commutative.  

By Definition \ref{def:treecatg}, the composite $\theta_V\phi$ corresponds to the decomposition:
\[\begin{split}\Vhat &= \Bigl[\Lunderline_1 \comp_y \bigl[U\bigl((\Cor_v)_{u\not=v\in\Vt(U)},\Cor_u\bigr)\bigr]\Bigr] \comp_u T\Bigl[T_t(\Cor_v)_{v\in\Vt(T_t)}\Bigr]_{t\in\Vt(T)}\\
&= \Uhat \comp_u T(T_t)_{t\in\Vt(T)}.\end{split}\]
On the other hand, by Definition \ref{def:treecatg} and \eqref{output-ext-phi}, the composite $\phihat\theta_T$ corresponds to the decomposition:
\[\begin{split}\Vhat &= \Bigl[\Cor_0 \comp_u \bigl[\Lunderline_1(\Uhat_{\out(u)},\Cor_y)\bigr]\Bigr] \comp_y T\bigl[\Cor_t(T_t)\bigr]_{t\in\Vt(T)}\\
&= \Uhat \comp_u T(T_t)_{t\in\Vt(T)}\end{split}\]
In the last equality above, we used the fact that \[\Uhat= \Lunderline_1(\Uhat_{\out(u)},\Cor_y).\]
Therefore, the morphisms $\theta_V\phi$ and $\phihat\theta_T$ are equal.
\end{proof}

\begin{example}[D\'{e}calage]
For the symmetric group operad $\S$, the symmetric tree category $\Treecats$ is equivalent to the symmetric dendroidal category $\Omega$ of Moerdijk-Weiss  \cite{mw07,mw09}.  In the setting of $\Omega$, Proposition \ref{output-ex-identity} above is a part of Proposition 4.8 in \cite{acm}, in which the output extension functor is denoted by \[D : \overline{\Omega} \to \overline{\Omega}\] and is called the \index{decalage@d\'{e}calage}\emph{d\'{e}calage}.  However, a detailed proof of the functoriality of the output extension functor, which is Proposition \ref{output-ext-functor} above, is not given there.\dqed
\end{example}

Next we relate the output extension functor to a constant functor.  Recall the $0$-corolla $\Cor_0$ in Example \ref{ex:corolla}.

\begin{proposition}\label{output-ex-constant}
There is a natural transformation\index{output extension!morphism from constant functor} \[\nicexy@C-.3cm{\Treecatug \ar@/^1pc/[rr]^-{C_0} \ar@/_1pc/[rr]_-{\Ex} & \xi\, \Downarrow & \Treecatug}\]
from the constant functor at $\Cor_0$,\[\nicexy{\Treecatug \ar[r]^-{C_0} & \Treecatug},\qquad C_0(T)=\Cor_0\] to the output extension functor $\Ex$, defined as the morphism \[\nicexy@C+1.5cm{C_0(T)=\Cor_0 \ar[r]^-{\xi_T\,=\, \bigl((\Cor_0,u); \That\bigr)} & \Cor_0(\That)=\That=\Ex(T) \in \Treecatug}\] for $T\in \Treecatug$.
\end{proposition}

\begin{proof}
Suppose \[\nicexy@C+3.5cm{(T,\id_0) \ar[r]^-{\phi\,=\,\bigl(((U,\id_0),u); (T_t,\sigma_t)_{t\in\Vt(T)}\bigr)} & (V,\id_0)}\in\Treecatug\] is a morphism as in Definition \ref{def:treecatg} with both $T$ and $V$ (and hence also $U$) closed.  We must show that the diagram \[\nicexy{C_0(T)=\Cor_0\ar@{=}[d] \ar[r]^-{\xi_T} & \That \ar[d]^{\phihat}\\ C_0(V)=\Cor_0 \ar[r]^-{\xi_V} & \Vhat}\] in $\Treecatug$ is commutative.  By the definition \eqref{output-ext-phi} of $\phihat$ and Definition \ref{def:treecatg}, the composite $\phihat\xi_T$ corresponds to the decomposition\[\begin{split}\Vhat&=
\Bigl[\Cor_0 \comp_u \bigl[\Cor_0(\Cor_0)\bigr]\Bigr] \comp_u \Cor_0\Bigl[\That\bigl((T_t)_{t\in\Vt(T)},\Uhat_{\out(u)}\bigr)\Bigr]\\
&=\Cor_0 \comp_u \Cor_0(\Vhat)\end{split}\] in which the last equality is from \eqref{phihat-decomp}.  The last line above also corresponds to the morphism $\xi_V$ by definition.
\end{proof}

Recall the nerve functor in \eqref{nerve-category}.  We are now ready for the main observation of this chapter.

\begin{theorem}\label{thm:treecatg-contractible}
Suppose $\G$ is an action operad with $\G(0)$ the trivial group.
Then the nerve of the category\index{G-tree category@$\G$-tree category!contractibility} $\Treecatg$ is contractible.
\end{theorem}

\begin{proof}
By Theorem \ref{treecatug-full-reflective} there is an adjunction $\Cl \dashv \eta$ between the categories $\Treecatg$ and $\Treecatug$.  By basic categorical homotopy theory (see, e.g., \cite{jardine} Section 1), the existence of this adjunction implies that the nerves of $\Treecatg$ and $\Treecatug$ are homotopy equivalent.  Therefore, it suffices to show that the full reflective subcategory $\Treecatug$ of closed $\G$-trees has a contractible nerve.  By Propositions \ref{output-ex-identity} and \ref{output-ex-constant}, there is a zig-zag of natural transformations 
\begin{equation}\label{id-ex-constant}
\nicexy{\Treecatug \ar@{-}[r] \ar@/_2pc/[rr]_-{C_0}^-{}="c" \ar@/^2pc/[rr]^-{\Id}_-{}="i" 
\ar@{}"i";"c" ^(.4){}="a" ^(.6){}="b"
\ar@{=>} "i";"a" _{\theta} \ar@{=>} "c";"b"^{\xi} & \Ex \ar[r] & \Treecatug}
\end{equation}
from the identity functor to the constant functor at $\Cor_0$.  This shows that the category $\Treecatug$ has a contractible nerve.
\end{proof}

\begin{remark}\label{rk:strict-test-cat}
By Corollaire 3.7 in \cite{cm11}, the diagram \eqref{id-ex-constant} shows that the category $\Treecatug$ of closed $\G$-trees is a Grothendieck\index{strict test category}\index{category!strict test} strict test category.\dqed
\end{remark}

\begin{example}
By Theorem \ref{thm:treecatg-contractible}, the $\G$-tree category $\Treecatg$ is contractible for each of the action operads $\P$, $\S$, $\R$, $\PR$, $\B$, $\PB$, $\Cac$, and $\PCac$ in Example \ref{ex:all-augmented-group-operads}.  Recall that there is an equivalence $\Treecats \simeq \Omega$ between the symmetric tree category and the Moerdijk-Weiss dendroidal category $\Omega$.  The  contractibility of the dendroidal category $\Omega$ is Theorem 4.10 in \cite{acm}.\dqed
\end{example}

\chapter{Generalized Reedy Structure}\label{ch:reedy}

In this chapter, we continue to study the categorical structure of the $\G$-tree category $\Treecatg$ in Definition \ref{def:treecatg} for an action operad $\G$.  The main result of this chapter is that $\Treecatg$ is both a dualizable generalized Reedy category and an EZ-category.

In Section \ref{sec:coface_codegen} we first show that each morphism in the $\G$-tree category $\Treecatg$ has a decomposition into codegeneracies, isomorphisms, and cofaces, similar to the familiar decomposition in the finite ordinal category $\Delta$.  Each codegeneracy in $\Treecatg$ corresponds to a $\G$-tree substitution of an exceptional edge into the domain.  Each inner coface corresponds to substituting a $\G$-tree with one internal edge into the domain.  Each outer coface corresponds to substituting the domain into a $\G$-tree with one internal edge.

In Section \ref{sec:dual_reedy} we observe that $\Treecatg$ is a dualizable generalized Reedy category, which is a generalization of a Reedy category due to Berger and Moerdijk that allows non-identity automorphisms.  In the dualizable generalized Reedy structure in $\Treecatg$, the degree function is given by the number of vertices.  The subcategory $\Treecatgplus$ is generated by isomorphisms, inner cofaces, and outer cofaces.  The subcategory $\Treecatgminus$ is generated by isomorphisms and codegeneracies.  

In Section \ref{sec:reedy} we state a model categorical consequence of the dualizable generalized Reedy structure on $\Treecatg$.  This result says that for each cofibrantly generated model category $\M$, the diagram categories $\M^{\Treecatg}$ and $\M^{\Treecatgop}$ both admit Reedy-type model structures.

In Section \ref{sec:ez} we extend the result in Section \ref{sec:dual_reedy} by showing that $\Treecatg$ is an EZ-category.  An EZ-category $\C$, where EZ stands for Eilenberg-Zilber, is a special kind of a dualizable generalized Reedy category with $\C^+$ the subcategory of monomorphisms and $\C^-$ the subcategory of split epimorphisms.

\section{Coface and Codegeneracy}\label{sec:coface_codegen}

We observed in Theorem \ref{thm:treecatg-contractible} that for each action operad $\G$ with $\G(0)$ the trivial group, the $\G$-tree category $\Treecatg$ is contractible, just like the finite ordinal category $\Delta$. In fact, the categories $\Treecatg$ and $\Delta$ have even more in common.  In this section, we show that each morphism in $\Treecatg$ has a codegeneracy-coface decomposition, similar to the epi-monic decomposition in $\Delta$; see, for example, Lemma 8.1.2 in \cite{weibel}.  In the next section, we will extend this observation to a dualizable generalized Reedy structure on the category $\Treecatg$.  Throughout this chapter, $\G$ is an arbitrary action operad.  We do not need to assume that $\G(0)$ is the trivial group.  

In the rest of this chapter, suppose \[\nicexy@C+3.5cm{(T,\sigma) \ar[r]^-{\phi\,=\,\bigl(((U,\sigma^U),u); (T_t,\sigma_t)_{t\in\Vt(T)}\bigr)} & (V,\sigma^V)}\in\Treecatg\] is an arbitrary morphism as in Definition \ref{def:treecatg}.  Recall that $\Cor_n$ denotes the $n$-corolla in Example \ref{ex:corolla}.  We begin with the definition of an inner coface.  

\begin{definition}\label{def:inner-coface}
An\index{inner coface}\index{G-tree category@$\G$-tree category!inner coface} \emph{inner coface} in $\Treecatg$ is a morphism $\phi$ such that:
\begin{itemize}
\item $U=\Cor_{|\inp(U)|}$.
\item There exists a vertex $s \in \Vt(T)$ such that $T_s$ has exactly one internal edge, called the \emph{internal edge corresponding to $\phi$}.
\item $T_t=\Cor_{|\inp(t)|}$ for $t\in\Vt(T)\setminus\{s\}$.
\end{itemize}
An inner coface $\phi$ is said to be\index{pure inner coface} \emph{pure} if 
\begin{itemize}
\item $\sigma^U = \id_{|\inp(U)|} \in \G(|\inp(U)|)$ and 
\item $\sigma_t=\id_{|\inp(t)|}\in\G(|\inp(t)|)$ for $t \in \Vt(T)\setminus\{s\}$.
\end{itemize}
\end{definition}

\begin{interpretation}
In the previous definition, the condition that $T_s$ has one internal edge is equivalent to the condition $|\Vt(T_s)|=2$.  The condition that $T_t$ is the $|\inp(t)|$-corolla is equivalent to the condition $|\Vt(T_t)|=1$.  It follows that \[|\Vt(T)|+1=|\Vt(V)|\] for an inner coface $\phi$.\dqed
\end{interpretation}

The following observation is a special case of the codegeneracy-coface decomposition in $\Treecatg$.

\begin{lemma}\label{inner-coface-iterated}
Suppose \[\phi : (T,\sigma) \to (V,\sigma^V) \in \Treecatg\] is a morphism such that:
\begin{itemize}
\item $(U,\sigma^U)=\bigl(\Cor_{|\inp(U)|},\id_{|\inp(U)|}\bigr)$.
\item There exists a vertex $s \in \Vt(T)$ such that $T_s$ has at least one internal edge.
\item $(T_t,\sigma_t)=\bigl(\Cor_{|\inp(t)|},\id_{|\inp(t)|}\bigr)$ for $t\in\Vt(T)\setminus\{s\}$.
\end{itemize}
Then $\phi$ decomposes as \[\phi=\phi_k \circ\cdots\circ\phi_1\in\Treecatg\] such that:
\begin{itemize}
\item $k$ is the number of internal edges in $T_s$.
\item $\phi_i$ is a pure inner coface for each $1 \leq i \leq k$.
\end{itemize}
\end{lemma}

\begin{proof}
This is proved by an induction on the number $k$ of internal edges in $T_s$, with the initial case $k=1$ being true by definition.  

If $T_s$ has at least two internal edges, we pick an arbitrary internal edge $e$ in it with initial vertex $a$ and terminal vertex $b$.  We will decompose $T_s$ using $e$ as follows.  Suppose $T_s^e$ is the planar tree obtained from $T_s$ by:
\begin{itemize}
\item removing $e\in\Ed(T_s)$;
\item redefining $\out(b)$ as the immediate predecessor of each $f \in \inp(a)$, i.e., combining $a$ and $b$ into a single vertex $b'$ in $T_s^e$ with \[\out(b')=\out(b)\andspace 
\inp(b')=\inp(a) \amalg \bigl[\inp(b)\setminus\{e\}\bigr];\]
\item defining the ordering \[\nicexy{\bigl\{1,\ldots,|\inp(a)|+|\inp(b)|-1\bigr\} \ar[r]^-{\ell_{b'}}_-{\cong} & \inp(b')}\] as \[\ell_{b'}(i)=\begin{cases} \ell_b(i) & \text{if $i < \ell_b^{-1}(e)$},\\
\ell_a\bigl(i-\ell_b^{-1}(e)+1\bigr) &\text{if $\ell_b^{-1}(e) \leq i \leq \ell_b^{-1}(e)-1+|\inp(a)|$},\\
\ell_b\bigl(i-|\inp(a)|+1\bigr) & \text{if $\ell_b^{-1}(e)+|\inp(a)| \leq i \leq |\inp(a)|+|\inp(b)|-1$}.\end{cases}\]
\end{itemize}
The rest of the planar tree structure is inherited from $T_s$.  Geometrically, $T_s^e$ is obtained from $T_s$ by shrinking away the internal edge $e$.  In particular, $T_s^e$ has one fewer internal edge than $T_s$.

Next suppose $T_e$ is the planar tree with \[\begin{split}\Ed(T_e) &= \out(b) \amalg \inp(b) \amalg \inp(a),\\ \out(T_e)&= \out(b),\\ \inp(T_e) &= \inp(a) \amalg \bigl[\inp(b)\setminus\{e\}\bigr],\end{split}\] and the planar tree structure inherited from $T$.  Note that $e$ is the only internal edge in $T_e$.  One can visualize $T_e$ as follows. 
\begin{center}\begin{tikzpicture}
\node[plain] (b) {$b$}; \node[plain, below=.5 of b] (a) {$a$}; 
\node[below=.08 of a] () {\footnotesize{$\inp(a)$}};
\draw[outputleg] (b) to node[pos=.6]{\footnotesize{$\out(b)$}} +(0,.8); 
\draw[thick] (a) to node{\footnotesize{$e$}} (b);
\draw[thick] (a) to +(-.6,-.6); \draw[thick] (a) to +(.6,-.6);
\draw[thick] (b) to +(-.6,-.6); \draw[thick] (b) to +(.6,-.6);
\end{tikzpicture}\end{center}
There is a planar tree decomposition \[T_s = T_s^e \comp_{b'} T_e.\]  This implies that there is a $\G$-tree decomposition 
\begin{equation}\label{tssigmas}
(T_s,\sigma_s) = (T_s^e,\sigma_s) \comp_{b'} (T_e,\id)
\end{equation} 
in which $\id=\id_{|\inp(a)|+|\inp(b)|-1}$.

By assumption the morphism $\phi$ corresponds to the $\G$-tree  decomposition \[(V,\sigma^V) = (T,\sigma) \comp_s (T_s,\sigma_s).\]  Substituting in \eqref{tssigmas} and using the vertical associativity of $\G$-tree substitution in Lemma \ref{lem:g-tree-sub-vertex}(7), we now have the $\G$-tree decomposition \[\begin{split}
(V,\sigma^V) &= (T,\sigma) \comp_s \Bigl[(T_s^e,\sigma_s) \comp_{b'} (T_e,\id)\Bigr]\\ 
&=\Bigl[(T,\sigma) \comp_s (T_s^e,\sigma_s)\Bigr] \comp_{b'} (T_e,\id).
\end{split}\]
The above equality corresponds to a decomposition $\phi=\phi_k\circ\phi'$, \[\nicexy{(T,\sigma) \ar[d]_-{\phi'} \ar[r]^-{\phi} & (V,\sigma^V)\\
(T,\sigma)\comp_s(T_s^e,\sigma_s) \ar[r]^-{\phi_k} & \Bigl[(T,\sigma) \comp_s (T_s^e,\sigma_s)\Bigr] \comp_{b'} (T_e,\id) \ar@{=}[u]}\] in $\Treecatg$ in which $\phi_k$ is a pure inner coface because $T_e$ has exactly one internal edge.  Since $T_s^e$ has one fewer internal edge than $T_s$, the induction hypothesis applies to the morphism $\phi'$, decomposing it into the composite of $k-1$ pure inner cofaces.
\end{proof}

Next we define outer cofaces.

\begin{definition}\label{def:outer-coface}
An\index{outer coface}\index{G-tree category@$\G$-tree category!outer coface} \emph{outer coface} in $\Treecatg$ is a morphism $\phi$ such that:
\begin{itemize}
\item $U$ has exactly one internal edge, called the \emph{internal edge corresponding to $\phi$}.
\item $T_t=\Cor_{|\inp(t)|}$ for $t\in\Vt(T)$.
\end{itemize}
An outer coface $\phi$ is said to be\index{pure outer coface} \emph{pure} if \[\sigma_t=\id_{|\inp(t)|}\] for $t \in \Vt(T)$.
\end{definition}

\begin{interpretation}
In Definition \ref{def:outer-coface} the condition that $U$ has exactly one internal edge is equivalent to $|\Vt(U)|=2$.  It follows that \[|\Vt(T)|+1=|\Vt(V)|\] for an outer coface $\phi$.\dqed
\end{interpretation}

\begin{example}\label{ex:outer-coface-corolla}
Suppose $\sigma\in\G(n)$ for some $n\geq 0$.  Then there are $n+1$ outer cofaces with target $(\Cor_n,\sigma)$.  Indeed, corresponding to the root of the $n$-corolla is the outer coface \[\nicexy@C+2.5cm{(\uparrow,\id_1) \ar[r]^-{\phi_0\,=\, \bigl(((\Cor_{n,0},\sigma),u); \varnothing\bigr)} & (\Cor_n,\sigma)} \in\Treecatg\] in which $\Cor_{n,0}$ is the planar tree
\begin{center}\begin{tikzpicture}
\node[plain] (u) {$u$}; \node[plain, below=.5 of u] (v) {$v$}; 
\node[below=.08 of v] () {$\cdots$};
\draw[outputleg] (u) to +(0,.8); \draw[thick] (v) to (u);
\draw[thick] (v) to node[swap]{\footnotesize{$1$}} +(-.6,-.6); 
\draw[thick] (v) to node{\footnotesize{$n$}} +(.6,-.6);
\end{tikzpicture}\end{center}
with one internal edge and $n$ inputs, all adjacent to the vertex $v$.

On the other hand, for each $1\leq i \leq n$, corresponding to the $i$th input of the $n$-corolla is the outer coface\[\nicexy@C+2.5cm{(\uparrow,\id_1) \ar[r]^-{\phi_i\,=\, \bigl(((\Cor_{n,i},\sigma),u); \varnothing\bigr)} & (\Cor_n,\sigma)} \in\Treecatg.\] Here $\Cor_{n,i}$ is the planar tree
\begin{center}\begin{tikzpicture}
\node[plain] (u) {$u$}; \node[plain, above=.5 of u] (v) {$v$}; 
\draw[outputleg] (v) to +(0,.8); \draw[thick] (u) to node{\footnotesize{$i$}}(v);
\draw[thick] (v) to node[swap]{\footnotesize{$1$}} +(-.6,-.6); 
\draw[thick] (v) to node{\footnotesize{$n$}} +(.6,-.6);
\draw[thick] (u) to +(0,-.8);
\end{tikzpicture}\end{center}
with one internal edge, which is the $i$th input of the vertex $v$, and $n$ inputs, all but one of which are adjacent to $v$.  Moreover, each $\phi_j$ for $0\leq j \leq n$ is a pure outer coface.\dqed
\end{example}

We will need an outer coface analogue of Lemma \ref{inner-coface-iterated}, for which we will use the following concepts.  Recall that the set of internal edges in a planar tree $U$ is denoted by $\Int(U)$ and that a vertex $z$ is the set  $\{\out(z)\} \amalg \inp(z)$, provided $\out(z)$ is not an input of $U$.

\begin{definition}\label{def:removable-vertex}
Suppose $U$ is a planar tree with at least one vertex.
\begin{enumerate}
\item A vertex $z$ in a planar tree $U$ is\index{removable vertex}\index{planar tree!removable vertex} \emph{removable} if there is exactly one internal edge of $U$ in $z$.
\item Suppose $z$ is a removable vertex in $U$, and $\{e\} = z\cap \Int(U)$.  Define $U_z$ as the planar tree obtained from $U$ by  
\begin{enumerate}[label=(\roman*)]
\item removing the edges $z\setminus\{e\}$;
\item defining $e$ as $\out(U_z)$ if $e \in \inp(z)$, or as an input of $U_z$ if $e=\out(z)$.
\end{enumerate}
\end{enumerate}
\end{definition}

\begin{definition}\label{def:ue}
Suppose $U$ is a planar tree.  For each $e \in \Int(U)$, define $U(e)$ as the planar tree with \[\Ed\bigl(U(e)\bigr) = \{\out(U)\}\amalg \{e\} \amalg \inp(U) \subseteq \Ed(U)\] and with the planar tree structure inherited from $U$.
\end{definition}

\begin{ginterpretation}
In Definition \ref{def:removable-vertex}, the planar tree $U_z$ is obtained from $U$ by removing the removable vertex $z$ but keeping the edge $e$ as the new root or an input in $U_z$, depending on whether $e \in \inp(z)$ or $e=\out(z)$.  In particular, we have  
\[\begin{split}\Vt(U_z) &= \Vt(U) \setminus \{z\},\\
\Int(U_z) &= \Int(U) \setminus\{e\}.\end{split}\]
In Definition \ref{def:ue}, the planar tree $U(e)$ is obtained from $U$ by shrinking away all the internal edges in $U$ except $e$.  Its only internal edge is $e$.  If, in addition, $e$ is the only internal edge in a removable vertex $z$ in $U$, then $U(e)$ has two vertices, one of which is $z$.  The other vertex in $U(e)$  has the same profile as $U_z$.\dqed
\end{ginterpretation}

\begin{lemma}\label{lem:removabe-vertex-decomp}
Suppose $z$ is a removable vertex in a planar tree $U$ with $\{e\} = z\cap \Int(U)$.  Then $U$ admits a planar tree decomposition \[U = U(e) \comp_y U_z\] in which $y \in \Vt(U(e))$ is the vertex other than $z$.
\end{lemma}

\begin{proof}
It follows from Definitions \ref{def:removable-vertex} and \ref{def:ue} that $U(e) \comp_y U_z$ and $U$ have the same sets of edges and the same planar tree structures.
\end{proof}

\begin{example}\label{ex:removable-vertex}
Suppose $U$ is the planar tree
\begin{center}\begin{tikzpicture}
\node[plain] (x) {$x$}; \node[plain, below=.5 of x] (w) {$w$};
\node[plain, below left=.6 of w] (u) {$u$}; 
\node[plain, below right=.6 of w] (v) {$v$};
\draw[outputleg] (x) to +(0,.8); 
\foreach \x in {x,v} {\draw[thick] (\x) to +(-.6,-.6); \draw[thick] (\x) to +(.6,-.6);}
\draw[thick] (w) to node{\footnotesize{$g$}} (x);
\draw[thick] (u) to node{\footnotesize{$e$}}(w);
\draw[thick] (v) to node[swap]{\footnotesize{$f$}} (w);
\draw[thick] (u) to +(0,-.6);
\end{tikzpicture}\end{center}
with internal edges $e$, $f$, and $g$.  Then the vertices $u$, $v$, and $x$ are removable, and $w$ is not.  In particular, the only internal edge of $U$ in the removable vertex $x$ is $g$, and the decomposition in Lemma \ref{lem:removabe-vertex-decomp} takes the form \[U= U(g) \comp_y U_x\] in which $U(g)$ and $U_x$ are the following planar trees.
\begin{center}\begin{tikzpicture}
\node[plain] (x) {$x$}; \node[above=.5 of x] {$U(g)$};
\node[plain, below=.5 of x] (y) {$y$};
\draw[outputleg] (x) to +(0,.8); 
\draw[thick] (y) to node{\footnotesize{$g$}} (x);
\foreach \x in{x,y} {\draw[thick] (\x) to +(-.6,-.6); \draw[thick] (\x) to +(.6,-.6);} 
\draw[thick] (y) to +(0,-.7); 
\node[plain, right=3 of x] (w') {$w$}; \node[above=.5 of w'] {$U_x$};
\node[plain, below left=.6 of w'] (u') {$u$}; 
\node[plain, below right=.6 of w'] (v') {$v$};
\draw[outputleg] (w') to node{\footnotesize{$g$}} +(0,.8); 
\draw[thick] (u') to node{\footnotesize{$e$}}(w');
\draw[thick] (v') to node[swap]{\footnotesize{$f$}} (w');
\draw[thick] (u') to +(0,-.6);
\draw[thick] (v') to +(-.6,-.6); \draw[thick] (v') to +(.6,-.6);
\end{tikzpicture}
\end{center}
Here $U(g)$ is obtained from $U$ by (i) shrinking away the internal edges $e$ and $f$ and (ii) combining $u$, $v$, and $w$ into a new vertex $y$.  The planar tree $U_x$ is obtained from $U$ by removing the three edges in $x$ not equal to $g$, which is now the root of $U_x$.\dqed
\end{example}

The following concept will be used to show the existence of certain removable vertices.

\begin{definition}\label{def:trail}
Suppose $U$ is a planar tree with at least one internal edge.  A\index{trail} \emph{trail} in $U$ is a pair \[\bigl((e_1,\ldots,e_n), j\bigr)\] such that the following four conditions hold.
\begin{enumerate}[label=(\roman*)]
\item The $e_i$'s are distinct internal edges in $U$, and $1 \leq j \leq n$.
\item For each $1 \leq i < j$, the terminal vertex of $e_i$ is the initial vertex of $e_{i+1}$.
\item If $j < n$, then the terminal vertex of $e_j$ is also the terminal vertex of $e_{j+1}$.
\item For each $j+1< i \leq n$, the terminal vertex of $e_i$ is the initial vertex of $e_{i-1}$.
\end{enumerate}
We call $n$ the \emph{length} of this trail.
\end{definition}

\begin{interpretation}
In the case $j=n$, a trail is a directed path in $U$.  In the case $j < n$, a trail is a concatenation of two directed paths in $U$ with the same terminal vertex, which is the common terminal vertex of $e_j$ and $e_{j+1}$.\dqed 
\end{interpretation}

\begin{example}\label{ex:trail}
Consider the planar tree $U$ in Example \ref{ex:removable-vertex}.  It has three trails with length $n=2$:
\begin{enumerate}
\item $\bigl((e,g),2\bigr)$.
\item $\bigl((f,g),2\bigr)$.
\item $\bigl((e,f),1\bigr)$.
\end{enumerate}
There is also the trail $\bigl((f,e),1\bigr)$, but we do not distinguish it from $\bigl((e,f),1\bigr)$.\dqed
\end{example}

\begin{lemma}\label{lem:removable-exist}
Suppose $U$ is a planar tree with at least one internal edge.  Then $U$ admits at least two removable vertices.  
\end{lemma}

\begin{proof}
Suppose $L=\bigl((e_1,\ldots,e_n), j\bigr)$ is a trail in $U$ with maximal length.  The initial vertex of $e_1$ is removable, since otherwise the trail can be lengthened, which would contradict the choice of $L$.  If $j=n$, then the terminal vertex of $e_n$ is removable by the maximality of $L$.  If $j<n$, then the initial vertex of $e_n$ is removable by the maximality of $L$.
\end{proof}

The following observation is another special case of the codegeneracy-coface decomposition in $\Treecatg$.

\begin{lemma}\label{outer-coface-iterated}
Suppose \[\phi : (T,\sigma) \to (V,\sigma^V) \in \Treecatg\] is a morphism such that:
\begin{itemize}
\item $U$ has $k\geq 1$ internal edges.
\item $(T_t,\sigma_t)=\bigl(\Cor_{|\inp(t)|},\id_{|\inp(t)|}\bigr)$ for $t\in\Vt(T)$.
\end{itemize}
Then $\phi$ decomposes as \[\phi=\phi_k \circ\cdots\circ\phi_1\in\Treecatg\] such that $\phi_i$ is a pure outer coface for each $1 \leq i \leq k$.
\end{lemma}

\begin{proof}
This is proved by an induction on $k$, with the initial case $k=1$ being true by definition.

Suppose $U$ has at least two internal edges.  By Lemma \ref{lem:removable-exist}, $U$ admits a removable vertex $z\not= u$.  The decomposition \[U = U(e) \comp_y U_z\] in Lemma \ref{lem:removabe-vertex-decomp} implies that there is a $\G$-tree decomposition \[(U,\sigma^U) = \bigl(U(e),\sigma^U\bigr) \comp_y \bigl(U_z,\id_{|\inp(U_z)|}\bigr) \in \Treecatg.\] By definition the morphism $\phi$ corresponds to the $\G$-tree decomposition \[(V,\sigma^V) = (U,\sigma^U) \comp_u (T,\sigma).\]  Combined with the previous decomposition of $(U,\sigma^U)$ and the vertical associativity of  $\G$-tree substitution in Lemma \ref{lem:g-tree-sub-vertex}(7), we now have the $\G$-tree decomposition \[\begin{split}
(V,\sigma^V) &= \Bigl[\bigl(U(e),\sigma^U\bigr) \comp_y \bigl(U_z,\id_{|\inp(U_z)|}\bigr)\Bigr]\comp_u (T,\sigma)\\
&= \bigl(U(e),\sigma^U\bigr) \comp_y \Bigl[\bigl(U_z,\id_{|\inp(U_z)|}\bigr) \comp_u (T,\sigma)\Bigr].
\end{split}\]
This equality corresponds to the decomposition $\phi=\phi_k\circ\phi'$,
\[\nicexy{(T,\sigma) \ar[d]_-{\phi'} \ar[r]^-{\phi} & (V,\sigma^V)\\
\bigl(U_z,\id_{|\inp(U_z)|}\bigr) \comp_u (T,\sigma) \ar[r]^-{\phi_k} & \bigl(U(e),\sigma^U\bigr) \comp_y \Bigl[\bigl(U_z,\id_{|\inp(U_z)|}\bigr)\comp_u (T,\sigma)\Bigr]\ar@{=}[u]}\] in $\Treecatg$ in which $\phi_k$ is a pure outer coface because $U(e)$ has exactly one internal edge.  Since $U_z$ has one fewer internal edge than $U$, the induction hypothesis applies to $\phi'$, decomposing it into the composite of $k-1$ pure outer cofaces.
\end{proof}

Next we define codegeneracies.  Recall the exceptional edge $\uparrow$ in Example \ref{ex:exceptional-edge}.

\begin{definition}\label{def:codegeneracy}
A\index{codegeneracy}\index{G-tree category@$\G$-tree category!codegeneracy} \emph{codegeneracy} in $\Treecatg$ is a morphism $\phi$ such that:
\begin{itemize}
\item $U=\Cor_{|\inp(U)|}$.
\item There exists a vertex $s$ in $T$ such that $T_s=~\uparrow$, called the \emph{exceptional edge corresponding to $\phi$}.
\item $T_t=\Cor_{|\inp(t)|}$ for $t\in\Vt(T)\setminus\{s\}$.
\end{itemize}
A codegeneracy $\phi$ is said to be \emph{pure} if
\begin{itemize}
\item $\sigma^U=\id_{|\inp(U)|}$ and 
\item $\sigma_t=\id_{|\inp(t)|}$ for $t \in \Vt(T)\setminus\{s\}$.
\end{itemize}
\end{definition}

\begin{remark}
In a codegeneracy $\phi : (T,\sigma) \to (V,\sigma^V)$, we have \[|\Vt(T)|=|\Vt(V)|+1.\]\dqed
\end{remark}

We are now ready for the codegeneracy-coface decomposition in $\Treecatg$.

\begin{theorem}\label{thm:codegen-coface-decomp}
Each morphism\index{G-tree category@$\G$-tree category!codegeneracy-coface decomposition} $\phi : (T,\sigma) \to (V,\sigma^V)$ in $\Treecatg$ admits a decomposition as \[\nicexy{(T,\sigma) \ar[rr]^-{\phi} \ar[d]_-{\phi_1} && (V,\sigma^V)\\ (X,\sigma^X) \ar[r]^-{\phi_2}_-{\cong} & (Y,\sigma^Y) \ar[r]^-{\phi_3} & (Z,\sigma^Z) \ar[u]_-{\phi_4}}\]in which:
\begin{itemize}
\item $\phi_1$ is the composite of $|P_0|$ pure codegeneracies, where \[P_0 = \bigl\{t\in \Vt(T) : T_t=\,\uparrow\bigr\}.\]
\item $\phi_2$ is an isomorphism.
\item $\phi_3$ is the composite of $\sum_{t\in P_{\geq 2}} \bigl(|\Vt(T_t)|-1\bigr)$ pure inner cofaces, where \[P_{\geq 2} = \bigl\{t\in \Vt(T) : |\Vt(T_t)|\geq 2 \bigr\}.\]
\item $\phi_4$ is the composite of $|\Vt(U)|-1$ pure outer cofaces if $|\Vt(U)|\geq 2$, or is an isomorphism if $|\Vt(U)|=1$.
\end{itemize}
\end{theorem}

\begin{proof}
We partition the set $\Vt(T)$ as \[\Vt(T) = P_0 \amalg P_1 \amalg P_{\geq 2}\] with $P_0$ and $P_{\geq 2}$ as stated above and \[P_1 = \bigl\{t \in \Vt(T) : T_t=\Cor_t\bigr\},\] where $\Cor_t = \Cor_{|\inp(t)|}$.  Now we define the morphisms $\phi_1$, $\phi_2$, $\phi_3$, and $\phi_4$ corresponding to the following $\G$-tree substitutions:
\[\begin{split}
(X,\sigma^X) &= (T,\sigma)\bigl((\uparrow,\sigma_t)_{t\in P_0}, (\Cor_t,\id_{|\inp(t)|})_{t\in\Vt(T)\setminus P_0}\bigr),\\
(Y,\sigma^Y) &= (X,\sigma^X)\bigl((\Cor_t,\sigma_t)_{t\in P_1}, (\Cor_t,\id_{|\inp(t)|})_{t\in P_{\geq 2}}\bigr),\\
(Z,\sigma^Z) &= (Y,\sigma^Y)\bigl(T_t,\sigma_t\bigr)_{t\in P_{\geq 2}},\\
(V,\sigma^V) &= (U,\sigma^U) \comp_u (Z,\sigma^Z).
\end{split}\]
The last equality is well-defined because \[(Z,\sigma^Z)= (T,\sigma)\bigl(T_t,\sigma_t\bigr)_{t\in\Vt(T)}\] by the horizontal associativity of $\G$-tree substitution in Lemma \ref{lem:g-tree-sub-vertex}(8), which also implies that $\phi_1$ is the composite of $|P_0|$ pure codegeneracies.  The morphism $\phi_2$ is an isomorphism because each $\sigma_t \in \G(|\inp(t)|)$ for $t\in P_1$ has an inverse.  The morphisms $\phi_3$ and $\phi_4$ have the stated properties by Lemma \ref{inner-coface-iterated} and Lemma \ref{outer-coface-iterated}, respectively.
\end{proof}

\begin{example}[Planar Dendroidal Category]\label{ex:planar-epi-monic}
Recall from Example \ref{ex:treecatp} that, for the planar group operad $\P$ in Example \ref{ex:trivial-group-operad}, there is an isomorphism of categories \[\Treecatp\cong\Omega_p,\] in which $\Omega_p$ is the planar dendroidal category in Definition 2.2.1 in \cite{mt}.  The codegeneracy-coface decomposition in $\Omega_p$, analogous to Theorem \ref{thm:codegen-coface-decomp}, is Lemma 2.2.2 in \cite{mt}.  Note that $\phi_2$ is the identity morphism in this case because each level $\P(n)$ is the trivial group.\dqed
\end{example}

\begin{example}[Dendroidal Category]\label{ex:symmetric-epi-monic}
Recall from Example \ref{ex:treecats} that, for the symmetric group operad $\S$ in Example \ref{ex:symmetric-group-operad}, there is an equivalence of categories \[\Treecats\simeq\Omega,\] in which $\Omega$ is the  dendroidal category in Definition 2.3.1 in \cite{mt} and Section 3 in \cite{mw07}.  The codegeneracy-coface decomposition in $\Omega$, analogous to Theorem \ref{thm:codegen-coface-decomp}, is Lemma 2.3.2 in \cite{mt}.\dqed
\end{example}

\begin{example}[Braided, Ribbon, and Cactus Tree Categories]\label{ex:br-epi-monic}
Suppose $\G$ is 
\begin{itemize}
\item the braid group operad $\B$ in Definition \ref{def:braid-group-operad}, 
\item the pure braid group operad $\PB$ in Definition \ref{def:pure-braid-group-operad}, 
\item the ribbon group operad $\R$ in Definition \ref{def:ribbon-group-operad}, 
\item the pure ribbon group operad $\PR$ in Definition \ref{def:pure-ribbon-group-operad},
\item the cactus group operad $\Cac$ in Definition \ref{def:cactus-group-operad}, or 
\item the pure cactus group operad $\PCac$ in Definition \ref{def:pure-cactus-group-operad}.
\end{itemize}
Then Theorem \ref{thm:codegen-coface-decomp} provides a codegeneracy-coface decomposition for each morphism in 
\begin{itemize}
\item the braided tree category $\Treecatb$, 
\item the pure braided tree category $\Treecatpb$, 
\item the ribbon tree category $\Treecatr$,
\item the pure ribbon tree category $\Treecatpr$,
\item the cactus tree category $\Treecatcac$, or
\item the pure cactus tree category $\Treecatpcac$.
\end{itemize}

Note that in the (pure) ribbon case, for each $t\in P_0$, the input equivariance \[\sigma_t\in\R(1)=\PR(1)\cong \bbZ\] can be non-trivial, unlike in the planar, symmetric, (pure) braided, and (pure) cactus cases where $\G(1)$ is the trivial group.\dqed
\end{example}

The following description of an isomorphism in $\Treecatg$ will be useful.

\begin{lemma}\label{treecatg-iso}
A morphism \[\phi : (T,\sigma) \to (V,\sigma^V) \in \Treecatg\] is an isomorphism if and only if
\begin{enumerate}
\item $U=\Cor_{|\inp(U)|}$ and
\item $T_t = \Cor_t$ for each $t\in\Vt(T)$.
\end{enumerate}
\end{lemma}

\begin{proof}
The \emph{if} direction follows from the fact that each $\G(n)$ is a group.  

For the other direction, if there exists a vertex $t$ in $T$ such that $T_t$ is not a corolla, then $T_t$ is
\begin{enumerate}[label=(\roman*)]
\item either an exceptional edge or
\item has at least one internal edge.
\end{enumerate}
In the first case, the map $\phi : \Ed(T)\to\Ed(V)$ is not injective.  In the second case, $\phi$ is not surjective on edge sets.  So in either case $\phi$ cannot be an isomorphism.  In other words, if $\phi$ is an isomorphism, then each $T_t=\Cor_t$, so we have the decomposition \[(V,\sigma)=(U,\sigma^U) \comp_u \Bigl[(T,\sigma)\bigl(\Cor_t,\sigma_t\bigr)_{t\in\Vt(T)}\Bigr].\]  Since $\phi$ is an isomorphism, this implies \[\Vt(T)=\Vt(V)=\Vt(T) \amalg \bigl[\Vt(U)\setminus\{u\}\bigr].\]  It follows that $U$ contains only the vertex $u$, so it is a corolla.
\end{proof}

\begin{example}\label{ex:non-iso}
By Lemma \ref{treecatg-iso}, if a morphism $\phi : (T,\sigma) \to (V,\sigma^V)$ in $\Treecatg$ is an isomorphism, then $|\Vt(T)|=|\Vt(V)|$.  However, this equality in the size of the vertex sets is \emph{not} enough to infer that $\phi$ is an isomorphism.  For instance, there are morphisms \[\nicexy@C+2cm{(\Cor_1,\id_1) \ar[r]^-{\bigl(((\Cor_1,\id_1),u);\,(\uparrow,\id_1)\bigr)} & (\uparrow,\id_1) \ar[d]^-{\out}\\ & (\Cor_1,\id_1)}\] in $\Treecatg$ in which the second morphism is the pure outer coface corresponding to the root of $\Cor_1$ in Example \ref{ex:outer-coface-corolla}.  In the composite of these two morphisms, both edges in the domain $\Cor_1$ are sent to the root of the codomain $\Cor_1$, so it is not an isomorphism.\dqed
\end{example}

\section{Dualizable Generalized Reedy Structure}\label{sec:dual_reedy}

A generalized Reedy category in the sense of \cite{berger-moerdijk-reedy} is a generalization of a classical Reedy category where non-identity automorphisms are allowed.  The purpose of this section is to show that, for each action operad $\G$, the $\G$-tree category $\Treecatg$ is a dualizable generalized Reedy category.  As a consequence, for each cofibrantly generated model category $\M$, the diagram categories $\M^{\Treecatg}$ and $\M^{\Treecatgop}$ both inherit from $\M$ Reedy-type model structures, which we will discuss in the next section.

Recall that a\index{wide subcategory} \emph{wide subcategory} of a category $\C$ is a subcategory that contains all the objects in $\C$.  The wide subcategory generated by all the isomorphisms in $\C$ is denoted by $\Iso(\C)$.  The following is Definition 1.1 in \cite{berger-moerdijk-reedy}.

\begin{definition}\label{def:generalized-reedy}
A\index{generalized Reedy category}\index{category!generalized Reedy} \emph{generalized Reedy structure} on a small category $\C$ is a triple \[(\C^+,\C^-,d)\] consisting of
\begin{enumerate}
\item wide subcategories $\C^+$ and $\C^-$ of $\C$ and 
\item a\index{degree function} \emph{degree function} $d : \Ob(\C) \to \bbN$ 
\end{enumerate}
that satisfies the following four axioms:
\begin{enumerate}[label=(\roman*)]
\item Every non-isomorphism in $\C^+$ (resp., $\C^-$) raises (resp., lowers) degree.  Every isomorphism in $\C$ preserves degree.
\item $\C^+\cap\C^-=\Iso(\C)$.
\item Every morphism $f$ in $\C$ admits a factorization as \[f=f^+f^-\] with $f^+\in\C^+$ and $f^-\in\C^-$ that is unique up to isomorphism.
\item If \[\theta f =f\] for some $\theta \in\Iso(\C)$ and $f\in\C^-$, then $\theta$ is an identity morphism.
\end{enumerate}
A generalized Reedy structure is\index{dualizable generalized Reedy category} \emph{dualizable} if the following additional axiom is satisfied:
\begin{enumerate}[label=(\roman*), resume]
\item If \[f\theta=f\] for some $\theta\in\Iso(\C)$ and $f\in\C^+$, then $\theta$ is an identity morphism.
\end{enumerate}
A \emph{(dualizable) generalized Reedy category} is a small category equipped with a (dualizable) generalized Reedy structure.
\end{definition}

\begin{interpretation}
The first three conditions in Definition \ref{def:generalized-reedy} generalize those in a Reedy category by allowing non-identity automorphisms.  Condition (iv) says that isomorphisms regard morphisms in $\C^-$ as epimorphisms.  Similarly, condition (v) says that isomorphisms regard morphisms in $\C^+$ as monomorphisms.\dqed
\end{interpretation}

\begin{example}[Reedy Categories]
A\index{Reedy category}\index{category!Reedy} \emph{Reedy category} is a tuple \[(\C,\C^+,\C^-,d)\] as in Definition \ref{def:generalized-reedy} that satisfies the following two axioms:
\begin{enumerate}[label=(\roman*)]
\item Every non-identity morphism in $\C^+$ (resp., $\C^-$) raises (resp., lowers) degree.
\item Every morphism in $\C$ admits a unique factorization as \[f=f^+f^-\] with $f^+\in\C^+$ and $f^-\in\C^-$.  
\end{enumerate}
In particular, every Reedy category is also a dualizable generalized Reedy category.  For example:
\begin{itemize}
\item The finite ordinal category $\Delta$, which is a Reedy category, is also a generalized Reedy category.  
\item The planar dendroidal category $\Omega_p$, which does not have any non-identity isomorphisms, is a Reedy category.\dqed
\end{itemize} 
\end{example}

\begin{example}
The opposite category $\C^{\op}$ of a dualizable generalized Reedy category $\C$ is also a dualizable generalized Reedy category, in which \[(\C^{\op})^+=(\C^-)^{\op} \andspace (\C^{\op})^-=(\C^+)^{\op}.\]\dqed
\end{example}

\begin{example}
The dendroidal category $\Omega$ is a dualizable generalized Reedy category by Example 5.3.3(v) in \cite{mt}.  We will prove a generalization of this fact below.\dqed
\end{example}

We now define the structure on the $\G$-tree category $\Treecatg$ that will be shown to be a generalized Reedy structure.

\begin{definition}\label{def:treecatg-reedy}
For an action operad $\G$:
\begin{enumerate}
\item Define \[\deg(T,\sigma)=|\Vt(T)|\] for each $\G$-tree $(T,\sigma)$.
\item Define $\Treecatgplus$ as the wide subcategory of $\Treecatg$ generated by isomorphisms, inner cofaces, and outer cofaces.
\item Define $\Treecatgminus$ as the wide subcategory of $\Treecatg$ generated by isomorphisms and codegeneracies.
\end{enumerate}
\end{definition}

As in Section \ref{sec:coface_codegen}, suppose \[\nicexy@C+3.5cm{(T,\sigma) \ar[r]^-{\phi\,=\,\bigl(((U,\sigma^U),u); (T_t,\sigma_t)_{t\in\Vt(T)}\bigr)} & (V,\sigma^V)}\in\Treecatg\] is an arbitrary morphism as in Definition \ref{def:treecatg}.  We will consider the induced map \[\phi : \Ed(T)\to\Ed(V)\] on edge sets.  

\begin{lemma}\label{treecatg-plus}
For a morphism $\phi \in \Treecatg$, the following conditions are equivalent:
\begin{enumerate}
\item $\phi \in \Treecatgplus$.
\item $T_t\not=\, \uparrow$ for each $t\in\Vt(T)$.
\item $\phi : \Ed(T)\to\Ed(V)$ is injective.
\end{enumerate}
\end{lemma}

\begin{proof}
The implication (2) $\Rightarrow$ (1) follows from Theorem \ref{thm:codegen-coface-decomp} with $P_0 = \varnothing$.  If $\phi$ is an isomorphism, or an inner coface, or an outer coface, then $\phi$ is injective on edge sets.  This yields the implication (1) $\Rightarrow$ (3).  Finally, if $T_t=\, \uparrow$ for some vertex $t$ in $T$, then $\phi$ is not injective on edge sets.  This yields the implication (3) $\Rightarrow$ (2).
\end{proof}

The following descriptions of morphisms in $\Treecatgminus$ are proved similarly.

\begin{lemma}\label{treecatg-minus}
For a morphism $\phi \in \Treecatg$, the following conditions are equivalent:
\begin{enumerate}
\item $\phi \in \Treecatgminus$.
\item $U=\Cor_{|\inp(U)|}$, and $T_t=\Cor_{t}$ or $\uparrow$ for each $t\in\Vt(T)$.
\item $\phi : \Ed(T)\to\Ed(V)$ is surjective.
\end{enumerate}
\end{lemma}

\begin{lemma}\label{treecatg-factorization}
Each morphism $\phi \in \Treecatg$ admits a factorization as \[\phi=\phi^+\phi^-\] with $\phi^+\in\Treecatgplus$ and $\phi^-\in\Treecatgminus$ that is unique up to isomorphism.
\end{lemma}

\begin{proof}
Existence follows from the decomposition \[\phi=(\phi_4\phi_3)(\phi_2\phi_1)\] in Theorem \ref{thm:codegen-coface-decomp}, where \[\phi^+=\phi_4\phi_3\in\Treecatgplus \andspace \phi^-=\phi_2\phi_1\in\Treecatgminus.\]  The uniqueness of this decomposition up to isomorphism follows from:
\begin{itemize}
\item Lemmas \ref{treecatg-iso}, \ref{treecatg-plus}, and \ref{treecatg-minus};
\item the uniqueness of the surjection-injection decomposition of each set map;
\item the fact that each $\G(n)$ is a group.
\end{itemize}
\end{proof}

\begin{theorem}\label{thm:treecatg-reedy}
For each action operad $\G$, the triple\index{G-tree category@$\G$-tree category!generalized Reedy structure}\index{generalized Reedy category!G-tree category@$\G$-tree category}\index{dualizable generalized Reedy category!G-tree category@$\G$-tree category} \[\bigl(\Treecatgplus,\Treecatgminus,\deg\bigr)\] is a dualizable generalized Reedy structure on $\Treecatg$.
\end{theorem}

\begin{proof}
We will use the numbering in Definition \ref{def:generalized-reedy}.  Conditions (i) and (ii) follow from Lemmas \ref{treecatg-iso}, \ref{treecatg-plus}, and \ref{treecatg-minus}.  Condition (iii) holds by Lemma \ref{treecatg-factorization}.  Condition (iv) holds by Lemmas \ref{treecatg-iso} and \ref{treecatg-minus}.  Condition (v) holds by Lemmas \ref{treecatg-iso} and \ref{treecatg-plus}.
\end{proof}

\begin{example}[Planar Tree Category is Reedy]
For the planar group operad $\P$ in Example \ref{ex:trivial-group-operad}, the planar tree category $\Treecatp$ is a Reedy category by Theorem \ref{thm:treecatg-reedy} because $\Treecatp$ has no non-identity isomorphisms.\dqed
\end{example}

\begin{example}[Symmetric, Braided, Ribbon, and Cactus Tree Categories]
For the other action operads in Example \ref{ex:all-augmented-group-operads}, the tree categories $\Treecats$, $\Treecatb$, $\Treecatpb$, $\Treecatr$, $\Treecatpr$, $\Treecatcac$, and $\Treecatpcac$ are dualizable generalized Reedy categories.  Moreover, the functors \[\nicexy{\Treecatp \ar@{=}[d] \ar[r]^-{\Treecat^{\iota}} & \Treecatpr \ar[r]^-{\Treecat^{\iota}} & \Treecatr \ar[r]^-{\Treecat^{\pi}}& \Treecats \ar@{=}[d]\\
\Treecatp \ar[r]^-{\Treecat^{\iota}} & \Treecatpb \ar[r]^-{\Treecat^{\iota}} & \Treecatb \ar[r]^-{\Treecat^{\pi}} & \Treecats\\
\Treecatp \ar@{=}[u] \ar[r]^-{\Treecat^{\iota}} & \Treecatpcac \ar[r]^-{\Treecat^{\iota}} & \Treecatcac \ar[r]^-{\Treecat^{\pi}} & \Treecats\ar@{=}[u]}\] preserve the generalized Reedy structures.\dqed
\end{example}

\section{Reedy-Type Model Structures}\label{sec:reedy}

The purpose of this section is to state the consequence of Theorem \ref{thm:treecatg-reedy} that says that the diagram categories indexed by $\Treecatg$ and $\Treecatgop$, based on any cofibrantly generated model category, admit Reedy-type model structures.  

Before stating the result, we first recall some basic concepts in model categories.   The concept of a model category is originally due to Quillen \cite{quillen}.  Roughly speaking, a model category is a category equipped with enough structure that allows one to make homotopical constructions.  The formulation of a model category below is due to \cite{may-ponto}.  The reader is referred to the references \cite{dwyer-spalinski}, \cite{hirschhorn}, and \cite{hovey} for more detailed discussion of model categories. 

\begin{definition}\label{def:lifting-property}
Suppose $\M$ is a category.
\begin{enumerate}
\item For morphisms $f : A \to B$ and $g : C \to D$ in $\M$, we write \index{left lifting property}\index{right lifting property}\label{notation:fboxslashg}$f \boxslash g$ if for each solid-arrow commutative diagram \[\nicexy{A \ar[d]_-{f} \ar[r] & C \ar[d]^-{g}\\ B \ar[r] \ar@{-->}[ur] & D}\] in $\M$, a dashed arrow exists that makes the entire diagram commutative.  
\item For a class $\cala$ of morphisms in $\M$, we define the classes of morphisms\label{notation:boxslasha}
\[\begin{split}
^{\boxslash}\!\cala &= \bigl\{f \in \M ~\vert~ f \boxslash a \text{ for all } a \in \cala\bigr\},\\
\cala^{\boxslash} &= \bigl\{g \in \M ~\vert~ a \boxslash g \text{ for all } a \in \cala\bigr\}.
\end{split}\]
\item We say that a pair $(\call,\calr)$ of classes of morphisms in $\M$ \index{functorially factors}\emph{functorially factors} $\M$ if each morphism $h$ in $\M$ has a functorial factorization $h=gf$ such that $f \in \call$ and $g \in \calr$.  
\item A \index{weak factorization system}\emph{weak factorization system} in $\M$ is a pair $(\call,\calr)$ of classes of morphisms in $\M$ such that 
\begin{enumerate}[label=(\roman*)]
\item $(\call,\calr)$ functorially factors $\M$, 
\item $\call = \> {}^{\boxslash}\!\calr$, and
\item $\calr = \call^{\boxslash}$.
\end{enumerate}
\end{enumerate}
\end{definition}

\begin{definition}\label{def:model-cat}
A \index{model category}\emph{model category} is a complete and cocomplete category $\M$ equipped with three classes of morphisms\label{notation:wcf} $(\calw, \calc,\calf)$, called \index{weak equivalence}\emph{weak equivalences}, \index{cofibration}\emph{cofibrations}, and \index{fibration}\emph{fibrations}, that satisfy the following two axioms:
\begin{description} 
\item[$2$-out-of-$3$] \index{2of3@$2$-out-of-$3$ property}For any morphisms $f$ and $g$ in $\M$ such that the composition $gf$ is defined, if any two of the three morphisms $f$, $g$, and $gf$ are in $\calw$, then so is the third.
\item[WFS] $(\calc, \calf \cap \calw)$ and $(\calc \cap \calw,\calf)$ are weak factorization systems.
\end{description}
In a model category:
\begin{enumerate}
\item An\index{acyclic fibration}\index{acyclic cofibration} \emph{acyclic (co)fibration} is a morphism that is both a (co)fibration and a weak equivalence.  
\item An object $A\in \M$ is \index{cofibrant object}\emph{cofibrant} if the unique morphism $\varnothing \to A$ from the initial object is a cofibration.
\item An object $A\in \M$ if \index{fibrant object}\emph{fibrant} if the unique morphism $A \to *$ to the terminal object is a fibration.
\end{enumerate}
\end{definition}

The reader is referred to \cite{pinter,shen} for discussion of \index{cardinal}cardinals and \index{ordinal}ordinals.  An ordinal is also regarded as a category with a unique morphism $i \to j$ if and only if $i \leq j$.   The concept of a small object $A$ in the following definition essentially means that a morphism from $A$ to the codomain of a sufficiently long composition must factor through some stage.

\begin{definition}\label{def:small-object}
Suppose $\C$ is a cocomplete category, and $\cali$ is a collection of morphisms in $\C$.
\begin{enumerate}
\item For an ordinal $\alpha$, an \emph{$\alpha$-sequence} in $\C$ is a colimit-preserving functor \[X : \alpha \to \C.\]  The induced morphism \[\nicexy{X_0 \ar[r] & \colimover{\beta<\alpha}\, X_\beta}\] is called a\index{transfinite composition} \emph{transfinite composition}.  
\item For an ordinal $\alpha$ and an $\alpha$-sequence $X$, if each morphism $X_{\beta}\to X_{\beta + 1}$ belongs to $\cali$ for each ordinal $\beta$ satisfying $\beta+1 < \alpha$, then we call the above morphism a \emph{transfinite composition of morphisms in $\cali$}.  
\item A \emph{relative $\cali$-cell complex}\index{relative cell complex} is a transfinite composition of pushouts of morphisms in $\cali$.  The collection of relative $\cali$-cell complexes is denoted by\label{not:cell} $\Cell(\cali)$.
\item For an object $A$ in $\C$ and a cardinal $\kappa$, we say that $A$ is \emph{$\kappa$-small relative to $\cali$} if for
\begin{itemize}
\item each regular cardinal $\alpha \geq \kappa$ and
\item each $\alpha$-sequence $X$ in $\C$ with each morphism $X_{\beta} \to X_{\beta+1}$ in $\cali$ for each ordinal $\beta$ satisfying $\beta+1 < \alpha$,
\end{itemize} 
the induced map of sets
\[\nicexy{\colimover{\beta<\alpha}\, \C(A,X_\beta) \ar[r] & \C\Bigl(A,\colimover{\beta<\alpha}\, X_\beta\Bigr)}\] is a bijection.  
\item We say that $A$ is \emph{small relative to $\cali$} if it is $\kappa$-small relative to $\cali$ for some cardinal $\kappa$. 
\item We say that $\cali$ \emph{permits the small object argument}\index{small object argument} if the domain of each morphism in $\cali$ is small relative to $\Cell(\cali)$.
\end{enumerate}
\end{definition}

\begin{definition}\label{def:cof-gen-model-cat}
A model category $(\M,\calw,\calc,\calf)$ is \index{cofibrantly generated model category}\emph{cofibrantly generated} if it is equipped with two sets $\cali$ and $\calj$ of morphisms such that the following three statements hold:
\begin{enumerate}[label=(\roman*)]
\item Both $\cali$ and $\calj$ permit the small object argument.
\item $\calf = \calj^{\boxslash}$.
\item $\calf \cap \calw = \cali^{\boxslash}$.
\end{enumerate}
\end{definition}

In particular, in a cofibrantly generated model category, fibrations are detected by the set $\calj$, and acyclic fibrations are detected by the set $\cali$.

\begin{example}\label{ex:model-cate}
Here are some basic examples of cofibrantly generated model categories.
\begin{enumerate}
\item The category $\CHau$ \index{space!model structure} of compactly generated weak Hausdorff spaces is a cofibrantly generated model category \cite{hovey} (Section 2.4) in which a weak equivalence is a weak homotopy equivalence, i.e., a map that induces isomorphisms on all homotopy groups for all choices of base points in the domain.  A fibration is a Serre fibration. 
\item The category $\Sset$ of \index{simplicial set!model structure}simplicial sets is a cofibrantly generated model category \cite{hovey} (Chapter 3) in which a weak equivalence is a morphism whose geometric realization is a weak homotopy equivalence.  A cofibration is a level-wise injection.
\item For a field $\fieldk$, the category $\Chaink$ is a \index{chain complex!model structure}cofibrantly generated model category \cite{hovey,quillen} with quasi-isomorphisms as weak equivalences, dimension-wise injections as cofibrations, and dimension-wise surjections as fibrations.
\item The category $\Cat$ of small \index{small category!model structure}categories is a cofibrantly generated model category \cite{rezk}, called the \index{folk model structure}\emph{folk model structure}.  A weak equivalence is an equivalence of categories, i.e., a functor that is full, faithful, and essentially surjective.  A cofibration is a functor that is injective on objects.\dqed
\end{enumerate}
\end{example}
 
\begin{example}\label{ex:cofgen-diagram}
For a cofibrantly generated model category $\M$ and a small category $\D$, the category $\M^{\D}$ of $\D$-diagrams in $\M$ inherits from $\M$ a cofibrantly generated model category structure \cite{hirschhorn} (11.6.1) with weak equivalences and fibrations defined entrywise in $\M$.  In particular, for a group $\Gamma$, regarded as a groupoid with one object, the category $\M^{\Gamma}$, whose objects are those in $\M$ equipped with a right $\Gamma$-action, admits a\index{projective model structure} \emph{projective model structure} with weak equivalences and fibrations defined in $\M$.\dqed
\end{example}

\begin{definition}\label{def:latching-matching}
Suppose $(\C,\C^+,\C^-,d)$ is a generalized Reedy category, and $\M$ is a complete and cocomplete category.  Suppose $c\in\Ob(\C)$, and $X : \C \to \M$ is a functor.
\begin{enumerate}
\item Define the category $\C^+(c)$ in which:
\begin{itemize}
\item Objects are non-invertible morphisms in $\C^+$ with codomain $c$.
\item A morphism $h$ from $f : a\to c$ to $g : b\to c$ in $\C^+(c)$ is a morphism $h : a \to b \in \C$ such that $gh=f$. 
\end{itemize}
\item Define the\index{latching object} \emph{latching object of $X$ at $c$} as the colimit taken in $\M$, \[L_c(X) = \colimover{f\in\C^+(c)} X_{\dom(f)} \in \M,\]
where $\dom(f)$ is the domain of $f$.
\item Define the category $\C^-(c)$ in which:
\begin{itemize}
\item Objects are non-invertible morphisms in $\C^-$ with domain $c$.
\item A morphism $h$ from $f : c\to a$ to $g : c\to b$ in $\C^-(c)$ is a morphism $h : a \to b \in \C$ such that $g=hf$.
\end{itemize} 
\item Define the\index{matching object} \emph{matching object of $X$ at $c$} as the limit taken in $\M$, \[M_c(X) = \limover{f\in\C^-(c)} X_{\codom(f)} \in \M,\]
where $\codom(f)$ is the codomain of $f$.
\end{enumerate}
\end{definition}

\begin{interpretation}
Roughly speaking, the latching object $L_c(X)$ consists of stuff that maps up to $X_c$, while the matching object $M_c(X)$ consists of stuff that $X_c$ maps down to.\dqed
\end{interpretation}

For an object $c\in\C$, $\Aut(c)$ denotes the automorphism group of $c$.

\begin{definition}\label{def:reedy-model-structure}
Suppose $(\C,\C^+,\C^-,d)$ is a generalized Reedy category, and $\M$ is a cofibrantly generated model category.  Suppose $\varphi : X \to Y\in\M^{\C}$.
\begin{enumerate}
\item We call $\varphi$ a \emph{Reedy weak equivalence} if the entry $\varphi_c : X_c\to Y_c$ is a weak equivalence in $\M$ for each object $c\in\C$.
\item We call $\varphi$ a \emph{Reedy cofibration} if for each object $c\in\C$, the induced morphism \[X_c \sqcupover{L_c(X)} L_c(Y) \to Y_c\] is a cofibration in $\M^{\Aut(c)}$.
\item We call $\varphi$ a \emph{Reedy fibration} if for each object $c\in\C$, the induced morphism \[X_c \to M_c(X)\timesover{M_c(Y)} Y_c\] is a fibration in $\M$.
\end{enumerate}
\end{definition}

\begin{corollary}\label{cor:generalized-reedy}
Suppose $\G$ is an action operad, and $\M$ is a cofibrantly generated model category.  Then the diagram categories\index{diagram category!model structure}\index{G-tree category@$\G$-tree category!Reedy-type model structure} $\M^{\Treecatg}$ and $\M^{\Treecatgop}$ are model categories with the Reedy weak equivalences, Reedy cofibrations, and Reedy fibrations in Definition \ref{def:reedy-model-structure}.
\end{corollary}

\begin{proof}
Theorem 1.6 in \cite{berger-moerdijk-reedy} implies that, for each generalized Reedy category $\D$, the diagram category $\M^{\D}$ is a model category with the classes of morphisms in Definition \ref{def:reedy-model-structure}.  Therefore, the result follows from Theorem \ref{thm:treecatg-reedy}, which implies that both $\Treecatg$ and $\Treecatgop$ are generalized Reedy categories.
\end{proof}

\begin{example}
Corollary \ref{cor:generalized-reedy} applies to the tree categories $\Treecatp$, $\Treecats$, $\Treecatb$, $\Treecatpb$, $\Treecatr$, $\Treecatpr$, $\Treecatcac$, and $\Treecatpcac$.  Note that the model structures on the diagram categories $\M^{\Treecatp}$ and $\M^{(\Treecatp)^{\op}}$ are the usual Reedy model structures in Theorem 5.2.5 in \cite{hovey}, since $\Treecatp$ and $(\Treecatp)^{\op}$ do not have any non-identity isomorphisms.\dqed
\end{example}

\section{Eilenberg-Zilber Structure}\label{sec:ez}

The purpose of this section is to show that, for each action operad $\G$, the $\G$-tree category $\Treecatg$ is an EZ-category in the sense of Berger and Moerdijk, where EZ stands for Eilenberg-Zilber.  

To state the definition of an EZ-category, first recall that a commutative square \[\nicexy{a \ar[d] \ar[r] & b\ar[d]\\ c\ar[r] & d}\] in a small category $\C$ is called an \emph{absolute pushout}\index{absolute pushout} if its image under the Yoneda embedding \[\nicexy{\Set^{\Cop}(-,a)\ar[d]\ar[r] & \Set^{\Cop}(-,b)\ar[d]\\ \Set^{\Cop}(-,c) \ar[r] & \Set^{\Cop}(-,d)}\] is a pushout in $\Set^{\Cop}$.  Also recall that a \index{split epimorphism}\emph{split epimorphism} in a category is an epimorphism $f$ that admits a section, i.e., a morphism $s$ such that $fs=\Id_{\codom(f)}$.  The following is Definition 6.7 in \cite{berger-moerdijk-reedy}.

\begin{definition}\label{def:ez}
An\index{EZ-category}\index{category!EZ} \emph{EZ-category} is a small category $\C$ equipped with a degree function $d : \Ob(\C)\to\bbN$ that satisfies the following three axioms:
\begin{enumerate}[label=(\roman*)]
\item A monomorphism preserves (resp., raises) the degree if and only if it is invertible (resp., not invertible).
\item Every morphism factors as a split epimorphism followed by a monomorphism.
\item Each pair of split epimorphisms with a common domain has an absolute pushout. 
\end{enumerate}
\end{definition}

\begin{example}
Each EZ-category $(\C,d)$ is a dualizable generalized Reedy category with:
\begin{itemize}
\item $\C^+$ the wide subcategory generated by the monomorphisms;
\item $\C^-$ the wide subcategory generated by the split epimorphisms.
\end{itemize}  
Both the finite ordinal category $\Delta$ and the symmetric dendroidal category $\Omega$ are EZ-categories by Example 6.8 in \cite{berger-moerdijk-reedy}.  We will prove a generalization of this fact below.\dqed
\end{example}

For an action operad $\G$, we will show that the category $\Treecatg$ is an EZ-category in several steps.  Recall the wide subcategories $\Treecatgplus$ and $\Treecatgminus$ of $\Treecatg$ in Definition \ref{def:treecatg-reedy}.

\begin{lemma}\label{treecatgplus-mono}
The subcategory $\Treecatgplus$ coincides with the subcategory of $\Treecatg$ containing all the monomorphisms.
\end{lemma}

\begin{proof}
In one direction, each isomorphism is trivially a monomorphism.  Also, an inner coface or an outer coface is a monomorphism by an inspection of the formula for the composite of two morphisms in Definition \ref{def:treecatg}.  So each morphism in $\Treecatgplus$ is a composite of monomorphisms, hence a monomorphism.

Conversely, for each monomorphism \[\phi : (T,\sigma) \to (V,\sigma^V)\] in $\Treecatg$, the induced map \[\phi : \Ed(T)\to\Ed(V)\] on edge sets cannot send two distinct edges in $T$ to the same edge in $V$. Therefore, $\phi \in \Treecatgplus$ by Lemma \ref{treecatg-plus}.
\end{proof}

\begin{lemma}\label{treecatgminus-splitepi}
The subcategory $\Treecatgminus$ coincides with the subcategory of $\Treecatg$ containing all the split epimorphisms.
\end{lemma}

\begin{proof}
First suppose \[\phi : (T,\sigma) \to (V,\sigma^V)\] is a split epimorphism.  Then the induced map \[\phi : \Ed(T)\to\Ed(V)\] on edge sets is surjective, so $\phi \in \Treecatgminus$ by Lemma \ref{treecatg-minus}.

Conversely, an isomorphism is trivially a split epimorphism.  If $\phi$ is a codegeneracy, then it has a section $s$, \[\nicexy{(T,\sigma) \ar[r]_-{\phi} & (V,\sigma^V) \ar@/_1pc/[l]_-{s}},\] because we can substitute in the planar tree $\Cor_{n,0}$ in Example \ref{ex:outer-coface-corolla} at the initial vertex of the image of the exceptional edge corresponding to the codegeneracy $\phi$.  The equality $\phi s = \Id_{(V,\sigma^V)}$ implies that $\phi$ is a split epimorphism.
\end{proof}

\begin{lemma}\label{treecatg-absolute-pushout}
In the category $\Treecatg$, each pair of split epimorphisms with a common domain has an absolute pushout. 
\end{lemma}

\begin{proof}
By Lemma \ref{treecatgminus-splitepi} split epimorphisms in $\Treecatg$ are precisely the morphisms in the subcategory $\Treecatgminus$ generated by isomorphisms and codegeneracies.  Therefore, it is enough to show that each pair of codegeneracies with a common domain has an absolute pushout.  Moreover, by the description of an isomorphism in $\Treecatg$ in Lemma \ref{treecatg-iso}, it is actually enough to consider a pair of \emph{pure} codegeneracies with a common domain,\[\nicexy@C+5cm{(T,\sigma) \ar[r]^{\phi_i\,=\,\bigl(((\Cor_{T},\id_{T}),u); (\Cor_{t},\id_{t})_{t_i\not=t\in\Vt(T)},\,(\uparrow,\sigma_i)\bigr)} & (V_i,\sigma^V_i)}\] for $i=1,2$ and $\sigma_i\in\G(1)$, where $(\uparrow,\sigma_i)$ is substituted at the vertex $t_i$ in $T$.  Here \[\Cor_T=\Cor_{|\inp(T)|} \andspace \Cor_t=\Cor_{|\inp(t)|},\] and similarly for $\id_T \in \G(|\inp(T)|)$ and $\id_t \in \G(|\inp(t)|)$.  For each $i=1,2$, since the underlying permutation $\overline{\sigma_i} \in S_1$ is the identity permutation, there is a decomposition \[V_i=T \comp_{t_i}\uparrow\] of planar trees.

If $t_1=t_2$, then $\phi_1$ and $\phi_2$ agree up to an isomorphism, so the pair $\{\phi_1,\phi_2\}$ has an absolute pushout.

Suppose $t_1\not= t_2$.  Then there is a commutative diagram
\begin{equation}\label{tv-absolution-pushout}
\nicexy{(T,\sigma)\ar[d]_-{\phi_1} \ar[r]^-{\phi_2} & (V_2,\sigma^V_2)\ar[d]^-{\phi_3}\\
(V_1,\sigma^V_1)\ar[r]^-{\phi_4} & (V,\sigma^V)}
\end{equation}
in $\Treecatg$ in which:
\begin{itemize}
\item The lower-right corner is defined as \[(V,\sigma^V) = (V_1,\sigma^V_1) \comp_{t_2} (\uparrow,\sigma_2) \in \Treecatg.\] 
\item $\phi_4$ is the pure codegeneracy corresponding to the previous decomposition of $(V,\sigma^V)$.
\item $\phi_3$ is a similar pure codegeneracy followed by an isomorphism.
\end{itemize}
By the argument in the second paragraph of the proof of Lemma \ref{treecatgminus-splitepi}, there exist morphisms $\phi_j'$ as in \[\nicexy@C+.5cm{(T,\sigma)\ar[d]_-{\phi_1} \ar[r]_-{\phi_2} & (V_2,\sigma^V_2) \ar@/_1pc/[l]_-{\phi_2'} \ar[d]_-{\phi_3}\\
(V_1,\sigma^V_1) \ar@/_1pc/[u]_-{\phi_1'} & 
(V,\sigma^V)\ar@/_1pc/[u]_-{\phi_3'}}\] 
that are sections of the corresponding morphisms (i.e., $\phi_j\phi_j'=\Id$ for $j=1,2,3$) such that the equality \[\phi_3'\phi_3 = \phi_2\phi_1'\phi_1\phi_2'\] holds.  By Lemma 3.1.7 in \cite{mt}, this implies that the diagram \eqref{tv-absolution-pushout} is an absolute pushout.
\end{proof}

\begin{theorem}\label{thm:treecatg-ez}
For each action operad $\G$, the category $\Treecatg$ is an\index{EZ-category!G-tree category@$\G$-tree category}\index{G-tree category@$\G$-tree category!EZ-category} EZ-category with the degree function in Definition \ref{def:treecatg-reedy}.
\end{theorem}

\begin{proof}
Using the numbering in Definition \ref{def:ez}, condition (i) follows from Lemmas \ref{treecatg-plus} and \ref{treecatgplus-mono}.  Condition (ii) follows from Lemmas \ref{treecatg-factorization}, \ref{treecatgplus-mono}, and \ref{treecatgminus-splitepi}.  Condition (iii) holds by Lemma \ref{treecatg-absolute-pushout}.
\end{proof}

\begin{example}
The tree categories $\Treecatp$, $\Treecats$, $\Treecatb$, $\Treecatpb$, $\Treecatr$, $\Treecatpr$, $\Treecatcac$, and $\Treecatpcac$ are EZ-categories.\dqed
\end{example}

\chapter{Realization-Nerve Adjunction for Group Operads}
\label{ch:group-nerve}

Throughout this chapter, $\G$ is an action operad.  The main purpose of this chapter is to study the $\G$-operad analogue of the realization-nerve adjunction between the category of simplicial sets and the category of small categories.  The $\G$-operad analogue is an adjunction between the category of $\Treecatg$-presheaves and the category of $\G$-operads in $\Set$.  

In Section \ref{sec:g-tree-set} we define the $\G$-realization functor and the $\G$-nerve functor.  Then we observe that the Yoneda embedding of the $\G$-tree category $\Treecatg$ factors through the $\G$-nerve functor.  

In Section \ref{sec:goperad-colimit} we observe that each $\G$-operad in $\Set$ is canonically the colimit of free $\G$-operads generated by $\G$-trees.  Using this result, in Section \ref{sec:gnerve-fully-faithful} we show that the $\G$-nerve functor is full and faithful.

In Section \ref{sec:gtreeset-monoidal} we equip the category of $\Treecatg$-presheaves with the structure of a symmetric monoidal closed category.  In Section \ref{sec:gnerve-symmetric-monoidal}, we observe that the $\G$-nerve functor is symmetric monoidal, and the $\G$-realization functor is strong symmetric monoidal.

In Section \ref{sec:nerve-change-g} we observe that the Yoneda embedding of $\Treecatg$, the $\G$-nerve functor, and the $\G$-realization functor behave well with respect to a change of the action operad $\G$.  In Section \ref{sec:presheaf-functor-monoidal} we show that the symmetric monoidal structure on the category of $\Treecatg$-presheaves has nice naturality properties with respect to a change of the action operad $\G$.

\section{Realization and Nerve}\label{sec:g-tree-set}

This section has two main purposes:
\begin{enumerate}
\item We first define the $\G$-realization-nerve adjunction.
\item Then we observe that the Yoneda embedding of $\Treecatg$ factors through the $\G$-nerve.  This fact is ultimately due to the fully faithfulness of the functor \[\G^- : \Treecatg \to \gopset\] established in Theorem \ref{gtree-gop-functor}.  
\end{enumerate}
We begin by defining the relevant functors.

\begin{definition}\label{def:gnerve}
Define the categories and functors in the diagram
\begin{equation}\label{realg-nerveg}
\nicexy@C+.4cm{\Treecatg \ar[d]_-{\Yonedag} \ar@{=}[r] & \Treecatg \ar[d]^-{\G^-}\\
\Settotreecatgop \ar@<2pt>[r]^-{\Realg} & \gopset \ar@<2pt>[l]^-{\Nerveg}}
\end{equation}
as follows.
\begin{enumerate}
\item $\G^-$ is the full and faithful functor \[\nicexy{\Treecatg\ar[r]^-{\G^-} & \gopset},\qquad \begin{cases}(T,\sigma) \mapsto \G^{(T,\sigma)},& \\ \phi \mapsto \G^{\phi}\end{cases}\] from the $\G$-tree category $\Treecatg$ to the category of $\G$-operads in $\Set$ in Theorem \ref{gtree-gop-functor}.
\item $\Settotreecatgop$ is\index{G-tree category@$\G$-tree category!presheaves} the category of $\Treecatg$-presheaves, i.e., of functors from $\Treecatgop$ to $\Set$.  Its morphisms are natural transformations between these functors.
\item $\Yonedag$ is the\label{not:yoneda} Yoneda embedding\index{Yoneda embedding} of $\Treecatg$, so \[\Yonedag(T,\sigma) = \Treecatg\bigl(-,(T,\sigma)\bigr) \in\Settotreecatgop\] for each $\G$-tree $(T,\sigma)$.  We call $\Yonedag(T,\sigma)$ the \emph{representable $\Treecatg$-presheaf}\index{representable presheaf} induced by $(T,\sigma)$.
\item $\Realg$\label{not:realg} is the Yoneda extension of $\G^-$, i.e., the left Kan extension of $\G^-$ along the Yoneda embedding $\Yonedag$.  It exists by Proposition \ref{gopm-bicomplete}.  There are coend formulas \[\begin{split}
\Realg(X) &\cong \int^{(T,\sigma)\in\Treecatg} \Bigl[\Settotreecatgop\bigl(\Yonedag(T,\sigma),X\bigr)\Bigr] \cdot \G^{(T,\sigma)}\\
&\cong \int^{(T,\sigma)\in\Treecatg} \coprodover{X(T,\sigma)} \G^{(T,\sigma)} \in \gopset
\end{split}\]
for $X \in \Settotreecatgop$.  We call $\Realg$ the \emph{$\G$-realization functor}.\index{G-realization@$\G$-realization}\index{realization!G-@$\G$-}
\item $\Nerveg$ is\label{not:nerveg} the functor such that \[\Nerveg(\colorc,\O) \in \Settotreecatgop\] is the $\Treecatg$-presheaf defined by \[\Nerveg(\colorc,\O)(T,\sigma) = \gopset\bigl(\G^{(T,\sigma)},(\colorc,\O)\bigr)\] for $(\colorc,\O)\in\gopset$ and $(T,\sigma)\in\Treecatg$.  We call $\Nerveg$ the \emph{$\G$-nerve functor}.\index{G-nerve@$\G$-nerve}\index{nerve!G-@$\G$-}
\end{enumerate}
\end{definition}

In the rest of this chapter, we study the $\G$-realization functor and the $\G$-nerve functor.  First we record the well-understood fact that the functors $\Realg$ and $\Nerveg$ form an adjunction.

\begin{lemma}\label{nerve-real-adjunction}
$\Realg$ is the left adjoint of $\Nerveg$.
\end{lemma}

\begin{proof}
Suppose $X\in\Settotreecatgop$ and $(\colorc,\O)\in\gopset$.  To save space, we will abbreviate $(\colorc,\O)$ to $\O$ below.   Then one computes as follows using (co)end calculus: 
\[\begin{split} \gopset\bigl(\Realg(X),\O\bigr) 
&\cong \gopset\Bigl(\int^{(T,\sigma)\in\Treecatg} \coprodover{X(T,\sigma)} \G^{(T,\sigma)}, \O\Bigr)\\
&\cong \int_{(T,\sigma)\in\Treecatg} \gopset\Bigl(\coprodover{X(T,\sigma)} \G^{(T,\sigma)}, \O\Bigr)\\
&\cong \int_{(T,\sigma)\in\Treecatg} \Set\Bigl(X(T,\sigma), \gopset(\G^{(T,\sigma)}, \O)\Bigr)\\
&= \int_{(T,\sigma)\in\Treecatg} \Set\Bigl(X(T,\sigma), \Nerveg(\O)(T,\sigma)\Bigr)\\
&\cong \Settotreecatgop\bigl(X,\Nerveg(\O)\bigr).
\end{split}\]
In the above calculation, $\int_{(T,\sigma)\in\Treecatg}$ denotes the end indexed by the category $\Treecatg$.  The last isomorphism follows from the usual interpretation of the set of natural transformations in terms of an end.  See, for example, Proposition 6.6.9 in \cite{borceux2} or IX.5 in \cite{maclane}.  
\end{proof}

We call $(\Realg,\Nerveg)$ the \emph{$\G$-realization-nerve adjunction}.  The reader is referred to Proposition 3.2 in \cite{loregian} for a discussion of the realization-nerve adjunction in a more general context.

The next observation says that the Yoneda embedding factors through the $\G$-nerve functor.

\begin{proposition}\label{nerve-factors-yoneda}
The diagram \[\nicexy@C+.4cm{\Treecatg \ar[d]_-{\Yonedag} \ar@{=}[r] & \Treecatg \ar[d]^-{\G^-}\\
\Settotreecatgop & \gopset \ar[l]_-{\Nerveg}}\]
is commutative up to a natural isomorphism.
\end{proposition}

\begin{proof}
For $\G$-trees $(T,\sigma)$ and $(T',\sigma')$, we have that\[\begin{split}
\Nerveg\G^{(T,\sigma)}(T',\sigma') 
&= \gopset\bigl(\G^{(T',\sigma')},\G^{(T,\sigma)}\bigr)\\
&\cong \Treecatg\bigl((T',\sigma'),(T,\sigma)\bigr)\\
&= \Yonedag(T,\sigma)(T',\sigma').
\end{split}\]
The middle isomorphism holds because the functor $\G^-$ is full and faithful by Theorem \ref{gtree-gop-functor}.
\end{proof}

\begin{example}[Planar and Symmetric Dendroidal Nerves]\label{ex:dendroidal-nerve}
When $\G$ is the symmetric group operad $\S$, the symmetric tree category $\Treecats$ is equivalent to the Moerdijk-Weiss symmetric dendroidal category $\Omega$, as mentioned in Example \ref{ex:treecats}.  In the setting of $\Omega$, the $\S$-nerve is called the\index{dendroidal nerve}\index{nerve!dendroidal} \emph{dendroidal nerve} and is denoted by $N_d$ in Example 3.1.4 in \cite{mt}.  It is noted there that the Yoneda embedding factors through the dendroidal nerve, which corresponds to the $\G=\S$ case of Proposition \ref{nerve-factors-yoneda}.  When $\G$ is the planar group operad $\P$, we have the planar analogues of these statements, which apply to the planar dendroidal category $\Omega_p \cong \Treecatp$.\dqed
\end{example}

\begin{example}[Braided Realization-Nerve Adjunction]\label{ex:braided-nerve}
When $\G$ is the braid group operad $\B$ in Definition \ref{def:braid-group-operad}, we call $\Realb$ the\index{realization!braided}\index{braided realization} \emph{braided realization functor}, $\Nerveb$ the\index{braided nerve}\index{nerve!braided} \emph{braided nerve functor}, and \[\nicexy@C+.4cm{\Settotreecatbop \ar@<2pt>[r]^-{\Realb} & \bopset \ar@<2pt>[l]^-{\Nerveb}}\] the \emph{braided realization-nerve adjunction}.  There is a natural isomorphism \[\Yonedab \cong \Nerveb\B^-\] by Proposition \ref{nerve-factors-yoneda}.  Similarly, for the pure braid group operad $\PB$ in Definition \ref{def:pure-braid-group-operad}, there are a \emph{pure braided realization-nerve adjunction} \[\nicexy@C+.4cm{\Settotreecatpbop \ar@<2pt>[r]^-{\Realpb} & \pbopset \ar@<2pt>[l]^-{\Nervepb}}\]  and a natural isomorphism $\Yonedapb \cong \Nervepb \PB^-$.\dqed
\end{example}

\begin{example}[Ribbon Realization-Nerve Adjunction]\label{ex:ribbon-nerve}
When $\G$ is the ribbon group operad $\R$ in Definition \ref{def:ribbon-group-operad}, we call $\Realr$ the\index{realization!ribbon}\index{ribbon realization} \emph{ribbon realization functor}, $\Nerver$ the\index{ribbon nerve}\index{nerve!ribbon} \emph{ribbon nerve functor}, and \[\nicexy@C+.4cm{\Settotreecatrop \ar@<2pt>[r]^-{\Realr} & \ropset \ar@<2pt>[l]^-{\Nerver}}\] the \emph{ribbon realization-nerve adjunction}.  There is a natural isomorphism \[\Yonedar \cong \Nerver\R^-\] by Proposition \ref{nerve-factors-yoneda}.  Similarly, for the pure ribbon group operad $\PR$ in Definition \ref{def:pure-ribbon-group-operad}, there are a \emph{pure ribbon realization-nerve adjunction} \[\nicexy@C+.4cm{\Settotreecatprop \ar@<2pt>[r]^-{\Realpr} & \propset \ar@<2pt>[l]^-{\Nervepr}}\] and a natural isomorphism $\Yonedapr \cong \Nervepr \PR^-$.\dqed
\end{example}

\begin{example}[Cactus Realization-Nerve Adjunction]\label{ex:cactus-nerve}
When $\G$ is the cactus group operad $\Cac$ in Definition \ref{def:cactus-group-operad}, we call $\Realcac$ the\index{realization!cactus}\index{cactus realization} \emph{cactus realization functor}, $\Nervecac$ the\index{cactus nerve}\index{nerve!cactus} \emph{cactus nerve functor}, and \[\nicexy@C+.4cm{\Settotreecatcacop \ar@<2pt>[r]^-{\Realcac} & \cacopset \ar@<2pt>[l]^-{\Nervecac}}\] the \emph{cactus realization-nerve adjunction}.  There is a natural isomorphism \[\Yonedacac \cong \Nervecac\Cac^-\] by Proposition \ref{nerve-factors-yoneda}.  Similarly, for the pure cactus group operad $\PCac$ in Definition \ref{def:pure-cactus-group-operad}, there are a \emph{pure cactus realization-nerve adjunction} \[\nicexy@C+.4cm{\Settotreecatpcacop \ar@<2pt>[r]^-{\Realpcac} & \pcacopset \ar@<2pt>[l]^-{\Nervepcac}}\] and a natural isomorphism $\Yonedapcac \cong \Nervepcac \PCac^-$.\dqed
\end{example}

\section{Group Operads as Colimits}\label{sec:goperad-colimit}

The purpose of this section is to show that every $\G$-operad in $\Set$ is canonically isomorphic to a colimit\index{G-operad@$\G$-operad!as colimit} of free $\G$-operads generated by $\G$-trees.  First we need the following definitions.

\begin{definition}\label{def:goperad-colim}
Suppose $\O$ is a $\colorc$-colored $\G$-operad in $\Set$.
\begin{enumerate}
\item Denote by $\Govero$ the category of $\G^-$ over $\O$.  Explicitly:
\begin{enumerate}[label=(\roman*)]
\item An object in $\Govero$ is a pair 
\begin{equation}\label{govero-object}
\Bigl((T,\sigma)\in\Treecatg; \narrowxy{\bigl(\Ed(T), \G^{(T,\sigma)}\bigr) \ar[r]^-{
\phi} & (\colorc,\O)}\in\gopset\Bigr)
\end{equation}
consisting of a $\G$-tree $(T,\sigma)$ and a morphism $\phi$ from the $\G$-operad freely generated by $(T,\sigma)$ to $(\colorc,\O)$.
\item A morphism \[\varphi : \bigl((T,\sigma);\phi\bigr)\to \bigl((T';\sigma'),\phi'\bigr)\] in $\Govero$ is a morphism $\varphi : (T,\sigma)\to (T',\sigma')$ in $\Treecatg$ such that the diagram 
\begin{equation}\label{govero-morphism}
\nicexy{\bigl(\Ed(T),\G^{(T,\sigma)}\bigr) \ar[r]^-{\phi} \ar[d]^{\G^{\varphi}} & (\colorc,\O) \ar@{=}[d]\\ 
\bigl(\Ed(T'),\G^{(T',\sigma')}\bigr) \ar[r]^-{\phi'} & (\colorc,\O)}
\end{equation}
in $\gopset$ commutes.
\end{enumerate} 
\item Define the functor\label{not:gpr} \[\Gpr : \Govero \to \gopset,\qquad \begin{cases}\Gpr\bigl((T,\sigma); \phi\bigr) = \bigl(\Ed(T),\G^{(T,\sigma)}\bigr),\\
\Gpr(\varphi) = \G^{\varphi}.\end{cases}\]
\end{enumerate}
\end{definition}

\begin{ginterpretation}\label{int:decoration-of-tree}
We think of the morphism $\phi\in\gopset$ in \eqref{govero-object} as an $\O$-decoration of the planar tree $T$.  Indeed, the map \[\phi : \Ed(T)\to\colorc\] is a $\colorc$-coloring of $T$.  Since the $\Ed(T)$-colored $\G$-operad $\G^{(T,\sigma)}$ is freely generated by the set of vertices in $T$, the morphism \[\phi : \G^{(T,\sigma)}\to\O\] of $\G$-operads is determined by the map \[\Vt(T) \ni t \mapsto \phi(t)\in\O\sbinom{\phi\,\out(t)}{\phi\,\inp(t)}\] that sends each vertex in $T$ to an element in $\O$ with the correct profile.  So an object $\bigl((T,\sigma);\phi\bigr)$ in $\Govero$ consists of a $\G$-tree $(T,\sigma)$ and a specific $\O$-decoration of it.  The commutative diagram \eqref{govero-morphism} expresses the idea that in $\Govero$ a morphism must be compatible with the $\O$-decorations of the domain and the codomain.  \dqed
\end{ginterpretation}

We will show that the colimit of the functor $\Gpr$ is isomorphic to the given $\G$-operad $\O$.  To do that, we will use the following morphisms.

\begin{definition}\label{def:morphisms-ab}
Suppose $\O$ is a $\colorc$-colored $\G$-operad in $\Set$.
\begin{enumerate}
\item Denote by 
\begin{equation}\label{morphism-a}
A : \colimover{\Govero}\, \Gpr \to (\colorc,\O) \in\gopset
\end{equation}
the natural morphism determined by the commutative diagrams \[\nicexy{\Gpr\bigl((T,\sigma); \phi\bigr) \ar@{=}[r] \ar[d]_-{\text{natural}} & \bigl(\Ed(T),\G^{(T,\sigma)}\bigr) \ar[d]^-{\phi}\\
\colimover{\Govero}\, \Gpr \ar[r]^-{A} & (\colorc,\O)}\] for $\bigl((T,\sigma); \phi\bigr) \in \Govero$.
\item Define the morphism 
\begin{equation}\label{morphism-b}
B : (\colorc,\O)\to \colimover{\Govero}\,\Gpr \in\gopset 
\end{equation}
as follows.
\begin{enumerate}[label=(\roman*)]
\item For each element $c\in\colorc$, $B(c)$ is the image of the unique element in the color set of $\G^{(\uparrow,\id_1)}$, namely $\Ed(\uparrow)=\{\uparrow\}$, under the natural morphism \[\nicexy{\Gpr\bigl((\uparrow,\id_1); \phi_c\bigr) \ar[r]^-{\omega} & \colimover{\Govero}\,\Gpr}\] in which the morphism \[\narrowxy{\bigl(\{\uparrow\},\G^{(\uparrow,\id_1)}\bigr) \ar[r]^-{\phi_c} & (\colorc,\O)}\in\gopset\] is determined by $\phi_c(\uparrow)=c$.
\item For each element $x \in \O\duc$ with $\duc\in\Profcc$ and $n=|\uc|$, there is a corresponding morphism \[\nicexy{\Bigl(\Ed(\Cor_{(\uc;d)}), \G^{(\Cor_{(\uc;d)},\id_n)}\Bigr) \ar[r]^-{\phi_x} & (\colorc,\O)} \in \gopset\] that sends the unique vertex $v$ in the $\duc$-corolla to $x$.  Then $B(x)$ is the image of $v$ under the natural morphism \[\nicexy{\Gpr\bigl((\Cor_{(\uc;d)},\id_n); \phi_x\bigr) \ar[r]^-{\omega} & \colimover{\Govero}\,\Gpr}.\]
\end{enumerate}
\end{enumerate}
\end{definition}

\begin{interpretation}
The colimit in the domain of the morphism $A$ in \eqref{morphism-a} is made up of objects of the form $\Gpr\bigl((T,\sigma); \phi\bigr)$, which is the $\G$-operad freely generated by a $\G$-tree $(T,\sigma)$.  We think of the morphism $A$ as an approximation of the $\colorc$-colored $\G$-operad $\O$ using $\G$-operads freely generated by $\G$-trees.  The morphism $B$ in \eqref{morphism-b} assigns:
\begin{itemize}
\item to each color $c\in\colorc$ the $c$-colored exceptional edge;
\item to each element $x\in\O$ the corolla whose vertex is decorated by $x$.\dqed
\end{itemize}
\end{interpretation}

While $A$ in \eqref{morphism-a} is a morphism of $\G$-operads by construction, it is not completely obvious why $B$ in \eqref{morphism-b} is a morphism of $\G$-operads.  We now check that this is indeed the case.

\begin{lemma}\label{B-is-morphism}
$B$ in \eqref{morphism-b} is indeed a morphism of $\G$-operads.
\end{lemma}

\begin{proof}
It suffices to show that $B$ is compatible with the $\G$-operad structure morphisms in Corollary \ref{goperad-gammat}.  Suppose $\O$ is a $\colorc$-colored $\G$-operad in $\Set$, and $(T,\sigma)$ is a $\colorc$-colored $\G$-tree with profile $\duc$ and $|\uc|=n$.  We write $\colim_{\Govero}\,\Gpr$ as $\colim\,\Gpr$, and write $\colord$ for the color set of $\colim\,\Gpr$.  We must show that the diagram
\begin{equation}\label{B-lambda-compatible}
\nicexy@C+.5cm{\O(T,\sigma)=\prodover{t\in\Vt(T)}\O\inoutt \ar[d]_-{\gamma^{\O}_{(T,\sigma)}} \ar[r]^-{\prodover{t\in\Vt(T)} B} & \prodover{t\in\Vt(T)} \bigl(\colim\,\Gpr\bigr)\Binoutt \ar[d]^-{\gamma^{\colim\,\Gpr}_{(BT,\sigma)}} \\
\O\duc \ar[r]^-{B} & \bigl(\colim\,\Gpr\bigr)\sbinom{Bd}{B\uc}}
\end{equation}
is commutative, in which:
\begin{itemize}
\item In the top row, the products are unordered.
\item $\gamma^{\O}_{(T,\sigma)}$ is the $\G$-operad structure morphism of $\O$ for the $\colorc$-colored $\G$-tree $(T,\sigma)$.
\item $BT$ is the $\colord$-colored planar tree obtained from $T$ by composing the $\colorc$-coloring of $T$ with the map $B : \colorc \to\colord$.
\item $\gamma^{\colim\,\Gpr}_{(BT,\sigma)}$ is the $\G$-operad structure morphism of $\colim\,\Gpr$ for the $\colord$-colored $\G$-tree $(BT,\sigma)$.
\end{itemize}

For each $t\in\Vt(T)$, pick an element $x_t\in\O\inoutt$.  Suppose \[y=\gamma^{\O}_{(T,\sigma)}\{x_t\}_{t\in\Vt(T)} \in \O\duc,\] and 
\begin{equation}\label{phi_xt}
\phi_{\{x_t\}_{t\in\Vt(T)}} : \G^{(T,\sigma)}\to\O \in\gopset
\end{equation} 
is the morphism that sends each generator $t\in\Vt(T)$ to $x_t\in\O\inoutt$.  Consider the morphisms
\begin{equation}\label{varphi_t}
\nicexy@C+2cm{\bigl(\Cor_{(\uc;d)},\id_n\bigr) \ar[r]^-{\varphi}_-{\bigl((\Cor_{(\uc;d)},\id_n),v; (T,\sigma)\bigr)} & (T,\sigma)&
\bigl(\Cor_{t},\id_t\bigr) \ar[l]_-{\varphi_t}^-{\bigl((T,\sigma),t; (\Cor_t,\id_t)\bigr)}}
\end{equation}
in $\Treecatg$ with 
\begin{itemize}
\item $v$ the unique vertex in the $\duc$-corolla,
\item $\Cor_t=\Cor_{|\inp(t)|}$, and
\item $\id_t=\id_{|\inp(t)|}$.
\end{itemize}  
There is a commutative diagram \[\nicexy@C+.5cm{\G^{(\Cor_{(\uc;d)},\id_n)} \ar[r]^-{\phi_y} \ar[d]_-{\G^{\varphi}} & \O\ar@{=}[d]\\ \G^{(T,\sigma)} \ar[r]^-{\phi_{\{x_t\}_{t\in\Vt(T)}}} & \O\\
\G^{(\Cor_t,\id_t)} \ar[u]^-{\G^{\varphi_t}} \ar[r]^-{\phi_{x_t}} & \O \ar@{=}[u]}\]in $\gopset$.   This implies that there is a commutative diagram 
\begin{equation}\label{gpr-varphi-omega}
\nicexy{\Gpr\Bigl((\Cor_{(\uc;d)},\id_n); \phi_y\Bigr)= \G^{(\Cor_{(\uc;d)},\id_n)} \ar[r]^-{\omega} \ar[d]_-{\G^{\varphi}} & \colim\,\Gpr \ar@{=}[d]\\ 
\Gpr\Bigl((T,\sigma); \phi_{\{x_t\}_{t\in\Vt(T)}}\Bigr)= \G^{(T,\sigma)} \ar[r]^-{\omega} & \colim\,\Gpr\ar@{=}[d]\\
\Gpr\Bigl((\Cor_t,\id_t);\phi_{x_t}\Bigr)= \G^{(\Cor_t,\id_t)} \ar[r]^-{\omega} \ar[u]^-{\G^{\varphi_t}} & \colim\,\Gpr}
\end{equation}
in $\gopset$.  

Writing $v_t$ for the unique vertex in $\Cor_t$, it follows from the previous diagram that there are equalities:
\[\begin{split}
B\Bigl(\gamma^{\O}_{(T,\sigma)}\{x_t\}_{t\in\Vt(T)}\Bigr)
&= B(y)= \omega(v)=\omega\G^{\varphi}(v)\\
&=\omega\Bigl(\gamma^{\G^{(T,\sigma)}}_{(T,\sigma)}\{t\}_{t\in\Vt(T)}\Bigr)\\
&=\gamma^{\colim\,\Gpr}_{(BT,\sigma)} \Bigl(\{\omega(t)\}_{t\in\Vt(T)}\Bigr)\\
&=\gamma^{\colim\,\Gpr}_{(BT,\sigma)} \Bigl(\{\omega(v_t)\}_{t\in\Vt(T)}\Bigr)\\
&=\gamma^{\colim\,\Gpr}_{(BT,\sigma)} \Bigl(\{B(x_t)\}_{t\in\Vt(T)}\Bigr).
\end{split}\]
This shows that the diagram \eqref{B-lambda-compatible} is commutative.
\end{proof}

The following result says that each $\G$-operad can be built from $\G$-operads that are freely generated by $\G$-trees.

\begin{theorem}\label{goperad-colimit}
Suppose $\O$ is a $\colorc$-colored $\G$-operad in $\Set$.  Then the morphisms
\[\nicexy{\colimover{\Govero}\, \Gpr \ar@<2pt>[r]^-{A} & (\colorc,\O) \ar@<2pt>[l]^-{B}}\] in $\gopset$ are inverses of each other.
\end{theorem}

\begin{proof}
The equality \[AB=\Id_{(\colorc,\O)}\] holds by the definitions of the morphisms $A$ and $B$.

To show that the composite $BA$ is also the identity morphism, suppose $(T,\sigma)$ is a $\colorc$-colored $\G$-tree.  Pick arbitrary elements $x_t\in\O\inoutt$ for each $t\in\Vt(T)$, and suppose \[\phi_{\{x_t\}_{t\in\Vt(T)}} : \G^{(T,\sigma)}\to\O\in\gopset\] is the morphism in \eqref{phi_xt}.  For each vertex $t$ in $T$, suppose \[\nicexy@C+2.5cm{\bigl(\Cor_t,\id_t\bigr) \ar[r]^-{\varphi_t\,=\, \bigl((T,\sigma),t; (\Cor_t,\id_t)\bigr)} & (T,\sigma)} \in \Treecatg\] is the morphism in \eqref{varphi_t}.  There is a commutative diagram
\[\nicexy{\Gpr\Bigl((\Cor_t,\id_t); \phi_{x_t}\Bigr) 
\ar `l[dd] `[dd]_(.7){\omega} [dd]
\ar[d]_-{\G^{\varphi_t}} \ar@{=}[rr] && \G^{(\Cor_t,\id_t)} \ar[dd]^-{\omega} \ar[ddl]|(.4){\phi_{x_t}}\\
\Gpr\Bigl((T,\sigma); \phi_{\{x_t\}_{t\in\Vt(T)}}\Bigr) \ar[d]_-{\omega} \ar@{=}[r] & \G^{(T,\sigma)} \ar[d]_-{\phi_{\{x_t\}_{t\in\Vt(T)}}} &\\
\colim\,\Gpr \ar[r]^(.4){A} & \O \ar[r]^-{B} & \colim\,\Gpr}\] in $\gopset$.  So the composite \[\nicexy{\Gpr\Bigl((T,\sigma); \phi_{\{x_t\}_{t\in\Vt(T)}}\Bigr) \ar[r]^-{BA\omega} & \colim\,\Gpr}\] applied to each generator $t\in\G^{(T,\sigma)}$ is simply $\omega(t)$.  Since the $\colorc$-colored $\G$-tree $(T,\sigma)$ and the elements $\{x_t\}\in \prod_{t\in\Vt(T)} \O\inoutt$ are arbitrary, this implies that $BA$ is the identity morphism on $\colim\,\Gpr$.
\end{proof}

\section{Nerve is Fully Faithful}\label{sec:gnerve-fully-faithful}

The purpose of this section is to observe that the $\G$-nerve functor is full and faithful. 

\begin{theorem}\label{nerveg-fully-faithful}
The counit\index{nerve!full and faithful} \[\eta : \Realg\Nerveg \iso\Id_{\gopset}\] of the $\G$-realization-nerve adjunction is a natural isomorphism.  In particular, the $\G$-nerve $\Nerveg$ is full and faithful.
\end{theorem}

\begin{proof}
Suppose $\O$ is a $\colorc$-colored $\G$-operad in $\Set$.  For each $\colorc$-colored $\G$-tree $(T,\sigma)$, a morphism \[\Yonedag(T,\sigma)\to \Nerveg(\colorc,\O)\in \Settotreecatgop\] is uniquely determined by a morphism \[\bigl(\Ed(T),\G^{(T,\sigma)}\bigr) \to (\colorc,\O)\in \gopset\] by Yoneda Lemma.  

The $\G$-nerve is a colimit of representables.  More precisely, the natural morphism 
\begin{equation}\label{nerve-colim-representables}
\nicexy{\colimover{\bigl((T,\sigma); \phi\bigr)\in\Govero} \Yonedag(T,\sigma) \ar[r]^-{\cong} & \Nerveg(\colorc,\O)}\in \Settotreecatgop
\end{equation}
is an isomorphism.  Applying the $\G$-realization functor, we obtain the natural isomorphisms
\[\nicexy{\Realg \Nerveg(\colorc,\O) & (\colorc,\O)\\
\Realg\Bigl(\colimover{\bigl((T,\sigma); \phi\bigr)\in\Govero} \Yonedag(T,\sigma)\Bigr) \ar[u]^-{\cong} & \colimover{\Govero}\,\Gpr \ar[u]_-{A}^-{\cong}\\ 
\colimover{\bigl((T,\sigma); \phi\bigr)\in\Govero} \Realg\Yonedag(T,\sigma) \ar[u]^-{\cong} \ar[r]^-{\cong} & \colimover{\bigl((T,\sigma); \phi\bigr)\in\Govero} \G^{(T,\sigma)}. \ar@{=}[u]}\]
In the above diagram:
\begin{itemize}
\item The lower left isomorphism follows from the fact in Lemma \ref{nerve-real-adjunction} that $\Realg$ is a left adjoint, which preserves colimits.  
\item The bottom horizontal isomorphism uses the fact that the Yoneda embedding is full and faithful, so the Yoneda extension $\Realg$ actually extends $\G^-$ up to a natural isomorphism by Corollary X.3.3 on page 239 in \cite{maclane}. 
\item The isomorphism $A$ is due to Theorem \ref{goperad-colimit}.  
\end{itemize}
This proves the first assertion.  

The second assertion--that $\Nerveg$ is full and faithful--is a formal consequence of the first assertion.  See, e.g., Proposition 1.3 on page 7 in \cite{gz} or \cite{maclane} IV.3 Theorem 1.
\end{proof}

\begin{example}[Planar and Symmetric Dendroidal Nerves]\label{ex:dendroidal-nerve-fullyfaithful}
When $\G$ is the\index{nerve!dendroidal}\index{dendroidal nerve} symmetric group operad $\S$, the symmetric tree category $\Treecats$ is equivalent to the Moerdijk-Weiss symmetric dendroidal category $\Omega$.  In the setting of $\Omega$, it is noted in Example 3.1.4 in \cite{mt} that the dendroidal nerve functor is full and faithful, which corresponds to the $\G=\S$ case of the second assertion in Theorem \ref{nerveg-fully-faithful}.  When $\G$ is the planar group operad $\P$, we have the planar analogue of this statement, which applies to the planar dendroidal category $\Omega_p \cong \Treecatp$.\dqed
\end{example}

\begin{example}[Braided Nerve]\label{ex:braided-nerve-fullyfaithful}
For\index{nerve!braided}\index{braided nerve} the braid group operad $\B$ in Definition \ref{def:braid-group-operad} and the pure braid group operad $\PB$ in Definition \ref{def:pure-braid-group-operad}, the braided nerve functor and the pure braided nerve functor,\[\nicexy{\Settotreecatbop & \bopset \ar[l]_-{\Nerveb}}\andspace \nicexy{\Settotreecatpbop & \pbopset \ar[l]_-{\Nervepb}},\] are full and faithful by Theorem \ref{nerveg-fully-faithful}.\dqed
\end{example}

\begin{example}[Ribbon Nerve]\label{ex:ribbon-nerve-fullyfaithful}
For\index{nerve!ribbon}\index{ribbon nerve} the ribbon group operad $\R$ in Definition \ref{def:ribbon-group-operad} and the pure ribbon group operad $\PR$ in Definition \ref{def:pure-ribbon-group-operad}, the ribbon nerve functor and the pure ribbon nerve functor \[\nicexy{\Settotreecatrop & \ropset \ar[l]_-{\Nerver}} \andspace \nicexy{\Settotreecatprop & \propset \ar[l]_-{\Nervepr}},\] are full and faithful by Theorem \ref{nerveg-fully-faithful}.\dqed
\end{example}

\begin{example}[Cactus Nerve]\label{ex:cactus-nerve-fullyfaithful}
For\index{nerve!cactus}\index{cactus nerve} the cactus group operad $\Cac$ in Definition \ref{def:cactus-group-operad} and the pure cactus group operad $\PCac$ in Definition \ref{def:pure-cactus-group-operad}, the cactus nerve functor and the pure cactus nerve functor \[\nicexy{\Settotreecatcacop & \cacopset \ar[l]_-{\Nervecac}} \andspace \nicexy{\Settotreecatpcacop & \pcacopset \ar[l]_-{\Nervepcac}},\] are full and faithful by Theorem \ref{nerveg-fully-faithful}.\dqed
\end{example}

\section{Symmetric Monoidal Structure on Presheaf Category}\label{sec:gtreeset-monoidal}

In this section, we observe that the category $\Settotreecatgop$ of $\Treecatg$-presheaves admits a symmetric monoidal closed structure that is closely related to the $\G$-tensor product of $\G$-operads.  In the next section, we will show that the $\G$-nerve functor is symmetric monoidal, and the $\G$-realization functor is strong symmetric monoidal.

\begin{definition}\label{def:treecatg-presheaf-monoidal}
Suppose $\G$ is an action operad, and $X,Y\in \Settotreecatgop$.
\begin{description}
\item[Monoidal Product] Define the object\index{monoidal product!of presheaves@of $\Treecatg$-presheaves} \[X\tensorsub{\G} Y = \colimover{(x;y)\,\in\, X(T,\sigma)\times Y(T',\sigma')} \Nerveg\bigl(\G^{(T,\sigma)}\tensorg \G^{(T',\sigma')}\bigr) \in \Settotreecatgop\] in which:
\begin{itemize}
\item $\tensorg$ is the symmetric monoidal product in $\gopset$ in Theorem \ref{goperad-symmetric-monoidal}.
\item $\G^{(T,\sigma)}$ is the $\Ed(T)$-colored $\G$-operad freely generated by the set of vertices in $T$ in Definition \ref{def:gtree-goperad}.
\item $\Nerveg : \gopset \to \Settotreecatgop$ is the $\G$-nerve functor in Definition \ref{realg-nerveg}.
\item The colimit is indexed by the category in which an object is a pair of morphisms \[\Bigl(\narrowxy{\Yonedag(T,\sigma)\ar[r]^-{x} & X}; \narrowxy{\Yonedag(T',\sigma')\ar[r]^-{y} & Y}\Bigr)\] in $\Settotreecatgop$ with $(T,\sigma),(T',\sigma')\in\Treecatg$, i.e., a pair of elements \[(x;y)\in X(T,\sigma)\times Y(T',\sigma').\]  A morphism from such an element $(x;y)$ to an element $(x_1;y_1)\in X(T_1,\sigma_1)\times Y(T'_1,\sigma'_1)$ is a pair of morphisms \[\Bigl(\narrowxy{(T,\sigma)\ar[r]^-{\phi} & (T_1,\sigma_1)}; \narrowxy{(T',\sigma')\ar[r]^-{\phi_1} & (T'_1,\sigma'_1)}\Bigr)\] in $\Treecatg$ such that \[X(\phi^{\op})(x_1)=x \andspace Y(\phi_1^{\op})(y_1)=y.\]
\end{itemize}
\item[Internal Hom] Define the object\label{not:homsubg} \[\Homsub{\G}(X,Y)\in\Settotreecatgop\] by \[\Homsub{\G}(X,Y)(T,\sigma) = \Settotreecatgop\Bigl(X\tensorsub{\G}\Yonedag(T,\sigma), Y\Bigr)\] for $(T,\sigma)\in\Treecatg$, in which $X\tensorsub{\G}\Yoneda(T,\sigma)$ is defined in the previous item.
\end{description}
\end{definition}

The next observation contains sample calculation involving representable $\Treecatg$-presheaves.

\begin{proposition}\label{representable-tensor}
Suppose $(T,\sigma)$ and $(T',\sigma')$ are in $\Treecatg$, and $\O\in\gopset$.  Then there are natural isomorphisms 
\[\begin{split}
\Yonedag(T,\sigma) \tensorsub{\G} \Yonedag(T',\sigma') &\cong \Nerveg\bigl(\G^{(T,\sigma)}\tensorg \G^{(T',\sigma')}\bigr),\\
\Nerveg(\O)\tensorsub{\G} \Yonedag(T',\sigma') &\cong \Nerveg\bigl(\O\tensorg \G^{(T',\sigma')}\bigr),\\
\Homsub{\G}\bigl(\Yonedag(T,\sigma),\Yonedag(T',\sigma')\bigr) &\cong
\gopset\bigl(\G^{(T,\sigma)}\tensorg \G^?,\G^{(T',\sigma')}\bigr)\end{split}\]
in $\Settotreecatgop$.
\end{proposition}

\begin{proof}
The first natural isomorphism is an immediate consequence of Yoneda Lemma.  The second isomorphism follows from the first one and the fact that an element in $\Nerveg(\O)(T,\sigma)$ is by definition a morphism $\G^{(T,\sigma)}\to\O$.

For the last natural isomorphism, suppose $(V,\tau)$ is in $\Treecatg$.  Then we compute as follows:
\[\begin{split}
&\Homsub{\G}\bigl(\Yonedag(T,\sigma),\Yonedag(T',\sigma')\bigr) (V,\tau)\\
&= \Settotreecatgop\Bigl(\Yonedag(T,\sigma) \tensorsub{\G} \Yonedag(V,\tau), \Yonedag(T',\sigma')\Bigr)\\
&\cong \Settotreecatgop\Bigl(\Nerveg\bigl(\G^{(T,\sigma)} \tensorg \G^{(V,\tau)}\bigr), \Nerveg\G^{(T',\sigma')}\Bigr)\\
&\cong \gopset\bigl(\G^{(T,\sigma)}\tensorg \G^{(V,\tau)},\G^{(T',\sigma')}\bigr).
\end{split}\]
The first isomorphism follows from the first natural isomorphism and Proposition \ref{nerve-factors-yoneda}.  The last isomorphism follows from Theorem \ref{nerveg-fully-faithful}, which says that the $\G$-nerve is full and faithful.
\end{proof}

\begin{theorem}\label{thm:presheaf-monoidal}
The tuple \[\Bigl(\Settotreecatgop, \tensorsub{\G}, \Yonedag(\uparrow,\id_1), \Homsub{\G}\Bigr)\] is a symmetric monoidal closed category.\index{G-tree category@$\G$-tree category!presheaves} 
\end{theorem}

\begin{proof}
It follows from the first two natural isomorphisms in Proposition \ref{representable-tensor} and the symmetry of $\tensorg$ that \[\Bigl(\Settotreecatgop, \tensorsub{\G}, \Yonedag(\uparrow,\id_1)\Bigr)\] is a symmetric monoidal category.  To see that it is closed with $\Homsub{\G}$, we must show that for $X,Y,Z \in \Settotreecatgop$, there is a natural bijection 
\begin{equation}\label{presheaf-tensor-hom}
\Settotreecatgop\bigl(X\tensorsub{\G} Y,Z\bigr) \cong \Settotreecatgop\bigl(X,\Homsub{\G}(Y,Z)\bigr).
\end{equation}
We may assume that $X$ is a representable $\Treecatg$-presheaf, i.e., $X=\Yonedag(T,\sigma)$ for some $\G$-tree $(T,\sigma)$, since every $\Treecatg$-presheaf is canonically a colimit of representable $\Treecatg$-presheaves and since $\tensorsub{\G}$ commutes with colimits on each side separately.  There are natural bijections
\[\begin{split}
\Settotreecatgop\bigl(\Yonedag(T,\sigma)\tensorsub{\G} Y, Z\bigr)
&\cong \Settotreecatgop\bigl(Y\tensorsub{\G} \Yonedag(T,\sigma),Z\bigr)\\
&= \Homsub{\G}(Y,Z)(T,\sigma)\\
&\cong \Settotreecatgop\bigl(\Yonedag(T,\sigma), \Homsub{\G}(Y,Z)\bigr),
\end{split}\]
which prove the required natural bijection.
\end{proof}

The following observation is an extension of the natural bijection in \eqref{presheaf-tensor-hom} to a natural isomorphism of $\Treecatg$-presheaves.

\begin{corollary}
For $X,Y,Z \in \Settotreecatgop$, there is a natural isomorphism 
\[\Homsub{\G}\bigl(X\tensorsub{\G} Y,Z\bigr) \cong \Homsub{\G}\bigl(X,\Homsub{\G}(Y,Z)\bigr) \in \Settotreecatgop.\]
\end{corollary}

\begin{proof}
For each $\G$-tree $(T,\sigma)$, there are natural bijections
\[\begin{split}
\Homsub{\G}\bigl(X\tensorsub{\G} Y,Z\bigr)(T,\sigma)
&= \Settotreecatgop\Bigl((X\tensorsub{\G} Y)\tensorsub{\G} \Yonedag(T,\sigma),Z\Bigr)\\
&\cong \Settotreecatgop\Bigl(\bigl(X\tensorsub{\G} \Yonedag(T,\sigma)\bigr)\tensorsub{\G} Y,Z\Bigr)\\
&\cong \Settotreecatgop\Bigl(X\tensorsub{\G} \Yonedag(T,\sigma), \Homsub{\G}(Y,Z)\Bigr)\\
&=\Homsub{\G}\bigl(X,\Homsub{\G}(Y,Z)\bigr)(T,\sigma).\end{split}\]
The first bijection follows from the symmetry of $\tensorsub{\G}$.  The second bijection follows from \eqref{presheaf-tensor-hom}.
\end{proof}

\section{Nerve is Symmetric Monoidal}
\label{sec:gnerve-symmetric-monoidal}

The purpose of this section is to observe that the $\G$-realization functor is \emph{strong} symmetric monoidal.  As a result, the $\G$-nerve functor is symmetric monoidal.

Recall from Theorem \ref{goperad-symmetric-monoidal} that $\gopset$ is a symmetric monoidal category.  Recall that a \emph{strong monoidal functor} is a monoidal functor \[(F,F_2,F_0) : (\C,\otimes^{\C},\tensorunit^{\C}) \to (\D,\otimes^{\D},\tensorunit^{\D})\] in which all the structure morphisms \[F_2(X,Y) : FX \otimes^{\D} FY \to F\bigl(X\otimes^{\C}Y\bigr) \andspace F_0 : \tensorunit^{\D}\to F\tensorunit^{\C}\] are isomorphisms.

\begin{theorem}\label{realization-strong-monoidal}
For each action operad $\G$, the $\G$-realization\index{realization!strong symmetric monoidal} functor \[\Realg : \Bigl(\Settotreecatgop, \tensorsub{\G}, \Yonedag(\uparrow,\id_1)\Bigr) \to \bigl(\gopset,\tensorg,\I\bigr)\] is strong symmetric monoidal.  In particular, the $\G$-nerve\index{nerve!symmetric monoidal} functor \[\Nerveg : \bigl(\gopset,\tensorg,\I\bigr)\to \Bigl(\Settotreecatgop, \tensorsub{\G}, \Yonedag(\uparrow,\id_1)\Bigr)\] is symmetric monoidal.
\end{theorem}

\begin{proof}
Since the $\G$-realization functor $\Realg$ is an extension of $\G^-$ up to a natural isomorphism, we have the isomorphism \[\Realg\Yonedag(\uparrow,\id_1)\cong \G^{(\uparrow,\id_1)} = \I^{\G}.\]  

It remains to show that there is a natural isomorphism 
\begin{equation}\label{real-xtensory}
\Realg\bigl(X\tensorsub{\G} Y\bigr) \cong (\Realg X) \tensorg (\Realg Y)\in \gopset
\end{equation}
for $X,Y\in\Settotreecatgop$.  We know that $\Realg$, being the left adjoint of $\Nerveg$, commutes with colimits, and both $\tensorg$ and $\tensorsub{\G}$ commute with colimits on each side.  Therefore, we may assume that $X$ and $Y$ are representable $\Treecatg$-presheaves, i.e., \[X=\Yonedag(T,\sigma) \andspace Y=\Yonedag(T',\sigma')\] for some $(T,\sigma)$, $(T',\sigma')\in\Treecatg$.  Then there are natural isomorphisms
\[\begin{split}
&\Realg\bigl(\Yonedag(T,\sigma) \tensorsub{\G}\Yonedag(T',\sigma')\bigr)\\
&\cong \Realg\Nerveg\bigl(\G^{(T,\sigma)} \tensorg \G^{(T',\sigma')}\bigr)\\
&\cong \G^{(T,\sigma)} \tensorg \G^{(T',\sigma')}\\
&\cong \bigl(\Realg\Yonedag(T,\sigma)\bigr) \tensorg \bigl(\Realg\Yonedag(T',\sigma')\bigr)
\end{split}\]
in $\gopset$.  The first isomorphism is from Proposition \ref{representable-tensor}.  The second isomorphism follows from Theorem \ref{nerveg-fully-faithful}.  The last isomorphism follows from the natural isomorphism $\Realg\Yonedag\cong \G^-$.

The assertion that $\Nerveg$ is a symmetric monoidal functor is a purely formal consequence of the fact that its left adjoint $\Realg$ is strong symmetric monoidal.
\end{proof}

\begin{corollary}\label{real-nerve-tensor-nerve}
For $\G$-operads $(\colorc,\O)$ and $(\colord,\P)$ in $\Set$, there is a natural isomorphism \[\Realg\Bigl(\Nerveg(\colorc,\O) \tensorsub{\G} \Nerveg(\colord,\P)\Bigr) \cong (\colorc,\O)\tensorg (\colord,\P).\]
\end{corollary}

\begin{proof}
This follows from the natural isomorphism $\Realg\Nerveg\cong\Id$ in Theorem \ref{nerveg-fully-faithful} and from \eqref{real-xtensory} by substituting in $\Nerveg(\colorc,\O)$ and $\Nerveg(\colord,\P)$ for $X$ and $Y$, respectively.
\end{proof}

\begin{example}
When $\G$ is the symmetric group operad $\S$, most of the results in this section can be found in Section 4.2 in \cite{mt} for the symmetric dendroidal category $\Omega$, which is equivalent to the symmetric tree category $\Treecats$.\dqed
\end{example}

\section{Change of Action Operads}\label{sec:nerve-change-g}

In this section, we study how the categories and the functors in the diagram \eqref{realg-nerveg}
\[\nicexy@C+.4cm{\Treecatg \ar[d]_-{\Yonedag} \ar@{=}[r] & \Treecatg \ar[d]^-{\G^-}\\
\Settotreecatgop \ar@<2pt>[r]^-{\Realg} & \gopset \ar@<2pt>[l]^-{\Nerveg}}\]
behave with respect to a change of the action operad $\G$.  Throughout this section, suppose $\varphi : \Gone \to\Gtwo$ is a morphism of action operads as in Definition \ref{def:augmented-group-operad-morphism}.

First we observe that the Yoneda embedding is compatible with a change of action operads.

\begin{proposition}\label{yoneda-varphi}\index{change of action operads!Yoneda embedding}
The diagram\index{Yoneda embedding!change of action operads} \[\nicexy{\Treecatgone \ar[d]_-{\Yonedagone} \ar[r]^-{\Treecat^{\varphi}} & \Treecatgtwo \ar[d]^-{\Yonedagtwo}\\ \Settotreecatgoneop \ar[r]^-{\varphi_!} & \Settotreecatgtwoop}\]
is commutative up to a natural isomorphism.
\end{proposition}

\begin{proof}
Suppose $(T,\sigma)$ is a $\Gone$-tree and $Z \in \Settotreecatgtwoop$.  There are natural bijections 
\[\begin{split}
&\Settotreecatgtwoop\Bigl(\varphi_!\Yonedagone(T,\sigma),Z\Bigr)\\
&\cong \Settotreecatgoneop\Bigl(\Yonedagone(T,\sigma),\varphi^*Z\Bigr)\\
&\cong (\varphi^*Z)(T,\sigma)\\
&= Z\bigl(\Treecat^{\varphi}(T,\sigma)\bigr)\\
&\cong \Settotreecatgtwoop\Bigl(\Yonedagtwo\Treecat^{\varphi}(T,\sigma), Z\Bigr)
\end{split}\]
in which the first isomorphism follows from the adjunction $\varphi_! \dashv\varphi^*$.  The second and the third isomorphisms are by Yoneda Lemma.   Since $(T,\sigma)$ and $Z$ are arbitrary, the result follows.
\end{proof}

Consider the diagram
\begin{equation}\label{nervegone-nervegtwo}
\nicexy@C+.5cm{\Settotreecatgoneop \ar@<2pt>[r]^-{\varphi_!} \ar@<-2pt>[d]_-{\Realgone} & \Settotreecatgtwoop \ar@<-2pt>[d]_-{\Realgtwo} \ar@<2pt>[l]^-{\varphi^*}\\
\goneopset \ar@<-2pt>[u]_-{\Nervegone} \ar@<2pt>[r]^-{\varphi_!} & \gtwoopset \ar@<2pt>[l]^-{\varphi^*} \ar@<-2pt>[u]_-{\Nervegtwo}}
\end{equation}
in which:
\begin{itemize}
\item The top functor $\varphi^*$ is the pullback induced by \[(\Treecat^{\varphi})^{\op} : \Treecatgoneop \to \Treecatgtwoop\] with $\Treecat^{\varphi}$ the functor in Proposition \ref{treecatg-functor}(1).
\item The top $\varphi_!$ is the left adjoint of $\varphi^*$ given by a left Kan extension.
\item The bottom adjunction $(\varphi_!,\varphi^*)$ is the one in Theorem \ref{goneopm-gtwoopm}.
\item $(\Realgone,\Nervegone)$ is the $\Gone$-realization-nerve adjunction, and $(\Realgtwo,\Nervegtwo)$ is the $\Gtwo$-realization-nerve adjunction.
\end{itemize}
Next we observe that the $\G$-realization functor and the $\G$-nerve functor are well-behaved with respect to a morphism of action operads.

\begin{theorem}\label{gnerve-changeofg}\index{change of action operads!G-realization-nerve@$\G$-realization-nerve}
Each morphism\index{nerve!change of action operads}\index{realization!change of action operads} $\varphi : \Gone\to\Gtwo$ of action operads yields  natural isomorphisms \[\Nervegone\varphi^*\cong\varphi^*\Nervegtwo \andspace \varphi_!\Realgone\cong\Realgtwo\varphi_!\] in the diagram \eqref{nervegone-nervegtwo}.
\end{theorem}

\begin{proof}
To prove the first natural isomorphism, suppose $\O$ is a $\colorc$-colored $\Gtwo$-operad in $\Set$, and $(T,\sigma)$ is a $\Gone$-tree.  By definition we have the set \[\Bigl(\Nervegone \varphi^*(\colorc,\O)\Bigr)(T,\sigma) = \goneopset\Bigl((\Gone)^{(T,\sigma)}, (\colorc,\varphi^*\O)\Bigr).\]  Recall that:
\begin{itemize}\item $(\Gone)^{(T,\sigma)}$ is the free $\Gone$-operad generated by the vertices in $T$.
\item The $\colorc$-colored $\Gone$-operad $\varphi^*\O$ has the same underlying entries as $\O$.
\end{itemize}
So a $\Gone$-operad morphism \[\phi : (\Gone)^{(T,\sigma)} \to (\colorc,\varphi^*\O)\] consists of
\begin{itemize} \item a map $\phi : \Ed(T)\to\colorc$ and
\item a choice of an element $\phi(t)\in\O\inoutphit$ for each $t\in\Vt(T)$.
\end{itemize}
For the other composite, we have the set \[\Bigl(\varphi^*\Nervegtwo(\colorc,\O)\Bigr)(T,\sigma)= \gtwoopset\Bigl((\Gtwo)^{(T,\varphi\sigma)}, (\colorc,\O)\Bigr).\]  Each morphism \[(\Gtwo)^{(T,\varphi\sigma)} \to (\colorc,\O)\] of $\Gtwo$-operads admits the same description as a morphism $\phi$ of $\Gone$-operads as above.  This proves the first natural isomorphism.

The second natural isomorphism follows from the first one by the uniqueness of left adjoints, since $\varphi_!\Realgone$ and $\Realgtwo\varphi_!$ are the left adjoints of $\Nervegone\varphi^*$ and $\varphi^*\Nervegtwo$, respectively.
\end{proof}

\section{Comparing Symmetric Monoidal Structures}\label{sec:presheaf-functor-monoidal}

In Theorem \ref{goperad-monoidal-functor} and Corollary \ref{goperad-comonoidal}, we showed that in the adjunction \[\nicexy{\Bigl(\goneopset,\tensorgone,\I\Bigr) \ar@<2pt>[r]^-{\varphi_!} & \Bigl(\gtwoopset,\tensorgtwo,\I\Bigr) \ar@<2pt>[l]^-{\varphi^*}},\] the left adjoint $\varphi_!$ is a symmetric comonoidal functor, and the right adjoint $\varphi^*$ is a symmetric monoidal functor.  The purpose of this section is to observe that the same is true for the top adjunction in the diagram \eqref{nervegone-nervegtwo}.

\begin{theorem}\label{treecatg-adjunction-comonoidal}\index{change of action operads!presheaf category}
In the adjunction \[\nicexy@C+.5cm{\Settotreecatgoneop \ar@<2pt>[r]^-{\varphi_!}  & \Settotreecatgtwoop \ar@<2pt>[l]^-{\varphi^*}}\] with the symmetric monoidal structures in Theorem \ref{thm:presheaf-monoidal}:
\begin{enumerate}
\item The left adjoint $\varphi_!$ is a symmetric comonoidal functor.
\item The right adjoint $\varphi^*$ is a symmetric monoidal functor.
\end{enumerate}
\end{theorem}

\begin{proof}
To see that $\varphi_!$ is a comonoidal functor, first note that for each $\Gone$-tree $(T,\sigma)$, there is a morphism \[\nicexy{\Yonedagone\bigl(\uparrow,\id_1^{\Gone}\bigr)(T,\sigma) \ar[r] & \Bigl(\varphi^*\Yonedagtwo\bigl(\uparrow, \id_1^{\Gtwo}\bigr)\Bigr)(T,\sigma)\\ & \Yonedagtwo\bigl(\uparrow,\id_1^{\Gtwo}\bigr)(T,\varphi\sigma) \ar@{=}[u] \ar@{=}[d]\\ 
\Treecatgone\Bigl((T,\sigma), (\uparrow,\id_1^{\Gone})\Bigr) \ar@{=}[uu] \ar[r]^-{\Treecat^{\varphi}} & \Treecatgtwo\Bigl((T,\varphi\sigma), (\uparrow,\id_1^{\Gtwo})\Bigr)}\] in which $\id_1^{\G^i} \in\G^i(1)$ is the multiplicative unit for $i=1,2$.  Since $(T,\sigma)$ is arbitrary, this defines a morphism 
\begin{equation}\label{varphistar-unit}
\nicexy{\Yonedagone\bigl(\uparrow,\id_1^{\Gone}\bigr)\ar[r]^-{\varphi^*_0} & \varphi^*\Yonedagtwo\bigl(\uparrow, \id_1^{\Gtwo}\bigr)}\in \Settotreecatgoneop.
\end{equation}
Its adjoint is the comonoidal structure morphism \[\nicexy{\varphi_!\Yonedagone\bigl(\uparrow,\id_1^{\Gone}\bigr) \ar[r]^-{\varphi_!^0} & \Yonedagtwo\bigl(\uparrow,\id_2^{\Gtwo}\bigr)} \in \Settotreecatgtwoop\] on the monoidal units.

To define the other comonoidal structure morphism \[\nicexy{\varphi_!\bigl(X\tensorsub{\Gone}Y\bigr) \ar[r]^-{\varphi_!^2} & \varphi_!X \tensorsub{\Gtwo} \varphi_!Y}\in \Settotreecatgtwoop,\] note that $\varphi_!$, being a left adjoint, commutes with colimits and that each $\tensorsub{\G^i}$ commutes with colimits on each side.  Since each $\Treecatgone$-presheaf is canonically a colimit of representable $\Treecatgone$-presheaves, we only need to define the comonoidal structure morphism $\varphi_!^2$ when \[X=\Yonedagone(T,\sigma) \andspace Y=\Yonedagone(T',\sigma')\] for some $\Gone$-trees $(T,\sigma)$ and $(T',\sigma')$.  In this case, $\varphi_!^2$ is the following composite of natural morphisms:
\begin{small}
\[\narrowxy{\varphi_!\Bigl(\Yonedagone(T,\sigma) \tensorsub{\Gone} \Yonedagone(T',\sigma')\Bigr) \ar[r]^-{\varphi_!^2} \ar[d]_-{\varphi_!(\cong)} & \bigl(\varphi_!\Yonedagone(T,\sigma)\bigr) \tensorsub{\Gtwo} \bigl(\varphi_!\Yonedagone(T',\sigma')\bigr)\\
\varphi_!\Nervegone\Bigl((\Gone)^{(T,\sigma)} \tensorgone (\Gone)^{(T',\sigma')}\Bigr) \ar[d]_-{\varphi_!\Nervegone(\omega,\omega)} & 
\Yonedagtwo(T,\varphi\sigma) \tensorsub{\Gtwo} \Yonedagtwo(T',\varphi\sigma') \ar[u]_-{\cong}\\
\varphi_!\Nervegone\Bigl(\varphi^*(\Gtwo)^{(T,\varphi\sigma)}\tensorgone \varphi^*(\Gtwo)^{(T',\varphi\sigma')}\Bigr) \ar[d]_-{\varphi_!\Nervegone(\varphi^*_2)} & 
\Nervegtwo\Bigl((\Gtwo)^{(T,\varphi\sigma)} \tensorgtwo (\Gtwo)^{(T',\varphi\sigma')}\Bigr) \ar[u]_-{\cong}\\
\varphi_!\Nervegone\varphi^*\Bigl((\Gtwo)^{(T,\varphi\sigma)} \tensorgtwo (\Gtwo)^{(T,\varphi\sigma)}\Bigr) \ar[r]^-{\varphi_!(\cong)} & 
\varphi_!\varphi^*\Nervegtwo\Bigl((\Gtwo)^{(T,\varphi\sigma)} \tensorgtwo (\Gtwo)^{(T,\varphi\sigma)}\Bigr) \ar[u]_-{\eta}}\]
\end{small}

In the above diagram:
\begin{itemize}
\item The upper-left isomorphism and the middle-right isomorphism are from Proposition \ref{representable-tensor}.
\item The upper-right isomorphism \[\varphi_!\Yonedagone \cong \Yonedagtwo\Treecat^{\varphi}\] is from Proposition \ref{yoneda-varphi}.
\item The bottom horizontal isomorphism \[\Nervegone\varphi^*\cong \varphi^*\Nervegtwo\] is from Theorem \ref{gnerve-changeofg}.
\item In the lower-left morphism, $\varphi^*_2$ is the monoidal structure morphism of $\varphi^*$ in Theorem \ref{goperad-monoidal-functor}.
\item The lower-right morphism \[\eta : \varphi_!\varphi^*\to\Id\] is the counit of the adjunction $\varphi_! \dashv\varphi^*$.
\item In the middle-left, the morphism \[\nicexy{(\Gone)^{(T,\sigma)} \ar[r]^-{\omega} & \varphi^*(\Gtwo)^{(T,\varphi\sigma)}}\in\goneopset\] is the identity map on $\Ed(T)$ on color sets. Each vertex in $T$, as a generator of $(\Gone)^{(T,\sigma)}$, is sent to itself as an element in $\varphi^*(\Gtwo)^{(T,\varphi\sigma)}$.  The other morphism named $\omega$ is defined in the same way with $T'$ in place of $T$.
\end{itemize}
To check that $(\varphi_!,\varphi_!^2,\varphi_!^0)$ is a symmetric comonoidal functor, it suffices to consider only representable $\Treecatgone$-presheaves.  For representables, the coassociativity of $\varphi_!^2$ boils down to the associativity of $\varphi^*_2$.

That the right adjoint $\varphi^*$ is a symmetric monoidal functor follows from the fact that its left adjoint $\varphi_!$ is a symmetric comonoidal functor and Theorem 1.2 in \cite{kelly-adjunction}.
\end{proof}

\begin{remark}
For the symmetric monoidal functor \[\nicexy{\Settotreecatgoneop & \Settotreecatgtwoop \ar[l]_-{\varphi^*}},\] the monoidal structure morphism $\varphi^*_0$ on the monoidal unit is \eqref{varphistar-unit}.  For $X,Y\in\Settotreecatgtwoop$, the monoidal structure morphism $\varphi^*_2$ is adjoint to the composite \[\nicexy{\varphi_!\Bigl(\varphi^*X \tensorsub{\Gone} \varphi^*Y\Bigr) \ar[r]^-{\varphi_!^2} & \varphi_!\varphi^*X \tensorsub{\Gtwo} \varphi_!\varphi^*Y \ar[r]^-{(\eta_X,\eta_Y)} & X\tensorsub{\Gtwo}Y}\] in $\Settotreecatgtwoop$.  Here each $\eta$ is the counit of the adjunction, and $\varphi^2_!$ is the comonoidal structure morphism of $\varphi_!$.\dqed
\end{remark}

\begin{example}\label{ex:realization-nerve-morphisms}
Proposition \ref{yoneda-varphi}, Theorem \ref{gnerve-changeofg}, and Theorem \ref{treecatg-adjunction-comonoidal} apply to the morphisms of action operads in Example \ref{ex:all-augmented-group-operads}.  For example, the morphisms \[\nicexy{\P\ar[r]^-{\iota} & \PCac\ar[r]^-{\iota} & \Cac\ar[r]^-{\pi} & \S}\] of action operads induce the diagram
\[\nicexy@C+.5cm{\Treecatp \ar[d]_-{\Yonedap} \ar[r]^-{\Treecat^{\iota}} & \Treecatpcac \ar[d]_-{\Yonedapcac} \ar[r]^-{\Treecat^{\iota}} & \Treecatcac \ar[d]_-{\Yonedacac} \ar[r]^-{\Treecat^{\pi}} & \Treecats \ar[d]_-{\Yonedas}\\
\Settotreecatpop \ar@<2pt>[r]^-{\iota_!} \ar@<-2pt>[d]_-{\Realp} & \Settotreecatpcacop \ar@<2pt>[r]^-{\iota_!} \ar@<-2pt>[d]_-{\Realpcac} \ar@<2pt>[l]^-{\iota^*} & \Settotreecatcacop \ar@<2pt>[r]^-{\pi_!} \ar@<-2pt>[d]_-{\Realcac} \ar@<2pt>[l]^-{\iota^*} & \Settotreecatsop \ar@<-2pt>[d]_-{\Reals} \ar@<2pt>[l]^-{\pi^*}\\
\popset \ar@<2pt>[r]^-{\iota_!} \ar@<-2pt>[u]_-{\Nervep} & \pcacopset \ar@<2pt>[r]^-{\iota_!} \ar@<-2pt>[u]_-{\Nervepcac} \ar@<2pt>[l]^-{\iota^*} & \cacopset \ar@<2pt>[r]^-{\pi_!} \ar@<-2pt>[u]_-{\Nervecac} \ar@<2pt>[l]^-{\iota^*} & \sopset \ar@<2pt>[l]^-{\pi^*} \ar@<-2pt>[u]_-{\Nerves}}\]
In the top half, there are natural isomorphisms \[\iota_!\Yonedap\cong \Yonedapcac\Treecat^{\iota},\quad \iota_!\Yonedapcac\cong \Yonedacac\Treecat^{\iota},\andspace \pi_!\Yonedacac\cong \Yonedas\Treecat^{\pi}\] by Proposition \ref{yoneda-varphi}.  

Furthermore, in each of the three squares in the bottom half, the right (resp., left) adjoint diagram is commutative up to a natural isomorphism by Theorem \ref{gnerve-changeofg}.  For example, in the bottom-right square, there are natural isomorphisms 
\[\Nervecac\pi^*\cong\pi^*\Nerves \andspace \pi_!\Realcac\cong\Reals\pi_!\] relating the cactus nerve (resp., realization) functor and the symmetric nerve (resp., realization) functor.

Using the ribbon group operad $\R$ or the braid group operad $\B$ instead, we obtain two other diagrams with similar properties.\dqed
\end{example}

\chapter{Nerve Theorem for Group Operads}\label{ch:infinity-goperad}

A \index{quasi-category}quasi-category, also known as an\index{infinity-category@$\infty$-category} $(\infty,1)$-category, is a simplicial set in which each inner horn admits a filler.  Quasi-categories were introduced by Boardman and Vogt \cite{boardman-vogt} as\index{weak Kan complexes} \emph{weak Kan complexes}.  They came into prominence in $\infty$-category theory through the work of Joyal  \cite{joyal} and Lurie \cite{lurie}.  A \emph{strict quasi-category} is a quasi-category in which each inner horn admits a unique filler.  A simplicial set is a strict quasi-category if and only if it satisfies the Segal condition \cite{grothendieck,segal}, which in turn is equivalent to being isomorphic to the nerve of some small category.  This part of $(\infty,1)$-category theory was extended to the dendroidal setting in \cite{cm11b,mt,mw07,mw09}, giving  equivalent descriptions of the dendroidal nerve in terms of strict $\infty$-symmetric operads and a symmetric operad analogue of the Segal condition.  The analogues for strict $\infty$-cyclic operads \cite{hry-cyclic}, strict $\infty$-properads, and strict $\infty$-wheeled properads \cite{hry15} are also true.

Suppose $\G$ is an action operad as in Definition \ref{def:augmented-group-operad}.  The are two main purposes of this chapter. 
\begin{enumerate}\item We introduce the $\G$-operad analogues of (strict) quasi-categories, called (strict) $\infty$-$\G$-operads.
\item We characterize the $\G$-nerve of a $\G$-operad as a $\Treecatg$-presheaf that satisfies a $\G$-operad analogue of the Segal condition, which in turn is equivalent to being a strict $\infty$-$\G$-operads.
\end{enumerate}

In Section \ref{sec:gnerve-segal} we first define a $\G$-operad analogue of the Segal condition.  Then we observe that the $\G$-nerve of each $\G$-operad satisfies the $\G$-Segal condition.  The $\G$-operad analogues of faces, horns, and (strict) quasi-categories are defined in Section \ref{sec:segal-infinity}.  A (strict) $\infty$-$\G$-operad is a $\Treecatg$-presheaf in which each inner horn, defined using inner and outer cofaces in $\Treecatg$, has a (unique) filler.  Then we observe that each $\Treecatg$-presheaf that satisfies the $\G$-Segal condition is automatically a strict $\infty$-$\G$-operad.  In particular, the $\G$-nerve of each $\G$-operad is a strict $\infty$-$\G$-operad.  

In Section \ref{sec:infinity-nerve} we show conversely that each strict $\infty$-$\G$-operad satisfies the $\G$-Segal condition.  In Section \ref{sec:nerve-characterization} we show that for a $\Treecatg$-presheaf, (i) being the $\G$-nerve of some $\G$-operad, (ii) being a strict $\infty$-$\G$-operad, and (iii) the $\G$-Segal condition are equivalent.
  
To avoid too much repetition, throughout the rest of this chapter, suppose \[\nicexy@C+3.5cm{(T,\sigma) \ar[r]^-{\phi\,=\,\bigl(((U,\sigma^U),u); (T_t,\sigma_t)_{t\in\Vt(T)}\bigr)} & (V,\sigma^V)}\] is a morphism in the $\G$-tree category $\Treecatg$ as in Definition \ref{def:treecatg}.

\section{Nerve Satisfies Segal Condition}\label{sec:gnerve-segal}

The purpose of this section is to prove that the $\G$-nerve of each $\G$-operad in $\Set$ satisfies a $\G$-operad analogue of the Segal condition, which is defined below.  

\begin{definition}\label{def:segal-condition}
Suppose $X\in\Settotreecatgop$ and $(T,\sigma)\in\Treecatg$ has at least one vertex.
\begin{enumerate}
\item Define the\index{corolla base} \emph{corolla base} of $X(T,\sigma)$, denoted \[X(T,\sigma)_1 = \Bigl(\prodover{t\in\Vt(T)} X(\Cor_t,\id_t) \Bigr)_{X(\uparrow,\id_1)},\] as the limit of the diagram consisting of the maps \[\nicexy{X(\Cor_u,\id_u) \ar[dr]^-{X(\phi_i)} && X(\Cor_t,\id_t) \ar[dl]_-{X(\phi_0)}\\ & X(\uparrow,\id_1) &}\] for each internal edge $e$ in $T$ with initial vertex $t$ and terminal vertex $u$, in which $e$ is the $i$th input of $u$.  Here $\Cor_t = \Cor_{|\inp(t)|}$ and $\id_t=\id_{|\inp(t)|}$.  The morphisms \[\nicexy{(\uparrow,\id_1) \ar[r]^-{\phi_0} & (\Cor_t,\id_t)} \andspace \nicexy{(\uparrow,\id_1) \ar[r]^-{\phi_i} & (\Cor_u,\id_u)}\] are the outer cofaces in Example \ref{ex:outer-coface-corolla} corresponding to the output of $\Cor_t$ and the $i$th input of $\Cor_u$, respectively.
\item The \emph{$\G$-Segal map}\index{G-Segal map@$\G$-Segal map} is the map 
\begin{equation}\label{g-segal-map} 
\nicexy{X(T,\sigma) \ar[r]^-{\chi_{(T,\sigma)}} & X(T,\sigma)_1}
\end{equation}
determined by the commutative diagrams \[\nicexy{X(T,\sigma) \ar[dr]_-{X(\xi_t)} \ar[r]^-{\chi_{(T,\sigma)}} & X(T,\sigma)_1 \ar[d]^-{\mathrm{natural}}\\ & X(\Cor_t,\id_t)}\] for $t\in\Vt(T)$, in which $\xi_t$ is the vertex inclusion 
\begin{equation}\label{vertex-inclusion}
\nicexy@C+3cm{(\Cor_t,\id_t) \ar[r]^-{\xi_t\,=\, \bigl(((T,\sigma),t); (\Cor_t,\id_t)\bigr)} & (T,\sigma) \in \Treecatg}.
\end{equation}
\item We say that $X$ satisfies the \emph{$\G$-Segal condition}\index{G-Segal condition@$\G$-Segal condition} if the $\G$-Segal map $\chi_{(T,\sigma)}$ is a bijection for each $(T,\sigma)\in\Treecatg$ with at least one vertex.
\end{enumerate}
\end{definition}

\begin{example}
When $\G$ is the action operad $\P$, $\S$, $\B$, $\PB$, $\R$, $\PR$, $\Cac$, or $\PCac$, we refer to the $\G$-Segal map/condition as the \index{planar Segal condition}planar, \index{symmetric Segal condition}symmetric, \index{braided Segal condition}braided, \index{pure braided Segal condition}pure braided, \index{ribbon Segal condition}ribbon, \index{pure ribbon Segal condition}pure ribbon, \index{cactus Segal condition}cactus, or \index{pure cactus Segal condition}pure cactus Segal map/condition.\dqed
\end{example}

The following description of the corolla base is immediate from the definition.

\begin{proposition}\label{corolla-base-description}
Each element in the corolla base of $X(T,\sigma)$ is a collection of elements \[\{x_t\}\in\prodover{t\in\Vt(T)} X(\Cor_t,\id_t),\] one for each vertex in $T$, that satisfies the following condition.  For each internal edge $e$ in $T$ with initial vertex $t$ and terminal vertex $u$, in which $e$ is the $i$th input of $u$, the equality \[X(\phi_i)(x_u) = X(\phi_0)(x_t) \in X(\uparrow,\id_1)\] holds.  
\end{proposition}

\begin{interpretation}\label{int:corolla-base}
One can think of an element $\{x_t\}$ in the corolla base of $X(T,\sigma)$ as a decoration of the planar tree $T$ by $X$, with each vertex $t \in \Vt(T)$ decorated by an element $x_t\in X(\Cor_t,\id_t)$.  The previous equality expresses the compatibility of the elements $\{x_t\}$ over the internal edges in $T$.  The $\G$-Segal condition says that each element in the corolla base of $X(T,\sigma)$ has a unique extension to an element in $X(T,\sigma)$ via the $\G$-Segal map.\dqed
\end{interpretation}

Next we provide a more conceptual interpretation of the $\G$-Segal map as a map  corepresented by a suitable morphism of $\Treecatg$-presheaves, for which we will use the following concept.

\begin{definition}\label{def:segal-core}
Suppose $(T,\sigma)\in\Treecatg$ with at least one vertex.  Define its\index{G-Segal core@$\G$-Segal core} \emph{$\G$-Segal core} \[\Segalg(T,\sigma)\in \Settotreecatgop\] as the union \[\Segalg(T,\sigma) = \bigcup_{t\in\Vt(T)} \face{\xi_t}(T,\sigma) \subseteq \Yonedag(T,\sigma)\] with
\begin{itemize}
\item each \[\nicexy{(\Cor_t,\id_t) \ar[r]^-{\xi_t} & (T,\sigma)} \in \Treecatg\] the vertex inclusion in \eqref{vertex-inclusion} and 
\item $\face{\xi_t}(T,\sigma)$ the image\label{not:facexit} \[\nicexy{\Image\Bigl(\Yonedag(\Cor_t,\id_t) \ar[r]^-{\xi_t \circ-} & \Yonedag(T,\sigma)\Bigr)}\in \Settotreecatgop\] induced by $\xi_t$.
\end{itemize}
We write\label{not:segalcor-inclusion} \[\nicexy{\Segalg(T,\sigma)\ar[r]^-{\iota_{(T,\sigma)}} & \Yonedag(T,\sigma)} \in \Settotreecatgop\] for the inclusion, called the \emph{$\G$-Segal core inclusion}.\index{G-Segal core@$\G$-Segal core!inclusion}
\end{definition}

The following observation says that the $\G$-Segal map is corepresented by the $\G$-Segal core inclusion.

\begin{proposition}\label{segal-core}
Suppose $X\in\Settotreecatgop$ and $(T,\sigma)\in\Treecatg$ with at least one vertex.  Then there is a commutative diagram \[\nicexy{X(T,\sigma) \ar[d]_-{\cong} \ar[r]^-{\chi_{(T,\sigma)}} & X(T,\sigma)_1 \ar[d]^-{\cong}\\ \Settotreecatgop\bigl(\Yonedag(T,\sigma),X\bigr) \ar[r]^-{\iota_{(T,\sigma)}^*} & \Settotreecatgop\bigl(\Segalg(T,\sigma),X\bigr)}\]
that is natural in $(T,\sigma)$.
\end{proposition}

\begin{proof}
The vertical isomorphisms follow from Yoneda Lemma, and the commutativity of the diagram follows by inspection.
\end{proof}

Next we observe that, to check the $\G$-Segal condition, it suffices to check a subset of the $\G$-Segal maps.

\begin{proposition}\label{segal-characterization}
Suppose $X\in \Settotreecatgop$.  Then $X$ satisfies the $\G$-Segal condition if and only if the $\G$-Segal map \[\nicexy@C+1cm{X(T,\id_{|\inp(T)|}) \ar[r]^-{\chi_{(T,\id_{|\inp(T)|})}} & X(T,\id_{|\inp(T)|})_1}\] is a bijection for each planar tree $T$ with at least two vertices. 
\end{proposition}

\begin{proof}
First note that if $(T,\sigma)$ has exactly one vertex, then the corresponding $\G$-Segal map is a bijection.  Indeed, $T$ must be a corolla, which has no internal edges, and the corolla base of $X(\Cor_{|\inp(T)|},\sigma)$ is $X(\Cor_{|\inp(T)|},\id_{|\inp(T)|})$.  The vertex inclusion $\xi_t \in \Treecatg$, where $t$ is the unique vertex in $T=\Cor_{|\inp(T)|}$, is an isomorphism.  It follows that the $\G$-Segal map \[\chi_{(\Cor_{|\inp(T)|},\id_{|\inp(T)|})} = X(\xi_t)\] is a bijection.  So to check the $\G$-Segal condition, it is enough to check that the $\G$-Segal maps $\chi_{(T,\sigma)}$ are bijections whenever $T$ has at least two vertices.

For a $\G$-tree $(T,\sigma)$ with at least two vertices, consider the isomorphism \[\nicexy@C+4cm{(T,\id_{|\inp(T)|}) \ar[r]^-{\phi\,=\, \Bigl(((\Cor_{|\inp(T)|},\sigma),u) ;(\Cor_t,\id_t)_{t\in\Vt(T)}\Bigr)}_-{\cong} & (T,\sigma)} \in \Treecatg.\]  Then the diagram \[\nicexy@C+1cm{X(T,\id_{|\inp(T)|}) \ar[r]^-{\chi_{(T,\id_{|\inp(T)|})}} & X(T,\id_{|\inp(T)|})_1 \ar@{=}[d]\\ X(T,\sigma) \ar[u]^-{\phi}_-{\cong} \ar[r]^-{\chi_{(T,\sigma)}} & X(T,\sigma)_1}\] is commutative.  So the bottom $\G$-Segal map $\chi_{(T,\sigma)}$ is a bijection if and only if the top $\G$-Segal map $\chi_{(T,\id_{|\inp(T)|})}$ is a bijection.
\end{proof}

Here is the main observation in this section.  Recall the $\G$-nerve in Definition \ref{def:gnerve}.

\begin{theorem}\label{nerve-is-segal}
For each action operad $\G$, the $\G$-nerve of each colored $\G$-operad in $\Set$ satisfies the $\G$-Segal condition.\index{G-Segal condition@$\G$-Segal condition!G-nerve@$\G$-nerve}\index{nerve!satisfies G-Segal condition@satisfies $\G$-Segal condition}
\end{theorem}

\begin{proof}
Suppose $\O$ is a $\colorc$-colored $\G$-operad.  By Proposition \ref{segal-characterization} it suffices to check that the $\G$-Segal map \[\nicexy@C+.5cm{\Nerveg(\colorc,\O)(T,\id_T) \ar[r]^-{\chi_{(T,\id_T)}} & \Nerveg(\colorc,\O)(T,\id_T)_1}\] is a bijection for each planar tree $T$ with at least two vertices, where $\id_T=\id_{|\inp(T)|} \in \G(|\inp(T)|)$.  Recall that $\G^{(T,\id_T)}$ is the $\Ed(T)$-colored $\G$-operad in $\Set$ freely generated by the set $\Vt(T)$ of vertices in $T$.  Each element in the set\[\Nerveg(\colorc,\O)(T,\id_T) = \gopset\bigl(\G^{(T,\id_T)},\O\bigr)\] consists of:
\begin{enumerate}[label=(\roman*)]
\item a map $\phi : \Ed(T)\to\colorc$;
\item an element $\phi(t)\in\O\sbinom{\phi\,\out(t)}{\phi\,\inp(t)}$ for each $t\in\Vt(T)$.
\end{enumerate}

The above data is equivalent to an element in the corolla base \[\Nerveg(\colorc,\O)(T,\id_T)_1 = \Bigl(\prod_{t\in\Vt(T)} \Nerveg(\colorc,\O)(\Cor_t,\id_t)\Bigr)_{\Nerveg(\colorc,\O)(\uparrow,\id_1)},\] as discussed in Interpretation \ref{int:corolla-base}, via the $\G$-Segal map $\chi_{(T,\id_T)}$.  Indeed, for each $t\in\Vt(T)$, each element \[\phi_t \in\Nerveg(\colorc,\O)(\Cor_t,\id_t) = \gopset\bigl(\G^{(\Cor_t,\id_t)},\O\bigr)\]consists of:
\begin{enumerate}[label=(\roman*)]
\item an element $\phi_t(e)\in\colorc$ for each edge $e\in t$;
\item an element $\phi_t \in \O\sbinom{\phi_t\,\out(t)}{\phi_t\,\inp(t)}$
\end{enumerate}
Similarly, each element in $\Nerveg(\colorc,\O)(\uparrow,\id_1)$ is a choice of an element in $\colorc$.  So each element in the corolla base $\Nerveg(\colorc,\O)(T,\id_T)_1$ is a collection of elements $\{\phi_t\in\O\}_{t\in\Vt(T)}$ as above such that for each internal edge, say $e$ with initial vertex $t$ and terminal vertex $u$, the output profile of $\phi_t$ coincides with the input profile of $\phi_u$ corresponding to $e$.
\end{proof}

\section{Segal Condition Implies Strict Infinity}\label{sec:segal-infinity}

In this section, for each action operad $\G$, we define the $\G$-operad analogues of faces, horns, and (strict) quasi-categories.  The main result says that every $\Treecatg$-presheaf that satisfies the $\G$-Segal condition is a strict $\infty$-$\G$-operad.  As a result, the $\G$-nerve of each colored $\G$-operad in $\Set$ is a strict $\infty$-$\G$-operad.  Recall the concepts of an inner coface and of an outer coface in Definitions \ref{def:inner-coface} and \ref{def:outer-coface}.

\begin{definition}
A\index{coface}\index{G-tree category@$\G$-tree category!coface} \emph{coface} is either an inner coface or an outer coface.  For a $\G$-tree $(V,\sigma^V)$, a \emph{coface of $(V,\sigma^V)$} is a coface with codomain $(V,\sigma^V)$.
\end{definition}

\begin{definition}\label{def:horn}
Suppose $\phi : (T,\sigma) \to (V,\sigma^V)$ is a coface of $(V,\sigma^V)$.
\begin{enumerate}
\item If $\phi' : (T',\sigma') \to (V,\sigma^V)$ is another coface of $(V,\sigma^V)$, then we say that $\phi$ is \emph{isomorphic} to $\phi'$ if there exists an isomorphism \[\theta : (T,\sigma)\iso(T',\sigma')\in\Treecatg\] such that the diagram \[\nicexy{(T,\sigma) \ar[r]^-{\phi} \ar[d]_-{\theta}^-{\cong} & (V,\sigma^V) \ar@{=}[d]\\ (T',\sigma') \ar[r]^-{\phi'} & (V,\sigma^V)}\] is commutative.  In this case, we write $\phi \cong \phi'$.
\item The\label{not:phiface} \emph{$\phi$-face}\index{face}\index{G-tree category@$\G$-tree category!face} of $\Yonedag(V,\sigma^V)$ is defined as the image \[\face{\phi}(V,\sigma^V) = \Image\Bigl(\nicexy{\Yonedag(T,\sigma) \ar[r]^-{\phi\circ -} & \Yonedag(V,\sigma^V)}\Bigr) \in \Settotreecatgop\] induced by $\phi$. 
\item The\label{not:phihorn} \emph{$\phi$-horn}\index{horn}\index{G-tree category@$\G$-tree category!horn} of $\Yonedag(V,\sigma^V)$ is defined as the union \[\horn{\phi}(V,\sigma^V) = \bigcup_{\phi'\not\cong\phi} \face{\phi'}(V,\sigma^V) \subseteq \Yonedag(V,\sigma^V)\in \Settotreecatgop\] indexed by all the cofaces of $(V,\sigma^V)$ not isomorphic to $\phi$.
\item A \emph{horn} of $\Yonedag(V,\sigma^V)$ is a $\phi$-horn for some coface $\phi$ of $(V,\sigma^V)$.  An \emph{inner horn} is the $\phi$-horn in which $\phi$ is an inner coface.
\item For $X\in \Settotreecatgop$, a \emph{horn of $X$} is a morphism of the form \[\nicexy{\horn{\phi}(V,\sigma^V) \ar[r] & X} \in \Settotreecatgop\] for some $(V,\sigma^V)\in\Treecatg$ and some coface $\phi$.  
\item An\index{inner horn} \emph{inner horn of $X$} is a horn of $X$ in which $\horn{\phi}(V,\sigma^V)$ is an inner horn.
\end{enumerate}
\end{definition}

Let us first provide an explicit description of a horn in a $\Treecatg$-presheaf.  To simplify the notation below, we will write a morphism in $\Treecatg$ and its image in $\Settotreecatgop$ under the Yoneda embedding using the same symbol.

\begin{lemma}\label{horn-description}
Suppose $X\in \Settotreecatgop$, $(V,\sigma^V) \in\Treecatg$, and $\phi$ is a coface of $(V,\sigma^V)$.  Then a horn \[f : \horn{\phi}(V,\sigma^V) \to X\in\Settotreecatgop\] is equivalent to a family of morphisms \[\Bigl\{\narrowxy{\Yonedag(T',\sigma')\ar[r]^-{x_{\phi'}} & X} : \text{cofaces $\narrowxy{(T',\sigma')\ar[r]^-{\phi'} & (V,\sigma^V)}$ with $\phi\not\cong\phi'$}\Bigr\},\] one for each coface of $(V,\sigma^V)$ not isomorphic to $\phi$, that satisfies the following condition.  Whenever
\begin{equation}\label{in-union-faces}
\nicexy{(Q,\sigma^Q) \ar[d]_-{\alpha'} \ar[r]^-{\alpha''} & (T'',\sigma'') \ar[d]^-{\phi''}\\ (T',\sigma') \ar[r]^-{\phi'} & (V,\sigma^V)}
\end{equation}
is a commutative diagram of cofaces in $\Treecatg$ with $\phi' \not\cong \phi \not\cong\phi''$, the diagram 
\begin{equation}\label{map-union-to-x}
\nicexy{\Yonedag(Q,\sigma^Q) \ar[d]_-{\alpha'} \ar[r]^-{\alpha''} & \Yonedag(T'',\sigma'') \ar[d]^-{x_{\phi''}}\\ \Yonedag(T',\sigma') \ar[r]^-{x_{\phi'}} & X}
\end{equation}
in $\Settotreecatgop$ is also commutative.
\end{lemma}

\begin{proof}
Given a horn $f$ and a coface $\phi' \not\cong \phi$, the morphism $x_{\phi'}$ is the composite
\begin{equation}\label{xphi-f}
\nicexy{\Yonedag(T',\sigma') \ar[r]^-{x_{\phi'}} \ar[d]_-{\phi'\circ-} & X\\ \face{\phi'}(V,\sigma^V) \ar@{>->}[r] & \horn{\phi}(V,\sigma^V) \ar[u]_-{f}}
\end{equation}
in $\Settotreecatgop$.  The commutative square \eqref{in-union-faces} means that \[\begin{split}
\phi'\alpha' = \phi''\alpha''& \in \Bigl(\face{\phi'}(V,\sigma^V) \cap \face{\phi''}(V,\sigma^V)\Bigr)(Q,\sigma^Q)\\
& \subseteq \Yonedag(V,\sigma^V)(Q,\sigma^Q).
\end{split}\]  
The square \eqref{map-union-to-x} is commutative because the morphism $f$ is well-defined on the intersection of the $\phi'$-face and the $\phi''$-face.  

Conversely, given a family of morphisms $\{x_{\phi'}\}$ as above, we define the morphism $f$ using the diagram \eqref{xphi-f}.  The diagrams \eqref{in-union-faces} and \eqref{map-union-to-x} and Lemma \ref{treecatg-factorization} ensure that $f$ is well-defined.
\end{proof}

We are now ready to define the $\G$-operad analogue of an $(\infty,1)$-category.

\begin{definition}\label{def:infinity-goperad}
Suppose $X\in\Settotreecatgop$.
\begin{enumerate}
\item We call $X$ an \emph{$\infty$-$\G$-operad}\index{G-operad@$\G$-operad!infinity}\index{infinity-G-operad@$\infty$-$\G$-operad} if for each inner horn $f$ of $X$, 
\begin{equation}\label{horn-filler}
\nicexy{\horn{\phi}(V,\sigma^V) \ar[r]^-{f} \ar@{>->}[d] & X\\ \Yonedag(V,\sigma^V) \ar@{-->}[ur] &}
\end{equation}
a dashed arrow, called an\index{inner horn filler} \emph{inner horn filler}, exists that makes the diagram commutative.
\item A\index{G-operad@$\G$-operad!strict infinity}\index{strict infinity-G-operad@strict $\infty$-$\G$-operad} \emph{strict $\infty$-$\G$-operad} is an $\infty$-$\G$-operad for which each inner horn filler in \eqref{horn-filler} is unique.
\end{enumerate}
\end{definition}

\begin{remark}
The simplicial set version of Definition \ref{def:infinity-goperad}(1) is due to Boardman and Vogt, who called the simplicial set version of the lifting condition \eqref{horn-filler} the\index{restricted Kan condition} \emph{restricted Kan condition}.  See Definition 4.8 on page 102 in \cite{boardman-vogt}.  The restricted Kan condition is the subset of the usual Kan lifting condition in which the omitted face is neither the top face nor the bottom face.  In more modern terminology, a simplicial set that satisfies the restricted Kan condition is called either a\index{quasi-category} \emph{quasi-category} or an\index{infinity-category@$\infty$-category} \emph{$(\infty,1)$-category}.  Roughly speaking, a quasi-category is a category in which the composition, identities, and associativity axiom have all been relaxed.  It is known that being a strict $(\infty,1)$-category, in which the inner horn fillers are unique, is equivalent to being the nerve of a small category.  See, e.g., Proposition 1.1.2.2 in \cite{lurie}.\dqed
\end{remark}

\begin{example}[Infinity Planar and Symmetric Operads]\label{ex:inf-psop}
When $\G$ is the planar group operad $\P$ or the symmetric group operad $\S$, Definition \ref{def:infinity-goperad} defines\index{infinity-planar operad@$\infty$-planar operad}\index{planar operad!infinity}\index{strict infinity-planar operad@strict $\infty$-planar operad} (strict) $\infty$-planar operads and (strict) $\infty$-symmetric operads.  For the symmetric dendroidal category $\Omega \simeq \Treecats$, (strict) $\infty$-symmetric operads\index{infinity-symmetric operad@$\infty$-symmetric operad}\index{symmetric operad!infinity}\index{strict infinity-symmetric operad@strict $\infty$-symmetric operad} were introduced in Section 5 in \cite{mw09} and were called (strict) inner Kan complexes.\dqed
\end{example}

\begin{example}[Infinity Braided, Ribbon, and Cactus Operads]\label{ex:inf-brop}
When $\G$ is the action operad $\B$, $\PB$, $\R$, $\PR$, $\Cac$, or $\PCac$, Definition \ref{def:infinity-goperad} defines\index{infinity-braided operad@$\infty$-braided operad}\index{braided operad!infinity}\index{strict infinity-braided operad@strict $\infty$-braided operad} \emph{(strict) $\infty$-braided operads}, \emph{(strict) $\infty$-pure braided operads}, \index{infinity-ribbon operad@$\infty$-ribbon operad}\index{ribbon operad!infinity}\index{strict infinity-ribbon operad@strict $\infty$-ribbon operad}\emph{(strict) $\infty$-ribbon operads}, \emph{(strict) $\infty$-pure ribbon operads}, \index{infinity-cactus operad@$\infty$-cactus operad}\index{cactus operad!infinity}\index{strict infinity-cactus operad@strict $\infty$-cactus operad}\emph{(strict) $\infty$-cactus operads}, and \emph{(strict) $\infty$-pure cactus operads}.   One should think of an $\infty$-braided/ribbon/cactus operad as a braided/ribbon/cactus operad in which the axioms are replaced by suitable homotopical analogues.  A special case of the main result of this chapter states that a $\Treecatb$/$\Treecatr$/$\Treecatcac$-presheaf is a strict $\infty$-braided/ribbon/cactus operad if and only if it satisfies the braided/ribbon/cactus Segal condition.  In particular, the braided/ribbon/cactus nerve of a braided/ribbon/cactus operad in $\Set$ is a strict $\infty$-braided/ribbon/cactus operad.  See Corollary \ref{nerve-is-strict} and Corollary \ref{segal=strict}.\dqed
\end{example}

Here is the main observation of this section.  Recall the $\G$-Segal condition in Definition \ref{def:segal-condition}.

\begin{theorem}\label{segal-is-strict}
Every $\Treecatg$-presheaf that satisfies the $\G$-Segal condition\index{G-Segal condition@$\G$-Segal condition!strict infinity-G-operad@strict $\infty$-$\G$-operad} is a strict $\infty$-$\G$-operad.
\end{theorem}

\begin{proof}
Suppose $X\in\Settotreecatgop$ satisfies the $\G$-Segal condition, and \[\nicexy{\horn{\phi}(V,\sigma^V) \ar[r]^-{f} & X}\in\Settotreecatgop\] is an inner horn of $X$.  We must show that $f$ has a unique extension $g$ that makes the diagram\[\nicexy{\horn{\phi}(V,\sigma^V) \ar[r]^-{f} \ar@{>->}[d] & X\\ \Yonedag(V,\sigma^V) \ar[ur]_g &}\]commutative.  Since $\phi$ is an inner coface, $V$ has at least one internal edge and at least two vertices.  By Lemma \ref{horn-description}, the inner horn $f$ is equivalent to a family of elements \[f=\bigl\{f_{\phi'} \in X(T',\sigma')\bigr\},\] one for each coface $\phi'$ of $(V,\sigma^V)$ not isomorphic to $\phi$, that are compatible as expressed in the diagram \eqref{map-union-to-x}.  

Since $\phi$ is an \emph{inner} coface of $(V,\sigma^V)$, for each vertex $v$ in $V$, there exists a factorization\[\nicexy{(\Cor_v,\id_v) \ar[d]_-{\xi_v'} \ar[r]^-{\xi_v} & (V,\sigma^V)\\ (T',\sigma')\ar[ur]_-{\phi'} &}\]of the vertex inclusion $\xi_v$ in \eqref{vertex-inclusion}, such that:
\begin{itemize}
\item $\phi'$ is an outer coface of $(V,\sigma^V)$.
\item $v\in\Vt(T')$ with vertex inclusion $\xi_v'$.
\end{itemize}
Note that the outer coface $\phi'$ is not uniquely determined by the vertex $v$; i.e., $\xi_v$ may factor through another outer coface of $(V,\sigma^V)$.  Nevertheless, the element \[f_v = X(\xi_v')(f_{\phi'}) \in X(\Cor_v,\id_v)\] is well-defined by the compatibility diagram \eqref{map-union-to-x}, which also implies that these elements are compatible over internal edges in $V$.

By the assumed $\G$-Segal condition on $X$, an element in $X(V,\sigma^V)$ is uniquely determined by its image in the corolla base \[X(V,\sigma^V)_1 = \Bigl(\prod_{v\in\Vt(V)} X(\Cor_v,\id_v)\Bigr)_{X(\uparrow,\id_1)}\] via the $\G$-Segal map $\chi_{(V,\sigma^V)}$.  So the collection \[g = \{f_v\}_{v\in\Vt(T)}\in X(V,\sigma^V)_1\] in the previous paragraph determines a unique element in $X(V,\sigma^V)$.  This establishes the uniqueness of an extension of $f$ and provides a candidate for such an extension.

To show that $g=\{f_v\}_{v\in\Vt(V)} \in X(V,\sigma^V)$ is indeed an extension of $f$, we need to show that for each coface $\phi'$ of $(V,\sigma^V)$ not isomorphic to $\phi$, we have that\[X(\phi')(g)=f_{\phi'} \in X(T',\sigma').\]  If $\phi'$ is an outer coface, then every vertex in $T'$ is also a vertex in $V$.  By the $\G$-Segal condition on $X$, the element $f_{\phi'}$ is uniquely determined by the images  
\[X(\xi_v')(f_{\phi'})=f_v\in X(\Cor_v,\id_v)\] for $v\in\Vt(T')\subseteq \Vt(V)$.  Since $g$ also restricts to these elements $f_v$'s via the vertex inclusions, we conclude that $X(\phi')(g)=f_{\phi'}$.

Next suppose $\phi'$ is an inner coface of $(V,\sigma^V)$ not isomorphic to $\phi$.  By the $\G$-Segal condition on $X$, the element $f_{\phi'}\in X(T',\sigma')$ is uniquely determined by the images \[\Bigl\{X(\xi_{t'})(f_{\phi'}) \in X(\Cor_{t'},\id_{t'})\Bigr\}_{t'\in\Vt(T')},\] where each morphism \[\nicexy{(\Cor_{t'},\id_{t'}) \ar[r]^-{\xi_{t'}} & (T',\sigma')}\in\Treecatg\] is the vertex inclusion for $t'$.  Since $\phi'$ is an inner coface of $(V,\sigma^V)$, all but one vertex in $T'$, say $u\in\Vt(T')$, belong to $V$.  

For each vertex $t'$ in $T'$ that is already in $V$, there is a commutative diagram\[\nicexy{(\Cor_{t'},\id_{t'}) \ar[d]_-{\xi_{t'}} \ar[r]^-{\xi_{t'}} & (T'',\sigma'') \ar[d]^-{\phi''} \\ (T',\sigma') \ar[r]^-{\phi'} & (V,\sigma^V)}\]in which:
\begin{itemize}
\item $\phi''$ is an outer coface of $(V,\sigma^V)$ with $t'\in\Vt(T'')$.
\item The top horizontal $\xi_{t'}$ is the vertex inclusion for $t'$ as a vertex in $T''$.
\end{itemize}
We now have:
\[\begin{split}
X(\xi_{t'})X(\phi')(g) &= X(\xi_{t'})X(\phi'')(g)\\
&= f_{t'}\\
&= X(\xi'_{t'})(f_{\phi'}).
\end{split}\]

Therefore, it remains to show that \[X(\xi_u)(f_{\phi'}) = X(\phi'\xi_u)(g)\in X(\Cor_u,\id_u),\] where $\phi'\xi_u$ is the composite \[\nicexy{(\Cor_u,\id_u) \ar[r]^-{\xi_u} & (T',\sigma') \ar[r]^-{\phi'} & (V,\sigma^V)} \in \Treecatg.\]  We can factor this composite as
\begin{equation}\label{xi-phi'}
\nicexy{(\Cor_u,\id_u)\ar[d]_-{\xi_u} \ar[r]^-{\varphi_0} & (T'_u,\sigma'_u) \ar[r]^-{\varphi_1} & (T'',\sigma'') \ar[d]^-{\phi''}\\ (T',\sigma') \ar[rr]^-{\phi'} && (V,\sigma^V)}
\end{equation}
in which:
\begin{itemize}
\item $(T'_u,\sigma'_u)$ is the inside $\G$-tree of $\phi'$ at the vertex $u$.
\item $\varphi_0$ is the inner coface corresponding to the $\G$-tree substitution \[(T'_u,\sigma'_u) = (\Cor_u,\id_u)(T'_u,\sigma'_u).\]
\item $\phi''$ is an outer coface of $(V,\sigma^V)$.
\item $\varphi_1$ is either a composite of outer cofaces if $|\Vt(V)|>2$, or an isomorphism if $|\Vt(V)|=2$. 
\end{itemize}

We now compute as follows:
\[\begin{split}
X(\phi'\xi_u)(g) &= X(\phi''\varphi_1\varphi_0)(g)\\
&= X(\varphi_1\varphi_0)X(\phi'')(g)\\
&= X(\varphi_1\varphi_0)(f_{\phi''})\\
&= X(\xi_u)(f_{\phi'}).
\end{split}\]
The first equality is from \eqref{xi-phi'}.  The second equality is by the functoriality of $X$.  The third equality is by the previous case and the fact that $\phi''$ is an outer coface of $(V,\sigma^V)$.  The last equality is by Lemma \ref{horn-description}.
\end{proof}

Combining Theorem \ref{nerve-is-segal} and Theorem \ref{segal-is-strict}, we obtain the following result.

\begin{corollary}\label{nerve-is-strict}
The $\G$-nerve of each colored $\G$-operad in $\Set$ is a strict $\infty$-$\G$-operad.\index{nerve!strict infinity-G-operad@strict $\infty$-$\G$-operad}
\end{corollary}

\begin{example}[Planar and Symmetric Dendroidal Nerves]
When\index{nerve!dendroidal}\index{dendroidal nerve} $\G$ is the symmetric group operad $\S$, the $\S$-nerve of each symmetric operad in $\Set$ is a strict $\infty$-symmetric operad.  Recall that the symmetric tree category $\Treecats$ is equivalent to the Moerdijk-Weiss symmetric dendroidal category $\Omega$.  In the dendroidal context, Corollary \ref{nerve-is-strict} recovers Proposition 5.3 in \cite{mw09}.  Similarly, if $\P$ is the planar group operad $\P$, then the $\P$-nerve of each planar operad in $\Set$ is a strict $\infty$-planar operad.\dqed
\end{example}

\begin{example}[Braided, Ribbon, and Cactus Nerves]
If $\G$ is the braid group operad $\B$, then the\index{braided nerve}\index{nerve!braided} braided nerve of each braided operad in $\Set$ is a strict $\infty$-braided operad.  Similarly, if $\G$ is the ribbon group operad $\R$, then the\index{ribbon nerve}\index{nerve!ribbon} ribbon nerve of each ribbon operad in $\Set$ is a strict $\infty$-ribbon operad.  Also, if $\G$ is the cactus group operad $\Cac$, then the\index{cactus nerve}\index{nerve!cactus} cactus nerve of each cactus operad in $\Set$ is a strict $\infty$-cactus operad.  There are also pure analogues for the pure braid group operad $\PB$, the pure ribbon group operad $\PR$, and the pure cactus group operad $\PCac$.\dqed
\end{example}

\section{Strict Infinity Implies Segal Condition}\label{sec:infinity-nerve}

The purpose of this section is to prove the following converse of Theorem \ref{segal-is-strict}.

\begin{theorem}\label{strict-is-segal}
Every\index{strict infinity-G-operad@strict $\infty$-$\G$-operad!G-Segal condition@$\G$-Segal condition} strict $\infty$-$\G$-operad satisfies the $\G$-Segal condition.
\end{theorem}

\begin{proof}
Suppose $X\in\Settotreecatgop$ is a strict $\infty$-$\G$-operad.  To show that it satisfies the $\G$-Segal condition, we use Proposition \ref{segal-characterization}.  All the input equivariances in the rest of this proof are identities, so we will omit writing them.  We must show that the $\G$-Segal map \[\nicexy{X(V) \ar[r]^-{\chi_V} & X(V)_1 = \Bigl(\prodover{v\in\Vt(V)} X(\Cor_v)\Bigr)_{X(\uparrow)}}\] is a bijection for each planar tree $V$ with at least two vertices.  In other words, we must show that each element $\{x_v\}_{v\in\Vt(V)} \in X(V)_1$ in the corolla base has a unique extension to an element in $X(V)$.  This is proved by an induction on $n=|\Vt(V)| \geq 2$.

If $V$ has exactly two vertices, then it has only one internal edge.  By Lemma \ref{horn-description} the two elements $\{x_v\}_{v\in\Vt(V)}$ define an inner horn \[\nicexy@C+.6cm{\horn{\phi}(V) \ar@{>->}[d] \ar[r]^-{\{x_v\}_{v\in\Vt(V)}} & X\\ \Yonedag(V) \ar@{-->}[ur]}\] of $X$.  By the assumption on $X$, this inner horn has a unique dashed extension, which is the desired unique extension to an element in $X(V)$.

Next suppose $|\Vt(V)|>2$, and $\{x_v\}_{v\in\Vt(V)} \in X(V)_1$ is an element in the corolla base.  We first claim that for each coface $\phi : T \to V$ of $V$, the collection $\{x_v\}$ has a unique extension to an element $x_T\in X(T)$.  To prove this claim, consider the two cases.
\begin{enumerate}
\item If $\phi$ is an outer coface, then the sub-collection $\{x_t\}_{t\in\Vt(T)}$ has a unique extension to an element $x_T \in X(T)$ by the induction hypothesis, since $\Vt(T)$ is a proper subset of $\Vt(V)$.
\item Suppose $\phi$ is an inner coface, and $e$ is the internal edge in $V$ corresponding to $\phi$ with initial vertex $e^0$ and terminal vertex $e^1$.  Consider the planar tree $V_e$ with:
\begin{itemize}
\item $\Vt(V_e) = \{e^0,e^1\}$;
\item only one internal edge $e$;
\item $V= T\comp_t V_e$, where $t\in \Vt(T)$ is the substitution vertex for $\phi$.  
\end{itemize}
As in the $n=2$ case, the pair of elements $\{x_{e^0},x_{e^1}\}$ has a unique extension to an element $x_{V_e} \in X(V_e)$.  For the inner coface \[\phi_e : \Cor_t \to V_e,\] consider the element 
\begin{equation}\label{xtxphi}
x_t = X(\phi_e)\left(x_{V_e}\right) \in X(\Cor_t).
\end{equation}
Note that \[\Vt(T) = \{t\} \amalg \bigl[\Vt(V)\setminus\{e^0,e^1\}\bigr].\]  By the induction hypothesis applied to $T$, the collection of elements 
\begin{equation}\label{xtxv}
\Bigl\{x_t, \{x_v\}_{v\in \Vt(V)\setminus\{e^0,e^1\}}\Bigr\}
\end{equation} has a unique extension to an element $x_T \in X(T)$.
\end{enumerate}
This proves the claim.

Now pick an arbitrary inner coface $\varphi : Q \to V$ of $V$.  The collection of elements \[\Phi = \Bigl\{x_{T}\in X(T) : \text{coface $\phi : T\to V$ with $\phi\not\cong\varphi$}\Bigr\},\] one for each coface of $V$ not isomorphic to $\varphi$, determines an inner horn \[\nicexy{\horn{\varphi}(V) \ar@{>->}[d] \ar[r]^-{\Phi} & X\\ \Yonedag(V) \ar@{-->}[ur]_-{y}}\] of $X$ by Lemma \ref{horn-description}.  By the assumption on $X$, the inner horn $\Phi$ has a unique dashed extension $y\in X(V)$.  We claim that this element $y$ is the desired unique extension of the original collection $\{x_v\}_{v\in \Vt(V)}$.  

Indeed, for each vertex $v$ in $V$, the corresponding vertex inclusion $\xi_v$ must factor through some outer coface $\phi$ of $V$, as in \[\nicexy{\Cor_v \ar[r]^-{\xi_v'} & T \ar[r]^-{\phi} & V \ar@{<-}`u[ll]`[ll]_-{\xi_v}[ll]},\] by Lemma \ref{outer-coface-iterated}.  So we have 
\[\begin{split} x_v &= X(\xi_v')(x_T)\\ &= X(\xi_v')X(\phi)(y)\\ &= X(\xi_v)(y)\in X(\Cor_v).
\end{split}\] Therefore, $x_v$ is a restriction of $y$.

To prove the uniqueness of $y$, suppose $z\in X(V)$ is another extension of the collection $\{x_v\}_{v\in\Vt(V)}$.  By the assumption on $X$, it is enough to show the equality \[X(\phi)(z) = x_T\in X(T)\] for each coface $\phi : T \to V$ of $V$, since this would imply that $z$ is also an extension of the inner horn $\Phi$, which is unique.
\begin{enumerate}
\item If $\phi$ is an outer coface, then $x_T$ is uniquely determined by the subset $\{x_t\}_{t\in\Vt(T)}$ of $\{x_v\}_{v\in\Vt(V)}$.  Since $X(\phi)(z)$ also restricts to $x_t$ for each $t\in \Vt(T)\subseteq \Vt(V)$, it follows that $X(\phi)(z) = x_T$.
\item Suppose $\phi$ is an inner coface.  The element $x_T$ is uniquely determined by the collection \eqref{xtxv}.  Therefore, it remains to show that \[X(\xi_t)X(\phi)(z)=x_t\in X(\Cor_t)\] with:
\begin{itemize}
\item $x_t$ the element in \eqref{xtxphi};
\item $\xi_t : \Cor_t \to T$ the vertex inclusion.
\end{itemize} 
The composite $\phi\xi_t$ factors as \[\nicexy{\Cor_t \ar[r]^-{\xi_t} \ar[d]_-{\phi_e} & T \ar[d]^-{\phi}\\ V_e \ar[r]^-{\varphi_e} & V}\] in which $\varphi_e$ is a composite of outer cofaces.   Since $x_{V_e}\in X(V_e)$ is the unique extension of the pair of elements $\{x_{e^0},x_{e^1}\} \subseteq \{x_v\}_{v\in\Vt(T)}$, we have \[x_{V_e} = X(\varphi_e)(z).\]  Therefore, we have \[\begin{split}x_t &= X(\phi_e)(x_{V_e})\\
&= X(\phi_e)X(\varphi_e)(z)\\
&= X(\xi_t)X(\phi)(z).
\end{split}\]
\end{enumerate}
We have shown that $X(\phi)(z)=x_T$, which proves the uniqueness of $y$.
\end{proof}

Combining Theorem \ref{segal-is-strict} and Theorem \ref{strict-is-segal}, we obtain the following characterization of strict $\infty$-$\G$-operads.

\begin{corollary}\label{segal=strict}
A $\Treecatg$-presheaf is a strict $\infty$-$\G$-operad if and only if it satisfies the $\G$-Segal condition.
\end{corollary}

\begin{remark}
The categorical analogue of Corollary \ref{segal=strict} is a basic fact in $\infty$-category theory.  Namely, a simplicial set is a strict quasi-category (i.e., every inner horn has a unique extension) if and only if it satisfies the Segal condition.  See, for example, Proposition 1.1.2.2 in \cite{lurie}.\dqed
\end{remark}

\begin{example}[Strict $\infty$-Braided/Ribbon/Cactus Operads]
When $\G$ is the braid group operad $\B$, Corollary \ref{segal=strict} says that a $\Treecatb$-presheaf is a strict $\infty$-braided operad if and only if it satisfies the braided Segal condition.  Similarly, when $\G$ is the ribbon group operad $\R$, a $\Treecatr$-presheaf is a strict $\infty$-ribbon operad if and only if it satisfies the ribbon Segal condition.  Also, when $\G$ is the cactus group operad $\Cac$, a $\Treecatcac$-presheaf is a strict $\infty$-cactus operad if and only if it satisfies the cactus Segal condition.   There are also pure analogues for the pure braid group operad $\PB$, the pure ribbon group operad $\PR$, and the pure cactus group operad $\PCac$.\dqed
\end{example}

\section{Characterization of the Nerve}\label{sec:nerve-characterization}

The purpose of this section is to prove the following characterizations of the $\G$-nerve for each action operad $\G$.

\begin{theorem}\label{nerve-characterization-theorem}
For\index{nerve!characterizations} each $X \in \Settotreecatgop$, the following statements are equivalent:
\begin{enumerate}
\item $X$ is isomorphic to the $\G$-nerve of some $\G$-operad in $\Set$.
\item $X$ satisfies the $\G$-Segal condition.
\item $X$ is a strict $\infty$-$\G$-operad.
\end{enumerate}
\end{theorem}

\begin{proof}
The implications (1) $\Rightarrow$ (2) $\Leftrightarrow$ (3) have already been established in Theorem \ref{nerve-is-segal} and Corollary \ref{segal=strict}.  It remains to establish the implication (3) $\Rightarrow$ (1), so suppose $X$ is a strict $\infty$-$\G$-operad.  We define an associated $\G$-operad $(\colorc,\O)$ in $\Set$ and observe that its $\G$-nerve is isomorphic to $X$.

First we define the color set \[\colorc = X\bigl(\uparrow,\id_1\bigr).\]  For each pair $\czerouc \in \Profcc$ with $\uc=(c_1,\ldots,c_n)$, we define the set \[\O\czerouc = \coprodover{\sigma\in\G(n)}\, \Bigl\{x \in X\bigl(\Cor_n,\sigma) :  X(\phi_0)(x) = c_0,\, X(\phi_{\sigmabar(j)})(x) = c_j \text{ for $1\leq j \leq n$}\Bigr\},\] where the morphisms \[\nicexy{(\uparrow,\id_1) \ar[r]^-{\phi_j} & (\Cor_n,\sigma)} \in \Treecatg\] for $0 \leq j \leq n$ are the outer cofaces in Example \ref{ex:outer-coface-corolla}.  Next we define a $\colorc$-colored $\G$-operad structure on $\O$ as in Proposition \ref{prop:g-operad-compi}.

The $\G$-sequence structure on $\O$ is induced by $X$ applied to the isomorphisms
\[\nicexy@C+2.5cm{(\Cor_n,\sigma) \ar[r]^-{\bigl(((\Cor_n,\tau),u) ; (\Cor_n,\id_n)\bigr)} & (\Cor_n,\sigma\tau)}\in \Treecatg\] for $\sigma,\tau\in \G(n)$.

For each $c\in X(\uparrow,\id_1)$, the $c$-colored unit in $\O$ is the element \[\operadunit_c = X(\phi)(c) \in X(\Cor_1,\id_1)\] with $\phi \in \Treecatg$ the codegeneracy \[\nicexy@C+2.4cm{(\Cor_1,\id_1) \ar[r]^-{\bigl(((\Cor_1,\id_1),u); (\uparrow,\id_1)\bigr)} & (\uparrow,\id_1)}.\]

To define the planar operad $\compi$-composition, first observe that it is enough to consider elements in $X(\Cor_n,\id_n)$, since the $\compi$-composition of general elements in $X(\Cor_n,\sigma)$ are then determined by the equivariance axiom in Proposition \ref{prop:g-operad-compi}.  So suppose given elements \[x \in X(\Cor_n,\id_n) \cap \O\duc \andspace y \in X(\Cor_m,\id_m) \cap \O\ciub\] for some $1 \leq i \leq n$.  Consider the planar tree $T$
\begin{center}\begin{tikzpicture}
\node[plain] (u) {$u$}; \node[plain, below=.8 of u] (v) {$v$}; 
\node[below=.08 of v] () {$\cdots$};
\draw[outputleg] (u) to +(0,.8); \draw[thick] (v) to node{\scriptsize{$i$}} (u);
\draw[thick] (u) to node[swap]{\scriptsize{$1$}} +(-.6,-.6); 
\draw[thick] (u) to node{\scriptsize{$n$}} +(.6,-.6);
\draw[thick] (v) to node[swap]{\scriptsize{$1$}} +(-.6,-.6); 
\draw[thick] (v) to node{\scriptsize{$m$}} +(.6,-.6);
\end{tikzpicture}\end{center}
with one internal edge, which is the $i$th input of the vertex $u$, $n$ inputs at $u$, and $m$ inputs at $v$.  Since \[X(\phi_i)(x) = c_i = X(\phi_0)(y) \in \colorc,\] it follows from Lemma \ref{horn-description} that $x$ and $y$ determine an inner horn $f$ of $X$,
\[\nicexy{\horn{\phi}(T,\id_{n+m-1}) \ar[r]^-{f} \ar@{>->}[d] & X\\ \Yonedag(T,\id_{n+m-1}) \ar@{-->}[ur]_-{g} &}\] in which $\phi \in \Treecatg$ is the inner coface \[\nicexy@C+4cm{(\Cor_{n+m-1}, \id_{n+m-1}) \ar[r]^-{\bigl(((\Cor_{n+m-1}, \id_{n+m-1}),u); (T,\id_{n+m-1})\bigr)} & (T,\id_{n+m-1})}.\]  Since $X$ is a strict $\infty$-$\G$-operad, there exists a unique dashed arrow extension $g$, which by Yoneda Lemma corresponds to an element $z \in X(T,\id_{n+m-1})$.  We now define \[x\compi y = X(\phi)(z) \in X(\Cor_{n+m-1}, \id_{n+m-1}).\]  Using the properties of planar tree substitution in Lemma \ref{lem:planar-tree-sub-vertex}, one can check that $x\compi y$ is well-defined and that $\O$ is a $\colorc$-colored $\G$-operad in $\Set$.

We now claim that there is an isomorphism \[X\cong \Nerveg(\colorc,\O).\]  Since $X$ is a strict $\infty$-$\G$-operad, it satisfies the $\G$-Segal condition by Theorem \ref{strict-is-segal}.  The $\G$-nerve $\Nerveg(\colorc,\O)$ also satisfies the $\G$-Segal condition by Theorem \ref{nerve-is-segal}.  Therefore, it is enough to establish the desired isomorphism on a $\G$-tree $(V,\sigma) \in \Treecatg$ in which the planar tree $V$ is either $\uparrow$ or a corolla.  In this case, the desired isomorphism holds by the construction of $\O$.  
\end{proof}

\chapter{Coherent Realization-Nerve and Infinity Group Operads}
\label{ch:hc-nerve}

Fix an action operad $\G$ as in Definition \ref{def:augmented-group-operad} and a commutative segment $(J,\mu,0,1,\epsilon)$ in the symmetric monoidal category $(\M,\otimes,\tensorunit)$ as in Definition \ref{def:segment}. The purpose of this chapter is to study a version of the $\G$-realization-nerve adjunction, called the \emph{coherent} $\G$-realization-nerve adjunction, that is based on the symmetric monoidal category $\M$ rather than $\Set$.  The $\G$-Boardman-Vogt construction, defined using the chosen commutative segment $J$, will play a crucial role.

In Section \ref{sec:hc-nerve} we define the coherent $\G$-realization functor $\Realghc$ and the coherent $\G$-nerve functor $\Nerveghc$.  Their definitions are similar to those of the $\G$-realization and the $\G$-nerve functors, with the $\G$-operad $\G^{(T,\sigma)}$ replaced by its $\G$-Boardman-Vogt construction $\Wg\iota_*\G^{(T,\sigma)}$ in $\M$.  Here $\iota : \Set \to\M$ is the canonical strong symmetric monoidal functor.  When $\G$ is the symmetric group operad $\S$, the coherent $\S$-nerve $\Nerveshc$ recovers the homotopy coherent dendroidal nerve of Moerdijk-Weiss \cite{mw09}.  While the coherent $\G$-nerve is defined using the $\G$-Boardman-Vogt construction, we show that it can also be computed using the planar Boardman-Vogt construction.

In Theorem \ref{nerveg-fully-faithful} we showed that the $\G$-realization of the $\G$-nerve is naturally isomorphic to the identity.  In Section \ref{sec:hcr-nerve} we obtain a coherent version of this result by showing that the coherent $\G$-realization of the $\G$-nerve is naturally isomorphic to the $\G$-Boardman-Vogt construction.  Using the result in the previous section, in Section \ref{sec:nerve-to-hcnerve} we compute the morphisms and the internal hom object from the $\G$-nerve of a $\G$-operad to the coherent $\G$-nerve of another $\G$-operad.

In Theorem \ref{gnerve-changeofg} we observed that the $\G$-realization and the $\G$-nerve functors behave nicely with respect to a change of the action operad.  In Section \ref{sec:hcr-hcn-change-goperad} we obtain a coherent version of this result for the coherent $\G$-realization functor and the coherent $\G$-nerve functor.

In Section \ref{sec:ex-bv} we provide an explicit description of the planar Boardman-Vogt construction of the planar operad freely generated by a planar tree.  This is a preliminary observation that is needed in the following section.  Recall from Corollary \ref{nerve-is-strict} that the $\G$-nerve of each $\G$-operad in $\Set$ is a strict $\infty$-$\G$-operad.  In Section \ref{sec:hcnerve-infinity-goperad} we obtain a coherent version of this result by showing that the coherent $\G$-nerve of each entrywise fibrant $\G$-operad is an $\infty$-$\G$-operad.

\section{Coherent Realization and Coherent Nerve}\label{sec:hc-nerve}

Before we define the coherent realization and the coherent nerve, let us first recall the following constructions.
\begin{enumerate}
\item  There is a diagram \[\nicexy{\Treecatg\ar[r]^-{\G^-} \ar[d]_-{\Yonedag} & \gopset\\ \Settotreecatgop &}\]in which:
\begin{itemize}
\item $\G^-$ is the full and faithful functor in Theorem \ref{gtree-gop-functor}.  For each $\G$-tree $(T,\sigma)$, $\G^{(T,\sigma)}$ is the $\Ed(T)$-colored $\G$-operad in $\Set$ freely generated by the set of vertices in $T$.
\item $\Yonedag$ is the Yoneda embedding in Definition \ref{def:gnerve} of the $\G$-tree category $\Treecatg$.
\end{itemize}
\item There is a strong symmetric monoidal functor $\iota : \Set \to \M$ that sends each set $X$ to $\iota(X) = \coprod_X \tensorunit$.  By Theorem \ref{change-base-category}(2) there is an induced functor \[\nicexy{\gopset \ar[r]^-{\iota_*} & \gopm}\] from the category of $\G$-operads in $\Set$ to the category of $\G$-operads in $\M$, which is entrywise the functor $\iota$.  For each $\O\in\gopset$, its image $\iota_*\O$ is also denoted by $\subm{\O}$.
\item By Theorem \ref{wg-augmentation} the $\G$-Boardman-Vogt construction gives a functor \[\nicexy{\gopm \ar[r]^-{\Wg} & \gopm}\] that is naturally augmented over the identity functor.  For each $\colorc$-colored $\G$-operad $\Q$ in $\M$, its $\G$-Boardman-Vogt construction $\Wg\Q$ is also a $\colorc$-colored $\G$-operad in $\M$ and is entrywise defined as the coend in \eqref{wgoduc}, where the chosen commutative segment in $\M$ is used.
\end{enumerate}

We are now ready to define the coherent versions of the $\G$-realization and the $\G$-nerve functors.

\begin{definition}\label{def:hcrealization-hcnerve}
Define the\index{realization!coherent}\index{coherent G-realization@coherent $\G$-realization} \emph{coherent $\G$-realization} $\Realghc$ and the\index{nerve!coherent}\index{coherent G-nerve@coherent $\G$-nerve} \emph{coherent $\G$-nerve} $\Nerveghc$,
\begin{equation}\label{hcr-hcn}
\nicexy@C-.3cm{\Treecatg \ar[d]_-{\Yonedag} \ar[r]^-{\G^-} & \gopset \ar[r]^-{\iota_*} & \gopm \ar[d]^-{\Wg}\\
\Settotreecatgop \ar@<2pt>[rr]^-{\Realghc} && \gopm \ar@<2pt>[ll]^-{\Nerveghc}}
\end{equation}
as follows.
\begin{enumerate}
\item $\Realghc$ is the Yoneda extension of the composite \[\nicexy@C+2cm{\Treecatg \ar[r]^-{\Wg\iota_*\G^- \,=\, \Wg\gdashsubm} & \gopm},\]i.e., the left Kan extension of $\Wg\gdashsubm$ along the Yoneda embedding $\Yonedag$.  It exists by Proposition \ref{gopm-bicomplete}.  There are coend formulas \[\begin{split}
\Realghc(X) &\cong \int^{(T,\sigma)\in\Treecatg} \Bigl[\Settotreecatgop\bigl(\Yonedag(T,\sigma),X\bigr)\Bigr] \cdot \Wg\gsubm^{(T,\sigma)}\\
&\cong \int^{(T,\sigma)\in\Treecatg} \coprodover{X(T,\sigma)} \Wg\gsubm^{(T,\sigma)} \in \gopm
\end{split}\]
for $X \in \Settotreecatgop$. 
\item $\Nerveghc$ is the functor such that \[\Nerveghc(\colorc,\O) \in \Settotreecatgop\] is the $\Treecatg$-presheaf defined by \[\Nerveghc(\colorc,\O)(T,\sigma) = \gopm\Bigl(\Wg\gsubm^{(T,\sigma)},(\colorc,\O)\Bigr)\] for $(\colorc,\O)\in\gopm$ and $(T,\sigma)\in\Treecatg$.  
\end{enumerate}
\end{definition}

In the rest of this chapter, we study the coherent $\G$-realization and the coherent $\G$-nerve.  First we observe that they form an adjunction.

\begin{lemma}\label{hcr-hcn-adjunction}
$\Realghc$ is the left adjoint of $\Nerveghc$.
\end{lemma}

\begin{proof}
This is almost identical to the proof of Lemma \ref{nerve-real-adjunction} with $\gopset$ and $\G^{(T,\sigma)}$ replaced by $\gopm$ and $\Wg\gsubm^{(T,\sigma)}$, respectively.
\end{proof}

We call $\Realghc \dashv \Nerveghc$ the \emph{coherent $\G$-realization-nerve adjunction}.  

\begin{example}[Coherent Planar and Symmetric Realization-Nerve]\label{ex:ps-hcrn}
When $\G$ is the symmetric group operad $\S$ in Example \ref{ex:symmetric-group-operad}, the symmetric tree category $\Treecats$ is equivalent to the Moerdijk-Weiss symmetric dendroidal category $\Omega$, as mentioned in Example \ref{ex:treecats}.  In the setting of $\Omega$, the coherent $\S$-nerve is called the\index{homotopy coherent dendroidal nerve}\index{dendroidal nerve!homotopy coherent} \emph{homotopy coherent dendroidal nerve} and is denoted by $hcN_d$ in \cite{mw09}.  Similarly, when $\G$ is the planar group operad $\P$ in Example \ref{ex:trivial-group-operad}, the planar tree category $\Treecatp$ is isomorphic to the Moerdijk-Weiss planar dendroidal category $\Omega_p$.  We have the analoguous coherent $\P$-realization-nerve adjunction\label{not:realphc} $\Realhc^{\P} \dashv \Nervehc^{\P}$.  However, the reader is reminded that the Boardman-Vogt construction in \cite{berger-moerdijk-bv} is a sequential colimit of pushouts, with each map dependent on the previous inductive stage.  On the other hand, our $\G$-Boardman-Vogt construction is defined in one step as a coend indexed by the substitution category.\dqed
\end{example}

\begin{example}[Coherent Braided Realization-Nerve]\label{ex:braided-hcrn}
For the braid group operad $\B$ in Definition \ref{def:braid-group-operad}, we call $\Realhc^{\B}$ and $\Nervehc^{\B}$ the\index{coherent braided realization}\index{braided realization!coherent} \emph{coherent braided realization} and the\index{coherent braided nerve}\index{braided nerve!coherent} \emph{coherent braided nerve}.  Similarly, for the pure braid group operad $\PB$ in Definition \ref{def:pure-braid-group-operad}, we call $\Realhc^{\PB}$ and $\Nervehc^{\PB}$ the \emph{coherent pure braided realization} and the \emph{coherent pure braided nerve}.\dqed
\end{example}

\begin{example}[Coherent Ribbon Realization-Nerve]\label{ex:ribbon-hcrn}
For the ribbon group operad $\R$ in Definition \ref{def:ribbon-group-operad}, we call $\Realhc^{\R}$ and $\Nervehc^{\R}$ the\index{coherent ribbon realization}\index{ribbon realization!coherent} \emph{coherent ribbon realization} and the\index{coherent ribbon nerve}\index{ribbon nerve!coherent} \emph{coherent ribbon nerve}.  Similarly, for the pure ribbon group operad $\PR$ in Definition \ref{def:pure-ribbon-group-operad}, we call $\Realhc^{\PR}$ and $\Nervehc^{\PR}$ the \emph{coherent pure ribbon realization} and the \emph{coherent pure ribbon nerve}.\dqed
\end{example}

\begin{example}[Coherent Cactus Realization-Nerve]\label{ex:cactus-hcrn}
For the cactus group operad $\Cac$ in Definition \ref{def:cactus-group-operad}, we call $\Realhc^{\Cac}$ and $\Nervehc^{\Cac}$ the\index{coherent cactus realization}\index{cactus realization!coherent} \emph{coherent cactus realization} and the\index{coherent cactus nerve}\index{cactus nerve!coherent} \emph{coherent cactus nerve}.  Similarly, for the pure cactus group operad $\PCac$ in Definition \ref{def:pure-cactus-group-operad}, we call $\Realhc^{\PCac}$ and $\Nervehc^{\PCac}$ the \emph{coherent pure cactus realization} and the \emph{coherent pure cactus nerve}.\dqed
\end{example}

Next we compute the coherent $\G$-nerve of a $\G$-operad in terms of the planar Boardman-Vogt construction.  The functors involved in the next result are in the diagram
\begin{equation}\label{hcnerve-computation-planar}
\nicexy{\Treecatp \ar[d]_-{\P^-} & \Treecatg \ar[d]^-{\G^-}\\
\popset \ar[d]_-{\iota_*} \ar[r]^-{\iota^{\G}_!} & \gopset \ar[d]^-{\iota_*}\\
\Poperadm \ar[d]_-{\Wp} \ar[r]^-{\iota^{\G}_!} & \gopm \ar[d]^-{\Wg}\\
\Poperadm \ar@<2pt>[r]^-{\iota^{\G}_!} & \gopm \ar@<2pt>[l]^-{(\iota^{\G})^*}}
\end{equation}
in which:
\begin{itemize}
\item $\P$ is the planar group operad in Example \ref{ex:trivial-group-operad}.
\item $\P^-$ and $\G^-$ are the functors in Theorem \ref{gtree-gop-functor}.
\item $\iota^{\G} : \P\to\G$ is the unique morphism from the planar group operad $\P$ to $\G$.
\item The pair $\iota^{\G}_! \dashv (\iota^{\G})^*$ is the adjunction in Theorem \ref{goneopm-gtwoopm} induced by the morphism $\iota^{\G}$.
\item $\Wp$ and $\Wg$ are the planar and $\G$-Boardman-Vogt constructions, respectively.
\item Each $\iota_*$ is the change-of-base functor in Theorem \ref{change-base-category}(2) induced by the strong symmetric monoidal functor $\iota : \Set \to \M$ defined by $\iota(X)=\coprod_X \tensorunit$.
\end{itemize}

\begin{proposition}\label{hcn-planar}
There is a natural bijection \[\bigl(\Nerveghc(\O)\bigr)(T,\sigma) \cong \Poperadm\Bigl(\Wp\iota_*\P^T, (\iota^{\G})^*\O\Bigr)\] for each $\colorc$-colored $\G$-operad $\O$ in $\M$ and each $\G$-tree $(T,\sigma)$.
\end{proposition}

\begin{proof}
Since the functor $\iota : \Set\to\M$ commutes with colimits, by Theorem \ref{varphi-F}(2) the middle square of the diagram \eqref{hcnerve-computation-planar} is commutative up to a natural isomorphism.  Moreover, in the bottom square of the diagram \eqref{hcnerve-computation-planar}, there is a natural isomorphism \[\iota^{\G}_!\Wp \cong \Wg\iota^{\G}_!\] by Theorem \ref{thm:w-left-adjoint}.  Combining the middle and the bottom squares of the diagram \eqref{hcnerve-computation-planar}, there is a natural isomorphism
\begin{equation}\label{w-iota}
\Wg\iota_*\iota^{\G}_! \cong \iota^{\G}_!\Wp\iota_*.
\end{equation}

The desired natural bijection follows from the following computation:
\[\begin{split}
\bigl(\Nerveghc(\O)\bigr)(T,\sigma) 
&= \gopm\Bigl(\Wg\iota_*\G^{(T,\sigma)},\O\Bigr)\\
&\cong \gopm\Bigl(\Wg\iota_*\iota^{\G}_!\P^T,\O\Bigr)\\
&\cong \gopm\Bigl(\iota^{\G}_!\Wp\iota_*\P^T,\O\Bigr)\\
&\cong \Poperadm\Bigl(\Wp\iota_*\P^T,(\iota^{\G})^*\O\Bigr).
\end{split}\]
The first isomorphism is from Corollary \ref{ptree-free-goperad}.  The second isomorphism is from \eqref{w-iota}.  The last isomorphism is from the adjunction $\iota^{\G}_! \dashv (\iota^{\G})^*$.
\end{proof}

\begin{interpretation}
In Proposition \ref{hcn-planar}:
\begin{enumerate}
\item The $\Ed(T)$-colored planar operad $\Wp\iota_*\P^T$ is the planar Boardman-Vogt construction of $\iota_*\P^T$, which is the image in $\M$ of the $\Ed(T)$-colored planar operad $\P^T$ in $\Set$ freely generated by the set of vertices in $T$.  In Theorem \ref{planar-bv} below, we will compute the planar operad $\Wp\iota_*\P^T$ explicitly.  As we will show there, each entry of $\Wp\iota_*\P^T$ is either
\begin{itemize}
\item the initial object in $\M$ if the corresponding entry of $\P^T$ is empty, or
\item $J^{\otimes |\Int(Q)|}$ if the corresponding entry of $\P^T$ is the singleton $\{Q\}$.
\end{itemize}
\item The object $(\iota^{\G})^*\O$ is the underlying $\colorc$-colored planar operad of the $\colorc$-colored $\G$-operad $\O$.  It is obtained from $\O$ by forgetting its $\G$-equivariant structure.  Its entries, colored operadic units, and $\compi$-composition remain the same.\dqed
\end{enumerate}
\end{interpretation}

\section{Coherent Realization of the Nerve}\label{sec:hcr-nerve}

The purpose of this section is to show that the coherent $\G$-realization of the $\G$-nerve is the $\G$-Boardman-Vogt construction.  In the following main result of this section, we refer to some of the functors in the diagram \eqref{hcr-hcn}.

\begin{theorem}\label{coherent-realization-nerve}
The diagram\index{nerve!coherent realization of}\index{coherent G-realization@coherent $\G$-realization!of nerve}
\begin{equation}\label{creal-nerve-wg}
\nicexy{\gopset \ar[d]_-{\Nerveg} \ar[r]^-{\iota_*} & \gopm \ar[d]^-{\Wg}\\
\Settotreecatgop \ar[r]^-{\Realghc} & \gopm}
\end{equation}
is commutative up to a natural isomorphism.  In other words, there is a natural isomorphism \[\Realghc \Nerveg(\O) \cong \Wg\iota_*\O \in \gopm\] for $\O\in\gopset$.
\end{theorem}

\begin{proof}
Pick an arbitrary $\colorc$-colored operad $\O$ in $\Set$.  Recall from \eqref{nerve-colim-representables} that the $\G$-nerve of $\O$ is naturally isomorphic to a colimit of representables
\[\nicexy{\colimover{\bigl((T,\sigma); \phi\bigr)\in\Govero} \Yonedag(T,\sigma) \ar[r]^-{\cong} & \Nerveg(\O)}\in \Settotreecatgop.\] The indexing category $\Govero$ of the colimit on the left is defined in Definition \ref{def:goperad-colim}.  An object in $\Govero$ is a pair $\bigl((T,\sigma);\phi\bigr)$ with
\begin{itemize}
\item $(T,\sigma)$ a $\G$-tree and 
\item $\phi : \G^{(T,\sigma)} \to \O$ a morphism of $\G$-operads in $\Set$. 
\end{itemize} 
Such a morphism $\phi$ may be regarded as an $\O$-decoration of the planar tree $T$, as discussed in \ref{int:decoration-of-tree}.

There is a natural isomorphism 
\begin{equation}\label{hcreal-nerve-colimt}
\Realghc \Nerveg(\O) \cong \colimover{\bigl((T,\sigma); \phi\bigr)\in\Govero} \Wg\iota_*\G^{(T,\sigma)}.
\end{equation}
Indeed, we have the natural isomorphisms
\[\begin{split}
\Realghc \Nerveg(\O) &\cong \Realghc\Bigl(\colimover{\bigl((T,\sigma); \phi\bigr)\in\Govero} \Yonedag(T,\sigma)\Bigr)\\
&\cong \colimover{\bigl((T,\sigma); \phi\bigr)\in\Govero} \Realghc\Yonedag(T,\sigma)\\
&\cong \colimover{\bigl((T,\sigma); \phi\bigr)\in\Govero} \Wg\iota_*\G^{(T,\sigma)}.
\end{split}\]
The first isomorphism is the coherent $\G$-realization applied to the isomorphism in \eqref{nerve-colim-representables}.  The second isomorphism uses the fact that 
the coherent $\G$-realization commutes with colimits, which holds because it is a left adjoint by Lemma \ref{hcr-hcn-adjunction}.  The third isomorphism comes from the definition of the coherent $\G$-realization as the Yoneda extension of the composite $\Wg\iota_*\G^-$ in the diagram \eqref{hcr-hcn}.  

The desired isomorphism now follows from \eqref{hcreal-nerve-colimt} and the coend definition of the $\G$-Boardman-Vogt construction in Section \ref{sec:boardman-vogt}.  More precisely, for each object $\bigl((T,\sigma);\phi\bigr)$ in $\Govero$, by the functoriality of the $\G$-Boardman-Vogt construction in Theorem \ref{wg-augmentation}(1), the morphism \[\nicexy{\G^{(T,\sigma)}\ar[r]^-{\phi} & \O} \in \gopset\] yields a morphism \[\nicexy@C+.5cm{\Wg\iota_*\G^{(T,\sigma)} \ar[r]^-{\Wg\iota_*\phi} & \Wg\iota_*\O} \in \gopm.\]  By the universal property of a colimit, these morphisms induce a morphism \[\nicexy{\Realghc \Nerveg(\O) \cong \colimover{\bigl((T,\sigma); \phi\bigr)\in\Govero} \Wg\iota_*\G^{(T,\sigma)} \ar[r]^-{A} & \Wg\iota_*\O}\in \gopm\] whose restriction to each object $\Wg\iota_*\G^{(T,\sigma)}$ is the morphism $\Wg\iota_*\phi$.  To see that the morphism $A$ is the desired isomorphism, we construct its inverse explicitly.

Observe that for each $\colorc$-colored $\G$-tree $(T,\sigma)$ with profile $\duc \in\Profcc$, there are isomorphisms 
\[\begin{split}
(\iota_*\O)(T,\sigma) &= \bigtensorover{t\in\Vt(T)} (\iota_*\O)\inoutt\\
&= \bigtensorover{t\in\Vt(T)} \Bigl(\coprodover{\O\inoutt} \tensorunit\Bigr)\\
&\cong \coprodover{\prodover{t\in\Vt(T)} \O\inoutt} \bigl(\bigotimes \tensorunit\bigr)\\
&\cong \coprodover{\prodover{t\in\Vt(T)} \O\inoutt} \tensorunit.
\end{split}\]  
Similarly, there is an isomorphism \[\bigl(\iota_*\G^{(T,\sigma)}\bigr)(T,\sigma) \cong \coprodover{\prodover{t\in\Vt(T)} \G^{(T,\sigma)}\inoutt} \tensorunit.\]  
Furthermore, each element $\{x_t\}\in\prod_{t\in\Vt(T)}\O\inoutt$ corresponds to a morphism \[\phi : \G^{(T,\sigma)}\to\O \in \gopset\] determined by the assignment \[\nicexy{\Vt(T) \ni t \ar@{|->}[r]^-{\phi} & x_t}\in \O\inoutt.\]  

This yields a morphism $B$ determined by the commutative diagram
\[\nicexy@C+1.5cm{\bigl(\Wg\iota_*\O\bigr)\duc \ar[r]^-{B} & \Bigl(\colimover{((T,\sigma);\phi)} \Wg\iota_*\G^{(T,\sigma)}\Bigr)\inout{T}\\
\J(T,\sigma)\otimes (\iota_*\O)(T,\sigma) \ar[u]^-{\eta_{(T,\sigma)}} & \bigl(\Wg\iota_*\G^{(T,\sigma)}\bigr)\inout{T} \ar[u]^-{((T,\sigma);\phi)}_-{\text{natural}}\\
\coprodover{\prodover{t\in\Vt(T)} \O\inoutt} \J(T,\sigma)\otimes \tensorunit \ar[u]^-{\cong} & \J(T,\sigma)\otimes \bigl(\iota_*\G^{(T,\sigma)}\bigr)(T,\sigma) \ar[u]_-{\eta_{(T,\sigma)}}\\
\J(T,\sigma)\otimes\tensorunit \ar[u]^-{\{x_t\}_{t\in\Vt(T)}}_-{\text{summand}} \ar[r]^-{\{(\Cor_t,\id_t)\}_{t\in\Vt(T)}}_-{\text{summand}} & \coprodover{\prodover{t\in\Vt(T)} \G^{(T,\sigma)}\inoutt} \J(T,\sigma)\otimes\tensorunit \ar[u]_-{\cong}.}\]
One can check that this defines a morphism \[\nicexy{\Wg\iota_*\O \ar[r]^-{B} & \colimover{\bigl((T,\sigma); \phi\bigr)\in\Govero} \Wg\iota_*\G^{(T,\sigma)} \cong \Realghc \Nerveg(\O)}\in \gopm\]that is inverse to the morphism $A$.
\end{proof}

\begin{example}[Coherent Symmetric Realization]
When\index{coherent realization!symmetric} $\G$ is the symmetric group operad $\S$ in Example \ref{ex:symmetric-group-operad}, the $\S$-tree category $\Treecats$ is equivalent to the symmetric dendroidal category $\Omega$ of Moerdijk-Weiss \cite{mw07,mw09}.  In this case, Theorem \ref{coherent-realization-nerve} recovers a result in Section 6.3 in \cite{mt}.\dqed
\end{example}

\begin{example}[Coherent Braided Realization]
When $\G$ is the braid group operad $\B$ in Definition \ref{def:braid-group-operad}, Theorem \ref{coherent-realization-nerve} says that the\index{coherent braided realization}\index{coherent realization!braided} coherent braided realization of the \index{braided nerve}\index{nerve!braided}braided nerve of a braided operad in $\Set$ is naturally isomorphic to its braided Boardman-Vogt construction.\dqed
\end{example}

\begin{example}[Coherent Ribbon Realization]
When $\G$ is the ribbon group operad $\R$ in Definition \ref{def:ribbon-group-operad}, Theorem \ref{coherent-realization-nerve} says that the\index{coherent ribbon realization}\index{coherent realization!ribbon} coherent ribbon realization of the\index{ribbon nerve}\index{nerve!ribbon} ribbon nerve of a ribbon operad in $\Set$ is naturally isomorphic to its ribbon Boardman-Vogt construction.\dqed
\end{example}

\begin{example}[Coherent Cactus Realization]
When $\G$ is the cactus group operad $\Cac$ in Definition \ref{def:cactus-group-operad}, Theorem \ref{coherent-realization-nerve} says that the\index{coherent cactus realization}\index{coherent realization!cactus} coherent cactus realization of the\index{cactus nerve}\index{nerve!cactus} cactus nerve of a cactus operad in $\Set$ is naturally isomorphic to its cactus Boardman-Vogt construction.\dqed
\end{example}

\section{Nerve to Coherent Nerve}\label{sec:nerve-to-hcnerve}

The purpose of this section is to compute
\begin{enumerate}[label=(\roman*)]
\item the set of morphisms and 
\item the internal hom object 
\end{enumerate}
from the $\G$-nerve to the coherent $\G$-nerve.  The first result below says that the set of morphisms from the $\G$-nerve to the coherent $\G$-nerve can be described in terms of the $\G$-Boardman-Vogt construction.  We continue to use the functors in the diagram \eqref{hcr-hcn}.  

\begin{corollary}\label{nerve-hcnerve}
There is a natural bijection \[\Settotreecatgop\Bigl(\Nerveg(\O),\Nerveghc(\iota_*\Q)\Bigr) \cong \gopm\bigl(\Wg\iota_*\O,\iota_*\Q\bigr)\] for $\O,\Q\in\gopset$.
\end{corollary}

\begin{proof}
This follows from (i) the coherent $\G$-realization-nerve adjunction in Lemma \ref{hcr-hcn-adjunction} and (ii) the natural isomorphism \[\Realghc \Nerveg(\O) \cong \Wg\iota_*\O\] in Theorem \ref{coherent-realization-nerve}.
\end{proof}

Recall from Theorem \ref{thm:presheaf-monoidal} that the presheaf category $\Settotreecatgop$ has a symmetric monoidal closed structure with internal hom object $\Homsub{\G}$.  Also recall from Theorem \ref{goperad-symmetric-monoidal} that the category $\gopset$ has a symmetric monoidal structure denoted by $\tensorg$.  The next observation strengthens the previous result to the internal hom object from the $\G$-nerve to the coherent $\G$-nerve.

\begin{theorem}\label{internal-hom-nerve-hcnerve}
There is a natural bijection \[\Homsub{\G}\Bigl(\Nerveg(\O),\Nerveghc(\iota_*\Q)\Bigr)(T,\sigma) \cong \gopm\Bigl(\Wg\iota_* \bigl(\O\tensorg \G^{(T,\sigma)}\bigr),\iota_*\Q\Bigr)\] for $\O,\Q\in\gopset$ and $(T,\sigma)\in\Treecatg$.
\end{theorem}

\begin{proof}
We compute as follows:
\[\begin{split}
&\Homsub{\G}\Bigl(\Nerveg(\O),\Nerveghc(\iota_*\Q)\Bigr)(T,\sigma)\\
&= \Settotreecatgop\Bigl(\Nerveg(\O)\tensorsub{\G} \Yonedag(T,\sigma), \Nerveghc(\iota_*\Q)\Bigr)\\
&\cong \Settotreecatgop\Bigl(\Nerveg\bigl(\O\tensorg \G^{(T,\sigma)}\bigr), \Nerveghc(\iota_*\Q)\Bigr)\\
&\cong \gopm\Bigl(\Wg\iota_*\bigl(\O\tensorg \G^{(T,\sigma)}\bigr),\iota_*\Q\Bigr).
\end{split}\]
The equality is the definition of the internal hom object $\Homsub{\G}$ in Definition \ref{def:treecatg-presheaf-monoidal}.  The first isomorphism uses the second isomorphism in Proposition \ref{representable-tensor}.  The last isomorphism follows from Corollary \ref{nerve-hcnerve}.
\end{proof}

\begin{example}
When $\G$ is the symmetric group operad $\S$ in Example \ref{ex:symmetric-group-operad},  Corollary \ref{nerve-hcnerve} and Theorem \ref{internal-hom-nerve-hcnerve} recover some results in Section 6.3 in \cite{mt} about the (homotopy coherent) dendroidal nerve.\dqed
\end{example}

\section{Change of Action Operads}\label{sec:hcr-hcn-change-goperad}

The purpose of this section is to show that the coherent $\G$-realization functor and the coherent $\G$-nerve functor are well-behaved with respect to a change of the action operad.  Fix a morphism $\varphi : \Gone\to\Gtwo$ of action operads as in Definition \ref{def:augmented-group-operad-morphism}.  Consider the diagram
\begin{equation}\label{hcnervegone-hcnervegtwo}
\nicexy@C+1.4cm{\Settotreecatgoneop \ar@<2pt>[r]^-{\varphi_!} \ar@<-2pt>[d]_-{\Realghcone} & \Settotreecatgtwoop \ar@<-2pt>[d]_-{\Realghctwo} \ar@<2pt>[l]^-{\varphi^*}\\
\goneopm \ar@<-2pt>[u]_-{\Nerveghcone} \ar@<2pt>[r]^-{\varphi_!} & \gtwoopm \ar@<2pt>[l]^-{\varphi^*} \ar@<-2pt>[u]_-{\Nerveghctwo}}
\end{equation}
in which:
\begin{itemize}
\item The top functor $\varphi^*$ is the pullback induced by \[(\Treecat^{\varphi})^{\op} : \Treecatgoneop \to \Treecatgtwoop\] with $\Treecat^{\varphi}$ the functor in Proposition \ref{treecatg-functor}(1).
\item The top $\varphi_!$ is the left adjoint of $\varphi^*$ given by a left Kan extension.
\item The bottom adjunction $\varphi_!\dashv\varphi^*$ is the one in Theorem \ref{goneopm-gtwoopm}.
\item $\Realghcone \dashv\Nerveghcone$ is the coherent $\Gone$-realization-nerve adjunction.
\item $\Realghctwo \dashv\Nerveghctwo$ is the coherent $\Gtwo$-realization-nerve adjunction.
\end{itemize}

The following observation is the coherent version of Theorem \ref{gnerve-changeofg}.

\begin{theorem}\label{hcgnerve-changeofg}\index{change of action operads!coherent G-realization-nerve@coherent $\G$-realization-nerve}
There are natural isomorphisms\index{coherent G-nerve@coherent $\G$-nerve!change of action operads}\index{coherent G-realization@coherent $\G$-realization!change of action operads} \[\Nerveghcone\varphi^*\cong\varphi^*\Nerveghctwo \andspace \varphi_!\Realghcone\cong\Realghctwo\varphi_!\] in the diagram \eqref{hcnervegone-hcnervegtwo}.
\end{theorem}

\begin{proof}
To prove the first natural isomorphism, suppose $\O$ is a $\colorc$-colored $\Gtwo$-operad in $\M$, and $(T,\sigma)$ is a $\Gone$-tree.  There are natural isomorphisms:
\[\begin{split}
\Bigl(\Nerveghcone \varphi^*(\O)\Bigr)(T,\sigma) 
&= \goneopm\Bigl(\Wgone\iota_*(\Gone)^{(T,\sigma)}, \varphi^*\O\Bigr)\\
&\cong \gtwoopm\Bigl(\varphi_!\Wgone\iota_*(\Gone)^{(T,\sigma)}, \O\Bigr)\\
&\cong \gtwoopm\Bigl(\Wgtwo\varphi_!\iota_*(\Gone)^{(T,\sigma)}, \O\Bigr)\\
&\cong \gtwoopm\Bigl(\Wgtwo\iota_*\varphi_!(\Gone)^{(T,\sigma)}, \O\Bigr)\\
&\cong \gtwoopm\Bigl(\Wgtwo\iota_*(\Gtwo)^{(T,\varphi\sigma)},\O\Bigr)\\
&=\Bigl(\Nerveghctwo(\O)\Bigr)(T,\varphi\sigma)\\
&= \Bigl(\varphi^*\Nerveghctwo(\O)\Bigr)(T,\sigma).
\end{split}\]  
The three equalities are from the definitions of the coherent $\Gone$-nerve, the coherent $\Gtwo$-nerve, and the two right adjoints $\varphi^*$.  From top to bottom:
\begin{itemize}
\item The first isomorphism is from the bottom adjunction $\varphi_! \dashv \varphi^*$.
\item The second isomorphism is from Theorem \ref{wphi-to-phiw}.
\item The third isomorphism is from Example \ref{ex:set-to-m}.
\item The fourth isomorphism is from Theorem \ref{gdash-left-adjoint}.
\end{itemize}
Since $\O$ and $(T,\sigma)$ are arbitrary, we have established the first natural isomorphism.

The second natural isomorphism follows from the first one and the uniqueness of left adjoints.
\end{proof}

\begin{example}
Theorem \ref{hcgnerve-changeofg} applies to all the morphisms of action operads in Example \ref{ex:all-augmented-group-operads}.  For instance, the morphisms of action operads \[\nicexy{\P \ar[r]^-{\iota} & \PCac \ar[r]^-{\iota} & \Cac \ar[r]^-{\pi} & \S}\] induce the diagram
\[\nicexy@C+.7cm{\Settotreecatpop \ar@<2pt>[r]^-{\iota_!} \ar@<-2pt>[d]_-{\Realphc} 
& \Settotreecatpcacop \ar@<-2pt>[d]_-{\Realhc^{\PCac}} \ar@<2pt>[l]^-{\iota^*} \ar@<2pt>[r]^-{\iota_!}  
& \Settotreecatcacop \ar@<-2pt>[d]_-{\Realhc^{\Cac}} \ar@<2pt>[l]^-{\iota^*} \ar@<2pt>[r]^-{\pi_!}
& \Settotreecatsop \ar@<-2pt>[d]_-{\Realshc} \ar@<2pt>[l]^-{\pi^*}\\
\Poperadm \ar@<-2pt>[u]_-{\Nervephc} \ar@<2pt>[r]^-{\iota_!} 
& \Pcacoperadm \ar@<2pt>[l]^-{\iota^*} \ar@<-2pt>[u]_-{\Nervehc^{\PCac}} \ar@<2pt>[r]^-{\iota_!} 
& \Cacoperadm \ar@<2pt>[l]^-{\iota^*} \ar@<-2pt>[u]_-{\Nervehc^{\Cac}} \ar@<2pt>[r]^-{\pi_!}
& \Soperadm. \ar@<2pt>[l]^-{\pi^*} \ar@<-2pt>[u]_-{\Nerveshc}}\]
In each of the three squares, the left adjoint diagram and the right adjoint diagram are both commutative.  Using the ribbon group operad $\R$ and the braid group operad $\B$ instead of the cactus operad $\Cac$, we obtain two other diagrams with similar properties.\dqed
\end{example}

\section{Planar BV Construction of Planar Trees}\label{sec:ex-bv}

In this section, we compute the planar Boardman-Vogt construction of the planar operad freely generated by a planar tree for a fixed commutative segment $(J,\mu,0,1,\epsilon)$ in the ambient symmetric monoidal category $(\M,\otimes,\tensorunit)$.  This result is needed in Section \ref{sec:hcnerve-infinity-goperad} to prove that, for each action operad $\G$, the coherent $\G$-nerve of an entrywise fibrant $\G$-operad is an $\infty$-$\G$-operad.

For each symmetric monoidal category $(\M,\otimes,\tensorunit)$, there is a strong symmetric monoidal functor $\iota : \Set \to \M$ defined by \[\iota(X) = \coprodover{X} \tensorunit\] for each set $X$.  As discussed in Example \ref{ex:set-to-m}, there is an induced functor \[\nicexy{\gopcset \ar[r]^-{\iota_*} & \goperadcm}\] from the category of $\colorc$-colored $\G$-operads in $\Set$ to the category of $\colorc$-colored $\G$-operads in $\M$.  When this functor is applied to $\O\in\gopcset$, we will also write $\iota_*\O$ as $\subm{\O} \in \goperadcm$.  In particular, for the $\Ed(T)$-colored planar operad $\P^T$ in $\Set$ freely generated by the set of vertices in a planar tree $T$ in Definition \ref{def:gtree-goperad}, we have \[\iota_*\P^T=\ptsubm.\] We will write $\varnothing^{\M}$ for an initial object in $\M$.

\begin{theorem}\label{planar-bv}
Suppose $T$ is a planar tree, and $\duc\in \profedtedt$.  Then\index{planar operad!Boardman-Vogt construction} 
\[\Wp\ptsubm\duc\cong \begin{cases} \varnothing^{\M} & \text{if  $\P^T\duc=\varnothing$},\\
J^{\otimes |\Int(Q)|} &\text{if $\P^T\duc$ contains a single element $Q$}.\end{cases}\]
\end{theorem}

\begin{proof}
Recall from Lemma \ref{ptree-poperad} that each entry of $\P^T$ is either
\begin{enumerate}[label=(\roman*)]
\item empty or
\item a single element, in which case it is the unique planar tree with the given profile and with vertices from $T$.
\end{enumerate}
In the former case, the substitution category $\PTree^{\Ed(T)}\duc$ is empty.  So by   definition \eqref{wgoduc} we have \[\Wp\ptsubm\duc = \int^{T\in \PTree^{\Ed(T)}\duc} \J(T)\otimes\ptsubm(T) = \varnothing^{\M}.\]

In the second case, when $\P^T\duc$ contains a single element $Q$, the substitution category $\PTree^{\Ed(T)}\duc$ is the terminal category $\{Q\}$ because there are no non-identity endomorphisms in the substitution category $\PTree$.  Each vertex $q$ in $Q$ is a vertex in $T$, so $\P^T(q)$ contains a single element.  The coend definition  \eqref{wgoduc} reduces to 
\[\begin{split}
\Wp\ptsubm\duc &= \int^{\{Q\}} \J(Q)\otimes\ptsubm(Q)\\ 
&= \J(Q)\otimes \ptsubm(Q)\\
&= J^{\otimes |\Int(Q)|} \otimes \bigtensorover{q\in\Vt(Q)} \ptsubm(q)\\
&= J^{\otimes |\Int(Q)|} \otimes \bigtensorover{q\in\Vt(Q)} \iota\bigl(\P^T(q)\bigr)\\
&= J^{\otimes |\Int(Q)|} \otimes \bigtensorover{q\in\Vt(Q)} \tensorunit\\
&\cong J^{\otimes |\Int(Q)|} .
\end{split}\]
This proves the second case.
\end{proof}

\begin{example}
For the $2$-level tree $T = T\bigl(\{\ub_j\};\uc;d\bigr)$ in Example \ref{ex:2-level-tree}, we have \[\Wp\ptsubm\dub = J^{\otimes |\uc|}.\]\dqed
\end{example}

\begin{example}
For the planar tree $K$ in Example \ref{ex:treesub}, we have 
\[\begin{split}
\Wp\P^{K}_{\smallm}\sbinom{c}{a,b} & \cong \Wp\P^{K}_{\smallm}\sbinom{e}{c,d} \cong \Wp\P^{K}_{\smallm}\sbinom{e}{a,b,f,g}\cong J,\\
\Wp\P^{K}_{\smallm}\sbinom{e}{a,b,f,d} &\cong \Wp\P^{K}_{\smallm}\sbinom{e}{a,b,g} \cong J^{\otimes 2},\\
\Wp\P^{K}_{\smallm}\sbinom{e}{a,b,d} &\cong J^{\otimes 3}.
\end{split}\]
All other entries of $\Wp\P^{K}_{\smallm}$ are the initial object in $\M$.\dqed
\end{example}

\section{Coherent Nerves are Infinity Group Operads}
\label{sec:hcnerve-infinity-goperad}

The purpose of this section is to show that, for each action operad $\G$, the coherent $\G$-nerve of each entrywise fibrant $\G$-operad is an $\infty$-$\G$-operad.  The reader may refer back to Section \ref{sec:reedy} for basic concepts of model categories.  

\begin{definition}\label{def:monoidal-model-cat}
Suppose $(\M,\otimes,\tensorunit)$ is a symmetric monoidal closed category that is also a model category.  We say that $\M$ is a\index{model category!monoidal} \emph{monoidal model category} \cite{schwede-shipley} (Definition 3.1) if it satisfies the following \index{pushout product axiom}\emph{pushout product axiom}:
\begin{quote} Given cofibrations $f : A \to B$ and $g : C \to D$ in $\M$, the pushout product $f \square g$ in the diagram \[\nicexy{A \otimes C \ar@{}[dr]|-{\mathrm{pushout}} \ar[r]^-{\Id_A \otimes g} \ar[d]_-{f \otimes \Id_C} & A \otimes D \ar[d] \ar[ddr]^-{f \otimes \Id_D}\\ B \otimes C \ar[r] \ar[drr]_-{\Id_B \otimes g} & Z \ar[dr]|-{f \square g} &\\ && B \otimes D}\] is a cofibration that is furthermore a weak equivalence if either $f$ or $g$ is also a weak equivalence.  Here \[Z = B \otimes C \sqcupover{A \otimes C} A \otimes D\] is the object of the pushout square.\end{quote}
\end{definition}

Note that in \cite{hovey} (Definition 4.2.6), a monoidal model category has an extra condition about the monoidal unit.

\begin{assumption}\label{assumption:monoidal-model}
Throughout this section, we assume the following.
\begin{enumerate}
\item $\M$ is a monoidal model category with $\tensorunit$ cofibrant.
\item $(J,\mu,0,1,\epsilon)$ is an \emph{interval} in $\M$ \cite{berger-moerdijk-bv}.  This means that $J$ is a commutative segment in $\M$ in the sense of Definition \ref{def:segment} such that:
\begin{enumerate}[label=(\roman*)]
\item The morphism $(0,1) : \tensorunit \amalg \tensorunit \to J$ is a cofibration in $\M$.
\item The counit $\epsilon : J \to \tensorunit$ is a weak equivalence in $\M$.
\end{enumerate}
\end{enumerate}
\end{assumption}
Under these assumptions, both morphisms $0,1 : \tensorunit \to J$ are acyclic cofibrations, and $J$ is cofibrant.  

\begin{example}
The above assumptions hold for the categories $\CHau$, $\Sset$, $\Chaink$, and $\Cat$ in Example \ref{ex:com-segment}.\dqed
\end{example}

Recall from Definition \ref{def:infinity-goperad} that an $\infty$-$\G$-operad is a $\Treecatg$-presheaf in which every inner horn admits a filler.  Also recall from Definition \ref{def:hcrealization-hcnerve} the coherent $\G$-nerve \[\nicexy@C+.5cm{\Settotreecatgop & \gopm \ar[l]_-{\Nerveghc}},\] which we computed in Proposition \ref{hcn-planar}. 

\begin{theorem}\label{cnerve-infinity-goperad}
Under Assumption \ref{assumption:monoidal-model}, suppose $\O$ is a $\colorc$-colored $\G$-operad in $\M$ such that each entry of $\O$ is a fibrant object.  Then its\index{coherent G-nerve@coherent $\G$-nerve!is infinity-G-operad@is $\infty$-$\G$-operad}\index{infinity-G-operad@$\infty$-$\G$-operad!from coherent nerve} coherent $\G$-nerve is an $\infty$-$\G$-operad.
\end{theorem}

\begin{proof}
Suppose given an inner horn $f$ of the coherent $\G$-nerve $\Nerveghc(\O)$,
\[\nicexy{\horn{\phi}(V,\sigma^V) \ar[r]^-{f} \ar@{>->}[d] & \Nerveghc(\O)
\\ \Yonedag(V,\sigma^V) \ar@{-->}[ur]_-{F} &}\] as in \eqref{horn-filler} for some inner coface \[\nicexy{(T,\sigma) \ar[r]^-{\phi} & (V,\sigma^V).}\] We must show that there exists a dashed arrow $F$ that makes the above diagram commutative.

By the descriptions of (i) a horn in Lemma \ref{horn-description} and (ii) the coherent $\G$-nerve in Proposition \ref{hcn-planar}, the given inner horn $f$ is equivalent to a family of compatible morphisms of planar operads\[\nicexy{\Wp\iota_*\P^{T'} \ar[r]^-{f_{\phi'}} & (\iota^{\G})^*\O} \in \Poperadm,\] one for each coface \[\nicexy{(T',\sigma')\ar[r]^-{\phi'} & (V,\sigma^V)}\] of $(V,\sigma^V)$ not isomorphic to $\phi$. By Proposition \ref{hcn-planar} again, the desired filler $F$ is equivalent to a morphism \[\nicexy{\Wp\iota_*\P^{V} \ar[r]^-{F} & (\iota^{\G})^*\O}\in \Poperadm\] that extends the given morphisms $\bigl\{f_{\phi'} : \phi' \not\cong \phi\bigr\}$.

To define the morphism $F$ on color sets, note that each $\P^{T'}$ is $\Ed(T')$-colored, and so is the planar Boardman-Vogt construction $\Wp\iota_*\P^{T'}$.  Similarly, the planar Boardman-Vogt construction $\Wp\iota_*\P^V$ is $\Ed(V)$-colored.  As $\phi'$ runs through the cofaces of $(V,\sigma^V)$ not isomorphic to the inner coface $\phi$, there is an equality \[\bigcup_{\text{$\phi\not\cong\phi'$ cofaces}}\Ed(T')=\Ed(V)\] on edge sets.  Therefore, the morphism $F$ on color sets is already defined by the morphisms $f_{\phi'}$ on color sets.

To define the entries of $F$, we use the descriptions of (i) the free $\Ed(T)$-colored planar operad $\P^T$ in Lemma \ref{ptree-poperad} and (ii) its planar Boardman-Vogt construction $\Wp\iota_*\P^T$ in Theorem \ref{planar-bv}.  Suppose \[\duc \in \profedvedv \withspace \P^V\duc \not=\varnothing,\] so $\P^V\duc = \{Q\}$ for a unique $\Ed(V)$-colored planar tree $Q$ with profile $\duc$ and with vertices in $V$.  There is a morphism \[\varphi : Q \to V \in\Treecatp\] corresponding to a planar tree substitution \[V = U \comp_u Q.\] There are two cases.

First, if $\duc \not= \sbinom{\out(V)}{\inp(V)}$, then $U$ has at least one internal edge, and $Q\not=V$.  By Lemma \ref{outer-coface-iterated} there is a factorization \[\nicexy{Q \ar[r]^-{\varphi'} & T' \ar[r]^-{\phi'} & V \ar@{<-}`u[ll]`[ll]_-{\varphi}[ll]}\]
in $\Treecatp$ in which:
\begin{itemize}
\item $\phi'$ is an outer coface.
\item $\varphi'$ is either the identity morphism or a finite composite of outer cofaces.
\end{itemize}
This implies that $\P^{T'}\duc = \{Q\}$, and 
\[\begin{split}
\Wp\iota_*\P^{V}\duc &= J^{\otimes|\Int(Q)|}\\
&= \Wp\iota_*\P^{T'}\duc.
\end{split}\]  
So for such pair $\duc$, the morphism $F$ must be defined as $f_{\phi'}$.

The only remaining case is $\duc= \sbinom{\out(V)}{\inp(V)}$.  This means that $U$ is a corolla and that $Q=V$.  For this entry, the desired morphism $F$ has the form \[\nicexy{\Wp\iota_*\P^V\sbinom{\out(V)}{\inp(V)} = J^{\otimes|\Int(V)|} \ar[r]^-{F} & (\iota^{\G})^*\O\sbinom{F\out(V)}{F\inp(V)}}\] that satisfies the following two conditions.
\begin{enumerate}[label=(\roman*)]
\item For each inner coface $\phi'$ of $(V,\sigma^V)$ not isomorphic to the given inner coface $\phi$, $T'$ has the same profile as $V$, and \[\Int(T') = \Int(V) \setminus \{e'\}\] for the internal edge $e'$ corresponding to $\phi'$.  Then the diagram
\[\nicexy{\Wp\iota_*\P^{T'}\sbinom{\out(T')}{\inp(T')} \cong J^{\otimes|\Int(T')|}\otimes \tensorunit \ar[r]^-{f_{\phi'}} \ar[d]_-{\Id \otimes 0} & (\iota^{\G})^*\O\sbinom{F\out(V)}{F\inp(V)}\\
\Wp\iota_*\P^V\sbinom{\out(V)}{\inp(V)} = J^{\otimes|\Int(V)|} \ar`r/3pt[ur]^-{F}[ur]
}\]
is commutative, where the morphism $0 : \tensorunit \to J$ corresponds to $e'$.
\item The morphism $F$ respects the $\compi$-composition relative to the other entries of $F$ that are defined in the previous case.  In view of the operadic structure in the $\G$-Boardman-Vogt construction in Definition \ref{def:wg-goperad-structure}, this condition determines the restriction of $F$ to
\[\nicexy@C+1cm{\underset{e'\in\Int(V)}{\bigcup} J^{\otimes |\Int(V)\setminus\{e'\}|} \otimes \tensorunit \ar[r]^-{\underset{e'\in\Int(V)}{\bigcup} \, \Id\otimes 1} & J^{\otimes |\Int(V)|}},\] in which each morphism $1 : \tensorunit \to J$ corresponds to $e'$.  The union symbol here denotes a colimit whose indexing category is the opposite category of a punctured cube.  More precisely, the indexing category $\C_V$ of this colimit has non-empty subsets of $\Int(V)$ as objects, and a morphism $A \to B$ is an inclusion $B \subseteq A \subseteq \Int(V)$.  The functor $\C_V\to\M$ that defines this colimit sends each non-empty subset $A \subseteq \Int(V)$ to the object \[J^{\otimes |\Int(V) \setminus A|} \otimes \tensorunit^{|A|}.\]  This functor sends a morphism $A \to B \in \C_V$ to the morphism given by:
\begin{itemize}
\item $1 : \tensorunit \to J$ on $A \setminus B$;
\item the identity morphism of either $J$ or $\tensorunit$ in all other tensor factors.
\end{itemize}
\end{enumerate}

The given inner coface $\phi$ corresponds to an internal edge $e$ in $V$.  The two conditions in the previous paragraph determine the morphism $\{f_{\phi'}\}$ below, which is defined by the given morphisms $f_{\phi'}$.
\[\nicexy@C-1.5cm{\Bigl[\underset{e\not=e'\in\Int(V)}{\bigcup} J^{\otimes |\Int(V)\setminus\{e'\}|} \otimes (\tensorunit \amalg \tensorunit)\Bigr] \cup \bigl[J^{\otimes|\Int(V)\setminus\{e\}|} \otimes \tensorunit\bigr] \ar[dd]_{\bigl[\underset{e\not=e'\in\Int(V)}{\bigcup} \Id \otimes (0,1)\bigr]}^{\cup\, [\Id\otimes 1]} &\\ 
& (\iota^{\G})^*\O\sbinom{F\out(V)}{F\inp(V)} \ar@{<-}`u/3pt[ul]_(.6){\{f_{\phi'}\}} [ul]\\
\Wp\iota_*\P^V\sbinom{\out(V)}{\inp(V)} = J^{\otimes|\Int(V)|} \ar@{-->}`r/3pt[ur]^-{F} [ur]&}\] 
The object in the upper-left corner is again a colimit indexed by the category $\C_V$ in case (ii) in the previous paragraph.  The functor that defines this colimit is modified as displayed.  More precisely, it sends a non-empty subset $A \subseteq \Int(V)$ to the object \[J^{\otimes |\Int(V) \setminus A|} \otimes (\tensorunit \amalg \tensorunit)^{A \setminus \{e\}} \otimes \tensorunit^{\otimes |A \cap \{e\}|}.\] This functor sends a morphism $A \to B \in \C_V$ to the morphism given by: 
\begin{itemize}
\item the cofibration $(0,1) : \tensorunit \amalg \tensorunit \to J$ for $(A\setminus \{e\}) \setminus (B \setminus \{e\})$;
\item the acyclic cofibration $1 : \tensorunit \to J$ for $(A \cap \{e\}) \setminus (B \cap \{e\})$;
\item the identity morphism of $J$, $\tensorunit \amalg \tensorunit$, or $\tensorunit$ in all other tensor factors.
\end{itemize}
Since $J$ is a cofibrant object in $\M$, the left vertical morphism in the previous diagram is an acyclic cofibration by repeated applications of the pushout product axiom.  Since each entry of $\O$ is a fibrant object, the dashed arrow $F$ exists that makes the entire diagram commutative.
\end{proof}

\begin{remark}
The analogue of Theorem \ref{cnerve-infinity-goperad} for simplicially enriched categories is a result in \cite{cp}.  It says that for each small category enriched in the category of simplicial sets, its homotopy coherent nerve is a quasi-category, i.e., a weak Kan complex in the sense of \cite{boardman-vogt}.\dqed
\end{remark}

\begin{example}[Every Object is Fibrant]
In Theorem \ref{cnerve-infinity-goperad} the $\G$-operad is assumed to have fibrant entries.  This is automatically true if every object in $\M$ is fibrant.  This is the case, for example, if $\M$ is the categories $\CHau$, $\Cat$, or $\Chaink$ in Example \ref{ex:com-segment}.  For such an ambient category $\M$, the coherent $\G$-nerve of each $\G$-operad in $\M$ is an $\infty$-$\G$-operad.\dqed
\end{example}

\begin{example}[Simplicial Sets]
For the category $\Sset$ of simplicial sets, fibrant objects are called\index{Kan complex} \emph{Kan complexes}.  They are the simplicial sets in which every horn admits an extension; see, for example, Section I.3 in \cite{goerss-jardine}.  The coherent $\G$-nerve of each colored $\G$-operad in $\Sset$ that is entrywise a Kan complex is an $\infty$-$\G$-operad. \dqed
\end{example}

\begin{example}[Coherent $\G$-Nerve of Translation Categories]
Recall from Definition \ref{def:translation-category} that the\index{translation category}\index{coherent G-nerve@coherent $\G$-nerve!of translation categories} translation category of a group $G$ is the small groupoid $E_G$ with object set $G$ and with a unique morphism between any two objects.  In Example \ref{ex:nerve-goperad} we mentioned that for each action operad $\G$, there is a one-colored $\G$-operad \[E_{\G} =\bigl\{E_{\G(n)}\bigr\}_{n\geq 0}\] in $\Cat$, in which the $n$th object is the translation category of the group $\G(n)$.  In the category $\Cat$ of all small categories with its folk model structure, every object is fibrant.  Since the one-colored $\G$-operad $E_{\G}$ is entrywise fibrant, by Theorem \ref{cnerve-infinity-goperad} its coherent $\G$-nerve is an $\infty$-$\G$-operad.\dqed
\end{example}

\begin{example}[Coherent Planar Nerve]
When $\G$ is the planar group operad $\P$ in Example \ref{ex:trivial-group-operad}, Theorem \ref{cnerve-infinity-goperad} says that the coherent planar nerve of an entrywise fibrant planar operad is an\index{infinity-planar operad@$\infty$-planar operad} $\infty$-planar operad.\dqed
\end{example}

\begin{example}[Coherent Symmetric Nerve]
For the symmetric group operad $\S$ in Example \ref{ex:symmetric-group-operad}, the symmetric tree category $\Treecats$ is equivalent to the symmetric dendroidal category $\Omega$ of Moerdijk-Weiss.  In this case, Theorem \ref{cnerve-infinity-goperad} recovers Theorem 7.1 in \cite{mw09}, which says that the homotopy coherent dendroidal nerve of an entrywise fibrant symmetric operad is an\index{infinity-symmetric operad@$\infty$-symmetric operad} $\infty$-symmetric operad.  For example, 
\begin{itemize}
\item the little $n$-cube operad $\C_n$ in Example \ref{ex:little-ncube},
\item the little $n$-disc operad $\D_n$ in Example \ref{ex:little-n-disc},
\item the framed little $n$-disc operad $\Df_n$ in Example \ref{ex:framed-little-n-disc},
\item the phylogenetic operad $\Phyl$ in Section \ref{sec:symmetric-operad-phylogenetic}, and
\item the planar tangle operad $\PTan$ in Definition \ref{def:planar-tangle-operad} 
\end{itemize}
are symmetric operads in $\CHau$, one-colored in the first four cases and $\SCir$-colored for $\PTan$.   Since every object in $\CHau$ is fibrant, the homotopy coherent dendroidal nerves of these symmetric operads are $\infty$-symmetric operads.\dqed
\end{example}

\begin{example}[Coherent Braided Nerve]
For the braid group operad $\B$ in Definition \ref{def:braid-group-operad}, Theorem \ref{cnerve-infinity-goperad} says that the coherent braided nerve of an entrywise fibrant braided operad is an\index{infinity-braided operad@$\infty$-braided operad} $\infty$-braided operad.  For example, \begin{itemize}
\item the universal cover $\Ctilde_2$ of the little $2$-cube operad in Example \ref{ex:universal-cover-c2} and
\item the universal cover $\Dtilde_2$ of the little $2$-disc operad in Example \ref{ex:universal-cover-d2}
\end{itemize} 
are one-colored braided operads in $\CHau$.  So their coherent braided nerves are $\infty$-braided operads.\dqed
\end{example}

\begin{example}[Coherent Ribbon Nerve]
For the ribbon group operad $\R$ in Definition \ref{def:ribbon-group-operad}, Theorem \ref{cnerve-infinity-goperad} says that the coherent ribbon nerve of an entrywise fibrant ribbon operad is an\index{infinity-ribbon operad@$\infty$-ribbon operad} $\infty$-ribbon operad.  For example, the universal cover $\Dftildetwo$ of the framed little $2$-disc operad in Example \ref{ex:universal-cover-fd2} is a one-colored ribbon operad in $\CHau$.  So its coherent ribbon nerve is an $\infty$-ribbon operad.\dqed
\end{example}

\begin{example}[Coherent Cactus Nerve]
For the cactus group operad $\Cac$ in Definition \ref{def:cactus-group-operad}, Theorem \ref{cnerve-infinity-goperad} says that the coherent cactus nerve of an entrywise fibrant cactus operad is an\index{infinity-cactus operad@$\infty$-cactus operad} $\infty$-cactus operad.  There is also a pure analogue for the pure cactus group operad $\PCac$ in Definition \ref{def:pure-cactus-group-operad}.\dqed
\end{example}

%% file: coherence_for_monoidal_categories.tex
\part{Coherence for Monoidal Categories with Group Equivariance}\label{part:monoidal-cat}

\chapter{Monoidal Categories}\label{ch:monoidal-cat}

A monoid is a set with a strictly associative multiplication and a two-sided multiplicative unit.  A monoidal category is a categorical analogue of a monoid with the strict axioms replaced by suitable natural isomorphisms.  In preparation for later discussion about $\G$-monoidal categories for an action operad $\G$, in this chapter we discuss the symmetric operads that encode monoidal categories, strict monoidal categories, monoidal functors, and strong monoidal functors, thereby providing an operadic view of these monoidal categorical concepts.

The definitions of a monoidal category and of a monoidal functor, as well as one version of the Coherence Theorem for monoidal categories, are recalled in Section \ref{sec:monoidal-cat-functor}.    A \emph{strict monoidal category} is a monoidal category whose structure morphisms are identity morphisms.  The structure morphisms of a monoidal functor are not required to be invertible; our monoidal functors are sometimes called \emph{lax monoidal functors} in the literature.  A monoidal functor with invertible structure morphisms is called a \emph{strong monoidal functor}; this is what Joyal and Street \cite{joyal-street} called a \emph{tensor functor}.  A \emph{strict monoidal functor} is a monoidal functor whose structure morphisms are identity morphisms.

In Section \ref{sec:monoidal-cat-operad} we define the monoidal category operad $\MCat$, which is a $1$-colored symmetric operad in the category $\Cat$ of small categories.  In Section \ref{sec:sym-operad-monoidal-cat} we show that $\MCat$-algebras are exactly small monoidal categories with general associativity isomorphism, left unit isomorphism, and right unit isomorphism.  In Section \ref{sec:strict-monoidal-functor} it is shown that $\MCat$-algebra morphisms are strict monoidal functors between general small monoidal categories.  In Section \ref{sec:mcatst} we discuss a one-colored symmetric operad in $\Cat$ whose algebras are small strict monoidal categories and whose algebra morphisms are strict monoidal functors.  

In Section \ref{sec:operad-monoidal-functor} we discuss a $2$-colored symmetric operad $\MFun$ in $\Cat$ whose algebras are monoidal functors.  More precisely, an $\MFun$-algebra is a triple $(A,B,F)$ in which (i) $A$ and $B$ are small monoidal categories, and (ii) $F : A \to B$ is a monoidal functor.   We also discuss another $2$-colored symmetric operad in $\Cat$ whose algebras are strong monoidal functors.

\section{Monoidal Categories and Monoidal Functors}
\label{sec:monoidal-cat-functor}

We begin by recalling the definitions of a monoidal category and of a monoidal functor and the Coherence Theorem for monoidal categories.  Our references for monoidal categories are \cite{joyal-street}, \cite{maclane-rice}, and \cite{maclane} (Chapters VII and XI).

\begin{definition}\label{def:monoidal-category}
A \index{monoidal category}\index{category!monoidal}\emph{monoidal category} is a tuple
\[(\M, \otimes, \tensorunit, \alpha, \lambda, \rho)\]
consisting of:
\begin{itemize}
\item a category $\M$;
\item a functor $\otimes : \M \times \M \to \M$\label{notation:monoidal-product} called the \index{monoidal product}\emph{monoidal product};
\item an object $\tensorunit\in\M$ called the \index{monoidal unit}\emph{monoidal unit};
\item a natural isomorphism
\begin{equation}\label{mon-cat-alpha}
\nicexy{(X \otimes Y) \otimes Z \ar[r]^-{\alpha_{X,Y,Z}}_-{\cong} &  X\otimes (Y \otimes Z)}
\end{equation}
for all objects $X,Y,Z \in \M$ called the \index{associativity isomorphism}\emph{associativity isomorphism};
\item natural isomorphisms
\begin{equation}\label{mon-cat-lambda}
\nicexy{\tensorunit \otimes X \ar[r]^-{\lambda_X}_-{\cong} & X} \andspace \nicexy{X \otimes \tensorunit \ar[r]^-{\rho_X}_-{\cong} & X}
\end{equation}
for all objects $X \in \M$ called the \emph{left unit isomorphism}\index{left unit isomorphism} and the \index{right unit isomorphism}\emph{right unit isomorphism}, respectively.
\end{itemize}
This data is required to satisfy the following two axioms.
\begin{description}
\item[Unity Axioms]
The diagram
\begin{equation}\label{monoidal-unit}
\nicexy{(X \otimes \tensorunit) \otimes Y \ar[d]_-{\rho \otimes \Id} \ar[r]^-{\alpha}
& X \otimes (\tensorunit \otimes Y) \ar[d]^-{\Id \otimes \lambda}\\ X \otimes Y \ar@{=}[r] & X \otimes Y}
\end{equation}
is commutative for all objects $X,Y \in \M$.  Moreover, the equality \[\nicexy{\lambda_{\tensorunit} = \rho_{\tensorunit} : \tensorunit \otimes \tensorunit \ar[r]^-{\cong} & \tensorunit}\] holds.
\item[Pentagon Axiom]
The pentagon\index{Pentagon Axiom}
\[\nicexy@C-1.5cm{& (W \otimes X) \otimes (Y \otimes Z) \ar[dr]^-{\alpha_{W,X,Y\otimes Z} }& \\
\bigl((W \otimes X) \otimes Y\bigr) \otimes Z \ar[ur]^-{\alpha_{W\otimes X,Y,Z}} \ar[d]_-{\alpha_{W,X,Y} \otimes \Id_Z}
&& W \otimes \bigl(X \otimes (Y \otimes Z)\bigr) \\
\bigl(W \otimes (X \otimes Y)\bigr) \otimes Z \ar[rr]^-{\alpha_{W,X\otimes Y,Z}} & \hspace{1cm} & W \otimes \bigl((X \otimes Y) \otimes Z\bigr)\ar[u]_-{\Id_W \otimes \alpha_{X,Y,Z}}}\]
is commutative for all objects $W,X,Y,Z \in \M$.
\end{description}
A \index{strict monoidal category}\index{monoidal category!strict}\emph{strict monoidal category} is a monoidal category in which $\alpha$, $\lambda$, and $\rho$ are all identity morphisms.
\end{definition}

\begin{remark}
Consider Definition \ref{def:monoidal-category}
\begin{enumerate}
\item In a strict monoidal category, an iterated monoidal product $a_1 \otimes \cdots \otimes a_n$ without any parentheses is well-defined.  
\item The axiom $\lambda_{\tensorunit} = \rho_{\tensorunit}$ in a monoidal category is actually redundant by \cite{kelly2}.  It is a consequence of the unity axiom \eqref{monoidal-unit} and the pentagon axiom. 
\item In a monoidal category, every diagram involving monoidal products of the associativity isomorphism, the unit isomorphisms, their inverses, and identity morphisms commutes.  This is one formulation of the Coherence Theorem for monoidal categories; see \cite{maclane} (VII.2 Corollary).  For example, in any monoidal category, the unity axiom \eqref{monoidal-unit} and the pentagon axiom imply the commutativity of the following unity diagrams.
\begin{equation}\label{moncat-other-unit-axioms}
\nicexy{(\tensorunit \otimes X) \otimes Y \ar[d]_-{\lambda \otimes \Id} \ar[r]^-{\alpha}
& \tensorunit \otimes (X\otimes Y) \ar[d]^-{\lambda}\\ X \otimes Y \ar@{=}[r]& X \otimes Y}\qquad
\nicexy{(X \otimes Y) \otimes \tensorunit \ar[d]_-{\rho} \ar[r]^-{\alpha}
& X \otimes (Y\otimes \tensorunit ) \ar[d]^-{\Id \otimes \rho}\\ X \otimes Y \ar@{=}[r]& X \otimes Y}
\end{equation} 
See\cite{maclane} VII.1 Eq. (9) and Exercise 1.\dqed
\end{enumerate}
\end{remark}

\begin{convention}\label{conv:empty-tensor}
In a monoidal category, an \index{empty tensor product}\emph{empty tensor product}, written as\label{notation:empty-tensor} $X^{\otimes 0}$ or $X^{\otimes \varnothing}$, means the monoidal unit $\tensorunit$.\dqed
\end{convention}

\begin{definition}\label{def:monoidal-functor}
Suppose $\M$ and $\N$ are monoidal categories.  A \index{monoidal functor}\index{functor!monoidal}\emph{monoidal functor} \[(F,F_2,F_0) : \M \to \N\] consists of the following data:
\begin{itemize}
\item a functor $F : \M \to \N$;
\item a natural transformation
\begin{equation}\label{monoidal-f2}
\nicexy{F(X) \otimes F(Y) \ar[r]^-{F_2} & F(X \otimes Y) \in \N,}
\end{equation}
where $X$ and $Y$ are objects in $\M$;
\item a morphism
\begin{equation}\label{monoidal-f0}
\nicexy{\tensorunit^{\N} \ar[r]^-{F_0} & F(\tensorunit^{\M}) \in \N,}
\end{equation}
where $\tensorunit^{\N}$ and $\tensorunit^\M$ are the monoidal units in $\N$ and $\M$, respectively.
\end{itemize} 
This data is required to satisfy the following three axioms.
\begin{description}
\item[Associativity]
The diagram\index{associativity!monoidal functor}
\begin{equation}\label{f2}
\nicexy{\bigl(F(X) \otimes F(Y)\bigr) \otimes F(Z) \ar[r]^-{\alpha^{\N}} \ar[d]_-{F_2 \otimes \Id} & F(X) \otimes \bigl(F(Y) \otimes F(Z)\bigr) \ar[d]^-{\Id \otimes F_2}\\
F(X \otimes Y) \otimes F(Z) \ar[d]_-{F_2} & F(X) \otimes F(Y \otimes Z) \ar[d]^-{F_2}\\
F\bigl((X \otimes Y) \otimes Z\bigr) \ar[r]^-{F(\alpha^{\M})} &
F\bigl(X \otimes (Y \otimes Z)\bigr)}
\end{equation}
is commutative for all objects $X,Y,Z \in \M$.
\item[Left Unity]
The diagram
\begin{equation}\label{f0-left}
\nicexy{\tensorunit_{\N} \otimes F(X) \ar[d]_-{F_0 \otimes \Id} \ar[r]^-{\lambda^{\N}}
& F(X) \\ F(\tensorunit_{\M}) \otimes F(X) \ar[r]^-{F_2} & F(\tensorunit_{\M} \otimes X)
\ar[u]_-{F(\lambda^{\M})}}
\end{equation}
is commutative for all objects $X \in \M$.
\item[Right Unity]
The diagram
\begin{equation}\label{f0-right}
\nicexy{F(X) \otimes \tensorunit_{\N} \ar[d]_-{\Id \otimes F_0} \ar[r]^-{\rho^{\N}}
& F(X) \\ F(X) \otimes F(\tensorunit_{\M}) \ar[r]^-{F_2} & F(X \otimes \tensorunit_{\M})
\ar[u]_-{F(\rho^{\M})}}
\end{equation}
is commutative for all objects $X \in \M$.
\end{description}
A \index{strong monoidal functor}\index{monoidal functor!strong}\emph{strong monoidal functor} is a monoidal functor in which the morphisms $F_0$ and $F_2$ are all isomorphisms.  A \index{strict monoidal functor}\index{monoidal functor!strict}\emph{strict monoidal functor} is a monoidal functor in which the morphisms $F_0$ and $F_2$ are all identity morphisms. 
\end{definition}

\begin{remark}
In a monoidal functor, the structure morphisms $F_2$ and $F_0$ are not required to be isomorphisms.  What we call a monoidal functor is sometimes called a \index{lax monoidal functor}\emph{lax monoidal functor} in the literature to emphasize the fact that $F_2$ and $F_0$ are not necessarily invertible.  A strong monoidal functor is what Joyal and Street \cite{joyal-street} called a\index{tensor functor} \emph{tensor functor}.\dqed
\end{remark}

The following coherence result for monoidal categories is due to Mac Lane.  See \cite{maclane-rice} (Section 3), \cite{maclane} (XI.3), or \cite{joyal-street} (Corollary 1.4) for a proof. 
 
\begin{theorem}[Coherence for Monoidal Categories]\label{maclane-thm}\index{strictification!monoidal category}\index{coherence!monoidal category}\index{monoidal category!coherence} 
For each monoidal category $\M$, there exist a strict monoidal category $\M_{\st}$ and an adjoint equivalence
\[\nicexy{\M \ar@<2pt>[r]^-{L} & \M_{\st} \ar@<2pt>[l]^-{R}}\]
with (i) both $L$ and $R$ strong monoidal functors and (ii) $RL=\Id_{\M}$.
\end{theorem}

\section{Monoidal Category Operad}\label{sec:monoidal-cat-operad}

The purpose of this section is to construct a one-colored symmetric operad in $\Cat$, called the monoidal category operad, whose algebras will be shown to be precisely monoidal categories with general associativity, left unit, and right unit isomorphisms.  The definition of the monoidal category operad uses the following bookkeeping devices.

\begin{definition}\label{def:nonass-monomial}
Fix a sequence $\{x_i\}_{i\geq 0}$ of variables.  Define\index{non-associative monomial} \emph{non-associative monomials} inductively as follows:
\begin{itemize}
\item Each $x_i$ for $i\geq 0$ is a non-associative monomial.
\item If $f$ and $g$ are non-associative monomials, then $(fg)$ is a non-associative monomial.
\end{itemize}
For $n\geq 0$, a\index{standard non-associative monomial}\index{non-associative monomial!standard} \emph{standard non-associative monomial of weight $n$} is a non-associative monomial in which:
\begin{enumerate}[label=(\roman*)]
\item Each $x_j$ appears exactly once for $1\leq j \leq n$.
\item No $x_k$ is involved for $k>n$.
\end{enumerate}
The associative monomial obtained from a standard non-associative monomial $p$ of weight $n$ by removing all the parentheses and all the $x_0$'s is denoted by $\pbar$ and is called a\index{standard monomial} \emph{standard monomial of weight $n$}.  It is also called the\index{underlying standard monomial} \emph{underlying standard monomial} of $p$.   A \emph{standard (non-associative) monomial} is a standard (non-associative) monomial of weight $n$ for some $n \geq 0$.
\end{definition}

\begin{example}
To simplify the notation, we will omit the outermost pair of parentheses.  
\begin{enumerate}
\item Standard non-associative monomials of weight $0$ are precisely the parenthesized words in one alphabet $\{x_0\}$, as in Example \ref{ex:mag-planar-operad}, of length $\geq 1$:
\[x_0,\quad x_0x_0,\quad (x_0x_0)x_0,\quad x_0(x_0x_0), \quad \bigl((x_0x_0)x_0\bigr)x_0, \ldots.\]
The only standard monomial of weight $0$ is $\varnothing$.
\item Standard non-associative monomials of weight $1$ include $x_1$, $x_0x_1$, $x_1x_0$, \[(x_0x_0)x_1,\, (x_0x_1)x_0,\, (x_1x_0)x_0,\, x_0(x_0x_1),\, x_0(x_1x_0),\, x_1(x_0x_0),\] and so forth.  The only standard monomial of weight $1$ is $x_1$.
\item Standard non-associative monomials of weight $2$ include $x_1x_2$, $x_2x_1$, \[(x_1x_2)x_0,\, (x_2x_1)x_0,\, (x_1x_0)x_2\, (x_2x_0)x_1,\, (x_0x_1)x_2,\, (x_0x_2)x_1,\] and so forth.  The only standard monomials of weight $2$ are $x_1x_2$ and $x_2x_1$.
\item Standard non-associative monomials of weight $3$ include \[(x_1x_2)x_3,\, x_2(x_1x_3),\,  x_0\bigl((x_1x_2)x_3\bigr),\, \bigl((x_1x_2)x_3\bigr)x_0,\, \bigl((x_0x_1)x_2\bigr)x_3,\, \bigl((x_1x_0)x_2\bigr)x_3,\] and so forth.  The only standard monomials of weight $3$ are: \[x_1x_2x_3,\, x_1x_3x_2,\, x_2x_1x_3,\, x_3x_1x_2,\, x_2x_3x_1, \andspace x_3x_2x_1.\]
\end{enumerate}
In general, standard non-associative monomials of weight $n \geq 1$ are obtained from standard non-associative monomials of weight $0$ and length $\geq n$ by replacing $n$ appearances of $x_0$'s with $\{x_1,\ldots,x_n\}$ in any order.  A standard monomial of weight $n$ is of the form $x_{\sigma(1)}\cdots x_{\sigma(n)}$ for a unique permutation $\sigma \in S_n$.\dqed
\end{example}

\begin{notation}\label{not:subscript-interval}
To simplify the presentation, we will use the abbreviation \[x_{[i,j]} = \begin{cases} \bigl(x_i,\ldots,x_j\bigr) & \text{if $i \leq j$},\\ \varnothing & \text{if $i > j$}\end{cases}\] for sequences in the variables, and similarly for other sequences of objects.\dqed 
\end{notation}

\begin{example} $x_{[2,6]} = (x_2,x_3,x_4,x_5,x_6)$, $x_{[0,0]} = x_0$, and $x_{[3,2]}=\varnothing$.\dqed
\end{example}

Recall that $\Cat$ is the category of small categories and functors between them.   A \emph{discrete category} is a category with only identity morphisms.  We now define the one-colored symmetric operad in $\Cat$ whose algebras are exactly small monoidal categories.  

\begin{definition}\label{def:mcat-operad}
Define the\index{monoidal category operad}\index{operad!monoidal category} \emph{monoidal category operad} $\MCat$ as the $1$-colored symmetric operad in $\Cat$ as follows.
\begin{description}
\item[Entries]
For each $n\geq 0$, $\MCat(n)$ is the groupoid whose object set is the set of standard non-associative monomials of weight $n$.  For $p,q\in \Ob(\MCat(n))$, it has the morphism set \[\MCat(n)(p,q) = \begin{cases} \{\ast\} & \text{if $\pbar=\qbar$},\\ \varnothing & \text{otherwise}.\end{cases}\]  In other words, there is a unique isomorphism from $p$ to $q$ if and only if they have the same underlying standard monomials.
\item[Equivariance]
The symmetric group $S_n$ acts on the objects in $\MCat(n)$ by permuting the entries $\{x_1,\ldots,x_n\}$ in each standard non-associative monomial of weight $n$.  This definition uniquely determines the action of $S_n$ on the morphism sets in $\MCat(n)$.
\item[Unit]
The operadic unit is the object $x_1 \in \MCat(1)$.
\item[Composition]
For $1\leq i \leq n$ and $m \geq 0$, the $\compi$-composition \[\nicexy{\MCat(n) \times \MCat(m) \ar[r]^-{\compi} & \MCat(n+m-1)}\] is the functor defined as \[p\bigl(x_{[0,n]}\bigr) \compi q\bigl(x_{[0,m]}\bigr) = p\Bigl(x_{[0,i-1]}, q\bigl(x_0,x_{[i,i+m-1]}\bigr), x_{[i+m,n+m-1]}\Bigr)\]
for objects $p\bigl(x_{[0,n]}\bigr)\in \MCat(n)$ and $q\bigl(x_{[0,m]}\bigr) \in \MCat(m)$.  This definition uniquely determines the $\compi$-composition on morphisms.
\end{description}
The finishes the definition of the monoidal category operad.
\end{definition}

\begin{lemma}\label{mcat-is-operad}
$\MCat$ is a $1$-colored symmetric operad in $\Cat$.
\end{lemma}

\begin{proof}
We need to check that $\MCat$ satisfies the associativity and unity axioms in Definition \ref{def:planar-operad-compi}  and the equivariance axiom in Proposition \ref{prop:symmetric-operad-compi}.  The unity axioms and the equivariance axiom are immediate from the definition.  

To check the horizontal associativity axiom, suppose $1 \leq i < j \leq n$ and $r\in \Ob(\MCat(l))$.  Each composite in the horizontal associativity diagram \eqref{compi-associativity}, when applied to \[(p,q,r) \in \MCat(n)\times \MCat(m)\times \MCat(l),\] gives \[p\Bigl(x_{[0,i-1]}, q\bigl(x_0,x_{[i,i+m-1]}\bigr), x_{[i+m,j+m-2]}, r\bigl(x_0,x_{[j+m-1,j+m+l-2]}\bigr), x_{[j+m+l-1,n+m+l-2]}\Bigr)\] in $\MCat(n+m+l-2)$, and similarly for morphisms.  This proves the horizontal associativity axiom.  

To check the vertical associativity axiom, suppose $1 \leq i \leq n$, $1 \leq j \leq m$, and $r\in \Ob(\MCat(l))$.  Each composite in the vertical associativity diagram \eqref{compi-associativity-two}, when applied to $(p,q,r)$ as above, gives\[p\Bigl(x_{[0,i-1]}, q\bigl(x_0,x_{[i,i+j-2]}, r\bigl(x_0,x_{[i+j-1,i+j+l-2]}\bigr), x_{[i+j+l-1,i+l+m-2]}\bigr), x_{[i+l+m-1,n+m+l-2]}\Bigr)\] in $\MCat(n+m+l-2)$, and similarly for morphisms.  This proves the vertical associativity axiom.  
\end{proof}

\begin{example}
The $\compi$-composition $p\compi q$ in $\MCat$ substitutes $q$ for the entry $x_i$ in $p$ and shifts the indices accordingly.  For example, the functor \[\nicexy{\MCat(2)\times\MCat(2) \ar[r]^-{\comp_2} & \MCat(3)}\] yields \[\Bigl\{\bigl((x_1x_0)x_2\bigr)x_0\Bigr\} \comp_2 \Bigl\{x_0\bigl((x_0(x_0x_2))x_1\bigr)\Bigr\} = \Bigl[(x_1x_0) \Bigl(x_0 \bigl((x_0(x_0x_3))x_2\bigr)\Bigr)\Bigr] x_0\]on objects.\dqed
\end{example}

Before we show that $\MCat$-algebras in $\Cat$ are precisely small monoidal categories with general associativity and unit isomorphisms, we first observe that $\MCat$ satisfies the operadic versions of the monoidal category axioms.

\begin{lemma}\label{mcat-monoidal-axioms}
Consider:
\begin{itemize}
\item the object $\mu = x_1x_2 \in \MCat(2)$;
\item the isomorphisms 
\[\begin{split}
& \nicexy{\mu\comp_1 \mu = (x_1x_2)x_3 \ar[r]^-{\alpha}_-{\cong} & x_1(x_2x_3) =\mu \comp_2\mu} \in \MCat(3),\\
& \nicexy{\mu \comp_1 x_0= x_0x_1 \ar[r]^-{\lambda}_-{\cong} & x_1}\in \MCat(1),\\
& \nicexy{\mu \comp_2 x_0= x_1x_0 \ar[r]^-{\rho}_-{\cong} & x_1}\in \MCat(1).
\end{split}\]
\end{itemize}
Then:
\begin{enumerate}
\item The Pentagon Axiom is satisfied in $\MCat$ in the sense that the diagram
\[\nicexy@C-1cm{& (x_1x_2)(x_3x_4) \ar[dr]^-{\alpha \comp_3 \Id_{\mu}} &\\
\bigl((x_1x_2)x_3\bigr)x_4 \ar[ur]^-{\alpha \comp_1 \Id_{\mu}} \ar[d]_-{\Id_{\mu}\comp_1 \alpha} && x_1\bigl(x_2(x_3x_4)\bigr)\\
\bigl(x_1(x_2x_3)\bigr)x_4 \ar[rr]^-{\alpha \comp_2 \Id_{\mu}} && x_1\bigl((x_2x_3)x_4\bigr) \ar[u]_-{\Id_{\mu} \comp_2 \alpha}}\]
in $\MCat(4)$ is commutative.
\item The Unity Axiom is satisfied in $\MCat$ in the sense that the diagram
\[\nicexy@C+.3cm{(x_1x_0)x_2 \ar[r]^-{\alpha \comp_2 \Id_{x_0}} \ar[d]_-{\Id_{\mu}\comp_1\rho} & x_1(x_0x_2) \ar[d]^-{\Id_{\mu}\comp_2\lambda}\\ x_1x_2 \ar@{=}[r] & x_1x_2}\] 
in $\MCat(2)$ is commutative.
\item The morphisms \[\nicexy@C+.3cm{x_0x_0 \ar[r]^-{\lambda\comp_1\Id_{x_0}} & x_0} \andspace \nicexy@C+.3cm{x_0x_0 \ar[r]^-{\rho\comp_1\Id_{x_0}} & x_0}\] in $\MCat(0)$ are equal. 
\end{enumerate}
\end{lemma}

\begin{proof}
All three assertions follow from the fact that each non-empty morphism set in $\MCat$ is a one-point set, so any two morphisms with the same domain and the same codomain are equal.
\end{proof}

To extend from a monoidal category structure to an $\MCat$-algebra structure, we also need to know that morphisms in $\MCat$ can be factored in a specific way, for which we use the following terminology.

\begin{definition}\label{def:mcat-terminology}
Suppose $p\in \Ob(\MCat(n))$.
\begin{enumerate}
\item Suppose $p$ involves the non-associative monomial $(fg)h$, where $f$, $g$, and $h$ are themselves non-associative monomials.  Suppose $q$ is obtained from $p$ by replacing $(fg)h$ by the non-associative monomial $f(gh)$.  The isomorphism \[\alpha : p \iso q \in \MCat(n)\] is called an \emph{associativity isomorphism}.
\item Suppose $p$ involves $x_0f$ for some non-associative monomial $f$, and $r \in \MCat(n)$ is obtained from $p$ by removing this instance of $x_0$ and the associated pair of parentheses around $x_0f$.  The isomorphism \[\lambda : p \iso r \in \MCat(n)\] is called a \emph{left unit isomorphism}.
\item Suppose $p$ involves $fx_0$ for some non-associative monomial $f$, and $s \in \MCat(n)$ is obtained from $p$ by removing this instance of $x_0$ and the associated pair of parentheses around $fx_0$.  The isomorphism \[\rho : p \iso s \in \MCat(n)\] is called a \emph{right unit isomorphism}.
\item Denote by $p_{\st} \in\Ob(\MCat(n))$ the non-associative monomial obtained from $p$ by removing all instances of $x_0$'s and their associated pairs of parentheses.  We call\label{not:pst} $p_{\st}$ the\index{strict non-associative monomial}\index{non-associative monomial!strict} \emph{strict non-associative monomial} associated to $p$.
\item Denote by\label{not:pnor} $p_{\nor}\in\Ob(\MCat(n))$ the non-associative monomial obtained from $p$ by insisting that every pair of parentheses start in the front.  We call $p_{\nor}$ the\index{normal non-associative monomial}\index{non-associative monomial!normal} \emph{normal non-associative monomial} associated to $p$.
\end{enumerate}
\end{definition}

\begin{example}
The following are examples of a left unit isomorphism, a right unit isomorphism, and an associativity isomorphism in $\MCat(3)$:
\[\nicexy{p=p_{\nor}\Bigl(\bigl((x_0x_3)x_2\bigr)x_0\Bigr)x_1 \ar[d]_-{\lambda} & x_3\bigl(x_2x_1).\\
\Bigl(\bigl(x_3x_2\bigr)x_0\Bigr)x_1 \ar[r]^-{\rho} & \bigl(x_3x_2\bigr)x_1=p_{\st} \ar[u]_-{\alpha}}\]\dqed
\end{example}

\begin{lemma}\label{mcat-factorization}
Every morphism $\varphi : p \to q \in \MCat(n)$ admits a unique factorization as
\[\nicexy{p \ar[r]^-{U} & p_{\st} \ar[r]^-{V} & q_{\st} \ar[r]^-{U'} & q \ar@{<-}`u[lll]`[lll]_-{\varphi}[lll]}\]
such that:
\begin{itemize}
\item $U$ is a composite of left unit isomorphisms and right unit isomorphisms.
\item $V$ is a composite of associativity isomorphisms and their inverses.
\item $U'$ is a composite of inverses of left unit isomorphisms and inverses of right unit isomorphisms.
\end{itemize}
\end{lemma}

\begin{proof}
Uniqueness follows from the fact that each non-empty morphism set in $\MCat$ is a one-point set.  

For the existence assertion, first note that there is a unique morphism \[U : p \to p_{\st}\] because they have the same underlying standard monomials.  That $U$ can be factored as stated follows by a finite induction on the set of $x_0$'s in $p$.  The same reasoning applies to $q$, and $U'$ is the inverse of the resulting composite.

For the morphism $V$, note that the assumed existence of the morphism $\varphi$ means that $p$ and $q$ have the same underlying standard monomials.  So $p_{\st}$ and $q_{\st}$ also have the same underlying standard monomials, and there is a unique morphism \[V : p_{\st}\to q_{\st}.\]  That $V$ can be factored as stated can be proved by the classic argument for the Coherence Theorem for monoidal categories in \cite{maclane} (VII.2 Theorem 1), which we briefly recall here for the reader's convenience.  By successively moving parentheses to the left, there is a unique isomorphism \[V_p : p_{\st} \to (p_{\st})_{\nor}\] that is a composite of inverses of associativity isomorphisms.  Likewise, there is a composite of inverses of associativity isomorphisms \[V_q : q_{\st} \to (q_{\st})_{\nor}.\]  Notice that \[(p_{\st})_{\nor} = (q_{\st})_{\nor},\] so we have $V = V_q^{-1}V_p$. 
\end{proof}

\section{Operadic Coherence for Monoidal Categories}
\label{sec:sym-operad-monoidal-cat}

The purpose of this section is to provide an operadic interpretation of the Coherence Theorem for monoidal categories.  More precisely, we show that algebras over the monoidal category operad $\MCat$ are precisely small monoidal categories with general associativity isomorphism, left unit isomorphism, and right unit isomorphism.  We will use the description of an algebra over a symmetric operad in Proposition \ref{prop:symmetric-operad-algebra-defs}.  We continue to use the notation in Lemma \ref{mcat-monoidal-axioms}.

\begin{theorem}\label{mcat-algebra}
For each small category $A$, the following two structures are equivalent:\index{coherence!monoidal category}\index{monoidal category!operadic interpretation}
\begin{enumerate}[label=(\roman*)]
\item A monoidal category structure on $A$.
\item An $\MCat$-algebra structure on $A$.
\end{enumerate}
\end{theorem}

\begin{proof}
Suppose $A$ is an $\MCat$-algebra with structure morphism \[\nicexy{\MCat(n)\times A^{\times n} \ar[r]^-{\theta} & A} \in \Cat\] for $n \geq 0$.  To define the associated monoidal category structure on $A$, consider the structure morphism when $n=0,1,2,3$.
\begin{description}
\item[Monoidal Unit] The structure morphism \[\nicexy{\MCat(0) \ar[r]^-{\theta} & A}\in \Cat\] yields the object \[\tensorunit^A = \theta(x_0)\in A.\]
\item[Monoidal Product] The structure morphism \[\nicexy{\MCat(2)\times A\times A  \ar[r]^-{\theta} & A} \in \Cat\] yields the functor \[\nicexy@C+1cm{A \times A \ar[r]^-{\otimes\,=\,\theta(\mu,-,-)} & A}\] when restricted to the object $\mu=x_1x_2\in \Ob(\MCat(2))$.
\item[Unit Isomorphisms] Using the left and right unit isomorphisms \[\nicexy{x_0x_1 \ar[r]^-{\lambda}_-{\cong} & x_1 & x_1x_0 \ar[l]_-{\rho}^-{\cong}}\] in $\MCat(1)$, the structure morphism \[\nicexy{\MCat(1)\times A \ar[r]^-{\theta} & A} \in \Cat\] yields the natural isomorphisms \[\nicexy@C+1cm{\tensorunit^A \otimes a \ar[r]^-{\lambda^A_a\,=\, \theta(\lambda,a)}_-{\cong} & a & a \otimes \tensorunit^A \ar[l]_-{\rho^A_a\,=\,\theta(\rho,a)}^-{\cong}}\in A\] for objects $a \in A$.
\item[Associativity Isomorphism] Using the associativity isomorphism \[\alpha : (x_1x_2)x_3 \iso x_1(x_2x_3) \in \MCat(3),\] the structure morphism \[\nicexy{\MCat(3)\times A^{\times 3} \ar[r]^-{\theta} & A} \in \Cat\] yields the natural isomorphism \[\nicexy@C+2cm{(a_1\otimes a_2)\otimes a_3 \ar[r]^-{\alpha^A_{a_1,a_2,a_3} \,=\, \theta(\alpha,a_1,a_2,a_3)}_-{\cong} & a_1 \otimes (a_2 \otimes a_3)} \in A\] for objects $a_1,a_2,a_3\in A$.
\end{description}
It follows from Lemma \ref{mcat-monoidal-axioms} that $\bigl(A,\otimes,\tensorunit^A,\alpha^A,\lambda^A,\rho^A)$ satisfies the unity axioms and the pentagon axiom of a monoidal category.

Conversely, suppose given a monoidal category structure $\bigl(A,\otimes,\tensorunit^A,\alpha^A,\lambda^A,\rho^A)$ on $A$.  We define the $\MCat$-algebra structure morphism \[\nicexy{\MCat(n)\times A^{\times n} \ar[r]^-{\theta} & A} \in \Cat\] as follows.  For objects $p=p\bigl(x_{[0,n]}\bigr)\in \MCat(n)$ and $a_i \in A$ for $1\leq i \leq n$, we define 
\begin{equation}\label{theta-p-ua}
\theta\bigl(p; a_1,\ldots,a_n\bigr) = p\bigl(\tensorunit^A,a_1,\ldots,a_n\bigr) \in A.
\end{equation} 
This is obtained from $p$ by
\begin{enumerate}[label=(\roman*)]
\item replacing each instance of $x_0$ by the monoidal unit $\tensorunit^A$,
\item replacing $x_i$ by $a_i$ for $1\leq i \leq n$, and
\item interpreting concatenation as the monoidal product $\otimes$.  
\end{enumerate}
For example, we have the object\[\theta\Bigl((x_0x_3)\bigl((x_1x_0)x_2\bigr); a_1,a_2,a_3\Bigr) = (\tensorunit^A \otimes a_3) \otimes \bigl((a_1 \otimes \tensorunit^A) \otimes a_2\bigr) \in A.\]  Similarly, if each $f_i$ is a morphism in $A$, then we define the morphism \[\theta\bigl(p; f_1,\ldots,f_n\bigr) = p\bigl(\tensorunit^A,f_1,\ldots,f_n\bigr) \in A.\] 

Next we define the structure morphism $\theta$ for morphisms in $\MCat(n)$.  Suppose \[\varphi : p \to q \in \MCat(n)\] is a morphism.  By Lemma \ref{mcat-factorization} it admits a unique decomposition as \[\varphi = U'VU\] in which: 
\begin{itemize}
\item $U$ is a composite of left unit isomorphisms and right unit isomorphisms.
\item $V$ is a composite of associativity isomorphisms and their inverses.
\item $U'$ is a composite of inverses of left unit isomorphisms and inverses of right unit isomorphisms.
\end{itemize}
For objects $a_i\in A$, we abbreviate the sequence $(a_1,\ldots,a_n)$ to $\ua$.  Then the morphism $\theta(\varphi;\ua)$ is the composite 
\begin{equation}\label{theta-varphi-ua}
\nicexy@C+.5cm{p\bigl(\tensorunit^A, \ua\bigr) \ar[r]^-{\theta(\varphi;\, \ua)} \ar[d]_-{\theta(U;\,\ua)} & q(\tensorunit^A, \ua\bigr)\\ p_{\st}(\ua) \ar[r]^-{\theta(V;\,\ua)} & q_{\st}(\ua) \ar[u]_-{\theta(U';\,\ua)}}
\end{equation} 
in $A$ defined as follows.
\begin{itemize}
\item $\theta(U;\ua)$ is the composite of
\begin{itemize}\item the left unit isomorphism $\lambda^A$ tensored with identity morphisms and
\item the right unit isomorphism $\rho^A$ tensored with identity morphisms,
\end{itemize}
corresponding to $U$.
\item $\theta(V;\ua)$ is the composite of
\begin{itemize}\item the associativity isomorphism $\alpha^A$ tensored with identity morphisms and
\item $(\alpha^A)^{-1}$ tensored with identity morphisms,
\end{itemize}
corresponding to $V$. 
\item $\theta(U';\ua)$ is the composite of
\begin{itemize}\item $(\lambda^A)^{-1}$ tensored with identity morphisms and
\item $(\rho^A)^{-1}$ tensored with identity morphisms, 
\end{itemize}
corresponding to $U'$.
\end{itemize}
For example, the composite 
\[\narrowxy{(x_3x_0)(x_1x_2) \ar[r]^-{\rho}_-{\cong} & x_3(x_1x_2) \ar[r]^-{\alpha^{-1}}_-{\cong} & (x_3x_1)x_2 \ar[r]^-{\lambda^{-1}}_-{\cong} & \bigl(x_0(x_3x_1)\bigr)x_2} \]in $\MCat(3)$ is sent by $\theta$ to the composite
\[\nicexy@C+.7cm{(a_3 \otimes \tensorunit^A) \otimes (a_1 \otimes a_2) \ar[d]_-{\rho^A_{a_3} \otimes \Id_{a_1\otimes a_2}} & \bigl(\tensorunit^A \otimes (a_3\otimes a_1)\bigr)\otimes a_2\\
a_3 \otimes (a_1\otimes a_2) \ar[r]^-{(\alpha^A_{a_3,a_1,a_2})^{-1}} & (a_3 \otimes a_1) \otimes a_2 \ar[u]_-{(\lambda^A_{a_3\otimes a_1})^{-1} \otimes \Id_{a_2}}}\]in $A$.  

In general, the morphism $\theta(\varphi;\ua)$ is a categorical composite of monoidal products of the associativity isomorphism $\alpha^A$, the left unit isomorphism $\lambda^A$, the right unit isomorphism $\rho^A$, their inverses, and identity morphisms.  To see that the morphism $\theta(\varphi;\ua)$ is well-defined, we need to show that it is independent of the choices involved in the decomposition $\varphi = U'VU$.  
\begin{itemize}
\item The morphism \[U : p \to p_{\st} \in \MCat(n)\] is a composite of left and right unit isomorphisms that successively remove instances of the variable $x_0$ in the standard non-associative monomial $p$.  Choosing a different presentation of $U$  involves removing the $x_0$'s in $p$ in a different order.  Since $A$ is a monoidal category, the unity axioms in Definition \ref{def:monoidal-category} and the additional unity diagrams \eqref{moncat-other-unit-axioms} imply that the morphism $\theta(U;\ua)$ is well-defined.  A similar argument applies to the morphism \[U' : q_{\st} \to q\in \MCat(n),\] which shows that the morphism $\theta(U';\ua)$ is well-defined.
\item The morphism \[V : p_{\st} \to q_{\st} \in \MCat(n)\] is a composite of associativity isomorphisms and their inverses that move the parentheses.  That the morphism $\theta(V;\ua)$ in $A$ is well-defined is shown by the proof of the Coherence Theorem for monoidal categories in \cite{maclane} (VII.2 Theorem 1).
\end{itemize}
An inspection shows that \[\theta : \MCat(n) \times A^{\times n} \to A\] is a functor.

It remains to check that $(A,\theta)$ is actually an $\MCat$-algebra, i.e., satisfies the associativity and unity axioms in Definition \ref{def:planar-operad-algebra-generating} and the equivariance axiom in Proposition \ref{prop:symmetric-operad-algebra-defs}.  These axioms for $(A,\theta)$ all follow immediately from the definition of $\theta$ ($=$ substituting objects in $A$ for the variables in $p$) and the definition of the operadic composition in $\MCat$ ($=$ substituting standard non-associative monomials for the variables).
\end{proof}

\begin{remark}
The proof of Theorem \ref{mcat-algebra} uses Mac Lane's Coherence Theorem for monoidal categories, so the former does not recover the latter.  On the other hand, in Chapter \ref{ch:coherence-gmonoidal} we will prove several versions of coherence for monoidal categories equipped with an action by an action operad $\G$, again assuming Mac Lane's Coherence Theorem for monoidal categories.  Restricting to appropriate action operads, we recover and unify coherence results for braided monoidal categories, symmetric monoidal categories, ribbon monoidal categories, and coboundary monoidal categories.\dqed
\end{remark}

Recall the symmetric circle product $\circs$ in \eqref{symmetric-circle-product}, from which symmetric operads and their algebras are defined.  Theorem \ref{mcat-algebra} leads to the following operadic interpretation of the free monoidal category.

\begin{corollary}\label{free-monoidal-category}\index{monoidal category!generated by a category}
For each small category $A$, the free monoidal category generated by $A$ is the category \[\MCat \circs A = \coprodover{n\geq 0}\, \MCat(n) \timesover{S_n} A^{\times n}.\]
\end{corollary}

\begin{remark}
In \cite{joyal-street} the free monoidal category is called the\index{free tensor category} \emph{free tensor category}.  The free monoidal category generated by a single object is constructed in \cite{maclane} (VII.2).\dqed
\end{remark}

\begin{remark}[Fresse's Parenthesized Permutation Operad]\label{rk:pap-operad}
Theorem \ref{mcat-algebra} says that $\MCat$-algebras are precisely small monoidal categories with general associativity isomorphisms and left/right unit isomorphisms.  In \cite{fresse-gt} Theorem 6.1.7(a-b) Fresse described a $1$-colored symmetric operad\label{not:fresse-pap} \textit{\textsf{PaP}} in $\Cat$, called the\index{parenthesized permutation operad} \emph{parenthesized permutation operad}, whose algebras are small monoidal categories with general associativity isomorphisms but without monoidal units.  To obtain Fresse's parenthesized permutation operad \textit{\textsf{PaP}}, we take the monoidal category operad $\MCat$ and forget about
\begin{enumerate}[label=(\roman*)]
\item the $0$th level $\MCat(0)$ and
\item all instances of $x_0$'s in $\MCat(n)$ for $n\geq 1$.
\end{enumerate}

A variation of \textit{\textsf{PaP}} is the $1$-colored symmetric operad\label{not:fresse-papplus} \textit{\textsf{PaP}}$_+$ in $\Cat$ in \cite{fresse-gt} Theorem 6.1.7(c), called the \emph{unitary parenthesized permutation operad}.  Its algebras are small monoidal categories with (i) general associativity isomorphisms and (ii)  monoidal units that are strict with respect to both the monoidal product and the associativity isomorphism, i.e., 
\[\begin{split}
\lambda_x &= \Id_x= \rho_x,\\
\alpha_{\tensorunit,x,y} &= \alpha_{x,\tensorunit,y} = \alpha_{x,y,\tensorunit} = \Id_{x\otimes y}.
\end{split}\] 
This symmetric operad is obtained from the parenthesized permutation operad \textit{\textsf{PaP}} by adjoining an object \textit{\textsf{PaP}}$_+(0) = *$.  

There is a morphism $\pi : \MCat \to$ \!\textit{\textsf{PaP}}$_+$ of one-colored symmetric operads in $\Cat$ that sends each standard non-associative monomial of weight $0$ to $* \in$ \textit{\textsf{PaP}}$_+(0)$.  For $n\geq 1$, $\pi$ sends each standard non-associative monomial $p$ of weight $n$ to the strict non-associative monomial $p_{\st}$ associated to $p$.\dqed
\end{remark}

\section{Strict Monoidal Functors as Algebra Morphisms}
\label{sec:strict-monoidal-functor}

Recall that a monoidal functor $F$ is \emph{strict} if the structure morphisms $F_2$ \eqref{monoidal-f2} and $F_0$ \eqref{monoidal-f0} are all identity morphisms.  The purpose of this section is to observe that $\MCat$-algebra morphisms are exactly strict monoidal functors.

\begin{theorem}\label{mcat-algebra-morphism}
A morphism \[F : (A,\theta^A) \to (B,\theta^B)\] of $\MCat$-algebras is equivalent to a strict monoidal functor $F : A\to B$.
\end{theorem}

\begin{proof}
The assumption on $F$ being an $\MCat$-algebra morphism means that it is a functor $F : A \to B$ such that the diagram
\begin{equation}\label{mcat-algebra-morphism-diagram}
\nicexy@C+.4cm{\MCat(n) \times A^{\times n} \ar[d]_-{\theta^A} \ar[r]^-{(\Id,F^{\times n})} & \MCat(n) \times B^{\times n} \ar[d]^-{\theta^B}\\ A \ar[r]^-{F} & B}
\end{equation}
in $\Cat$ is commutative for $n\geq 0$.  
\begin{itemize}
\item When $n=0$, this diagram, when applied to the object $x_0\in \MCat(0)$, implies the equality \[F(\tensorunit^A) = \tensorunit^B.\]  
\item When $n=2$, the above diagram, when applied to the object \[\mu=x_1x_2\in \MCat(2),\] implies the equality \[F(a_1\otimes a_2) = F(a_1)\otimes F(a_2)\] for objects $a_1,a_2\in A$.  
\item When $n=1$, the above diagram, when applied to the left and right unit isomorphisms \[\nicexy{x_0x_1\ar[r]^-{\lambda}_-{\cong} & x_1 & x_1x_0 \ar[l]_-{\rho}^-{\cong}}\] in $\MCat(1)$, yields the equalities \[F(\lambda^A_a)  = \lambda^B_{F(a)} \andspace F(\rho^A_a) = \rho^B_{F(a)}\] for objects $a\in A$.
\item When $n=3$, the above diagram, when applied to the associativity isomorphism \[\nicexy{(x_1x_2)x_3 \ar[r]^-{\alpha}_-{\cong} & x_1(x_2x_3)} \in \MCat(3),\] yields the equality \[F\bigl(\alpha^A_{a_1,a_2,a_3}\bigr) = \alpha^B_{F(a_1), F(a_2), F(a_3)}\] for objects $a_1,a_2,a_3 \in A$.
\end{itemize}
Therefore, $F$ has the structure of a strict monoidal functor $F : A \to B$.

Conversely, suppose $F : A \to B$ is a strict monoidal functor.  To see that the diagram \eqref{mcat-algebra-morphism-diagram} is commutative, suppose
\begin{itemize}
\item $p \in \Ob(\MCat(n))$ and 
\item $\ua = (a_1,\ldots,a_n) \in A^{\times n}$.  
\end{itemize}
Using \eqref{theta-p-ua} and the assumption that $F$ is a strict monoidal functor, there are equalities
\[\begin{split}
F\bigl(\theta^A(p;\ua)\bigr) &= F\bigl(p(\tensorunit^A; \ua)\bigr)\\
&= p\bigl(F\tensorunit^A; Fa_1,\ldots, Fa_n\bigr)\\
&= p\bigl(\tensorunit^B; Fa_1,\ldots, Fa_n\bigr)\\
&= \theta^B\bigl(p; Fa_1, \ldots, Fa_n\bigr).
\end{split}\]
This shows that the diagram \eqref{mcat-algebra-morphism-diagram} is commutative on objects.  The same argument, with each $a_i$ replaced by a morphism in $A$, shows that the diagram \eqref{mcat-algebra-morphism-diagram} is commutative on morphisms in $A$.  On the other hand, suppose $p$ is replaced by a morphism $\varphi : p \to q$ in $\MCat(n)$.  We use the decomposition
\[\nicexy{\theta^A(\varphi; \ua) = \theta^A(U';\ua) \circ \theta^A(V;\ua) \circ \theta^A(U;\ua)}\]
in \eqref{theta-varphi-ua} and the assumption that $F$ is a strict monoidal functor to infer that the diagram \eqref{mcat-algebra-morphism-diagram} is commutative on morphisms in $\MCat$.
\end{proof}

Theorem \ref{mcat-algebra} and Theorem \ref{mcat-algebra-morphism} together mean that $\MCat$ is the one-colored symmetric operad in $\Cat$ whose algebras are small monoidal categories with general associativity and unit isomorphisms, and whose algebra morphisms are strict monoidal functors.

\section{Strict Monoidal Category Operad}\label{sec:mcatst}

Recall that a strict monoidal category is a monoidal category in which the associativity isomorphism, the left unit isomorphism, and the right unit isomorphism are all identity morphisms.  The purpose of this section is to discuss the one-colored symmetric operad in $\Cat$ whose algebras are strict monoidal categories.  This is needed later when we discuss $\G$-monoidal categories, both general and strict, for an action operad $\G$.

Recall from Definition \ref{def:nonass-monomial} that a \emph{standard monomial of weight $n$} is an associative monomial of the form \[x_{\sigma(1)}\cdots x_{\sigma(n)}\] for some permutation $\sigma \in S_n$.  Also recall that a\index{discrete category}\index{category!discrete} \emph{discrete category} is a category in which there are no non-identity morphisms.  We now define the strict version of the monoidal category operad $\MCat$.

\begin{definition}\label{def:mcatst}
Define the\index{strict monoidal category!operad}\index{monoidal category operad!strict}\index{operad!strict monoidal category} \emph{strict monoidal category operad} $\MCatst$ as the $1$-colored symmetric operad in $\Cat$ as follows.
\begin{description}
\item[Entries]
For each $n\geq 0$, $\MCatst(n)$ is the discrete category whose object set is the set of standard monomials of weight $n$.
\item[Equivariance]
The symmetric group $S_n$ acts on the objects in $\MCatst(n)$ by permuting the variables $\{x_1,\ldots,x_n\}$.
\item[Unit]
The operadic unit is the object $x_1 \in \MCat(1)$.
\item[Composition]
For $1\leq i \leq n$ and $m \geq 0$, the $\compi$-composition \[\nicexy{\MCatst(n) \times \MCatst(m) \ar[r]^-{\compi} & \MCatst(n+m-1)}\] is the functor defined as \[p\bigl(x_{[1,n]}\bigr) \compi q\bigl(x_{[1,m]}\bigr) = p\Bigl(x_{[1,i-1]}, q\bigl(x_{[i,i+m-1]}\bigr), x_{[i+m,n+m-1]}\Bigr)\]
for standard monomials $p\bigl(x_{[1,n]}\bigr)\in \MCatst(n)$ and $q\bigl(x_{[1,m]}\bigr) \in \MCatst(m)$.
\end{description}
\end{definition}

\begin{example}
In $\MCatst$:
\begin{itemize}
\item $\MCatst(0)$ has only one object $\varnothing$.
\item $\MCatst(1)$ has only one object $x_1$.
\item $\MCatst(2)$ has two objects $x_1x_2$ and $x_2x_1$.
\item $\MCatst(3)$ has six objects: $x_1x_2x_3$, $x_1x_3x_2$, $x_2x_1x_3$, $x_3x_1x_2$, $x_2x_3x_1$, and $x_3x_2x_1$.
\end{itemize}
In general, the object set in $\MCatst(n)$ is canonically isomorphic to the set $S_n$ of permutations, with the standard monomial $x_{\sigma(1)}\cdots x_{\sigma(n)}$ identified with the permutation $\sigma \in S_n$.  Under this identification, $\MCatst(n)$ isomorphic to the translation category of $S_n$; see Definition \ref{def:translation-category}.\dqed
\end{example}

The following result is the strict version of Lemma \ref{mcat-is-operad}, whose proof is reused here essentially verbatim.

\begin{lemma}\label{mcatst-is-operad}
$\MCatst$ is a $1$-colored symmetric operad in $\Cat$.
\end{lemma}

The following observation is the strict version of Theorem \ref{mcat-algebra}

\begin{theorem}\label{mcatst-algebra}\index{strict monoidal category!operadic interpretation}
For each small category $A$, the following two structures are equivalent:
\begin{enumerate}[label=(\roman*)]
\item A strict monoidal category structure on $A$.
\item An $\MCatst$-algebra structure on $A$.
\end{enumerate}
\end{theorem}

\begin{proof}
We reuse the proof of Theorem \ref{mcat-algebra} by interpreting $x_0$ as the empty sequence $\varnothing$ and forgetting all the parentheses.
\end{proof}

Theorem \ref{mcatst-algebra} leads to the following operadic interpretation of the free strict monoidal category.

\begin{corollary}\label{free-strict-monoidal-category}
For each small category $A$, the free strict monoidal category generated by $A$ is the category \[\MCatst \circs A = \coprodover{n\geq 0}\, \MCatst(n) \timesover{S_n} A^{\times n}.\]
\end{corollary}

The following result is proved by reusing the proof of Theorem \ref{mcat-algebra-morphism}.

\begin{proposition}\label{mcatst-algebra-morphism}
A morphism \[F : (A,\theta^A) \to (B,\theta^B)\] of $\MCatst$-algebras is equivalent to a strict monoidal functor $F : A\to B$.
\end{proposition}

Theorem \ref{mcatst-algebra} and Proposition \ref{mcatst-algebra-morphism} together mean that $\MCatst$ is the one-colored symmetric operad in $\Cat$ whose algebras are small strict monoidal categories and whose algebra morphisms are strict monoidal functors.

\begin{remark}
The strict monoidal category operad $\MCatst$ is the same as Fresse's\label{not:copplus} \textit{\textsf{CoP}}$_+$, called the \emph{unitary colored permutation operad} in \cite{fresse-gt} Section 6.1.6.  Theorem \ref{mcatst-algebra} is essentially the same as \cite{fresse-gt} Theorem 6.1.8.\dqed
\end{remark}

The following is an exercise in using the morphism part of Proposition \ref{def:g-operad-generating} for the symmetric group operad $\S$.

\begin{lemma}\label{mcat-mcatst}
There is a morphism of $1$-colored symmetric operads \[\pi : \MCat \to \MCatst\] in $\Cat$ that sends each standard non-associative monomial $p \in \Ob(\MCat(n))$ to its underlying standard monomial $\pbar\in \Ob(\MCatst(n))$.
\end{lemma}

The following is Theorem \ref{changing-goperad} (2-3) for the symmetric group operad $\S$ and the morphism $\pi : \MCat \to \MCatst$ of one-colored symmetric operads in $\Cat$.

\begin{corollary}\label{mcatalg-mcatstalg}
The morphism $\pi : \MCat \to \MCatst$ induces an adjunction
\[\nicexy{\algm\bigl(\MCat\bigr) \ar@<2pt>[r]^-{\pi_!} & \algm\bigl(\MCatst\bigr) \ar@<2pt>[l]^-{\pi^*}}\] between the category of $\MCat$-algebras and the category of $\MCatst$-algebras.
\end{corollary}

\begin{interpretation}
Interpreting the adjunction $\pi_! \dashv \pi^*$ using Theorem \ref{mcat-algebra} and Theorem \ref{mcatst-algebra}, the right adjoint $\pi^*$ sends a strict monoidal category to itself, which is a monoidal category.  The left adjoint $\pi_!$ strictifies each monoidal category to a strict monoidal category by forcing the associativity isomorphism and the unit isomorphisms to be identity morphisms.\dqed
\end{interpretation}

\section{Monoidal Functor Operad}\label{sec:operad-monoidal-functor}

Recall from Theorem \ref{mcat-algebra} that $\MCat$ is a $1$-colored symmetric operad in $\Cat$ whose algebras are exactly small monoidal categories with general associativity isomorphism and left/right unit isomorphisms.  On the other hand, by Theorem \ref{mcat-algebra-morphism}, $\MCat$-algebra morphisms are strict monoidal functors, in which the structure morphisms $F_2$ and $F_0$ are identity morphisms.  

The purpose of this section is to describe a $2$-colored symmetric operad in $\Cat$ for which an algebra is a triple $(A,B,F)$ such that:
\begin{itemize}\item $A$ and $B$ are general small monoidal categories.
\item $F : A\to B$ is a general monoidal functor, in which the structure morphisms $F_2$ and $F_0$ are not required to be isomorphisms.  
\end{itemize}
We will also discuss a variation of this $2$-colored symmetric operad in $\Cat$ for which an algebra is a triple $(A,B,F)$ as above but with $F$ a strong monoidal functor, in which the structure morphisms $F_2$ and $F_0$ are isomorphisms.

Recall from Example \ref{ex:underlying-object-operad} that, for each action operad $\G$, the underlying object functor \[\Ob : \Cat \to \Set\] sends each $\colorc$-colored $\G$-operad in $\Cat$ to its underlying object $\colorc$-colored $\G$-operad in $\Set$.  We will use this below when $\G$ is the symmetric group operad $\S$.  Specifically, we will define a $2$-colored symmetric operad $\MFun$ in $\Cat$ by first defining the underlying object $2$-colored symmetric operad $\Ob\bigl(\MFun\bigr)$ in $\Set$.

Also recall from Definition \ref{def:nonass-monomial} that for a standard non-associative monomial $p$ of weight $n$, its underlying standard monomial is denoted by $\pbar$.

\begin{notation}
In Definition \ref{def:monoidal-functor-operad} below, we use the color set $\{s,t\}$, with $s$ and $t$ standing for \emph{source} and \emph{target}, respectively.  
\begin{enumerate}[label=(\roman*)]
\item For $n \geq 0$ and $u\in \{s,t\}$, we write \[u^n = (u,\ldots,u)\] for the sequence with $n$ copies of $u$.
\item In addition to the variables $\{x_i\}_{i \geq 0}$ in Definition \ref{def:nonass-monomial}, which will be used for the source, we pick another sequence $\{x_i'\}_{i\geq 0}$ of variables, which will be used for the target.
\end{enumerate}
\end{notation}

\begin{definition}\label{def:monoidal-functor-operad}
Define the\index{monoidal functor operad}\index{operad!monoidal functor} \emph{monoidal functor operad} $\MFun$ as the $\{s,t\}$-colored symmetric operad in $\Cat$ as follows.
\begin{description}
\item[Objects] The underlying object $\{s,t\}$-colored symmetric operad $\Ob(\MFun)$ in $\Set$ is operadically generated by the elements:
\begin{enumerate}[label=(\roman*)]
\item $x_0 \in \Ob\left(\MFun\sbinom{s}{\varnothing}\right)$ and $\mu \in \Ob\left(\MFun\sbinom{s}{s,s}\right)$;
\item $x_0' \in \Ob\left(\MFun\sbinom{t}{\varnothing}\right)$ and $\mu' \in \Ob\left(\MFun\sbinom{t}{t,t}\right)$;
\item $f \in \Ob\left(\MFun\sbinom{t}{s}\right)$.
\end{enumerate}
Denote the $s$-colored and the $t$-colored operadic units by, respectively, \[x_1\in \Ob\left(\MFun\sbinom{s}{s}\right) \andspace x_1'\in \Ob\left(\MFun\sbinom{t}{t}\right).\]  
\item[Notation]
We use non-associative monomial notation for these objects with \[\mu=x_1x_2 \andspace \mu'=x_1'x_2'.\]  For example, we have
\[\begin{split}
&\mu(1,2) = x_2x_1 \in \Ob\left(\MFun\sbinom{s}{s,s}\right),\\
&\mu\comp_1 x_0 = x_0x_1 \in \Ob\left(\MFun\sbinom{s}{s}\right),\\
&\mu \comp_1 \mu = (x_1x_2)x_3, \quad \mu \comp_2 \mu = x_1(x_2x_3)  \in \Ob\left(\MFun\sbinom{s}{s,s,s}\right),\\
&f\comp_1 x_0 = f(x_0) \in \Ob\left(\MFun\sbinom{t}{\varnothing}\right),\\
&f\comp_1 \mu = f(x_1x_2) \in \Ob\left(\MFun\sbinom{t}{s,s}\right),\\
&\bigl(\mu' \comp_1 f\bigr) \comp_2 f = f(x_1)f(x_2) \in \Ob\left(\MFun\sbinom{t}{s,s}\right),\\
&\bigl(\mu'\comp_1 x_0')\comp_2 f = x_0'f(x_1) \in \Ob\left(\MFun\sbinom{t}{s}\right),\\
&\bigl(((\mu'\comp_1\mu') \comp_1 f)\comp_2 f\bigr)\comp_3 f = \bigl(f(x_1)f(x_2)\bigr)f(x_3) \in \Ob\left(\MFun\sbinom{t}{s,s,s}\right),\\
&\bigl(\mu' \comp_2 (f\comp_1 \mu)\bigr) \comp_1 f = f(x_1) f(x_2x_3) \in \Ob\left(\MFun\sbinom{t}{s,s,s}\right),\\
&f\comp_1\bigl(\mu \comp_2\mu\bigr) = f\bigl(x_1(x_2x_3)\bigr) \in \Ob\left(\MFun\sbinom{t}{s,s,s}\right), 
\end{split}\]
and so forth.
\item[Morphisms] The morphism sets in the categories in $\MFun$ are categorically and operadically generated by:
\begin{enumerate}[label=(\roman*)]
\item an isomorphism \[*_s^{p,q} \in \MFun\sbinom{s}{s^n}(p,q)\] for each pair of standard non-associative monomials $p,q$ of weight $n \geq 0$ in $\{x_0,\ldots,x_n\}$ with $\pbar=\qbar$;
\item an isomorphism \[*_t^{p,q} \in \MFun\sbinom{t}{t^n}(p,q)\] for each pair of standard non-associative monomials $p,q$ of weight $n \geq 0$ in $\{x_0',\ldots,x_n'\}$ with $\pbar=\qbar$;
\item a morphism \[f_2 \in \MFun\sbinom{t}{s,s}\bigl(f(x_1)f(x_2), f(x_1x_2)\bigr);\]
\item a morphism \[f_0 \in \MFun\sbinom{t}{\varnothing}\bigl(x_0',f(x_0)\bigr).\]
\end{enumerate}
We will sometimes omit $p,q$ in $*_s^{p,q}$ and $*_t^{p,q}$ if there is no danger of confusion.
\item[Relations] The above data is subject to the following five types of relations.
\begin{description}
\item[Source] The assignments \[\begin{cases}* \mapsto s \in \{s,t\},\\ \Ob\bigl(\MCat(n)\bigr) \ni p \mapsto p \in \Ob\bigl(\MFun\sbinom{s}{s^n}\bigr),\\ \MCat(n)(p,q) \ni * \mapsto *_s^{p,q} \in \MFun\sbinom{s}{s^n}(p,q)\end{cases}\] define a morphism 
\begin{equation}\label{iota-source}
\iota_s : \MCat \to \MFun
\end{equation}
of symmetric operads in $\Cat$ in the sense of Definition \ref{def:g-operads}.
\item[Target] The assignments \[\begin{cases}* \mapsto t \in \{s,t\},\\ \Ob\bigl(\MCat(n)\bigr) \ni p(x_{[0,n]}) \mapsto p(x'_{[0,n]}) \in \Ob\bigl(\MFun\sbinom{t}{t^n}\bigr), \\ \MCat(n)(p,q) \ni * \mapsto *_t^{p,q} \in \MFun\sbinom{t}{t^n}\bigl(p(x'_{[0,n]}), q(x'_{[0,n]})\bigr)\end{cases}\] 
define a morphism 
\begin{equation}\label{iota-target}
\iota_t : \MCat \to \MFun
\end{equation} 
of symmetric operads in $\Cat$.
\item[Associativity] Writing
\[\begin{split}
\alpha &= *_s \in \MFun\sbinom{s}{s,s,s}\bigl((x_1x_2)x_3,x_1(x_2x_3)\bigr),\\
\alpha' &= *_t \in \MFun\sbinom{t}{t,t,t}\bigl((x_1'x_2')x_3',x_1'(x_2'x_3')\bigr),
\end{split}\]
the diagram
\begin{equation}\label{mfun-associativity}
\nicexy@C+.5cm{\bigl(f(x_1)f(x_2)\bigr)f(x_3) \ar[d]_-{\mu'(f_2,f)} \ar[r]^-{\alpha'(f,f,f)} & f(x_1)\bigl(f(x_2)f(x_3)\bigr) \ar[d]^-{\mu'(f,f_2)}\\
f(x_1x_2)f(x_3) \ar[d]_-{f_2\comp_1\mu} & f(x_1)f(x_2x_3) \ar[d]^-{f_2\comp_2\mu}\\
f\bigl((x_1x_2)x_3\bigr) \ar[r]^-{f(\alpha)} & f\bigl(x_1(x_2x_3)\bigr)}
\end{equation}
in $\MFun\sbinom{t}{s,s,s}$ is commutative. 
\item[Left Unity] Writing
\[\begin{split}
\lambda &= *_s \in \MFun\sbinom{s}{s}\bigl(x_0x_1,x_1\bigr),\\
\lambda' &= *_t \in \MFun\sbinom{t}{t}\bigl(x_0'x_1',x_1'\bigr),
\end{split}\]
the diagram
\begin{equation}\label{mfun-left-unity}
\nicexy@C+.3cm{x_0'f(x_1) \ar[d]_-{\mu'(f_0,f)} \ar[r]^-{\lambda'(f)} & f(x_1)\\
f(x_0)f(x_1) \ar[r]^-{f_2\comp_1 x_0} & f(x_0x_1) \ar[u]_-{f(\lambda)}}
\end{equation}
in $\MFun\sbinom{t}{s}$ is commutative.
\item[Right Unity] Writing
\[\begin{split}
\rho &= *_s \in \MFun\sbinom{s}{s}\bigl(x_1x_0,x_1\bigr),\\
\rho' &= *_t \in \MFun\sbinom{t}{t}\bigl(x_1'x_0',x_1'\bigr),
\end{split}\]
the diagram
\begin{equation}\label{mfun-right-unity}
\nicexy@C+.3cm{f(x_1)x_0' \ar[d]_-{\mu'(f,f_0)} \ar[r]^-{\rho'(f)} & f(x_1)\\
f(x_1)f(x_0) \ar[r]^-{f_2\comp_2 x_0} & f(x_1x_0) \ar[u]_-{f(\rho)}}
\end{equation}
in $\MFun\sbinom{t}{s}$ is commutative.
\end{description}
\end{description}
This completes the definition of the monoidal functor operad $\MFun$.
\end{definition}

We now observe that $\MFun$ is the symmetric operad that precisely captures the concept of a general monoidal functor, in which the structure morphisms are not necessarily invertible.  We will use the descriptions of an algebra over a symmetric operad in both Proposition \ref{prop:symmetric-operad-algebra-defs} and Example \ref{ex:endomorphism-operad}.

\begin{theorem}\label{mfun-algebra}\index{monoidal functor operad}\index{operad!monoidal functor}
An $\MFun$-algebra is exactly a triple $(A,B,F)$ such that:
\begin{itemize}
\item $A$ and $B$ are small monoidal categories.
\item $F : A \to B$ is a monoidal functor.
\end{itemize}  
\end{theorem}
 
\begin{proof}
Suppose $(A,\theta)$ is an $\MFun$-algebra.  Since the color set of $\MFun$ is $\{s,t\}$, it follows that $A \in \Cat^{\{s,t\}}$.  In other words, $A$ consists of two small categories $A_s$ and $A_t$.  By pre-composing with the morphism \[\iota_s : \MCat \to \MFun\] in \eqref{iota-source} of symmetric operads in $\Cat$, $(A,\theta)$ induces an $\MCat$-algebra \[(A_s,\theta_s = \theta\iota_s).\]  By Theorem \ref{mcat-algebra}, this $\MCat$-algebra structure on $A_s$ is exactly a monoidal category structure on $A_s$.  

Similarly, by pre-composing with the morphism \[\iota_t : \MCat\to \MFun\] in \eqref{iota-target} of symmetric operads in $\Cat$, $(A,\theta)$ induces an $\MCat$-algebra \[(A_t,\theta_t = \theta\iota_t),\] which by Theorem \ref{mcat-algebra} is equivalent to a monoidal category structure on $A_t$.  In the next paragraph, we use the notations in the proof of Theorem \ref{mcat-algebra} to denote the monoidal category structures in $A_s$ and $A_t$.
 
The $\MFun$-algebra structure morphism \[\nicexy{\MFun\sbinom{t}{s} \times A_s \ar[r]^-{\theta} & A_t}\in \Cat,\] when applied to the generating object $f \in \Ob\left(\MFun\sbinom{t}{s}\right)$, yields a functor \[\nicexy@C+.7cm{A_s \ar[r]^-{F \,=\, \theta(f,-)} & A_t}.\]  Moreover, the $\MFun$-algebra structure morphism \[\nicexy{\MFun\sbinom{t}{s,s}\times A_s \times A_s \ar[r]^-{\theta} & A_t},\] when applied to the generating morphism \[f_2 \in \MFun\sbinom{t}{s,s}\bigl(f(x_1)f(x_2), f(x_1x_2)\bigr),\]yields a natural transformation\[\nicexy@C+1cm{F(a_1)\otimes F(a_2) \ar[r]^-{F_2\,=\,\theta(f_2,a_1,a_2)} & F(a_1\otimes a_2)} \in A_t\] for objects $a_1,a_2\in A_s$.  Furthermore, the $\MFun$-algebra structure morphism \[\nicexy{\MFun\sbinom{t}{\varnothing} \ar[r]^-{\theta} & A_t},\]when applied to the generating morphism \[f_0 \in \MFun\sbinom{t}{\varnothing}\bigl(x_0',f(x_0)\bigr),\] yields a morphism \[\nicexy@C+.5cm{\tensorunit^{A_t} \ar[r]^-{F_0\,=\, \theta(f_0)} & F(\tensorunit^{A_s})} \in A_t.\] The associativity axiom \eqref{f2}, the left unity axiom \eqref{f0-left}, and the right unity axiom \eqref{f0-right} for a monoidal functor are satisfied by the triple $(F,F_2,F_0)$ by the associativity relation \eqref{mfun-associativity}, the left unity relation \eqref{mfun-left-unity}, and the right unity relation \eqref{mfun-right-unity}, respectively, in the monoidal functor operad $\MFun$.  So \[(F,F_2,F_0) : A_s \to A_t\] is a monoidal functor.
 
Conversely, suppose $F : A \to B$ is a monoidal functor between small monoidal categories $A$ and $B$.  From this data, we define an $\MFun$-algebra $\bigl(\{A,B\},\theta\bigr)$ as follows.  
\begin{description}
\item[Underlying Categories] The pair of categories $\{A,B\}$ is an object in $\Cat^{\{s,t\}}$, which is the underlying object of the desired $\MFun$-algebra.
\item[Objects] To define the $\MFun$-algebra structure morphism $\theta$, we first consider the generating objects in $\MFun$.
\begin{itemize}
\item For the generating object $x_0 \in \Ob\left(\MFun\sbinom{s}{\varnothing}\right)$, we define $\theta(x_0)$ as the monoidal unit $\tensorunit^A$ in $A$.
\item For the generating object $x_0' \in \Ob\left(\MFun\sbinom{t}{\varnothing}\right)$, we define $\theta(x_0')$ as the monoidal unit $\tensorunit^B$ in $B$.
\item For the generating object $\mu \in \Ob\left(\MFun\sbinom{s}{s,s}\right)$, we define \[\theta(\mu,-,-) :  A \times A \to A\] as the monoidal product $\otimes^A$ in $A$.
\item For the generating object $\mu' \in \Ob\left(\MFun\sbinom{t}{t,t}\right)$, we define \[\theta(\mu',-,-) :  B \times B \to B\] as the monoidal product $\otimes^B$ in $B$.
\item For the generating object $f \in \Ob\left(\MFun\sbinom{t}{s}\right)$, we define \[\theta(f,-) : A \to B\] as the functor $F$.
\end{itemize}
\item[Morphisms] Next we consider the generating morphisms in $\MFun$.
\begin{itemize}
\item For the generating isomorphism \[*_s \in \MFun\sbinom{s}{s^n}(p,q),\] we define $\theta(*_s)$ as the unique natural isomorphism \[\nicexy{p\bigl(\tensorunit^A,a_1,\ldots,a_n\bigr) \ar[r]^-{\cong} & q\bigl(\tensorunit^A,a_1,\ldots,a_n\bigr)} \in A\] for objects $a_1,\ldots,a_n\in A$, that is a categorical composite of monoidal products of identity morphisms, the associativity isomorphism $\alpha^A$, the left unit isomorphism $\lambda^A$, the right unit isomorphism $\rho^A$, and their inverses.  Indeed, since $p$ and $q$ are standard non-associative monomials of weight $n$ with the same underlying standard monomial, such a unique natural isomorphism is guaranteed to exist by the Coherence Theorem for monoidal categories in \cite{maclane} (VII.2 Corollary).
\item For the generating isomorphism \[*_t \in \MFun\sbinom{t}{t^n}(p,q),\] we similarly define $\theta(*_t)$ as the unique natural isomorphism \[\nicexy{p\bigl(\tensorunit^B,b_1,\ldots,b_n\bigr) \ar[r]^-{\cong} & q\bigl(\tensorunit^B,b_1,\ldots,b_n\bigr)} \in B\] for objects $b_1,\ldots,b_n\in B$, that is made up of identity morphisms, $\alpha^B$, $\lambda^B$, $\rho^B$, and their inverses.
\item For the generating morphism \[f_2 \in \MFun\sbinom{t}{s,s}\bigl(f(x_1)f(x_2), f(x_1x_2)\bigr),\] we define $\theta(f_2)$ as the natural transformation \[\nicexy{F(a_1)\otimes F(a_2) \ar[r]^-{F_2} & F\bigl(a_1\otimes a_2\bigr)} \in B\] for objects $a_1,a_2\in A$, which is part of the structure of the given monoidal functor $F$.
\item For the generating morphism \[f_0 \in \MFun\sbinom{t}{\varnothing}\bigl(x_0',f(x_0)\bigr),\] we define $\theta(f_0)$ as the morphism \[\nicexy{\tensorunit^B \ar[r]^-{F_0} & F(\tensorunit^A)} \in B,\]which is also part of the structure of the given monoidal functor $F$.
\end{itemize}
\end{description}

To see that the previous paragraph defines an $\MFun$-algebra $\bigl(\{A,B\},\theta\bigr)$, it suffices to show that the relations in Definition \ref{def:monoidal-functor-operad} are preserved.
\begin{itemize}
\item The source relation is preserved by Theorem \ref{mcat-algebra}, which implies that the monoidal category structure on $A$ is the same as an $\MCat$-algebra structure on $A$.
\item Similarly, the target relation is preserved by Theorem \ref{mcat-algebra} applied to the monoidal category structure on $B$.
\item The associativity relation \eqref{mfun-associativity}, the left unity relation \eqref{mfun-left-unity}, and the right unity relation \eqref{mfun-right-unity} are preserved because $F$ is a monoidal functor.
\end{itemize} 
Therefore, each monoidal functor between small monoidal categories defines an $\MFun$-algebra.
\end{proof}

Recall that a \emph{strong} monoidal functor is a monoidal functor $F$ in which the structure morphisms $F_2$ and $F_0$ are invertible.  Strong monoidal functors are what Joyal and Street \cite{joyal-street} called \emph{tensor functors}.  We end this section with the description of a variation of the monoidal functor operad that captures the concept of a strong monoidal functor.

\begin{definition}\label{def:strong-monoidal-functor-operad}
The\index{operad!strong monoidal functor}\index{monoidal functor operad!strong} \emph{strong monoidal functor operad} $\MFunstg$ is the $\{s,t\}$-colored symmetric operad in $\Cat$ with the same definition as the monoidal functor operad $\MFun$, except for:
\begin{enumerate}[label=(\roman*)]
\item two more generating morphisms
\[\begin{split} f_2^{-1} &\in \MFunstg\sbinom{t}{s,s}\bigl(f(x_1x_2), f(x_1)f(x_2)\bigr),\\
f_0^{-1} &\in \MFunstg\sbinom{t}{\varnothing}\bigl(f(x_0), x_0'\bigr);\end{split}\]
\item four more relations \[\begin{cases} f_2^{-1}f_2 = \Id_{f(x_1)f(x_2)},\\ f_2f_2^{-1} = \Id_{f(x_1x_2)},\end{cases}\qquad \begin{cases} f_0^{-1}f_0 = \Id_{x_0'},\\ f_0f_0^{-1} = \Id_{f(x_0)}.\end{cases}\] 
\end{enumerate}
\end{definition}

The following result is the strong monoidal functor analogue of Theorem \ref{mfun-algebra}, whose proof can be reused here with only cosmetic changes.

\begin{theorem}\label{mfunstg-algebra}\index{operad!strong monoidal functor}\index{strong monoidal functor!operad}
An $\MFunstg$-algebra is exactly a triple $(A,B,F)$ such that:
\begin{itemize}
\item $A$ and $B$ are small monoidal categories.
\item $F : A \to B$ is a strong monoidal functor.
\end{itemize}  
\end{theorem}

\chapter{$\G$-Monoidal Categories}\label{ch:gmonoidal-cat}

Fix an action operad $(\G,\omega)$ as in Definition \ref{def:augmented-group-operad}.  In this chapter, we introduce $\G$-monoidal categories as general small monoidal categories in which the action operad $\G$ acts on iterated monoidal products.  As we will see later, $\G$-monoidal categories specialize to various types of small monoidal categories with extra structures, including:
\begin{enumerate}[label=(\roman*)]
\item monoidal categories when $\G$ is the planar group operad $\P$ in Example \ref{ex:trivial-group-operad};
\item symmetric monoidal categories when $\G$ is the symmetric group operad $\S$ in Example \ref{ex:symmetric-group-operad};
\item braided monoidal categories when $\G$ is the braid group operad $\B$ in Definition \ref{def:braid-group-operad};
\item ribbon monoidal categories, which are called \emph{balanced tensor categories} in \cite{joyal-street}, when $\G$ is the ribbon group operad $\R$ in Definition \ref{def:ribbon-group-operad};
\item coboundary monoidal categories in the sense of Drinfel'd \cite{drinfeld} when $\G$ is the cactus group operad $\Cac$ in Definition \ref{def:cactus-group-operad}.
\end{enumerate}
Therefore, the concept of a $\G$-monoidal category provides a unifying framework for proving results about small monoidal categories with extra structures.

In Section \ref{sec:gmonoidal-cat} we define a $1$-colored symmetric operad $\MCatg$ in $\Cat$, called the \emph{$\G$-monoidal category operad}, whose algebras are defined as $\G$-monoidal categories.  The $\G$-monoidal category operad is an extension of the monoidal category operad $\MCat$ in Definition \ref{def:mcat-operad} by allowing the elements in $\G$ to act on the objects, which are standard non-associative monomials.  In fact, restricting the definition to the case $\G=\P$ yields $\MCatp$, which coincides with the monoidal category operad $\MCat$.  Therefore, a $\P$-monoidal category structure is the same as a monoidal category structure.

In Section \ref{sec:coherence-gmonoidal-cat} we observe that a $\G$-monoidal category structure on a small category is exactly a monoidal category structure together with a suitable $\G$-action on iterated monoidal products.  We regard this as a coherence result for $\G$-monoidal categories, which are originally defined as algebras over the $\G$-monoidal category operad $\MCatg$.  We emphasize that a typical $\G$-monoidal category has a general underlying monoidal category, whose associativity isomorphism, left unit isomorphism, and right unit isomorphism are not required to be identity morphisms.

A $\G$-monoidal category is \emph{strict} if its underlying monoidal category is strict, i.e., if the associativity and left/right unit isomorphisms are identity morphisms.  In Section \ref{sec:strict-gmonoidal-cat} we construct the $1$-colored symmetric operad $\MCatgst$ in $\Cat$ whose algebras are exactly strict $\G$-monoidal categories.  In \cite{gurski} Gurski defined a categorical Borel construction corresponding to $\G$, whose algebras are strict $\G$-monoidal categories.  Therefore, the $\G$-monoidal category operad $\MCatg$ may be regarded as a generalization of Gurski's categorical Borel construction by allowing non-strict monoidal categories.

In Section \ref{sec:gmonoidal-functor} we define a (strong/strict) $\G$-monoidal functor between $\G$-monoidal categories as a (strong/strict) monoidal functor between the underlying monoidal categories that is compatible with the $\G$-monoidal structure isomorphisms.  For the planar operad $\P$, a (strong/strict) $\P$-monoidal functor is exactly a (strong/strict) monoidal functor.  We observe that an $\MCatg$-algebra morphism is a strict $\G$-monoidal functor.

In Section \ref{sec:operad-gmonoidal-functor} we describe the $2$-colored symmetric operad $\MFung$ whose algebras are general $\G$-monoidal functors between $\G$-monoidal categories.  This is the $\G$-analogue of the monoidal functor operad $\MFun$ in Definition \ref{def:monoidal-functor-operad}.  There is also a strong analogue $\MFungstg$ whose algebras are strong $\G$-monoidal functors between $\G$-monoidal categories.  This is the $\G$-analogue of the strong monoidal functor operad $\MFunstg$ in Definition \ref{def:strong-monoidal-functor-operad}.

\section{$\G$-Monoidal Category Operad}\label{sec:gmonoidal-cat}

The purpose of this section is to define $\G$-monoidal categories as algebras over a suitable symmetric operad in $\Cat$.  First we define the categories that constitute this symmetric operad.  

For an action operad $\G$ with augmentation $\omega : \G \to \S$, recall that the underlying permutation of each element $g\in \G(n)$ is denoted by \[\gbar = \omega(g)\in S_n.\]  For a standard non-associative monomial $p$ of weight $n \geq 0$ as in Definition \ref{def:nonass-monomial}, recall that $\pbar$ denotes its underlying standard monomial of weight $n$.  It is obtained from $p$ by removing all instances of the variable $x_0$ and all pairs of parentheses.  Permutations in $S_n$ act on standard monomials by permuting the variables, so 
\[\begin{split}
\tau\pbar\bigl(x_1,\ldots,x_n\bigr) &= \pbar\bigl(x_{\tau^{-1}(1)},\ldots, x_{\tau^{-1}(n)}\bigr),\\
\pbar\bigl(x_1,\ldots,x_n\bigr)\tau &= \pbar\bigl(x_{\tau(1)},\ldots, x_{\tau(n)}\bigr)\end{split}\] 
for $\tau \in S_n$.  Each standard monomial $p$ has the form 
\begin{equation}\label{pxtaup}
p\bigl(x_1,\ldots,x_n\bigr) = (x_1 \cdots x_n)\tau_p
\end{equation}
for some unique permutation $\tau_p \in S_n$.

\begin{definition}\label{def:mcatg-underlying-categories}
For $n \geq 0$, define the groupoid $\MCatg(n)$ as follows.
\begin{description}
\item[Objects] The set of objects in $\MCatg(n)$ is the set of standard non-associative monomials in $\{x_0,\ldots,x_n\}$ of weight $n$.
\item[Morphisms] For $p,q\in \Ob\bigl(\MCatg(n)\bigr)$, we define the morphism set
\[\MCatg(n)(p,q) = \Bigl\{g \in \G(n) : \gbar = \tau_{\qbar} \tau_{\pbar}^{-1}\in S_n\Bigr\}\] with $\tau_{\pbar}$ and $\tau_{\qbar}$ as in \eqref{pxtaup}.
\item[Identity] The identity morphism of each object $p \in\Ob\bigl(\MCatg(n)\bigr)$ is the multiplicative unit $\id_n\in \G(n)$.
\item[Composition] The categorical composition in $\MCatg(n)$ is given by multiplication in the group $\G(n)$.
\end{description}
\end{definition}

\begin{example}[Planar Group Operad]\label{ex:mcatpn}
If $\G$ is the planar group operad $\P$ in Example \ref{ex:trivial-group-operad}, then $\MCatp(n)$ is the groupoid $\MCat(n)$ in Definition \ref{def:mcat-operad}, since $\P(n) = \{\id_n\}$ is the trivial group.\dqed
\end{example}

\begin{example}[Symmetric Group Operad]\label{ex:mcatsn}
If $\G$ is the symmetric group operad $\S$ in Example \ref{ex:symmetric-group-operad}, then for each pair of standard non-associative monomials $p,q$ of weight $n$, the morphism set $\MCats(n)(p,q)$ contains the unique permutation $\tau_{\qbar} \tau_{\pbar}^{-1}$.\dqed
\end{example}

\begin{example}[Braid Group Operad]\label{ex:mcatbn}
If $\G$ is the braid group operad $\B$ in Definition \ref{def:braid-group-operad}, then for each pair of standard non-associative monomials $p,q$ of weight $n$, the morphism set $\MCatb(n)(p,q)$ contains all the braids with underlying permutation $\tau_{\qbar} \tau_{\pbar}^{-1}$.  This morphism set can be non-canonically identified with the underlying set of the pure braid group $PB_n$, which is the kernel of the underlying permutation map $\pi : B_n \to S_n$.\dqed 
\end{example}

\begin{example}[Ribbon Group Operad]\label{ex:mcatrn}
If $\G$ is the ribbon group operad $\R$ in Definition \ref{def:ribbon-group-operad}, then for each pair of standard non-associative monomials $p,q$ of weight $n$, the morphism set $\MCatr(n)(p,q)$ contains all the ribbons with underlying permutation $\tau_{\qbar} \tau_{\pbar}^{-1}$.  This morphism set can be non-canonically identified with the underlying set of the pure ribbon group $PR_n$, which is the kernel of the underlying permutation map $\pi : R_n \to S_n$.\dqed 
\end{example}

\begin{example}[Cactus Group Operad]\label{ex:mcatcacn}
If $\G$ is the cactus group operad $\Cac$ in Definition \ref{def:cactus-group-operad}, then for each pair of standard non-associative monomials $p,q$ of weight $n$, the morphism set $\MCatcac(n)(p,q)$ contains all the cacti with underlying permutation $\tau_{\qbar} \tau_{\pbar}^{-1}$.  This morphism set can be non-canonically identified with the underlying set of the pure cactus group $PCac_n$, which is the kernel of the underlying permutation map $\pi : Cac_n \to S_n$.\dqed 
\end{example}

We now define a symmetric operad with the groupoids $\MCatg(n)$ as constituent categories.

\begin{definition}\label{def:mcatg-operad}\index{monoidal category operad!G-@$\G$-}\index{G-monoidal category operad@$\G$-monoidal category operad}\index{G-monoidal category@$\G$-monoidal category!operad}\index{operad!G-monoidal category@$\G$-monoidal category}
Define the \emph{$\G$-monoidal category operad} $\MCatg$ as the $1$-colored symmetric operad in $\Cat$ as follows.
\begin{description}
\item[Entries] For $n\geq 0$ the entry $\MCatg(n)$ is the groupoid in Definition \ref{def:mcatg-underlying-categories}.
\item[Unit] The operadic unit is the object $x_1 \in \MCatg(1)$.
\item[Equivariance] For each permutation $\sigma\in S_n$, the functor \[\nicexy{\MCatg(n)\ar[r]^-{\sigma} & \MCatg(n)}\] sends each standard non-associative monomial $p = p\bigl(x_{[0,n]}\bigr) \in \MCatg(n)$ to \[p\sigma = p\bigl(x_0,x_{\sigma(1)},\ldots,x_{\sigma(n)}\bigr)\] for $\sigma \in S_n$.  For a morphism $g\in \MCatg(n)(p,q)$, the image \[g\sigma \in \MCatg(n)\bigl(p\sigma, q\sigma\bigr)\] is the same element $g\in \G(n)$.
\item[Composition] For $1 \leq i \leq n$ and $m \geq 0$, the $\compi$-composition \[\nicexy{\MCatg(n) \times \MCatg(m) \ar[r]^-{\compi} & \MCatg(n+m-1)} \in \Cat\] on objects is given by \[p\bigl(x_{[0,n]}\bigr) \compi q\bigl(x_{[0,m]}\bigr) = p\Bigl(x_{[0,i-1]}, q\bigl(x_0,x_{[i,i+m-1]}\bigr), x_{[i+m,n+m-1]}\Bigr)\] $p \in \Ob\bigl(\MCatg(n)\bigr)$ and $q \in \Ob\bigl(\MCatg(m)\bigr)$.  On morphisms it sends a pair \[(g,h) \in \MCatg(n)(p,q)\times \MCatg(m)(r,s)\] to \[g\compi h \in \MCatg(n+m-1)\bigl(p\compi r, q\compi s\bigr),\] which is the $\compi$-composition in the action operad $\G$.  This is well-defined because the augmentation $\omega : \G \to \S$ preserves $\compi$, so \[\overline{g\compi h} = \gbar \compi \overline{h}.\]
\end{description}
This finishes the definition of the $\G$-monoidal category operad.
\end{definition}

\begin{interpretation}\label{int:mcatg-eq-welldef}
In the above definition, to see that the equivariant structure is well-defined on  morphisms, we use the notation in Definition \ref{def:mcatg-underlying-categories}.  Then we have the equalities
\[\begin{split}
\tau_{\overline{q\sigma}}\tau_{\overline{p\sigma}}^{-1} 
&= (\tau_{\qbar}\sigma)(\tau_{\pbar}\sigma)^{-1}\\
&= \tau_{\qbar}\tau_{\pbar}^{-1}.
\end{split}\]
So if $g \in \G(n)$ belongs to the morphism set $\MCatg(n)(p,q)$, then it also belongs to the morphism set $\MCatg(n)(p\sigma,q\sigma)$ for $\sigma \in S_n$.\dqed
\end{interpretation}

\begin{example}
For the planar group operad $\P$ in Example \ref{ex:trivial-group-operad}, the above definition of $\MCatp$ coincides with that of the monoidal category operad $\MCat$ in Definition \ref{def:mcat-operad}.\dqed
\end{example}

\begin{lemma}
For each action operad $\G$, $\MCatg$ is a $1$-colored symmetric operad in $\Cat$.
\end{lemma}

\begin{proof}
At the object set level, $\MCatg$ coincides with the monoidal category operad $\MCat$ in Definition \ref{def:mcat-operad}.  At the morphism level, the symmetric operad axioms follow from an inspection, using the fact that $\G$ is an action operad.
\end{proof}

\begin{definition}\label{def:g-monoidal-category}
For each action operad $\G$, a \emph{$\G$-monoidal category} is defined as an $\MCatg$-algebra.
\end{definition}

Since $\MCatg$ is a $1$-colored symmetric operad in $\Cat$, each $\G$-monoidal category has an  underlying small category.

\begin{corollary}\label{mcatp-algebra}
For each small category $A$, a $\P$-monoidal category structure and a monoidal category structure are the same thing.
\end{corollary}

\begin{proof}
This follows from Theorem \ref{mcat-algebra} and the fact that $\MCatp=\MCat$.
\end{proof}

To understand the structure of $\G$-monoidal categories, we first make the following observation about the $\G$-monoidal category operad with respect to changing the action operad.

\begin{proposition}\label{mcatg-change-g}\index{change of action operads!G-monoidal category operad@$\G$-monoidal category operad}
Suppose $\varphi : \Gone \to \Gtwo$ is a morphism of action operads in the sense of Definition \ref{def:augmented-group-operad-morphism}.  Then:
\begin{enumerate}
\item For each $n\geq 0$, there is an induced functor \[\nicexy{\MCatgone(n) \ar[r]^-{\varphi_*} & \MCatgtwo(n)}\] that is the identity function on objects and is induced by $\varphi : \Gone(n) \to \Gtwo(n)$ on morphisms.
\item The level-wise morphisms $\varphi_*$ form an induced morphism \[\nicexy{\MCatgone \ar[r]^-{\varphi_*} & \MCatgtwo}\] of $1$-colored symmetric operads in $\Cat$.
\end{enumerate}
\end{proposition}

\begin{proof}
The first assertion follows from the fact that $\varphi : \Gone(n) \to \Gtwo(n)$ is a group homomorphism that preserves underlying permutations.  The second assertion follows from the assumption that $\varphi$ is a morphism of $1$-colored planar operads in $\Set$ that preserves underlying permutations.
\end{proof}

\begin{proposition}\label{gmonoidal-cat-underlying-cat}
For each action operad $\G$, pulling back along the initial morphism $\iota : \P \to\G$ gives each $\G$-monoidal category the structure of a monoidal category.
\end{proposition}

\begin{proof}
By Proposition \ref{mcatg-change-g}(2), there is an induced morphism \[\nicexy{\MCat = \MCatp \ar[r]^-{\iota^*} & \MCatg}\] of $1$-colored symmetric operads in $\Cat$.  By Theorem \ref{changing-goperad}(2), each $\MCatg$-algebra pulls back via $\iota^*$ to an $\MCat$-algebra, which is the same thing as a monoidal category structure by Theorem \ref{mcat-algebra}.
\end{proof}

We call the monoidal category structure in Proposition \ref{gmonoidal-cat-underlying-cat} the \emph{underlying monoidal category} of a $\G$-monoidal category.

\section{Coherence for $\G$-Monoidal Categories, I}
\label{sec:coherence-gmonoidal-cat}

For an action operad $\G$, recall that a $\G$-monoidal category is defined as an $\MCatg$-algebra for the $1$-colored symmetric operad $\MCatg$ in $\Cat$ in Definition \ref{def:mcatg-operad}.  The purpose of this section is to provide the first version of a coherence result for $\G$-monoidal categories.  It subsumes Theorem \ref{mcat-algebra}, which provides an operadic interpretation of the Coherence Theorem for monoidal categories.  

In the next result, we will use Notation \ref{not:subscript-interval} for sequences of objects.

\begin{theorem}\label{gmonoidal-coherence-1}\index{coherence!G-monoidal category@$\G$-monoidal category}\index{G-monoidal category@$\G$-monoidal category!coherence}
For each small category $A$, a $\G$-monoidal category structure on $A$ is equivalent to a pair \[\Bigl(\bigl(A,\otimes,\tensorunit,\alpha,\lambda,\rho\bigr), (-)_*\Bigr)\] in which:
\begin{enumerate}
\item $\bigl(A,\otimes,\tensorunit,\alpha,\lambda,\rho\bigr)$ is a monoidal category structure on $A$.
\item $(-)_*$ assigns to each isomorphism $\sigma \in \MCatg(n)(p,q)$, for $n\geq 0$ and $p,q\in \Ob\bigl(\MCatg(n)\bigr)$, an isomorphism 
\begin{equation}\label{sigmastar-p-to-q}
\nicexy{p\bigl(\tensorunit,a_1,\ldots,a_n\bigr) \ar[r]^-{\sigma_*}_-{\cong} & q\bigl(\tensorunit,a_1,\ldots,a_n\bigr)}\in A,
\end{equation}
in which both sides are interpreted as iterated monoidal products involving the objects $a_1,\ldots,a_n\in A$ and copies of $\tensorunit$ as in \eqref{theta-p-ua}.  
\end{enumerate}
This isomorphism is natural in the objects $a_1,\ldots,a_n$ and satisfy the following five conditions.
\begin{description}
\item[Monoidal Structure Isomorphisms] $\alpha$, $\lambda$, and $\rho$ are given by 
\begin{equation}\label{alpha-lambda-rho}
\begin{split}
\alpha &= \Bigl(\nicexy{(x_1x_2)x_3\ar[r]^-{\id_3} & x_1(x_2x_3) \in \MCatg(3)}\Bigr)_*,\\
\lambda &= \Bigl(\nicexy{x_0x_1 \ar[r]^-{\id_1} & x_1 \in \MCatg(1)}\Bigr)_*,\\
\rho &= \Bigl(\nicexy{x_1x_0 \ar[r]^-{\id_1} & x_1 \in \MCatg(1)}\Bigr)_*.
\end{split}
\end{equation}
\item[Categorical Identity] For each object $p$ in $\MCatg(n)$, the equality
\begin{equation}\label{gmon-coherence-unit}
(\operadunit_p)_* = \Id_{p(\tensorunit,a_1,\ldots,a_n)}
\end{equation}
holds in $A$, in which $\operadunit_p \in \MCatg(n)(p,p)$ is the identity morphism of $p$, i.e., $\id_n\in \G(n)$.
\item[Categorical Composition] For each pair of isomorphisms \[(\pi,\sigma) \in \MCatg(n)(q,r) \times \MCatg(n)(p,q),\] the diagram
\begin{equation}\label{gmon-coherence-cat-comp}
\nicexy{p\bigl(\tensorunit,a_1,\ldots,a_n\bigr) \ar[r]^-{\sigma_*} \ar@(d,l)[dr]_-{(\pi\sigma)_*} & q\bigl(\tensorunit,a_1,\ldots,a_n\bigr) \ar[d]^-{\pi_*}\\
& r\bigl(\tensorunit,a_1,\ldots,a_n\bigr)}
\end{equation}
in $A$ is commutative.
\item[Equivariance] For
\begin{itemize}
\item $\sigma \in \MCatg(n)(p,q)$ and
\item $\tau \in S_n$ with $\sigma\tau \in \MCatg(n)(p\tau,q\tau)$, 
\end{itemize}
the isomorphism
\[\nicexy{(p\tau)\bigl(\tensorunit, a_{\tau(1)},\ldots,a_{\tau(n)}\bigr) \ar[r]^-{(\sigma\tau)_*} & (q\tau)\bigl(\tensorunit, a_{\tau(1)},\ldots,a_{\tau(n)}\bigr)}\]
is equal to $\sigma_*$ in \eqref{sigmastar-p-to-q}.
\item[Operadic Composition] For each pair of isomorphisms \[(\sigma,\tau) \in \MCatg(n)(p,q)\times \MCatg(m)(r,s)\] and $1\leq i \leq n$, the diagram
\begin{equation}\label{gmon-coherence-op-comp}
\begin{small}\nicexy@C-3.8cm{(p\compi r)\bigl(\tensorunit, a_{[1,n+m-1]}\bigr) \ar@{=}[d] \ar[rr]^-{(\sigma\compi\tau)_*} && (q\compi s)\bigl(\tensorunit, a_{[1,n+m-1]}\bigr) \ar@{=}[d] \\
p\bigl(\tensorunit, a_{[1,i-1]}, r(\tensorunit, a_{[i,i+m-1]}), a_{[i+m,n+m-1]}\bigr) \ar@(d,dl)[dr]_-{p \compi \tau_*} &&
q\bigl(\tensorunit, a_{[1,i-1]}, s(\tensorunit, a_{[i,i+m-1]}), a_{[i+m,n+m-1]}\bigr)\\
& p\bigl(\tensorunit, a_{[1,i-1]}, s(\tensorunit, a_{[i,i+m-1]}), a_{[i+m,n+m-1]}\bigr) \ar@(dr,d)[ur]_-{\sigma_*}&}\end{small}
\end{equation}
in $A$ is commutative.  Here $p \compi \tau_*$ is the monoidal product of the isomorphism \[\nicexy{r\bigl(\tensorunit, a_i,\ldots,a_{i+m-1}\bigr) \ar[r]^-{\tau_*}_-{\cong} & s\bigl(\tensorunit, a_i,\ldots,a_{i+m-1}\bigr)}\] with the identity morphisms of $\tensorunit$ and of $a_j$ for $1\leq j \leq i-1$ and $i+m \leq j \leq n+m-1$. 
\end{description}
\end{theorem}

\begin{proof}
Suppose $(A,\theta)$ is an $\MCatg$-algebra.  The $\MCatg$-algebra structure morphism
\begin{equation}\label{mcatg-algebra-structure-map}
\nicexy{\MCatg(n)\times A^{\times n} \ar[r]^-{\theta} & A} \in \Cat,
\end{equation}
when applied to each isomorphism $\sigma \in \MCatg(n)(p,q)$, yields the natural isomorphism $\sigma_*$ in \eqref{sigmastar-p-to-q}.  By Proposition \ref{gmonoidal-cat-underlying-cat} the composite $\theta\iota^*$ defines an underlying monoidal category structure on $A$, where $\iota : \P \to \G$ is the initial morphism to the action operad $\G$.  The associativity isomorphism $\alpha$, the left unit isomorphism $\lambda$, and the right unit isomorphism $\rho$ all arise in the form $(-)_*$ as stated above.  We refer the reader back to the proof of Theorem \ref{mcat-algebra} for the explicit identification of an $\MCat$-algebra structure and a monoidal category structure on $A$.  

The equality \eqref{gmon-coherence-unit} and the commutativity of the diagram \eqref{gmon-coherence-cat-comp} hold because $\theta$ is a functor.  The equivariance condition holds by \eqref{operad-algebra-eq}, which is the equivariance axiom of $(A,\theta)$ as an algebra over the symmetric operad $\MCatg$.  Similarly, the diagram \eqref{gmon-coherence-op-comp} is the associativity diagram \eqref{planar-operad-algebra-associativity} for the $\MCatg$-algebra $(A,\theta)$  restricted to the $\compi$-composition in $\MCatg$.

Conversely, suppose $A$ is equipped with a monoidal category structure and natural isomorphisms $\sigma_*$ in \eqref{sigmastar-p-to-q} that satisfy the stated conditions.  We define the functor $\theta$ in \eqref{mcatg-algebra-structure-map} as follows.
\begin{itemize}
\item At the object set level, $\theta$ is defined as the $\MCat$-algebra structure on $A$, which by Theorem \ref{mcat-algebra} is the same as the monoidal category structure on $A$.
\item For each isomorphism $\sigma \in \MCatg(n)(p,q)$, we define $\theta(\sigma;-)$ as the natural isomorphism $\sigma_*$ in \eqref{sigmastar-p-to-q}.  The categorical identity axiom \eqref{gmon-coherence-unit} and the categorical composition axiom \eqref{gmon-coherence-cat-comp} ensure that $\theta$ is indeed a functor. 
\end{itemize}
It remains to observe that the morphism $\theta$ satisfies the associativity axiom \eqref{planar-operad-algebra-associativity}, the unity axiom \eqref{planar-operad-algebra-unity}, and the equivariance axiom \eqref{operad-algebra-eq} of an algebra over the symmetric operad $\MCatg$.  

The associativity axiom \eqref{planar-operad-algebra-associativity} holds by:
\begin{enumerate}[label=(\roman*)]
\item the assumed commutativity of the diagram \eqref{gmon-coherence-op-comp};
\item a finite induction using \eqref{compi-to-gamma}, which expresses the operadic composition $\gamma$ in terms of the $\compi$-compositions.
\end{enumerate}
The unity axiom \eqref{planar-operad-algebra-unity} is satisfied because the underlying monoidal category of $A$ as an $\MCat$-algebra satisfies the unity axiom.  Finally, the equivariance axiom \eqref{operad-algebra-eq} follows from the assumed equivariance condition in the statement of the Theorem.  
\end{proof}

We call $\sigma_*$ in \eqref{sigmastar-p-to-q} the \emph{$\G$-monoidal structure isomorphism}.

\begin{interpretation}
Theorem \ref{gmonoidal-coherence-1} says that a $\G$-monoidal category is a monoidal category equipped with a compatible action by the action operad $\G$ on iterated monoidal products.  In later chapters, we will observe that $\G$-monoidal categories are symmetric monoidal categories, braided monoidal categories, ribbon monoidal categories, or coboundary monoidal categories when $\G$ is the symmetric group operad $\S$, the braid group operad $\B$, the ribbon group operad $\R$, or the cactus group operad $\Cac$, respectively.\dqed
\end{interpretation}

\section{Strict $\G$-Monoidal Categories}\label{sec:strict-gmonoidal-cat}

This section provides a strict version of Theorem \ref{gmonoidal-coherence-1}.  First we define a strict version of the $\G$-monoidal category operad.

\begin{definition}\label{def:mcatgst-operad}\index{strict G-monoidal category@strict $\G$-monoidal category!operad}\index{operad!strict G-monoidal category@strict $\G$-monoidal category}
The \emph{strict $\G$-monoidal category operad} $\MCatgst$ is the $1$-colored symmetric operad in $\Cat$ obtained from the $\G$-monoidal category operad $\MCatg$ in Definition \ref{def:mcatg-operad} by forgetting all instances of the variable $x_0$ and all parentheses in the standard non-associative monomials.  A \emph{strict $\G$-monoidal category} is defined as an $\MCatgst$-algebra.
\end{definition}

\begin{interpretation}
The set of objects in the groupoid $\MCatgst(n)$ consists of all the standard monomials of weight $n$ in $\{x_1,\ldots,x_n\}$, with the symmetric group action given by permutation of the variables.  The $\compi$-composition is given by substituting a standard monomial for the $i$th variable and shifting the indices accordingly.  For standard monomials $p,q$ of weight $n$, the morphism set is \[\MCatgst(n)(p,q) = \Bigl\{g \in \G(n) : \gbar = \tau_q \tau_p^{-1}\in S_n \Bigr\}\] with $\tau_p$ and $\tau_q$ as in \eqref{pxtaup}.  The identity morphism of each object in $\MCatgst(n)$ is $\id_n\in \G(n)$.  The categorical composition in $\MCatgst(n)$ is given by multiplication in the group $\G(n)$.

In addition to being the strict analogue of the $\G$-monoidal category operad, the symmetric operad $\MCatgst$ is also the $\G$-analogue of the strict monoidal category operad $\MCatst$ in Definition \ref{def:mcatst}.  In fact, at the object set level, $\MCatgst$ agrees with the strict monoidal category operad $\MCatst$, whose algebras are strict monoidal categories.  Recall that a monoidal category is \emph{strict} if the associativity isomorphism, the left unit isomorphism, and the right unit isomorphism are all identity morphisms.\dqed
\end{interpretation}

\begin{theorem}\label{mcatgst-algebra}
For each small category $A$, the following two structures are equivalent:
\begin{enumerate}
\item An $\MCatgst$-algebra structure on $A$.
\item An $\MCatg$-algebra structure on $A$ whose underlying monoidal category is strict.
\end{enumerate}
\end{theorem}

\begin{proof}
We simply reuse the proof of Theorem \ref{gmonoidal-coherence-1} by forgetting all instances of $x_0$'s and all parentheses in standard non-associative monomials.  Instead of Theorem \ref{mcat-algebra}, here we use Theorem \ref{mcatst-algebra}, which says that a strict monoidal category structure is equivalent to an $\MCatst$-algebra structure on $A$.
\end{proof}

\begin{interpretation}
Theorem \ref{mcatgst-algebra} says that a strict $\G$-monoidal category is exactly a $\G$-monoidal category with an underlying strict monoidal category.  Along with Theorem \ref{gmonoidal-coherence-1}, this means that a strict $\G$-monoidal category is a strict monoidal category equipped with a compatible action by the action operad $\G$ on iterated monoidal products.\dqed
\end{interpretation}

\begin{remark}[Gurski's Categorical Borel Construction]\label{rk:gurski-borel}
In \cite{gurski} Gurski constructed a\index{categorical Borel construction} categorical Borel construction \[\Cat \ni A \mapsto \coprod_{n\geq 0} (E_{\G(n)})\timesover{\G(n)} A^{\times n},\] in which $E_{\G(n)}$ is the translation category of $\G(n)$.  Theorem \ref{mcatgst-algebra} above, which characterizes strict $\G$-monoidal categories as strict monoidal categories with extra structures, is close to Corollary 3.12 in \cite{gurski}.  In this context, we may regard the $\G$-monoidal category operad $\MCatg$ as a more general version of Gurski's categorical Borel construction that allows non-strict monoidal categories.  Below we will prove a coherence result that says that $\G$-monoidal categories can always be strictified to strict $\G$-monoidal categories, extending Mac Lane's Coherence Theorem for monoidal categories.\dqed
\end{remark}

For each standard non-associative monomial $p$, recall that $\pbar$ denotes the associated standard monomial, which is obtained from $p$ by forgetting all instances of the variable $x_0$ and all parentheses.

\begin{proposition}\label{mcatgst-change-g}
Suppose $\G$ is an action operad, and $\varphi : \Gone \to \Gtwo$ is a morphism of action operads.  Then:
\begin{enumerate}
\item For each $n\geq 0$, there is an induced functor\index{change of action operads!strict G-monoidal category operad@strict $\G$-monoidal category operad} \[\nicexy{\MCatgonest(n) \ar[r]^-{\varphi_*} & \MCatgtwost(n)}\] that is the identity function on objects and is induced by $\varphi : \Gone(n) \to \Gtwo(n)$ on morphisms.
\item The level-wise morphisms $\varphi_*$ form an induced morphism \[\nicexy{\MCatgonest \ar[r]^-{\varphi_*} & \MCatgtwost}\] of $1$-colored symmetric operads in $\Cat$.
\item There is a morphism\index{G-monoidal category operad@$\G$-monoidal category operad!to strict G-monoidal category operad@to strict $\G$-monoidal category operad} \[\nicexy{\MCatg \ar[r]^-{\pi} & \MCatgst}\] of $1$-colored symmetric operads in $\Cat$ that is given by the assignments
\[\begin{cases}
\Ob\bigl(\MCatg(n)\bigr)\ni p \mapsto \pbar \in \Ob\bigl(\MCatgst(n)\bigr),\\
\MCatg(n)(p,q) \ni \sigma \mapsto \sigma \in \MCatgst(n)(\pbar,\qbar).\end{cases}\] 
\item The diagram 
\[\nicexy{\MCatgone \ar[d]_-{\pi} \ar[r]^-{\varphi_*} & \MCatgtwo \ar[d]^-{\pi}\\
\MCatgonest \ar[r]^-{\varphi_*} & \MCatgtwost}\]
of $1$-colored symmetric operads in $\Cat$ is commutative.
\end{enumerate}
\end{proposition}

\begin{proof}
The first two assertions are the strict analogue of Proposition \ref{mcatg-change-g}, whose proof is reused here without change.  The last two assertions follow by a simple inspection.
\end{proof}

\begin{example}
The morphism $\pi : \MCatg \to \MCatgst$ of $1$-colored symmetric operads in $\Cat$ induces an adjunction \[\nicexy{\alg_{\Cat}\bigl(\MCatg\bigr) \ar@<2pt>[r]^-{\pi_!} & \alg_{\Cat}\bigl(\MCatgst\bigr) \ar@<2pt>[l]^-{\pi^*}}\] in which the right adjoint $\pi^*$ is the forgetful functor from strict $\G$-monoidal categories to $\G$-monoidal categories.  The left adjoint $\pi_!$ sends each $\G$-monoidal category to a strict $\G$-monoidal category.\dqed
\end{example}

\section{$\G$-Monoidal Functors}\label{sec:gmonoidal-functor}

In preparation for further discussion of coherence for $\G$-monoidal categories, in this section we discuss the $\G$-analogue of monoidal functors.  We will show that $\MCatg$-algebra morphisms and $\MCatgst$-algebra morphisms are strict $\G$-monoidal functors.  

Recall from Definition \ref{def:monoidal-functor} the notion of a monoidal functor between two monoidal categories.  Also recall from Theorem \ref{gmonoidal-coherence-1} that each $\G$-monoidal category has an underlying monoidal category and structure isomorphisms $\sigma_*$ as in \eqref{sigmastar-p-to-q}.  

\begin{notation}
For a monoidal functor $(F,F_2,F_0)$, a composite of morphisms, each of which is a monoidal product of $F_2$ and identity morphisms, is also denoted by $F_2$. 
\end{notation}

\begin{definition}\label{def:gmonoidal-functor}\index{G-monoidal functor@$\G$-monoidal functor}\index{functor!G-monoidal@$\G$-monoidal}
Suppose $A$ and $B$ are $\G$-monoidal categories.  A \emph{$\G$-monoidal functor} $F : A \to B$ is a monoidal functor $(F,F_2,F_0) : A \to B$ between the underlying monoidal categories such that the diagram
\begin{equation}\label{gmonoidal-functor-diagram}
\nicexy{p\bigl(\tensorunit^{B}, Fa_1,\ldots,Fa_n\bigr) \ar[d]_-{p(F_0,\Id_{Fa_1},\ldots,\Id_{Fa_n})} \ar[r]^-{\sigma_*}_-{\cong} & q\bigl(\tensorunit^{B}, Fa_1,\ldots,Fa_n\bigr) \ar[d]^-{q(F_0,\Id_{Fa_1},\ldots,\Id_{Fa_n})}\\
p\bigl(F\tensorunit^{A}, Fa_1,\ldots,Fa_n\bigr)  \ar[d]_-{F_2} & q\bigl(F\tensorunit^{A}, Fa_1,\ldots,Fa_n\bigr) \ar[d]^-{F_2}\\
F\Bigl[p\bigl(\tensorunit^{A}, a_1,\ldots, a_n\bigr)\Bigr] \ar[r]^-{F(\sigma_*)}_-{\cong} & F\Bigl[q\bigl(\tensorunit^{A}, a_1,\ldots, a_n\bigr)\Bigr]}
\end{equation}
is commutative for each morphism $\sigma \in \MCatg(n)(p,q)$ with $n \geq 0$ and objects $a_1,\ldots,a_n\in A$.  A $\G$-monoidal functor is said to be \emph{strong} (resp., \emph{strict}) if the underlying monoidal functor is strong (resp., strict).
\end{definition}

\begin{example}[$\P$-Monoidal Functors]
For the planar group operad $\P$ in Example \ref{ex:trivial-group-operad}, a $\P$-monoidal functor is exactly a monoidal functor.  Indeed, since each $\P(n) = \{\id_n\}$ is the trivial group, $\sigma$ must be $\id_n$.  So $p$ and $q$ have the same underlying standard monomial.  In other words, $p$ and $q$ only differ in their parenthesizations.  Both $(\id_n)_*$ in the diagram \eqref{gmonoidal-functor-diagram} are the unique natural isomorphisms guaranteed by the Coherence Theorem on monoidal categories \cite{maclane} (VII.2 Corollary).  So the commutativity of the diagram \eqref{gmonoidal-functor-diagram} follows from the monoidal functor axioms of $F$.  Similarly, a strong/strict $\P$-monoidal functor is exactly a strong/strict monoidal functor.\dqed
\end{example}

\begin{example}[Strict $\G$-Monoidal Functor]
For a strict $\G$-monoidal functor \[(F,\Id,\Id) : A \to B\] between $\G$-monoidal categories, the diagram \eqref{gmonoidal-functor-diagram} becomes the diagram
\begin{equation}\label{strict-gmon-functor}
\nicexy{p\bigl(\tensorunit^{B}, Fa_1,\ldots,Fa_n\bigr) \ar@{=}[d] \ar[r]^-{\sigma_*}_-{\cong} & q\bigl(\tensorunit^{B}, Fa_1,\ldots,Fa_n\bigr) \ar@{=}[d]\\
p\bigl(F\tensorunit^{A}, Fa_1,\ldots,Fa_n\bigr)  \ar@{=}[d] & q\bigl(F\tensorunit^{A}, Fa_1,\ldots,Fa_n\bigr) \ar@{=}[d]\\
F\Bigl[p\bigl(\tensorunit^{A}, a_1,\ldots, a_n\bigr)\Bigr] \ar[r]^-{F(\sigma_*)}_-{\cong} & F\Bigl[q\bigl(\tensorunit^{A}, a_1,\ldots, a_n\bigr)\Bigr].}
\end{equation}
In other words, $F$ sends the $\G$-monoidal structure isomorphism $\sigma_*$ in $A$ to the $\G$-monoidal structure isomorphism $\sigma_*$ in $B$.\dqed
\end{example}

The following observation is the $\G$-analogue of Theorem \ref{mcat-algebra-morphism}.

\begin{theorem}\label{mcatg-algebra-morphism}
Suppose $A$ and $B$ are $\G$-monoidal categories.  Then a morphism $F : A \to B$ of $\MCatg$-algebras is exactly a strict $\G$-monoidal functor.
\end{theorem}

\begin{proof}
The initial morphism $\iota : \P \to \G$ to the action operad $\G$ induces the forgetful functor \[\nicexy{\alg_{\Cat}\bigl(\MCatp\bigr) & \alg_{\Cat}\bigl(\MCatg\bigr) \ar[l]_-{\iota^*}}.\]  Since $\MCat=\MCatp$, each morphism \[F : (A,\theta^A) \to (B,\theta^B)\] of $\MCatg$-algebras has an underlying $\MCat$-algebra morphism, which is the same thing as a strict monoidal functor by Theorem \ref{mcat-algebra-morphism}.  The diagram \eqref{strict-gmon-functor} above is the commutative diagram \eqref{algebra-map-compatibility}
\begin{equation}\label{mcatg-algebra-map}
\nicexy@C+.4cm{\MCatg(n) \times A^{\times n} \ar[r]^-{(\Id, F^{\times n})} \ar[d]_-{\theta^A} & \MCatg(n) \times B^{\times n} \ar[d]^-{\theta^B}\\ A \ar[r]^-{F} & B}
\end{equation} 
applied to the morphism $\sigma \in \MCatg(n)(p,q)$.  So $F$ is a strict $\G$-monoidal functor.

Conversely, suppose $F : A \to B$ is a strict $\G$-monoidal functor.  In other words, $F$ is a strict monoidal functor between the underlying monoidal categories such that the diagram \eqref{strict-gmon-functor} is commutative for all morphisms $\sigma$ in $\MCatg$.  But then the diagram \eqref{mcatg-algebra-map} is commutative, which means that $F$ is a morphism of $\MCatg$-algebras.
\end{proof}

We close this section with the strict analogue of Proposition \ref{mcatg-algebra-morphism}, whose proof is reused here essentially without change.

\begin{proposition}\label{mcatgst-algebra-morphism}
Suppose $A$ and $B$ are strict $\G$-monoidal categories.  Then a morphism $F : A \to B$ of $\MCatgst$-algebras is exactly a strict $\G$-monoidal functor.
\end{proposition}

In summary, $\MCatg$ (resp., $\MCatgst$) is the $1$-colored symmetric operad in $\Cat$ whose algebras are (strict) $\G$-monoidal categories by definition and whose algebra morphisms are strict $\G$-monoidal functors.

\section{$\G$-Monoidal Functor Operad}\label{sec:operad-gmonoidal-functor}

Recall from Definition \ref{def:gmonoidal-functor} that, for an action operad $\G$, a $\G$-monoidal functor between two $\G$-monoidal categories is a monoidal functor between the underlying monoidal categories that is compatible with the $\G$-monoidal structures in the sense that the diagram \eqref{gmonoidal-functor-diagram} is commutative.  A $\G$-monoidal functor is \emph{strict} (resp., \emph{strong}) if the underlying monoidal functor is a strict (resp., strong) monoidal functor.  Proposition \ref{mcatg-algebra-morphism} says that $\MCatg$-algebra morphisms are strict $\G$-monoidal functors, not general $\G$-monoidal functors.

The purpose of this section is to describe a colored symmetric operad whose algebras are general $\G$-monoidal functors.  More precisely, we construct a $2$-colored symmetric operad $\MFung$ in $\Cat$ for which an algebra is a triple $(A,B,F)$ such that:
\begin{itemize}
\item $A$ and $B$ are $\G$-monoidal categories.
\item $F : A \to B$ is a $\G$-monoidal functor.
\end{itemize}
We will also discuss a variation for which an algebra is as above, but with $F$ a strong $\G$-monoidal functor.

The $2$-colored symmetric operad $\MFung$ is the $\G$-analogue of the monoidal functor operad $\MFun$ in Definition \ref{def:monoidal-functor-operad}, whose notations will be reused in this section.  Recall that an isomorphism $\sigma \in \MCatg(n)(p,q)$ is an element $g \in \G(n)$ with underlying permutation \[\gbar = \tau_{\qbar} \tau_{\pbar}^{-1}\in S_n,\] where $\tau_{\pbar}$ and $\tau_{\qbar}$ are as in \eqref{pxtaup}.

\begin{definition}\label{def:gmonoidal-functor-operad}\index{G-monoidal functor@$\G$-monoidal functor!operad}\index{operad!G-monoidal functor@$\G$-monoidal functor}
Define the \emph{$\G$-monoidal functor operad} $\MFung$ as the $\{s,t\}$-colored symmetric operad in $\Cat$ as follows.
\begin{description}
\item[Objects] The underlying object $\{s,t\}$-colored symmetric operad $\Ob(\MFung)$ in $\Set$ is equal to $\Ob(\MFun)$.  In other words, it is operadically generated by the elements:
\begin{enumerate}[label=(\roman*)]
\item $x_0 \in \Ob\left(\MFung\sbinom{s}{\varnothing}\right)$ and $\mu \in \Ob\left(\MFung\sbinom{s}{s,s}\right)$;
\item $x_0' \in \Ob\left(\MFung\sbinom{t}{\varnothing}\right)$ and $\mu' \in \Ob\left(\MFung\sbinom{t}{t,t}\right)$;
\item $f \in \Ob\left(\MFung\sbinom{t}{s}\right)$.
\end{enumerate}
The $s$-colored and the $t$-colored operadic units are denoted by, respectively, \[x_1\in \Ob\left(\MFung\sbinom{s}{s}\right) \andspace x_1'\in \Ob\left(\MFung\sbinom{t}{t}\right).\]  
\item[Morphisms] The morphism sets in the categories in $\MFung$ are categorically and operadically generated by:
\begin{enumerate}[label=(\roman*)]
\item an isomorphism \[\sigma_s \in \MFung\sbinom{s}{s^n}(p,q)\] for each isomorphism $\sigma \in \MCatg(n)(p,q)$;
\item an isomorphism \[\sigma_t \in \MFung\sbinom{t}{t^n}\bigl(p(x'_{[0,n]}), q(x'_{[0,n]})\bigr)\] for each isomorphism $\sigma \in \MCatg(n)(p,q)$;
\item a morphism \[f_2 \in \MFung\sbinom{t}{s,s}\bigl(f(x_1)f(x_2), f(x_1x_2)\bigr);\]
\item a morphism \[f_0 \in \MFung\sbinom{t}{\varnothing}\bigl(x_0',f(x_0)\bigr).\]
\end{enumerate}
\item[Relations] The above data is subject to the following four types of relations.
\begin{description}
\item[Source] The assignments \[\begin{cases}* \mapsto s \in \{s,t\},\\ \Ob\bigl(\MCatg(n)\bigr) \ni p \mapsto p \in \Ob\bigl(\MFung\sbinom{s}{s^n}\bigr),\\ \MCatg(n)(p,q) \ni \sigma \mapsto \sigma_s \in \MFun\sbinom{s}{s^n}(p,q)\end{cases}\] define a morphism 
\begin{equation}\label{iota-source-g}
\iota_s : \MCatg \to \MFung
\end{equation}
of symmetric operads in $\Cat$ in the sense of Definition \ref{def:g-operads}.
\item[Target] The assignments \[\begin{cases}* \mapsto t \in \{s,t\},\\ \Ob\bigl(\MCatg(n)\bigr) \ni p(x_{[0,n]}) \mapsto p(x'_{[0,n]}) \in \Ob\bigl(\MFung\sbinom{t}{t^n}\bigr), \\ \MCatg(n)(p,q) \ni \sigma \mapsto \sigma_t \in \MFung\sbinom{t}{t^n}\bigl(p(x'_{[0,n]}), q(x'_{[0,n]})\bigr)\end{cases}\] 
define a morphism 
\begin{equation}\label{iota-target-g}
\iota_t : \MCatg \to \MFung
\end{equation} 
of symmetric operads in $\Cat$.
\item[Monoidal Functor] The associativity relation \eqref{mfun-associativity}, the left unity relation \eqref{mfun-left-unity}, and the right unity relation \eqref{mfun-right-unity} are satisfied.
\item[$\G$-Monoidal Structure] The diagram
\begin{equation}\label{gmonoidal-functor-relation}
\nicexy{p\bigl(x_0', f(x_1),\ldots,f(x_n)\bigr) \ar[d]_-{p(f_0,f(x_1),\ldots,f(x_n))} \ar[r]^-{\sigma_t(f,\ldots,f)} & q\bigl(x_0', f(x_1),\ldots,f(x_n)\bigr) \ar[d]^-{q(f_0,f(x_1),\ldots,f(x_n))}\\
p\bigl(f(x_0),f(x_1),\ldots,f(x_n)\bigr) \ar[d]_-{\text{iteration of}\,f_2} & q\bigl(f(x_0),f(x_1),\ldots,f(x_n)\bigr) \ar[d]^-{\text{iteration of}\,f_2}\\
f\Bigl[p\bigl(x_0, x_1,\ldots, x_n\bigr)\Bigr] \ar[r]^-{f(\sigma_s)} & f\Bigl[q\bigl(x_0, x_1,\ldots, x_n\bigr)\Bigr]}
\end{equation}
in $\MFung\sbinom{t}{s^n}$ is commutative for each isomorphism $\sigma \in \MCatg(n)(p,q)$ with $n \geq 0$.
\end{description}
\end{description}
This completes the definition of the $\G$-monoidal functor operad $\MFung$.
\end{definition}

We now observe that $\MFung$ is the symmetric operad that precisely captures the concept of a general $\G$-monoidal functor, in which the structure morphisms are not necessarily invertible.

\begin{theorem}\label{mfung-algebra}
An $\MFung$-algebra is exactly a triple $(A,B,F)$ such that:
\begin{itemize}
\item $A$ and $B$ are $\G$-monoidal categories.
\item $F : A \to B$ is a $\G$-monoidal functor.
\end{itemize}  
\end{theorem}

\begin{proof}
This is a cosmetic modification of the proof of Theorem \ref{mfun-algebra}, which is the special case when $\G$ is the planar group operad $\P$.  For an $\MFung$-algebra $(A,\theta)$, the restrictions of the structure morphism $\theta$ via the morphisms $\iota_s$ in \eqref{iota-source-g} and $\iota_t$ in \eqref{iota-target-g} give the $\G$-monoidal categories $A_s$ and $A_t$. Proceeding as in the proof of Theorem \ref{mfun-algebra}, we obtain a monoidal functor \[\nicexy@C+.5cm{A_s \ar[r]^-{F\, =\, \theta(f,-)} & A_t}.\]  The $\G$-monoidal structure relation \eqref{gmonoidal-functor-relation}, when interpreted in the $\MFung$-algebra $(A,\theta)$, ensures that the monoidal functor $F$ renders the diagram \eqref{gmonoidal-functor-diagram} commutative.  So $F$ is a $\G$-monoidal functor.  

Conversely, given a $\G$-monoidal functor $F : A \to B$ between $\G$-monoidal categories, $F$ is in particular a monoidal functor between small monoidal categories.  We define an $\MFung$-algebra $\bigl(\{A,B\},\theta\bigr)$ exactly as in the proof of Theorem \ref{mfun-algebra}.  For an isomorphism $\sigma \in \MCatg(n)(p,q)$, we define the structure isomorphism
\[\nicexy{p\bigl(\tensorunit,a_1,\ldots,a_n\bigr) \ar[r]^-{\theta(\sigma)}_-{\cong} & q\bigl(\tensorunit,a_1,\ldots,a_n\bigr)}\in A_s\] as the $\G$-monoidal structure isomorphism $\sigma_*$ in \eqref{sigmastar-p-to-q} for $A_s$, and similarly for $A_t$.  The $\G$-monoidal structure relation \eqref{gmonoidal-functor-relation} is preserved precisely because $F$ makes the diagram \eqref{gmonoidal-functor-diagram} commutative.
\end{proof}

We close this section with the strong analogue of the $\G$-monoidal functor operad.  Recall that a $\G$-monoidal functor is said to be \emph{strong} if the underlying monoidal functor is strong, i.e., the structure morphisms $F_2$ and $F_0$ are isomorphisms.

\begin{definition}\label{def:strong-gmonoidal-functor-operad}\index{operad!strong G-monoidal functor@strong $\G$-monoidal functor}
The \emph{strong $\G$-monoidal functor operad} $\MFungstg$ is the $\{s,t\}$-colored symmetric operad in $\Cat$ with the same definition as the $\G$-monoidal functor operad $\MFung$, except for:
\begin{enumerate}[label=(\roman*)]
\item two more generating morphisms
\[\begin{split} f_2^{-1} &\in \MFungstg\sbinom{t}{s,s}\bigl(f(x_1x_2), f(x_1)f(x_2)\bigr),\\
f_0^{-1} &\in \MFungstg\sbinom{t}{\varnothing}\bigl(f(x_0), x_0'\bigr);\end{split}\]
\item four more relations \[\begin{cases} f_2^{-1}f_2 = \Id_{f(x_1)f(x_2)},\\ f_2f_2^{-1} = \Id_{f(x_1x_2)},\end{cases}\qquad \begin{cases} f_0^{-1}f_0 = \Id_{x_0'},\\ f_0f_0^{-1} = \Id_{f(x_0)}.\end{cases}\] 
\end{enumerate}
\end{definition}

The following result is the strong $\G$-monoidal functor analogue of Theorem \ref{mfung-algebra}, whose proof can be reused here with only minor changes.

\begin{theorem}\label{mfungstg-algebra}
An $\MFungstg$-algebra is exactly a triple $(A,B,F)$ such that:
\begin{itemize}
\item $A$ and $B$ are $\G$-monoidal categories.
\item $F : A \to B$ is a strong $\G$-monoidal functor.
\end{itemize}  
\end{theorem}

Finally, we note that there are morphisms relating various $\G$-monoidal functor operads.

\begin{proposition}\label{mfung-mfungstg-morphism}
Suppose $\G$ is an action operad, and $\varphi : \Gone\to\Gtwo$ is a morphism of action operads.
\begin{enumerate}
\item There is a morphism \[\nicexy{\MFung \ar[r]^-{\pi} & \MFungstg}\] of $\{s,t\}$-colored symmetric operads in $\Cat$ that is the identity assignment on generating objects and generating morphisms.
\item There is a commutative diagram \[\nicexy{\MFunof{\Gone} \ar[r]^-{\varphi_*} \ar[d]_-{\pi} & \MFunof{\Gtwo} \ar[d]^-{\pi}\\
\MFunstgof{\Gone} \ar[r]^-{\varphi_*} & \MFunstgof{\Gtwo}}\] 
of $\{s,t\}$-colored symmetric operads in $\Cat$, in which:
\begin{itemize}
\item Each $\varphi_*$ is the identity assignment on generating objects
\item The top $\varphi_*$ on generating morphisms is given by 
\[\begin{cases} \MFunof{\Gone}\sbinom{s}{s^n}(p,q) \ni \sigma_s \mapsto \varphi(\sigma_s) \in \MFunof{\Gtwo}\sbinom{s}{s^n}(p,q),\\
\MFunof{\Gone}\sbinom{t}{t^n}(p,q) \ni \sigma_t \mapsto \varphi(\sigma_t) \in \MFunof{\Gtwo}\sbinom{t}{t^n}(p,q),\\
\MFunof{\Gone}\sbinom{t}{s,s} \ni f_2 \mapsto f_2 \in \MFunof{\Gtwo}\sbinom{t}{s,s},\\
\MFunof{\Gone}\sbinom{t}{\varnothing} \ni f_0 \mapsto f_0 \in \MFunof{\Gtwo}\sbinom{t}{\varnothing}.\\
\end{cases}\]
\item The bottom $\varphi_*$ on generating morphisms is the same as the top $\varphi_*$ together with the assignment
\[\begin{cases} \MFunstgof{\Gone}\sbinom{t}{s,s} \ni f_2^{-1} \mapsto f_2^{-1} \in \MFunstgof{\Gtwo}\sbinom{t}{s,s},\\
\MFunstgof{\Gone}\sbinom{t}{\varnothing} \ni f_0^{-1} \mapsto f_0^{-1} \in \MFunstgof{\Gtwo}\sbinom{t}{\varnothing}.
\end{cases}\]
\end{itemize}
\end{enumerate}
\end{proposition}

\begin{proof}
Both assertions follow from a simple inspection.
\end{proof}

\chapter{Coherence for $\G$-Monoidal Categories}
\label{ch:coherence-gmonoidal}

The purpose of this chapter is to prove several forms of coherence for $\G$-monoidal categories for an action operad $\G$ as in Definition \ref{def:augmented-group-operad}.  As motivation, recall the following statements, which are all called the Coherence Theorem for monoidal categories in the literature.
\begin{enumerate}
\item Every monoidal category is adjoint equivalent to a strict monoidal category via strong monoidal functors \cite{maclane} (XI.3).
\item For each category, there is a strict monoidal equivalence from its free monoidal category to its free strict monoidal category \cite{joyal-street} (Theorem 1.2).
\item In the free monoidal category generated by a set of objects, every diagram is commutative \cite{maclane} (VII.2).
\item An explicit description of the free monoidal category generated by a category is also referred to as the Coherence Theorem for monoidal categories.
\end{enumerate}
We will prove the $\G$-analogues of these statements.  Together with the following chapters, this discussion of $\G$-monoidal categorical coherence subsumes various coherence statements in the literature for monoidal categories with extra structures.

In Section \ref{sec:coherence-gmonoidal-cat2} we observe that every $\G$-monoidal category is adjoint equivalent to a strict $\G$-monoidal category via strong $\G$-monoidal functors.  In Section \ref{sec:coherence-gmonoidal-cat3} we show that for each small category, there is a canonical strict $\G$-monoidal equivalence from its free $\G$-monoidal category to its free strict $\G$-monoidal category.  

In Section \ref{sec:coherence-gmonoidal-cat4} we first provide an explicit description of the free strict $\G$-monoidal category generated by a small category.  Using this description, next we prove that, in the free $\G$-monoidal category generated by a set of objects, a diagram is commutative if and only if composites with the same domain and the same codomain have the same image in the free strict $\G$-monoidal category generated by one object.

In Section \ref{sec:free-gmon-cat} we provide an explicit description of the free $\G$-monoidal category generated by a small category.  Using this description, along with its strict analogue, we provide alternative proofs of the coherence results in Section \ref{sec:coherence-gmonoidal-cat3} and Section \ref{sec:coherence-gmonoidal-cat4}.

\section{Strictification of $\G$-Monoidal Categories}
\label{sec:coherence-gmonoidal-cat2}

A version of the Coherence Theorem for monoidal categories, which we recalled in Theorem \ref{maclane-thm} above, says that each monoidal category is adjoint equivalent to a strict monoidal category via strong monoidal functors.  The purpose of this section is to prove the following $\G$-analogue of the Coherence Theorem for monoidal categories.

\begin{theorem}\label{gmonoidal-coherence2}\index{strictification!G-monoidal category@$\G$-monoidal category}\index{coherence!G-monoidal category@$\G$-monoidal category}\index{G-monoidal category@$\G$-monoidal category!strictification}
Suppose $A$ is a $\G$-monoidal category.  Then there exists an adjoint equivalence \[\nicexy{A \ar@<2pt>[r]^-{L} & A_{\st} \ar@<2pt>[l]^-{R}}\] in which:
\begin{enumerate}[label=(\roman*)]
\item $A_{\st}$ is a strict $\G$-monoidal category.
\item Both $L$ and $R$ are strong $\G$-monoidal functors.
\end{enumerate} 
\end{theorem}

\begin{proof}
We first apply Mac Lane's Coherence Theorem \ref{maclane-thm} to the underlying monoidal category of $A$ to obtain an adjoint equivalence \[\nicexy{A \ar@<2pt>[r]^-{L} & A_{\st} \ar@<2pt>[l]^-{R}}\] in which:
\begin{enumerate}[label=(\roman*)]
\item $A_{\st}$ is a strict monoidal category.
\item Both $(L,L_2,L_0)$ and $(R,R_2,R_0)$ are strong monoidal functors.
\end{enumerate}
Next we extend the strict monoidal category structure on $A_{\st}$ to a strict $\G$-monoidal category structure.  By Theorem \ref{mcatgst-algebra}, this is equivalent to extending $A_{\st}$ to an $\MCatg$-algebra.  Using the description of $\G$-monoidal categories in Theorem \ref{gmonoidal-coherence-1}, we need to:
\begin{enumerate}
\item Define the $\G$-monoidal structure isomorphism $\sigma_*$ in \eqref{sigmastar-p-to-q} for $A_{\st}$.
\item Check its naturality and the conditions \eqref{alpha-lambda-rho}--\eqref{gmon-coherence-op-comp}.
\end{enumerate}

Below we will write:
\begin{itemize}\item $\tensorunit$ and $\tensorunit^{\st}$ for the monoidal units in $A$ and in $A_{\st}$, respectively;
\item $\eta : LR \iso \Id_{A_{\st}}$ for the counit of the adjunction $L \dashv R$.
\end{itemize}
Suppose $\sigma \in \MCatg(n)(p,q)$ is an isomorphism for some $n \geq 0$, and $b_1,\ldots,b_n$ are objects in $A_{\st}$.  We define $\sigma_*$ in $A_{\st}$, which we denote by $\sigma_*^{\st}$, as the following composite of isomorphisms:
\begin{equation}\label{astrict-sigmastar}
\nicexy{p\bigl(\tensorunit^{\st}, b_1,\ldots,b_n\bigr) \ar[d]_-{p(L_0,\eta^{-1},\ldots,\eta^{-1})} \ar[r]^-{\sigma_*^{\st}} & q\bigl(\tensorunit^{\st}, b_1,\ldots,b_n\bigr)\\
p\bigl(L\tensorunit, LRb_1,\ldots,LRb_n\bigr) \ar[d]_-{L_2} & q\bigl(L\tensorunit,LR b_1,\ldots,LRb_n\bigr) \ar[u]_-{q(L_0^{-1},\eta,\ldots,\eta)}\\
L\Bigl[p\bigl(\tensorunit, Rb_1,\ldots,Rb_n\bigr)\Bigr] \ar[r]^-{L(\sigma_*)} & L\Bigl[q\bigl(\tensorunit, Rb_1,\ldots,Rb_n\bigr)\Bigr]. \ar[u]_-{L_2^{-1}}}
\end{equation}
The $\sigma_*$ in the bottom row of the previous diagram comes from the $\G$-monoidal category structure on $A$.  The naturality and the conditions \eqref{alpha-lambda-rho}--\eqref{gmon-coherence-op-comp} for $\sigma_*^{\st}$ follow from those for $\sigma_*$ in $A$ and the fact that $L$ is a strong monoidal functor.  Therefore, we have extended the strict monoidal category $A_{\st}$ to a strict $\G$-monoidal category.

It remains to check that the strong monoidal functors $L$ and $R$ are $\G$-monoidal functors, i.e., that the diagram \eqref{gmonoidal-functor-diagram} is commutative for both $L$ and $R$.  To save space in the diagram below, we will write:
\begin{itemize}
\item $\ua=(a_1,\ldots,a_n)\in A^{\times n}$, $L\ua=(La_1,\ldots,La_n) \in A_{\st}^{\times n}$, and similarly for $RL\ua$ and $LRL\ua$;
\item $f^n=(f,\ldots,f)$ with $n$ copies of $f\in \{\Id,\eta,\eta^{-1}, \epsilon, \epsilon^{-1}\}$.  
\end{itemize}
The diagram \eqref{gmonoidal-functor-diagram} for the left adjoint $L$ is the outermost diagram below.
\[\nicexy{p\bigl(\tensorunit^{\st},L\ua\bigr) \ar[dd]_-{p(L_0,\Id^n)} \ar[dr]|-{p(L_0,(\eta^{-1})^n)} \ar[rrr]^-{\sigma_*^{\st}} &&& q\bigl(\tensorunit^{\st},L\ua\bigr) \ar[dd]^-{q(L_0,\Id^n)}\\
& p\bigl(L\tensorunit, LRL\ua\bigr) \ar[dd]^-{L_2} & q\bigl(L\tensorunit, LRL\ua\bigr) \ar[ur]|-{q(L_0^{-1},\eta^n)} \ar[dr]|-{q(\Id,L(\epsilon^{-1})^n)}\\
p\bigl(L\tensorunit, L\ua\bigr) \ar[dd]_-{L_2} \ar[ur]|-{p(\Id,L\epsilon^n)} &&& q\bigl(L\tensorunit, L\ua\bigr) \ar[dd]^-{L_2}\\
& L\Bigl[p\bigl(\tensorunit, RL\ua\bigr)\Bigr] \ar[r]^-{L(\sigma_*)}  & L\Bigl[q\bigl(\tensorunit, RL\ua\bigr)\Bigr] \ar[uu]^-{L_2^{-1}} \ar[dr]|-{Lq(\Id,(\epsilon^{-1})^n)}&\\
L\Bigl[p\bigl(\tensorunit, \ua\bigr)\Bigr] \ar[ur]|-{Lp(\Id,\epsilon^n)} \ar[rrr]^-{L(\sigma_*)} &&& L\Bigl[q\bigl(\tensorunit,\ua\bigr)\Bigr]}\]
In the diagram above:
\begin{itemize}
\item The sub-diagram in the middle is commutative by the definition of $\sigma_*^{\st}$.  
\item The upper-left and upper-right triangles are commutative by the \index{triangle identity}triangle identity for the adjunction $L \dashv R$ in \cite{maclane} (IV Theorem 1(ii)), which says that \[\eta_L \circ L\epsilon = \Id_L.\]  
\item The left and right quadrilaterals are commutative by the naturality of $L_2$.
\item The bottom sub-diagram is commutative by the naturality of $\sigma_*$. 
\end{itemize}
This shows that $L$ is a strong $\G$-monoidal functor.  The proof that the right adjoint $R$ is a strong $\G$-monoidal functor is similar.
\end{proof}

\section{Strictification of Free $\G$-Monoidal Categories}
\label{sec:coherence-gmonoidal-cat3}

Another version of the Coherence Theorem for monoidal categories says that for each category $\C$, there is a strict monoidal equivalence from the free monoidal category generated by $\C$ to the free strict monoidal category generated by $\C$.  The purpose of this section is to prove a $\G$-analogue of this version of the Coherence Theorem for monoidal categories.

To prepare for this coherence result, recall from Definition \ref{def:mcatg-operad} that $\MCatg$ is the $1$-colored symmetric operad in $\Cat$ whose algebras are $\G$-monoidal categories.  There is a free-forgetful adjunction \[\nicexy@C+.8cm{\Cat \ar@<2pt>[r]^-{\MCatg \circs -} & \alg_{\Cat}\bigl(\MCatg\bigr) \ar@<2pt>[l]^-{U}},\] in which $\circs$ is the one-colored symmetric circle product in \eqref{symmetric-circle-product}.  So for each small category $A$ and each $\G$-monoidal category $B$, there is a bijection
\begin{equation}\label{mcatg-free-forget}
\Cat(A,UB) \cong \alg_{\Cat}(\MCatg)\Bigl(\MCatg(A),B\Bigr),
\end{equation}
in which \[\MCatg(A) = \MCatg \circs A.\]

The strict analogue $\MCatgst$ in Definition \ref{def:mcatgst-operad} is the $1$-colored symmetric operad in $\Cat$ whose algebras are strict $\G$-monoidal categories.  By Theorem \ref{mcatgst-algebra} each strict $\G$-monoidal category is in particular a $\G$-monoidal category.  For each small category $A$, the bijection \eqref{mcatg-free-forget} gives an $\MCatg$-algebra morphism 
\begin{equation}\label{lambdag-a}
\nicexy{\MCatg \circs A = \MCatg(A) \ar[r]^-{\Lambdag_A} & \MCatgst(A) = \MCatgst \circs A}
\end{equation}
that corresponds to the composite functor $\iotag_A$ in the diagram \[\nicexy@C+.5cm{A \ar[d]_-{\iotag_{A,1}} \ar[r]^-{\iotag_A} & \MCatgst(A)\\
\MCatgst(1) \times A \ar[r]^-{\text{include}} & \coprodover{n\geq 0}\, \MCatgst(n) \timesover{S_n} A^{\times n}. \ar@{=}[u]}\] The functor $\iotag_{A,1}$ is defined by the assignments \[\begin{cases}\Ob(A) \ni a \mapsto (x_1; a) \in \Ob\bigl(\MCatgst(1) \times A\bigr),\\
A(a_1,a_2) \ni f \mapsto (\Id_{x_1}; f) \in \MCatgst(1)\times A.\end{cases}\]  

Also note that  there is an equality\[\Lambdag_A = \pi \circs A\] with \[\pi : \MCatg \to \MCatgst\] the morphism of $1$-colored symmetric operads in Proposition \ref{mcatgst-change-g}(3).  By Theorem \ref{mcatg-algebra-morphism}, the $\MCatg$-algebra morphism $\Lambdag_A$ is a strict $\G$-monoidal functor from the free $\G$-monoidal category $\MCatg(A)$ generated by $A$ to the free strict $\G$-monoidal category $\MCatgst(A)$ generated by $A$.  We now observe that this is an equivalence of categories as well.

\begin{theorem}\label{gmonoidal-coherence3}\index{strictification!free G-monoidal category@free $\G$-monoidal category}\index{coherence!G-monoidal category@$\G$-monoidal category}
For each small category $A$, the strict $\G$-monoidal functor \[\nicexy{\MCatg(A) \ar[r]^-{\Lambdag_A} & \MCatgst(A)}\] in \eqref{lambdag-a} is an equivalence of categories.
\end{theorem}

\begin{proof}
The initial morphism $\iota : \P\to\G$ to the action operad $\G$ induces a commutative diagram
\begin{equation}\label{mcatp-mcatg-diagram}
\nicexy{\MCat(A) = \MCatp(A) \ar[r]^-{\iota_*} \ar[d]_-{\Lambdap_A} & \MCatg(A) \ar[d]^-{\Lambdag_A}\\  \MCatst(A) = \MCatpst(A) \ar[r]^-{\iota_*} & \MCatgst(A)}
\end{equation}
in $\Cat$ by Proposition \ref{mcatgst-change-g}(4).  For the planar group operad $\P$, the morphism $\Lambdap_A$ is the strict monoidal \emph{equivalence} from the free monoidal category $\MCat(A)$ generated by $A$ to the free strict monoidal category $\MCatst(A)$ generated by $A$ in  \cite{joyal-street} Theorem 1.2.  Since the top $\iota_*$ is the identity function on objects, to show that $\Lambdag_A$ is an equivalence of categories, it is enough to show that the diagram \eqref{mcatp-mcatg-diagram} is a pushout in $\Cat$.

To show that the diagram \eqref{mcatp-mcatg-diagram} is a pushout, consider a commutative diagram of solid-arrows\[\nicexy{\MCat(A) \ar[r]^-{\iota_*} \ar[d]_-{\Lambdap_A} & \MCatg(A) \ar[d]^-{\Lambdag_A} \ar@(r,u)[ddr]^-{f'} &\\  \MCatst(A) \ar[r]^-{\iota_*} \ar@(d,l)[drr]^-{f} & \MCatgst(A) \ar@{-->}[dr]^-{g} &\\ && B}\] in $\Cat$.  We must show that there exists a unique functor \[\nicexy{\MCatgst(A) \ar[r]^-{g} & B}\] that extends both $f$ and $f'$.  Since the bottom functor $\iota_*$ in the diagram \eqref{mcatp-mcatg-diagram} is the identity function on objects, we must define $g$ to be $f$ on objects.  

The functor $\Lambdag_A$ sends each morphism \[\varphi = \Bigl(\narrowxy{p\ar[r]^-{\sigma} & q}; h_1,\ldots,h_n\Bigr) \in \MCatg(n) \timesover{S_n} A^{\times n} \subseteq \MCatg(A)\] to the morphism \[\overline{\varphi} = \Bigl(\narrowxy{\pbar\ar[r]^-{\sigma} & \qbar}; h_1,\ldots,h_n\Bigr) \in \MCatgst(n) \timesover{S_n} A^{\times n} \subseteq \MCatgst(A).\]  Therefore, we must define \[g(\overline{\varphi}) = f'(\varphi).\] Since $\Lambdag_A$ is surjective on morphisms, this provides a candidate for the functor $g$ and also establishes its uniqueness.  To see that $g$ is well-defined on morphisms, suppose \[\varphi' = \Bigl(\narrowxy{p'\ar[r]^-{\sigma} & q'}; h_1,\ldots,h_n\Bigr) \in \MCatg(n) \timesover{S_n} A^{\times n} \subseteq \MCatg(A)\] is another lift of $\overline{\varphi}$ back to $\MCatg(A)$.  Since $\pbar' = \pbar$ and $\qbar'=\qbar$ as standard monomials of weight $n$, there is a commutative diagram \[\nicexy{p \ar[d]_-{\id_n} \ar[r]^-{\sigma} & q \ar[d]^-{\id_n}\\ p' \ar[r]^-{\sigma} & q'}\] in $\MCatg(n)$, where $\id_n$ is the multiplicative unit in $\G(n)$.  Together with the given condition \[f \Lambdap_A = f'\iota_*,\] this implies that \[f'(\varphi) = f'(\varphi') \in B,\] which means that $g$ is well-defined on morphisms as well.

The given condition, $f \Lambdap_A = f'\iota_*$, also implies that the above definition of $g$ defines a functor and satisfies the equalities \[g \iota_*=f \andspace g\Lambdag_A = f'.\]  This proves that the diagram \eqref{mcatp-mcatg-diagram} is a pushout in $\Cat$.
\end{proof}

\section{Free Strict $\G$-Monoidal Categories}
\label{sec:coherence-gmonoidal-cat4}

Another way to state the Coherence Theorem for monoidal categories is that in the free monoidal category generated by a set of objects, every diagram is commutative.  The purpose of this section is to prove a $\G$-analogue of this version of the Coherence Theorem for monoidal categories.

The strict $\G$-monoidal category operad $\MCatgst$ induces a free-forgetful adjunction \[\nicexy@C+.8cm{\Cat \ar@<2pt>[r]^-{\MCatgst(-)} & \alg_{\Cat}\bigl(\MCatgst\bigr) \ar@<2pt>[l]^-{U}},\] in which \[\MCatgst(-) = \MCatgst \circs -\] with $\circs$ the one-colored symmetric circle product in \eqref{symmetric-circle-product}.   Let us first describe explicitly this free strict $\G$-monoidal category functor.  For an element $\sigma \in \G(n)$, recall that $\sigmabar \in S_n$ denotes its underlying permutation.

\begin{theorem}\label{free-strict-gmon-category}\index{strict G-monoidal category@strict $\G$-monoidal category!free}
Consider the free strict $\G$-monoidal category $\MCatgst(A)$ generated by a small category $A$.  Then the following statements hold.
\begin{description}
\item[Objects] The objects in $\MCatgst(A)$ are the finite, possibly empty, sequences in $\Ob(A)$.
\item[Morphisms] For objects $\ua=(a_1,\ldots,a_n) \in A^{\times n}$ and $\ua'=(a_1',\ldots,a_m') \in A^{\times m}$, the morphism set is \[\MCatgst(A)(\ua;\ua') =\begin{cases} \varnothing & \text{if $m\not=n$};\\
\Bigl\{\bigl(\sigma; \{f_i\}_{i=1}^n\bigr) \in \G(n) \times \overset{n}{\prodover{i=1}} A\bigl(a_i,a'_{\sigmabar(i)}\bigr) \Bigr\} & \text{if $m=n$}.\end{cases}\]
\item[Identity] The identity morphism of an object $\ua = (a_1,\ldots,a_n)$ is \[\Bigl(\id_n\in \G(n); \{\Id_{a_i}\}_{i=1}^n\Bigr).\]
\item[Composition] Categorical composition is induced by that in $A$ and the group multiplication in each $\G(n)$.
\item[Strict Monoidal Structure] The monoidal product on objects is given by concatenation \[\ua \otimes \ua' = \bigl(a_1,\ldots,a_n,a_1',\ldots,a_m'\bigr),\] with the empty sequence as the monoidal unit.  The monoidal product of two morphisms is given by \[\bigl(\sigma; \{f_i\}_{i=1}^n\bigr) \otimes \bigr(\sigma'; \{f'_j\}_{j=1}^m\bigr) = \Bigl(\gammag\bigl(\id_2; \sigma,\sigma'\bigr); (f_1,\ldots,f_n,f'_1,\ldots,f'_m)\Bigr),\] in which \[\nicexy{\G(2)\times \G(n) \times \G(m) \ar[r]^-{\gammag} & \G(n+m)}\] is the operadic composition in $\G$.
\item[$\G$-Monoidal Structure] Suppose given:
\begin{itemize}
\item an isomorphism $\sigma \in \MCatgst(n)(p,q)$ with $p = (x_1\cdots x_n)\tau_p$ for some unique permutation $\tau_p\in S_n$;
\item objects $\ua^i = \bigl(a^i_1,\ldots,a^i_{k_i}\bigr) \in A^{\times k_i}$ for $1\leq i \leq n$.  
\end{itemize}
Then the structure isomorphism \eqref{sigmastar-p-to-q} \[\nicexy{p\big(\ua^1,\ldots,\ua^n\bigr) \ar[r]^-{\sigma_*}_-{\cong} & q\bigl(\ua^1,\ldots,\ua^n\bigr)} \in \MCatgst(A)\] is given by the pair \[\sigma_* = \Bigl(\gammag\bigl(\sigma;\id_{k_{\tau_p(1)}},\ldots,\id_{k_{\tau_p(n)}}\bigr); \{\Id_{a^{\tau_p(i)}_j}\}^{1\leq i \leq n}_{1 \leq j \leq k_{\tau_p(i)}}\Bigr),\] in which \[\nicexy{\G(n)\times \G(k_1)\times \cdots \times \G(k_n) \ar[r]^-{\gammag} & \G(k_1+\cdots+k_n)}\] is the operadic composition in $\G$.
\end{description}
\end{theorem}

\begin{proof}
The definition of the symmetric circle product gives the equality 
\begin{equation}\label{mcatgstofa}
\MCatgst(A) = \coprodover{n\geq 0}\, \MCatgst(n) \timesover{S_n} A^{\times n} \in \Cat.
\end{equation}  
The objects in $\MCatgst(n)$ are standard monomials in $\{x_i\}_{i=1}^n$ of weight $n$, so each object has the form \[\bigl(x_1\cdots x_n\bigr)\tau\] for a unique permutation $\tau \in S_n$.  Along with \eqref{mcatgstofa}, this yields the stated description for the objects in $\MCatgst(A)$.

For standard monomials $p$ and $q$ of weight $n$, there is the morphism set \[\MCatgst(n)(p,q) = \Bigl\{\sigma \in \G(n) : \sigmabar = \tau_q\tau_p^{-1}\Bigr\}.\]  Together with the previous paragraph, this yields the stated descriptions for the morphism sets, identity morphisms, and the categorical composition in $\MCatgst(A)$.

The description of the strict monoidal structure follows from the proof of Theorem \ref{mcat-algebra}, suitably adapted to $\MCatst$.  The formula for the $\G$-monoidal structure isomorphism $\sigma_*$ comes from the fact that $\MCatgst(A)$ is the free $\MCatgst$-algebra of $A$, so the $\MCatgst$-action is induced by the symmetric operad structure of $\MCatgst$ using the factors $\MCatgst(n)$ in \eqref{mcatgstofa}.
\end{proof}

\begin{example}[Free Strict Monoidal Categories]\label{ex:free-strict-monoidal-cat}\index{strict monoidal category!free}
Theorem \ref{free-strict-gmon-category} for $\G=\P$, the planar group operad in Example \ref{ex:trivial-group-operad}, gives an explicit description of the free strict monoidal category $\MCatst(A)$ generated by a small category $A$.  Its objects are finite sequences of objects in $A$ with morphism sets \[\MCatst(A)(\ua;\ua') = \begin{cases} \varnothing & \text{if $|\ua|\not=|\ua'|$},\\
\overset{|\ua|}{\prodover{i=1}} A(a_i,a_i') & \text{if $|\ua|=|\ua'|$}.\end{cases}\] Composition and identity morphisms are inherited from those in $A$.  The strict monoidal structure is given by concatenation.\dqed
\end{example}

\begin{example}[Free Strict $\G$-Monoidal Category of One Object]
\label{ex:free-strict-gmon-one-element}
Every set can be regarded as the object set of a discrete category, in which there are no non-identity morphisms.  Applying the functor $\MCatgst(-)$ to a set yields the free strict $\G$-monoidal category generated by that set of objects.  For example, the free strict $\G$-monoidal category $\MCatgst(*)$ generated by one element has object set \[\Ob\bigl(\MCatgst(*)\bigr) = \bigl\{0,1,2,\ldots\bigr\},\] i.e., the set of non-negative integers.  The morphism sets are \[\MCatgst(*)(m;n) = \begin{cases} \varnothing & \text{if $m\not=n$};\\ \G(n) & \text{if $m=n$}.\end{cases}\] Categorical composition is given by multiplication in the groups $\G(n)$, with the multiplicative unit $\id_n\in \G(n)$ as the identity morphism of the object $n$.  

The strict monoidal structure is given by:
\begin{itemize}
\item addition \[m \otimes n = m+n\] on objects, with $0$ as the monoidal unit;
\item the operadic composition in $\G$\[\sigma \otimes \sigma' = \gammag\bigl(\id_2; \sigma, \sigma')\] for morphisms $(\sigma,\sigma') \in \G(n)\times\G(m)$.
\end{itemize}
Suppose given an isomorphism $\sigma \in \MCatgst(n)(p,q)$ with $p = (x_1\cdots x_n)\tau_p$ for some unique permutation $\tau_p\in S_n$, and $k_1,\ldots,k_n \geq 0$ with $k=k_1+\cdots+k_n$.  Then the $\G$-monoidal structure isomorphism \[\nicexy{k=p\big(k_1,\ldots,k_n\bigr) \ar[r]^-{\sigma_*}_-{\cong} & q\bigl(k_1,\ldots,k_n\bigr) = k} \in \MCatgst(*)(k;k)\] is given by \[\sigma_* = \gammag\bigl(\sigma;\id_{k_{\tau_p(1)}},\ldots,\id_{k_{\tau_p(n)}}\bigr) \in \G(k),\] which is an operadic composition in the action operad $\G$.\dqed
\end{example}

\begin{corollary}\label{free-strict-gmon-set}
For each set $A$ of objects, the unique function $A \to *$ induces a faithful functor \[\nicexy{\MCatgst(A) \ar[r]^-{\pi} & \MCatgst(*)}.\]
\end{corollary}

\begin{proof}
By Theorem \ref{free-strict-gmon-category}, for each pair of objects $\ua,\ua'$ in $\MCatgst(A)$, the morphism set is \[\MCatgst(A)(\ua;\ua') = \begin{cases} \varnothing & \text{if $|\ua|\not=|\ua'|$},\\ \Bigl\{\sigma \in \G(n) : \sigmabar\ua=\ua'\Bigr\} & \text{if $|\ua|=|\ua'|$}.\end{cases}\]  The functor $\pi$ is faithful because it takes each morphism $\sigma \in \MCatgst(A)(\ua;\ua')$ to the morphism $\sigma\in \MCatgst(*)(|\ua|;|\ua|) = \G(|\ua|)$.
\end{proof}

The main observation in this section is the following coherence result for $\G$-monoidal categories.

\begin{theorem}\label{gmonoidal-coherence4}\index{coherence!G-monoidal category@$\G$-monoidal category}\index{G-monoidal category@$\G$-monoidal category!coherence}
For each set $A$ of objects, the composite functor \[\nicexy{\MCatg(A) \ar[r]^-{\Lambdag_A} & \MCatgst(A) \ar[r]^-{\pi} & \MCatgst(*)}\] is faithful.
\end{theorem}

\begin{proof}
The functor $\Lambdag_A$ is an equivalence of categories by Theorem \ref{gmonoidal-coherence3}, so it is faithful.  The functor $\pi$ is faithful by Corollary \ref{free-strict-gmon-set}.
\end{proof}

\begin{corollary}\label{gmonoidal-coherence-underlying}
For each set $A$ of objects, a diagram in $\MCatg(A)$ is commutative if and only if composites with the same domain and the same codomain have the same $\pi\Lambdag_A$-image in $\MCatgst(*)$.
\end{corollary}

\begin{example}[Free Monoidal Category of a Set of Objects]\label{ex:free-moncat-set}
For the planar group operad $\P$, $\MCatp(A) = \MCat(A)$ is the free monoidal category generated by $A$.  The free strict monoidal category generated by one object $\MCatst(*)$ has no non-identity morphisms.  For a set $A$, the faithfulness of the functor \[\nicexy{\MCat(A) \ar[r]^-{\pi\Lambdap_A} & \MCatst(*)}\] implies that each diagram in the free monoidal category $\MCat(A)$ generated by $A$ is commutative.  This is one formulation of the Coherence Theorem for monoidal categories; see, for example, \cite{maclane} (VII.2) or \cite{joyal-street} (Corollary 1.6).\dqed  
\end{example}

\section{Free $\G$-Monoidal Categories}\label{sec:free-gmon-cat}

The purpose of this section is to provide an explicit description of the free $\G$-monoidal category generated by a small category.  As a consequence, we obtain alternative proofs of Theorem \ref{gmonoidal-coherence3} and Theorem \ref{gmonoidal-coherence4}.  We will use the description of a $\G$-monoidal category in Theorem \ref{gmonoidal-coherence-1}.

\begin{theorem}\label{free-gmon-category}\index{coherence!G-monoidal category@$\G$-monoidal category}\index{G-monoidal category@$\G$-monoidal category!coherence}\index{G-monoidal category@$\G$-monoidal category!free}
Consider the free $\G$-monoidal category $\MCatg(A)$ generated by a small category $A$.  Then the following statements hold.
\begin{description}
\item[Objects] Each object in $\MCatg(A)$ is a standard non-associative monomial in the objects $\{\tensorunit, a_1,\ldots,a_n\}$, i.e., has the form \[p\bigl(\tensorunit, a_{[1,n]}\bigr) = p\bigl(\tensorunit,a_1,\ldots,a_n\bigr)\] with:
\begin{itemize}
\item $p(x_0,x_1,\ldots,x_n)$ a standard non-associative monomial of weight $n$ in which the variables $\{x_1,\ldots,x_n\}$ appear in this order from left to right;
\item each $a_i$ an object in $A$;
\item $\tensorunit$ a newly adjoined object not from $A$.
\end{itemize}
\item[Morphisms] The morphism set is 
\[\begin{split}
&\MCatg(A)\Bigl(p\bigl(\tensorunit,a_{[1,n]}\bigr); q\bigl(\tensorunit,a'_{[1,m]}\bigr)\Bigr)\\
&= \begin{cases} \varnothing & \text{if $m\not=n$};\\
\Bigl\{\bigl(\sigma; \{f_i\}_{i=1}^n\bigr) \in \G(n) \times \overset{n}{\prodover{i=1}} A\bigl(a_i,a'_{\sigmabar(i)}\bigr) \Bigr\} & \text{if $m=n$}.\end{cases}\end{split}\]
\item[Identity] The identity morphism of an object $p\bigl(\tensorunit,a_{[1,n]}\bigr)$ is the pair \[\Bigl(\id_n; \{\Id_{a_i}\}_{i=1}^n\Bigr)\in \G(n) \times \overset{n}{\prodover{i=1}} A\bigl(a_i,a_i\bigr).\]
\item[Composition] Categorical composition is induced by that in $A$ and the group multiplication in the $\G(n)$.
\item[Monoidal Unit] The monoidal unit is the object $\tensorunit$.  
\item[Monoidal Product] The monoidal product on objects is given by concatenation, \[p\bigl(\tensorunit,a_{[1,n]}\bigr) \otimes q\bigl(\tensorunit,a'_{[1,m]}\bigr) = p\bigl(\tensorunit,a_{[1,n]}\bigr) q\bigl(\tensorunit,a'_{[1,m]}\bigr).\]  The monoidal product of two morphisms is given by \[\bigl(\sigma; \{f_i\}_{i=1}^n\bigr) \otimes \bigr(\sigma'; \{f'_j\}_{j=1}^m\bigr) = \Bigl(\gammag\bigl(\id_2; \sigma,\sigma'\bigr); (f_1,\ldots,f_n,f'_1,\ldots,f'_m)\Bigr),\] in which \[\nicexy{\G(2)\times \G(n) \times \G(m) \ar[r]^-{\gammag} & \G(n+m)}\] is the operadic composition in $\G$.
\item[Associativity Isomorphisms]
The associativity isomorphism 
\[\begin{small}\nicexy@C-.4cm{\Bigl[p\bigl(\tensorunit,a_{[1,n]}\bigr) \otimes q\bigl(\tensorunit,a'_{[1,m]}\bigr)\Bigr] \otimes r\bigl(\tensorunit,a''_{[1,l]}\bigr) \ar[r]^-{\alpha}_-{\cong} & p\bigl(\tensorunit,a_{[1,n]}\bigr) \otimes \Bigl[q\bigl(\tensorunit,a'_{[1,m]}\bigr) \otimes r\bigl(\tensorunit,a''_{[1,l]}\bigr)\Bigr]}\end{small}\] 
is \[\Bigl(\id_{n+m+l}\in \G(n+m+l); \{\Id_{a_i}\}_{i=1}^n, \{\Id_{a'_j}\}_{j=1}^m, \{\Id_{a''_k}\}_{k=1}^l\Bigr).\]
\item[Unit Isomorphisms] The left unit isomorphism $\lambda$ and the right unit isomorphism $\rho$, \[\nicexy{\tensorunit \otimes p\bigl(\tensorunit,a_{[1,n]}\bigr) \ar[r]^-{\lambda}_-{\cong} & p\bigl(\tensorunit,a_{[1,n]}\bigr) & p\bigl(\tensorunit,a_{[1,n]}\bigr) \otimes \tensorunit \ar[l]_-{\rho}^-{\cong}},\] are given by 
\[\begin{split}
\lambda &= \Bigl(\id_n\in \G(n); \{\Id_{a_i}\}_{i=1}^n\Bigr),\\
\rho &= \Bigl(\id_n\in \G(n); \{\Id_{a_i}\}_{i=1}^n\Bigr).
\end{split}\]
\item[$\G$-Monoidal Structure] Suppose given:
\begin{itemize}
\item an isomorphism $\sigma \in \MCatg(n)(p,q)$ with $\pbar = (x_1\cdots x_n)\tau_{\pbar}$ for some unique permutation $\tau_{\pbar}\in S_n$;
\item objects $p_i\bigl(\tensorunit,\ua^i\bigr)$ in $\MCatg(A)$ for $1\leq i \leq n$ with $\ua^i = \bigl(a^i_1,\ldots,a^i_{k_i}\bigr) \in A^{\times k_i}$.  
\end{itemize} 
Then the $\G$-monoidal structure isomorphism \eqref{sigmastar-p-to-q} \[\nicexy{p\Big(\tensorunit,p_1\bigl(\tensorunit,\ua^1\bigr), \ldots, p_n\bigl(\tensorunit,\ua^n\bigr)\Bigr) \ar[r]^-{\sigma_*}_-{\cong} & q\Bigl(\tensorunit, p_1\bigl(\tensorunit,\ua^1\bigr), \ldots, p_n\bigl(\tensorunit,\ua^n\bigr)\Bigr)}\] in $\MCatg(A)$ is given by the pair \[\sigma_* = \Bigl(\gammag\bigl(\sigma;\id_{k_{\tau_p(1)}},\ldots,\id_{k_{\tau_p(n)}}\bigr); \{\Id_{a^{\tau_p(i)}_j}\}^{1\leq i \leq n}_{1 \leq j \leq k_{\tau_p(i)}}\Bigr),\] in which \[\nicexy{\G(n)\times \G(k_1)\times \cdots \times \G(k_n) \ar[r]^-{\gammag} & \G(k_1+\cdots+k_n)}\] is the operadic composition in $\G$.
\end{description}
\end{theorem}

\begin{proof}
The definition of the symmetric circle product gives the coproduct
\begin{equation}\label{mcatgofa}
\MCatg(A) = \coprodover{n\geq 0}\, \MCatg(n) \timesover{S_n} A^{\times n} \in \Cat.
\end{equation}  
The objects in $\MCatg(n)$ are standard non-associative monomials in $\{x_i\}_{i=0}^n$ of weight $n$.  Each standard non-associative monomial of weight $n$ can be written uniquely as $p\tau$, in which $\tau \in S_n$, and $p$ is a standard non-associative monomial of weight $n$ in which the variables $\{x_1,\ldots,x_n\}$ appear in this order from left to right.  Along with \eqref{mcatgofa}, this yields the stated description of the objects in $\MCatg(A)$.

For standard non-associative monomials $p$ and $q$ of weight $n$, there is the morphism set \[\MCatg(n)(p,q) = \Bigl\{\sigma \in \G(n) : \sigmabar = \tau_{\qbar}\tau_{\pbar}^{-1}\Bigr\},\] in which $\pbar$ and $\qbar$ are the standard monomials of weight $n$ associated to $p$ and $q$, respectively.  Together with the previous paragraph, this yields the stated descriptions of the morphism sets, identity morphisms, and the categorical composition in $\MCatg(A)$.

The description of the monoidal structure follows from the proof of Theorem \ref{mcat-algebra}.  The formula for the $\G$-monoidal structure isomorphism $\sigma_*$ comes from the fact that $\MCatg(A)$ is the free $\MCatg$-algebra of $A$, so the $\MCatg$-action is induced by the symmetric operad structure of $\MCatg$ using the factors $\MCatg(n)$ in \eqref{mcatgofa}.
\end{proof}

\begin{example}[Free Monoidal Categories]\label{ex:free-monoidal-cat}\index{monoidal category!free}
Theorem \ref{free-gmon-category} for the planar group operad $\P$ gives an explicit description of the free monoidal category $\MCat(A)$ generated by a small category $A$.  Its objects are standard non-associative monomials in the objects $\{\tensorunit,a_1,\ldots,a_n\}$ with $n\geq 0$ and each $a_i$ an object in $A$, in which the objects $\{a_1,\ldots,a_n\}$ appear in this order from left to right.  There is the morphism set \[\MCat(A)\Bigl(p\bigl(\tensorunit,a_{[1,n]}\bigr); q\bigl(\tensorunit,a'_{[1,m]}\bigr)\Bigr)= \begin{cases} \varnothing & \text{if $m\not=n$};\\
\overset{n}{\prodover{i=1}} A\bigl(a_i,a'_i\bigr) & \text{if $m=n$}.\end{cases}\]  Identity morphisms and composition are inherited from those in $A$.  The monoidal product is given by concatenation.  The associativity isomorphism, the left unit isomorphism, and the right unit isomorphism are all induced by identity morphisms in $A$.\dqed
\end{example}

\begin{example}[Free $\G$-Monoidal Category on One Object]\label{ex:free-gmon-one-object}
When applied to the discrete category with one object $*$, Theorem \ref{free-gmon-category} gives an explicit description of the free $\G$-monoidal category $\MCatg(*)$ generated by one object.  Its objects are standard non-associative monomials in $\{\tensorunit, *,\ldots,*\}$ with $n\geq 0$ copies of the object $*$.  The morphism set between two such objects is empty if they have different numbers of copies of $*$, and is $\G(n)$ if they both have $n$ copies of $*$.  Identity morphisms and composition are inherited from the groups $\G(n)$.  The monoidal product and the $\G$-monoidal structure are induced by the planar operad structure of $\G$.  The associativity isomorphism, the left unit isomorphism, and the right unit isomorphism are all induced by the units in the groups $\G(n)$.\dqed
\end{example}

Observe that the explicit descriptions of the free strict $\G$-monoidal category $\MCatgst(A)$ in Theorem \ref{free-strict-gmon-category} and of the free $\G$-monoidal category $\MCatg(A)$ in Theorem \ref{free-gmon-category} only depend on
\begin{enumerate}[label=(\roman*)]
\item Theorem \ref{mcat-algebra}, which is an operadic interpretation of Mac Lane's Coherence Theorem for monoidal categories, and 
\item Theorem \ref{gmonoidal-coherence-1}, which is the first form of coherence for $\G$-monoidal categories.
\end{enumerate}
We now use these explicit descriptions of $\MCatgst(A)$ and $\MCatg(A)$ to give alternative proofs of Theorem \ref{gmonoidal-coherence3} and Theorem \ref{gmonoidal-coherence4}.

\begin{theorem}[$=$ Theorem \ref{gmonoidal-coherence3}]\label{gmonoidal-coherence3-alt}
For each small category $A$, the strict $\G$-monoidal functor \[\nicexy{\MCatg(A) \ar[r]^-{\Lambdag_A} & \MCatgst(A)}\] in \eqref{lambdag-a} is an equivalence of categories.
\end{theorem}

\begin{proof}
We check directly that this functor is essentially surjective on objects, and is full and faithful on morphism sets.  There is an equality \[\Lambdag_A = \pi \circs A\] with \[\pi : \MCatg \to \MCatgst\] the morphism of $1$-colored symmetric operads in Proposition \ref{mcatgst-change-g}(3).  In terms of the descriptions of $\MCatgst(A)$ in Theorem \ref{free-strict-gmon-category} and of $\MCatg(A)$ in Theorem \ref{free-gmon-category}, $\pi \circs A$ sends an object $p\bigl(\tensorunit,\ua\bigr)$ in $\MCatg(A)$ to the object $\ua$ in $\MCatgst(A)$.  So this functor is surjective on object sets.

For each pair of objects $p\bigl(\tensorunit, \ua\bigr)$ and $q\bigl(\tensorunit, \ua'\bigr)$ in $\MCatg(A)$, if $|\ua|\not=|\ua'|$, then we have the empty morphism sets \[\MCatg(A)\Bigl(p\bigl(\tensorunit, \ua\bigr); q\bigl(\tensorunit, \ua'\bigr)\Bigr) = \varnothing = \MCatgst(A)\bigl(\ua;\ua'\bigr).\]  On the other hand, if $|\ua|=|\ua'|=n$, then both morphism sets are equal to \[\Bigl\{\bigl(\sigma; \{f_i\}_{i=1}^n\bigr) \in \G(n) \times \overset{n}{\prodover{i=1}} A\bigl(a_i,a'_{\sigmabar(i)}\bigr) \Bigr\}\] with $\Lambdag_A$ the identity assignment.  So the functor $\Lambdag_A$ is full and faithful.
\end{proof}

\begin{remark}
The benefit of the previous proof is that it does not use the fact from \cite{joyal-street} (Theorem 1.2) that the strict monoidal functor \[\nicexy{\MCat(A) \ar[r]^-{\Lambdap_A} & \MCatst(A)},\] from the free monoidal category generated by $A$ to the free strict monoidal category generated by $A$, is an equivalence of categories.  On the other hand, the original proof of Theorem \ref{gmonoidal-coherence3} is more conceptual and does not require an explicit computation of $\MCatg(A)$ and $\MCatgst(A)$.\dqed
\end{remark}

\begin{theorem}[$=$ Theorem \ref{gmonoidal-coherence4}]
\label{gmonoidal-coherence4-alt}
For each set $A$ of objects, the composite functor \[\nicexy{\MCatg(A) \ar[r]^-{\Lambdag_A} & \MCatgst(A) \ar[r]^-{\pi} & \MCatgst(*)}\] is faithful.
\end{theorem}

\begin{proof}
Since $A$ is a discrete category, for each pair of objects $p\bigl(\tensorunit, \ua\bigr)$ and $q\bigl(\tensorunit, \ua'\bigr)$ in $\MCatg(A)$ with $|\ua|=|\ua'|=n$, there is the morphism set \[\MCatg(A)\Bigl(p\bigl(\tensorunit,\ua\bigr); q\bigl(\tensorunit,\ua'\bigr)\Bigr) = \Bigl\{\sigma \in \G(n) : \sigmabar \ua= \ua'\Bigr\}.\]  Using the description of $\MCatgst(*)$ in Example \ref{ex:free-strict-gmon-one-element}, the functor $\pi\Lambdag_A$ sends such a morphism $\sigma$ to $\sigma \in \G(n)$.  So $\pi\Lambdag_A$ is injective on each morphism set.
\end{proof}

\chapter{Braided and Symmetric Monoidal Categories}
\label{ch:braided-symmetric-monoidal-cat}

Recall from Definition \ref{def:monoidal-category} that a monoidal category is a category equipped with a multiplication $\otimes$ and a unit object $\tensorunit$, as well as an associativity isomorphism, a left unit isomorphism, and a right unit isomorphism that satisfy the unity axioms and the pentagon axiom.  A braided monoidal category is a monoidal category equipped with a natural braiding \[\nicexy{X\otimes Y \ar[r]^-{\xi_{X,Y}}_-{\cong} & Y \otimes X}\] that satisfies some appropriate axioms.  Braided monoidal categories arise naturally in the representation theory of quantum groups and in knot invariants; see, for example, \cite{chari,drinfeld-quantum,kassel,savage,street}.  A symmetric monoidal category is a braided monoidal category whose braiding is an involution.  

In this chapter, we observe that the concept of a $\G$-monoidal category in Definition \ref{def:g-monoidal-category} restricts to (i) a braided monoidal category when $\G$ is the braid group operad $\B$ and (ii) a symmetric monoidal category when $\G$ is the symmetric group operad $\S$.  The analogous statements for $\G$-monoidal functors in Definition \ref{def:gmonoidal-functor} are also true.  As a result, the coherence results for $\G$-monoidal categories in Chapter \ref{ch:coherence-gmonoidal} restrict to coherence results for braided monoidal categories and for symmetric monoidal categories.

In Section \ref{sec:braided-monoidal-cat} we show that a $\B$-monoidal category is exactly a small braided monoidal category with general associativity isomorphism, left/right unit isomorphisms, and braiding.  The $\B$-monoidal category operad $\MCatb$, whose algebras are by definition $\B$-monoidal categories, is more general than Fresse's unitary parenthesized braid operad, whose algebras are small braided monoidal categories with a strict monoidal unit.  See Remark \ref{rk:fresse-pab} for a more detailed discussion of the relationship between the $\B$-monoidal category operad and Fresse's unitary parenthesized braid operad.  It is also observed that a strict $\B$-monoidal category is exactly a small strict braided monoidal category.

In Section \ref{sec:coherence-bmon-functor} we observe that a $\B$-monoidal functor is exactly a braided monoidal functor.  We note that our braided monoidal functors are \emph{lax} in general; i.e., their structure morphisms are not required to be isomorphisms.  Braided monoidal functors with invertible structure morphisms are said to be \emph{strong}; these are what Joyal and Street \cite{joyal-street} called \emph{braided tensor functors}.  It is also observed that a strong/strict $\B$-monoidal functor is exactly a strong/strict braided monoidal functor.

In Section \ref{sec:braided-monoidal-coherence} we recover known coherence results for braided monoidal categories by restricting the coherence results in Chapter \ref{ch:coherence-gmonoidal} to the braid group operad $\B$.  This is possible by the results in Section \ref{sec:braided-monoidal-cat} and Section \ref{sec:coherence-bmon-functor}.  The first coherence result says that every small braided monoidal category can be strictified via an adjoint equivalence of strong braided monoidal functors.  The free analogue of this coherence result says that the free braided monoidal category generated by a small category can be strictified to its free strict braided monoidal category via a strict braided monoidal equivalence.  Another coherence result says that, in the free braided monoidal category generated by a set of objects, a diagram is commutative if and only if composites with the same (co)domain have the same underlying braids.

In the remaining sections of this chapter, we obtain analogous results for symmetric monoidal categories.  In Section \ref{sec:symmetric-monoidal-cat} we observe that (strict) $\S$-monoidal categories, where $\S$ is the symmetric group operad, are exactly small (strict) symmetric monoidal categories.  In Section \ref{sec:symmetric-monoidal-functor} we show that (strong/strict) $\S$-monoidal functors are exactly (strong/strict) symmetric monoidal functors.  In Section \ref{sec:symmetric-monoidal-coherence} we recover known coherence results for symmetric monoidal categories by restricting the coherence results in Chapter \ref{ch:coherence-gmonoidal} to the symmetric group operad $\S$.

\section{Braided Monoidal Categories are $\B$-Monoidal Categories}
\label{sec:braided-monoidal-cat}

The purpose of this section is to observe that for the braid group operad $\B$ in Definition \ref{def:braid-group-operad}, $\B$-monoidal categories in the sense of Definition \ref{def:g-monoidal-category} (i.e., \label{not:mcatb-operad}$\MCatb$-algebras) are exactly small braided monoidal categories in the usual sense.  What we call a braided monoidal category is what Joyal and Street \cite{joyal-street} called a \emph{braided tensor category}.  

\begin{definition}\label{def:braided-monoidal-category}\index{braided monoidal category}\index{monoidal category!braided}
A \emph{braided monoidal category} is a pair $(\M,\xi)$ in which:
\begin{itemize}
\item $(\M,\otimes,\tensorunit,\alpha,\lambda,\rho)$ is a monoidal category as in Definition \ref{def:monoidal-category}.
\item $\xi$ is a natural isomorphism 
\begin{equation}\label{braiding-isomorphism}
\nicexy{X \otimes Y \ar[r]^-{\xi_{X,Y}}_-{\cong} & Y \otimes X}
\end{equation}
for objects $X,Y \in \M$, called the \index{braiding}\emph{braiding}.
\end{itemize}
This data is required to satisfy the following axioms.
\begin{description}
\item[Unit Axiom] The diagram 
\begin{equation}\label{braiding-unit}
\nicexy{X \otimes \tensorunit \ar[d]_-{\rho} \ar[r]^-{\xi_{X,\tensorunit}}
& \tensorunit \otimes X \ar[d]^-{\lambda}\\ X \ar@{=}[r] & X}
\end{equation}
is commutative for all objects $X\in \M$.
\item[Hexagon Axioms] The following two hexagon diagrams\index{Hexagon Axioms} are required to be commutative for objects $X,Y,Z \in \M$.
\begin{equation}\label{hexagon-b1}
\nicexy@C-40pt{& (Y \otimes X) \otimes Z\ar[rr]^-{\alpha} & \hspace{1.5cm}
& Y \otimes (X \otimes Z) \ar[dr]^-{\Id_Y \otimes \xi_{X,Z}}&\\ 
(X \otimes Y) \otimes Z \ar[ur]^-{\xi_{X,Y}\otimes \Id_Z} \ar[dr]_-{\alpha} &&&& Y \otimes (Z\otimes X)\\
& X \otimes (Y \otimes Z) \ar[rr]^-{\xi_{X,Y\otimes Z}} & \hspace{1.5cm}
& (Y \otimes Z) \otimes X \ar[ur]_-{\alpha} &}
\end{equation}
\begin{equation}\label{hexagon-b2}
\nicexy@C-40pt{& X \otimes (Z\otimes Y) \ar[rr]^-{\alpha^{-1}} & \hspace{1.5cm}
& (X \otimes Z) \otimes Y \ar[dr]^-{\xi_{X,Z} \otimes \Id_Y}&\\ 
X \otimes (Y \otimes Z) \ar[ur]^-{\Id_X \otimes \xi_{Y,Z}} \ar[dr]_-{\alpha^{-1}} &&&& (Z \otimes X) \otimes Y\\
& (X \otimes Y) \otimes Z \ar[rr]^-{\xi_{X\otimes Y, Z}} & \hspace{1.5cm}
& Z \otimes (X \otimes Y) \ar[ur]_-{\alpha^{-1}} &}
\end{equation}
\end{description}
A braided monoidal category is said to be\index{strict braided monoidal category} \emph{strict} if the underlying monoidal category is strict.
\end{definition}

\begin{interpretation}
One way to make sense of the two hexagon diagrams \eqref{hexagon-b1} and \eqref{hexagon-b2} is to interpret them as braids (i.e., elements in the braid group $B_3$), with the braiding $\xi$ interpreted as the generator $s_1$ in the braid group $B_2$.  The reader may wish to refer to Section \ref{sec:braid-groups} for a brief discussion of braids.  

For example, ignoring the associativity isomorphism $\alpha$ for the moment, the first hexagon diagram \eqref{hexagon-b1} corresponds to the braid
\begin{center}\begin{tikzpicture}[xscale=.7, yscale=.4]
\foreach \x in {0,2,4} {\coordinate (A\x) at (\x,-1.5); \coordinate (B\x) at (\x,5);}
\draw[thick] (A0) to[out=90, in=270] (B4);
\foreach \x in {2,4} {\pgfmathsetmacro\xminus{\x -2}
\draw[line width=7pt, white] (A\x) to [out=90, in=270] (B\xminus);
\draw[thick] (A\x) to [out=90, in=270] (B\xminus);}
\draw[->] (-.5,-1.5) -- (4.5,-1.5); \node at (0,-2) {\scriptsize{$X$}};
\node at (2,-2) {\scriptsize{$Y$}}; \node at (4,-2) {\scriptsize{$Z$}};
\end{tikzpicture}\end{center}
in $B_3$, in which the two strings labeled by $Y$ and $Z$ cross over the string labeled by $X$.  Similarly, the second hexagon diagram \eqref{hexagon-b2} corresponds to the braid
\begin{center}\begin{tikzpicture}[xscale=.7, yscale=.4]
\foreach \x in {0,2,4} {\coordinate (A\x) at (\x,-1.5); \coordinate (B\x) at (\x,5);}
\draw[thick] (A0) to[out=90, in=270] (B2); \draw[thick] (A2) to[out=90, in=270] (B4);
\draw[line width=7pt, white] (A4) to [out=90, in=270] (B0);
\draw[thick] (A4) to [out=90, in=270] (B0);
\draw[->] (-.5,-1.5) -- (4.5,-1.5); \node at (0,-2) {\scriptsize{$X$}};
\node at (2,-2) {\scriptsize{$Y$}}; \node at (4,-2) {\scriptsize{$Z$}};
\end{tikzpicture}\end{center}
in which the string labeled by $Z$ crosses over the two strings labeled by $Y$ and $X$.\dqed
\end{interpretation}

Next is the main result of this section, in which $\B$-monoidal categories (i.e., $\MCatb$-algebras) are identified with small braided monoidal categories.

\begin{theorem}\label{bmonoidal-braided-monoidal}\index{braided monoidal category!operadic interpretation}
For each small category, the following two structures are equivalent:
\begin{enumerate}[label=(\roman*)]
\item A $\B$-monoidal category structure for the braid group operad $\B$.
\item A braided monoidal category structure.
\end{enumerate}
\end{theorem}

\begin{proof}
Suppose $(A,\xi)$ is a small braided monoidal category with underlying monoidal category $(A,\otimes,\tensorunit,\alpha,\lambda,\rho)$ and braiding $\xi$.  To show that $(A,\xi)$ defines a $\B$-monoidal category, we use Theorem \ref{gmonoidal-coherence-1} with $\G=\B$.  In other words, for each isomorphism $\sigma \in \MCatb(n)(p,q)$, we need to 
\begin{enumerate}
\item define the $\B$-monoidal structure isomorphism
\[\nicexy{p\bigl(\tensorunit,\ua\bigr)= p\bigl(\tensorunit,a_1,\ldots,a_n\bigr) \ar[r]^-{\sigma_*}_-{\cong} & q\bigl(\tensorunit,a_1,\ldots,a_n\bigr) = q\bigl(\tensorunit,\ua\bigr)}\in A\] in \eqref{sigmastar-p-to-q} for objects $a_1,\ldots,a_n\in A$, and
\item check the five conditions there. 
\end{enumerate} 
In the rest of this proof, we call an isomorphism in $A$ \emph{canonical} if it is a categorical composite of iterated monoidal products of identity morphisms, $\alpha$, $\lambda$, $\rho$, $\xi$, and their inverses.  We will define the $\B$-monoidal structure isomorphism $\sigma_*$ as a canonical isomorphism.

First we claim that there exists at least one canonical isomorphism from the object $p\bigl(\tensorunit,\ua\bigr)$ to the object $q\bigl(\tensorunit,\ua\bigr)$.  To see this, note that the braid $\sigma \in B_n$ satisfies $\sigmabar= \tau_{\qbar}\tau_{\pbar}^{-1}$ as in Definition \ref{def:mcatg-underlying-categories}.  This implies the equality
\begin{equation}\label{sigmabar-pst-qst}
\sigmabar= \tau_{\overline{q_{\st}}}\tau_{\overline{p_{\st}}}^{-1} \in S_n,
\end{equation}
where $p_{\st}$ is the standard non-associative monomial obtained from $p$ by removing all instances of $x_0$'s and likewise for $q_{\st}$.  The braid $\sigma\in B_n$ can be written non-uniquely as a finite product 
\begin{equation}\label{sigma-decomp-in-Bn}
\sigma = s_{i_k}^{\epsilon_k} \cdots s_{i_1}^{\epsilon_1} \in B_n
\end{equation} 
with each $s_{i_j} \in B_n$ a generating braid, $1\leq i_j \leq n-1$, and each $\epsilon_j\in \{1,-1\}$.  

We define the $\B$-monoidal structure isomorphism $\sigma_*$ as the composite 
\begin{equation}\label{sigmastar-braided-monoidal}
\nicexy{p\bigl(\tensorunit,\ua\bigr) \ar[rrr]^-{\sigma_*}_-{\cong} \ar[d]_-{u_p} &&& q\bigl(\tensorunit,\ua\bigr)\\
p_{\st}(\ua) \ar[d]_{\alpha_1} &&& q_{\st}(\ua) \ar[u]_-{u_q^{-1}}\\
p_{\st}^1(\ua) \ar[d]_-{c_{i_1}^{\epsilon_1}} &&& p_{\st}^{2k}(\ua) \ar[u]_-{\alpha_{k+1}} \\
p_{\st}^2(\ua) \ar[r]^-{\alpha_2} & \cdots \ar[r]^-{c_{i_{k-1}}^{\epsilon_{k-1}}} & p_{\st}^{2k-2}(\ua) \ar[r]^-{\alpha_k} & p_{\st}^{2k-1}(\ua) \ar[u]_-{c_{i_k}^{\epsilon_k}}}
\end{equation}
of canonical isomorphisms in $A$ defined as follows.
\begin{itemize}
\item $u_p$ and $u_q$ are the unique canonical isomorphisms involving only $\Id$, $\lambda$, and $\rho$ that are guaranteed by the Coherence Theorem for monoidal categories \cite{maclane} (VII.2 Corollary).  Each of these two morphisms removes all copies of the monoidal unit $\tensorunit$ corresponding to the $x_0$'s.
\item $p_{\st}^1$ is obtained from $p_{\st}$ by rearranging the parentheses such that the $i_1$th and the $(i_1+1)$st variables, counting from the left, are within an innermost pair of parentheses.  
\item The isomorphism $\alpha_1$ is the unique canonical isomorphism involving only $\Id$ and $\alpha^{\pm 1}$ guaranteed by the Coherence Theorem for monoidal categories. 
\item The canonical isomorphism $c_{i_1}^{\epsilon_1}$ is the monoidal product
\begin{equation}\label{ci-epsilon}
c_{i_1}^{\epsilon_1} = \underbrace{\Id \otimes \cdots \otimes \Id}_{i_1-1} \otimes \xi^{\epsilon_1} \otimes \underbrace{\Id \otimes \cdots \otimes \Id}_{n-i_1-1}
\end{equation}
corresponding to the braid group generator $s_{i_1}^{\epsilon_1} \in B_n$ in \eqref{sigma-decomp-in-Bn} with:
\begin{itemize}
\item $i_1-1$ (resp., $n-i_1-1$) instances of the identity morphism to the left (resp., right) of $\xi^{\epsilon_1}$;
\item the same parenthesization as $p_{\st}^1(\ua)$;
\item $p_{\st}^2$ obtained from $p_{\st}^1$ by switching the $i_1$th and the $(i_1+1)$st variables, counting from the left.
\end{itemize}
\item The canonical isomorphisms \[\nicexy{p_{\st}^2(\ua) \ar[r]^-{\alpha_2} & p_{\st}^3(\ua) \ar[r]^-{c_{i_2}^{\epsilon_2}} & p_{\st}^4(\ua)}\] are defined as in the previous case $(c_{i_1}^{\epsilon_1}, \alpha_1)$, starting with $p_{\st}^2(\ua)$ and using the index $i_2$ instead of $i_1$.  Inductively, this defines $\alpha_j$ and $c_{i_j}^{\epsilon_j}$ for $1\leq j \leq k$.
\item Because of \eqref{sigmabar-pst-qst}, $p_{\st}^{2k}$ and $q_{\st}$ differ only by the arrangement of the parentheses.  The isomorphism $\alpha_{k+1}$ is once again the unique canonical isomorphism involving only $\Id$ and $\alpha^{\pm 1}$ guaranteed by the Coherence Theorem for monoidal categories. 
\end{itemize}
By construction, $\sigma_*$ is a canonical isomorphism, which is natural in the objects $a_1,\ldots,a_n$.  

Next we check that $\sigma_*$ is independent of the choices in its construction.  For a given decomposition \eqref{sigma-decomp-in-Bn} of $\sigma\in B_n$ in terms of the braid group generators, the canonical isomorphism $\sigma_*$ in \eqref{sigmastar-braided-monoidal} is unique by the Coherence Theorem for monoidal categories, the unit axiom for the braiding $\xi$, and the two hexagon diagrams \eqref{hexagon-b1} and \eqref{hexagon-b2}.  It remains to check that, for each braid relation in Definition \ref{def:braid-group}, the two sides induce the same canonical isomorphism.  For the same reason as two sentences before, for this part of the proof, we may even (i) assume that $A$ is a \emph{strict} monoidal category and (ii) ignore the monoidal unit.

For the first braid relation, suppose $1\leq i,j \leq n-1$ with $|i-j|\geq 2$.  Then both $s_is_j$ and $s_js_i$ induce the canonical isomorphism
\[\underbrace{\Id \otimes \cdots \otimes \Id}_{i-1} \otimes \xi \otimes \underbrace{\Id \otimes \cdots \otimes \Id}_{j-i-2} \otimes \xi \otimes \underbrace{\Id \otimes \cdots \otimes \Id}_{n-j-1}\]if $j-i\geq 2$, and similarly if $i-j\geq 2$.

The second braid relation says \[s_is_{i+1}s_i = s_{i+1}s_i s_{i+1}\] for $1 \leq i \leq n-2$.  In the induced canonical isomorphisms, we are only concerned about the three consecutive objects starting from the $i$th object, since the identity morphisms are used in all other entries.  To simplify the diagram below, we call the three relevant, consecutive objects $x$, $y$, and $z$ (with $x$ being the $i$th object), respectively, and omit displaying the other objects and their identity morphisms.  We must show that the outermost diagram in the diagram
\[\nicexy@C+.5cm@R+.4cm{x \otimes y \otimes z \ar[d]_-{\Id_x \otimes \xi_{y,z}} \ar[r]^-{\xi_{x,y} \otimes \Id_z} \ar[dr]^-{\xi_{x,y\otimes z}} & y \otimes x \otimes z \ar[d]^-{\Id_y \otimes \xi_{x,z}}\\
x \otimes z \otimes y \ar[d]_-{\xi_{x,z}\otimes \Id_y} \ar[dr]^-{\xi_{x,z\otimes y}} & y \otimes z \otimes x \ar[d]^-{\xi_{y,z}\otimes \Id_x}\\
z \otimes x \otimes y \ar[r]^-{\Id_{z}\otimes \xi_{x,y}} & z \otimes y \otimes x}\]
is commutative.  The two triangles are commutative by the hexagon diagram \eqref{hexagon-b1}.  The middle parallelogram is commutative by the naturality of the braiding $\xi$.  We have shown that the canonical isomorphism $\sigma_*$ is well-defined.

Next we check the five conditions in Theorem \ref{gmonoidal-coherence-1} with $\G=\B$.
\begin{enumerate}[label=(\roman*)]
\item The axiom \eqref{alpha-lambda-rho} says that the associativity isomorphism in $A$ is $(\id_3)_*$ and that the left/right unit isomorphism is $(\id_1)_*$.  This is true by the construction of $\sigma_*$.
\item Similarly, the categorical identity axiom \eqref{gmon-coherence-unit} holds by construction.
\item For the categorical composition axiom \eqref{gmon-coherence-cat-comp}, as above we may assume that $A$ is a strict monoidal category and ignore the monoidal unit.  With this reduction, the axiom \eqref{gmon-coherence-cat-comp} holds simply because, if $\pi\in B_n$ admits a decomposition as \[\pi = s_{m_l}^{\varepsilon_l} \cdots s_{m_1}^{\varepsilon_1} \in B_n\] with $1\leq m_1,\ldots,m_l \leq n-1$ and $\varepsilon_1,\ldots, \varepsilon_l\in\{1,-1\}$, then $\pi\sigma$ admits a decomposition as \[\pi\sigma = s_{m_l}^{\varepsilon_l} \cdots s_{m_1}^{\varepsilon_1} s_{i_k}^{\epsilon_k} \cdots s_{i_1}^{\epsilon_1}\in B_n.\]
\item The equivariance condition holds because, for $\sigma \in \MCatb(n)(p,q)$ and a permutation $\tau \in S_n$, the morphism $\sigma\tau \in \MCatb(n)(p\tau,q\tau)$ is given by the same braid $\sigma$ by definition.  So $\sigma_*$ and $(\sigma\tau)_*$ are the same morphism.
\item For the operadic composition axiom \eqref{gmon-coherence-op-comp}, we once again apply the same reduction as in case (iii).  The $\compi$-composition of braids \eqref{braided-compi} is given by the product \[\sigma \compi \tau = \sigma\langle \underbrace{1,\ldots,1}_{i-1},m,\underbrace{1,\ldots,1}_{n-i}\rangle \cdot \Bigl(\id_1^{\oplus i-1}\oplus \tau\oplus \id_1^{\oplus n-i}\Bigr) \in B_{n+m-1}\] of a block braid induced by $\sigma$ and a direct sum braid induced by $\tau$.  So the diagram \eqref{gmon-coherence-op-comp} is commutative by the categorical composition axiom \eqref{gmon-coherence-cat-comp}. 
\end{enumerate}
This shows that $\bigl(A,(-)_*\bigr)$ is a $\B$-monoidal category.

Conversely, suppose $A$ is an $\B$-monoidal category, i.e., an $\MCatb$-algebra.  By Theorem \ref{gmonoidal-coherence-1}, $A$ has an underlying monoidal category.  We define the braiding \eqref{braiding-isomorphism} as 
\begin{equation}\label{onetwostar-braided}
\nicexy@C+1cm{a_1 \otimes a_2 \ar[r]^-{\xi_{a_1,a_2}\,=\, (s_1)_*} & a_2 \otimes a_1} \in A
\end{equation} 
in which \[s_1\in \MCatb(2)\bigl(x_1x_2,x_2x_1\bigr)\] with $s_1$ the generating braid in $B_2$.  Now we check the axioms for a braided monoidal category in Definition \ref{def:braided-monoidal-category}.
\begin{enumerate}[label=(\roman*)]
\item The unit axiom \eqref{braiding-unit} is satisfied because $B_1$ is the trivial group, so $\MCatb(1)(x_1x_0,x_1)$ contains a unique morphism.  
\item The first hexagon diagram \eqref{hexagon-b1} is commutative by:
\begin{itemize}\item the equality
\[\bigl(\id_1 \oplus s_1\bigr) \cdot \id_3 \cdot \bigl(s_1\oplus \id_1\bigr) = 
\id_3 \cdot \bigl(s_1 \comp_2 \id_2\bigr) \cdot \id_3 \in B_3;\]
\item the axioms \eqref{alpha-lambda-rho}, \eqref{gmon-coherence-cat-comp}, and \eqref{gmon-coherence-op-comp} in Theorem \ref{gmonoidal-coherence-1}.  
\end{itemize}
\item The second hexagon diagram \eqref{hexagon-b2} follows similarly from the equality
\[\bigl(s_1\oplus \id_1\bigr) \cdot \id_3 \cdot \bigl(\id_1 \oplus s_1\bigr) = 
\id_3 \cdot \bigl(s_1 \comp_1 \id_2\bigr) \cdot \id_3 \in B_3\] and Theorem \ref{gmonoidal-coherence-1}.  
\end{enumerate}
Therefore, $(A,\xi)$ is a braided monoidal category.
\end{proof}

\begin{remark}[Coherence for Braided Monoidal Categories]
\label{rk:braided-monoidal-coherence}
The proof of Theorem \ref{bmonoidal-braided-monoidal} above is modelled after the proof of the Coherence Theorem for \emph{symmetric} monoidal categories in \cite{maclane} (XI.1 Theorem 1).  Note that the previous proof did \emph{not} use the Coherence Theorem for braided monoidal categories in \cite{joyal-street} Corollary 2.6.  It only uses the Coherence Theorem for monoidal categories.  This is important in Section \ref{sec:braided-monoidal-coherence} when we recover known coherence results for braided monoidal categories by applying Theorem \ref{bmonoidal-braided-monoidal} and restricting the coherence results in Chapter \ref{ch:coherence-gmonoidal} to the braid group operad $\B$.  In other words, there is no circularity in our claim that coherence results for $\G$-monoidal categories restrict to coherence results for braided monoidal categories.\dqed
\end{remark}

\begin{remark}[Fresse's Parenthesized Braid Operad]\label{rk:fresse-pab}
Fresse in \cite{fresse-gt} (Theorem 6.2.4) described a $1$-colored symmetric operad \textit{\textsf{PaB}} in $\Cat$, called the\index{parenthesized braid operad}\index{operad!parenthesized braid} \emph{parenthesized braid operad}, whose algebras are small braided monoidal categories without monoidal units.  The parenthesized braid operad can be obtained from the $\B$-monoidal category operad $\MCatb$ by forgetting
\begin{enumerate}[label=(\roman*)]
\item the $0$th level $\MCatb(0)$ and
\item all instances of $x_0$'s in $\MCatb(n)$ for $n\geq 1$.
\end{enumerate}
Fresse proved that the parenthesized braid operad is equivalent to the fundamental groupoid operad of the little $2$-disc operad.

A variation of \textit{\textsf{PaB}} is the $1$-colored symmetric operad\label{not:pabplus} \textit{\textsf{PaB}}$_+$ in $\Cat$, called the \emph{unitary parenthesized braid operad}.  Its algebras are small braided monoidal categories whose monoidal units are strict with respect to the monoidal product, the associativity isomorphism, the left and right unit isomorphisms, and the braiding, i.e., 
\[\begin{split}
\lambda_x &= \rho_x = \xi_{\tensorunit,x} = \xi_{x,\tensorunit} = \Id_x ,\\
\alpha_{\tensorunit,x,y} &= \alpha_{x,\tensorunit,y} = \alpha_{x,y,\tensorunit} = \Id_{x\otimes y}.
\end{split}\]
This symmetric operad is obtained from the parenthesized braid operad \textit{\textsf{PaB}} by adjoining an object at level $0$, i.e., \textit{\textsf{PaB}}$_+(0) = *$.  

There is a morphism $\pi : \MCatb \to$ \!\!\textit{\textsf{PaB}}$_+$ of $1$-colored symmetric operads in $\Cat$ that sends each standard non-associative monomial of weight $0$ to $*$ in \textit{\textsf{PaB}}$_+(0)$.  For $n\geq 1$, $\pi$ sends each standard non-associative monomial $p$ of weight $n$ to the strict non-associative monomial $p_{\st}$, as in Definition \ref{def:mcat-terminology}, associated to $p$.\dqed
\end{remark}

Next is the strict analogue of Theorem \ref{bmonoidal-braided-monoidal}.  We will use the strict $\G$-monoidal category operad $\MCatgst$ in Definition \ref{def:mcatgst-operad} when $\G$ is the braid group operad $\B$.

\begin{corollary}\label{strict-bmonoidal-braided-monoidal}\index{strict braided monoidal category!operadic interpretation}
For each small category, the following two structures are equivalent:
\begin{enumerate}[label=(\roman*)]
\item An $\MCatbst$-algebra structure.
\item A strict braided monoidal category structure.
\end{enumerate}
\end{corollary}

\begin{proof}
By Theorem \ref{mcatgst-algebra} applied to $\G=\B$, an $\MCatbst$-algebra structure on $A$ is equal to an $\MCatb$-algebra structure on $A$ whose underlying monoidal category is strict.  By Theorem \ref{bmonoidal-braided-monoidal}, an $\MCatb$-algebra structure is equal to a braided monoidal category structure.  Therefore, an $\MCatbst$-algebra structure is exactly a braided monoidal category structure whose underlying monoidal category is strict.
\end{proof}

\begin{remark}[Earlier and Related Incarnations]\label{rk:fresse-cob-operad}
The strict $\B$-monoidal category operad $\MCatbst$, which encodes small strict braided monoidal categories by Corollary \ref{strict-bmonoidal-braided-monoidal}, has appeared in the literature before.
\begin{enumerate}
\item $\MCatbst$ is equal to Fresse's unitary colored braid operad\label{not:cobplus} \textit{\textsf{CoB}}$_+$.  Corollary \ref{strict-bmonoidal-braided-monoidal} is essentially Theorem 6.2.6 in \cite{fresse-gt}.
\item In Wahl's thesis \cite{wahl} the strict $\B$-monoidal category operad  $\MCatbst$ is denoted by $\mathcal{C}^{\beta}$, and Corollary \ref{strict-bmonoidal-braided-monoidal} recovers \cite{wahl} Theorem 1.4.4.
\item In Remark \ref{rk:gurski-borel} we mentioned Gurski's \cite{gurski} categorical Borel construction \[\Cat \ni A \mapsto \coprod_{n\geq 0} (E_{\G(n)})\timesover{\G(n)} A^{\times n}\] for an action operad $\G$, in which $E_{\G(n)}$ is the translation category of $\G(n)$.  Although the categorical Borel construction is technically different from the strict $\B$-monoidal category operad, they encode the same structure.  Indeed, when $\G$ is the braid group operad $\B$, Gurski in \cite{gurski} (Example 4.2) observed that the algebras over his categorical Borel construction are small strict braided monoidal categories.\dqed
\end{enumerate}
\end{remark}

\section{Braided Monoidal Functors are $\B$-Monoidal Functors}
\label{sec:coherence-bmon-functor}

The purpose of this section is to observe that, for the braid group operad $\B$, a $\B$-monoidal functor in the sense of Definition \ref{def:gmonoidal-functor} is exactly a braided monoidal functor in the usual sense.  The strong and strict analogues are also true.  As a consequence, the $\B$-monoidal category operad $\MCatb$ is the symmetric operad in $\Cat$ that encodes braided monoidal categories and strict braided monoidal functors.  Moreover, the $\B$-monoidal functor operad $\MFunb$ is the $2$-colored symmetric operad in $\Cat$ that encodes general braided monoidal functors between small braided monoidal categories.

Let us first recall the definition of a braided monoidal functor.

\begin{definition}\label{def:braided-monoidal-functor}
Suppose $\M$ and $\N$ are braided monoidal categories.  A \index{functor!braided monoidal}\index{braided monoidal functor}\index{monoidal functor!braided}\emph{braided monoidal functor} $(F,F_2,F_0) : \M \to \N$
is a monoidal functor between the underlying monoidal categories that is compatible with the braidings, in the sense that the diagram
\begin{equation}\label{monoidal-functor-braided}
\nicexy@C+.5cm{F(X) \otimes F(Y) \ar[d]_-{F_2} \ar[r]^-{\xi_{FX,FY}}_-{\cong} & F(Y) \otimes F(X) \ar[d]^-{F_2} \\ F(X \otimes Y) \ar[r]^-{F\xi_{X,Y}}_-{\cong} & F(Y \otimes X)}
\end{equation}
is commutative for all objects $X,Y \in \M$.  A braided monoidal functor is said to be \emph{strong} (resp., \emph{strict}) if the underlying monoidal functor is so.
\end{definition}

\begin{remark}A strong braided monoidal functor--i.e., a braided monoidal functor in which $F_2$ and $F_0$ are isomorphisms--is what Joyal and Street \cite{joyal-street} called a\index{braided tensor functor} \emph{braided tensor functor}.\dqed
\end{remark}

The main result of this section identifies braided monoidal functors with $\B$-monoidal functors for the braid group operad $\B$.

\begin{theorem}\label{bmonoidal-functor}\index{braided monoidal functor!operadic interpretation}
Suppose $F : \M \to \N$ is a functor between small braided monoidal categories.  Then the following two structures on $F$ are equivalent:
\begin{enumerate}[label=(\roman*)]
\item A braided monoidal functor structure.
\item A $\B$-monoidal functor structure.
\end{enumerate}
\end{theorem}

\begin{proof}
First note that $\M$ and $\N$, being small braided monoidal categories, are $\B$-monoidal categories by Theorem \ref{bmonoidal-braided-monoidal}.  The braidings in both $\M$ and $\N$ are given by $(s_1)_*$ for the isomorphism \[s_1 \in \MCatb(2)\bigl(x_1x_2,x_2x_1\bigr)\] with $s_1$ the generating braid in the braid group $B_2$.

Suppose $F$ is a $\B$-monoidal functor in the sense of Definition \ref{def:gmonoidal-functor}.  In other words, $(F,F_2,F_0) : \M\to\N$ is a monoidal functor such that the diagram \eqref{gmonoidal-functor-diagram} is commutative for each isomorphism $\sigma \in \MCatb(n)(p,q)$.  For the isomorphism $s_1 \in \MCatb(2)$ as in the previous paragraph, the commutative diagram \eqref{gmonoidal-functor-diagram} for a $\B$-monoidal functor becomes the commutative diagram \eqref{monoidal-functor-braided} for a braided monoidal functor.

Conversely, suppose $(F,F_2,F_0) : \M\to\N$ is a braided monoidal functor.  For each isomorphism $\sigma \in \MCatb(n)(p,q)$, we need to show that the diagram \eqref{gmonoidal-functor-diagram} for a $\B$-monoidal functor is commutative.  Since $F$ is a monoidal functor that is compatible with the braidings, the diagram \eqref{gmonoidal-functor-diagram} is commutative when $\sigma$ is the isomorphism:
\begin{itemize}
\item $\id_3 \in \MCatb(3)\bigl((x_1x_2)x_3,x_1(x_2x_3)\bigr)$, which yields the associativity isomorphisms in $\M$ and $\N$;
\item $\id_1 \in \MCatb(1)\bigl(x_0x_1,x_1\bigr)$, which yields the left unit isomorphisms in $\M$ and $\N$;
\item $\id_1 \in \MCatb(1)\bigl(x_1x_0,x_1\bigr)$, which yields the right unit isomorphisms in $\M$ and $\N$;
\item $s_1 \in \MCatb(2)\bigl(x_1x_2,x_2x_1\bigr)$, which yields the braidings in $\M$ and $\N$.
\end{itemize}

For a general isomorphism $\sigma \in \MCatb(n)(p,q)$, in each of $\M$ and $\N$, the $\B$-monoidal structure isomorphism \eqref{sigmastar-braided-monoidal} \[\nicexy{p\bigl(\tensorunit,a_1,\ldots,a_n\bigr) \ar[r]^-{\sigma_*}_-{\cong} & q\bigl(\tensorunit,a_1,\ldots,a_n\bigr)}\] is defined as the categorical composition of monoidal products of identity morphisms, the associativity isomorphism, the left unit isomorphism, the right unit isomorphism, the braiding, and their inverses.  Therefore, for each $\sigma$, the diagram \eqref{gmonoidal-functor-diagram} for a $\B$-monoidal functor factors into a ladder diagram, in which each square is commutative by the previous paragraph.
\end{proof}

\begin{interpretation}
Theorem \ref{bmonoidal-functor} asserts that for a monoidal functor $F : \M \to \N$ between small braided monoidal categories, in order for the diagram \eqref{gmonoidal-functor-diagram} to be commutative for each isomorphism in the $\B$-monoidal category operad $\MCatb$, it is necessary and sufficient that the single diagram \eqref{monoidal-functor-braided} be commutative.  Therefore, it is a coherence result for braided monoidal functors.\dqed
\end{interpretation}

Recall from Definition \ref{def:gmonoidal-functor} that a $\G$-monoidal functor is said to be \emph{strong} (resp., \emph{strict}) if the underlying monoidal functor is so.  The next observation is an immediate consequence of Theorem \ref{bmonoidal-functor}.

\begin{corollary}\label{strong-bmonoidal-functor}
Suppose $F : \M \to \N$ is a functor between small braided monoidal categories.  Then the following two structures on $F$ are equivalent:
\begin{enumerate}[label=(\roman*)]
\item A strong (resp., strict) braided monoidal functor structure.
\item A strong (resp., strict) $\B$-monoidal functor structure.
\end{enumerate}
\end{corollary}

We can now identify the algebras and the algebra morphisms of the $\B$-monoidal category operad $\MCatb$ with well-known concepts.

\begin{corollary}\label{mcatb-algebra}\index{braided monoidal category!operad}\index{operad!braided monoidal category}
$\MCatb$ is a $1$-colored symmetric operad in $\Cat$ whose:
\begin{itemize}
\item algebras are small braided monoidal categories;
\item algebra morphisms are strict braided monoidal functors.
\end{itemize}
\end{corollary}

\begin{proof}
This follows from Theorem \ref{bmonoidal-braided-monoidal}, the $\G=\B$ special case of Theorem \ref{mcatg-algebra-morphism}, and Corollary \ref{strong-bmonoidal-functor}.
\end{proof}

Next is the strict analogue.

\begin{corollary}\label{mcatbst-algebra}\index{strict braided monoidal category!operad}\index{operad!strict braided monoidal category}
$\MCatbst$ is a $1$-colored symmetric operad in $\Cat$ whose:
\begin{itemize}
\item algebras are small strict braided monoidal categories;
\item algebra morphisms are strict braided monoidal functors.
\end{itemize}
\end{corollary}

\begin{proof}
This follows from Corollary \ref{strict-bmonoidal-braided-monoidal}, the $\G=\B$ special case of Proposition \ref{mcatgst-algebra-morphism}, and Corollary \ref{strong-bmonoidal-functor}.
\end{proof}

For an action operad $\G$, recall from Definition \ref{def:gmonoidal-functor-operad} the $\G$-monoidal functor operad $\MFung$.

\begin{corollary}\label{mfunb-algebra}\index{braided monoidal functor!operad}\index{operad!braided monoidal functor}
An $\MFunb$-algebra is exactly a triple $(A,B,F)$ such that:
\begin{itemize}
\item $A$ and $B$ are small braided monoidal categories.
\item $F : A \to B$ is a braided monoidal functor.
\end{itemize}
\end{corollary}

\begin{proof}
This follows from Theorem \ref{mfung-algebra} applied to the braid group operad $\B$, Theorem \ref{bmonoidal-braided-monoidal}, and Theorem \ref{bmonoidal-functor}.
\end{proof}

Next is the strong analogue involving the strong $\G$-monoidal functor operad $\MFungstg$ in Definition \ref{def:strong-gmonoidal-functor-operad}.

\begin{corollary}\label{mfunbstg-algebra}\index{operad!strong braided monoidal functor}
An $\MFunbstg$-algebra is exactly a triple $(A,B,F)$ such that:
\begin{itemize}
\item $A$ and $B$ are small braided monoidal categories.
\item $F : A \to B$ is a strong braided monoidal functor.
\end{itemize}
\end{corollary}

\begin{proof}
This follows from Theorem \ref{mfungstg-algebra} applied to the braid group operad $\B$, Theorem \ref{bmonoidal-braided-monoidal}, and Corollary \ref{strong-bmonoidal-functor}.
\end{proof}

\begin{remark}[Idrissi's Parenthesized Permutations and Braids Operad]
Somewhat analogous to the strong $\B$-monoidal functor operad $\MFunbstg$ is  Idrissi's \cite{idrissi} \index{operad!parenthesized permutations and braids} \emph{operad of parenthesized permutations and braids} $\mathsf{PaPB}$.  It is a $2$-colored symmetric operad in $\Cat$ for which an algebra is a triple $(A,B,F)$ such that:
\begin{itemize}
\item $A$ is a small monoidal category without monoidal unit.
\item $B$ is a small braided monoidal category without monoidal unit.
\item $F$ is a strong braided monoidal functor from $B$ to the\index{Drinfel'd center} Drinfel'd center of $A$.
\end{itemize}
There is also a unitary extension $\mathsf{PaPB}_+$ for which an algebra is such a triple in which:
\begin{itemize}
\item $A$ and $B$ have monoidal units that are strict with respect to the monoidal products, the associativity isomorphisms, the left and right unit isomorphisms, and the braiding in $B$, as in Remark \ref{rk:fresse-pab}.
\item $F$ strictly preserves the monoidal units.
\end{itemize}
Idrissi proved that the operad of parenthesized permutations and braids is equivalent to the fundamental groupoid operad of Voronov's \index{Swiss-Cheese operad}\index{operad!Swiss-Cheese}Swiss-Cheese operad \cite{voronov}.  The Swiss-Cheese operad is closely related to \index{open-closed homotopy algebra}open-closed homotopy algebras \cite{hoefel,kajiura-stasheff}.\dqed
\end{remark}

\section{Coherence for Braided Monoidal Categories}
\label{sec:braided-monoidal-coherence}

The purpose of this section is to record several coherence statements for braided monoidal categories by restricting the coherence results for $\G$-monoidal categories in Chapter \ref{ch:coherence-gmonoidal} to the braid group operad $\B$ in Definition \ref{def:braid-group-operad}. 

In the previous two sections, we proved the following results.
\begin{enumerate}
\item In Theorem \ref{bmonoidal-braided-monoidal} we proved that $\B$-monoidal categories and small braided monoidal categories are the same thing.  As we pointed out in Remark \ref{rk:braided-monoidal-coherence}, the proof of Theorem \ref{bmonoidal-braided-monoidal} only uses the Coherence Theorem for monoidal categories and \emph{not} the Coherence Theorem for braided monoidal categories.   
\item The strict version is Corollary \ref{strict-bmonoidal-braided-monoidal}, which says that strict $\B$-monoidal categories and small strict braided monoidal categories are the same thing.  
\item Moreover, in Theorem \ref{bmonoidal-functor} we established the equivalence between $\B$-monoidal functors and braided monoidal functors, with the strong/strict versions in Corollary \ref{strong-bmonoidal-functor}. 
\end{enumerate}

Therefore, the coherence results in Chapter \ref{ch:coherence-gmonoidal} applies to braided monoidal categories and braided monoidal functors when $\G$ is the braid group operad $\B$.  First is the $\G=\B$ special case of Theorem \ref{gmonoidal-coherence2}, which says that braided monoidal categories can always be strictified.

\begin{theorem}\label{braided-monoidal-strictification}\index{braided monoidal category!coherence}\index{strictification!braided monoidal category}\index{coherence!braided monoidal category}
Suppose $A$ is a small braided monoidal category.  Then there exists an adjoint equivalence \[\nicexy{A \ar@<2pt>[r]^-{L} & A_{\st} \ar@<2pt>[l]^-{R}}\] in which:
\begin{enumerate}[label=(\roman*)]
\item $A_{\st}$ is a strict braided monoidal category.
\item Both $L$ and $R$ are strong braided monoidal functors.
\end{enumerate} 
\end{theorem}

Next is the $\G=\B$ special case of Theorem \ref{gmonoidal-coherence3}, which provides a strictification of a free braided monoidal category.  It is called the Coherence for Braided Tensor Categories in \cite{joyal-street} Theorem 2.5.

\begin{theorem}\label{free-braided-monoidal-strictification}
For each small category $A$, the strict braided monoidal functor \[\nicexy{\MCatb(A) \ar[r]^-{\Lambdab_A} & \MCatbst(A)}\] in \eqref{lambdag-a} is an equivalence of categories, where $\MCatb(A)$ (resp., $\MCatbst(A)$) is the free (strict) braided monoidal category generated by $A$.
\end{theorem}

The following result is the $\G=\B$ special case of Theorem \ref{gmonoidal-coherence4}, which says that the canonical functor from the free braided monoidal category generated by a set of objects to the free strict braided monoidal category generated by one object is faithful.

\begin{theorem}\label{bmonoidal-coherence4}
For each set $A$ of objects, the composite functor \[\nicexy{\MCatb(A) \ar[r]^-{\Lambdab_A} & \MCatbst(A) \ar[r]^-{\pi} & \MCatbst(*)}\] is faithful.
\end{theorem}

Next is the $\G=\B$ special case of Corollary \ref{gmonoidal-coherence-underlying}, which recovers \cite{joyal-street} Corollary 2.6.

\begin{corollary}\label{braided-monoidal-coherence-underlying}
For each set $A$ of objects, a diagram in the free braided monoidal category $\MCatb(A)$ generated by $A$ is commutative if and only if composites with the same domain and the same codomain have the same $\pi\Lambdab_A$-image in $\MCatbst(*)$.
\end{corollary}

\begin{interpretation}
Recall from Example \ref{ex:free-strict-gmon-one-element}, restricted to the braid group operad $\B$, that the free strict braided monoidal category generated by one object, $\MCatbst(*)$, has objects the non-negative integers $\{0,1,2,\ldots\}$.  With $B_n$ denoting the $n$th braid group, the morphism sets are \[\MCatbst(*)(m;n) = \begin{cases} \varnothing & \text{if $m\not=n$};\\ B_n & \text{if $m=n$}.\end{cases}\] Using the explicit description of the composite $\pi\Lambdab_A$ in the proof of Theorem \ref{gmonoidal-coherence4-alt}, it makes sense to call \[\pi\Lambdab_A(\sigma) =\sigma \in B_n = \MCatbst(*)(n;n)\] the \emph{underlying braid} of each morphism $\sigma \in \MCatb(A)$ for a set $A$.  Therefore, Corollary \ref{braided-monoidal-coherence-underlying} says that, in the free braided monoidal category $\MCatb(A)$ generated by a set $A$ of objects, a diagram is commutative if and only if composites with the same (co)domain have the same underlying braid.\dqed
\end{interpretation}

\begin{remark}[Free Strict Braided Monoidal Categories]
An explicit description of the free strict braided monoidal category $\MCatbst(A)$ generated by a small category $A$ is given in Theorem \ref{free-strict-gmon-category} by restricting to the case $\G=\B$.  This description coincides with the one given in \cite{joyal-street} Example 2.1 and Proposition 2.2.\dqed
\end{remark}

\begin{remark}[Free Braided Monoidal Categories]
Similarly, an explicit description of the free braided monoidal category $\MCatb(A)$ generated by a small category $A$ is given in Theorem \ref{free-gmon-category}.  The $\B$-monoidal category operad $\MCatb$ is called the \emph{operad of unital parenthesized braids} in \cite{horel} Construction 6.3.  Citing the proof in \cite{fresse-gt} 6.2.7, in \cite{horel} Remark 6.5 it is briefly mentioned that the operad of unital parenthesized braids encodes small braided monoidal categories, i.e., Theorem \ref{bmonoidal-braided-monoidal}.\dqed
\end{remark}

\section{Symmetric Monoidal Categories are $\S$-Monoidal Categories}
\label{sec:symmetric-monoidal-cat}

The purpose of this section is to observe that for the symmetric group operad $\S$ in Example \ref{ex:symmetric-group-operad}, $\S$-monoidal categories in the sense of Definition \ref{def:g-monoidal-category} (i.e., $\MCats$-algebras) are exactly small symmetric monoidal categories in the usual sense.  Let us first recall the definition of a symmetric monoidal category from \cite{maclane} (VII.7).

\begin{definition}\label{def:symmetric-monoidal-category}
A \index{category!symmetric monoidal}\index{symmetric monoidal category}\index{monoidal category!symmetric}\emph{symmetric monoidal category} is a pair $\left(\M, \xi\right)$ in which:
\begin{itemize}
\item $\M = (\M,\otimes,\tensorunit,\alpha,\lambda,\rho)$ is a monoidal category as in Definition \ref{def:monoidal-category}.
\item $\xi$ is a natural isomorphism\label{notation:symmetry-iso}
\begin{equation}\label{symmetry-isomorphism}
\nicexy{X \otimes Y \ar[r]^-{\xi_{X,Y}}_-{\cong} & Y \otimes X}
\end{equation}
for objects $X,Y \in \M$, called the \index{symmetry isomorphism}\emph{symmetry isomorphism}.
\end{itemize}
This data is required to satisfy the following three axioms.
\begin{description}
\item[Symmetry Axiom]
The diagram
\begin{equation}\label{monoidal-symmetry-axiom}
\nicexy{X \otimes Y \ar[r]^-{\xi_{X,Y}} \ar@{=}[dr] & Y \otimes X \ar[d]^-{\xi_{Y,X}}
\\ & X \otimes Y}
\end{equation}
is commutative for all objects $X,Y \in \M$.
\item[Unit Axiom]
The diagram
\begin{equation}\label{symmetry-unit}
\nicexy{X \otimes \tensorunit \ar[d]_-{\rho} \ar[r]^-{\xi_{X,\tensorunit}}
& \tensorunit \otimes X \ar[d]^-{\lambda}\\ X \ar@{=}[r] & X}
\end{equation}
is commutative for all objects $X \in \M$.
\item[Hexagon Axiom]
The diagram\index{Hexagon Axioms!for symmetric monoidal categories}
\begin{equation}
\label{hexagon-axiom}
\nicexy@C-40pt{& X \otimes (Z \otimes Y) \ar[rr]^-{\Id_X \otimes \xi_{Z,Y}} & \hspace{1.8cm}
& X \otimes (Y \otimes Z) \ar[dr]^-{\alpha^{-1}}&\\ (X \otimes Z) \otimes Y 
\ar[ur]^-{\alpha} \ar[dr]_-{\xi_{X\otimes Z,Y}} &&&& (X \otimes Y) \otimes Z\\
& Y \otimes (X \otimes Z) \ar[rr]^-{\alpha^{-1}} & \hspace{1.5cm}
& (Y \otimes X) \otimes Z \ar[ur]_-{\xi_{Y,X} \otimes \Id_Z} &}
\end{equation}
is commutative for all objects $X,Y, Z \in \M$.
\end{description}
\end{definition}

\begin{remark}\label{rk:symmetric-braided-moncat}
A symmetric monoidal category is a special kind of a braided monoidal category.  Indeed, if the braiding in a braided monoidal category satisfies the symmetry axiom \eqref{monoidal-symmetry-axiom}, $\xi_{X,Y}^{-1} = \xi_{Y,X}$, then the first hexagon diagram \eqref{hexagon-b1} implies the second hexagon diagram \eqref{hexagon-b2}, and vice versa.  Furthermore, each of these hexagon diagrams is equivalent to the one in \eqref{hexagon-axiom}.  In other words, a symmetric monoidal category is a braided monoidal category whose braiding satisfies the symmetry axiom \eqref{monoidal-symmetry-axiom}.\dqed
\end{remark}

Recall the $\G$-monoidal category operad $\MCatg$ from Definition \ref{def:mcatg-operad}, whose algebras are defined as $\G$-monoidal categories, for an action operad $\G$.  In the next observation, we restrict to the symmetric group operad $\S$ and use the $\S$-monoidal category operad\label{not:mcats} $\MCats$, whose algebras are by definition $\S$-monoidal categories.  Theorem \ref{gmonoidal-coherence-1}, which describes a $\G$-monoidal category structure as a monoidal category structure together with $\G$-monoidal structure isomorphisms $\sigma_*$ in \eqref{sigmastar-p-to-q}, will be used in the proof below.

\begin{theorem}\label{smonoidal-symmetric-monoidal}\index{symmetric monoidal category!operadic interpretation}
For each small category, the following two structures are equivalent:
\begin{enumerate}[label=(\roman*)]
\item An $\S$-monoidal category structure for the symmetric group operad $\S$.
\item A symmetric monoidal category structure.
\end{enumerate}
\end{theorem}

\begin{proof}
We recycle the proof of Theorem \ref{bmonoidal-braided-monoidal} and only point out the necessary modifications.  In that proof, we first replace the braid group operad $\B$ and the word \emph{braided} with the symmetric group operad $\S$ and the word \emph{symmetric}, respectively.  The decomposition \eqref{sigma-decomp-in-Bn} in the braid group is replaced by a decomposition \[\sigma = \tau_{i_k} \cdots \tau_{i_i} \in S_n\] with each $\tau_{i_j} = (i_j,i_j+1) \in S_n$ an adjacent transposition that permutes the two indicated consecutive objects.  This makes sense because the adjacent transpositions $\bigl\{\tau_1,\ldots,\tau_{n-1}\bigr\}$ generate the symmetric group $S_n$.  For each $1\leq j \leq k$, the canonical isomorphism $c_{i_j}$ in \eqref{ci-epsilon} now takes the form
\[c_{i_j} = \underbrace{\Id \otimes \cdots \otimes \Id}_{i_j-1} \otimes \xi \otimes \underbrace{\Id \otimes \cdots \otimes \Id}_{n-i_j-1}.\]

To check that the canonical isomorphism $\sigma_*$ is well-defined, we now use the fact that the symmetric group $S_n$ is generated by $\bigl\{\tau_1,\ldots,\tau_{n-1}\bigr\}$ and is subject to the braid relations in Definition \ref{def:braid-group} and the involution relation, $\tau_i^2 = \id$ for $1\leq i \leq n-1$, which corresponds to the symmetry axiom \eqref{monoidal-symmetry-axiom}.

In the second half of the proof, we replace the braid group generator $s_1 \in B_2$ with the non-trivial permutation $\tau_1 = (1,2) \in S_2$, which yields the unique morphism in $\MCats(2)\bigl(x_1x_2,x_2x_1\bigr)$.
\end{proof}

\begin{remark}[Fresse's Parenthesized Symmetry Operad]\label{rk:fresse-pas}
Fresse in \cite{fresse-gt} (Theorem 6.3.2) described a $1$-colored symmetric operad \textit{\textsf{PaS}} in $\Cat$, called the\index{parenthesized symmetry operad}\index{operad!parenthesized symmetry} \emph{parenthesized symmetry operad}, whose algebras are small symmetric monoidal categories without monoidal units.  The parenthesized symmetry operad can be obtained from the $\S$-monoidal category operad $\MCats$ by forgetting
\begin{enumerate}[label=(\roman*)]
\item the $0$th level $\MCats(0)$ and
\item all instances of $x_0$'s in $\MCats(n)$ for $n\geq 1$.
\end{enumerate}

A variation of \textit{\textsf{PaS}} is the $1$-colored symmetric operad\label{not:pasplus} \textit{\textsf{PaS}}$_+$ in $\Cat$, called the \emph{unitary parenthesized symmetry operad}.  Its algebras are small symmetric monoidal categories whose monoidal units are strict with respect to the monoidal product, the associativity isomorphism, the left and right unit isomorphisms, and the symmetry isomorphism, i.e., 
\[\begin{split}
\lambda_x &= \rho_x = \xi_{\tensorunit,x} = \xi_{x,\tensorunit} = \Id_x ,\\
\alpha_{\tensorunit,x,y} &= \alpha_{x,\tensorunit,y} = \alpha_{x,y,\tensorunit} = \Id_{x\otimes y}.
\end{split}\]
This symmetric operad is obtained from the parenthesized symmetry operad \textit{\textsf{PaS}} by adjoining an object at level $0$, i.e., \textit{\textsf{PaS}}$_+(0) = *$.  

There is a morphism $\pi : \MCats \to$ \!\!\textit{\textsf{PaS}}$_+$ of $1$-colored symmetric operads in $\Cat$ that sends each standard non-associative monomial of weight $0$ to $*$ in \textit{\textsf{PaS}}$_+(0)$.  For $n\geq 1$, $\pi$ sends each standard non-associative monomial $p$ of weight $n$ to the strict non-associative monomial $p_{\st}$, as in Definition \ref{def:mcat-terminology}, associated to $p$.\dqed
\end{remark}

We end this section with the strict analogue of Theorem \ref{smonoidal-symmetric-monoidal}.  Recall the strict $\G$-monoidal category operad $\MCatgst$ in Definition \ref{def:mcatgst-operad}.  We will use $\MCatgst$ when $\G$ is the symmetric group operad $\S$ in the following result.

\begin{corollary}\label{strict-smonoidal-symmetric-monoidal}
For each small category, the following two structures are equivalent:
\begin{enumerate}[label=(\roman*)]
\item An $\MCatsst$-algebra structure.
\item A strict symmetric monoidal category structure.
\end{enumerate}
\end{corollary}

\begin{proof}
By Theorem \ref{mcatgst-algebra} applied to $\G=\S$, an $\MCatsst$-algebra structure on $A$ is equal to an $\MCats$-algebra structure on $A$ whose underlying monoidal category is strict.  By Theorem \ref{smonoidal-symmetric-monoidal}, an $\MCats$-algebra structure is equal to a symmetric monoidal category structure.  Therefore, an $\MCatsst$-algebra structure is exactly a symmetric monoidal category structure whose underlying monoidal category is strict.
\end{proof}

\begin{remark}[Fresse's Colored Symmetry Operad]\label{rk:fresse-cos-operad}
Corollary \ref{strict-smonoidal-symmetric-monoidal} says that $\MCatsst$ is the $1$-colored symmetric operad in $\Cat$ that encodes small strict symmetric monoidal categories.  The symmetric operad $\MCatsst$ is equal to Fresse's unitary colored symmetry operad in \cite{fresse-gt} (Theorem 6.3.3), denoted\label{not:cosplus} \textit{\textsf{CoS}}$_+$.\dqed
\end{remark}

\begin{remark}[Gurski's Categorical Borel Construction]\label{rk:gurski-borel-sym}\index{categorical Borel construction}
In Remark \ref{rk:gurski-borel} we mentioned Gurski's \cite{gurski} categorical Borel construction \[\Cat \ni A \mapsto \coprod_{n\geq 0} (E_{\G(n)})\timesover{\G(n)} A^{\times n}\] for an action operad $\G$, in which $E_{\G(n)}$ is the translation category of $\G(n)$.  When $\G$ is the symmetric group operad $\S$, Gurski in \cite{gurski} (Example 4.1) observed that the algebras over his categorical Borel construction are small strict symmetric monoidal categories.\dqed
\end{remark}

\section{Symmetric Monoidal Functors are $\S$-Monoidal Functors}
\label{sec:symmetric-monoidal-functor}

The purpose of this section is to observe that, for the symmetric group operad $\S$, an $\S$-monoidal functor in the sense of Definition \ref{def:gmonoidal-functor} is exactly a symmetric monoidal functor in the usual sense.  The strong and strict analogues are also true.  As a consequence, the $\S$-monoidal category operad $\MCats$ is the symmetric operad in $\Cat$ that encodes small symmetric monoidal categories and strict symmetric monoidal functors.  Moreover, the $\S$-monoidal functor operad $\MFuns$ is the $2$-colored symmetric operad in $\Cat$ that encodes general symmetric monoidal functors between small symmetric monoidal categories.

Let us first recall the definition of a symmetric monoidal functor.

\begin{definition}\label{def:symmetric-monoidal-functor}
Suppose $\M$ and $\N$ are symmetric monoidal categories.  A \index{functor!symmetric monoidal}\index{symmetric monoidal functor}\index{monoidal functor!symmetric}\emph{symmetric monoidal functor} $(F,F_2,F_0) : \M \to \N$
is a monoidal functor between the underlying monoidal categories that is compatible with the symmetry isomorphisms, in the sense that the diagram
\begin{equation}\label{monoidal-functor-symmetry}
\nicexy@C+.5cm{F(X) \otimes F(Y) \ar[d]_-{F_2} \ar[r]^-{\xi_{FX,FY}}_-{\cong} & F(Y) \otimes F(X) \ar[d]^-{F_2} \\ F(X \otimes Y) \ar[r]^-{F\xi_{X,Y}}_-{\cong} & F(Y \otimes X)}
\end{equation}
is commutative for all objects $X,Y \in \M$.  A symmetric monoidal functor is said to be \emph{strong} (resp., \emph{strict}) if the underlying monoidal functor is so.
\end{definition}

\begin{remark} The definition of a symmetric monoidal functor is obtained from that of a braided monoidal functor in Definition \ref{def:braided-monoidal-functor} by replacing the word \emph{braided} with the word \emph{symmetric}.\dqed
\end{remark}  

\begin{example}
Suppose $(\M,\otimes,\tensorunit)$ is a symmetric monoidal category with all set-indexed coproducts.  Then the functor
\[\nicexy{\Set \ar[r]^-{\iota} & \M},\quad \iota(X) = \coprodover{x \in X} \tensorunit\]
is a strong symmetric monoidal functor.\dqed
\end{example}

\begin{example}
The singular chain functor\index{singular chain functor}\index{functor!singular chain} \[\CHau \to \Chainz\] is a symmetric monoidal functor \cite{massey} (XI.3).\dqed
\end{example}

Next is the main result of this section, which identifies symmetric monoidal functors with $\S$-monoidal functors for the symmetric group operad $\S$.

\begin{theorem}\label{smonoidal-functor}\index{symmetric monoidal functor!operadic interpretation}
Suppose $F : \M \to \N$ is a functor between small symmetric monoidal categories.  Then the following two structures on $F$ are equivalent:
\begin{enumerate}[label=(\roman*)]
\item A symmetric monoidal functor structure.
\item An $\S$-monoidal functor structure.
\end{enumerate}
\end{theorem}

\begin{proof}
We reuse the proof of Theorem \ref{bmonoidal-functor} by replacing:
\begin{itemize}
\item the braid group operad $\B$ by the symmetric group operad $\S$;
\item the words \emph{braided} and \emph{braiding} by the words \emph{symmetric} and \emph{symmetry isomorphism}, respectively;
\item the isomorphism $s_1\in \MCatb(2)$ by the isomorphism  \[(1,2) \in \MCats(2)\bigl(x_1x_2,x_2x_1\bigr)\] with $(1,2)$ the non-trivial permutation in $S_2$.
\end{itemize}
\end{proof}

\begin{interpretation}
Theorem \ref{smonoidal-functor} is a coherence result for symmetric monoidal functors.  Indeed, it asserts that for a monoidal functor $F : \M \to \N$ between small symmetric monoidal categories, in order for the diagram \eqref{gmonoidal-functor-diagram} to be commutative for each isomorphism in the $\S$-monoidal category operad $\MCats$, it is necessary and sufficient that the single diagram \eqref{monoidal-functor-symmetry} be commutative.\dqed
\end{interpretation}

Recall from Definition \ref{def:gmonoidal-functor} that a $\G$-monoidal functor is said to be \emph{strong} (resp., \emph{strict}) if the underlying monoidal functor is so.  The next observation is an immediate consequence of Theorem \ref{smonoidal-functor}.

\begin{corollary}\label{strong-smonoidal-functor}
Suppose $F : \M \to \N$ is a functor between small symmetric monoidal categories.  Then the following two structures on $F$ are equivalent:
\begin{enumerate}[label=(\roman*)]
\item A strong (resp., strict) symmetric monoidal functor structure.
\item A strong (resp., strict) $\S$-monoidal functor structure.
\end{enumerate}
\end{corollary}

We can now identify the algebras and the algebra morphisms of the $\S$-monoidal category operad $\MCats$ with well-known concepts.

\begin{corollary}\label{mcats-algebra}\index{symmetric monoidal category!operad}\index{operad!symmetric monoidal category}
$\MCats$ is a $1$-colored symmetric operad in $\Cat$ whose:
\begin{itemize}
\item algebras are small symmetric monoidal categories;
\item algebra morphisms are strict symmetric monoidal functors.
\end{itemize}
\end{corollary}

\begin{proof}
The assertion about algebras is Theorem \ref{smonoidal-symmetric-monoidal}.  The assertion about algebra morphisms follows from (i) the $\G=\S$ special case of Theorem \ref{mcatg-algebra-morphism}, which says that $\MCats$-algebra morphisms are strict $\S$-monoidal functors, and (ii) Corollary \ref{strong-smonoidal-functor}.
\end{proof}

Next is the strict analogue.

\begin{corollary}\label{mcatsst-algebra}\index{strict symmetric monoidal category!operad}\index{operad!strict symmetric monoidal category}
$\MCatsst$ is a $1$-colored symmetric operad in $\Cat$ whose:
\begin{itemize}
\item algebras are small strict symmetric monoidal categories;
\item algebra morphisms are strict symmetric monoidal functors.
\end{itemize}
\end{corollary}

\begin{proof}
The assertion about algebras is Corollary \ref{strict-smonoidal-symmetric-monoidal}.  The assertion about algebra morphisms follows from (i) the $\G=\S$ special case of Proposition \ref{mcatgst-algebra-morphism}, which says that $\MCatsst$-algebra morphisms are strict $\S$-monoidal functors, and (ii) Corollary \ref{strong-smonoidal-functor}.
\end{proof}

For an action operad $\G$, recall from Definition \ref{def:gmonoidal-functor-operad} the $\G$-monoidal functor operad $\MFung$.

\begin{corollary}\label{mfuns-algebra}\index{symmetric monoidal functor!operad}\index{operad!symmetric monoidal functor}
An $\MFuns$-algebra is exactly a triple $(A,B,F)$ such that:
\begin{itemize}
\item $A$ and $B$ are small symmetric monoidal categories.
\item $F : A \to B$ is a symmetric monoidal functor.
\end{itemize}
\end{corollary}

\begin{proof}
Theorem \ref{mfung-algebra} applied to the symmetric group operad $\S$ says that an $\MFuns$-algebra is a triple $(A,B,F)$ such that (i) $A$ and $B$ are $\S$-monoidal categories and (ii) $F : A\to B$ is an $\S$-monoidal functor.  The assertion now follows from Theorem \ref{smonoidal-symmetric-monoidal} and Theorem \ref{smonoidal-functor}.
\end{proof}

Next is the strong analogue involving the strong $\G$-monoidal functor operad $\MFungstg$ in Definition \ref{def:strong-gmonoidal-functor-operad}.

\begin{corollary}\label{mfunsstg-algebra}\index{operad!strong symmetric monoidal functor}
An $\MFunsstg$-algebra is exactly a triple $(A,B,F)$ such that:
\begin{itemize}
\item $A$ and $B$ are small symmetric monoidal categories.
\item $F : A \to B$ is a strong symmetric monoidal functor.
\end{itemize}
\end{corollary}

\begin{proof}
Theorem \ref{mfungstg-algebra} applied to the symmetric group operad $\S$ says that an $\MFunsstg$-algebra is a triple $(A,B,F)$ such that (i) $A$ and $B$ are $\S$-monoidal categories and (ii) $F : A\to B$ is a strong $\S$-monoidal functor.  The assertion now follows from Theorem \ref{smonoidal-symmetric-monoidal} and Corollary \ref{strong-smonoidal-functor}.
\end{proof}

\section{Coherence for Symmetric Monoidal Categories}
\label{sec:symmetric-monoidal-coherence}

The purpose of this section is to record various coherence statements for symmetric monoidal categories by restricting the results in Chapter \ref{ch:coherence-gmonoidal} to the symmetric group operad $\S$ in Example \ref{ex:symmetric-group-operad}.

Recall the following statements from the previous two sections.
\begin{enumerate}
\item In Theorem \ref{smonoidal-symmetric-monoidal} we proved that $\S$-monoidal categories and small symmetric monoidal categories are the same thing.  
\item The strict version is Corollary \ref{strict-smonoidal-symmetric-monoidal}, which says that strict $\S$-monoidal categories and small strict symmetric monoidal categories are the same thing.  
\item In Theorem \ref{smonoidal-functor} we established the equivalence between $\S$-monoidal functors and symmetric monoidal functors, with the strong/strict versions in Corollary \ref{strong-smonoidal-functor}.  
\end{enumerate}

Therefore, the coherence results in Chapter \ref{ch:coherence-gmonoidal} applies to symmetric monoidal categories and symmetric monoidal functors when $\G$ is the symmetric group operad $\S$.  First is the $\G=\S$ special case of Theorem \ref{gmonoidal-coherence2}, which says that symmetric monoidal categories can always be strictified.  This is the symmetric version of the Coherence Theorem for monoidal categories in \cite{maclane} (XI.3 Theorem 1).

\begin{theorem}\label{symmetric-monoidal-strictification}\index{strictification!symmetric monoidal category}\index{coherence!symmetric monoidal category}\index{symmetric monoidal category!coherence}
Suppose $A$ is a small symmetric monoidal category.  Then there exists an adjoint equivalence \[\nicexy{A \ar@<2pt>[r]^-{L} & A_{\st} \ar@<2pt>[l]^-{R}}\] in which:
\begin{enumerate}[label=(\roman*)]
\item $A_{\st}$ is a strict symmetric monoidal category.
\item Both $L$ and $R$ are strong symmetric monoidal functors.
\end{enumerate} 
\end{theorem}

Next is the $\G=\S$ special case of Theorem \ref{gmonoidal-coherence3}, which provides a strictification of a free symmetric monoidal category.

\begin{theorem}\label{free-symmetric-monoidal-strictification}
For each small category $A$, the strict symmetric monoidal functor \[\nicexy{\MCats(A) \ar[r]^-{\Lambdas_A} & \MCatsst(A)}\] in \eqref{lambdag-a} is an equivalence of categories, where $\MCats(A)$ (resp., $\MCatsst(A)$) is the free (strict) symmetric monoidal category generated by $A$.
\end{theorem}

The following result is the $\G=\S$ special case of Theorem \ref{gmonoidal-coherence4}, which says that the canonical functor from the free symmetric monoidal category generated by a set of objects to the free strict symmetric monoidal category generated by one object is faithful.

\begin{theorem}\label{smonoidal-coherence4}
For each set $A$, the composite functor \[\nicexy{\MCats(A) \ar[r]^-{\Lambdas_A} & \MCatsst(A) \ar[r]^-{\pi} & \MCatsst(*)}\] is faithful.
\end{theorem}

Next is the $\G=\S$ special case of Corollary \ref{gmonoidal-coherence-underlying}.

\begin{corollary}\label{symmetric-monoidal-coherence-underlying}
For each set $A$ of objects, a diagram in the free symmetric monoidal category $\MCats(A)$ generated by $A$ is commutative if and only if  composites with the same domain and the same codomain have the same $\pi\Lambdas_A$-image in $\MCatsst(*)$.
\end{corollary}

\begin{interpretation}
Recall from Example \ref{ex:free-strict-gmon-one-element}, restricted to the symmetric group operad $\S$, that the free strict symmetric monoidal category generated by one object, $\MCatsst(*)$, has objects the non-negative integers $\{0,1,2,\ldots\}$.  With $S_n$ denoting the $n$th symmetric group, the morphism sets are \[\MCatsst(*)(m;n) = \begin{cases} \varnothing & \text{if $m\not=n$};\\ S_n & \text{if $m=n$}.\end{cases}\] Using the explicit description of the composite $\pi\Lambdas_A$ in the proof of Theorem \ref{gmonoidal-coherence4-alt}, it makes sense to call \[\pi\Lambdas_A(\sigma) =\sigma \in S_n = \MCatsst(*)(n;n)\] the \emph{underlying permutation} of each morphism $\sigma \in \MCats(A)$ for a set $A$.  Therefore, Corollary \ref{symmetric-monoidal-coherence-underlying} says that, in the free symmetric monoidal category $\MCats(A)$ generated by a set $A$ of objects, a diagram is commutative if and only if  composites with the same (co)domain have the same underlying permutation.\dqed
\end{interpretation}

\begin{remark}[Free (Strict) Symmetric Monoidal Categories]
Explicit descriptions of the free strict symmetric monoidal category $\MCatsst(A)$ and of the free symmetric monoidal category $\MCats(A)$ generated by a small category $A$ are given in Theorem \ref{free-strict-gmon-category} and Theorem \ref{free-gmon-category}, respectively, by restricting to the case $\G=\S$.\dqed
\end{remark}

\chapter{Ribbon Monoidal Categories}
\label{ch:ribbon-moncat}

Recall from Definition \ref{def:braided-monoidal-category} that a braided monoidal category is a monoidal category $\M$ equipped with a natural braiding \[\nicexy{X\otimes Y \ar[r]^-{\xi_{X,Y}}_-{\cong} & Y \otimes X}\] that is  compatible with the associativity isomorphism, the left unit isomorphism, and the right unit isomorphism, as expressed by the two hexagon axioms and the unit axiom.  This chapter is about ribbon monoidal categories, which are extensions of braided monoidal categories, and their coherence.  Ribbon monoidal categories are called \emph{balanced tensor categories} in \cite{joyal-street} and \emph{ribbon braided monoidal categories} in \cite{wahl}.  Just like braided monoidal categories, ribbon monoidal categories naturally arise in the representation theory of quantum groups; see \cite{chari,kassel,street}.

A \emph{ribbon monoidal category} is a braided monoidal category that is furthermore equipped with a natural isomorphism \[\nicexy{X \ar[r]^-{\theta_X}_-{\cong} & X},\] called the \emph{twist}, that is compatible with the monoidal unit and the braiding.  The compatibility with the monoidal unit says that the twist applied to the monoidal unit is equal to the identity morphism.  The compatibility with the braiding says that the twist applied to a monoidal product $X\otimes Y$ is equal to the composite of the braiding, the twist on each tensor factor, and the braiding again.  

Recall the ribbon group operad $\R$ in Definition \ref{def:ribbon-group-operad}.  Algebras over the $\R$-monoidal category operad $\MCatr$ in Definition \ref{def:mcatg-operad} are called $\R$-monoidal categories.  In Section \ref{sec:ribbon-monoidal-cat} we show that an $\R$-monoidal category is exactly a small ribbon monoidal category with general associativity isomorphism, left/right unit isomorphisms, braiding, and twist.  The strict variant of this result is also true and has appeared in the literature before.  In fact, the strict $\R$-monoidal category operad $\MCatrst$ is equal to Wahl's ribbon categorical operad in \cite{wahl}; see Remark \ref{rk:wahl-operad}. Moreover, Gurski's \cite{gurski} categorical Borel construction, when applied to the ribbon group operad $\R$, has as algebras small strict ribbon monoidal categories.  

In Section \ref{sec:ribbon-monoidal-functor} we observe that an $\R$-monoidal functor is exactly a ribbon monoidal functor, which is a braided monoidal functor that is compatible with the twists.  The structure morphisms of our ribbon monoidal functors are not invertible in general.  Ribbon monoidal functors with invertible structure morphisms are said to be \emph{strong}; these are what Joyal and Street \cite{joyal-street} called \emph{balanced  tensor functors}.  We also observe that a strong/strict $\R$-monoidal functor is exactly a strong/strict ribbon monoidal functor.

In Section \ref{sec:ribbon-monoidal-cat-coherence} we recover known coherence results for ribbon monoidal categories by restricting the coherence results in Chapter \ref{ch:coherence-gmonoidal} to the ribbon group operad $\R$.  The first coherence result says that every small ribbon monoidal category can be strictified via an adjoint equivalence of strong ribbon monoidal functors.  The free analogue of this coherence result says that the free ribbon monoidal category generated by a small category can be strictified to its free strict ribbon monoidal category via a strict ribbon monoidal equivalence.  Another coherence result says that, in the free ribbon monoidal category generated by a set of objects, a diagram is commutative if and only if composites with the same (co)domain have the same underlying ribbon.

\section{Ribbon Monoidal Categories are $\R$-Monoidal Categories}
\label{sec:ribbon-monoidal-cat}

The purpose of this section is to show that for the ribbon group operad $\R$ in Definition \ref{def:ribbon-group-operad}, $\R$-monoidal categories as in Definition \ref{def:g-monoidal-category} (i.e., $\MCatr$-algebras) are exactly small ribbon monoidal categories.  The strict version is also true.  Let us first recall the definition of a ribbon monoidal category, which is what Joyal and Street \cite{joyal-street} called a\index{balanced tensor category} \emph{balanced tensor category}.  

\begin{definition}\label{def:ribbon-monoidal-category}\index{ribbon monoidal category}\index{monoidal category!ribbon}
A \emph{ribbon monoidal category} is a pair $(\M,\theta)$ in which:
\begin{itemize}
\item $(\M,\otimes,\tensorunit,\alpha,\lambda,\rho,\xi)$ is a braided monoidal category as in Definition \ref{def:braided-monoidal-category}.
\item $\theta$ is a natural isomorphism 
\begin{equation}\label{twist-ribbon-moncat}
\nicexy{X \ar[r]^-{\theta_X}_-{\cong} & X}
\end{equation}
for objects $X \in \M$, called the \emph{twist}.\index{twist}
\end{itemize}
This data is required to satisfy the following two axioms.
\begin{description}
\item[Identity Axiom] $\theta_{\tensorunit} = \Id_{\tensorunit} : \tensorunit \iso \tensorunit$.
\item[Pair Twisting Axiom] The diagram\index{Pair Twisting Axiom}
\begin{equation}\label{pair-twist-diagram}
\nicexy@C+.3cm{X \otimes Y \ar[r]^-{\theta_{X\otimes Y}} \ar[d]_-{\xi_{X,Y}} & X \otimes Y\\
Y \otimes X \ar[r]^-{\theta_Y \otimes \theta_X} & Y \otimes X \ar[u]_-{\xi_{Y,X}}}
\end{equation}
is commutative for all objects $X,Y\in \M$.
\end{description}
A ribbon monoidal category is \emph{strict}\index{strict ribbon monoidal category} if the underlying monoidal category is strict.
\end{definition}

\begin{interpretation}\label{int:ribbon-moncat-axioms}
One may visualize the twist $\theta_X$ as a strip
\begin{center}\begin{tikzpicture}[xscale=1.2, yscale=.6, thick]
\node at (-.5,.5) {$\theta_X = $};
\foreach \y in {0,.5}{\draw (0,\y) to[out=90,in=270] (.2,\y+.5);
\draw[shorten >=1pt, shorten <=1pt, line width=5pt, white] (.2,\y) to[out=90,in=270] (0,\y+.5);
\draw (.2,\y) to[out=90,in=270] (0,\y+.5);}
\foreach \x in {0,1} \draw (0,\x) -- (.2,\x); 
\foreach \x in {-.4,1.4} \node at (.1,\x) {\footnotesize{$X$}};
\end{tikzpicture}\end{center}
labeled by an object $X$ with a full $2\pi$ twist.  This twisted strip is the generator $r_1$ in the first ribbon group $R_1$; see Example \ref{ex:generator-rn}.  

The pair twisting axiom \eqref{pair-twist-diagram} may be visualized as the picture below.
\begin{center}\begin{tikzpicture}[xscale=.7, yscale=.5, thick]
\draw [gray!50, dotted, line width=2pt] (1,0) -- (1,4);
\foreach \y in {0,2}{\draw (0,\y) to[out=90,in=270] (2,\y+2);
\draw[shorten >=1pt, shorten <=1pt, line width=5pt, white] (2,\y) to[out=90,in=270] (0,\y+2);eqref
\draw (2,\y) to[out=90,in=270] (0,\y+2);}
\foreach \x in {0,4} \draw (0,\x) -- (2,\x); 
\node at (3.3,2) {\Huge{$=$}};
\foreach \x in {-.4,4.4} {\node at (.5,\x) {\footnotesize{$X$}};
\node at (1.5,\x) {\footnotesize{$Y$}};}
\end{tikzpicture}\qquad
\begin{tikzpicture}[xscale=1.2, yscale=.5, thick]
\foreach \x in {1,1.2}{\draw (\x,0) -- (\x+1,1);} 
\draw[shorten >=.4cm, shorten <=.4cm, line width=8pt, white] (2.1,0) -- (1.1,1);
\foreach \x in {2,2.2}{\draw (\x,0) -- (\x-1,1);}
\foreach \y in {1,2}{\draw (1,\y) to[out=90,in=270] (1.2,\y+1);
\draw[shorten >=1pt, shorten <=1pt, line width=5pt, white] (1.2,\y) to[out=90,in=270] (1,\y+1); 
\draw (1.2,\y) to[out=90,in=270] (1,\y+1);}
\foreach \y in {1,2}{\draw (2,\y) to[out=90,in=270] (2.2,\y+1);
\draw[shorten >=1pt, shorten <=1pt, line width=5pt, white] (2.2,\y) to[out=90,in=270] (2,\y+1); 
\draw (2.2,\y) to[out=90,in=270] (2,\y+1);}
\foreach \x in {1,1.2}{\draw (\x,3) -- (\x+1,4);} 
\draw[shorten >=.4cm, shorten <=.4cm, line width=8pt, white] (2.1,3) -- (1.1,4);
\foreach \x in {2,2.2}{\draw (\x,3) -- (\x-1,4);}
\foreach \y in {0,4}{\draw (1,\y) -- (1.2,\y); \draw (2,\y) -- (2.2,\y);}
\foreach \x in {-.4,4.4} {\node at (1.1,\x) {\footnotesize{$X$}};
\node at (2.1,\x) {\footnotesize{$Y$}};}
\end{tikzpicture}\end{center}
The left-hand side represents the twist $\theta_{X\otimes Y}$, with two side-by-side strips labeled by $X$ and $Y$ undergoing a full $2\pi$ twist.  Considering the effects on the two individual strips, this can also be achieved as the picture on the right-hand side, which represents the composite
\[\xi_{Y,X} \circ \bigl(\theta_Y \otimes \theta_X\bigr) \circ \xi_{X,Y}\]
in the pair twisting axiom \eqref{pair-twist-diagram}.\dqed
\end{interpretation}

Next is the main result of this section, in which $\R$-monoidal categories\label{not:mcatr} (i.e., $\MCatr$-algebras) are identified with small ribbon monoidal categories.

\begin{theorem}\label{rmonoidal-ribbon-monoidal}\index{ribbon monoidal category!operadic interpretation}
For each small category, the following two structures are equivalent:
\begin{enumerate}[label=(\roman*)]
\item An $\R$-monoidal category structure for the ribbon group operad $\R$.
\item A ribbon monoidal category structure.
\end{enumerate}
\end{theorem}

\begin{proof}
We recycle the proof of Theorem \ref{bmonoidal-braided-monoidal} and point out the necessary modifications.  First we replace the braid group operad $\B$ with the ribbon group operad $\R$.

In the first half of the proof, we suppose $(A,\theta)$ is a ribbon monoidal category, and we want to define an $\R$-monoidal category structure on it.  For each isomorphism $\sigma \in \MCatr(n)(p,q)$, we consider a decomposition
\[\sigma = r_{i_k}^{\epsilon_k} \cdots r_{i_1}^{\epsilon_1} \in R_n\]
with each $r_{i_j} \in R_n$ a generating ribbon as in Interpretation \ref{int:ribbon-group-generators}, $1\leq i_j \leq n$, and each $\epsilon_j\in \{1,-1\}$.  Recall that the generating ribbons $\{r_1,\ldots,r_{n-1}\}$ satisfy the braid relations in Definition \ref{def:braid-group}.  In the rest of the proof, for $n\geq 2$ and $1\leq i \leq n-1$, the generating ribbon $r_i \in R_n$ plays the role of the generating braid $s_i \in B_n$ in the proof of Theorem \ref{bmonoidal-braided-monoidal}.

A canonical isomorphism is also allowed to involve the twist $\theta$ and its inverse.  In the definition of the canonical isomorphism $\sigma_*$ in \eqref{sigmastar-braided-monoidal}, whenever $i_j=n$ in the above decomposition of $\sigma$ for some $1\leq j \leq k$, we modify its construction as follows.
\begin{itemize}
\item $\alpha_j : p_{\st}^{2j-2} \to p_{\st}^{2j-1}$ is the identity morphism.
\item $c_{i_j}^{\epsilon_j}$ in \eqref{ci-epsilon} is now the monoidal product
\[c_{i_j}^{\epsilon_j} = \underbrace{\Id \otimes \cdots \otimes \Id}_{n-1} \otimes \theta^{\epsilon_1}.\]
\end{itemize}
To check that $\sigma_*$ is well-defined, we also need to check that the two sides of the generating relation \eqref{ribbon-new-relation} \[r_{n-1}r_nr_{n-1}r_n = r_nr_{n-1}r_nr_{n-1} \in R_n\] induce the same canonical isomorphism.  Since we are only concerned about two consecutive objects, with the identity morphisms in all other entries, it is sufficient to show that the outermost diagram in the diagram
\[\nicexy@C+.3cm{x\otimes y \ar[d]_-{\xi_{x,y}} \ar[r]^-{\Id_x\otimes \theta_y} & x \otimes y \ar[r]^-{\xi_{x,y}} & y \otimes x \ar[d]^-{\Id_y \otimes \theta_x}\\
y \otimes x \ar[d]_-{\Id_y \otimes \theta_x} \ar[rr]|-{\theta_y\otimes \theta_x} \ar[urr]|-{\theta_y \otimes \Id_x} && y \otimes x \ar[d]^-{\xi_{y,x}}\\
y\otimes x \ar[r]_-{\xi_{y,x}} \ar[urr]|-{\theta_y \otimes \Id_x} & x \otimes y \ar[r]_-{\Id_x \otimes \theta_y} & x \otimes y}\]
is commutative.  The upper-left and lower-right triangles are commutative by the naturality of the braiding $\xi$.  The other two triangles are commutative by the naturality of the monoidal product.

In the second half of the proof, given an $\R$-monoidal category $A$, we first equip $A$ with the structure of a braided monoidal category as in the second half of the proof of Theorem \ref{bmonoidal-braided-monoidal}, with the generating ribbon $r_1\in R_2$ playing the role of the generating braid $s_1 \in B_2$.  Next we define the twist in $A$ as the $\R$-monoidal structure isomorphism \[\nicexy@C+.5cm{x \ar[r]^-{\theta_x \,=\, (r_1)_*} & x} \in A\] in which \[r_1\in \MCatr(1)(x_1,x_1)\] with $r_1$ the generating ribbon in $R_1$. 

Now we check the axioms of a ribbon monoidal category in $A$.
\begin{itemize}
\item The identity axiom--i.e., $\theta_{\tensorunit} = \Id_{\tensorunit}$--is satisfied because the morphism set $\MCatr(0)(x_0,x_0)$ contains only the identity morphism, since the $0$th ribbon group $R_0$ is the trivial group. 
\item The pair twisting axiom \eqref{pair-twist-diagram} follows from the equality \[r_1 \comp_1 \id_2 = r_1 \cdot \bigl(r_1 \oplus r_1\bigr) \cdot r_1 \in R_2,\] which is true because the two sides are equal as braid on two strips as explained in Interpretation \ref{int:ribbon-moncat-axioms}.  In the above equality, in the direct sum $r_1\oplus r_1$ and on the left-hand side, each copy of $r_1 \in R_1$ is the generating ribbon, which corresponds to the twist in $A$.  The other two $r_1$'s on the right-hand side are the generating ribbon $r_1 \in R_2$, which yields the braiding in $A$.
\end{itemize} 
Therefore, $A$ has the structure of a ribbon monoidal category.
\end{proof}

Next is the strict analogue of Theorem \ref{rmonoidal-ribbon-monoidal}.  We will use the strict $\G$-monoidal category operad $\MCatgst$ in Definition \ref{def:mcatgst-operad} when $\G$ is the ribbon group operad $\R$.

\begin{corollary}\label{strict-rmonoidal-ribbon-monoidal}\index{strict ribbon monoidal category!operadic interpretation}
For each small category, the following two structures are equivalent:
\begin{enumerate}[label=(\roman*)]
\item An $\MCatrst$-algebra structure.
\item A strict ribbon monoidal category structure.
\end{enumerate}
\end{corollary}

\begin{proof}
By Theorem \ref{mcatgst-algebra} applied to $\G=\R$, an $\MCatrst$-algebra structure on $A$ is equal to an $\MCatr$-algebra structure on $A$ whose underlying monoidal category is strict.  By Theorem \ref{rmonoidal-ribbon-monoidal}, an $\MCatr$-algebra structure is equal to a ribbon monoidal category structure.  Therefore, an $\MCatrst$-algebra structure is exactly a ribbon monoidal category structure whose underlying monoidal category is strict.
\end{proof}

\begin{remark}[Earlier and Related Incarnations]\label{rk:wahl-operad}
The strict $\R$-monoidal category operad $\MCatrst$ in Corollary \ref{strict-rmonoidal-ribbon-monoidal}, which encodes small strict ribbon monoidal categories, has appeared in the literature before.
\begin{enumerate}
\item $\MCatrst$ is equal to Wahl's \emph{ribbon categorical operad} $\mathcal{C}^{R\beta}$ in \cite{wahl} Example 1.2.9.  Corollary \ref{strict-rmonoidal-ribbon-monoidal} recovers \cite{wahl} Theorem 1.4.7.
\item When $\G$ is the ribbon group operad $\R$, Gurski in \cite{gurski} (Example 4.2) mentioned that the algebras over his categorical Borel construction are small strict ribbon monoidal categories.\dqed
\end{enumerate}
\end{remark}

\section{Ribbon Monoidal Functors are $\R$-Monoidal Functors}
\label{sec:ribbon-monoidal-functor}

The purpose of this section is to observe that, for the ribbon group operad $\R$, an $\R$-monoidal functor in the sense of Definition \ref{def:gmonoidal-functor} is exactly a ribbon monoidal functor in the usual sense.  The strong and strict analogues are also true.  As a consequence, the $\R$-monoidal category operad $\MCatr$ is the symmetric operad in $\Cat$ that encodes ribbon monoidal categories and strict ribbon monoidal functors.  Moreover, the $\R$-monoidal functor operad $\MFunr$ is the $2$-colored symmetric operad in $\Cat$ that encodes general ribbon monoidal functors between small ribbon monoidal categories.

Let us first recall the definition of a ribbon monoidal functor.

\begin{definition}\label{def:ribbon-monoidal-functor}
Suppose $\M$ and $\N$ are ribbon monoidal categories.  A \index{functor!ribbon monoidal}\index{ribbon monoidal functor}\index{monoidal functor!ribbon}\emph{ribbon monoidal functor} $(F,F_2,F_0) : \M \to \N$
is a monoidal functor between the underlying monoidal categories that is compatible with:
\begin{enumerate}[label=(\roman*)]
\item the braidings in the sense of \eqref{monoidal-functor-braided};
\item the twists in the sense of the equality
\[\nicexy{F(\theta_X) = \theta_{F(X)} : F(X) \ar[r] & F(X)}\]
for objects $X \in \M$. 
\end{enumerate}
A ribbon monoidal functor is said to be \emph{strong} (resp., \emph{strict}) if the underlying monoidal functor is so.
\end{definition}

\begin{remark} A ribbon monoidal functor forgets to a braided monoidal functor between the underlying braided monoidal categories.  A strong ribbon monoidal functor--i.e., a ribbon monoidal functor with $F_2$ and $F_0$ invertible--is what Joyal and Street \cite{joyal-street} called a\index{balanced tensor functor} \emph{balanced tensor functor}.\dqed
\end{remark}

The main result of this section identifies ribbon monoidal functors with $\R$-monoidal functors for the ribbon group operad $\R$.

\begin{theorem}\label{rmonoidal-functor}\index{ribbon monoidal functor!operadic interpretation}
Suppose $F : \M \to \N$ is a functor between small ribbon monoidal categories.  Then the following two structures on $F$ are equivalent:
\begin{enumerate}[label=(\roman*)]
\item A ribbon monoidal functor structure.
\item An $\R$-monoidal functor structure.
\end{enumerate}
\end{theorem}

\begin{proof}
We reuse the proof of Theorem \ref{bmonoidal-functor} and point out the necessary modifications.  For an $\R$-monoidal functor $F$, the commutative diagram \eqref{gmonoidal-functor-diagram} for the isomorphism \[r_1 \in \MCatr(2)\bigl(x_1x_2,x_2x_1\bigr)\] gives the commutative diagram \eqref{monoidal-functor-braided}.  On the other hand, for the isomorphism \[r_1\in \MCatr(1)(x_1,x_1),\] the commutative diagram \eqref{gmonoidal-functor-diagram} gives the equality $F(\theta_X) = \theta_{F(X)}$.

In the second half of the proof, we note that each ribbon monoidal functor is compatible with the associativity isomorphisms, the left/right unit isomorphisms, the braidings, and the twists.  For each isomorphism $\sigma \in \MCatr$, by the construction of $\sigma_*$ in the proof of Theorem \ref{rmonoidal-ribbon-monoidal}, the diagram \eqref{gmonoidal-functor-diagram} factors into a ladder diagram, in which each square is commutative by the previous sentence.
\end{proof}

\begin{interpretation}
Theorem \ref{rmonoidal-functor} asserts that for a monoidal functor $F : \M \to \N$ between small ribbon monoidal categories, in order for the diagram \eqref{gmonoidal-functor-diagram} to be commutative for each isomorphism in the $\R$-monoidal category operad $\MCatr$, it is necessary and sufficient that (i) the diagram \eqref{monoidal-functor-braided} be commutative and (ii) $F(\theta_X) = \theta_{F(X)}$ for each object $X\in\M$.  So it is a coherence result for ribbon monoidal functors.\dqed
\end{interpretation}

Recall from Definition \ref{def:gmonoidal-functor} that a $\G$-monoidal functor is said to be \emph{strong} (resp., \emph{strict}) if the underlying monoidal functor is so.  The next observation is an immediate consequence of Theorem \ref{rmonoidal-functor}.

\begin{corollary}\label{strong-rmonoidal-functor}
For a functor $F : \M \to \N$ between small ribbon monoidal categories, the following two structures on $F$ are equivalent:
\begin{enumerate}[label=(\roman*)]
\item A strong (resp., strict) ribbon monoidal functor structure.
\item A strong (resp., strict) $\R$-monoidal functor structure.
\end{enumerate}
\end{corollary}

We can now identify the algebras and the algebra morphisms of the $\R$-monoidal category operad $\MCatr$ with well-known concepts.

\begin{corollary}\label{mcatr-algebra}\index{ribbon monoidal category!operad}\index{operad!ribbon monoidal category}
$\MCatr$ is a $1$-colored symmetric operad in $\Cat$ whose:
\begin{itemize}
\item algebras are small ribbon monoidal categories;
\item algebra morphisms are strict ribbon monoidal functors.
\end{itemize}
\end{corollary}

\begin{proof}
This follows from Theorem \ref{rmonoidal-ribbon-monoidal}, the $\G=\R$ special case of Theorem \ref{mcatg-algebra-morphism}, and Corollary \ref{strong-rmonoidal-functor}.
\end{proof}

Next is the strict analogue.

\begin{corollary}\label{mcatrst-algebra}\index{strict ribbon monoidal category!operad}\index{operad!strict ribbon monoidal category}
$\MCatrst$ is a $1$-colored symmetric operad in $\Cat$ whose:
\begin{itemize}
\item algebras are small strict ribbon monoidal categories;
\item algebra morphisms are strict ribbon monoidal functors.
\end{itemize}
\end{corollary}

\begin{proof}
This follows from Corollary \ref{strict-rmonoidal-ribbon-monoidal}, the $\G=\R$ special case of Proposition \ref{mcatgst-algebra-morphism}, and Corollary \ref{strong-rmonoidal-functor}.
\end{proof}

For an action operad $\G$, recall from Definition \ref{def:gmonoidal-functor-operad} the $\G$-monoidal functor operad $\MFung$.

\begin{corollary}\label{mfunr-algebra}\index{ribbon monoidal functor!operad}\index{operad!ribbon monoidal functor}
An $\MFunr$-algebra is exactly a triple $(A,B,F)$ such that:
\begin{itemize}
\item $A$ and $B$ are small ribbon monoidal categories.
\item $F : A \to B$ is a ribbon monoidal functor.
\end{itemize}
\end{corollary}

\begin{proof}
This follows from Theorem \ref{mfung-algebra} applied to the ribbon group operad $\R$, Theorem \ref{rmonoidal-ribbon-monoidal}, and Theorem \ref{rmonoidal-functor}.
\end{proof}

Next is the strong analogue involving the strong $\G$-monoidal functor operad $\MFungstg$ in Definition \ref{def:strong-gmonoidal-functor-operad}.

\begin{corollary}\label{mfunrstg-algebra}\index{operad!strong ribbon monoidal functor}
An $\MFunrstg$-algebra is exactly a triple $(A,B,F)$ such that:
\begin{itemize}
\item $A$ and $B$ are small ribbon monoidal categories.
\item $F : A \to B$ is a strong ribbon monoidal functor.
\end{itemize}
\end{corollary}

\begin{proof}
This follows from Theorem \ref{mfungstg-algebra} applied to the ribbon group operad $\R$, Theorem \ref{rmonoidal-ribbon-monoidal}, and Corollary \ref{strong-rmonoidal-functor}.
\end{proof}

\section{Coherence for Ribbon Monoidal Categories}
\label{sec:ribbon-monoidal-cat-coherence}

The purpose of this section is to record several coherence statements for ribbon monoidal categories by restricting the coherence results for $\G$-monoidal categories in Chapter \ref{ch:coherence-gmonoidal} to the ribbon group operad $\R$ in Definition \ref{def:ribbon-group-operad}. 

In the previous two sections, we proved the following results.
\begin{enumerate}
\item In Theorem \ref{rmonoidal-ribbon-monoidal} we observed that $\R$-monoidal categories and small ribbon monoidal categories are the same thing.  
\item The strict version is Corollary \ref{strict-rmonoidal-ribbon-monoidal}, which says that strict $\R$-monoidal categories and small strict ribbon monoidal categories are the same thing.  
\item In Theorem \ref{rmonoidal-functor} we established the equivalence between $\R$-monoidal functors and ribbon monoidal functors, with the strong/strict versions in Corollary \ref{strong-rmonoidal-functor}. 
\end{enumerate}

Therefore, the coherence results in Chapter \ref{ch:coherence-gmonoidal} applies to ribbon monoidal categories and ribbon monoidal functors when $\G$ is the ribbon group operad $\R$.  First is the $\G=\R$ special case of Theorem \ref{gmonoidal-coherence2}, which says that ribbon monoidal categories can always be strictified.

\begin{theorem}\label{ribbon-monoidal-strictification}\index{strictification!ribbon monoidal category}\index{coherence!ribbon monoidal category}\index{ribbon monoidal category!coherence}
Suppose $A$ is a small ribbon monoidal category.  Then there exists an adjoint equivalence \[\nicexy{A \ar@<2pt>[r]^-{L} & A_{\st} \ar@<2pt>[l]^-{R}}\] in which:
\begin{enumerate}[label=(\roman*)]
\item $A_{\st}$ is a strict ribbon monoidal category.
\item Both $L$ and $R$ are strong ribbon monoidal functors.
\end{enumerate} 
\end{theorem}

Next is the $\G=\R$ special case of Theorem \ref{gmonoidal-coherence3}, which provides a strictification of a free ribbon monoidal category.  It is called the Coherence for Balanced Tensor Categories in \cite{joyal-street} Theorem 6.2.

\begin{theorem}\label{free-ribbon-monoidal-strictification}
For each small category $A$, the strict ribbon monoidal functor \[\nicexy{\MCatr(A) \ar[r]^-{\Lambdar_A} & \MCatrst(A)}\] in \eqref{lambdag-a} is an equivalence of categories, where $\MCatr(A)$ (resp., $\MCatrst(A)$) is the free (strict) ribbon monoidal category generated by $A$.
\end{theorem}

The following result is the $\G=\R$ special case of Theorem \ref{gmonoidal-coherence4}, which says that the canonical functor from the free ribbon monoidal category generated by a set of objects to the free strict ribbon monoidal category generated by one object is faithful.

\begin{theorem}\label{rmonoidal-coherence4}
For each set $A$ of objects, the composite functor \[\nicexy{\MCatr(A) \ar[r]^-{\Lambdar_A} & \MCatrst(A) \ar[r]^-{\pi} & \MCatrst(*)}\] is faithful.
\end{theorem}

Next is the $\G=\R$ special case of Corollary \ref{gmonoidal-coherence-underlying}, which recovers \cite{joyal-street} Corollary 6.3.

\begin{corollary}\label{ribbon-monoidal-coherence-underlying}
For each set $A$ of objects, a diagram in the free ribbon monoidal category $\MCatr(A)$ generated by $A$ is commutative if and only if composites with the same domain and the same codomain have the same $\pi\Lambdar_A$-image in $\MCatrst(*)$.
\end{corollary}

\begin{interpretation}
Recall from Example \ref{ex:free-strict-gmon-one-element}, restricted to the ribbon group operad $\R$, that the free strict ribbon monoidal category generated by one object, $\MCatrst(*)$, has objects the non-negative integers $\{0,1,2,\ldots\}$.  With $R_n$ denoting the $n$th ribbon group, the morphism sets are \[\MCatrst(*)(m;n) = \begin{cases} \varnothing & \text{if $m\not=n$};\\ R_n & \text{if $m=n$}.\end{cases}\] Using the explicit description of the composite $\pi\Lambdar_A$ in the proof of Theorem \ref{gmonoidal-coherence4-alt}, it makes sense to call \[\pi\Lambdar_A(\sigma) =\sigma \in R_n = \MCatrst(*)(n;n)\] the \emph{underlying ribbon} of each morphism $\sigma \in \MCatr(A)$ for a set $A$.  Therefore, Corollary \ref{ribbon-monoidal-coherence-underlying} says that, in the free ribbon monoidal category $\MCatr(A)$ generated by a set $A$ of objects, a diagram is commutative if and only if composites with the same (co)domain have the same underlying ribbon.\dqed
\end{interpretation}

\begin{remark}[Free (Strict) Ribbon Monoidal Categories]
An explicit description of the free strict ribbon monoidal category $\MCatrst(A)$ generated by a small category $A$ is given in Theorem \ref{free-strict-gmon-category} by restricting to the case $\G=\R$.  This description coincides with the one given in \cite{joyal-street} Example 6.4 and Proposition 6.1.  Similarly, an explicit description of the free ribbon monoidal category $\MCatr(A)$ generated by a small category $A$ is given in Theorem \ref{free-gmon-category}. \dqed
\end{remark}

\chapter{Coboundary Monoidal Categories}
\label{ch:coboundary-monoidal-category}

This chapter is about another variation of braided monoidal categories called coboundary monoidal categories, and their coherence.  Whereas a ribbon monoidal category is a braided monoidal category with the extra structure of a twist, a coboundary monoidal category does not have an underlying braided monoidal category.  A \emph{coboundary monoidal category} is a monoidal category equipped with a natural isomorphism \[\nicexy{X \otimes Y \ar[r]^-{\xi_{X,Y}}_-{\cong} & Y \otimes X},\] called the \emph{commutor}, that satisfies the symmetry axiom and the unit axiom of a symmetric monoidal category, and that is compatible with the associativity isomorphism in the form of the \emph{cactus axiom} \eqref{cactus-axiom}.

Coboundary monoidal categories were first introduced by Drinfel'd \cite{drinfeld}, who showed that the category of representations of a coboundary Hopf algebra is a coboundary monoidal category.  Furthermore, Henrique and Kamnitzer \cite{hen-kam} (Theorem 6) showed that the category of crystals of a finite dimensional complex reductive Lie algebra is a coboundary monoidal category.  The reader is referred to \cite{savage} for a survey of coboundary monoidal categories in the representation theory of quantum groups and crystals. 

In Section \ref{sec:coboundary-cat} we prove a few preliminary results about coboundary monoidal categories that are needed in the subsequent section.  These preliminary results say that each coboundary monoidal category is equipped with natural isomorphisms acting on multiple monoidal products that behave like the generating cacti in the cactus groups.  In particular, they satisfy categorical analogues of the involution axiom and the containment axiom in Definition \ref{def:cactus-group}.  

Recall from Definition \ref{def:cactus-group-operad} the cactus group operad $\Cac$.  Algebras over the $\Cac$-monoidal category operad $\MCatcac$ are called $\Cac$-monoidal categories.  In  Section \ref{sec:coboundary-monoidal-cacmonoidal} we observe that $\Cac$-monoidal categories are exactly small coboundary monoidal categories with general associativity isomorphism, left/right unit isomorphisms, and commutor.  The strict version is also true.  In Section \ref{sec:coboundary-monoidal-functor} we identify (strong/strict) $\Cac$-monoidal functors with (strong/strict) coboundary monoidal functors, which are defined exactly like braided monoidal functors.  

In Section \ref{sec:coboundary-monoidal-cat-coherence} we obtain several coherence results for coboundary monoidal categories by restricting the coherence results in Chapter \ref{ch:coherence-gmonoidal} to the cactus group operad $\Cac$.   The first coherence result says that every small coboundary monoidal category can be strictified via an adjoint equivalence consisting of strong coboundary monoidal functors.  The free analogue of this coherence result says that the free coboundary monoidal category generated by a small category can be strictified to its free strict coboundary monoidal category via a strict coboundary monoidal equivalence.  Another coherence result says that, in the free coboundary monoidal category generated by a set of objects, a diagram is commutative if and only if composites with the same (co)domain have the same underlying cactus.

\section{Natural Actions on Coboundary Monoidal Categories}
\label{sec:coboundary-cat}

The purposes of this section are (i) to recall the definition of a coboundary monoidal category and (ii) to obtain some preliminary results about them that are needed in the proof in the next section that identifies small coboundary monoidal categories with $\Cac$-monoidal categories.  These results, Propositions \ref{ck-involution} and \ref{ck-containment}, say that each coboundary monoidal category is equipped with natural isomorphisms acting on multiple monoidal products that behave like the generating cacti in the cactus groups.  Let us first recall the definition of a coboundary monoidal category, called a \emph{coboundary category} in \cite{drinfeld,hen-kam}. 

\begin{definition}\label{def:coboundary-monoidal-category}\index{coboundary monoidal category}\index{monoidal category!coboundary}
A \emph{coboundary monoidal category} is a pair $(\M,\xi)$ in which:
\begin{itemize}
\item $(\M,\otimes,\tensorunit,\alpha,\lambda,\rho)$ is a monoidal category as in Definition \ref{def:monoidal-category}.
\item $\theta$ is a natural isomorphism 
\begin{equation}\label{commutor}
\nicexy{X\otimes Y \ar[r]^-{\xi_{X,Y}}_-{\cong} & Y \otimes X}
\end{equation}
for objects $X,Y \in \M$, called the \emph{commutor}.\index{commutor}
\end{itemize}
This data is required to satisfy the following three axioms.
\begin{description}
\item[Symmetry Axiom] The diagram
\begin{equation}\label{coboundary-symmetry-axiom}
\nicexy{X \otimes Y \ar[r]^-{\xi_{X,Y}} \ar@{=}[dr] & Y \otimes X \ar[d]^-{\xi_{Y,X}}\\ & X \otimes Y}
\end{equation}
is commutative for all objects $X,Y \in \M$.
\item[Unit Axiom] The diagram 
\[\nicexy{X \otimes \tensorunit \ar[d]_-{\rho} \ar[r]^-{\xi_{X,\tensorunit}}
& \tensorunit \otimes X \ar[d]^-{\lambda}\\ X \ar@{=}[r] & X}\]
is commutative for all objects $X\in \M$.
\item[Cactus Axiom] The diagram\index{Cactus Axiom}
\begin{equation}\label{cactus-axiom}
\nicexy@C+.5cm{(X \otimes Y) \otimes Z \ar[r]^-{\xi_{X,Y}\otimes \Id_Z} \ar[d]_-{\alpha} & (Y\otimes X) \otimes Z \ar[r]^-{\xi_{Y\otimes X,Z}} & Z \otimes (Y \otimes X) \\
X \otimes (Y \otimes Z) \ar[r]^-{\Id_X\otimes \xi_{Y,Z}} 
& X \otimes (Z \otimes Y) \ar[r]^-{\xi_{X,Z\otimes Y}} 
& (Z \otimes Y) \otimes X \ar[u]_-{\alpha}}
\end{equation}
is commutative for all objects $X,Y,Z \in \M$.
\end{description}
A coboundary monoidal category is \emph{strict}\index{strict coboundary monoidal category} if the underlying monoidal category is strict.
\end{definition}

\begin{interpretation}
Regarding each object as an interval, one can think of the commutor $\xi_{X,Y}$ as permuting two consecutive intervals $X$ and $Y$.  Ignoring the associativity isomorphism for the moment, the cactus axiom \eqref{cactus-axiom} says that, to reverse the order of three consecutive intervals $X$, $Y$, and $Z$, one can perform either one of the following two sequences of moves: 
\begin{enumerate}[label=(\roman*)]
\item Starting at the back, first permute the intervals $Y$ and $Z$, and then permute $X$ and the combined interval $ZY$.
\item Starting on the left, first permute the intervals $X$ and $Y$, and then permute the combined interval $YX$ and $Z$.
\end{enumerate}
With the associativity isomorphisms taken into account, the cactus axiom says that moving the parentheses using $\alpha$ and then performing the moves (i) are the same as performing the moves (ii) and then moving the parentheses using $\alpha^{-1}$.\dqed
\end{interpretation}

To define the cactus group action in a coboundary monoidal category in the next section, we need a few technical results, for which we first introduce some notations.

\begin{definition}\label{def:normalized-monoidal-product}
Suppose $A$ is a monoidal category, and $a_1,\ldots,a_n$ are objects in $A$.
\begin{enumerate}
\item Define inductively the object\label{not:left-normalized}\index{normalized monoidal product} \[\bigl(a_1\otimes \cdots \otimes a_n\bigr)_l = \begin{cases} a_1 & \text{if $n=1$},\\
\Bigl(\bigl(a_1\otimes \cdots \otimes a_{n-1}\bigr)_l \Bigr)\otimes a_n & \text{if $n \geq 2$},\end{cases}\] called a \emph{left normalized monoidal product}.
\item Define inductively the object \[\bigl(a_1\otimes \cdots \otimes a_n\bigr)_r= \begin{cases} a_1 & \text{if $n=1$},\\
a_1 \otimes \Bigl(\bigl(a_2\otimes \cdots \otimes a_{n}\bigr)_r\Bigr)& \text{if $n \geq 2$},\end{cases}\] called a \emph{right normalized monoidal product}.
\end{enumerate}
Furthermore, we apply the left/right normalized terminology to variables in a standard non-associative monomial in the obvious way.
\end{definition}

\begin{example}
We have the left normalized monoidal products
\[\begin{split}
(a_1 \otimes a_2 \otimes a_3)_l &= \bigl(a_1\otimes a_2\bigr) \otimes a_3,\\
(a_1 \otimes a_2 \otimes a_3 \otimes a_4)_l &= \Bigl(\bigl(a_1\otimes a_2\bigr) \otimes a_3\Bigr) \otimes a_4,\end{split}\]
in which every pair of parentheses starts in the front.  Similarly, in a right normalized monoidal product, such as \[(a_1 \otimes a_2 \otimes a_3 \otimes a_4)_r = a_1\otimes \Bigl(a_2 \otimes \bigl(a_3\otimes a_4\bigr)\Bigr),\] the right half of every pair of parentheses is at the end.  Also note that \[(a_1\otimes a_2)_l = a_1 \otimes a_2 = (a_1\otimes a_2)_r.\]  In the standard non-associative monomial without $x_0$'s, \[\bigl((x_4x_1)x_2\bigr)\bigl(x_5(x_3x_6)\bigr),\] the variables $x_4, x_1$, and $x_2$ are in a left normalized form, and the variables $x_5,x_3$, and $x_6$ are in a right normalized form.\dqed
\end{example}

\begin{convention}
In the rest of this chapter, to save space we sometimes omit writing the symbol $\otimes$ for the monoidal product, so $a_1a_2$ means $a_1 \otimes a_2$. \dqed
\end{convention}

\begin{definition}\label{def:coboundary-cat-cpq}
Suppose $(A,\xi)$ is a coboundary monoidal category, and $k \geq 2$.   We define the\index{interval-reversing isomorphism} \emph{interval-reversing isomorphism} $c^{(k)}$ in $A$ as the composite
\[\nicexy@C-1cm{\bigl(a_1 a_{2}\cdots a_{k-1} a_k\bigr)_r \ar[d]_-{\Id^{k-2} \xi_{a_{k-1},a_k}} \ar[rr]^-{c^{(k)}_{a_1,\ldots,a_k}}_-{\cong} && \bigl(a_k  a_{k-1}  \cdots  a_{2}a_1\bigr)_l\\
\bigl(a_1 \cdots a_{k-2} (a_ka_{k-1})_l\bigr)_r \ar[d]_-{\Id^{k-3}\xi_{a_{k-2}, (a_ka_{k-1})_l}} && \bigl(a_1(a_k a_{k-1}\cdots a_{2})_l\bigr)_r \ar[u]_-{\xi_{a_1,(a_k a_{k-1}\cdots a_{2})_l}}\\
\bigl(a_1\cdots a_{k-3} (a_k a_{k-1}a_{k-2})_l\bigr)_r \ar`d/3pt[dr]_(.7){\Id^{k-4} \xi_{a_{k-3}, (a_k a_{k-1}a_{k-2})_l}} [dr]
&& \bigl(a_1 a_{2} (a_k a_{k-1}\cdots a_{3})_l\bigr)_r \ar[u]_-{\Id_{a_1} \xi_{a_{2}, (a_k a_{k-1}\cdots a_{3})_l}}\\
& \cdots \ar`r/3pt[ur][ur]_(.3){\Id^2 \xi_{a_{3}, (a_k a_{k-1}\cdots a_{4})_l}} &}\]
for objects $a_1,\ldots,a_k\in A$.  Algebraically, $c^{(k)}$ is the composite
\begin{equation}\label{ck-algebraic}
c^{(k)}_{a_1,\ldots,a_k} = \bigcircle_{j=2}^k \,\Bigl[\Id^{k-j} \otimes \xi_{a_{k-j+1},\, (a_k\cdots a_{k-j+2})_l} \Bigr]
\end{equation}
with the big circle denoting categorical composition and with the constituents written from right to left as the index increases.  We also define $c^{(1)}= \Id$.  
\end{definition}

Note that each $c^{(k)}$ is a natural isomorphism.

\begin{example}
The interval-reversing isomorphism $c^{(5)}$ is the composite
\[\nicexy@C-1.7cm{a_1\Bigl(a_2\bigl(a_3(a_4a_5)\bigr)\Bigr) \ar[rr]^-{c^{(5)}} \ar[d]_-{\Id_{a_1}\Id_{a_2}\Id_{a_3} \xi_{a_4,a_5}} && \Bigl(\bigl((a_5a_4)a_3\bigr)a_2\Bigr)a_1\\
a_1\Bigl(a_2\bigl(a_3(a_5a_4)\bigr)\Bigr) \ar`d/3pt[dr]_(.7){\Id_{a_1}\Id_{a_2} \xi_{a_3, a_5a_4}} [dr] && a_1\Bigl(\bigl((a_5a_4)a_3\bigr)a_2\Bigr) \ar[u]_-{\xi_{a_1,((a_5a_4)a_3)a_2}}\\
& a_1\Bigl(a_2\bigl((a_5a_4)a_3\bigr)\Bigr) \ar`r/3pt[ur][ur]_(.3){\Id_{a_1}\xi_{a_2,(a_5a_4)a_3}} &}\]
of four isomorphisms.\dqed
\end{example}

\begin{interpretation}
The interval-reversing isomorphism $c^{(k)}$ is a categorical analogue of the generating cactus $s^{(k)}_{1,k}$ in the cactus group $Cac_k$ in Definition \ref{def:cactus-group}.  Half of Theorem \ref{cacmonoidal-coboundary-monoidal} below involves showing that each coboundary monoidal category is a $\Cac$-monoidal category.  The action of each cactus on iterated monoidal products will be defined using these interval-reversing isomorphisms.  In the rest of this section, we show that the interval-reversing isomorphisms satisfy relations similar to the generating relations in the cactus groups.\dqed
\end{interpretation}

\begin{convention}
To simplify the presentation, for the rest of this section only, we consider \emph{strict} coboundary monoidal categories, in which the associativity isomorphism $\alpha$, the left unit isomorphism $\lambda$, and the right unit isomorphism $\rho$ are identity morphisms.  With this assumption, we can omit parentheses for multiple monoidal products, except when we want to clarify some expressions.  This is purely for notational convenience, since we can incorporate $\alpha$, $\lambda$, and $\rho$ using the coboundary monoidal category axioms.\dqed
\end{convention}

\begin{notation}\label{not:coboundary-monoidal}
For the rest of this section, suppose $a_1,a_2,\ldots$ are objects in a \emph{strict} coboundary monoidal category $(A,\xi)$.  For $1\leq i, j$, we define the objects
\[\begin{split}
a_{[i,j]} &= \begin{cases} a_ia_{i+1}\cdots a_j &\text{if $i\leq j$},\\ \tensorunit & \text{if $i>j$},\end{cases}\\
a_{[j,i]\downarrow} &= \begin{cases} a_ja_{j-1} \cdots a_i & \text{if $j \geq i$},\\ \tensorunit & \text{if $j<i$}.\end{cases}
\end{split}\]
In $a_{[j,i]\downarrow}$, from left to right, the indices decrease from $j$ down to $i$ if $j \geq i$.  For example, we have \[a_{[4,2]\downarrow} = a_4a_3a_2.\]  We also extend the above notations to subscripts of the interval-reversing isomorphisms.  For example, for $i \leq j$, we have
\[\begin{split}
c^{(j-i+1)}_{a_{[i,j]}} &= c^{(j-i+1)}_{a_i,\ldots,a_j} : a_{[i,j]} \iso a_{[j,i]\downarrow},\\
c^{(j-i+1)}_{a_{[j,i]\downarrow}} &= c^{(j-i+1)}_{a_j,a_{j-1},\ldots,a_i} : a_{[j,i]\downarrow} \iso a_{[i,j]}.
\end{split}\]  
If we want to use the objects $a_{[i,j]}$ or $a_{[j,i]\downarrow}$ as a single subscript in $c^{(k)}$, then we put a pair of parentheses around it.  For example, we have
\[c^{(5)}_{a_1,(a_{[4,2]\downarrow}),a_{[5,7]}} = c^{(5)}_{a_1,(a_4a_3a_2),a_5,a_6,a_7} : a_1(a_4a_3a_2)a_5a_6a_7 \iso a_7a_6a_5(a_4a_3a_2)a_1,\] in which the second subscript is the object $a_{[4,2]\downarrow} = a_4a_3a_2$.\dqed
\end{notation}

Our first observation is that the interval-reversing isomorphism $c^{(k)}$ can be decomposed as $c^{(k-i)}$ at the back followed by $c^{(i+1)}$.

\begin{lemma}\label{ck-altdef}
For $1 \leq i \leq k-2$, the diagram
\[\nicexy@C+1cm{a_{[1,k]} \ar@{=}[d] \ar[rr]^-{c^{(k)}_{a_{[1,k]}}} && a_{[k,1]\downarrow}\\
a_{[1,i]} a_{[i+1,k]} \ar[r]^-{\Id^i \otimes c^{(k-i)}_{a_{[i+1,k]}}} & 
a_{[1,i]} a_{[k,i+1]\downarrow} \ar[r]^-{c^{(i+1)}_{a_{[1,i]}, (a_{[k,i+1]\downarrow})}} & a_{[k,i+1]\downarrow}a_{[i,1]\downarrow} \ar@{=}[u]}\]
is commutative.
\end{lemma} 

\begin{proof}
We first split the algebraic definition of $c^{(k)}$ in \eqref{ck-algebraic} as the composite
\[\biggl[\bigcircle_{j=k-i+1}^k \,\Id^{k-j} \otimes \xi_{a_{k-j+1},\, (a_{[k,k-j+2]\downarrow})}\biggr] \circ \biggl[\bigcircle_{j=2}^{k-i}\, \Id^{k-j} \otimes \xi_{a_{k-j+1},\, (a_{[k,k-j+2]\downarrow})}\biggr].\]  In the right factor above, we factor out $\Id^i \otimes -$ using the naturality of the monoidal product.  In the left factor above, we shift the index $j$ by $k-i-1$.  Then the above decomposition for $c^{(k)}$ becomes
\[\begin{split}
& \biggl[\bigcircle_{j=2}^{i+1} \,\Id^{i+1-j} \otimes \xi_{a_{i+2-j},\, (a_{[k,i+3-j]\downarrow})}\biggr] \circ
\biggl[\Id^i \otimes \Bigl[\bigcircle_{j=2}^{k-i} \,\Id^{k-i-j} \otimes \xi_{a_{k-j+1},\, (a_{[k,k-j+2]\downarrow})}\Bigr] \biggr]\\
&= c^{(i+1)}_{a_{[1,i]}, (a_{[k,i+1]\downarrow})} \circ \bigl(\Id^i \otimes c^{(k-i)}_{a_{[i+1,k]}}\bigr),
\end{split}\]
where in the last equality we used \eqref{ck-algebraic} twice.
\end{proof}

\begin{example}\label{ck-ck-1-xi}
The $i=1$ case of Lemma \ref{ck-altdef} is the commutative diagram
\[\nicexy@C+1cm{a_{[1,k]} \ar[rr]^-{c^{(k)}_{a_{[1,k]}}} && a_{[k,1]\downarrow}\\
a_{1}a_{[2,k]} \ar@{=}[u] \ar[r]^-{\Id_{a_1} \otimes c^{(k-1)}_{a_{[2,k]}}} & a_{1}a_{[k,2]\downarrow} \ar[r]^-{\xi_{a_1,(a_{[k,2]\downarrow})}} & a_{[k,2]\downarrow} a_1 \ar@{=}[u]}\] that expresses $c^{(k)}_{a_1,\ldots,a_k}$ as the interval-reversing isomorphism $c^{(k-1)}_{a_2,\ldots,a_k}$ at the back followed by the commutor $\xi_{a_1,(a_k\cdots a_2)}$.\dqed
\end{example}

\begin{example}\label{ck-xi-ck-1}
The $1\leq i=k-2$ case of Lemma \ref{ck-altdef} is the commutative diagram
\[\nicexy@C+1cm{a_{[1,k]} \ar[rr]^-{c^{(k)}_{a_{[1,k]}}} && a_{[k,1]\downarrow}\\
a_{[1,k-2]}a_{k-1}a_k \ar@{=}[u] \ar[r]^-{\Id^{k-2} \otimes \xi_{a_{k-1},a_k}} & a_{[1,k-2]}(a_ka_{k-1}) \ar[r]^-{c^{(k-1)}_{a_{[1,k-2]},(a_ka_{k-1})}} & a_ka_{k-1} a_{[k-2,1]\downarrow} \ar@{=}[u]}\] that expresses $c^{(k)}_{a_1,\ldots,a_k}$ as the commutor $\xi_{a_{k-1},a_k}$ at the back followed by the interval-reversing isomorphism $c^{(k-1)}_{a_1,\ldots,a_{k-2},(a_ka_{k-1})}$.\dqed
\end{example}

The following observation says that the interval-reversing isomorphism $c^{(k)}$ can also be expressed as the commutor $\xi$ at the \emph{front} followed by $c^{(k-1)}$.  This observation should be compared with Example \ref{ck-xi-ck-1}.

\begin{lemma}\label{ck-xi-infront}
For $k\geq 3$, the diagram
\[\nicexy@C+1cm{a_{[1,k]} \ar[rr]^-{c^{(k)}_{a_{[1,k]}}} && a_{[k,1]\downarrow}\\
(a_1a_2)a_{[3,k]} \ar@{=}[u] \ar[r]^-{\xi_{a_1,a_2} \otimes \Id^{k-2}} & (a_2a_1)a_{[3,k]} \ar[r]^-{c^{(k-1)}_{(a_2a_1),a_{[3,k]}}} & a_{[k,3]\downarrow} (a_2a_1) \ar@{=}[u]}\] 
is commutative.
\end{lemma}

\begin{proof}
We proceed by induction on $k$.  The first case $k=3$ is the cactus axiom \eqref{cactus-axiom}, since (i) $c^{(2)} = \xi$ and (ii) by definition $c^{(3)}_{a_1,a_2,a_3}$ is the composite
\[\nicexy@R-.3cm{a_1(a_2a_3) \ar`d/3pt[dr]_(.7){\Id_{a_1}\otimes \xi_{a_2,a_3}} [dr] \ar[rr]^-{c^{(3)}} && a_3a_2a_1\\
& a_1(a_3a_2) \ar`r/3pt[ur][ur]_(.3){\xi_{a_1,(a_3a_2)}}}.\]

For the induction step with $k\geq 4$, we must show that the outermost diagram in the diagram
\[\nicexy@C-.5cm@R+.5cm{a_{[1,k]} \ar`d/4pt[dddr]_-{\xi_{a_1,a_2} \otimes \Id^{k-2}} [dddr] \ar[dr]|-{\Id^{k-2} \xi_{a_{k-1},a_k}} \ar[rr]^-{c^{(k)}_{a_{[1,k]}}} && a_{[k,1]\downarrow}\\
& a_{[1,k-2]}(a_ka_{k-1}) \ar[d]_-{\xi_{a_1,a_2}\Id^{k-2}} \ar[ur]|-{c^{(k-1)}} &\\
& (a_2a_1)a_{[3,k-2]}(a_ka_{k-1}) \ar@(r,dl)[uur]|-{c^{(k-2)}} &\\
& (a_2a_1)a_{[3,k]} \ar[u]^-{\Id^{k-2}\xi_{a_{k-1},a_k}} \ar`r/4pt[uuur][uuur]_-{c^{(k-1)}_{(a_2a_1),a_{[3,k]}}} &}\] 
is commutative, where
\[c^{(k-1)} = c^{(k-1)}_{a_{[1,k-2]},(a_ka_{k-1})} \andspace 
c^{(k-2)} = c^{(k-2)}_{(a_2a_1), a_{[3,k-2]}, (a_ka_{k-1})}.\]
In the above diagram:
\begin{itemize}
\item The left trapezoid is commutative by the naturality of the monoidal product.
\item The right trapezoid and the top triangle are commutative by the decomposition of $c^{(k)}$ in Example \ref{ck-xi-ck-1}.
\item The remaining middle-right sub-diagram is commutative by the induction hypothesis.
\end{itemize}
This finishes the induction.
\end{proof}

The next observation says that the interval-reversing isomorphism $c^{(k)}$ can be expressed as $c^{(k-1)}$ at the \emph{front} followed by the commutor $\xi$.  This should be compared with Example \ref{ck-ck-1-xi}.

\begin{lemma}\label{ck-1infront}
For $k \geq 3$, the diagram
\[\nicexy@C+1cm{a_{[1,k]} \ar[rr]^-{c^{(k)}_{a_{[1,k]}}} && a_{[k,1]\downarrow}\\
a_{[1,k-1]}a_k \ar@{=}[u] \ar[r]^-{c^{(k-1)}_{a_{[1,k-1]}} \otimes \Id_{a_k}} & a_{[k-1,1]\downarrow} a_k \ar[r]^-{\xi_{(a_{[k-1,1]\downarrow}), a_k}} & a_ka_{[k-1,1]\downarrow} \ar@{=}[u]}\] 
is commutative.
\end{lemma}

\begin{proof}
We proceed by induction on $k$.  The initial case $k=3$ is true by the cactus axiom \eqref{cactus-axiom}, the definition of $c^{(3)}$, and the fact that $c^{(2)} = \xi$.  

For the induction step, suppose $k\geq 4$.  Using the decomposition of $c^{(k)}$ in Example \ref{ck-ck-1-xi}, we must show that the outermost diagram in the diagram
\[\nicexy@C-.7cm{a_{[1,k]} \ar[rrr]^-{c^{(k-1)}_{a_{[1,k-1]}} \otimes \Id_{a_k}} \ar[drr]|-{\xi_{a_1,a_2}\Id^{k-2}} \ar[ddr]|-{\Id^2 c^{(k-2)}_{a_{[3,k]}}} \ar[dddd]_-{\Id_{a_1}\otimes c^{(k-1)}_{a_{[2,k]}}} &&& (a_{[k-1,1]\downarrow}) a_k \ar[dd]^-{\xi_{(a_{[k-1,1]\downarrow}), a_k}}\\
&& (a_2a_1)a_{[3,k]} \ar[dd]|-{\Id^2 c^{(k-2)}_{a_{[3,k]}}} \ar[ur]|-{c^{(k-2)} \Id_{a_k}} \ar[dr]|-{c^{(k-1)}} &\\
& (a_1a_2)a_{[k,3]\downarrow} \ar[ddl]|-{\Id_{a_1}\xi_{a_2,a_{[k,3]\downarrow}}} \ar[dr]|-{\xi_{a_1,a_2} \Id^{k-2}} && a_{[k,1]\downarrow} \ar@{=}[dd]\\
&& (a_2a_1)a_{[k,3]\downarrow} \ar[dr]|-{\xi_{(a_2a_1),(a_{[k,3]\downarrow})}} &\\
a_1a_{[k,2]\downarrow} \ar[rrr]^-{\xi_{a_1,(a_{[k,2]\downarrow})}} &&& a_{[k,1]\downarrow}}\]
is commutative, in which \[c^{(k-2)} = c^{(k-2)}_{(a_2a_1),a_{[3,k-1]}} \andspace c^{(k-1)}=c^{(k-1)}_{(a_2a_1),a_{[3,k]}}\] in the upper-right triangle.
\begin{itemize}
\item The left triangle and the lower-right parallelogram are commutative by the decomposition of $c^{(k)}$ in Example \ref{ck-ck-1-xi}, applied to $c^{(k-1)}$ here.
\item The bottom triangle is commutative by the cactus axiom \eqref{cactus-axiom}.
\item The upper-right triangle is commutative by the induction hypothesis.
\item The top triangle is commutative by Lemma \ref{ck-xi-infront}, applied to $c^{(k-1)}$ here.
\item The middle quadrilateral is commutative by the naturality of the monoidal product.
\end{itemize}
This finishes the induction.
\end{proof}

The following observation says that the interval-reversing isomorphism $c^{(k)}$ is an involution.  It is a categorical analogue of the involution axiom in the cactus group in Definition \ref{def:cactus-group}.

\begin{proposition}\label{ck-involution}
For $k\geq 2$, the composite
\[\nicexy@C+.3cm{a_{[1,k]} \ar[r]^-{c^{(k)}_{a_{[1,k]}}} & a_{[k,1]\downarrow} \ar[r]^-{c^{(k)}_{a_{[k,1]\downarrow}}} & a_{[1,k]}}\]
is the identity morphism.
\end{proposition}

\begin{proof}
We proceed by induction on $k$.  The initial case $k=2$ is true by the symmetry axiom \eqref{coboundary-symmetry-axiom} and that $c^{(2)} = \xi$.

For the induction step, suppose $k \geq 3$.  Using the symmetry axiom \eqref{coboundary-symmetry-axiom}, it is enough to show that the outermost diagram in the diagram
\[\nicexy@C+.3cm@R+.4cm{a_{[1,k]} \ar[dd]_-{\xi_{a_1,(a_{[2,k]})}} \ar[dr]|-{\Id_{a_1}c^{(k-1)}_{a_{[2,k]}}} \ar[rr]^-{c^{(k)}_{a_{[1,k]}}} && a_{[k,1]\downarrow}  \ar[r]^-{c^{(k)}_{a_{[k,1]\downarrow}}} & a_{[1,k]}\\
& a_1a_{[k,2]\downarrow} \ar[r]^-{\xi_{a_1,(a_{[k,2]\downarrow})}} & a_{[k,2]\downarrow} a_1 \ar@{=}[u] \ar[dr]|-{c^{(k-1)}_{a_{[k,2]\downarrow}} \Id_{a_1}} &\\
a_{[2,k]}a_1 \ar[urr]|-{c^{(k-1)}_{a_{[2,k]}} \Id_{a_1}}
\ar@{=}[rrr] &&& a_{[2,k]}a_1 \ar[uu]_-{\xi_{(a_{[2,k]}), a_1}}}\]
is commutative.
\begin{itemize}
\item The top left trapezoid is commutative by Example \ref{ck-ck-1-xi}.
\item The right trapezoid is commutative by Lemma \ref{ck-1infront}.
\item The bottom triangle is commutative by the induction hypothesis.
\item The remaining quadrilateral on the left is commutative by the naturality of the commutor $\xi$.
\end{itemize}
This finishes the induction.
\end{proof}

We will also need a categorical analogue of the containment axiom in the cactus group, for which we first need the next two preliminary observations.  The next observation is another redefinition of $c^{(k)}$.

\begin{lemma}\label{cikminusm}
For $1 \leq i < k$, the diagram
\[\nicexy@C-1.2cm{a_{[1,i]} a_{[i+1,k]} 
\ar[rr]^-{c^{(k)}_{a_{[1,k]}}}
\ar`d/4pt[dr]_(.8){c^{(i+1)}_{a_{[1,i]}, (a_{[i+1,k]})}}[dr] 
&& a_{[k,i+1]\downarrow} a_{[i,1]\downarrow}\\
& a_{[i+1,k]} a_{[i,1]\downarrow} \ar`r/4pt[ur][ur]_(.2){c^{(k-i)}_{a_{[i+1,k]}} \otimes \Id^{i}} &}\]
is commutative.
\end{lemma}

\begin{proof}
The diagram in the statement is the outermost diagram below.
\[\nicexy@C-1.2cm{
a_{[1,i]} a_{[i+1,k]} \ar`d/4pt[ddr]_(.75){c^{(i+1)}} [ddr] 
\ar[dr]|-{\Id^{i} c^{(k-i)}} \ar[rr]^-{c^{(k)}}
&& a_{[k,i+1]\downarrow} a_{[i,1]\downarrow} \\
& a_{[1,i]} \bigl(a_{[k,i+1]\downarrow}\bigr) \ar[ur]|-{c^{(i+1)}} &\\ 
& a_{[i+1,k]} a_{[i,1]\downarrow} \ar`r/4pt[uur][uur]_(.25){c^{(k-i)} \Id^{i}} &}\]
In the diagram above, the triangle is commutative by Lemma \ref{ck-altdef}.  The remaining sub-diagram is commutative by the naturality of $c^{(i+1)}$.
\end{proof}

\begin{lemma}\label{cimone}
For $1\leq i < j$, the diagram
\[\nicexy@C+1.5cm@R+.7cm{a_{[1,j+1]} = a_{[1,j]} a_{j+1} 
\ar[d]_-{\xi_{(a_{[1,j]}), a_{j+1}}} \ar[r]^-{c^{(j+1)}_{a_{[1,j+1]}}} 
& a_{j+1} a_{[j,1]\downarrow} = a_{[j+1,1]\downarrow}\\
a_{j+1} a_{[1,i]} a_{[i+1,j]} \ar[r]^-{\Id \otimes c^{(i+1)}_{a_{[1,i]}, (a_{[i+1,j]})}} 
& a_{j+1} a_{[i+1,j]} a_{[i,1]\downarrow}  \ar[u]_-{\Id \otimes c^{(j-i)}_{a_{[i+1,j]}} \otimes \Id^i}}\]
is commutative.
\end{lemma}

\begin{proof}
The diagram in the statement is the outermost diagram below.
\[\nicexy@C-1cm@R+.2cm{a_{[1,j]} a_{j+1} \ar[dd]_-{\xi} \ar[dr]|-{c^{(j)}_{a_{[1,j]}} \Id} \ar[rr]^-{c^{(j+1)}} 
&& a_{j+1} a_{[j,1]\downarrow} \\
& \bigl(a_{[j,1]\downarrow}\bigr) (a_{j+1}) \ar[dr]|-{\xi} &\\
a_{j+1} a_{[1,i]} a_{[i+1,j]} \ar`d/4pt[dr]_(.7){\Id\otimes c^{(i+1)}} [dr] 
\ar[rr]^-{\Id \otimes c^{(j)}}
&& a_{j+1} a_{[j,1]\downarrow} \ar@{=}[uu]\\
& a_{j+1} a_{[i+1,j]} a_{[i,1]\downarrow} \ar`r/4pt[ur][ur]_(.3){\Id \otimes c^{(j-i)} \otimes \Id^i} &}\]
In the diagram above:
\begin{itemize}
\item The subscripts for the commutor $\xi$ in the interior are explained by the parentheses in its domain.
\item The upper-right triangle is commutative by Lemma \ref{ck-1infront}.
\item The upper-left triangle is commutative by the naturality of the commutor $\xi$.
\item The remaining bottom sub-diagram is commutative by Lemma \ref{cikminusm}.
\end{itemize}
This finishes the proof.
\end{proof}

The next observation is a categorical analogue of the containment axiom in the cactus group in Definition \ref{def:cactus-group}.

\begin{proposition}\label{ck-containment}
For $0 \leq i \leq j-2 \leq k-2$, the diagram
\[\nicexy@C+2cm@R+.7cm{a_{[1,k]} \ar[d]|-{\Id^i \otimes c^{(j-i)}_{a_{[i+1,j]}} \otimes \Id^{k-j}} \ar[r]^-{c^{(k)}_{a_{[1,k]}}} & a_{[k,1]\downarrow} \ar[d]|-{\Id^{k-j} \otimes c^{(j-i)}_{a_{[j,i+1]\downarrow}} \otimes \Id^{i}}\\
a_{[1,i]} a_{[j,i+1]\downarrow} a_{[j+1,k]} \ar[r]^-{c^{(k)}_{a_{[1,i]}, a_{[j,i+1]\downarrow}, a_{[j+1,k]}}} & a_{[k,j+1]\downarrow} a_{[i+1,j]} a_{[i,1]\downarrow}}\]
is commutative.
\end{proposition}

\begin{proof}
If $j=k$, then the stated diagram is the outermost diagram below.
\[\nicexy@R+.7cm{a_{[1,k]} \ar[d]_-{\Id^i \otimes c^{(k-i)}_{a_{[i+1,k]}}} \ar[r]^-{c^{(k)}_{a_{[1,k]}}} & a_{[k,1]\downarrow} \ar@{=}[r] & a_{[k,i+1]\downarrow} a_{[i,1]\downarrow} \ar[d]^-{c^{(k-i)}_{a_{[k,i+1]\downarrow}} \otimes \Id^{i}}\\
a_{[1,i]} a_{[k,i+1]\downarrow} \ar[urr]|-{c^{(i+1)}_{a_{[1,i]}, (a_{[k,i+1]\downarrow})}}
\ar[rr]_-{c^{(k)}_{a_{[1,i]}, a_{[k,i+1]\downarrow}}} && a_{[i+1,k]} a_{[i,1]\downarrow}}\]
The upper-left triangle is commutative by Lemma \ref{ck-altdef}.  The lower-right triangle is commutative by Lemma \ref{cikminusm}.

Next suppose $j < k$.  Then the stated diagram is the outermost diagram below.
\[\begin{small}\nicexy@C-1.2cm{a_{[1,k]} \ar[dddddd]|-{\Id^i  c^{(j-i)} \Id^{k-j}} \ar[dddr]|-{c^{(j)}\Id^{k-j}} \ar[dr]|-{\Id^j c^{(k-j)}} \ar[rrr]^-{c^{(k)}} &&& a_{[k,1]\downarrow} \ar[dddddd]|-{\Id^{k-j} c^{(j-i)} \Id^{i}}\\
& a_{[1,j]} a_{[k,j+1]\downarrow} \ar[dr]|-{c^{(j)}\Id^{k-j}} \ar[urr]|-{c^{(j+1)}}  && \\
&& a_{[j,1]\downarrow} a_{[k,j+1]\downarrow} \ar[uur]|-{\xi} & \\
& a_{[j,1]\downarrow} a_{[j+1,k]} \ar[ur]|-{\Id^j c^{(k-j)}} &&\\
&& a_{[k,j+1]\downarrow} a_{[1,i]} a_{[j,i+1]\downarrow} \ar[uuuur]|(.2){\Id^{k-j} c^{(i+1)}} &\\
& a_{[1,i]} a_{[j,i+1]\downarrow} a_{[k,j+1]\downarrow} \ar[uuur]|-{c^{(i+1)} \Id^{k-j}} \ar[ur]|-{\xi} \ar[drr]|-{c^{(j+1)}} &&\\
a_{[1,i]} a_{[j,i+1]\downarrow} a_{[j+1,k]} \ar[uuur]|-{c^{(i+1)} \Id^{k-j}} \ar[ur]|-{\Id^j c^{(k-j)}} \ar[rrr]^-{c^{(k)}} &&& a_{[k,j+1]\downarrow} a_{[i+1,j]} a_{[i,1]\downarrow}}\end{small}\]
In the above diagram:
\begin{itemize}
\item In the interior, the subscripts of each instance of the commutor $\xi$ are explained by its domain and codomain.  The interval-reversing isomorphisms are given by
\[\begin{split}
c^{(j)} &= c^{(j)}_{a_{[1,j]}},\\
c^{(j+1)} &= c^{(j+1)}_{a_{[1,j]}, (a_{[k,j+1]\downarrow})} \orspace 
c^{(j+1)}_{a_{[1,i]}, a_{[j,i+1]\downarrow}, (a_{[k,j+1]\downarrow})},\\
c^{(i+1)} &= c^{(i+1)}_{a_{[1,i]}, (a_{[j,i+1]\downarrow})},\\
c^{(k-j)} &= c^{(k-j)}_{a_{[j+1,k]}}.
\end{split}\] 
\item The top, left, and bottom triangles are commutative by Lemma \ref{ck-altdef}.
\item The right quadrilateral is commutative by Lemma \ref{cimone}.
\item Among the remaining four sub-diagrams in the interior, the triangle and the parallelogram on the left are commutative by the naturality of the monoidal product.
\item The remaining upper-right triangle is commutative by Lemma \ref{ck-1infront}.
\item The remaining parallelogram on the right is commutative by the naturality of the commutor $\xi$.
\end{itemize}
This finishes the proof.
\end{proof}

\begin{remark}
Proposition \ref{ck-involution} and Proposition \ref{ck-containment} are essentially Lemma 3(iv) and Lemma 4 in \cite{hen-kam}, which are the key results for identifying the cactus groups as the groups that act naturally on multiple monoidal products in coboundary monoidal categories.  However, our presentation in this section is more elementary and has more detail, which we hope is useful to the reader.\dqed
\end{remark}

\section{Coboundary Monoidal Categories are $\Cac$-Monoidal Categories}
\label{sec:coboundary-monoidal-cacmonoidal}

The purpose of this section is to show that for the cactus group operad $\Cac$ in Definition \ref{def:cactus-group-operad}, $\Cac$-monoidal categories as in Definition \ref{def:g-monoidal-category}\label{not:mcatcac} (i.e., $\MCatcac$-algebras) are exactly small coboundary monoidal categories.  The strict version is also true.  In \cite{hen-kam} (Theorem 7) Henriques and Kamnitzer proved that the cactus group $Cac_n$ acts on $n$-fold monoidal products in a coboundary monoidal category, and that these natural actions are exactly the ones generated by the commutor.  Theorem \ref{cacmonoidal-coboundary-monoidal} below provides another proof of \cite{hen-kam} (Theorem 7) and also provides an operadic interpretation.

\begin{theorem}\label{cacmonoidal-coboundary-monoidal}\index{coboundary monoidal category!operadic interpretation}
For each small category, the following two structures are equivalent:
\begin{enumerate}[label=(\roman*)]
\item A $\Cac$-monoidal category structure for the cactus group operad $\Cac$.
\item A coboundary monoidal category structure.
\end{enumerate}
\end{theorem}

\begin{proof}
This proof follows the same outline as the proof of Theorem \ref{bmonoidal-braided-monoidal}, but the details are different.  Suppose $(A,\xi)$ is a small coboundary monoidal category.  We must show that this yields an $\MCatcac$-algebra structure on $A$.  Each isomorphism $\sigma \in \MCatcac(n)(p,q)$ is an element in the $n$th cactus group $\Cac_n$ in Definition \ref{def:cactus-group}, such that $\sigmabar = \tau_{\qbar}\tau_{\pbar}^{-1}$ as in Definition \ref{def:mcatg-underlying-categories}.  There exists a decomposition 
\begin{equation}\label{sigma-cac-decomp}
\sigma = s^{(n)}_{p_k,q_k} \cdots s^{(n)}_{p_1,q_1} \in Cac_n
\end{equation}
with each $s^{(n)}_{p_j,q_j} \in Cac_n$ a generating cactus and $1\leq p_j < q_j \leq n$.  We call an isomorphism in $A$ \emph{canonical} if it is a categorical composite of iterated monoidal products of identity morphisms, $\alpha$, $\lambda$, $\rho$, $\xi$, and their inverses.  

We define the $\Cac$-monoidal structure isomorphism $\sigma_*$ in $A$ as the composite 
\begin{equation}\label{sigmastar-coboundary-monoidal}
\nicexy@C+.3cm{p\bigl(\tensorunit,\ua\bigr) \ar[rrr]^-{\sigma_*}_-{\cong} \ar[d]_-{u_p} &&& q\bigl(\tensorunit,\ua\bigr)\\
p_{\st}(\ua) \ar[d]_{\alpha_1} &&& q_{\st}(\ua) \ar[u]_-{u_q^{-1}}\\
p_{\st}^1(\ua) \ar[d]_-{c_{p_1,q_1}^{(n)}} &&& p_{\st}^{2k}(\ua) \ar[u]_-{\alpha_{k+1}} \\
p_{\st}^2(\ua) \ar[r]^-{\alpha_2} & \cdots \ar[r]^-{c_{p_{k-1},q_{k-1}}^{(n)}}& p_{\st}^{2k-2}(\ua) \ar[r]^-{\alpha_k} & p_{\st}^{2k-1}(\ua) \ar[u]_-{c_{p_k,q_k}^{(n)}}}
\end{equation}
of canonical isomorphisms defined as follows.
\begin{itemize}
\item $p_{\st}$ and $q_{\st}$ are obtained from $p$ and $q$, respectively, by removing all instances of $x_0$'s.
\item $u_p$ and $u_q$ are the unique canonical isomorphisms involving only $\Id$, $\lambda$, and $\rho$ that are guaranteed by the Coherence Theorem for monoidal categories \cite{maclane} (VII.2 Corollary).  Each of these two morphisms removes all copies of the monoidal unit $\tensorunit$ corresponding to the $x_0$'s.
\item $p_{\st}^1$ is obtained from $p_{\st}$ by rearranging the parentheses such that the $p_1$th to the $q_1$th variables, counting from the left, are arranged in right normalized form.  
\item The isomorphism $\alpha_1$ is the unique canonical isomorphism involving only $\Id$ and $\alpha^{\pm 1}$ guaranteed by the Coherence Theorem for monoidal categories. 
\item The canonical isomorphism $c_{p_1,q_1}^{(n)}$ is the monoidal product \[c_{p_1,q_1}^{(n)} = \Id^{\otimes p_1-1} \otimes c^{(q_1-p_1+1)} \otimes \Id^{\otimes n-q_1}\] such that:
\begin{itemize}
\item $c^{(q_1-p_1+1)}$ is the interval-reversing isomorphism in Definition \ref{def:coboundary-cat-cpq}, applied to the $p_1$th to the $q_1$th tensor factors in $p_{\st}^1(\ua)$, counting from the left.
\item $c_{p_1,q_1}^{(n)}$ has the same parenthesization as $p_{\st}^1(\ua)$.
\item $p_{\st}^2$ is obtained from $p_{\st}^1$ by:
\begin{itemize}
\item rearranging the parentheses such that the $p_1$th to the $q_1$th variables, counting from the left, are arranged in left normalized form;
\item reversing the left-to-right order of these $q_1-p_1+1$ variables.
\end{itemize}
\end{itemize}
\item The next pair of canonical isomorphisms \[\nicexy{p^2_{\st}(\ua) \ar[r]^-{\alpha_2} & p_{\st}^3(\ua) \ar[r]^-{c^{(n)}_{p_2,q_2}} & p_{\st}^4(\ua)}\] are defined as in the previous pair $(\alpha_1,c^{(n)}_{p_1,q_1})$, starting with $p_{\st}^2(\ua)$  and replacing $p_1<q_1$ with $p_2<q_2$.  Inductively, this defines $\alpha_j$ and $c^{(n)}_{p_j,q_j}$ for $1\leq j \leq k$.  
\item By construction, $p_{\st}^{2k}$ and $q_{\st}$ differ only by the arrangement of the parentheses.  The isomorphism $\alpha_{k+1}$ is once again the unique canonical isomorphism involving only $\Id$ and $\alpha^{\pm 1}$ guaranteed by the Coherence Theorem for monoidal categories. 
\end{itemize}
By construction, $\sigma_*$ is a canonical isomorphism, which is natural in the objects $a_1,\ldots,a_n$.  

To show that $\sigma_*$ is well-defined, by the Coherence Theorem for monoidal categories and the coboundary monoidal category axioms, we may (i) assume that $A$ is a strict monoidal category and (ii) ignore the monoidal unit.  We need to show that, for each generating axiom in the cactus group $Cac_n$ in Definition \ref{def:cactus-group}, the two sides induce the same canonical isomorphism.
\begin{enumerate}[label=(\roman*)]
\item For the involution axiom, it is enough to show that \[c^{(k)} \circ  c^{(k)} = \Id.\] This is true by Proposition \ref{ck-involution}.
\item For the disjointness axiom, suppose $p<q$ and $r<t$ are disjoint.  Then both sides of the disjointness axiom induce the canonical isomorphism \[\Id^{\otimes p-1} \otimes c^{(q-p+1)} \otimes \Id^{\otimes r-q-1} \otimes c^{(t-r+1)} \otimes \Id^{\otimes n-t}\] if $q<r$, and similarly if $t < p$.
\item For the containment axiom, suppose $p<q$ contains $r<t$.  Then it is enough to show that 
\[\begin{split}
& c^{(q-p+1)} \circ \bigl(\Id^{\otimes r-p} \otimes c^{(t-r+1)} \otimes \Id^{\otimes q-t}\bigr)\\
& = \bigl(\Id^{\otimes q-t} \otimes c^{(t-r+1)} \otimes \Id^{\otimes r-p}\bigr) \circ c^{(q-p+1)}.
\end{split}\] 
This is trivially true if both $p=r$ and $q=t$.  On the other hand, if either $p<r$ or $t<q$, or both, then the above equality holds by Proposition \ref{ck-containment}.
\end{enumerate}
The proof that the five conditions in Theorem \ref{gmonoidal-coherence-1} with $\G=\Cac$ are satisfied is the same as in the proof of Theorem \ref{bmonoidal-braided-monoidal}.

Conversely, each $\Cac$-monoidal category $A$--i.e., an $\MCatcac$-algebra--has an underlying monoidal category by Theorem \ref{gmonoidal-coherence-1}.  We define the commutor $\xi$ as \[\nicexy@C+1cm{a_1 \otimes a_2 \ar[r]^-{\xi_{a_1,a_2}\,=\, (s^{(2)}_{1,2})_*} & a_2 \otimes a_1} \in A\]
in which \[s^{(2)}_{1,2} \in \MCatcac(2)\bigl(x_1x_2,x_2x_1\bigr)\] with $s^{(2)}_{1,2} \in Cac_2$ the generating cactus.  Now we check the coboundary monoidal category axioms in Definition \ref{def:coboundary-monoidal-category}, for which we use the axioms in the cactus group in Definition \ref{def:cactus-group}.
\begin{enumerate}[label=(\roman*)]
\item The symmetry axiom, $\xi_{Y,X}\circ \xi_{X,Y} = \Id_{X\otimes Y}$, is satisfied by the involution axiom in $Cac_2$, i.e., \[s^{(2)}_{1,2} s^{(2)}_{1,2} = \id_2.\]
\item The unit axiom is satisfied because $Cac_1$ is the trivial group.
\item The cactus axiom \eqref{cactus-axiom} is satisfied by the equality
\[s^{(2)}_{1,2} \comp_1 s^{(2)}_{1,2} = s^{(2)}_{1,2} \comp_2 s^{(2)}_{1,2} \in Cac_3,\] which follows from the following computation in the cactus group operad $\Cac$:
\[\begin{split}
s^{(2)}_{1,2} \comp_1 s^{(2)}_{1,2} 
&= \gammacac\bigl(s^{(2)}_{1,2}; s^{(2)}_{1,2}, \id_1\bigr)\\
&= s^{(2)}_{1,2}\langle 2,1\rangle \cdot \bigl(s^{(2)}_{1,2} \oplus \id_1\bigr)\\
&= \bigl(s^{(3)}_{1,3} \cdot s^{(3)}_{1,2}\bigr) \cdot s^{(3)}_{1,2}\\
&= s^{(3)}_{1,3}\\
&= \bigl(s^{(3)}_{1,3} \cdot s^{(3)}_{2,3}\bigr) \cdot s^{(3)}_{2,3}\\
&= s^{(2)}_{1,2}\langle 1,2\rangle \cdot \bigl(\id_1 \oplus s^{(2)}_{1,2}\bigr)\\
&= \gammacac\bigl(s^{(2)}_{1,2}; \id_1,s^{(2)}_{1,2}\bigr)\\
&= s^{(2)}_{1,2} \comp_2 s^{(2)}_{1,2}.
\end{split}\]
\end{enumerate}
In the above computation, from top to bottom:
\begin{itemize}
\item The first equality is the definition of $\comp_1$ in terms of the operadic composition $\gammacac$ as in \eqref{compi-def}.
\item The second equality is the definition of the operadic composition $\gammacac$ in Definition \ref{def:cactus-operad-comp}.
\item The third equality uses Lemma \ref{block-cactus-properties}(1) on the block cactus $s^{(2)}_{1,2}\langle 2,1\rangle$ and Definition \ref{def:cactus-direct-sum} on the direct sum cactus $s^{(2)}_{1,2} \oplus \id_1$.
\item The fourth equality uses the involution axiom in $Cac_3$.
\item The remaining equalities follow as in the first four equalities in the reverse order.
\end{itemize}
This shows that $(A,\xi)$ is a coboundary monoidal category.
\end{proof}

Next is the strict analogue of Theorem \ref{cacmonoidal-coboundary-monoidal}.  We will use the strict $\G$-monoidal category operad $\MCatgst$ in Definition \ref{def:mcatgst-operad} when $\G$ is the cactus group operad $\Cac$.

\begin{corollary}\label{strict-cacmonoidal-coboundary-monoidal}\index{strict coboundary monoidal category!operadic interpretation}
For each small category, the following two structures are equivalent:
\begin{enumerate}[label=(\roman*)]
\item An $\MCatcacst$-algebra structure.
\item A strict coboundary monoidal category structure.
\end{enumerate}
\end{corollary}

\begin{proof}
By Theorem \ref{mcatgst-algebra} applied to $\G=\Cac$, an $\MCatcacst$-algebra structure on $A$ is equal to an $\MCatcac$-algebra structure on $A$ whose underlying monoidal category is strict.  By Theorem \ref{cacmonoidal-coboundary-monoidal}, an $\MCatcac$-algebra structure is equal to a coboundary monoidal category structure.  Therefore, an $\MCatcacst$-algebra structure is exactly a coboundary monoidal category structure whose underlying monoidal category is strict.
\end{proof}

\begin{remark}[Gurski's Categorical Borel Construction]\label{rk:gurski-cactus}
When $\G$ is the cactus group operad $\Cac$, Gurski in \cite{gurski} (Corollary 4.14) showed that the algebras over his\index{categorical Borel construction} categorical Borel construction are small strict coboundary monoidal categories.  So we may regard the $\Cac$-monoidal category operad $\MCatcac$, whose algebras have general underlying small monoidal categories, as a generalization of Gurski's categorical Borel construction for the cactus group operad $\Cac$.\dqed
\end{remark}

\section{Coboundary Monoidal Functors are $\Cac$-Monoidal Functors}
\label{sec:coboundary-monoidal-functor}

The purpose of this section is to observe that, for the cactus group operad $\Cac$, a $\Cac$-monoidal functor in the sense of Definition \ref{def:gmonoidal-functor} is exactly a cactus monoidal functor, which is defined exactly like a braided monoidal functor.  The strong and strict analogues are also true.  As a consequence, the $\Cac$-monoidal category operad $\MCatcac$ is the symmetric operad in $\Cat$ that encodes coboundary monoidal categories and strict coboundary monoidal functors.  Moreover, the $\Cac$-monoidal functor operad $\MFuncac$ is the $2$-colored symmetric operad in $\Cat$ that encodes general coboundary monoidal functors between small coboundary monoidal categories.

\begin{definition}\label{def:coboundary-monoidal-functor}
Suppose $\M$ and $\N$ are coboundary monoidal categories.  A \index{functor!coboundary monoidal}\index{coboundary monoidal functor}\index{monoidal functor!coboundary}\emph{coboundary monoidal functor} $(F,F_2,F_0) : \M \to \N$
is a monoidal functor between the underlying monoidal categories that is compatible with the commutors, in the sense that the diagram
\begin{equation}\label{monoidal-functor-coboundary}
\nicexy@C+.5cm{F(X) \otimes F(Y) \ar[d]_-{F_2} \ar[r]^-{\xi_{FX,FY}}_-{\cong} & F(Y) \otimes F(X) \ar[d]^-{F_2} \\ F(X \otimes Y) \ar[r]^-{F\xi_{X,Y}}_-{\cong} & F(Y \otimes X)}
\end{equation}
is commutative for all objects $X,Y \in \M$.  A coboundary monoidal functor is said to be \emph{strong} (resp., \emph{strict}) if the underlying monoidal functor is so.
\end{definition}

\begin{remark}The definition of a coboundary monoidal functor is the same as that of a braided monoidal functor, except for replacing the words \emph{braided} and \emph{braidings} with the words \emph{coboundary} and \emph{commutors}.\dqed
\end{remark}

\begin{remark}
What we call a strong coboundary monoidal functor is what Gurski \cite{gurski} called a \emph{coboundary functor}.  This is a strong monoidal functor--i.e., a monoidal functor $F$ with invertible structure morphisms $F_0$ and $F_2$--such that the diagram \eqref{monoidal-functor-coboundary} is commutative.\dqed
\end{remark}

The main result of this section identifies coboundary monoidal functors with $\Cac$-monoidal functors for the cactus group operad $\Cac$.

\begin{theorem}\label{cacmonoidal-functor}\index{coboundary monoidal functor!operadic interpretation}
Suppose $F : \M \to \N$ is a functor between small coboundary monoidal categories.  Then the following two structures on $F$ are equivalent:
\begin{enumerate}[label=(\roman*)]
\item A coboundary monoidal functor structure.
\item A $\Cac$-monoidal functor structure.
\end{enumerate}
\end{theorem}

\begin{proof}
We simply recycle the proof of Theorem \ref{bmonoidal-functor} by replacing:
\begin{itemize}
\item the braid group operad $\B$ with the cactus group operad $\Cac$;
\item the generating braid $s_1\in B_2$ with the generating cactus $s^{(2)}_{1,2}\in Cac_2$;
\item the words \emph{braided} and \emph{braiding} with the words \emph{coboundary} and \emph{commutor}, respectively;
\item Theorem \ref{bmonoidal-braided-monoidal} with Theorem \ref{cacmonoidal-coboundary-monoidal}.
\end{itemize}
\end{proof}

\begin{interpretation}
Theorem \ref{cacmonoidal-functor} is a coherence result for coboundary monoidal functors.  Indeed, it asserts that for a monoidal functor $F : \M \to \N$ between small coboundary monoidal categories, in order for the diagram \eqref{gmonoidal-functor-diagram} to be commutative for each isomorphism in the $\Cac$-monoidal category operad $\MCatcac$, it is necessary and sufficient that the single diagram \eqref{monoidal-functor-coboundary} be commutative.\dqed
\end{interpretation}

Recall from Definition \ref{def:gmonoidal-functor} that a $\G$-monoidal functor is said to be \emph{strong} (resp., \emph{strict}) if the underlying monoidal functor is so.  The next observation is an immediate consequence of Theorem \ref{cacmonoidal-functor}.

\begin{corollary}\label{strong-cacmonoidal-functor}
Suppose $F : \M \to \N$ is a functor between small coboundary monoidal categories.  Then the following two structures on $F$ are equivalent:
\begin{enumerate}[label=(\roman*)]
\item A strong (resp., strict) coboundary monoidal functor structure.
\item A strong (resp., strict) $\Cac$-monoidal functor structure.
\end{enumerate}
\end{corollary}

\begin{corollary}\label{mcatcac-algebra}\index{coboundary monoidal category!operad}\index{operad!coboundary monoidal category}
$\MCatcac$ is a $1$-colored symmetric operad in $\Cat$ whose:
\begin{itemize}
\item algebras are small coboundary monoidal categories;
\item algebra morphisms are strict coboundary monoidal functors.
\end{itemize}
\end{corollary}

\begin{proof}
This follows from Theorem \ref{cacmonoidal-coboundary-monoidal}, the $\G=\Cac$ special case of Theorem \ref{mcatg-algebra-morphism}, and Corollary \ref{strong-cacmonoidal-functor}.
\end{proof}

Next is the strict analogue.

\begin{corollary}\label{mcatcacst-algebra}\index{strict coboundary monoidal category!operad}\index{operad!strict coboundary monoidal category}
$\MCatcacst$ is a $1$-colored symmetric operad in $\Cat$ whose:
\begin{itemize}
\item algebras are small strict coboundary monoidal categories;
\item algebra morphisms are strict coboundary monoidal functors.
\end{itemize}
\end{corollary}

\begin{proof}
This follows from Corollary \ref{strict-cacmonoidal-coboundary-monoidal}, the $\G=\Cac$ special case of Proposition \ref{mcatgst-algebra-morphism}, and Corollary \ref{strong-cacmonoidal-functor}.
\end{proof}

For an action operad $\G$, recall from Definition \ref{def:gmonoidal-functor-operad} the $\G$-monoidal functor operad $\MFung$.

\begin{corollary}\label{mfuncac-algebra}\index{coboundary monoidal functor!operad}\index{operad!coboundary monoidal functor}
An $\MFuncac$-algebra is exactly a triple $(A,B,F)$ such that:
\begin{itemize}
\item $A$ and $B$ are small coboundary monoidal categories.
\item $F : A \to B$ is a coboundary monoidal functor.
\end{itemize}
\end{corollary}

\begin{proof}
This follows from Theorem \ref{mfung-algebra} applied to the cactus group operad $\Cac$, Theorem \ref{cacmonoidal-coboundary-monoidal}, and Theorem \ref{cacmonoidal-functor}.
\end{proof}

Next is the strong analogue involving the strong $\G$-monoidal functor operad $\MFungstg$ in Definition \ref{def:strong-gmonoidal-functor-operad}.

\begin{corollary}\label{mfuncacstg-algebra}\index{operad!strong coboundary monoidal functor}
An $\MFuncacstg$-algebra is exactly a triple $(A,B,F)$ such that:
\begin{itemize}
\item $A$ and $B$ are small coboundary monoidal categories.
\item $F : A \to B$ is a strong coboundary monoidal functor.
\end{itemize}
\end{corollary}

\begin{proof}
This follows from Theorem \ref{mfungstg-algebra} applied to the cactus group operad $\Cac$, Theorem \ref{cacmonoidal-coboundary-monoidal}, and Corollary \ref{strong-cacmonoidal-functor}.
\end{proof}

\section{Coherence for Coboundary Monoidal Categories}
\label{sec:coboundary-monoidal-cat-coherence}

The purpose of this section is to record several coherence results for coboundary monoidal categories by restricting the coherence results for $\G$-monoidal categories in Chapter \ref{ch:coherence-gmonoidal} to the cactus group operad $\Cac$ in Definition \ref{def:braid-group-operad}. 

In the previous two sections, we proved the following results.
\begin{enumerate}
\item In Theorem \ref{cacmonoidal-coboundary-monoidal} we observed that $\Cac$-monoidal categories and small coboundary monoidal categories are the same thing.
\item The strict version is Corollary \ref{strict-cacmonoidal-coboundary-monoidal}, which says that strict $\Cac$-monoidal categories and small strict coboundary monoidal categories are the same thing.  
\item In Theorem \ref{cacmonoidal-functor} we established the equivalence between $\Cac$-monoidal functors and coboundary monoidal functors, with the strong/strict versions in Corollary \ref{strong-cacmonoidal-functor}. 
\end{enumerate}

Therefore, the coherence results in Chapter \ref{ch:coherence-gmonoidal} applies to coboundary monoidal categories and coboundary monoidal functors when $\G$ is the cactus group operad $\Cac$.  First is the $\G=\Cac$ special case of Theorem \ref{gmonoidal-coherence2}, which says that coboundary monoidal categories can always be strictified.

\begin{theorem}\label{coboundary-monoidal-strictification}\index{strictification!coboundary monoidal category}\index{coherence!coboundary monoidal category}\index{coboundary monoidal category!coherence}
Suppose $A$ is a small coboundary monoidal category.  Then there exists an adjoint equivalence \[\nicexy{A \ar@<2pt>[r]^-{L} & A_{\st} \ar@<2pt>[l]^-{R}}\] in which:
\begin{enumerate}[label=(\roman*)]
\item $A_{\st}$ is a strict coboundary monoidal category.
\item Both $L$ and $R$ are strong coboundary monoidal functors.
\end{enumerate} 
\end{theorem}

\begin{remark}
Theorem \ref{coboundary-monoidal-strictification} is essentially \cite{gurski} Proposition 4.7.\dqed
\end{remark}

Next is the $\G=\Cac$ special case of Theorem \ref{gmonoidal-coherence3}, which provides a strictification of a free coboundary monoidal category.

\begin{theorem}\label{free-coboundary-monoidal-strictification}
For each small category $A$, the strict coboundary monoidal functor \[\nicexy{\MCatcac(A) \ar[r]^-{\Lambdacac_A} & \MCatcacst(A)}\] in \eqref{lambdag-a} is an equivalence of categories, where $\MCatcac(A)$ (resp., $\MCatcacst(A)$) is the free (strict) coboundary monoidal category generated by $A$.
\end{theorem}

The following result is the $\G=\Cac$ special case of Theorem \ref{gmonoidal-coherence4}.  It says that the canonical functor from the free coboundary monoidal category generated by a set of objects to the free strict coboundary monoidal category generated by one object is faithful.

\begin{theorem}\label{cacmonoidal-coherence4}
For each set $A$ of objects, the composite functor \[\nicexy{\MCatcac(A) \ar[r]^-{\Lambdacac_A} & \MCatcacst(A) \ar[r]^-{\pi} & \MCatcacst(*)}\] is faithful.
\end{theorem}

Next is the $\G=\Cac$ special case of Corollary \ref{gmonoidal-coherence-underlying}.

\begin{corollary}\label{coboundary-monoidal-coherence-underlying}
For each set $A$ of objects, a diagram in the free coboundary monoidal category $\MCatcac(A)$ generated by $A$ is commutative if and only if composites with the same domain and the same codomain have the same $\pi\Lambdacac_A$-image in $\MCatcacst(*)$.
\end{corollary}

\begin{interpretation}
Recall from Example \ref{ex:free-strict-gmon-one-element}, restricted to the cactus group operad $\Cac$, that the free strict coboundary monoidal category generated by one object, $\MCatcacst(*)$, has objects the non-negative integers $\{0,1,2,\ldots\}$.  The morphism sets are \[\MCatcacst(*)(m;n) = \begin{cases} \varnothing & \text{if $m\not=n$};\\ Cac_n & \text{if $m=n$}.\end{cases}\] Using the explicit description of the composite $\pi\Lambdacac_A$ in the proof of Theorem \ref{gmonoidal-coherence4-alt}, it makes sense to call \[\pi\Lambdacac_A(\sigma) =\sigma \in Cac_n = \MCatcacst(*)(n;n)\] the \emph{underlying cactus} of each morphism $\sigma \in \MCatcac(A)$ for a set $A$.  Therefore, Corollary \ref{coboundary-monoidal-coherence-underlying} says that, in the free coboundary monoidal category $\MCatcac(A)$ generated by a set $A$ of objects, a diagram is commutative if and only if composites with the same (co)domain have the same underlying cactus.\dqed
\end{interpretation}

\begin{remark}[Free (Strict) Coboundary Monoidal Categories]
An explicit description of the free strict coboundary monoidal category $\MCatcacst(A)$ generated by a small category $A$ is given in Theorem \ref{free-strict-gmon-category} by restricting to the case $\G=\Cac$.  Similarly, an explicit description of the free coboundary monoidal category $\MCatcac(A)$ generated by a small category $A$ is given in Theorem \ref{free-gmon-category}.\dqed
\end{remark}